\documentclass[twoside]{article}

\usepackage{tocbibind}
\usepackage{fancyhdr,bm,mathtools,graphicx}
\usepackage{amssymb,amsmath}
\usepackage{marvosym}
\usepackage{mathstyle}
\usepackage{hieroglf}
\usepackage[T1]{fontenc}
\usepackage{eufrak}
\usepackage{romanbar}
\usepackage[dvipsnames]{xcolor}
\usepackage{ifthen}
\usepackage{tikz-cd}
\tikzcdset{
arrow style=tikz,
diagrams={>={Stealth[scale=1]}}
}
\usepackage[T1]{fontenc}
\usepackage{textcomp}

\usepackage{stmaryrd}
\usepackage[mathscr]{euscript}
\usepackage{chngcntr}
\usepackage[style]{fncychap}
\usepackage{latexsym}

\def\Ref#1{[{equation~\ref{#1}}]}

\newcounter{chapterr}
\counterwithin{section}{chapterr}

\newcounter{definition}
\counterwithin{definition}{chapterr}
\def\definition
{\refstepcounter{definition}\vskip0.4\baselineskip\noindent{{\textbf{Definition~}}}{$\bf\arabic{chapterr}.$}{$\bf\arabic{definition}$}:\hskip 0.4\parindent}
\newcounter{theorem}
\counterwithin{theorem}{chapterr}
\def\theorem
{\refstepcounter{theorem}\vskip0.4\baselineskip\noindent{{\textbf{Theorem~}}}{$\bf\arabic{chapterr}.$}{$\bf\arabic{theorem}$}:\hskip 0.4\parindent}
\newcounter{symbol}
\counterwithin{symbol}{chapterr}

\newcounter{QMaxiom}
\counterwithin{QMaxiom}{chapterr}

\newcounter{corollary}
\counterwithin{corollary}{chapterr}
\def\corollary
{\refstepcounter{corollary}\vskip0.4\baselineskip\noindent{{\textbf{Corollary~}}}{$\bf\arabic{chapterr}.$}{$\bf\arabic{corollary}$}:\hskip 0.4\parindent}
\newcounter{lemma}
\counterwithin{lemma}{chapterr}

\newcounter{proposition}
\counterwithin{proposition}{chapterr}

\newcounter{remark}
\counterwithin{remark}{chapterr}

\def\proof
{\vskip0.3\baselineskip\noindent{\textbf{proof}}\hskip.08\baselineskip:\hskip 0.4\parindent}

\newcounter{fixed}
\counterwithin{fixed}{chapterr}

\def\prooff
{\vskip0.3\baselineskip\noindent{\textbf{proof}}:\hskip 0.4\parindent}

\def\caution{\textasteriskcentered\hskip0.25\baselineskip}

\def\Ref#1{(\ref{#1})}
\def\refthm#1{($\ref{#1}$)}
\def\V{\mathbb{V}}
\def\F{{\mathcal F}}
\def\C{\mathbb{C}}

\def\R{\mathbb{R}}

\def\S{{\rm S}}

\def\U{{\rm U}}
\def\M{{\rm M}}
\def\X{{\rm X}}

\def\f{{\rm f}}
\def\h{{\rm h}}
\def\ff{{\rm F}}

\def\N{\mathbb{N}}
\def\Z{\mathbb{Z}}

\def\vector{\overrightarrow}

\newcommand{\p}[1]{{#1^{\prime}}}

\def\refthm#1{($\ref{#1}$)}
\def\F{\mathbb{F}}
\def\C{\mathbb{C}}
\def\R{\mathbb{R}}

\def\M{{\rm M}}

\def\f{{\rm f}}
\def\ff{{\rm F}}

\def\card{{\textsf{card}}}
\def\endef{{\flushright{\noindent$\blacksquare$\\}}\noindent}
\def\endthm{{\flushright{\noindent$\square$\\}}\noindent}
\def\endcor{{\flushright{\noindent$\square$\\}}\noindent}

\def\endp{}

\def\then{\Rightarrow}
\def\thenn{\Leftrightarrow}

\def\vthenn{\Updownarrow}
\def\vthen{\Downarrow}
\def\iffdef{:\Leftrightarrow}
\def\({\left(}
\def\){\right)}
\def\[{\left[}
\def\]{\right]}
\newcommand{\bpair}[2]{\negthinspace\left(#1,\thinspace #2\right)}
\newcommand{\bbpair}[2]{\negthinspace\left(#1;\thinspace #2\right)}
\newcommand{\bbsingle}[1]{\negthinspace\left({#1}\right)}

\newcommand{\CSs}[1]{{\mathcal{P}}\(#1\)}

\newcommand{\Card}[1]{\left| #1\right|}
\newcommand{\CarD}[1]{\card\(#1\)}
\def\cardeq{\overset{\underset{\mathrm{card}}{}}{=}}

\newcommand{\union}[1]{\bigcup#1}
\newcommand{\intersection}[1]{\bigcap#1}
\newcommand{\Intersection}[3]{\bigcap_{#1\in#2}#3}
\newcommand{\Cinc}[1]{{\mathsf{Inc}}_{#1}}
\newcommand{\Union}[3]{\bigcup_{#1\in#2}#3}
\newcommand{\fUnion}[4]{\bigcup_{#1=#2}^{#3}#4}
\newcommand{\disunion}[1]{\bigsqcup{#1}}
\newcommand{\Dproduct}[3]{\prod_{{#1}\in{#2}}{#3}}
\newcommand{\cmp}[2]{#1\circ#2}
\newcommand{\Func}[2]{{\textsf{F}}\opair{#1}{#2}}
\newcommand{\surFunc}[2]{{{\mathfrak{s}}\textsf{F}}\opair{#1}{#2}}
\newcommand{\IF}[2]{{\mathsf{{\mathfrak{b}}\hphantom{}F}}\opair{#1}{#2}}
\newcommand{\InF}[2]{{\mathfrak{i}}\hphantom{}{\mathsf{F}}\opair{#1}{#2}}
\newcommand{\Zfamily}[2]{{\seta{#1}}_{#2}}
\def\indexf{\imath}
\newcommand{\domain}[1]{{\mathsf{dom}}\bbsingle{#1}}
\newcommand{\codomain}[1]{{\mathsf{img}}\bbsingle{#1}}
\newcommand{\funcimage}[1]{{\mathsf{img}}\bbsingle{#1}}
\newcommand{\Injection}[2]{{\mathrm{Inj}}_{{#1}\to{#2}}}
\newcommand{\funcprod}[2]{{#1}{\overset{\underset{\mathrm{\mathfrak{f}}}{}}{\times}}{#2}}
\newcommand{\EqR}[1]{{\mathsf{EqR}}\bbsingle{#1}}
\newcommand{\EqClass}[2]{{#1}/{#2}}
\newcommand{\pEqclass}[2]{\[{#1}\]_{#2}}
\newcommand{\PEqclass}[2]{{\mathsf{EqC}}\bbpair{#1}{#2}}
\newcommand{\eqrel}[1]{{\mathscr{R}}_{#1}}
\newcommand{\Cprod}[2]{{#1}\times{#2}}
\newcommand{\psCprod}[2]{{#1}\ltimes{#2}}
\newcommand{\iCprod}[2]{{#1}\Join{#2}}
\newcommand{\Cproduct}[1]{\prod_{}{#1}}
\newcommand{\CProduct}[2]{\prod_{#1}{#2}}
\newcommand{\OR}[2]{#1\thinspace\lor\thinspace#2}
\newcommand{\AND}[2]{#1\thinspace\land\thinspace#2}
\newcommand*{\suchthat}{\;\ifnum\currentgrouptype=16 \middle\fi|\;}
\newcommand{\Foreach}[2]{\forall\thinspace#1\in#2:\thickspace}
\newcommand{\Exists}[2]{\exists\thinspace#1\in#2\thickspace|\thickspace}
\newcommand{\Existss}[2]{\exists\thinspace#1\subseteq#2\thickspace|\thickspace}
\newcommand{\Existssis}[2]{\(\exists\)\thinspace#1\subseteq#2\thickspace|\thickspace}
\newcommand{\Existsu}[2]{\exists!\thinspace#1\in#2\thickspace|\thickspace}
\newcommand{\Existsis}[2]{\(\exists\)\thinspace#1\in#2\thickspace|\thickspace}
\newcommand{\Existsuis}[2]{\(\exists!\)\thinspace#1\in#2\thickspace|\thickspace}

\newcommand{\defset}[3]{\left\{#1\in#2\thickspace:\thickspace#3\right\}}
\newcommand{\defsets}[3]{\left\{#1\subseteq#2\thickspace:\thickspace#3\right\}}
\newcommand{\defSet}[2]{\left\{#1\thickspace:\thickspace#2\right\}}
\newcommand{\negation}[1]{\neg{#1}}
\def\eqdef{\overset{\underset{\mathrm{def}}{}}{=}}
\def\indef{\overset{\underset{\mathrm{def}}{}}{\in}}
\newcommand{\resd}[1]{{\mathsf{resD}}_{#1}}
\newcommand{\rescd}[1]{{\mathsf{resC\negthinspace D}}_{#1}}
\newcommand{\res}[1]{{\mathfrak{res}}_{#1}}
\newcommand{\finv}[1]{{#1}^{-1}}
\def\Bb{{\rm B}}

\def\M{{\rm M}}

\def\index{{\mathscr{I}}}
\def\Zp{\Z^{+}}

\def\Zpz{\Z^{\geq 0}}

\def\empty{\varnothing}
\def\Xt{{\mathbb{X}}}
\def\Yt{{\mathbb{Y}}}
\newcommand{\opair}[2]{\(#1,\thinspace #2\)}

\newcommand{\topology}[1]{{\mathsf{T}}_{#1}}
\newcommand{\Foreachs}[2]{\forall\thinspace#1\subseteq#2\thickspace:\thickspace}

\newcommand{\Fclosed}[2]{{\textsf C}_{#1}\(#2\)}
\newcommand{\compl}[2]{#1\setminus#2}
\def\Tt{{\mathscr T}}
\def\Y{{\mathrm Y}}
\def\y{y}
\newcommand{\seta}[1]{\left\{#1\right\}}
\newcommand{\func}[2]{#1\(#2\)}

\newcommand{\pfunc}[3]{#1^{#3}\(#2\)}
\newcommand{\power}[2]{#1^{#2}}
\def\base{{\mathsf B}}
\newcommand{\baseof}[1]{{\mathsf B}_{#1}}
\def\subbase{{\mathcal S}}
\newcommand{\Cbase}[1]{{\mathscr B}\[#1\]}
\newcommand{\Cbases}[1]{{\widetilde{\mathscr B}}\[#1\]}
\def\x{x}
\def\const{c}
\def\U{U}
\def\V{V}
\def\W{W}
\def\B{B}
\def\S{S}
\def\Ss{{\mathsf{S}}}
\def\sC{\Gamma}
\def\sCi{\Lambda}
\def\sCii{\Omega}
\def\f{\eta}
\def\ff{\gamma}
\newcommand{\topgen}[2]{\tau_{#1}\[#2\]}
\newcommand{\basegen}[2]{\base_{#1}\[#2\]}
\newcommand{\covers}[1]{{\mathsf{COV}}\bbsingle{#1}}
\def\subtop{\overset{\underset{\mathrm{top}}{}}{\leq}}
\newcommand{\sspower}[2]{\left. #1\right|_{#2}}
\newcommand{\stopology}[2]{\left. #1\right|_{#2}}
\newcommand{\sbase}[2]{\left. #1\right|_{#2}}
\def\asubset{A}
\newcommand{\Asubset}[1]{A_{#1}}
\def\bsubset{B}

\def\csubset{C}
\newcommand{\Csubset}[1]{C_{#1}}
\def\point{p}

\newcommand{\Sint}[2]{{\mathsf{Sint}}_{#1}\[#2\]}
\newcommand{\Sext}[2]{{\mathsf{Sext}}_{#1}\[#2\]}
\newcommand{\Sbound}[2]{{\mathsf{Sbnd}}_{#1}\[#2\]}
\newcommand{\Slim}[2]{{\mathsf{Slim}}_{#1}\[#2\]}
\newcommand{\Sadh}[2]{{\mathsf{Sadh}}_{#1}\[#2\]}
\newcommand{\Siso}[2]{{\mathsf{Siso}}_{#1}\[#2\]}
\newcommand{\nei}[1]{{\mathsf{N}}_{#1}}
\newcommand{\cnei}[1]{{\mathfrak{c}}{\mathsf{N}}_{#1}}
\newcommand{\Int}[1]{{\mathsf{int}}_{#1}}
\newcommand{\Ext}[1]{{\mathsf{ext}}_{#1}}
\newcommand{\Cl}[1]{{\mathsf{Cl}}_{#1}}
\newcommand{\Fr}[1]{{\mathsf{Fr}}_{#1}}

\newcommand{\Fdense}[2]{{\mathsf{D}}_{#1}\(#2\)}
\newcommand{\Fnwdense}[2]{{\mathsf{NWD}}_{#1}\(#2\)}
\newcommand{\Ctops}[1]{{\mathsf{Top}}\bbsingle{#1}}
\newcommand{\sCtops}[2]{{\mathsf{Top_{#1}}}\bbsingle{#2}}
\newcommand{\CtopsK}[1]{{\mathsf{Top_{0}}}\bbsingle{#1}}
\newcommand{\CtopsF}[1]{{\mathsf{Top_{1}}}\bbsingle{#1}}
\newcommand{\CtopsH}[1]{{\mathsf{Top_{2}}}\bbsingle{#1}}
\newcommand{\singletonTS}[1]{\mathbb{TS}_{#1}}
\def\cf{f}
\def\cg{g}
\def\om{g}
\def\hf{h}
\def\em{e}
\newcommand{\image}[1]{#1^{\rightarrow}}
\newcommand{\pimage}[1]{#1^{\leftarrow}}
\newcommand{\CF}[2]{{\mathsf{C\negthinspace F}}\bpair{#1}{#2}}
\newcommand{\HOF}[2]{{\mathsf{H\negthinspace F}}\bpair{#1}{#2}}
\newcommand{\OM}[2]{{\mathsf{O\negthinspace M}}\bpair{#1}{#2}}
\newcommand{\CM}[2]{{\mathsf{C\negthinspace M}}\bpair{#1}{#2}}
\newcommand{\idf}[1]{{\mathbf 1}_{#1}}
\newcommand{\function}[3]{#1:\thinspace#2\to#3}
\newcommand{\incf}[2]{\iota:\thinspace{#1}\hookrightarrow{#2}}
\def\atopology{\tau}
\newcommand{\binary}[2]{#1,\thinspace#2}
\def\acover{{\mathscr C}}
\def\apartition{{\mathscr P}}
\def\covelm{S}
\def\aconnectedset{\csubset}
\newcommand{\Fcov}[1]{{\mathsf{FCOV}}\bbsingle{#1}}
\newcommand{\Ocov}[1]{{\mathfrak{O}}{\mathsf{cov}}\bbsingle{#1}}
\newcommand{\Ccov}[1]{{\mathsf{CCOV}}\bbsingle{#1}}
\newcommand{\Lfcov}[1]{{\mathsf{LfCOV}}\bbsingle{#1}}
\newcommand{\fcover}[1]{{\mathsf{cov}}_{#1}}
\newcommand{\homeomorphic}[2]{{#1}\bumpeq{#2}}
\newcommand{\EM}[2]{{\mathsf{EM}}\bpair{#1}{#2}}
\newcommand{\embedded}[2]{#1\olessthan#2}
\newcommand{\embeq}[4]{{#1}\stackrel{\mathrm{emb}}{\equiv}_{\opair{#3}{#4}}{#2}}
\newcommand{\Cpart}[1]{{\mathsf{Par}}\bbsingle{#1}}
\newcommand{\Opart}[1]{{\mathsf{OPar}}\bbsingle{#1}}
\newcommand{\Clpart}[1]{{\mathsf{CPar}}\bbsingle{#1}}
\newcommand{\Ncpart}[1]{{\mathsf{NcPar}}\bbsingle{#1}}
\newcommand{\connecteds}[1]{{\mathcal{CND}}\bbsingle{#1}}
\newcommand{\maxcon}[1]{{\mathcal{MC}}\bbsingle{#1}}
\newcommand{\mcp}[1]{{\mathsf{MCp}}_{#1}}
\newcommand{\mcR}[2]{{{{\mathscr{MC}}{\mathsf{R}}}_{\rm{#2}}}\llbracket{#1}\rrbracket}
\newcommand{\sepA}[1]{{\mathbf{T}}_{#1}}
\newcommand{\Dist}[1]{{\mathsf{Dis}}\bbsingle{#1}}
\newcommand{\InDist}[1]{{\mathsf{IDis}}\bbsingle{#1}}
\newcommand{\Sep}[1]{{\mathsf{Sep}}\bbsingle{#1}}
\newcommand{\Lib}[1]{{\mathsf{Lib}}\bbsingle{#1}}
\newcommand{\ocov}[1]{{\mathfrak{o}}{\mathsf{cov}}_{#1}}
\newcommand{\focov}[1]{{{\mathsf{f}}\mathfrak{o}}{\mathsf{cov}}_{#1}}
\newcommand{\fOcov}[1]{{\mathsf{f}}{\mathfrak{O}}{\mathsf{cov}}\bbsingle{#1}}
\newcommand{\compacts}[1]{{\mathcal{CMP}}\bbsingle{#1}}
\newcommand{\acompactset}[1]{{\mathrm{K}}_{#1}}
\newcommand{\somecompacts}[1]{{\mathcal{K}}_{#1}}

\newcommand{\alexT}[1]{{\mathfrak{alex}}_{#1}}
\newcommand{\alexTS}[1]{{#1}^{\star}}
\newcommand{\collection}[1]{{\mathcal{S}}_{#1}}
\newcommand{\aset}[1]{{\mathrm{S}}_{#1}}
\newcommand{\Lim}[2]{{\mathsf{Lim}}_{#1}\bbsingle{#2}}

\newcommand{\Limits}[4]{{\underset{\overset{\opair{#1}{#2}}{}}{{\mathsf{Limit}}}}\bbpair{#3}{#4}}
\newcommand{\linterval}[2]{[{#1},\thinspace{#2})}
\newcommand{\VV}[1]{\V_{#1}}
\newcommand{\VVS}[1]{\VS_{#1}}
\newcommand{\WW}[1]{{\mathrm{W}}_{#1}}
\newcommand{\WS}[1]{{\mathbb{W}}_{#1}}
\def\pp{\p{\point}}
\def\VS{{\mathbb{V}}}
\newcommand{\NVS}[1]{{\mathscr{V}}_{#1}}
\newcommand{\NWS}[1]{{\mathscr{W}}_{#1}}
\newcommand{\vsum}[1]{+_{#1}}
\newcommand{\spro}[1]{\cdot_{#1}}
\newcommand{\Lin}[2]{{\mathsf{L}}\bpair{#1}{#2}}
\newcommand{\VFunc}[2]{{\boldsymbol{\mathfrak{F}}}\bpair{#1}{#2}}
\newcommand{\VLin}[2]{{\mathbb{L}}\bpair{#1}{#2}}
\newcommand{\zerov}[1]{{\mathbf{0}}_{#1}}

\newcommand{\triple}[3]{\opair{#1}{\binary{#2}{#3}}}
\newcommand{\tuple}[4]{\opair{#1}{\binary{#2}{\binary{#3}{#4}}}}
\newcommand{\mtuple}[2]{\(\suc{#1}{#2}\)}

\newcommand{\mult}[2]{\suc{#1}{#2}}
\def\fsum{+_{\ff}}
\def\fpro{\times_{\ff}}
\def\FF{{\mathbb{F}}}
\newcommand{\vecspaces}[2]{{\mathbf{Vec}}\bpair{#1}{#2}}
\newcommand{\anorm}[1]{\lVert\thinspace\rVert_{#1}}
\newcommand{\norm}[2]{\left\lVert{#1}\right\rVert_{#2}}
\newcommand{\norms}[1]{{\mathsf{Nr}}\bbsingle{#1}}
\newcommand{\setnorms}[2]{{\widetilde{\mathsf{Nr}}_{#1}}\bbsingle{#2}}
\def\v{v}
\def\u{u}
\def\w{w}
\newcommand{\vv}[1]{v_{#1}}
\newcommand{\ww}[1]{w_{#1}}
\newcommand{\uu}[1]{u_{#1}}

\newcommand{\pnt}[1]{p_{#1}}
\newcommand{\abs}[1]{\left|{#1}\right|}
\newcommand{\normdist}[2]{{\mathcal{N}}\negthinspace{\mathrm{d}}^{#1}_{#2}}
\newcommand{\normtop}[2]{{\mathcal{N}}\negthinspace{\mathsf{T}}^{#1}_{#2}}
\newcommand{\metrictop}[1]{{\mathcal{M}}{\mathsf{T}}_{#1}}
\newcommand{\normedtopspace}[1]{{\mathbb{T}}\bbsingle{#1}}
\newcommand{\prodnormedtopspace}[1]{{\mathbf{P}}\negthinspace{\mathbb{T}}\bbsingle{#1}}
\newcommand{\normeqR}[1]{{\mathfrak{eq}\mathsf{Nr}}\bbsingle{#1}}
\newcommand{\ball}[1]{{\mathfrak{B}}_{#1}}
\newcommand{\mball}[2]{{\mathcal{B}}_{\opair{#1}{#2}}}
\newcommand{\cball}[1]{{\mathfrak{CB}}_{#1}}
\newcommand{\mcball}[2]{{\mathcal{CB}}_{\opair{#1}{#2}}}
\newcommand{\LinC}[2]{{\mathsf{L\negthinspace C}}\bpair{#1}{#2}}
\newcommand{\VLinC}[2]{{\mathbb{L\negthinspace C}}\bpair{#1}{#2}}
\newcommand{\VLinCnorm}[2]{\anorm{\tiny{\vector{{#1}{#2}}}}}
\newcommand{\VLinCnormf}[3]{\norm{#1}{\tiny{\vector{{#2}{#3}}}}}
\newcommand{\NVLinC}[2]{{\overline{\mathbb{L\negthinspace C}}}\bpair{#1}{#2}}
\def\lcf{T}
\newcommand{\conv}[2]{{\mathbf{Conv}}_{#1}\bbsingle{#2}}
\newcommand{\seq}[1]{s_{#1}}
\newcommand{\seqq}[1]{t_{#1}}
\newcommand{\slimit}[4]{\overset{#1}{\lim_{{#2}\rightarrow\infty}}{#3}={#4}}
\newcommand{\Con}[1]{{\mathsf{ConS}}\bbsingle{#1}}
\newcommand{\Cauchy}[1]{{\mathsf{CauchyS}}\bbsingle{#1}}
\newcommand{\dualNVS}[1]{{#1}^{\star}}

\def\NF{\overline{\FF}}
\newcommand{\Vdim}[1]{{\mathfrak{dim}}_{#1}}
\newcommand{\normsum}[2]{{#1}\dotplus{#2}}
\newcommand{\suc}[2]{{#1},\ldots,{#2}}
\newcommand{\NVSiso}[2]{{\mathsf{NVSiso}}\bpair{#1}{#2}}
\newcommand{\ml}[1]{\mu_{#1}}
\newcommand{\findex}[1]{\index_{#1}}
\newcommand{\sindex}[3]{\seta{#1}_{{#2}\in{#3}}}
\newcommand{\alphaa}[1]{\alpha_{#1}}
\newcommand{\Vprodinj}[1]{{\boldsymbol{\mathfrak{I}}}_{#1}}
\newcommand{\MLin}[2]{{\mathsf{M\negthinspace L}}\bbpair{#1}{#2}}
\newcommand{\prodnormtop}[2]{{\mathbf{P}}\negthinspace{\mathcal{N}}\negthinspace{\mathsf{T}}^{#1}_{#2}}
\newcommand{\prodtop}[1]{{\mathbf{P}}{\mathsf{T}}_{#1}}
\newcommand{\Times}[2]{{#1}\times\ldots\times{#2}}
\newcommand{\MCF}[2]{{\mathsf{M\negthinspace C\negthinspace F}}\bbpair{#1}{#2}}
\newcommand{\MLC}[2]{{\mathsf{M\negthinspace L\negthinspace C}}\bbpair{#1}{#2}}

\newcommand{\VMLC}[2]{{\mathbb{M\negthinspace L\negthinspace C}}\bpair{#1}{#2}}
\newcommand{\VMLCnorm}[2]{\anorm{\tiny{\vector{{#1}{#2}}}}}
\newcommand{\VMLCnormf}[3]{\norm{#1}{\tiny{\vector{{#2}{#3}}}}}
\newcommand{\NVMLC}[2]{{\overline{\mathbb{M\negthinspace L\negthinspace C}}}\bpair{#1}{#2}}
\newcommand{\Banachnorms}[1]{{\mathsf{BanachNr}}\bbsingle{#1}}

\newcommand{\vbase}[1]{{\boldsymbol{e}}^{#1}}

\newcommand{\vdim}[1]{d_{#1}}
\newcommand{\com}[3]{{\mathsf{com}}_{#1}\negthickspace\opair{#2}{#3}\negthinspace}

\def\Zpz{\Z^{\geq 0}}

\def\Rp{\R^{+}}

\def\Rpz{\R^{\geq 0}}

\newcommand{\ssum}[3]{\sum_{#1}^{#2}{#3}}
\newcommand{\identity}[1]{{\bf{Id}}_{#1}}
\newcommand{\indexedcollection}[1]{\mathtt{#1}}

\newcommand{\Proj}[2]{{\mathrm{pr}}^{#2}_{#1}}
\newcommand{\topprod}[2]{{#1}{\overset{\underset{\mathrm{\mathfrak{t}}}{}}{\times}}{#2}}
\newcommand{\producttop}[1]{{\mathsf{ProdTop}}\(#1\)}
\newcommand{\quotienttop}[2]{{\mathsf{QTop}}\bbpair{#1}{#2}}
\newcommand{\topq}[2]{{#1}/{#2}}
\newcommand{\metricspace}[1]{{\mathbb{M}}_{#1}}
\newcommand{\metricball}[3]{{\mathsf{Ball}}_{#1}\opair{#2}{#3}}
\newcommand{\Metrictop}[1]{{\mathsf{MetrTop}}\(#1\)}
\def\refthm#1{[{\bf theorem~\ref{#1}}]}
\def\refdef#1{[{\bf definition~\ref{#1}}]}
\def\refcor#1{[{\bf corollary~\ref{#1}}]}

\def\varfill{\dotfill}

\def\toclevel@section{1}\def\toclevel@subection{2}
\addtocontents{toc}{\string\let\string\l@section\string\l@subsection}
\def\toclevel@subsection{2}\def\toclevel@subsubection{3}
\addtocontents{toc}{\string\let\string\l@subsection\string\l@subsubsection}
\usepackage{tocloft}
\newcommand{\chapteR}[1]{\cleardoublepage
{\refstepcounter{chapterr}\vskip\baselineskip\centering{\fontsize{21}{21}\selectfont${\bf\Roman{chapterr}}$
\vskip0.6\baselineskip{\fontsize{21}{21}\selectfont{\textbf{#1}}}}
\vskip 5.9\baselineskip}
\addcontentsline{toc}{0}
{\protect\vskip0.5\baselineskip\noindent\bf{\Roman{chapterr}\hskip0.5\baselineskip#1\hspace{\fill}}}\par
\fancyhead[LO]{\ifthenelse{\value{chapterr}=0}{#1}{$\bf\Roman{chapterr}$\hskip0.7\baselineskip#1}}
}
\newcommand{\Bibliography}[1]{\vskip0.5\baselineskip\centering{\huge\bf{References}}\vskip \baselineskip
\addcontentsline{toc}{0}
{\protect\vskip0.5\baselineskip\noindent\bf{References}\hspace{\fill}}\par
\fancyhead[LO,RE]{\bf{References}#1}
}

\def\refthm#1{[{theorem~\ref{#1}}]}

\def\refdef#1{[{definition~\ref{#1}}]}
\def\refcor#1{[{corollary~\ref{#1}}]}

\newcommand{\mathleft}{\@fleqntrue\@mathmargin0pt}

\renewcommand{\sectionmark}[1]{\ifthenelse{\value{section}=0}{\markright{#1}{}}
{\markright{${\arabic{section}}$ #1}{}}}
\makeatletter
\def\cleardoublepage{\clearpage\if@twoside \ifodd\c@page\else
\hbox{}
\vspace*{\fill}
\vspace{\fill}
\thispagestyle{empty}
\newpage
\if@twocolumn\hbox{}\newpage\fi\fi\fi}
\makeatother
\fancypagestyle{plain}{
\fancyhf{}
\fancyhead{}
\fancyhead[R]{\thepage}
\fancyhead[L]{\leftmark}
}
\newcommand{\newsymp}[1]{{#1}\equiv}
\newcommand{\newsymb}[1]{{#1}\equiv}

\newcommand{\SET}[1]{{\mathscr{S}}_{#1}}
\newcommand{\FUNCTION}[1]{{{f}}_{#1}}

\newcommand{\prop}[1]{{\mathfrak{p}}\llparenthesis{#1}\rrparenthesis}
\def\Prop{\mathfrak{p}}
\newcommand{\propos}[1]{{\mathbf{\mathfrak{p}}}_{#1}}
\def\dummy{\centerdot}
\usepackage[hang,flushmargin]{footmisc}
\renewcommand{\footnoterule}{%
  \kern 20pt
  \hrule width \textwidth height 0.5pt
  \kern 5pt
}
\def\quotl{``}
\def\quotr{"}
\usepackage[inner=3.8cm,outer=3cm,bottom=3.9cm,twoside]{geometry}
\begin{document}
\thispagestyle{empty}
%\rule[0pt]{1.5pt}{200pt}
%\section{}
\noindent
{\\ \\\textbf{\fontsize{40}{40}\selectfont
{\textsf{Point-Set Topology}}}}
%\\[0.8\baselineskip]{\textbf{\fontsize{20}{20}\selectfont
%{\textsf{of Finite Rank}}}}
\\[8\baselineskip]
\noindent
{\fontsize{21}{21}\selectfont
{\textsf{Farzad Shahi}}
}
\\[11\baselineskip]
\noindent
{\fontsize{11}{11}\selectfont
{\textrm{Version:} 1.00}
}
\vfill\hfill
$\underline{\Huge{\textsf{\bf F}}\negthickspace\negthickspace\negthinspace{{\rotatebox{90}{\textsf{\bf S}}}}}$
%{{\Huge{\Bat}}}
%\setlength{\unitlength}{2mm}
%\begin{picture}(30,20)
%\linethickness{0.075mm}
%\multiput(-100,0)(1,0){1}%
%{\line(0,1){20}}
%\end{picture}
%%%%%%%%%%%%%%%%%%%%%%%%%%%%%%%%%%%%%%%%%
\newpage
\thispagestyle{empty}
\noindent
{\fontsize{9.4}{9.4}\selectfont
{\underline{{\bf\textsf{Title:}} Point-Set Topology}}}\\
{\fontsize{9.4}{9.4}\selectfont
{\underline{{\bf\textsf{Author:}} Farzad Shahi}}}\\
{{\bf\textsf{email}}}:~{\texttt{shahi.farzad@gmail.com}}
\vskip 4\baselineskip
\noindent
{\fontsize{9.4}{9.4}\selectfont
{\underline{{\bf\textsf{Version:}} 1.00}}
}\\
\noindent
{\fontsize{9.4}{9.4}\selectfont
{\bf\textsf{2022}}}
\vskip 4\baselineskip
\noindent
{\fontsize{9.4}{9.4}\selectfont
{\bf\textsf{Typesetting:}} By the author, using \TeX}
\vskip 5\baselineskip
\noindent
{\fontsize{9.4}{9.4}\selectfont
{\bf\textsf{Abstract:}}
This is a review of the fundamental concepts of general topology.}
%%%%%%%%%%%%%%%%%%%%%%%%%%%%%%%%%%%%%%%%%%
%\tableofcontents
\newpage
\thispagestyle{empty}
%\Thecontents{}
\section*{\fontsize{21}{21}\selectfont\bf{Contents}}
\addtocontents{toc}{\protect\setcounter{tocdepth}{-1}}
\tableofcontents
\addtocontents{toc}{\protect\setcounter{tocdepth}{3}}
%\tableofcontents
\newpage
%%%%%%%%%%%%%%%%%%%%%%%%%%%%%%%%%%%%%%%%%%%%%%%%%%%%%%%%%%
\chapteR{
Mathematical Symbols
}
\thispagestyle{fancy}
\section*{
Set-theory
}
%%%%%%%%%%%%%%%%%%%%%%%%%%%%%%%%%%%%%%%%%%%%%%%%%%%%%%%%%%%%%%%%%%%%%%%%%%%%%%%%%%%%%%%%
$\newsymb{\empty}$
empty-set
\varfill
$\empty$\\
%%%%%%%%%%%%%%%%%%%%%%%%%%%%%%%%%%%%%%%%%%%%%%%%%%%%%%%%%%%%%%%%%%%%%%%%%%%%%%%%%%%%%%%%
$\newsymp{\SET{1}=\SET{2}}$
$\SET{1}$
equals
$\SET{2}$.
\varfill
$=$\\
%%%%%%%%%%%%%%%%%%%%%%%%%%%%%%%%%%%%%%%%%%%%%%%%%%%%%%%%%%%%%%%%%%%%%%%%%%%%%%%%%%%%%%%%
%%%%%%%%%%%%%%%%%%%%%%%%%%%%%%%%%%%%%%%%%%%%%%%%%%%%%%%%%%%%%%%%%%%%%%%%%%%%%%%%%%%%%%%%
$\newsymp{\SET{1}\in\SET{2}}$
$\SET{1}$
is an element of
$\SET{2}$.
\varfill
$\in$\\
%%%%%%%%%%%%%%%%%%%%%%%%%%%%%%%%%%%%%%%%%%%%%%%%%%%%%%%%%%%%%%%%%%%%%%%%%%%%%%%%%%%%%%%%
$\newsymp{\SET{1}\ni\SET{2}}$
$\SET{1}$
contains
$\SET{2}$.
\varfill
$\ni$\\
%%%%%%%%%%%%%%%%%%%%%%%%%%%%%%%%%%%%%%%%%%%%%%%%%%%%%%%%%%%%%%%%%%%%%%%%%%%%%%%%%%%%%%%%
$\newsymb{\seta{\binary{\SET{1}}{\SET{2}}}}$
the set composed of
$\SET{1}$
and
$\SET{2}$
\varfill
$\seta{\binary{\dummy}{\dummy}}$\\
%%%%%%%%%%%%%%%%%%%%%%%%%%%%%%%%%%%%%%%%%%%%%%%%%%%%%%%%%%%%%%%%%%%%%%%%%%%%%%%%%%%%%%%%
$\newsymb{\defset{\SET{1}}{\SET{2}}{\prop{\SET{1}}}}$
all elements of
$\SET{2}$
having the property
$\Prop$
\varfill
$\defset{\dummy}{\dummy}{\dummy}$\\
%%%%%%%%%%%%%%%%%%%%%%%%%%%%%%%%%%%%%%%%%%%%%%%%%%%%%%%%%%%%%%%%%%%%%%%%%%%%%%%%%%%%%%%%
$\newsymb{\union{\SET{}}}$
union of all elements of
$\SET{}$
\varfill
$\bigcup$\\
%%%%%%%%%%%%%%%%%%%%%%%%%%%%%%%%%%%%%%%%%%%%%%%%%%%%%%%%%%%%%%%%%%%%%%%%%%%%%%%%%%%%%%%%
$\newsymb{\intersection{\SET{}}}$
intersection of all elements of
$\SET{}$
\varfill
$\bigcap$\\
%%%%%%%%%%%%%%%%%%%%%%%%%%%%%%%%%%%%%%%%%%%%%%%%%%%%%%%%%%%%%%%%%%%%%%%%%%%%%%%%%%%%%%%%
$\newsymb{\CSs{\SET{}}}$
power-set of
$\SET{}$
\varfill
$\CSs{\dummy}$\\
%%%%%%%%%%%%%%%%%%%%%%%%%%%%%%%%%%%%%%%%%%%%%%%%%%%%%%%%%%%%%%%%%%%%%%%%%%%%%%%%%%%%%%%%
$\newsymb{\Dproduct{\alpha}{\index}{\SET{\alpha}}}$
Cartesian-product of the collection of indexed sets
${\seta{\SET{\alpha}}}_{\alpha\in\index}$
\varfill
$\prod$\\
%%%%%%%%%%%%%%%%%%%%%%%%%%%%%%%%%%%%%%%%%%%%%%%%%%%%%%%%%%%%%%%%%%%%%%%%%%%%%%%%%%%%%%%%
$\newsymp{\SET{1}\subseteq\SET{2}}$
$\SET{1}$
is a subset of
$\SET{2}$.
\varfill
$\subseteq$\\
%%%%%%%%%%%%%%%%%%%%%%%%%%%%%%%%%%%%%%%%%%%%%%%%%%%%%%%%%%%%%%%%%%%%%%%%%%%%%%%%%%%%%%%%
$\newsymp{\SET{1}\supseteq\SET{2}}$
$\SET{1}$
includes
$\SET{2}$.
\varfill
$\supseteq$\\
%%%%%%%%%%%%%%%%%%%%%%%%%%%%%%%%%%%%%%%%%%%%%%%%%%%%%%%%%%%%%%%%%%%%%%%%%%%%%%%%%%%%%%%%
$\newsymp{\SET{1}\subset\SET{2}}$
$\SET{1}$
is a proper subset of
$\SET{2}$.
\varfill
$\subset$\\
%%%%%%%%%%%%%%%%%%%%%%%%%%%%%%%%%%%%%%%%%%%%%%%%%%%%%%%%%%%%%%%%%%%%%%%%%%%%%%%%%%%%%%%%
$\newsymp{\SET{1}\supset\SET{2}}$
$\SET{1}$
properly includes
$\SET{2}$.
\varfill
$\supset$\\
%%%%%%%%%%%%%%%%%%%%%%%%%%%%%%%%%%%%%%%%%%%%%%%%%%%%%%%%%%%%%%%%%%%%%%%%%%%%%%%%%%%%%%%%
$\newsymb{\SET{1}\cup\SET{2}}$
union of
$\SET{1}$
and
$\SET{2}$
\varfill
$\cup$\\
%%%%%%%%%%%%%%%%%%%%%%%%%%%%%%%%%%%%%%%%%%%%%%%%%%%%%%%%%%%%%%%%%%%%%%%%%%%%%%%%%%%%%%%%
$\newsymb{\SET{1}\cap\SET{2}}$
intersection of
$\SET{1}$
and
$\SET{2}$
\varfill
$\cap$\\
%%%%%%%%%%%%%%%%%%%%%%%%%%%%%%%%%%%%%%%%%%%%%%%%%%%%%%%%%%%%%%%%%%%%%%%%%%%%%%%%%%%%%%%%
$\newsymb{\SET{1}\times\SET{2}}$
Cartesian-product of
$\SET{1}$
and
$\SET{2}$
\varfill
$\times$\\
%%%%%%%%%%%%%%%%%%%%%%%%%%%%%%%%%%%%%%%%%%%%%%%%%%%%%%%%%%%%%%%%%%%%%%%%%%%%%%%%%%%%%%%%
$\newsymb{\compl{\SET{1}}{\SET{2}}}$
the relative complement of
$\SET{2}$
with respect to
$\SET{1}$
\varfill
$\setminus$\\
%%%%%%%%%%%%%%%%%%%%%%%%%%%%%%%%%%%%%%%%%%%%%%%%%%%%%%%%%%%%%%%%%%%%%%%%%%%%%%%%%%%%%%%%
$\newsymb{\func{\FUNCTION{}}{\SET{}}}$
value of the function
$\FUNCTION{}$
at
$\SET{}$
\varfill
$\func{\dummy}{\dummy}$\\
%%%%%%%%%%%%%%%%%%%%%%%%%%%%%%%%%%%%%%%%%%%%%%%%%%%%%%%%%%%%%%%%%%%%%%%%%%%%%%%%%%%%%%%%
$\newsymb{\domain{\FUNCTION{}}}$
domain of the function
$\FUNCTION{}$
\varfill
$\domain{\dummy}$\\
%%%%%%%%%%%%%%%%%%%%%%%%%%%%%%%%%%%%%%%%%%%%%%%%%%%%%%%%%%%%%%%%%%%%%%%%%%%%%%%%%%%%%%%%
$\newsymb{\codomain{\FUNCTION{}}}$
codomain of the function
$\FUNCTION{}$
\varfill
$\codomain{\dummy}$\\
%%%%%%%%%%%%%%%%%%%%%%%%%%%%%%%%%%%%%%%%%%%%%%%%%%%%%%%%%%%%%%%%%%%%%%%%%%%%%%%%%%%%%%%%
$\newsymb{\funcimage{\FUNCTION{}}}$
image of the function
$\FUNCTION{}$
\varfill
$\funcimage{\dummy}$\\
%%%%%%%%%%%%%%%%%%%%%%%%%%%%%%%%%%%%%%%%%%%%%%%%%%%%%%%%%%%%%%%%%%%%%%%%%%%%%%%%%%%%%%%%
$\newsymb{\image{\FUNCTION{}}}$
image-map of the function
$\FUNCTION{}$
\varfill
$\image{\dummy}$\\
%%%%%%%%%%%%%%%%%%%%%%%%%%%%%%%%%%%%%%%%%%%%%%%%%%%%%%%%%%%%%%%%%%%%%%%%%%%%%%%%%%%%%%%%
$\newsymb{\pimage{\FUNCTION{}}}$
inverse-image-map of the function
$\FUNCTION{}$
\varfill
$\pimage{\dummy}$\\
%%%%%%%%%%%%%%%%%%%%%%%%%%%%%%%%%%%%%%%%%%%%%%%%%%%%%%%%%%%%%%%%%%%%%%%%%%%%%%%%%%%%%%%%
$\newsymb{\resd{\FUNCTION{}}}$
domain-restriction-map of the function
$\FUNCTION{}$
\varfill
$\resd{\dummy}$\\
%%%%%%%%%%%%%%%%%%%%%%%%%%%%%%%%%%%%%%%%%%%%%%%%%%%%%%%%%%%%%%%%%%%%%%%%%%%%%%%%%%%%%%%%
$\newsymb{\rescd{\FUNCTION{}}}$
codomain-restriction-map of the function
$\FUNCTION{}$
\varfill
$\rescd{\dummy}$\\
%%%%%%%%%%%%%%%%%%%%%%%%%%%%%%%%%%%%%%%%%%%%%%%%%%%%%%%%%%%%%%%%%%%%%%%%%%%%%%%%%%%%%%%%
$\newsymb{\func{\res{\FUNCTION{}}}{\SET{}}}$
domain-restriction and codomain-restriction of\\ the function
$\FUNCTION{}$ to $\SET{}$ and $\func{\image{\FUNCTION{}}}{\SET{}}$, respectively
\varfill
$\res{\dummy}$\\
%%%%%%%%%%%%%%%%%%%%%%%%%%%%%%%%%%%%%%%%%%%%%%%%%%%%%%%%%%%%%%%%%%%%%%%%%%%%%%%%%%%%%%%%
$\newsymb{\Func{\SET{1}}{\SET{2}}}$
the set of all maps from
$\SET{1}$
to
$\SET{2}$
\varfill
$\Func{\dummy}{\dummy}$\\
%%%%%%%%%%%%%%%%%%%%%%%%%%%%%%%%%%%%%%%%%%%%%%%%%%%%%%%%%%%%%%%%%%%%%%%%%%%%%%%%%%%%%%%%
$\newsymb{\IF{\SET{1}}{\SET{2}}}$
the set of all bijective functions from
$\SET{1}$
to
$\SET{2}$
\varfill
$\IF{\dummy}{\dummy}$\\
%%%%%%%%%%%%%%%%%%%%%%%%%%%%%%%%%%%%%%%%%%%%%%%%%%%%%%%%%%%%%%%%%%%%%%%%%%%%%%%%%%%%%%%%
$\newsymb{\finv{\FUNCTION{}}}$
the inverse mapping of the bijective function $\FUNCTION{}$
\varfill
$\finv{\dummy}$\\
%%%%%%%%%%%%%%%%%%%%%%%%%%%%%%%%%%%%%%%%%%%%%%%%%%%%%%%%%%%%%%%%%%%%%%%%%%%%%%%%%%%%%%%%
$\newsymb{\surFunc{\SET{1}}{\SET{2}}}$
the set of all surjective functions from
$\SET{1}$
to
$\SET{2}$
\varfill
$\surFunc{\dummy}{\dummy}$\\
%%%%%%%%%%%%%%%%%%%%%%%%%%%%%%%%%%%%%%%%%%%%%%%%%%%%%%%%%%%%%%%%%%%%%%%%%%%%%%%%%%%%%%%%
$\newsymb{\cmp{\FUNCTION{1}}{\FUNCTION{2}}}$
composition of the function
$\FUNCTION{1}$
with the function
$\FUNCTION{2}$
\varfill
$\cmp{}{}$\\
%%%%%%%%%%%%%%%%%%%%%%%%%%%%%%%%%%%%%%%%%%%%%%%%%%%%%%%%%%%%%%%%%%%%%%%%%%%%%%%%%%%%%%%%
$\newsymb{\Injection{\SET{1}}{\SET{2}}}$
the injection-mapping of the set $\SET{1}$ into the set $\SET{2}$
\varfill
$\Injection{\dummy}{\dummy}$\\
%%%%%%%%%%%%%%%%%%%%%%%%%%%%%%%%%%%%%%%%%%%%%%%%%%%%%%%%%%%%%%%%%%%%%%%%%%%%%%%%%%%%%%%%
$\newsymb{\funcprod{\cf_1}{\cf_2}}$
the function-product of the function $\cf_1$ and $\cf_2$
\varfill
$\Cprod{\dummy}{\dummy}$\\
%%%%%%%%%%%%%%%%%%%%%%%%%%%%%%%%%%%%%%%%%%%%%%%%%%%%%%%%%%%%%%%%%%%%%%%%%%%%%%%%%%%%%%%%
$\newsymb{\EqR{\SET{}}}$
the set of all equivalence relations on the set
$\SET{}$
\varfill
$\EqR{\dummy}$\\
%%%%%%%%%%%%%%%%%%%%%%%%%%%%%%%%%%%%%%%%%%%%%%%%%%%%%%%%%%%%%%%%%%%%%%%%%%%%%%%%%%%%%%%%
$\newsymb{\EqClass{\SET{1}}{\SET{2}}}$
quotient-set of
$\SET{1}$
by the equivalence-relation
$\SET{1}$
\varfill
$\EqClass{\dummy}{\dummy}$\\
%%%%%%%%%%%%%%%%%%%%%%%%%%%%%%%%%%%%%%%%%%%%%%%%%%%%%%%%%%%%%%%%%%%%%%%%%%%%%%%%%%%%%%%%
$\newsymb{\pEqclass{\SET{1}}{\SET{2}}}$
equivalence-class of
$\SET{1}$
by the equivalence-relation
$\SET{2}$
\varfill
$\pEqclass{\dummy}{\dummy}$\\
%%%%%%%%%%%%%%%%%%%%%%%%%%%%%%%%%%%%%%%%%%%%%%%%%%%%%%%%%%%%%%%%%%%%%%%%%%%%%%%%%%%%%%%%
$\newsymb{\PEqclass{\SET{1}}{\SET{2}}}$
equivalence-class of
$\SET{2}$
by the equivalence-relation
$\SET{1}$
\varfill
$\PEqclass{\dummy}{\dummy}$\\
%%%%%%%%%%%%%%%%%%%%%%%%%%%%%%%%%%%%%%%%%%%%%%%%%%%%%%%%%%%%%%%%%%%%%%%%%%%%%%%%%%%%%%%%
$\newsymp{\Card{\SET{1}}\cardeq\Card{\SET{2}}}$
$\IF{\SET{1}}{\SET{2}}$
is non-empty.
\varfill
$\Card{\dummy}\cardeq\Card{\dummy}$\\
%%%%%%%%%%%%%%%%%%%%%%%%%%%%%%%%%%%%%%%%%%%%%%%%%%%%%%%%%%%%%%%%%%%%%%%%%%%%%%%%%%%%%%%%
$\newsymb{\CarD{\SET{}}}$
cardinality of
$\SET{}$
\varfill
$\CarD{\dummy}$
%%%%%%%%%%%%%%%%%%%%%%%%%%%%%%%%%%%%%%%%%%%%%%%%%%%%%%%%%%%%%%%%%%%%%%%%%%%%%%%%%%%%%%%%
%%%%%%%%%%%%%%%%%%%%%%%%%%%%%%%%%%%%%%%%%%%%%%%%%%%%%%%%%%%%%%%%%%%%%%%%%%%%%%%%%%%%%%%%
%%%%%%%%%%%%%%%%%%%%%%%%%%%%%%%%%%%%%%%%%%%%%%%%%%%%%%%%%%%%%%%%%%%%%%%%%%%%%%%%%%%%%%%%
%%%%%%%%%%%%%%%%%%%%%%%%%%%%%%%%%%%%%%%%%%%%%%%%%%%%%%%%%%%%%%%%%%%%%%%%%%%%%%%%%%%%%%%%
%%%%%%%%%%%%%%%%%%%%%%%%%%%%%%%%%%%%%%%%%%%%%%%%%%%%%%%%%%%%%%%%%%%%%%%%%%%%%%%%%%%%%%%%
\section*{
Logic
}
%%%%%%%%%%%%%%%%%%%%%%%%%%%%%%%%%%%%%%%%%%%%%%%%%%%%%%%%%%%%%%%%%%%%%%%%%%%%%%%%%%%%%%%%
$\newsymp{\AND{\propos{1}}{\propos{2}}}$
$\propos{1}$
and
$\propos{2}$.
\varfill
$\AND{}{}$\\
%%%%%%%%%%%%%%%%%%%%%%%%%%%%%%%%%%%%%%%%%%%%%%%%%%%%%%%%%%%%%%%%%%%%%%%%%%%%%%%%%%%%%%%%
$\newsymp{\OR{\propos{1}}{\propos{2}}}$
$\propos{1}$
or
$\propos{2}$.
\varfill
$\OR{}{}$\\
%%%%%%%%%%%%%%%%%%%%%%%%%%%%%%%%%%%%%%%%%%%%%%%%%%%%%%%%%%%%%%%%%%%%%%%%%%%%%%%%%%%%%%%%
$\newsymp{{\propos{1}}\then{\propos{2}}}$
if
$\propos{1}$,
then
$\propos{2}$.
\varfill
$\then$\\
%%%%%%%%%%%%%%%%%%%%%%%%%%%%%%%%%%%%%%%%%%%%%%%%%%%%%%%%%%%%%%%%%%%%%%%%%%%%%%%%%%%%%%%%
$\newsymp{{\propos{1}}\thenn{\propos{2}}}$
$\propos{1}$,
if-and-only-if
$\propos{2}$.
\varfill
$\thenn$\\
%%%%%%%%%%%%%%%%%%%%%%%%%%%%%%%%%%%%%%%%%%%%%%%%%%%%%%%%%%%%%%%%%%%%%%%%%%%%%%%%%%%%%%%%
$\newsymp{\negation{\propos{}}}$
$\propos{1}$,
negation of
$\propos{}$.
\varfill
$\negation{}$\\
%%%%%%%%%%%%%%%%%%%%%%%%%%%%%%%%%%%%%%%%%%%%%%%%%%%%%%%%%%%%%%%%%%%%%%%%%%%%%%%%%%%%%%%%
$\newsymp{\Foreach{\SET{1}}{\SET{2}}{\prop{\SET{1}}}}$
for every
$\SET{2}$
in
$\SET{1}$,
$\prop{\SET{1}}$.
\varfill
$\Foreach{\dummy}{\dummy}\dummy$\\
%%%%%%%%%%%%%%%%%%%%%%%%%%%%%%%%%%%%%%%%%%%%%%%%%%%%%%%%%%%%%%%%%%%%%%%%%%%%%%%%%%%%%%%%
$\newsymp{\Exists{\SET{1}}{\SET{2}}{\prop{\SET{1}}}}$
exists
$\SET{2}$
in
$\SET{1}$ such that
$\prop{\SET{1}}$.
\varfill
$\Foreach{\dummy}{\dummy}\dummy$\\
%%%%%%%%%%%%%%%%%%%%%%%%%%%%%%%%%%%%%%%%%%%%%%%%%%%%%%%%%%%%%%%%%%%%%%%%%%%%%%%%%%%%%%%%
%%%%%%%%%%%%%%%%%%%%%%%%%%%%%%%%%%%%%%%%%%%%%%%%%%%%%%%%%%%%%%%%%%%%%%%%%%%%%%%%%%%%%%%%
%%%%%%%%%%%%%%%%%%%%%%%%%%%%%%%%%%%%%%%%%%%%%%%%%%%%%%%%%%%%%%%%%%%%%%%%%%%%%%%%%%%%%%%%
%%%%%%%%%%%%%%%%%%%%%%%%%%%%%%%%%%%%%%%%%%%%%%%%%%%%%%%%%%%%%%%%%%%%%%%%%%%%%%%%%%%%%%%%
%%%%%%%%%%%%%%%%%%%%%%%%%%%%%%%%%%%%%%%%%%%%%%%%%%%%%%%%%%%%%%%%%%%%%%%%%%%%%%%%%%%%%%%%
%%%%%%%%%%%%%%%%%%%%%%%%%%%%%%%%%%%%%%%%%%%%%%%%%%%%%%%%%%%%%%%%%%%%%%%%%%%%%%%%%%%%%%%%
%%%%%%%%%%%%%%%%%%%%%%%%%%%%%%%%%%%%%%%%%%%%%%%%%%%%%%%%%%%%%%%%%%%%%%%%%%%%%%%%%%%%%%%%
\section*{Mathematical Environments}
$\newsymp{\blacksquare}$
end of definition
\varfill
$\blacksquare$\\
%%%%%%%%%%%%%%%%%%%%%%%%%%%%%%%%%%%%%%%%%%%%%%%%%%%%%%%%%%%%%%%%%%%%%%%%%%%%%%%%%%%%%%%%
$\newsymp{\square}$
end of theorem, lemma, proposition, or corollary
\varfill
$\square$\\
%%%%%%%%%%%%%%%%%%%%%%%%%%%%%%%%%%%%%%%%%%%%%%%%%%%%%%%%%%%%%%%%%%%%%%%%%%%%%%%%%%%%%%%%
$\newsymp{\Diamond}$
end of the introduction of new fixed objects
\varfill
$\Diamond$
%%%%%%%%%%%%%%%%%%%%%%%%%%%%%%%%%%%%%%%%%%%%%%%%%%%%%%%%%%%%%%%%%%%%%%%%%%%%%%%%%%%%%%%%
%%%%%%%%%%%%%%%%%%%%%%%%%%%%%%%%%%%%%%%%%%%%%%%%%%%%%%%%%%%%%%%%%%%%%%%%%%%%%%%%%%%%%%%%%%%%%%%%%%%%%%%%%%%%%%%%%%%%%%%%%%%%%%%%%%%%%%%%%%%%%%%%%%%%%%%%%%%%%%%%%%%%%%%%%%%%%%%%
%%%%%%%%%%%%%%%%%%%%%%%%%%%%%%%%%%%%%%%%%%%%%%%%%%%%%%%%%%%%%%%%%%%%%%%%%%%%%%%%%%%%%%%%%%%%%%%%%%%%%%%%%%%%%%%%%%%%%%%%%%%%%%%%%%%%%%%%%%%%%%%%%%%%%%%%%%%%%%%%%%%%%%%%%%%%%%%%
%%%%%%%%%%%%%%%%%%%%%%%%%%%%%%%%%%%%%%%%%%%%%%%%%%%%%%%%%%%%%%%%%%%%%%%%%%%%%%%%%%%%%%%%%%%%%%%%%%%%%%%%%%%%%%%%%%%%%%%%%%%%%%%%%%%%%%%%%%%%%%%%%%%%%%%%%%%%%%%%%%%%%%%%%%%%%%%%
%%%%%%%%%%%%%%%%%%%%%%%%%%%%%%%%%%%%%%%%%%%%%%%%%%%%%%%%%%%%%%%%%%%%%%%%%%%%%%%%%%%%%%%%%%%%%%%%%%%%%%%%%%%%%%%%%%%%%%%%%%%%%%%%%%%%%%%%%%%%%%%%%%%%%%%%%%%%%%%%%%%%%%%%%%%%%%%%
\chapteR{
Topological Spaces
}
\thispagestyle{fancy}
\section{
Basic Structure of a Topological Space
}
\definition\label{deftopologicalspace}
$\X$
is taken as a set, and
$\topology{}$
as an element of
$\CSs{\CSs{\X}}$
(a subset of $\CSs{\X}$).
$\Xt=\opair{\X}{\topology{}}$
is defined to be a $\quotl$topological-space$\quotr$, iff it possesses these properties.
\begin{itemize}
\item[${\textbf{\textsf{T1}}}$]
%\hfill
$\empty\in\topology{}$.
\item[${\textbf{\textsf{T2}}}$]
%\hfill
$\X\in\topology{}$.
\item[${\textbf{\textsf{T3}}}$]
%\hfill
$\Foreach{\opair{\U_{1}}{\U_{2}}}{\(\topology{}\times\topology{}\)}\U_{1}\cap\U_{2}\in\topology{}$.
\item[${\textbf{\textsf{T4}}}$]
%\hfill
$\displaystyle
\Foreachs{\sC}{\topology{}}
\(\union{\sC}\)\in\topology{}$.
\end{itemize}
\endef
%%%%%%%%%%%%%%%%%%%%%%%%%%%%%%%%%%%%%%%%%%%%%%%%%%%%%%%%%%%%%%%%%%%%%%%%%%%%%%%%%%%%%%%%%%%%%%%%%%%%%%%%%%%%%%%%%%%%%%%%%%%%%%%%%%%%
\definition
$\Xt=\opair{\X}{\topology{}}$
is taken as a topological-space, and
$\U$
as a subset of
$\X$.
\begin{itemize}
\item
$\topology{}$
is called a $\quotl$topology on $\X$$\quotr$.
\item
Every
$\point$
in
$\X$,
is called a $\quotl$point of $\Xt$$\quotr$.
\item
$\U$
is called an $\quotl$open set of $\Xt$$\quotr$, iff
$\U\in\topology{}$.
\item
$\U$
is called a $\quotl$closed set of $\Xt$$\quotr$, iff
$\(\X\setminus\U\)\in\topology{}$.
\item
$\U$
is called a $\quotl$clopen set of $\Xt$$\quotr$, iff
$\U\in\topology{}$,
and
$\(\X\setminus\U\)\in\topology{}$.
\item
For every
$\asubset$
in
$\CSs{\X}$,
$\U$
is called a $\quotl$neighbourhood of $\asubset$ in $\Xt$$\quotr$, iff
$\U\in\topology{}$,
and
$\asubset\subseteq\U$.
\end{itemize}
\endef
%%%%%%%%%%%%%%%%%%%%%%%%%%%%%%%%%%%%%%%%%%%%%%%%%%%%%%%%%%%%%%%%%%%%%%%%%%
\definition\label{defclassofalltopologies}
$\X$
is taken as a set.
The set of all topologies on $\X$ is denoted by $\Ctops{\X}$.
\begin{equation}
\Ctops{\X}:=\defset{\topology{}}{\CSs{\CSs{\X}}}
{\[\opair{\X}{\topology{}}~is~a~topological~space.\]}.
\end{equation}
\endef
%%%%%%%%%%%%%%%%%%%%%%%%%%%%%%%%%%%%%%%%%%%%%%%%%%%%%%%%%%%%%%%%%%%%%%%%%%
\definition\label{defnbdclassofsets}
$\Xt=\opair{\X}{\topology{}}$ is taken as a topological-space
$\nei{\Xt}$ is defined to be this operator.
\begin{itemize}
\item[${\textbf{\textsf{Nei1}}}$]
\hfill
$\nei{\Xt}\indef\Func{\CSs{\X}}{\CSs{\CSs{\X}}}.$
\item[${\textbf{\textsf{Nei2}}}$]
\hfill
$\Foreach{\asubset}{\CSs{\X}}
\func{\nei{\Xt}}{\asubset}\eqdef\defset{\U}{\topology{}}{\asubset\subseteq\U}.$
\end{itemize}
In other words,
$\nei{\Xt}$
is defined to be the operator
from
$\CSs{\X}$
to
$\CSs{\CSs{\X}}$ that for every
$\asubset$
in
$\CSs{\X}$,
$\func{\nei{\Xt}}{\asubset}$
is the set of all neighbourhoods of
$\asubset$
in
$\Xt$.
\endef
%%%%%%%%%%%%%%%%%%%%%%%%%%%%%%%%%%%%%%%%%%%%%%%%%%%%%%%%%%%%
%%%%%%%%%%%%%%%%%%%%%%%%%%%%%%%%%%%%%%%%%%%%%%%%%%%%%%%%%%%%%%%%%%%%%%%%%%
\definition\label{deffamilyofclosedsets}
$\Xt=\opair{\X}{\topology{}}$ is taken as a topological-space.
The set of aa closed sets of $\Xt$
$\Xt$
is denoted by
$\Fclosed{\X}{\topology{}}$.
\begin{equation}
\Fclosed{\X}{\topology{}}:=\defset{\V}{\CSs{\X}}{\[\(\X\setminus\V\)\in\topology{}\]}.
\end{equation}
\endef
%%%%%%%%%%%%%%%%%%%%%%%%%%%%%%%%%%%%%%%%%%%%%%%%%%%%%%%%%%%%%%%%%%%%%%%%%%
\theorem\label{thmclosedsets}
$\Xt=\opair{\X}{\topology{}}$ is taken as a topological-space.
\begin{itemize}
\item[${\textbf{\textsf{t1}}}$]
%\hfill
$\Fclosed{\X}{\topology{}}\subseteq\CSs{\X}.$
\item[${\textbf{\textsf{t2}}}$]
%\hfill
$\empty\in\Fclosed{\X}{\topology{}}.$
\item[${\textbf{\textsf{t3}}}$]
%\hfill
$\X\in\Fclosed{\X}{\topology{}}.$
\item[${\textbf{\textsf{t4}}}$]
%\hfill
$\Foreach{\opair{\U_{1}}{\U_{2}}}{\Fclosed{\X}{\topology{}}\times\Fclosed{\X}{\topology{}}}
\(\U_{1}\cup\U_{2}\)\in\Fclosed{\X}{\topology{}}.$
\item[${\textbf{\textsf{t5}}}$]
%\hfill
$\displaystyle
\Foreachs{\sC}{\Fclosed{\X}{\topology{}}}\(\bigcap_{\U\in\sC}\U\)\in\Fclosed{\X}{\topology{}}.$
\end{itemize}
\prooff
\begin{itemize}
\item[${\textbf{\textsf{pt1}}}$]
According to
\refdef{deffamilyofclosedsets},
it is obvious.
\item[${\textbf{\textsf{pt2}}}$]
Considering the fact that,
$\empty\in\CSs{\X}$
and
\begin{align}
\X\setminus\empty&=\X\cr
&\in\topology{},
\end{align}
\refdef{deffamilyofclosedsets}
implies,
\begin{equation}
\empty\in\Fclosed{\X}{\topology{}}.
\end{equation}
\item[${\textbf{\textsf{pt3}}}$]
Considering that
$\X\in\CSs{\X}$,
and
\begin{align}
\X\setminus\X&=\empty\cr
&\in\topology{},
\end{align}
\refdef{deffamilyofclosedsets}
implies,
\begin{equation}
\X\in\Fclosed{\X}{\topology{}}.
\end{equation}
\item[${\textbf{\textsf{pt4}}}$]
Each $\U_1$
and
$\U_2$ is taken as e closed set of $\Xt$. So,
\begin{align}
&\(\U_1,\U_2\)\in\CSs{\X}\times\CSs{\X},\label{thmclosedsetspt4eq1}\\
&\(\X\setminus\U_1,\X\setminus\U_2\)\in\topology{}\times\topology{}.\label{thmclosedsetspt4eq2}
\end{align}
\Ref{thmclosedsetspt4eq2}
and
\refdef{deftopologicalspace}
imply,
\begin{equation}
\(\X\setminus\U_1\)\cap\(\X\setminus\U_2\)\in\topology{}.
\end{equation}
accordingly, by considering that,
\begin{align}
\X\setminus\(\U_1\cup\U_2\)=\(\X\setminus\U_1\)\cap\(\X\setminus\U_2\),
\end{align}
it is seen that,
\begin{equation}
\X\setminus\(\U_1\cup\U_2\)\in\topology{}.
\end{equation}
Additionally, according to
\Ref{thmclosedsetspt4eq1}
it is clear that
\begin{equation}
\U_1\cup\U_2\in\CSs{\X}.
\end{equation}
Therefore, according to
\refdef{deffamilyofclosedsets}
\begin{equation}
\U_1\cup\U_2\in\Fclosed{\X}{\topology{}}.
\end{equation}
\item[${\textbf{\textsf{pt5}}}$]
$\sC$
is taken as a collection of closed sets of.
$\Xt$
\begin{equation}
\sC\subseteq\Fclosed{\X}{\topology{}}.
\end{equation}
So,
\begin{align}
&\Foreach{\U}{\sC}\U\in\CSs{\X},\label{thmclosedsetspt5eq1}\\
&\Foreach{\U}{\sC}\X\setminus\U\in\topology{}.\label{thmclosedsetspt5eq2}
\end{align}
\Ref{thmclosedsetspt5eq2}
and
\refdef{deftopologicalspace}
imply,
\begin{equation}
\bigcup_{\U\in\sC}\(\X\setminus\U\)\in\topology{}.
\end{equation}
Therefore, considering that,
\begin{align}
\X\bigg\backslash\(\bigcap_{\U\in\sC}\U\)=\bigcup_{\U\in\sC}\(\X\setminus\U\),
\end{align}
it is seen that,
\begin{equation}
\X\bigg\backslash\(\bigcap_{\U\in\sC}\U\)\in\topology{}.
\end{equation}
Additionally, according to
\Ref{thmclosedsetspt5eq1},
it is clear that,
\begin{equation}
\bigcap_{\U\in\sC}\U\in\CSs{\X}.
\end{equation}
As a result of these, and
\refdef{deffamilyofclosedsets},
\begin{equation}
\bigcap_{\U\in\sC}\U\in\Fclosed{\X}{\topology{}}.
\end{equation}
\end{itemize}
\endthm
%%%%%%%%%%%%%%%%%%%%%%%%%%%%%%%%%%%%%%%%%%%%%%%%%%%%%%%%%%%%%%%%%%%%%%%%%%
\definition\label{defcnbdclassofsets}
$\Xt=\opair{\X}{\topology{}}$ is taken as a topological-space.
$\cnei{\Xt}$
is defined to be this operator.
\begin{itemize}
\item[${\textbf{\textsf{CNei1}}}$]
\hfill
$\cnei{\Xt}\indef\Func{\CSs{\X}}{\CSs{\CSs{\X}}}.$
\item[${\textbf{\textsf{CNei2}}}$]
\hfill
$\Foreach{\asubset}{\CSs{\X}}
\func{\cnei{\Xt}}{\asubset}\eqdef\defset{\U}{\Fclosed{\X}{\topology{}}}{\asubset\subseteq\U}.$
\end{itemize}
That is,
$\cnei{\Xt}$
is defined to be the operator from
$\CSs{\X}$
to
$\CSs{\CSs{\X}}$,
so that for every
$\asubset$
in
$\CSs{\X}$,
$\func{\cnei{\Xt}}{\asubset}$
is the set of all closed sets of
$\Xt$
that include
$\asubset$.\\
%\begin{itemize}
%\item
For every
$\asubset$
in
$\CSs{\X}$,
each element of
$\func{\cnei{\Xt}}{\asubset}$
is called a $\quotl$closed neighbourhood of
$\asubset$
in (the topological-space)
$\Xt$$\quotr$.
%\end{itemize}
\endef
%%%%%%%%%%%%%%%%%%%%%%%%%%%%%%%%%%%%%%%%%%%%%%%%%%%%%%%%%%%%%%%%%%%%%%%%%%%%%%%%%%%%%%%%%%%%%%%%%%%%%%%%%%
\theorem\label{thmintersectionoftopologies}
$\X$ is taken as a set, and
$\Tt$ as a non-empty collection of topologies on $\X$. That is, $\Tt$
is taken as a non-empty subset of
$\Ctops{\X}$
(or an element of
$\compl{\CSs{\Ctops{\X}}}{\seta{\empty}}$)
$\displaystyle\opair{\X}{\intersection{\Tt}}$
is a topological-space.
\begin{equation}
\intersection{\Tt}\in\Ctops{\X}.
\end{equation}
\proof
Considering that for every
$\topology{}$
in
$\Tt$,
$\(\X,\topology{}\)$
is a topological-space, and according to
\refdef{deftopologicalspace},
\begin{align}
&\Foreach{\topology{}}{\Tt}\(\empty\in\topology{}\),\label{thmintersectionoftopologiesp1}\\
&\Foreach{\topology{}}{\Tt}\(\X\in\topology{}\),\label{thmintersectionoftopologiesp2}\\
&\Foreach{\topology{}}{\Tt}
\[\Foreach{\opair{\U_{1}}{\U_{2}}}{\topology{}\times\topology{}}
\U_{1}\cap\U_{2}\in\topology{}\],\label{thmintersectionoftopologiesp3}\\
&\Foreach{\topology{}}{\Tt}\[
\Foreachs{\sC}{\topology{}}\(\bigcup_{\U\in\sC}\U\)\in\topology{}\].\label{thmintersectionoftopologiesp4}
\end{align}
\begin{itemize}
\item[${\textbf{\textsf{p1}}}$]
According to
\Ref{thmintersectionoftopologiesp1},
\begin{equation}
\empty\in\(\bigcap_{\topology{}\in\Tt}\topology{}\).
\end{equation}
\item[${\textbf{\textsf{p2}}}$]
According to
\Ref{thmintersectionoftopologiesp2},
\begin{equation}
\X\in\(\bigcap_{\topology{}\in\Tt}\topology{}\).
\end{equation}
\item[${\textbf{\textsf{p3}}}$]
$\opair{\U_1}{\U_2}$
is taken as an arbitrary element of
$\displaystyle\(\bigcap_{\topology{}\in\Tt}\topology{}\times\bigcap_{\topology{}\in\Tt}\topology{}\)$
It is obvious that,
\begin{equation}
\Foreach{\topology{}}{\Tt}\opair{\U_1}{\U_2}\in\topology{}\times\topology{}.
\end{equation}
This and
\Ref{thmintersectionoftopologiesp3}
imply that,
\begin{equation}
\Foreach{\topology{}}{\Tt}\U_1\cap\U_2\in\topology{},
\end{equation}
and accordingly,
\begin{equation}
\U_1\cap\U_2\in\(\bigcap_{\topology{}\in\Tt}\topology{}\).
\end{equation}
\item[${\textbf{\textsf{p4}}}$]
$\sC$
is taken as a collection of sets in
$\displaystyle\bigcap_{\topology{}\in\Tt}\topology{}$
That is,
\begin{equation}
\sC\subseteq\(\bigcap_{\topology{}\in\Tt}\topology{}\).
\end{equation}
So, considering that,
\begin{equation}
\Foreach{\topology{}}{\Tt}\(\bigcap_{\topology{}\in\Tt}\topology{}\)\subseteq\topology{},
\end{equation}
it is clear that,
\begin{equation}
\Foreach{\topology{}}{\Tt}\sC\subseteq\topology{}.
\end{equation}
This and \Ref{thmintersectionoftopologiesp4}
imply,
\begin{equation}
\Foreach{\topology{}}{\Tt}\(\bigcup_{\U\in\sC}\U\)\in\topology{},
\end{equation}
and accordingly,
\begin{equation}
\(\bigcup_{\U\in\sC}\U\)\in\(\bigcap_{\topology{}\in\Tt}\topology{}\).
\end{equation}
\end{itemize}
Thus, according to
\refdef{deftopologicalspace},
it is clear that
$\displaystyle\opair{\X}{\bigcap_{\topology{}\in\Tt}\topology{}}$
is a topological-space.
\endthm
%%%%%%%%%%%%%%%%%%%%%%%%%%%%%%%%%%%%%%%%%%%%%%%%%%%%%%%%%%%%%%%%%%%%%%%%%%%%%%%%%%%%%%%%%%%%%%%%%%%%%%%%%%%%%%%%%%%%
\definition\label{deffinersoarsertopology}
$\X$ is taken as a set, and each
$\topology{}$
and
$\p{\topology{}}$
as a topology on
$\X$.
\begin{itemize}
\item
It is said that $\quotl$$\p{\topology{}}$ is finer that $\topology{}$$\quotr$,
or $\quotl$$\topology{}$ is coarser than $\p{\topology{}}$$\quotr$, iff
$\topology{}\subseteq\p{\topology{}}$.
\item
It is said that
$\quotl$$\p{\topology{}}$ is strictly finer than $\p{\topology{}}$$\quotr$,
or $\quotl$$\p{\topology{}}$ is strictly coarser than$\p{\topology{}}$$\quotr$, iff
$\topology{}\subset\p{\topology{}}$.
\end{itemize}
\endef
%%%%%%%%%%%%%%%%%%%%%%%%%%%%%%%%%%%%%%%%%%%%%%%%%%%%%%%%%%%%%%%%%%%%%%%%%%%%%%%%%%%%%%%%%%%%%%%%%%%%%%%%%%%%%%%%%%%%
\theorem\label{thmclosedsetsofcoarsertopology}
$\X$ is taken as a set, and each
$\topology{}$
and
$\p{\topology{}}$
as a topology on
$\X$. If
$\topology{}$
is coarser than
$\p{\topology{}}$, then
every closed set of the topological-space
$\opair{\X}{\topology{}}$
is a closed set of the topological-space
$\opair{\X}{\p{\topology{}}}$.
That is,
\begin{equation}
\bigg(\topology{}\subseteq\p{\topology{}}\bigg)\then
\bigg(\Fclosed{\X}{\topology{}}\subseteq\Fclosed{\X}{\p{\topology{}}}\bigg).
\end{equation}
\prooff
It is supposed that,
\begin{equation}
\topology{}\subseteq\p{\topology{}}.
\end{equation}
So, considering that,
\begin{equation}
\Foreach{\U}{\Fclosed{\X}{\topology{}}}
\(\compl{\X}{\U}\)\in\topology{},
\end{equation}
it is clear that
\begin{equation}
\Foreach{\U}{\Fclosed{\X}{\topology{}}}
\(\compl{\X}{\U}\)\in\p{\topology{}},
\end{equation}
and hence according to
\refdef{deffamilyofclosedsets},
\begin{equation}
\Foreach{\U}{\Fclosed{\X}{\topology{}}}
\U\in\Fclosed{\X}{\p{\topology{}}}.
\end{equation}
\endthm
\endef
%%%%%%%%%%%%%%%%%%%%%%%%%%%%%%%%%%%%%%%%%%%%%%%%%%%%%%%%%%%%%%%%%%%%%%%%%%%%%%%%%%%%%%%%%%%%%%%%%%%%%%%%%%%%%%%%%%%%
\definition
$\Xt=\opair{\X}{\topology{}}$
is taken as a topological-space.
\begin{itemize}
\item
$\opair{\X}{\topology{}}$
is called a $\quotl$finite topological-space$\quotr$, iff
$\X$
is a finite set.
\item
$\opair{\X}{\topology{}}$
is called an $\quotl$infinite topological-space$\quotr$, iff
$\X$
is an infinite set.
\item
$\opair{\X}{\topology{}}$
is called a $\quotl$countably-infinite topolofical-space$\quotr$, iff
$\X$
is a countably-infinite set.
\item
$\opair{\X}{\topology{}}$
is called an $\quotl$uncountable topological-space$\quotr$, iff
$\X$
is an uncountable set.
\end{itemize}
\endef
\theorem\label{thm(in)discretetopology}
$\X$
is taken as a set.
\begin{itemize}
\item
$\(\X,\CSs{\X}\)$
is a topological-space. That is,
\begin{equation}
\CSs{\X}\in\Ctops{\X}.
\end{equation}
\item
$\(\X,\left\{\empty,\X\right\}\)$
is a topological-space. That is,
\begin{equation}
\seta{\binary{\empty}{\X}}\in\Ctops{\X}.
\end{equation}
\end{itemize}
\prooff
It is obvious.
\endthm
\definition\label{def(in)discretetopology}
$\X$
is taken as a set.
\begin{itemize}
\item
$\(\X,\CSs{\X}\)$
is called a $\quotl$discrete topological-space$\quotr$. Additionally,
$\CSs{\X}$
is referred to as the $\quotl$discrete topology on $\X$$\quotr$.
\item
$\(\X,\left\{\empty,\X\right\}\)$
is called a $\quotl$indiscrete topological-space$\quotr$. Additionally,
$\CSs{\X}$
is referred to as the $\quotl$indiscrete topology on $\X$$\quotr$.
\end{itemize}
\endef
%%%%%%%%%%%%%%%%%%%%%%%%%%%%%%%%%%%%%%%%%%%%%%%%%%%%%%%%%%%%%%%%%%%%%%%%%%%%%%%%%%%%%%%%%%%%%%%%%%%%%%%%%%%%%%%%%%%%
\theorem\label{thmdiscretetopologyclopensets}
$\X$ is taken as a set.
Every subset of $\X$ is clopen set of the discrete topological-space $\opair{\X}{\CSs{\X}}$. That is,
\begin{equation}
\Foreach{\asubset}{\CSs{\X}}
\asubset\in\Fclosed{\X}{\CSs{\X}}.
\end{equation}
\prooff
It is obvious.
\endthm
%%%%%%%%%%%%%%%%%%%%%%%%%%%%%%%%%%%%%%%%%%%%%%%%%%%%%%%%%%%%%%%%%%%%%%%%%%%%%%%%%%%%%%%%%%%%%%%%%%%%%%%%%%%%%%%%%%%%
\theorem\label{thmNandSdiscretetopology}
$\X$
is taken as a set
$\topology{}$
is taken as a topology on $\X$ (an element of $\Ctops{\X}$).
$\topology{}$
the discrete topology on $\X$, if-and-only-if
every singleton subset of $\X$
is an open set of $\opair{\X}{\topology{}}$. That is,
\begin{equation}
\[\topology{}=\CSs{\X}\]\thenn
\[\Foreach{\x}{\X}\seta{\x}\in\topology{}\].
\end{equation}
\prooff
According to
\refdef{deftopologicalspace}
and
\refthm{thm(in)discretetopology},
it is clear.
\endthm
%%%%%%%%%%%%%%%%%%%%%%%%%%%%%%%%%%%%%%%%%%%%%%%%%%%%%%%%%%%%%%%%%%%%%%%%%%%%%%%%%%%%%%%%%%%%%%%%%%%%%%%%%%%%%%%%%%%%
\theorem\label{thmoneelementdif}
$\opair{\X}{\topology{}}$
is taken as a topological-space, and
each $\Y$ and $\seta{a}$
are taken as a set, such that
\begin{equation}\label{thmoneelementdifeq1}
\Y=\X\cup\seta{a}.
\end{equation}
$\defSet{\seta{a}\cup\U}{\U\in\topology{}}\cup\seta{\empty}$
is a topology on $\Y$.
\proof
According to
\refdef{deftopologicalspace},
\begin{align}
&\X\in\topology{},\label{thmoneelementdifp01}\\
&\Foreach{\opair{\U_{1}}{\U_{2}}}{\topology{}\times\topology{}}\U_{1}\cap\U_{2}\in\topology{},\label{thmoneelementdifp02}\\
&\Foreachs{\sC}{\topology{}}\(\bigcup_{\U\in\sC}\U\)\in\topology{}\label{thmoneelementdifp03}.
\end{align}
$\topology{}^{\prime}$
is defined as,
\begin{equation}\label{thmoneelementdifp04}
\topology{}^{\prime}:=\defSet{\seta{a}\cup\U}{\U\in\topology{}}\cup\seta{\empty}.
\end{equation}
\begin{itemize}
\item[${\textbf{\textsf{p0}}}$]
It is clear that
\begin{equation}
\topology{}^{\prime}\subseteq\CSs{\Y}.
\end{equation}
\item[${\textbf{\textsf{p1}}}$]
It is clear that
\begin{equation}
\empty\in\topology{}^{\prime}.
\end{equation}
Additionally, considering that,
\begin{equation}
\X\in\topology{},
\end{equation}
according to \Ref{thmoneelementdifp04},
it is seen that,
\begin{equation}
\seta{a}\cup\X\in\topology{}^{\prime},
\end{equation}
and hence according to \Ref{thmoneelementdifeq1},
\begin{equation}
\Y\in\topology{}^{\prime}
\end{equation}
\item[${\textbf{\textsf{p2}}}$]
$\opair{\V_1}{\V_2}$
is taken as an arbitrary element of
$\topology{}^{\prime}\times\topology{}^{\prime}$
According to \Ref{thmoneelementdifp04},
there exists an element of
$\topology{}\times\topology{}$, like $\opair{\U_1}{\U_2}$, such that,
\begin{align}
\V_1=\seta{a}\cup\U_1,\label{thmoneelementdifp21}\\
\V_2=\seta{a}\cup\U_2.\label{thmoneelementdifp22}
\end{align}
Therefore,
\begin{align}\label{thmoneelementdifp23}
\V_1\cap\V_2&=\(\seta{a}\cup\U_1\)\cap\(\seta{a}\cup\U_2\)\cr
&=\seta{a}\cup\(\U_1\cap\U_2\).
\end{align}
In addition, according to \refdef{deftopologicalspace}, and considering that,
\begin{equation}\label{thmoneelementdifp24}
\opair{\U_1}{\U_2}\in\topology{}\times\topology{},
\end{equation}
it is clear that,
\begin{equation}\label{thmoneelementdifp25}
\U_1\cap\U_2\in\topology{}.
\end{equation}
As a result of this, \Ref{thmoneelementdifp04}, and \Ref{thmoneelementdifp23},
\begin{equation}\label{thmoneelementdifp26}
\V_1\cap\V_2\in\topology{}^{\prime}\times\topology{}^{\prime}.
\end{equation}
\item[${\textbf{\textsf{p3}}}$]
$\sC$
is taken as a collection of sets in $\topology{}^{\prime}$. That is,
\begin{equation}\label{thmoneelementdifp27}
\sC\subseteq\topology{}^{\prime}.
\end{equation}
According to
\Ref{thmoneelementdifp04},
\begin{equation}\label{thmoneelementdifp28}
\Foreach{\V}{\sC}\[\Exists{\U}{\topology{}}\V=\seta{a}\cup\U\].
\end{equation}
$f$
is taken as a function that,
\begin{align}
&f\in\Func{\sC}{\topology{}},\label{thmoneelementdifp29}\\
&\Foreach{\V}{\sC}\V=\seta{a}\cup\func{f}{\V}.\label{thmoneelementdifp210}
\end{align}
Thus,
\begin{align}\label{thmoneelementdifp211}
\bigcup_{\V\in\sC}\V&=\bigcup_{\V\in\sC}\[\seta{a}\cup\func{f}{\V}\]\cr
&=\seta{a}\cup\[\bigcup_{\V\in\sC}\func{f}{\V}\].
\end{align}
According to
\refdef{deftopologicalspace},
and by considering
\begin{equation}\label{thmoneelementdifp212}
\Foreach{\V}{\sC}\func{f}{\V}\in\topology{},
\end{equation}
it is clear that,
\begin{equation}\label{thmoneelementdifp213}
\[\bigcup_{\V\in\sC}\func{f}{\V}\]\in\topology{}.
\end{equation}
\Ref{thmoneelementdifp04}, \Ref{thmoneelementdifp211}, and \Ref{thmoneelementdifp213} imply,

\begin{equation}
\(\bigcup_{\V\in\sC}\V\)\in\topology{}^{\prime}.
\end{equation}
\end{itemize}
Therefore, based on \refdef{deftopologicalspace},
$\(\Y,\topology{}^{\prime}\)$
is a topological-space
\endthm
%%%%%%%%%%%%%%%%%%%%%%%%%%%%%%%%%%%%%%%%%%%%%%%%%%%%%%%%%%%%%%%%%%%%%%%%%%%%%%%%%%%%%%%%%%%%%%%%%%%%%%%%%%%%%%%%%%%%%%%%%%%%%%%%%%%%%%%%%%%%%%%%%%%%%%%%%
\subsection{
One-point Topological Space
}
\theorem\label{thmsingletontopology0}
$\x$
is taken as a set.
$\opair{\seta{\x}}{\CSs{\seta{\x}}}$
is both a discrete and an indiscrete topological-space.
\proof
According to \refdef{def(in)discretetopology},
and considering that,
$\CSs{\seta{\x}}=\seta{\binary{\empty}{\seta{\x}}}$,
it is clear.
\endthm
%%%%%%%%%%%%%%%%%%%%%%%%%%%%%%%%%%%%%%%%%%%%%%%%%%%%%%%%%%%%%%%%%%%%%%%%%%%%%%%%%%%%%%%%%%%%%%%%%%%%%%%%%%%%%%%%%%%%%%%%%%%%%%%%%%%%%%%%%%%%%%%%%%%%%%%%%
\theorem\label{thmsingletontopology}
$\x$
is taken as a set. The power-set of
$\seta{\x}$, that is the set $\seta{\binary{\empty}{\seta{\x}}}$,
is the only topology on $\seta{\x}$. That is,
\begin{align}\label{thmsingletontopologyeq1}
\Ctops{\seta{\x}}&=\seta{\CSs{\seta{\x}}}\cr
&=\seta{\seta{\binary{\empty}{\seta{\x}}}}.
\end{align}
\prooff
According to
\refdef{deftopologicalspace},
\begin{equation}
\Foreach{\topology{}}{\Ctops{\seta{\x}}}
\topology{}\subseteq\CSs{\seta{\x}},
\end{equation}
and
\begin{align}
\Foreach{\topology{}}{\Ctops{\seta{\x}}}
\CSs{\seta{\x}}&=\seta{\binary{\empty}{\seta{\x}}}\cr
&\subseteq\topology{}.
\end{align}
Thus,
\begin{equation}
\Foreach{\topology{}}{\Ctops{\seta{\x}}}
\topology{}=\CSs{\seta{\x}},
\end{equation}
that is,
\Ref{thmsingletontopologyeq1}.
\endthm
%%%%%%%%%%%%%%%%%%%%%%%%%%%%%%%%%%%%%%%%%%%%%%%%%%%%%%%%%%%%%%%%%%%%%%%%%%%%%%%%%%%%%%%%%%%%%%%%%%%%%%%%%%%%%%%%%%%%%%%%%%%%%%%%%%%%%%%%%%%%%%%%%%%%%%%%%
\definition\label{defsingletontopology}
$\x$
is taken as a set.
\begin{equation}
\singletonTS{\x}:=\opair{\seta{\x}}{\CSs{\seta{\x}}}.
\end{equation}
$\singletonTS{\x}$
is referred to as the $\quotl$one-point topological-space of the point $\x$$\quotr$.
\endef
%%%%%%%%%%%%%%%%%%%%%%%%%%%%%%%%%%%%%%%%%%%%%%%%%%%%%%%%%%%%%%%%%%%%%%%%%%%%%%%%%%%%%%%%%%%%%%%%%%%%%%%%%%%%%%%%%%%%%%%%%%%%%%%%%%%%%%%%%%%%%%%%%%%%%%%%%
%%%%%%%%%%%%%%%%%%%%%%%%%%%%%%%%%%%%%%%%%%%%%%%%%%%%%%%%%%%%%%%%%%%%%%%%%%%%%%%%%%%%%%%%%%%%%%%%%%%%%%%%%%%%%%%%%%%%%%%%%%%%%%%%%%%%%%%%%%%%%%%%%%%%%%%%%
%%%%%%%%%%%%%%%%%%%%%%%%%%%%%%%%%%%%%%%%%%%%%%%%%%%%%%%%%%%%%%%%%%%%%%%%%%%%%%%%%%%%%%%%%%%%%%%%%%%%%%%%%%%%%%%%%%%%%%%%%%%%%%%%%%%%%%%%%%%%%%%%%%%%%%%%%
%%%%%%%%%%%%%%%%%%%%%%%%%%%%%%%%%%%%%%%%%%%%%%%%%%%%%%%%%%%%%%%%%%%%%%%%%%%%%%%%%%%%%%%%%%%%%%%%%%%%%%%%%%%%%%%%%%%%%%%%%%%%%%%%%%%%%%%%%%%%%%%%%%%%%%%%%
%%%%%%%%%%%%%%%%%%%%%%%%%%%%%%%%%%%%%%%%%%%%%%%%%%%%%%%%%%%%%%%%%%%%%%%%%%%%%%%%%%%%%%%%%%%%%%%%%%%%%%%%%%%%%%%%%%%%%%%%%%%%%%%%%%%%%%%%%%%%%%%%%%%%%%%%%
\section{
Bases and Sub-bases of Topological Spaces
}
\subsection{
Base
}
\definition\label{defbase}
$\Xt=\opair{\X}{\topology{}}$
is taken as a topological-space, and
$\base$
as an element of
$\CSs{\CSs{\X}}$
(a collection of subsets of $\X$)
$\base$
is called a $\quotl$base for the topological-space $\Xt$$\quotr$, iff
$\base$
is a subset of
$\topology{}$
(a collection of open sets of $\Xt$), and every element of $\topology{}$
(every open set of $\Xt$)
equals the union of the elements of a subset of $\base$. That is,
$\base$
is called a $\quotl$base for the topological-space $\opair{\X}{\topology{}}$$\quotr$, iff these properties
are satisfied.
\begin{itemize}
\item[${\textbf{\textsf{B1}}}$]
%\hfill
$\base\subseteq\topology{}.$
\item[${\textbf{\textsf{B2}}}$]
%\hfill
$\displaystyle
\Foreach{\U}{\topology{}}\[\Existss{\sC}{\base}\U=\(\union{\sC}\)\].$
\end{itemize}
\endef
%%%%%%%%%%%%%%%%%%%%%%%%%%%%%%%%%%%%%%%%%%%%%%%%%%%%%%%%%%%%%%%%%%%%%%%%%%%%%%%%%%%%%%%%%%%%%%%%%%%%%%%
\theorem\label{thmSandNconditionsofbase0}
$\Xt=\opair{\X}{\topology{}}$
is taken as a topological-space, and $\base$ as a subset of $\topology{}$.
$\base$ is a base for $\Xt$, if-and-only-if,
\begin{equation}\label{thmSandNconditionsofbase0eq1}
\topology{}=\defset{\U}{\CSs{\X}}{\[\Exists{\sC}{\CSs{\base}}\U=\(\union{\sC}\)\]}.
\end{equation}
That is,
$\base$ is a base for $\Xt$, if-and-only-if the collection of all open sets of $\Xt$
equals the collection of all unions of elements of $\CSs{\base}$.
\proof
\begin{itemize}
\item[${\textbf{\textsf{p1}}}$]
It is assumed that,
$\base$ is a base for $\Xt$.
According to \refdef{defbase},
\begin{equation}\label{thmSandNconditionsofbase0p11}
\Foreach{\U}{\topology{}}\[\Existss{\Bb}{\base}\U=\bigcup_{\V\in\Bb}\V\].
\end{equation}
Accordingly,
\begin{equation}\label{thmSandNconditionsofbase0p12}
\topology{}\subseteq\defset{\U}{\CSs{\X}}{\[\Exists{\sC}{\CSs{\base}}\U=\bigcup_{\V\in\sC}\V\]}.
\end{equation}
Additionally, considering that,
$\base\subseteq\topology{}$,
it is clear that,
\begin{equation}\label{thmSandNconditionsofbase0p13}
\Foreach{\sC}{\CSs{\base}}\sC\subseteq\topology{},
\end{equation}
and hence according to \refdef{deftopologicalspace},
\begin{equation}\label{thmSandNconditionsofbase0p14}
\Foreach{\sC}{\CSs{\base}}\(\bigcup_{\V\in\sC}\V\)\in\topology{}.
\end{equation}
Thus,
\begin{equation}\label{thmSandNconditionsofbase0p15}
\defset{\U}{\CSs{\X}}{\[\Exists{\sC}{\CSs{\base}}\U=\bigcup_{\V\in\sC}\V\]}\subseteq\topology{}.
\end{equation}
\Ref{thmSandNconditionsofbase0p12} and \Ref{thmSandNconditionsofbase0p15} imply
\Ref{thmSandNconditionsofbase0eq1}.
\item[${\textbf{\textsf{p2}}}$]
If \Ref{thmSandNconditionsofbase0eq1} is true, then based on \refdef{defbase},
it is clear that $\base$ is a base for $\Xt$.
\end{itemize}
\endthm
%%%%%%%%%%%%%%%%%%%%%%%%%%%%%%%%%%%%%%%%%%%%%%%%%%%%%%%%%%%%%%%%%%%%%%%%%%%%%%%%%%%%%%%%%%%%%%%%%%%%%%%
\theorem\label{thmtrivialbase}
$\opair{\X}{\topology{}}$
is taken as a topological-space.
$\topology{}$
is a base for the topological-space $\opair{\X}{\topology{}}$
\proof
It is trivial.
\endthm
%%%%%%%%%%%%%%%%%%%%%%%%%%%%%%%%%%%%%%%%%%%%%%%%%%%%%%%%%%%%%%%%%%%%%%%%%%%%%%%%%%%%%%%%
\definition\label{defcollectionofbases}
$\Xt=\opair{\X}{\topology{}}$
is taken as a topological-space. The collection of all bases for $\opair{\X}{\topology{}}$
is denoted by $\Cbase{\Xt}$.
\endef
%%%%%%%%%%%%%%%%%%%%%%%%%%%%%%%%%%%%%%%%%%%%%%%%%%%%%%%%%%%%%%%%%%%%%%%%%%%%%%%%%%%%%%%%
\theorem\label{thmSandNconditionsofbase1}
$\Xt=\opair{\X}{\topology{}}$
is taken as a topological-space, and $\base$ as a subset of $\CSs{\X}$.
$\base$ is a base for $\Xt$, if-and-only-if
$\base$
is a collection of open sets of $\Xt$ that for every
$\U$
in
$\topology{}$,
and for every $\x$ in $\U$,
there exists an element of $\base$ that contains $\x$ and is a subset of $\U$. That is,
\begin{align}
\(\base\in\Cbase{\Xt}\)\thenn
\[\base\subseteq\topology{},~\(\Foreach{\U}{\topology{}}\left\{\Foreach{\x}{\U}\[\Exists{\B}{\base}\(\x\in\B,~\B\subseteq\U\)\]\right\}\)\].
\end{align}
\prooff
\begin{itemize}
\item[${\textbf{\textsf{p1}}}$]
$\base$
is taken as a base for the topological-space $\opair{\X}{\topology{}}$.
\begin{itemize}
\item[${\textbf{\textsf{p1-1}}}$]
According to \refdef{defbase},
it is clear that,
\begin{equation}
\base\subseteq\topology{}.
\end{equation}
\endp
\item[${\textbf{\textsf{p1-2}}}$]
$\U$
is taken as an arbitrary element of
$\topology{}$
(an open set of $\Xt$).
According to \refdef{defbase},
there exists a subset of
$\base$, like $\Bb$, such that,
\begin{equation}\label{thmSandNconditionsofbase1p1-2eq1}
\U=\bigcup_{\V\in\Bb}\V.
\end{equation}
\begin{itemize}
\item[${\textbf{\textsf{p1-2-1}}}$]
$\x$
is taken as an arbitrary element of $\U$. According to \Ref{thmSandNconditionsofbase1p1-2eq1},
\begin{equation}\label{thmSandNconditionsofbase1p1-2-1eq1}
\x\in\bigcup_{\V\in\Bb}\V,
\end{equation}
and hence there exists an element of
$\Bb$, like $\V_{\x}$, that,
\begin{equation}\label{thmSandNconditionsofbase1p1-2-1eq2}
\x\in\V_{\x}.
\end{equation}
Considering that,
\begin{equation}\label{thmSandNconditionsofbase1p1-2-1eq3}
\V_{\x}\in\Bb\subseteq\base,
\end{equation}
it is clear that,
\begin{equation}\label{thmSandNconditionsofbase1p1-2-1eq4}
\V_{\x}\in\base.
\end{equation}
\Ref{thmSandNconditionsofbase1p1-2eq1}
and
\Ref{thmSandNconditionsofbase1p1-2-1eq3}
imply,
\begin{equation}
\V_{\x}\subseteq\U.
\end{equation}
\endp
\end{itemize}
So,
\begin{equation}
\Foreach{\x}{\U}\[\Exists{\B}{\base}\(\x\in\B,~\B\subseteq\U\)\].
\end{equation}
\endp
\end{itemize}
\endp
\item[${\textbf{\textsf{p2}}}$]
$\base$
is taken as such a collection that,
\begin{align}\label{thmSandNconditionsofbase1p2-2eq1}
&\base\subseteq\topology{},\\
&\Foreach{\U}{\topology{}}\left\{\Foreach{\x}{\U}\[\Exists{\B}{\base}\(\x\in\B,~\B\subseteq\U\)\]\right\}.
\end{align}
\begin{itemize}
\item[${\textbf{\textsf{p2-1}}}$]
$\U$
is taken as an arbitrary element of $\topology{}$.
$\f$
is taken as such a function that,
\begin{align}
&\f\in\Func{\U}{\base},\\
&\Foreach{\x}{\U}\x\in\func{\f}{\x}\subseteq\U.
\end{align}
These simply imply that,
\begin{equation}
\U=\bigcup_{\x\in\U}\func{\f}{\x}.
\end{equation}
Therefore, by defining
$\sCi$
as,
\begin{equation}
\sCi:=\defSet{\func{\f}{\x}}{\x\in\U},
\end{equation}
it is clear that
\begin{align}
&\sCi\subseteq\base,\\
&\U=\bigcup_{\V\in\sCi}\V.
\end{align}
\endp
\end{itemize}
So,
\begin{equation}
\Foreach{\U}{\topology{}}\[\Existss{\sC}{\base}\U=\bigcup_{\V\in\sC}\V\].
\end{equation}
Thus according to this and \Ref{thmSandNconditionsofbase1p2-2eq1}, based on \refdef{defbase},
it can be seen that,
$\base$
is a base for $\Xt$.
\endp
\end{itemize}
\endthm
%%%%%%%%%%%%%%%%%%%%%%%%%%%%%%%%%%%%%%%%%%%%%%%%%%%%%%%%%%%%%%%%%%%%%%%%%%%%%%%%%%%%%%%%%%%%%%%%%%%%%%%%%%%%%%%%%
\theorem\label{thmopensetpointsbase}
$\Xt=\opair{\X}{\topology{}}$
is taken as a topological-space, and
$\base$ as a base for $\Xt$, and $\U$ as a subset of $\X$.
$\U$ is an open set of $\Xt$, if-and-only-if for every $\x$ in $\U$,
there exists an element of $\base$ that contains $\x$, and is a subset of $\U$.
\begin{equation}
\(\U\in\topology{}\)\thenn
\left\{\Foreach{\x}{\U}\[\Exists{\B}{\base}\(\x\in\B,~\B\subseteq\U\)\]\right\}.
\end{equation}
\proof
\begin{itemize}
\item[${\textbf{\textsf{p1}}}$]
According to,
\refthm{thmSandNconditionsofbase1},
it is clear that,
\begin{equation}
\(\U\in\topology{}\)\then
\left\{\Foreach{\x}{\U}\[\Exists{\B}{\base}\(\x\in\B,~\B\subseteq\U\)\]\right\}.
\end{equation}
\item[${\textbf{\textsf{p2}}}$]
$\U$
is taken as such a set that,
\begin{equation}
\Foreach{\x}{\U}\[\Exists{\B}{\base}\(\x\in\B,~\B\subseteq\U\)\].
\end{equation}
$\f$
is taken as such a function that,
\begin{align}
&\f\in\Func{\U}{\base},\\
&\Foreach{\x}{\U}\x\in\func{\f}{\x}\subseteq\U.
\end{align}
Therefore,
\begin{equation}
\U=\bigcup_{\x\in\U}\func{\f}{\x}.
\end{equation}
So, by defining $\sCi$ as,
\begin{equation}
\sCi:=\defSet{\func{\f}{\x}}{\x\in\U},
\end{equation}
it is clear that
\begin{align}
&\sCi\subseteq\base,\\
&\U=\bigcup_{\V\in\sCi}\V.
\end{align}
Therefore, based on \refthm{thmSandNconditionsofbase0}, it is clear that,
\begin{equation}
\U\in\topology{}.
\end{equation}
\end{itemize}
\endthm
%%%%%%%%%%%%%%%%%%%%%%%%%%%%%%%%%%%%%%%%%%%%%%%%%%%%%%%%%%%%%%%%%%%%%%%%%%%%%%%%%%%%%%%%%%%%%%%%%%%%%%%%
%%%%%%%%%%%%%%%%%%%%%%%%%%%%%%%%%%%%%%%%%%%%%%%%%%%%%%%%%%%%%%%%%%%%%%%%%%%%%%%%%%%%%%%%%%%%%%%%%%%%%%%%%%%%%%%%%
%%%%%%%%%%%%%%%%%%%%%%%%%%%%%%%%%%%%%%%%%%%%%%%%%%%%%%%%%%%%%%%%%%%%%%%%%%%%%%%%%%%%%%%%%%%%%%%%%%%%%%%%%%%%%%%%%
\theorem\label{thmopensetpointstopology}
$\Xt=\opair{\X}{\topology{}}$
is taken as a topological-space
\begin{equation}
\Foreach{\U}{\CSs{\X}}
\[\(\U\in\topology{}\)\thenn
\left\{\Foreach{\x}{\U}\[\Exists{\V}{\topology{}}\(\x\in\V,~\V\subseteq\U\)\]\right\}\].
\end{equation}
In other words,
\begin{equation}
\Foreach{\U}{\CSs{\X}}
\[\(\U\in\topology{}\)\thenn
\left\{\Foreach{\x}{\U}\[\Exists{\V}{\func{\nei{\Xt}}{\seta{\x}}}\(\V\subseteq\U\)\]\right\}\].
\end{equation}
\prooff
According to,
\refdef{defnbdclassofsets},
\refthm{thmtrivialbase},
and
\refthm{thmopensetpointsbase},
it is clear.
\endthm
%%%%%%%%%%%%%%%%%%%%%%%%%%%%%%%%%%%%%%%%%%%%%%%%%%%%%%%%%%%%%%%%%%%%%%%%%%%%%%%%%%%%%%%%%%%%%%%%%%%%%%%%%%%%%%%%%
%%%%%%%%%%%%%%%%%%%%%%%%%%%%%%%%%%%%%%%%%%%%%%%%%%%%%%%%%%%%%%%%%%%%%%%%%%%%%%%%%%%%%%%%%%%%%%%%%%%%%%%%%%%%%%%%%
\theorem\label{thmbaseNcondition1}
$\Xt=\opair{\X}{\topology{}}$
is taken as a topological-space, and $\base$ as a bases for $\Xt$.
\begin{align}
&\X=\bigcup_{\B\in\base}\B,\label{thmbaseNcondition1eq1}\\
&\Foreach{\(\B_1,\B_2\)}{\base\times\base}\(\Existss{\sC}{\base}\B_1\cap\B_2=\bigcup_{\V\in\sC}\V\).
\end{align}
\prooff
\begin{itemize}
\item[${\textbf{\textsf{p1}}}$]
Considering that
$\X\in\topology{}$,
and according to \refdef{defbase},
there exists a subset of $\base$, like $\sCi$m that,
\begin{equation}
\X=\bigcup_{\V\in\sCi}\V.
\end{equation}
Considering that,
$\sCi\subseteq\base$,
it is clear that,
\begin{equation}
\(\bigcup_{\V\in\sCi}\V\)\subseteq\(\bigcup_{\B\in\base}\B\).
\end{equation}
Therefore,
\begin{equation}
\X\subseteq\(\bigcup_{\B\in\base}\B\).
\end{equation}
Additionally, considering that,
$\base\subseteq\CSs{\X}$,
it is clear that,
\begin{equation}
\(\bigcup_{\B\in\base}\B\)\subseteq\X.
\end{equation}
Therefore, \Ref{thmbaseNcondition1eq1} is achieved.
\endp
\item[${\textbf{\textsf{p2}}}$]
Each $\B_1$ and $\B_2$
is taken as an element of $\base$.
Considering that, $\base$
is a base for $\Xt$, clearly $\base\subseteq\topology{}$,
and accordingly,
\begin{equation}
\opair{\B_1}{\B_2}\in\topology{}\times\topology{}.
\end{equation}
Hence according to \refdef{deftopologicalspace},
\begin{equation}
\B_1\cap\B_2\in\topology{},
\end{equation}
and accordingly, based on \refdef{defbase},
\begin{equation}
\Existss{\sC}{\base}\B_1\cap\B_2=\bigcup_{\V\in\sC}\V.
\end{equation}
\endp
\end{itemize}
\endthm
%%%%%%%%%%%%%%%%%%%%%%%%%%%%%%%%%%%%%%%%%%%%%%%%%%%%%%%%%%%%%%%%%%%%%%%%%%%%%%%%%%%%%%%%%%%%%%%%%%%%%%%%%%%%%%%%%
%%%%%%%%%%%%%%%%%%%%%%%%%%%%%%%%%%%%%%%%%%%%%%%%%%%%%%%%%%%%%%%%%%%%%%%%%%%%%%%%%%%%%%%%%%%%%%%%%%%%%%%%%%%%%%%%%
\theorem\label{thmbaseScondition1}
$\X$
is taken as a set, and $\base$ as such a subset of $\CSs{\X}$
(collection of subsets of $\X$) that,
\begin{align}
%\begin{split}
&\X=\bigcup_{\B\in\base}\B,\label{thmbaseScondition1eq1}\\
&\Foreach{\(\B_1,\B_2\)}{\base\times\base}\(\Existss{\sC}{\base}\B_1\cap\B_2=\bigcup_{\V\in\sC}\V\).\label{thmbaseScondition1eq2}
%\end{split}
\end{align}
$\defset{\U}{\CSs{\X}}{\[\Exists{\sC}{\CSs{\base}}\U=\bigcup_{\V\in\sC}\V\]}$
is a topology on $\X$, and $\base$ is a base for the topological-space
$\opair{\X}{\defset{\U}{\CSs{\X}}{\[\Exists{\sC}{\CSs{\base}}\U=\bigcup_{\V\in\sC}\V\]}}$.
\proof
$\topology{}$
is defined as,
\begin{equation}\label{thmbaseScondition1p}
\topology{}:=\defset{\U}{\CSs{\X}}{\[\Exists{\sC}{\CSs{\base}}\U=\bigcup_{\V\in\sC}\V\]}.
\end{equation}
\begin{itemize}
\item[${\textbf{\textsf{p0}}}$]
It is clear that,
\begin{equation}
\topology{}\subseteq\CSs{\X}.
\end{equation}
\endp
\item[${\textbf{\textsf{p1}}}$]
Considering that,
\begin{align}
&\empty\in\CSs{\X},\\
&\empty\in\CSs{\base},\\
&\empty=\bigcup_{\V\in\empty}\V,
\end{align}
it is clear that,
\begin{equation}
\empty\in\topology{}.
\end{equation}
\endp
\item[${\textbf{\textsf{p2}}}$]
Since
$\base\in\CSs{\base}$,
\Ref{thmbaseScondition1eq1}
and
\Ref{thmbaseScondition1p}
imply,
\begin{equation}
\X\in\topology{}.
\end{equation}
\endp
\item[${\textbf{\textsf{p3}}}$]
$\sCi$
is taken as an arbitrary subset of $\topology{}$.
According to \Ref{thmbaseScondition1p},
\begin{equation}\label{thmbaseScondition1p31}
\Foreach{\U}{\sCi}\[\Exists{\sC}{\CSs{\base}}\U=\bigcup_{\V\in\sC}\V\].
\end{equation}
$\f$ is taken as such a function that,
\begin{align}
&\f\in\Func{\sCi}{\CSs{\base}},\label{thmbaseScondition1p32}\\
&\Foreach{\U}{\sCi}\U=\bigcup_{\V\in\func{\f}{\U}}\V.\label{thmbaseScondition1p33}
\end{align}
Thus,
\begin{align}\label{thmbaseScondition1p34}
\bigcup_{\U\in\sCi}\U=\bigcup_{\U\in\sCi}\(\bigcup_{\V\in\func{\f}{\U}}\V\),
\end{align}
and accordingly, by defining $\Phi$ as,
\begin{equation}\label{thmbaseScondition1p35}
\Phi:=\defset{\V}{\CSs{\X}}{\[\V\in\func{\f}{\U},~\U\in\sCi\]},
\end{equation}
it can be seen that,
\begin{equation}\label{thmbaseScondition1p36}
\bigcup_{\U\in\sCi}\U=\bigcup_{\V\in\Phi}\V.
\end{equation}
\Ref{thmbaseScondition1p32} and \Ref{thmbaseScondition1p35}
imply,
\begin{equation}\label{thmbaseScondition1p37}
\Phi\in\CSs{\base}.
\end{equation}
\Ref{thmbaseScondition1p},
\Ref{thmbaseScondition1p36}, and
\Ref{thmbaseScondition1p37}
imply,
\begin{equation}
\(\bigcup_{\U\in\sCi}\U\)\in\topology{}.
\end{equation}
\endp
\item[${\textbf{\textsf{p4}}}$]
Each $\U_1$
and
$\U_2$
is taken as an element of $\topology{}$.
Each $\sC_1$ and $\sC_2$ is taken as such a subset of $\base$ that,
\begin{align}
\U_1&=\bigcup_{\V\in\sC_1}\V,\label{thmbaseScondition1p41}\\
\U_2&=\bigcup_{\V\in\sC_2}\V.\label{thmbaseScondition1p42}
\end{align}
Therefore,
\begin{align}\label{thmbaseScondition1p43}
\U_1\cap\U_2&=\(\bigcup_{\V\in\sC_1}\V\)\cap\(\bigcup_{\V\in\sC_2}\V\)\cr
&=\bigcup_{\V\in\sC_1}\[\bigcup_{\V^{\prime}\in\sC_2}\(\V\cap\V^{\prime}\)\],
\end{align}
and accordingly, by defining $\Psi$ as,
\begin{equation}\label{thmbaseScondition1p44}
\Psi:=\defset{\U}{\CSs{\X}}{\[\Exists{\opair{\V}{\V^{\prime}}}{\sC_1\times\sC_2}\U=\V\cap\V^{\prime}\]},
\end{equation}
it is seen that,
\begin{equation}\label{thmbaseScondition1p45}
\U_1\cap\U_2=\bigcup_{\U\in\Psi}\U.
\end{equation}
Since,
\begin{align}
\sC_1&\subseteq\base,\label{thmbaseScondition1p46}\\
\sC_2&\subseteq\base,\label{thmbaseScondition1p47}
\end{align}
it is clear that,
\begin{equation}\label{thmbaseScondition1p48}
\Foreach{\opair{\V}{\V^{\prime}}}{\sC_1\times\sC_2}\opair{\V}{\V^{\prime}}\in\base\times\base.
\end{equation}
\Ref{thmbaseScondition1eq2},
\Ref{thmbaseScondition1p44},
and
\Ref{thmbaseScondition1p48}
imply,
\begin{equation}\label{thmbaseScondition1p49}
\Foreach{\U}{\Psi}\[\Exists{\sCii}{\CSs{\base}}\U=\bigcup_{\V\in\sCii}\V\].
\end{equation}
$\ff$ is taken as such a function that,
\begin{align}
&\ff\in\Func{\Psi}{\CSs{\base}},\label{thmbaseScondition1p410}\\
&\Foreach{\U}{\Psi}\U=\bigcup_{\V\in\func{\ff}{\U}}\V.\label{thmbaseScondition1p411}
\end{align}
Hence,
\begin{equation}\label{thmbaseScondition1p412}
\bigcup_{\U\in\Psi}\U=\bigcup_{\U\in\Psi}\(\bigcup_{\V\in\func{\ff}{\U}}\V\),
\end{equation}
and accordingly, by defining $\Delta$ as,
\begin{equation}\label{thmbaseScondition1p413}
\Delta:=\defset{\V}{\CSs{\X}}{\[\V\in\func{\ff}{\U},~\U\in\Psi\]},
\end{equation}
it is seen that,
\begin{equation}\label{thmbaseScondition1p414}
\bigcup_{\U\in\Psi}\U=\bigcup_{\V\in\Delta}\V.
\end{equation}
\Ref{thmbaseScondition1p410}
and
\Ref{thmbaseScondition1p413}
imply,
\begin{equation}\label{thmbaseScondition1p415}
\Delta\in\CSs{\base}.
\end{equation}
\Ref{thmbaseScondition1p},
\Ref{thmbaseScondition1p414}
and
\Ref{thmbaseScondition1p415}
imply,
\begin{equation}
\(\bigcup_{\U\in\Psi}\U\)\in\topology{}.
\end{equation}
This, together with \Ref{thmbaseScondition1p45} imply that,
\begin{equation}
\(\B_1\cap\B_2\)\in\topology{}.
\end{equation}
\endp
\end{itemize}
Therefore, based on \refdef{deftopologicalspace},
$\opair{\X}{\topology{}}$
is a topological-space.
Additionally, it is obvious that,
\begin{equation}
\base\subseteq\topology{}.
\end{equation}
From this and \Ref{thmbaseScondition1p}, based on \refdef{defbase}
it can be deduced that,
$\base$
is a base for $\opair{\X}{\topology{}}$.
\endthm
%%%%%%%%%%%%%%%%%%%%%%%%%%%%%%%%%%%%%%%%%%%%%%%%%%%%%%%
%%%%%%%%%%%%%%%%%%%%%%%%%%%%%%%%%%%%%%%%%%%%%%%%%%%%%%%%%%%%%%%%%%%%%%%%%%%%%%%%%%%%%%%%%%%%%%%%%%%%%%%
\definition\label{defCbases}
$\X$
is taken as a set.
$\Cbases{\X}$
is defined as,
\begin{align}
\(\base\in\Cbases{\X}\)\thenn
\left\{
\begin{aligned}
&\base\in\CSs{\CSs{\X}},\\
&\X=\bigcup_{\B\in\base}\B,\\
&\Foreach{\(\B_1,\B_2\)}{\base\times\base}\(\Existss{\sC}{\base}\B_1\cap\B_2=\bigcup_{\V\in\sC}\V\).
\end{aligned}
\right.
\end{align}
Every
$\base$
in
$\Cbases{\X}$
is called a $\quotl$base for a topology on $\X$$\quotr$.
\endef
%%%%%%%%%%%%%%%%%%%%%%%%%%%%%%%%%%%%%%%%%%%%%%%%%%%%%%%%%%%%%%%%%%%%%%%%%%%%%%%%%%%%%%%%%%%%%%%%%%%%%%%
\definition\label{defgeneratedtopology}
$\X$
is taken as a set.
\begin{equation}
\Foreach{\base}{\Cbases{\X}}
\topgen{\X}{\base}:=\defset{\U}{\CSs{\X}}{\[\Exists{\sC}{\CSs{\base}}\U=\bigcup_{\V\in\sC}\V\]}.
\end{equation}
$\topgen{\X}{\base}$
is called the $\quotl$derived topology on $\X$ from $\base$$\quotr$ or the
$\quotl$topology on $\X$ generated by $\base$$\quotr$.
\endef
%%%%%%%%%%%%%%%%%%%%%%%%%%%%%%%%%%%%%%%%%%%%%%%%%%%%%%%%%%%%%%%%%%%%%%%%%%%%%%%%%%%%%%%%%%%%%%%%%%%%%%%
\corollary\label{corcbase1}
$\X$
is taken as a set. If $\base$ is a base for a topology on $\X$, then $\base$
is a base for the topological-space $\opair{\X}{\topgen{\X}{\base}}$. That is,
\begin{equation}
\(\base\in\Cbases{\X}\)\then\left\{\base\in\Cbase{\opair{\X}{\topgen{\X}{\base}}}\right\}.
\end{equation}
\endcor
%%%%%%%%%%%%%%%%%%%%%%%%%%%%%%%%%%%%%%%%%%%%%%%%%%%%%%%%%%%%%%%%%%%%%%%%%%%%%%%%%%%%%%%%%%%%%%%%%%%%%%%
\corollary\label{corcbase2}
$\Xt=\opair{\X}{\topology{}}$
is taken as a topological-space. If $\base$ is a base for $\Xt$, then $\base$
is a base for a topology on $\X$, and $\topology{}$ is the same as the derived topology on $\X$
from $\base$. That is,
\begin{align}
\(\base\in\Cbase{\Xt}\)\then
\left\{
\begin{aligned}
&\base\in\Cbases{\X},\\
&\topology{}=\topgen{\X}{\base}.
\end{aligned}
\right.
\end{align}
\endcor
%%%%%%%%%%%%%%%%%%%%%%%%%%%%%%%%%%%%%%%%%%%%%%%%%%%%%%%%%%%%%%%%%%%%%%%%%%%%%%%%%%%%%%%%%%%%%%%%%%%%%%%%%%%
%%%%%%%%%%%%%%%%%%%%%%%%%%%%%%%%%%%%%%%%%%%%%%%%%%%%%%%%%%%%%%%%%%%%%%%%%%%%%%%%%%%%%%%%%%%%%%%%%%%%%%%%%%%
\theorem\label{thmbaserelations1}
$\X$
is taken as a set, and each $\base_1$ and $\base_2$ as an element of $\Cbases{\X}$.
\begin{align}
\(\topgen{\X}{\base_1}\subseteq\topgen{\X}{\base_2}\)\thenn
\left\{\Foreach{\B_1}{\base_1}\[\Foreach{\x}{\B_1}\(\Exists{\B_2}{\base_2}\x\in\B_2\subseteq\B_1\)\]\right\}.
\end{align}
\prooff
\begin{itemize}
\item[${\textbf{\textsf{p1}}}$]
It is assumed that,
$\topgen{\X}{\base_1}\subseteq\topgen{\X}{\base_2}$.
\begin{itemize}
\item[${\textbf{\textsf{p1-1}}}$]
$\B_1$
is taken as an arbitrary element of
$\base_1$
Since,
$\topgen{\X}{\base_1}\subseteq\topgen{\X}{\base_2}$, and,
\begin{equation}
\base_1\subseteq\topgen{\X}{\base_1},
\end{equation}
it is clear that,
\begin{equation}
\base_1\subseteq\topgen{\X}{\base_2},
\end{equation}
and accordingly,
\begin{equation}
\B_1\in\topgen{\X}{\base_2}.
\end{equation}
So since $\base_2$ is a base for $\topgen{\X}{\base_2}$, according to \refthm{thmopensetpointsbase},
\begin{equation}
\Foreach{\x}{\B_1}\(\Exists{\B_2}{\base_2}\x\in\B_2\subseteq\B_1\).
\end{equation}
\endp
\end{itemize}
\item[${\textbf{\textsf{p2}}}$]
It is assumed that,
\begin{equation}\label{thmbaserelations1p2eq1}
\Foreach{\B_1}{\base_1}\[\Foreach{\x}{\B_1}\(\Exists{\B_2}{\base_2}\x\in\B_2\subseteq\B_1\)\].
\end{equation}
\begin{itemize}
\item[${\textbf{\textsf{p2-1}}}$]
$\U$
is taken as an arbitrary element of
$\topgen{\X}{\base_1}$
Since $\base_1$ is a base for $\topgen{\X}{\base_1}$, according to \refdef{defbase},
there exists a subset of $\base_1$, like $\sC$, such that,
\begin{equation}\label{thmbaserelations1p2-1eq1}
\U=\bigcup_{\V\in\sC}\V.
\end{equation}
Since $\sC\subseteq\base_1$, \Ref{thmbaserelations1p2eq1} implies,
\begin{equation}
\Foreach{\V}{\sC}\[\Foreach{\x}{\V}\(\Exists{\B_2}{\base_2}\x\in\B_2\subseteq\V\)\].
\end{equation}
For every $\V$ in $\sC$,
$\f_{\V}$
is taken as such a function that,
\begin{align}
&\f_{\V}\in\Func{\V}{\base_2},\\
&\Foreach{\x}{\V}\x\in\func{\f_{\V}}{\x}\subseteq\V.
\end{align}
Hence it can be easily seen that,
\begin{equation}
\V=\bigcup_{\x\in\V}\func{\f_{\V}}{\x}.
\end{equation}
This and \Ref{thmbaserelations1p2-1eq1} imply,
\begin{equation}
\U=\bigcup_{\V\in\sC}\[\bigcup_{\x\in\V}\func{\f_{\V}}{\x}\],
\end{equation}
and accordingly, by defining $\sCi$ as,
\begin{equation}
\sCi:=\defset{\W}{\CSs{\X}}{\Exists{\V}{\sCi}\[\Exists{\x}{\V}\W=\func{\f_{\V}}{\x}\]},
\end{equation}
it is seen that,
\begin{equation}
\U=\bigcup_{\W\in\sCi}\W.
\end{equation}
It is obvious that,
\begin{equation}
\sCi\subseteq\base_2.
\end{equation}
So, since $\base_2$ is a base for $\topgen{\X}{\base_2}$, base on \refthm{thmopensetpointsbase},
\begin{equation}
\U\in\topgen{\X}{\base_2}.
\end{equation}
\endp
\end{itemize}
Therefore,
\begin{equation}
\topgen{\X}{\base_1}\subseteq\topgen{\X}{\base_2}.
\end{equation}
\endp
\end{itemize}
\endthm
%%%%%%%%%%%%%%%%%%%%%%%%%%%%%%%%%%%%%%%%%%%%%%%%%%%%%%%%%%%%%%%%%%%%%%%%%%%%%%%%%%%%%%%%%%%%%%%%%%%%%%%%%%
\theorem\label{thmbaserelations2}
$\X$
is taken as a set, and each $\base_1$ and $\base_2$
as an element of $\Cbases{\X}$.
\begin{equation}
\topgen{\X}{\base_1}=\topgen{\X}{\base_2},
\end{equation}
if-and-only-if
\begin{align}
\left\{
\begin{aligned}
&\left\{\Foreach{\B_1}{\base_1}\[\Foreach{\x}{\B_1}\(\Exists{\B_2}{\base_2}\x\in\B_2\subseteq\B_1\)\]\right\},\\
&\left\{\Foreach{\B_2}{\base_2}\[\Foreach{\x}{\B_2}\(\Exists{\B_1}{\base_1}\x\in\B_1\subseteq\B_2\)\]\right\}.
\end{aligned}
\right.
\end{align}
\prooff
According to \refthm{thmbaserelations1},
and considering that,
\begin{align}
\(\topgen{\X}{\base_1}=\topgen{\X}{\base_2}\)\thenn
\(\topgen{\X}{\base_1}\subseteq\topgen{\X}{\base_2},~
\topgen{\X}{\base_2}\subseteq\topgen{\X}{\base_1}\),
\end{align}
it is clear.
\endthm
%%%%%%%%%%%%%%%%%%%%%%%%%%%%%%%%%%%%%%%%%%%%%%%%%%%%%%%%%%%%%%%%%%%%%%%%%%%%%%%%%%%%%%%%%%%%%%%%%%%%%%%%%%
\subsection{
Sub-base
}
\definition\label{defsubbase}
$\Xt=\opair{\X}{\topology{}}$
is taken as a topological-space.
$\subbase$ is called a $\quotl$sub-base of the topological-space $\Xt$$\quotr$, iff
it possesses these properties.
\begin{itemize}
\item[${\textbf{\textsf{SB1}}}$]
%\hfill
$\subbase\subseteq\topology{}.$
\item[${\textbf{\textsf{SB2}}}$]
%\hfill
$\displaystyle\Exists{\base}{\Cbase{\Xt}}
\left\{\Foreach{\B}{\base}\[\Existss{\sCii}{\subbase}
\(\CarD{\sCii}\in\Zp,~\B=\bigcap_{\V\in\sCii}\V\)\]\right\}.$
\end{itemize}
\endef
%%%%%%%%%%%%%%%%%%%%%%%%%%%%%%%%%%%%%%%%%%%%%%%%%%%%%%%%%%%%%%%%%%%%%%%%%%%%%%%%%%%%%%
\theorem\label{thmNconditionofsubbase1}
$\Xt=\opair{\X}{\topology{}}$
is taken as a topological-space, and $\subbase$ as a sub-base of $\Xt$.
$\subbase$ covers the set $\X$. That is,
\begin{equation}
\(\bigcup_{\S\in\subbase}\S\)=\X.
\end{equation}
\prooff
According to,
\refdef{defsubbase},
There exists a base for $\Xt$, like $\base$, such that,
\begin{equation}\label{thmNconditionofsubbase1p1}
\Foreach{\B}{\base}\[\Existss{\sCii}{\subbase}
\(\CarD{\sCii}\in\Zp,~\B=\bigcap_{\V\in\sCii}\V\)\].
\end{equation}
\begin{itemize}
\item[${\textbf{\textsf{p1}}}$]
Since
$\subbase\subseteq\topology{}$,
and
$\topology{}\subseteq\CSs{\X}$,
it is clear that
$\subbase\subseteq\CSs{\X}$,
and hence,
\begin{equation}
\(\bigcup_{\S\in\subbase}\S\)\subseteq\X.
\end{equation}
\endp
\end{itemize}
\begin{itemize}
\item[${\textbf{\textsf{p2}}}$]
${}$
\begin{itemize}
\item[${\textbf{\textsf{p2-1}}}$]
$\x$
is taken as an arbitrary element of $\X$. Since,
$\base$ is a base for $\Xt$, according to \refthm{thmbaseNcondition1},
\begin{equation}
\X=\bigcup_{\B\in\base}\B.
\end{equation}
Thus, there exists an element of $\base$, like $\B_{\x}$, such that,
\begin{equation}
\x\in\B_{\x}.
\end{equation}
Additionally, according to \Ref{thmNconditionofsubbase1p1},
there exists a non-empty and finite subset of $\subbase$, like $\sCii_{\x}$,such that,
\begin{equation}
\B_{\x}=\bigcap_{\V\in\sCii_{\x}}\V.
\end{equation}
Therefore,
\begin{equation}
\x\in\bigcap_{\V\in\sCii_{\x}}\V,
\end{equation}
and accordingly,
\begin{equation}
\Foreach{\V}{\sCii_{\x}}x\in\V.
\end{equation}
Since
$\sCii_{\x}\subseteq\subbase$,
\begin{equation}
\Foreach{\V}{\sCii_{\x}}\V\in\subbase.
\end{equation}
Therefore,
\begin{equation}
\Exists{\S}{\subbase}x\in\S,
\end{equation}
and hence,
\begin{equation}
\x\in\bigcup_{\S\in\subbase}\S.
\end{equation}
\endp
\end{itemize}
So,
\begin{equation}
\X\subseteq\(\bigcup_{\S\in\subbase}\S\).
\end{equation}
\endp
\end{itemize}
\endthm
%%%%%%%%%%%%%%%%%%%%%%%%%%%%%%%%%%%%%%%%%%%%%%%%%%%%%%%%%%%%%%%%%%%%%%%%%%%%%%%%%%%%%%%%%%%%%%%%%
\theorem\label{thmsubbaseScondition1}
$\X$
is taken as a set, and
$\subbase$
as such a subset of
$\CSs{\X}$
that covers
$\X$.
That is,
\begin{equation}\label{thmsubbaseScondition1eq1}
\(\bigcup_{\S\in\subbase}\S\)=\X.
\end{equation}
Then,
\begin{align}\label{thmsubbaseScondition1eq2}
\defset{\B}{\CSs{\X}}{\[\Exists{\sCii}{\CSs{\subbase}}
\(\CarD{\sCii}\in\Zp,~\B=\bigcap_{\V\in\sCii}\V\)\]}
\in\Cbases{\X}.
\end{align}
\prooff
$\base$
is defined as,
\begin{equation}\label{thmsubbaseScondition1peq1}
\base:=\defset{\B}{\CSs{\X}}{\[\Exists{\sCii}{\CSs{\subbase}}
\(\CarD{\sCii}\in\Zp,~\B=\bigcap_{\V\in\sCii}\V\)\]}.
\end{equation}
\begin{itemize}
\item[${\textbf{\textsf{p1}}}$]
For every
$\S$
in
$\subbase$,
\begin{align}
\S&\in\CSs{\X},\\
\left\{\S\right\}&\in\CSs{\subbase},\\
\CarD{\left\{\S\right\}}&=1\in\Zp,\\
\S&=\bigcap_{\V\in\left\{\S\right\}}\V.
\end{align}
Thus according to
\Ref{thmsubbaseScondition1peq1},
\begin{equation}
\Foreach{\S}{\subbase}\S\in\base,
\end{equation}
which means,
\begin{equation}
\subbase\subseteq\base.
\end{equation}
Thus,
\begin{equation}
\(\bigcup_{\S\in\subbase}\S\)\subseteq\(\bigcup_{\B\in\base}\B\).
\end{equation}
According to this and
\Ref{thmsubbaseScondition1eq1},
and considering that
$\base\subseteq\CSs{\X}$,
it is clear that,
\begin{equation}
\X=\(\bigcup_{\B\in\base}\B\).
\end{equation}
\endp
\item[${\textbf{\textsf{p2}}}$]
${}$
\begin{itemize}
\item[${\textbf{\textsf{p2-1}}}$]
Each $\B_1$
and
$\B_2$
is taken as an arbitrary element of
$\base$.\\
$\sCii_1$
and
$\sCii_2$
are taken as such subsets of
$\subbase$
that,
\begin{align}
&\CarD{\sCii_{j}}\in\Zp,~j\in\left\{1,2\right\},\\
&\B_{j}=\(\bigcap_{\V\in\sCii_{j}}\V\),~j\in\left\{1,2\right\}.
\end{align}
So,
\begin{align}
\(\sCii_1\cup\sCii_2\)&\in\CSs{\subbase},\\
\CarD{\sCii_1\cup\sCii_2}&\in\Zp,
\end{align}
and
\begin{align}
\B_1\cap\B_2&=\(\bigcap_{\V\in\sCii_{1}}\V\)\cap\(\bigcap_{\V\in\sCii_{2}}\V\)\cr
&=\bigcap_{\V\in\(\sCii_1\cup\sCii_2\)}\V.
\end{align}
According to
\Ref{thmsubbaseScondition1eq1},
these imply,
\begin{equation}
\(\B_1\cap\B_2\)\in\base.
\end{equation}
Thus,
\begin{equation}
\left\{\B_1\cap\B_2\right\}\subseteq\base.
\end{equation}
Additionally, it is clear that,
\begin{equation}
\(\B_1\cap\B_2\)=\(\bigcup_{\V\in\left\{\B_1\cap\B_2\right\}}\V\).
\end{equation}
\endp
\end{itemize}
Therefore,
\begin{equation}
\Foreach{\opair{\B_1}{\B_2}}{\base\times\base}\Existss{\sC}{\base}\(\B_1\cap\B_2=\bigcup_{\U\in\sC}\U\).
\end{equation}
\endp
\end{itemize}
These together with
\refdef{defCbases},
imply
\Ref{thmsubbaseScondition1eq2}.
\endthm
%%%%%%%%%%%%%%%%%%%%%%%%%%%%%%%%%%%%%%%%%%%%%%%%%%%%%%%%%%%%%%%%%%%%%%%%%%%%%%%%%%%%%%%%%%%%%%%%%%%%%%%%%%%%%%%%%%%%%%%%%
\definition
$\X$
is taken as a set.
The set of all covers of
$\X$
is denoted by
$\quotl$$\covers{\X}$$\quotr$.
That is,
$\covers{\X}$
is defined as,
\begin{itemize}
\item[${\textbf{\textsf{COV1}}}$]
%\hfill
$\covers{\X}\subseteq\CSs{\CSs{\X}}.$
\item[${\textbf{\textsf{COV2}}}$]
%hfill
$\displaystyle\(\subbase\in\covers{\X}\)\thenn\[\(\bigcup_{\S\in\subbase}\S\)=\X\].$
\end{itemize}
\endef
%%%%%%%%%%%%%%%%%%%%%%%%%%%%%%%%%%%%%%%%%%%%%%%%%%%%%%%%%%%%%%%%%%%%%%%%%%%%%%%%%%%%%%%%%%%%%%%%%%%%%%%%%%%%%%%%%%%%%%%%%
\definition
$\X$
is taken as a set.
For every
$\subbase$
in
$\covers{\X}$
(Every cover of $\X$, like $\subbase$),
\begin{align}
\basegen{\X}{\subbase}:=\defset{\B}{\CSs{\X}}{\[\Exists{\sCii}{\CSs{\subbase}}
\(\CarD{\sCii}\in\Zp,~\B=\bigcap_{\V\in\sCii}\V\)\]}.
\end{align}
\endef
%%%%%%%%%%%%%%%%%%%%%%%%%%%%%%%%%%%%%%%%%%%%%%%%%%%%%%%%%%%%%%%%%%%%%%%%%%%%%%%%%%%%%%%%%%%%%%%%%%
\corollary
$\X$
is taken as a set, and
$\subbase$
as a cover of
$\X$.
($\subbase\in\covers{\X}$).
\begin{equation}
\basegen{\X}{\subbase}\in\Cbases{\X},
\end{equation}
and
$\subbase$
is a sub-base of the topological-space
$\opair{\X}{\topgen{\X}{\basegen{\X}{\subbase}}}$.
\endcor
%%%%%%%%%%%%%%%%%%%%%%%%%%%%%%%%%%%%%%%%%%%%%%%%%%%%%%%%%%%%%%%%%%%%%%%%%%%%%%%%%%%%%%%%%%%%%%%%%%%%%%%%%%%%%%%%%%%%%%%%%%%%%%%%%%%%%%%%%%%%%%%%%%%%%%%%%
%%%%%%%%%%%%%%%%%%%%%%%%%%%%%%%%%%%%%%%%%%%%%%%%%%%%%%%%%%%%%%%%%%%%%%%%%%%%%%%%%%%%%%%%%%%%%%%%%%%%%%%%%%%%%%%%%%%%%%%%%%%%%%%%%%%%%%%%%%%%%%%%%%%%%%%%%
%%%%%%%%%%%%%%%%%%%%%%%%%%%%%%%%%%%%%%%%%%%%%%%%%%%%%%%%%%%%%%%%%%%%%%%%%%%%%%%%%%%%%%%%%%%%%%%%%%%%%%%%%%%%%%%%%%%%%%%%%%%%%%%%%%%%%%%%%%%%%%%%%%%%%%%%%
%%%%%%%%%%%%%%%%%%%%%%%%%%%%%%%%%%%%%%%%%%%%%%%%%%%%%%%%%%%%%%%%%%%%%%%%%%%%%%%%%%%%%%%%%%%%%%%%%%%%%%%%%%%%%%%%%%%%%%%%%%%%%%%%%%%%%%%%%%%%%%%%%%%%%%%%%
%%%%%%%%%%%%%%%%%%%%%%%%%%%%%%%%%%%%%%%%%%%%%%%%%%%%%%%%%%%%%%%%%%%%%%%%%%%%%%%%%%%%%%%%%%%%%%%%%%%%%%%%%%%%%%%%%%%%%%%%%%%%%%%%%%%%%%%%%%%%%%%%%%%%%%%%%
\section{
Topological Subspaces
}
\definition\label{defsubspacetopology1}
$\X$
is taken as a set,
$\Ss$
as an element of
$\CSs{\CSs{\X}}$
(a collection of subsets of $\X$), and
$\Y$
as a subset of
$\X$.
$\sspower{\Ss}{\Y}$
is defined to be the collection of subsets of
$\Y$,
each element of which equals the intersection of
$\Y$
with an element of
$\Ss$.
That is,
\begin{equation}
\sspower{\Ss}{\Y}:=
\defset{\V}{\CSs{\Y}}{\[\Exists{\U}{\Ss}\(\V=\Y\cap\U\)\]}.
\end{equation}
\endef
%%%%%%%%%%%%%%%%%%%%%%%%%%%%%%%%%%%%%%%%%%%%%%%%%%%%%%%%%%%%%%%%%%%%%%%%%%%%%%%%%%%%%%%%%%%%
\theorem\label{thmsubspacetopology}
$\Xt=\opair{\X}{\topology{}}$
is taken as a topological-space, and
$\Y$
as a subset of
$\X$.
$\opair{\Y}{\stopology{\topology{}}{\Y}}$
is a topological-space.
\proof
According to
\refdef{deftopologicalspace},
\begin{align}
\left\{
\begin{aligned}
&\empty\in\topology{},\\
&\X\in\topology{},\\
&\Foreach{\opair{\U_1}{\U_2}}{\topology{}\times\topology{}}\(\U_1\cap\U_2\)\in\topology{},\\
&\Foreachs{\sC}{\topology{}}\(\bigcup_{\U\in\sC}\U\)\in\topology{}.
\end{aligned}
\right.
\end{align}
Additionally, considering that
$\Y\subseteq\X$,
it is clear that,
$\Y=\Y\cap\X$.
\begin{itemize}
\item[${\textbf{\textsf{p1}}}$]
It is obvious that,
\begin{equation}
\stopology{\topology{}}{\Y}\subseteq\CSs{\Y}.
\end{equation}
\endp
\item[${\textbf{\textsf{p2}}}$]
Considering that,
\begin{align}
\empty&\in\topology{},\\
\empty&=\Y\cap\empty,
\end{align}
According to
\refdef{defsubspacetopology1},
it is clear that,
\begin{equation}
\empty\in\stopology{\topology{}}{\Y}.
\end{equation}
\endp
\item[${\textbf{\textsf{p3}}}$]
Considering that,
\begin{align}
\X&\in\topology{},\\
\Y&=\Y\cap\X,
\end{align}
and according to,
\refdef{defsubspacetopology1},
it is clear that,
\begin{equation}
\Y\in\stopology{\topology{}}{\Y}.
\end{equation}
\endp
\item[${\textbf{\textsf{p4}}}$]
$\opair{\V_1}{\V_2}$
is taken as an arbitrary element of
$\stopology{\topology{}}{\Y}\times\stopology{\topology{}}{\Y}$.\\
Each $\U_1$
and
$\U_2$
is taken as such an element of
$\topology{}$
that,
\begin{align}
\V_1=\Y\cap\U_1,\\
\V_2=\Y\cap\U_2.
\end{align}
So,
\begin{equation}
\(\V_1\cap\V_2\)=\Y\cap\(\U_1\cap\U_2\).
\end{equation}
Considering that,
$\opair{\U_1}{\U_2}\in\topology{}\times\topology{}$,
it is clear that,
\begin{equation}
\(\U_1\cap\U_2\)\in\topology{}.
\end{equation}
These and
\refdef{defsubspacetopology1}
imply
\begin{equation}
\(\V_1\cap\V_2\)\in\stopology{\topology{}}{\Y}.
\end{equation}
\endp
\item[${\textbf{\textsf{p5}}}$]
$\sCi$
is taken as a subset of
$\stopology{\topology{}}{\Y}$.\\
$\f$
is taken as such a function that,
\begin{align}
&\f\in\Func{\sCi}{\topology{}},\label{thmsubspacetopologyp5eq1}\\
&\Foreach{\V}{\sCi}\V=\[\Y\cap\func{\f}{\V}\].\label{thmsubspacetopologyp5eq2}
\end{align}
Thus,
\begin{align}\label{thmsubspacetopologyp5eq3}
\bigcup_{\V\in\sCi}\V&=\bigcup_{\V\in\sCi}\[\Y\cap\func{\f}{\V}\]\cr
&=\Y\cap\[\bigcup_{\V\in\sCi}\func{\f}{\V}\].
\end{align}
Considering that,
\begin{equation}\label{thmsubspacetopologyp5eq4}
\Foreach{\V}{\sCi}\func{\f}{\V}\in\topology{},
\end{equation}
it can be seen that,
\begin{equation}\label{thmsubspacetopologyp5eq5}
\[\bigcup_{\V\in\sCi}\func{\f}{\V}\]\in\topology{}.
\end{equation}
\Ref{thmsubspacetopologyp5eq3},
\Ref{thmsubspacetopologyp5eq5},
and
\refdef{defsubspacetopology1}
imply
\begin{equation}
\(\bigcup_{\V\in\sCi}\V\)\in\stopology{\topology{}}{\Y}.
\end{equation}
\endp
\end{itemize}
Therefore, according to
\refdef{deftopologicalspace},
it is clear that
$\opair{\Y}{\stopology{\topology{}}{\Y}}$
is a topological-space.
\endthm
%%%%%%%%%%%%%%%%%%%%%%%%%%%%%%%%%%%%%%%%%%%%%%%%%%%%%%%%%%%%%%%%%%%%%%%%%%%%%%%%%%%%
\definition
$\opair{\X}{\topology{}}$
is taken as a topological-space.
\begin{itemize}
\item
$\opair{\Y}{\topology{}^{\prime}}$
is called a $\quotl$topological-subspace of $\opair{\X}{\topology{}}$$\quotr$,
and it is written\\
$\opair{\Y}{\topology{}^{\prime}}\subtop\opair{\X}{\topology{}}$,
if-and-only-if
\begin{align}
\left\{
\begin{aligned}
&\Y\subseteq\X,\\
&\topology{}^{\prime}=\stopology{\topology{}}{\Y}.
\end{aligned}
\right.
\end{align}
\item
For every
$\Y$
in
$\CSs{\X}$,
$\stopology{\topology{}}{\Y}$
is called the $\quotl$induced topology on $\Y$ from $\topology{}$$\quotr$.
\end{itemize}
\endef
%%%%%%%%%%%%%%%%%%%%%%%%%%%%%%%%%%%%%%%%%%%%%%%%%%%%%%%%%%%%%%%%%%%%%%%%%%%%%%%%%%%%%%%%%%%%%
\theorem\label{thmsubspaceclosedsets}
$\Xt=\opair{\X}{\topology{}}$
is taken as a topological-space, and
$\Y$
as a subset of
$\X$.
\begin{equation}
\[\V\in\Fclosed{\Y}{\stopology{\topology{}}{\Y}}\]\thenn
\[\V\in\CSs{\Y},~\Exists{\U}{\Fclosed{\X}{\topology{}}}\V=\Y\cap\U\].
\end{equation}
In other words,
\begin{equation}
\Fclosed{\Y}{\stopology{\topology{}}{\Y}}=
\defset{\V}{\CSs{\Y}}{\[\Exists{\U}{\Fclosed{\X}{\topology{}}}\V=\Y\cap\U\]}.
\end{equation}
\prooff
\begin{itemize}
\item[${\textbf{\textsf{p1}}}$]
$\V$
is taken as an arbitrary element of
$\Fclosed{\Y}{\stopology{\topology{}}{\Y}}$.
According to
\refdef{deffamilyofclosedsets},
\begin{align}
\V&\subseteq\Y,\\
\(\compl{\Y}{\V}\)&\in\stopology{\topology{}}{\Y}.
\end{align}
Hence according to
\refdef{defsubspacetopology1},
\begin{equation}
\Exists{\W}{\topology{}}\(\compl{\Y}{\V}\)=\Y\cap\W.
\end{equation}
$\W$
is taken as such an element of
$\topology{}$
that,
\begin{equation}
\(\compl{\Y}{\V}\)=\Y\cap\W.
\end{equation}
Thus, considering that
$\W\subseteq\X$
and
$\V\subseteq\Y\subseteq\X$,
it is clear that
\begin{equation}
\V=\Y\cap\(\compl{\X}{\W}\).
\end{equation}
Considering that
$\W\in\topology{}$,
and according to
\refdef{deffamilyofclosedsets},
it is clear that
\begin{equation}
\(\compl{\X}{\W}\)\in\Fclosed{\X}{\topology{}}.
\end{equation}
These imply,
$\Exists{\U}{\Fclosed{\X}{\topology{}}}\V=\Y\cap\U$.
\endp
\item[${\textbf{\textsf{p2}}}$]
$\V^{\prime}$
is taken as such an element of
$\CSs{\Y}$
that
\begin{equation}
\Exists{\U}{\Fclosed{\X}{\topology{}}}\V^{\prime}=\Y\cap\U.
\end{equation}
$\W^{\prime}$
is taken as such an element of
$\Fclosed{\X}{\topology{}}$
that
\begin{equation}
\V^{\prime}=\Y\cap\W^{\prime}
\end{equation}
So, considering that
$\Y\subseteq\X$,
and
$\W^{\prime}\subseteq\X$,
it is clear that
\begin{equation}
\(\compl{\Y}{\V^{\prime}}\)=\Y\cap\(\compl{\X}{\W^{\prime}}\).
\end{equation}
Considering that
$\W^{\prime}\in\Fclosed{\X}{\topology{}}$,
it is clear that
\begin{equation}
\(\compl{\X}{\W^{\prime}}\)\in\topology{}.
\end{equation}
These imply, according to
\refdef{defsubspacetopology1},
\begin{equation}
\(\compl{\Y}{\V^{\prime}}\)\in\stopology{\topology{}}{\Y},
\end{equation}
and accordingly, considering that
$\V^{\prime}\subseteq\Y$,
and according to
\refdef{deffamilyofclosedsets},
\begin{equation}
\V^{\prime}\in\Fclosed{\Y}{\stopology{\topology{}}{\Y}}.
\end{equation}
\endp
\end{itemize}
\endthm
%%%%%%%%%%%%%%%%%%%%%%%%%%%%%%%%%%%%%%%%%%%%%%%%%%%%%%%%%%%%%%%%%%%%%%%%%%%%%
\theorem\label{thmsubspacetopologybase}
$\Xt=\opair{\X}{\topology{}}$
is taken as a topological-space, and
$\Y$
as a subset of
$\X$. If
$\base$
is a base for
$\Xt$, then
$\sbase{\base}{\Y}$
is a base for the topological-space
$\opair{\Y}{\stopology{\topology{}}{\Y}}$. Tha is,
\begin{equation}\label{thmsubspacetopologybaseeq1}
\Foreach{\base}{\Cbase{\opair{\X}{\topology{}}}}\sbase{\base}{\Y}\in\Cbase{\opair{\Y}{\stopology{\topology{}}{\Y}}}.
\end{equation}
\prooff
\begin{itemize}
\item[${\textbf{\textsf{p}}}$]
$\base$
is taken as a base for
$\opair{\X}{\topology{}}$. According to
\refthm{thmopensetpointsbase},
\begin{equation}\label{thmsubspacetopologybasep1}
\(\U\in\topology{}\)\then
\left\{\Foreach{\x}{\U}\[\Exists{\B}{\base}\(\x\in\B,~\B\subseteq\U\)\]\right\}.
\end{equation}
Additionally, according to
\refdef{defsubspacetopology1},
\begin{align}
\stopology{\topology{}}{\Y}&=\defset{\V}{\CSs{\Y}}{\Exists{\U}{\topology{}}\V=\Y\cap\U},\label{thmsubspacetopologybasep2}\\
\sbase{\base}{\Y}&=\defset{\V}{\CSs{\Y}}{\Exists{\U}{\base}\V=\Y\cap\U}.\label{thmsubspacetopologybasep3}
\end{align}
\begin{itemize}
\item[${\textbf{\textsf{p1}}}$]
According to
\refdef{defbase},
\begin{equation}\label{thmsubspacetopologybasep1eq1}
\base\subseteq\topology{}.
\end{equation}
\Ref{thmsubspacetopologybasep2}
and
\Ref{thmsubspacetopologybasep3}
imply
\begin{equation}\label{thmsubspacetopologybasep1eq2}
\sbase{\base}{\Y}\subseteq\stopology{\topology{}}{\Y}.
\end{equation}
\endp
\item[${\textbf{\textsf{p2}}}$]
$\V$
is taken as an arbitrary element of
$\stopology{\topology{}}{\Y}$. According to
\Ref{thmsubspacetopologybasep2},
\begin{equation}\label{thmsubspacetopologybasep2eq1}
\Exists{\U}{\topology{}}\V=\Y\cap\U.
\end{equation}
$\U$
is taken as such an element of
$\topology{}$ that,
\begin{equation}\label{thmsubspacetopologybasep2eq2}
\V=\Y\cap\U.
\end{equation}
\begin{itemize}
\item[${\textbf{\textsf{p2-1}}}$]
$\x$
is taken as an arbitrary element of
$\V$. According to
\Ref{thmsubspacetopologybasep2eq2},
\begin{align}
\x&\in\Y,\label{thmsubspacetopologybasep2-1eq1}\\
\x&\in\U.\label{thmsubspacetopologybasep2-1eq2}
\end{align}
This and
\Ref{thmsubspacetopologybasep1}
imply,
\begin{equation}\label{thmsubspacetopologybasep2-1eq3}
\Exists{\B}{\base}\(\x\in\B,~\B\subseteq\U\).
\end{equation}
$\B$
is taken as such an element of
$\base$
taht
\begin{align}
\x&\in\B,\label{thmsubspacetopologybasep2-1eq4}\\
\B&\subseteq\U.\label{thmsubspacetopologybasep2-1eq5}
\end{align}
According to
\Ref{thmsubspacetopologybasep2-1eq1}
and
\Ref{thmsubspacetopologybasep2-1eq4},
\begin{equation}\label{thmsubspacetopologybasep2-1eq6}
\x\in\(\Y\cap\B\).
\end{equation}
Additionally, considering that
$\B\in\base$,
and according to
\Ref{thmsubspacetopologybasep3},
\begin{equation}\label{thmsubspacetopologybasep2-1eq7}
\(\Y\cap\B\)\in\sbase{\base}{\Y}.
\end{equation}
According to
\Ref{thmsubspacetopologybasep2-1eq5},
\begin{equation}
\(\Y\cap\B\)\subseteq\(\Y\cap\U\).
\end{equation}
Therefore,
\begin{equation}
\Exists{\B^{\prime}}{\sbase{\base}{\Y}}
\(\x\in\B^{\prime}\subseteq\sbase{\base}{\Y}\).
\end{equation}
\endp
\end{itemize}
\endp
\end{itemize}
So,
\begin{equation}
\Foreach{\V}{\stopology{\topology{}}{\Y}}
\left\{\Foreach{\x}{\V}\[\Exists{\B^{\prime}}{\sbase{\base}{\Y}}
\(\x\in\B^{\prime}\subseteq\sbase{\base}{\Y}\)\]\right\}.
\end{equation}
\Ref{thmsubspacetopologybasep1eq2}
and
\refthm{thmSandNconditionsofbase1}
imply
\begin{equation}
\sbase{\base}{\Y}\in\Cbase{\opair{\Y}{\stopology{\topology{}}{\Y}}}.
\end{equation}
\endp
\end{itemize}
\endthm
%%%%%%%%%%%%%%%%%%%%%%%%%%%%%%%%%%%%%%%%%%%%%%%%%%%%%%%%%%%%%%%%%%%%%%%%%%%%%%%%%%%
\theorem
$\X$
is taken as a set, and
$\Y$
as a subset of
$\X$. If $\base$
is a base for a topology on $\X$, then
$\sbase{\base}{\Y}$
is a base for a topology on $\Y$. That is,
\begin{equation}
\Foreach{\base}{\Cbases{\X}}\(\sbase{\base}{\Y}\in\Cbases{\Y}\).
\end{equation}
\prooff
$\base$
is taken as an arbitrary element of $\Cbases{\X}$.
According to \refcor{corcbase1},
$\base$
is a base for the topological-space $\opair{\X}{\topgen{\X}{\base}}$.
Hence according to \refthm{thmsubspacetopologybase},
$\sbase{\base}{\Y}$
is a base for the topological-space $\opair{\Y}{\stopology{\topgen{\X}{\base}}{\Y}}$.
Thus according to \refcor{corcbase2},
$\sbase{\base}{\Y}$
is a base for a topology on $\Y$.
\endthm
%%%%%%%%%%%%%%%%%%%%%%%%%%%%%%%%%%%%%%%%%%%%%%%%%%%%%%%%%%%%%%%%%%%%%%%%%%%%
\theorem\label{thmsubspacetopologytransitivity}
$\opair{\X}{\topology{}}$
is taken as a topological-space, and each $\Y$ and $\Y^{\prime}$
as a subset of $\X$ such that
$\Y^{\prime}\subseteq\Y$.
the induced topology on $\Y^{\prime}$ from $\stopology{\topology{}}{\Y}$
equals the induced topology on $\Y^{\prime}$ from $\topology{}$.
That is,
\begin{equation}\label{thmsubspacetopologytransitivityeq1}
\stopology{\stopology{\topology{}}{\Y}}{\Y^{\prime}}=\stopology{\topology{}}{\Y^{\prime}}.
\end{equation}
\prooff
According to
\refdef{defsubspacetopology1},
\begin{align}
\stopology{\topology{}}{\Y}&=\defset{\V}{\CSs{\Y}}{\[\Exists{\U}{\topology{}}\(\V=\Y\cap\U\)\]},\label{thmsubspacetopologytransitivitypeq1}\\
\stopology{\stopology{\topology{}}{\Y}}{\Y^{\prime}}&=
\defset{\V}{\CSs{\Y^{\prime}}}{\[\Exists{\U}{\stopology{\topology{}}{\Y}}\(\V=\Y^{\prime}\cap\U\)\]},\label{thmsubspacetopologytransitivitypeq2}\\
\stopology{\topology{}}{\Y^{\prime}}&=
\defset{\V}{\CSs{\Y^{\prime}}}{\[\Exists{\U}{\topology{}}\(\V=\Y^{\prime}\cap\U\)\]}.\label{thmsubspacetopologytransitivitypeq3}
\end{align}
Additionally, considering that
$\Y^{\prime}\subseteq\Y$,
\begin{equation}\label{thmsubspacetopologytransitivitypeq4}
\(\Y^{\prime}\cap\Y\)=\Y^{\prime}.
\end{equation}
\begin{itemize}
\item[${\textbf{\textsf{p1}}}$]
$\V$
is taken as an arbitrary element of
$\stopology{\topology{}}{\Y^{\prime}}$.\\
$\U$
is taken as such an element of $\topology{}$ that,
\begin{equation}
\V=\Y^{\prime}\cap\U.
\end{equation}
So, according to \Ref{thmsubspacetopologytransitivitypeq4},
\begin{equation}
\V=\Y^{\prime}\cap\(\Y\cap\U\).
\end{equation}
Considering that,
$\U\in\topology{}$, and according to
\Ref{thmsubspacetopologytransitivitypeq1},
\begin{equation}
\(\Y\cap\U\)\in\stopology{\topology{}}{\Y},
\end{equation}
These ans \Ref{thmsubspacetopologytransitivitypeq2} imply,
\begin{equation}
\V\in\stopology{\stopology{\topology{}}{\Y}}{\Y^{\prime}}.
\end{equation}
\endp
\item[${\textbf{\textsf{p2}}}$]
$\W$
is taken as an arbitrary element of
$\stopology{\stopology{\topology{}}{\Y}}{\Y^{\prime}}$.
$\U_1$
is taken as such an element of $\stopology{\topology{}}{\Y}$ that,
\begin{equation}
\W=\Y^{\prime}\cap\U_1.
\end{equation}
$\U_2$
is taken as such an element of $\topology{}$ that,
\begin{equation}
\U_1=\Y\cap\U_2.
\end{equation}
Therefore, according to \Ref{thmsubspacetopologytransitivitypeq4},
\begin{align}
\W&=\Y^{\prime}\cap\(\Y\cap\U\)\cr
&=\(\Y'\cap\Y\)\cap\U\cr
&=\Y^{\prime}\cap\U.
\end{align}
Hence considering that,
$\U\in\topology{}$,
and according to
\Ref{thmsubspacetopologytransitivitypeq3},
\begin{equation}
\W\in\stopology{\topology{}}{\Y^{\prime}}.
\end{equation}
\endp
\end{itemize}
\endthm
%%%%%%%%%%%%%%%%%%%%%%%%%%%%%%%%%%%%%%%%%%%%%%%%%%%%%%%%%%%%%%%%%
\corollary
$\opair{\X}{\topology{}}$
is taken as a topological-space, and
$\Y_1$
and
$\Y_2$
are taken as subsets of $\X$.
If
$\opair{\Y_2}{\topology{2}}$
is a topological-subspace of
$\opair{\Y_1}{\topology{1}}$,
and
$\opair{\Y_1}{\topology{1}}$
is a topological-subspace of
$\opair{\X}{\topology{}}$,
then
$\opair{\Y_2}{\topology{2}}$
is a topological-subspace of
$\opair{\X}{\topology{}}$.
That is,
\begin{equation}
\[\opair{\Y_2}{\topology{2}}\subtop\opair{\Y_1}{\topology{1}},~
\opair{\Y_1}{\topology{1}}\subtop\opair{\X}{\topology{}}\]\then
\[\opair{\Y_2}{\topology{2}}\subtop\opair{\X}{\topology{}}\].
\end{equation}
\endcor
%%%%%%%%%%%%%%%%%%%%%%%%%%%%%%%%%%%%%%%%%%%%%%%%%%%%%%%%%%%%%%%%%%%%%%%%%%%
\theorem\label{thmsubsubopen}
$\opair{\X}{\topology{}}$
is taken as a topological-space,
$\Y$
as a subset of
$\X$,
and
$\V$
as a subset of
$\Y$.
If
$\V$
is an open set of
$\opair{\Y}{\stopology{\topology{}}{\Y}}$,
and
$\Y$
is an open set of
$\opair{\X}{\topology{}}$,
then
$\V$
is an open set of
$\opair{\X}{\topology{}}$.
That is,
\begin{equation}
\[\V\in\stopology{\topology{}}{\Y},~\Y\in\topology{}\]\then
\(\V\in\topology{}\).
\end{equation}
\prooff
It is assumed that,
\begin{align}
\V&\in\stopology{\topology{}}{\Y},\label{thmsubsubopenp1}\\
\Y&\in\topology{}.\label{thmsubsubopenp2}
\end{align}
According to
\Ref{thmsubsubopenp1},
and
\refdef{defsubspacetopology1},
\begin{equation}\label{thmsubsubopenp3}
\Exists{\U}{\topology{}}\V=\Y\cap\U.
\end{equation}
This and
\Ref{thmsubsubopenp2}, and \refdef{deftopologicalspace} imply,
\begin{equation}
\V\in\topology{}.
\end{equation}
\endthm
%%%%%%%%%%%%%%%%%%%%%%%%%%%%%%%%%%%%%%%%%%%%%%%%%%%%%%%%%%%%%%%%%%%%%%%%%%%%%%%%%%%%%
\theorem\label{thmsubsubclosed}
$\opair{\X}{\topology{}}$
is taken as a topological-space,
$\Y$
as a subset of
$\X$,
and
$\V$
as a subset of
$\Y$.
If
$\V$
is a closed set of
$\opair{\Y}{\stopology{\topology{}}{\Y}}$,
and
$\Y$
is a closed set of
$\opair{\X}{\topology{}}$,
then
$\V$
is a closed set of
$\opair{\X}{\topology{}}$.
That is,
\begin{equation}
\[\V\in\Fclosed{\Y}{\stopology{\topology{}}{\Y}},~\Y\in\Fclosed{\X}{\topology{}}\]\then
\(\V\in\Fclosed{\X}{\topology{}}\).
\end{equation}
\prooff
It is assumed that,
\begin{align}
\V&\in\Fclosed{\Y}{\stopology{\topology{}}{\Y}},\label{thmsubsubclosedp1}\\
\Y&\in\Fclosed{\X}{\topology{}}.\label{thmsubsubclosedp2}
\end{align}
According to
\Ref{thmsubsubclosedp1},
and
\refthm{thmsubspaceclosedsets},
\begin{equation}\label{thmsubsubclosedp3}
\Exists{\U}{\Fclosed{\X}{\topology{}}}\V=\Y\cap\U.
\end{equation}
This,
\Ref{thmsubsubclosedp2},
and \refthm{thmclosedsets} imply that,
\begin{equation}
\V\in\Fclosed{\X}{\topology{}}.
\end{equation}
\endthm
%%%%%%%%%%%%%%%%%%%%%%%%%%%%%%%%%%%%%%%%%%%%%%%%%%%%%%%%%%%%%%%%%%%%%%%%%%%%%%%%%%%%%%%%%%%%%%%%%%%%%%%%%%%%%%%%%%%%%%%%%%%%%%%%%%%%%%%%%%%%%%%%%%%%%%%%%
\theorem\label{thmbaseofopensubspace}
$\opair{\X}{\topology{}}$
is taken as a topological-space, and
$\base$
as an element of
$\Cbase{\Xt}$
(a base for the topological-space $\Xt$).
For every open set of
$\Xt$,
like $\U$,
the set of all such elements of $\base$ included in $\U$,
is a base for the topological-space
$\opair{\U}{\stopology{\topology{}}{\U}}$. That is,
\begin{equation}
\Foreach{\U}{\topology{}}
\bigg(\defset{\B}{\base}{\B\subseteq\U}\in
\Cbase{\opair{\U}{\stopology{\topology{}}{\U}}}\bigg).
\end{equation}
\prooff
$\U$
is taken as an arbitrary element of
$\topology{}$.
Considering that
$\base\subseteq\topology{}$
(\refdef{defbase}),
\begin{equation}\label{thmbaseofopensubspacepeq1}
\Foreach{\asubset}{\defset{\B}{\base}{\B\subseteq\U}}
\asubset\in\topology{},
\end{equation}
and hence according to
\refdef{defsubspacetopology1},
\begin{equation}\label{thmbaseofopensubspacepeq2}
\Foreach{\asubset}{\defset{\B}{\base}{\B\subseteq\U}}
\(\asubset\cap\U\)\in\stopology{\topology{}}{\U}.
\end{equation}
In addition, it is clear that,
\begin{equation}\label{thmbaseofopensubspacepeq3}
\Foreach{\asubset}{\defset{\B}{\base}{\B\subseteq\U}}
\(\asubset\cap\U\)=\asubset.
\end{equation}
According to
\Ref{thmbaseofopensubspacepeq2} and
\Ref{thmbaseofopensubspacepeq3},
\begin{equation}\label{thmbaseofopensubspacepeq4}
\defset{\B}{\base}{\B\subseteq\U}\subseteq
\stopology{\topology{}}{\U}.
\end{equation}
\begin{itemize}
\item[${\textbf{\textsf{p1}}}$]
$\V$
is taken as an arbitrary element of
$\stopology{\topology{}}{\U}$
(an open set of the topological-space
$\opair{\U}{\stopology{\topology{}}{\U}}$)
Then, considering that
$\U$
is an open set of
$\Xt$,
\refthm{thmsubsubopen},
implies
$\V$
is also an open set of the topological-space
$\Xt$.
\begin{equation}\label{thmbaseofopensubspacep1eq1}
\V\in\topology{}.
\end{equation}
Thus according to
\refdef{defbase},
\begin{equation}\label{thmbaseofopensubspacep1eq2}
\Existsis{\sC}{\CSs{\base}}
\V=\(\union{\sC}\).
\end{equation}
Hence it is clear that,
\begin{equation}\label{thmbaseofopensubspacep1eq3}
\Foreach{\asubset}{\sC}
\bigg(\asubset\in\base,~\asubset\subseteq\V\bigg),
\end{equation}
and hence considering that,
$\V\subseteq\U$,
\begin{equation}\label{thmbaseofopensubspacep1eq4}
\sC\subseteq\defset{\B}{\base}{\B\subseteq\U}.
\end{equation}
\Ref{thmbaseofopensubspacep1eq2}
and
\Ref{thmbaseofopensubspacep1eq4}
imply,
\begin{equation}
\Existsis{\sC}{\CSs{\defset{\B}{\base}{\B\subseteq\U}}}\V=\(\union{\sC}\).
\end{equation}
\endp
\end{itemize}
Therefore,
\begin{equation}\label{thmbaseofopensubspacepeq5}
\Foreach{\V}{\stopology{\topology{}}{\U}}
\bigg[\Existsis{\sC}{\CSs{\defset{\B}{\base}{\B\subseteq\U}}}\V=\(\union{\sC}\)\bigg].
\end{equation}
According to
\refdef{defbase},
\Ref{thmbaseofopensubspacepeq4}
and
\Ref{thmbaseofopensubspacepeq5}
imply
$\defset{\B}{\base}{\B\subseteq\U}$
is a base for
$\opair{\U}{\stopology{\topology{}}{\U}}$. That is,
\begin{equation}
\defset{\B}{\base}{\B\subseteq\U}\in
\Cbase{\opair{\U}{\stopology{\topology{}}{\U}}}.
\end{equation}
\endthm
%thmSandNconditionsofbase0
%%%%%%%%%%%%%%%%%%%%%%%%%%%%%%%%%%%%%%%%%%%%%%%%%%%%%%%%%%%%%%%%%%%%%%%%%%%%%%%%%%%%%%%%%%%%%%%%%%%%%%%%%%%%%%%%%%%%%%%%%%%%%%%%%%%%%%%%%%%%%%%%%%%%%%%%%
\theorem\label{thmsubspacesofdiscretespace}
$\X$
is taken as a set.
For every subset
$\asubset$
of
$\X$,
$\opair{\asubset}{\CSs{\asubset}}$
is a topological-subspace of the discrete topological-space
$\opair{\X}{\CSs{\X}}$.
\begin{equation}
\Foreach{\asubset}{\CSs{\X}}
\opair{\asubset}{\CSs{\asubset}}\subtop\opair{\X}{\CSs{\X}}.
\end{equation}
In other words,
\begin{equation}
\Foreach{\asubset}{\CSs{\X}}
\stopology{\CSs{\X}}{\asubset}=\CSs{\asubset}.
\end{equation}
$\caution$
This means, every topological-subspace of a discrete topological-space is a discrete topological-space.
\proof
$\asubset$
is taken as an arbitrary element of
$\CSs{\X}$.
According to,
\refdef{defsubspacetopology1},
\begin{align}
\stopology{\CSs{\X}}{\asubset}&=
\defset{\V}{\CSs{\asubset}}
{\[\Exists{\U}{\CSs{\X}}\V=\U\cap\asubset\]}\cr
&=\CSs{\asubset}.
\end{align}
\endthm
%%%%%%%%%%%%%%%%%%%%%%%%%%%%%%%%%%%%%%%%%%%%%%%%%%%%%%%%%%%%%%%%%%%%%%%%%%%%%%%%%%%%%%%%%%%%%%%%%%%%%%%%%%%%%%%%%%%%%%%%%%%%%%%%%%%%%%%%%%%%%%%%%%%%%%%%%
%%%%%%%%%%%%%%%%%%%%%%%%%%%%%%%%%%%%%%%%%%%%%%%%%%%%%%%%%%%%%%%%%%%%%%%%%%%%%%%%%%%%%%%%%%%%%%%%%%%%%%%%%%%%%%%%%%%%%%%%%%%%%%%%%%%%%%%%%%%%%%%%%%%%%%%%%
%%%%%%%%%%%%%%%%%%%%%%%%%%%%%%%%%%%%%%%%%%%%%%%%%%%%%%%%%%%%%%%%%%%%%%%%%%%%%%%%%%%%%%%%%%%%%%%%%%%%%%%%%%%%%%%%%%%%%%%%%%%%%%%%%%%%%%%%%%%%%%%%%%%%%%%%%
\subsection{
Singleton Topological Subspaces}
\theorem\label{thmsingletonsubspacetopology}
$\opair{\X}{\topology{}}$
is taken as a topological-space.
For every
$\x$
in
$\X$,
the topological-subspace
$\opair{\seta{\x}}{\stopology{\topology{}}{\seta{\x}}}$
of
$\Xt$
equals
$\singletonTS{\x}$
(The singleton topological-space of the point $\x$). That is,
\begin{align}
\Foreach{\x}{\X}
\stopology{\topology{}}{\seta{\x}}&=\CSs{\seta{\x}}\cr
&=\seta{\binary{\empty}{\seta{\x}}}.
\end{align}
\prooff
$\x$
is taken as an arbitrary element of
$\X$
According to,
\refdef{defsubspacetopology1},
\begin{equation}
\stopology{\topology{}}{\seta{\x}}\subseteq\CSs{\seta{\x}}.
\end{equation}
In addition, according to,
\refdef{deftopologicalspace}
($\empty$
and
$\X$
are elements of
$\topology{}$),
and
\refdef{defsubspacetopology1},
\begin{align}
\CSs{\seta{\x}}&=\seta{\binary{\empty}{\seta{\x}}}\cr
&=\seta{\binary{\empty\cap\seta{\x}}{\X\cap\seta{\x}}}\cr
&\subseteq\stopology{\topology{}}{\seta{\x}}.
\end{align}
Thus,
\begin{equation}
\stopology{\topology{}}{\seta{\x}}=\CSs{\seta{\x}}.
\end{equation}
\endthm
%%%%%%%%%%%%%%%%%%%%%%%%%%%%%%%%%%%%%%%%%%%%%%%%%%%%%%%%%%%%%%%%%%%%%%%%%%%%%%%%%%%%%%%%%%%%%%%%%%%%%%%%%%%%%%%%%%%%%%%%%%%%%%%%%%%%%%%%%%%%%%%%%%%%%%%%%
%thmclosedsets
%thmsubspaceclosedsets
%%%%%%%%%%%%%%%%%%%%%%%%%%%%%%%%%%%%%%%%%%%%%%%%%%%%%%%%%%%%%%%%%%%%%%%%%%%%%%%%%%%%%%%%%%%%%%%%%%%%%%%%%%%%%%%%%%%%%%%%%%%%%%%%%%%%%%%%%%%%%%%%%%%%%%%%%
%%%%%%%%%%%%%%%%%%%%%%%%%%%%%%%%%%%%%%%%%%%%%%%%%%%%%%%%%%%%%%%%%%%%%%%%%%%%%%%%%%%%%%%%%%%%%%%%%%%%%%%%%%%%%%%%%%%%%%%%%%%%%%%%%%%%%%%%%%%%%%%%%%%%%%%%%
%%%%%%%%%%%%%%%%%%%%%%%%%%%%%%%%%%%%%%%%%%%%%%%%%%%%%%%%%%%%%%%%%%%%%%%%%%%%%%%%%%%%%%%%%%%%%%%%%%%%%%%%%%%%%%%%%%%%%%%%%%%%%%%%%%%%%%%%%%%%%%%%%%%%%%%%%
%%%%%%%%%%%%%%%%%%%%%%%%%%%%%%%%%%%%%%%%%%%%%%%%%%%%%%%%%%%%%%%%%%%%%%%%%%%%%%%%%%%%%%%%%%%%%%%%%%%%%%%%%%%%%%%%%%%%%%%%%%%%%%%%%%%%%%%%%%%%%%%%%%%%%%%%%
%%%%%%%%%%%%%%%%%%%%%%%%%%%%%%%%%%%%%%%%%%%%%%%%%%%%%%%%%%%%%%%%%%%%%%%%%%%%%%%%%%%%%%%%%%%%%%%%%%%%%%%%%%%%%%%%%%%%%%%%%%%%%%%%%%%%%%%%%%%%%%%%%%%%%%%%%
\section{
Point and Set Positions
}
\definition\label{definteriorpoint}
$\Xt=\opair{\X}{\topology{}}$
is taken as a topological-space, and
$\asubset$
as a subset of
$\X$.
\begin{align}
\Sint{\Xt}{\asubset}:&=
\defset{\point}{\X}{\(\Exists{\U}{\func{\nei{\Xt}}{\left\{\point\right\}}}\U\subseteq\asubset\)}\cr
&=\defset{\point}{\X}{\(\Exists{\U}{\[\topology{}\cap\CSs{\asubset}\]}\point\in\U\)}.
\end{align}
$\x$
is referred to as an $\quotl$interior-point of $\asubset$ in the topological-space $\Xt$$\quotr$
iff
$\x\in\Sint{\Xt}{\asubset}$.
\endef
%%%%%%%%%%%%%%%%%%%%%%%%%%%%%%%%%%%%%%%%%%%%%%%%%%%%%%%%%%%%
\definition\label{defexteriorpoint}
$\Xt=\opair{\X}{\topology{}}$
is taken as a topological-space, and
$\asubset$
as a subset of
$\X$.
\begin{align}
\Sext{\Xt}{\asubset}:&=\defset{\point}{\X}
{\(\Exists{\U}{\func{\nei{\Xt}}{\left\{\point\right\}}}\U\cap\asubset=\empty\)}\cr
&=\defset{\point}{\X}{\(\Exists{\U}{\[\topology{}\cap\CSs{\compl{\X}{\asubset}}\]}\point\in\U\)}.
\end{align}
$\x$
is reffered to as an $\quotl$exterior-point of $\asubset$ in the topological-space $\Xt$$\quotr$
iff
$\x\in\Sext{\Xt}{\asubset}$.
\endef
%%%%%%%%%%%%%%%%%%%%%%%%%%%%%%%%%%%%%%%%%%%%%%%%%%%%%%%%%%%%
%%%%%%%%%%%%%%%%%%%%%%%%%%%%%%%%%%%%%%%%%%%%%%%%%%%%%%%%%%%%
\definition\label{defboundarypoint}
$\Xt=\opair{\X}{\topology{}}$
is taken as a topological-space, and
$\asubset$
as a subset of
$\X$.
\begin{align}
\Sbound{\Xt}{\asubset}:=\defset{\point}{\X}
{\(\Foreach{\U}{\func{\nei{\Xt}}{\left\{\point\right\}}}
\[\U\cap\asubset\neq\empty,~\U\cap\(\compl{\X}{\asubset}\)\neq\empty\]\)}.
\end{align}
$\x$
is reffered to as a $\quotl$boundary-point of $\asubset$ in the topological-space $\Xt$$\quotr$
iff
$\x\in\Sbound{\Xt}{\asubset}$.
\endef
%%%%%%%%%%%%%%%%%%%%%%%%%%%%%%%%%%%%%%%%%%%%%%%%%%%%%%%%%%%%
\definition\label{defadherentpoint}
$\Xt=\opair{\X}{\topology{}}$
is taken as a topological-space, and
$\asubset$
as a subset of
$\X$.
\begin{equation}
\Sadh{\Xt}{\asubset}:=\defset{\point}{\X}
{\[\Foreach{\U}{\func{\nei{\Xt}}{\left\{\point\right\}}}
\U\cap\asubset\neq\empty\]}.
\end{equation}
$\x$
is referred to as an $\quotl$adherent-point of $\asubset$ in the topological-space $\Xt$$\quotr$
iff
$\x\in\Sadh{\Xt}{\asubset}$.
\endef
%%%%%%%%%%%%%%%%%%%%%%%%%%%%%%%%%%%%%%%%%%%%%%%%%%%%%%%%%%%%
\definition\label{deflimitpoint}
$\Xt=\opair{\X}{\topology{}}$
is taken as a topological-space, and
$\asubset$
as a subset of
$\X$.
\begin{equation}
\Slim{\Xt}{\asubset}:=\defset{\point}{\X}
{\[\Foreach{\U}{\func{\nei{\Xt}}{\left\{\point\right\}}}
\U\cap\(\compl{\asubset}{\left\{\point\right\}}\)\neq\empty\]}.
\end{equation}
$\x$
is referred to as a $\quotl$limit-point of $\asubset$ in the topological-space $\Xt$$\quotr$
iff
$\x\in\Slim{\Xt}{\asubset}$.
\endef
%%%%%%%%%%%%%%%%%%%%%%%%%%%%%%%%%%%%%%%%%%%%%%%%%%%%%%%%%%%%
\definition\label{defisolatedpoint}
$\Xt=\opair{\X}{\topology{}}$
is taken as a topological-space, and
$\asubset$
as a subset of
$\X$.
\begin{equation}
\Siso{\Xt}{\asubset}:=\defset{\point}{\X}
{\[\Exists{\U}{\func{\nei{\Xt}}{\left\{\point\right\}}}
\U\cap\asubset=\seta{\point}\]}.
\end{equation}
$\x$
is referred to as an $\quotl$isolated-point of $\asubset$ in the topological-space $\Xt$$\quotr$
iff
$\x\in\Siso{\Xt}{\asubset}$.
\endef
%%%%%%%%%%%%%%%%%%%%%%%%%%%%%%%%%%%%%%%%%%%%%%%%%%%%%%%%%%%%
\definition\label{defdenseinitselfset}
$\Xt=\opair{\X}{\topology{}}$
is taken as a topological-space, and
$\asubset$
a subset of
$\X$.
$\asubset$
is referred to as a $\quotl$dense-in-itself set of $\Xt$$\quotr$
iff
\begin{equation}
\Siso{\Xt}{\asubset}=\empty.
\end{equation}
\endef
%%%%%%%%%%%%%%%%%%%%%%%%%%%%%%%%%%%%%%%%%%%%%%%%%%%%%%%%%%%%
\corollary\label{corSintSext}
$\Xt=\opair{\X}{\topology{}}$
is taken as a topological-space, and
$\asubset$
a subset of
$\X$.
$\point$
is an exterior-point of
$\asubset$
in the topological-space
$\Xt$, if-and-only-if
$\point$
is an interior-point of
$\compl{\X}{\asubset}$
in the topological-space
$\Xt$. That is,
\begin{equation}
\Sext{\Xt}{\asubset}=\Sint{\Xt}{\compl{\X}{\asubset}}.
\end{equation}
\endcor
%%%%%%%%%%%%%%%%%%%%%%%%%%%%%%%%%%%%%%%%%%%%%%%%%%%%%%%%%%%%
\corollary\label{corSextSadh}
$\Xt=\opair{\X}{\topology{}}$
is taken as a topological-space, and
$\asubset$
a subset of
$\X$.
For every,
$\point$
in
$\X$,
$\point$
is an adherent point of
$\asubset$
in the topological-space
$\Xt$, if-and-only-if
$\point$
is not an exterior-point of
$\asubset$
in the topological-space
$\Xt$. That is
\begin{equation}
\Sadh{\Xt}{\asubset}=\compl{\X}{\(\Sext{\Xt}{\asubset}\)}.
\end{equation}
\endcor
%%%%%%%%%%%%%%%%%%%%%%%%%%%%%%%%%%%%%%%%%%%%%%%%%%%%%%%%%%%%
\corollary\label{correlationofpointsetpositions}
$\Xt=\opair{\X}{\topology{}}$
is taken as a topological-space, and
$\asubset$
a subset of
$\X$.
\begin{itemize}
\item
For every
$\point$
in
$\X$,
if
$\point$
is an interior-point of
$\asubset$
in the topological-space
$\Xt$,
then
$\point$
is an adherent point of
$\asubset$
in the topological-space
$\Xt$. That is,
\begin{equation}
\Sint{\Xt}{\asubset}\subseteq\Sadh{\Xt}{\asubset}.
\end{equation}
\item
For every
$\point$
in
$\X$,
if
$\point$
is a boundary-point of
$\asubset$
in the topological-space
$\Xt$, then
$\point$
is an adherent-point of
$\asubset$
in the topological-space
$\Xt$. That is,
\begin{equation}
\Sbound{\Xt}{\asubset}\subseteq\Sadh{\Xt}{\asubset}.
\end{equation}
\item
For every
$\point$
in
$\X$,
if
$\point$
is a limit-point of
$\asubset$
in
$\Xt$, then
$\point$
is an adherent-point of
$\asubset$
in
$\Xt$. That is,
\begin{equation}
\Slim{\Xt}{\asubset}\subseteq\Sadh{\Xt}{\asubset}.
\end{equation}
\item
For every
$\point$
in
$\X$,
if
$\point$
is an adherent point of
$\asubset$
in
$\Xt$,
and
$\point$
is not in
$\asubset$, then
$\point$
is a limit-point of
$\asubset$
in
$\Xt$. That is
\begin{equation}
\Sadh{\Xt}{\asubset}\cap\(\compl{\X}{\asubset}\)\subseteq\Slim{\Xt}{\asubset}.
\end{equation}
\end{itemize}
\endcor
%%%%%%%%%%%%%%%%%%%%%%%%%%%%%%%%%%%%%%%%%%%%%%%%%%%%%%%%%%%%%%%%%%%%%%%%%%%%%%%%%%%%%
\corollary\label{corlimitpointsandisolatedpointsaredistinct}
$\Xt=\opair{\X}{\topology{}}$
is taken as a topological-space, and
$\asubset$
a subset of
$\X$.
\begin{equation}
\Slim{\Xt}{\asubset}\cap\Siso{\Xt}{\asubset}=\empty.
\end{equation}
In addition,
\begin{equation}
\(\asubset\cap\Slim{\Xt}{\asubset}\)\sqcup\Siso{\Xt}{\asubset}=\asubset.
\end{equation}
\endcor
%%%%%%%%%%%%%%%%%%%%%%%%%%%%%%%%%%%%%%%%%%%%%%%%%%%%%%%%%%%%%%%%%%%%%%%%%%%%%%%%%%%%%
%%%%%%%%%%%%%%%%%%%%%%%%%%%%%%%%%%%%%%%%%%%%%%%%%%%%%%%%%%%%%%%%%%%%%%%%%%%%%%%%%%%%%
%%%%%%%%%%%%%%%%%%%%%%%%%%%%%%%%%%%%%%%%%%%%%%%%%%%%%%%%%%%%%%%%%%%%%%%%%%%%%%%%%%%%%
\subsection{
Interior of a Subset of a Topological Space}
\theorem\label{thmpreint0}
$\Xt=\opair{\X}{\topology{}}$
is taken as a topological-space, and
$\asubset$
as a subset of
$\X$.
There exists only one open set of $\Xt$
in
$\CSs{\asubset}$
that includes every open set of $\Xt$
in
$\CSs{\asubset}$. That is,
\begin{equation}
\Existsu{\U}{\[\topology{}\cap\CSs{\asubset}\]}
\left\{\Foreach{\V}{\[\topology{}\cap\CSs{\asubset}\]}\V\subseteq\U\right\}.
\end{equation}
More precisely,
\begin{align}
\left\{\U\in\[\topology{}\cap\CSs{\asubset}\],~
\Foreach{\V}{\[\topology{}\cap\CSs{\asubset}\]}\U\subseteq\V\right\}\thenn
\left\{\U=\union{\[\topology{}\cap\CSs{\asubset}\]}\right\}.
\end{align}
\prooff
The set
$\U$ is defined as,
\begin{equation}\label{thmpreint0peq1}
\U:=\union{\[\topology{}\cap\CSs{\asubset}\]}.
\end{equation}
\begin{itemize}
\item[${\textbf{\textsf{p1}}}$]
Considering that
$\[\topology{}\cap\CSs{\asubset}\]\subseteq\topology{}$,
and according to
\refdef{deftopologicalspace},
it is clear that
\begin{equation}
\union{\[\topology{}\cap\CSs{\asubset}\]}\in\topology{}.
\end{equation}
In addition, considering that
$\[\topology{}\cap\CSs{\asubset}\]\subseteq\CSs{\asubset}$,
it is clear that
\begin{equation}
\union{\[\topology{}\cap\CSs{\asubset}\]}\in\CSs{\asubset}.
\end{equation}
According to these,
\begin{equation}
\union{\[\topology{}\cap\CSs{\asubset}\]}\in\[\topology{}\cap\CSs{\asubset}\].
\end{equation}
it is trivial that,
\begin{equation}
\Foreach{\V}{\[\topology{}\cap\CSs{\asubset}\]}
\V\subseteq\union{\[\topology{}\cap\CSs{\asubset}\]},
\end{equation}
Therefore,
\begin{equation}
\Existsis{\U}{\[\topology{}\cap\CSs{\asubset}\]}
\[\Foreach{\V}{\[\topology{}\cap\CSs{\asubset}\]}\V\subseteq\U\].
\end{equation}
\endp
\item[${\textbf{\textsf{p2}}}$]
$\U^{\prime}$
is taken as such an element of
$\[\topology{}\cap\CSs{\asubset}\]$
that,
\begin{equation}
\Foreach{\V}{\[\topology{}\cap\CSs{\asubset}\]}\V\subseteq\U^{\prime}.
\end{equation}
Then, it is trivial that,
\begin{equation}
\left\{
\begin{aligned}
\U^{\prime}\subseteq\U,\\
\U\subseteq\U^{\prime},
\end{aligned}
\right.
\end{equation}
which means,
\begin{equation}
\U^{\prime}=\U.
\end{equation}
\endp
\end{itemize}
\endthm
%%%%%%%%%%%%%%%%%%%%%%%%%%%%%%%%%%%%%%%%%%%%%%%%%%%%%%%%%%%%%%%%%%%%%%%%%%%%%%%%%%%%
\definition\label{defintofset}
$\Xt=\opair{\X}{\topology{}}$
is taken as a topological space. The mapping
$\Int{\Xt}$
is defined as this.
\begin{itemize}
\item[${\textbf{\textsf{int1}}}$]
\hfill
$\Int{\Xt}\indef\Func{\CSs{\X}}{\CSs{\X}}.$
\item[${\textbf{\textsf{int2}}}$]
\hfill
$\Foreach{\asubset}{\CSs{\X}}
\func{\Int{\Xt}}{\asubset}\eqdef\union{\[\topology{}\cap\CSs{\asubset}\]}.$
\end{itemize}
For every
$\asubset$
in
$\CSs{\X}$,
$\func{\Int{\Xt}}{\asubset}$
is referred to as the $\quotl$interior of $\asubset$ in the topological-space $\Xt$$\quotr$.
In other words, for every
$\asubset$
in
$\CSs{\X}$,
the union of all open sets of
$\Xt$
that are subsets of
$\asubset$,
is reffered to as the $\quotl$interior of $\asubset$
in the topological-space $\Xt$$\quotr$.
\endef
%%%%%%%%%%%%%%%%%%%%%%%%%%%%%%%%%%%%%%%%%%%%%%%%%%%%%%%%%%%%%%%%%%%%%%%%%%%%%%%%%%%%%
\corollary\label{corintofset0}
$\Xt=\opair{\X}{\topology{}}$
is taken as a topological-space.
\begin{align}
\Foreach{\asubset}{\CSs{\X}}\left\{
\begin{aligned}
&\func{\Int{\Xt}}{\asubset}\in\topology{},\\
&\func{\Int{\Xt}}{\asubset}\subseteq\asubset,\\
&\Foreach{\U}{\[\topology{}\cap\CSs{\asubset}\]}\U\subseteq\func{\Int{\Xt}}{\asubset}.
\end{aligned}
\right.
\end{align}
In other words, for every
$\asubset$
in
$\CSs{\X}$,
$\func{\Int{\Xt}}{\asubset}$
is an open set of $\Xt$,
a subset of $\asubset$
including every open set of $\Xt$,
and a subset of $\asubset$.
In addition, for every
$\asubset$
in
$\X$,
such subset of
$\X$
is unique.
\endcor
%%%%%%%%%%%%%%%%%%%%%%%%%%%%%%%%%%%%%%%%%%%%%%%%%%%%%%%%%%%%%%%%%%%%%%%%%%%%%%%%%%%%%
\theorem\label{thmintofemptyset}
$\Xt=\opair{\X}{\topology{}}$
is taken as a topological-space.
\begin{equation}
\func{\Int{\Xt}}{\empty}=\empty.
\end{equation}
\prooff
Considering that,
\begin{align}
\CSs{\empty}&=\seta{\empty},\\
\empty&\in\topology{},
\end{align}
it is clear that,
\begin{equation}
\[\topology{}\cap\CSs{\empty}\]=\seta{\empty}.
\end{equation}
So, according to \refdef{defintofset},
\begin{align}
\func{\Int{\Xt}}{\empty}&=\union{\[\topology{}\cap\CSs{\empty}\]}\cr
&=\union{\seta{\empty}}\cr
&=\empty.
\end{align}
\endthm
%%%%%%%%%%%%%%%%%%%%%%%%%%%%%%%%%%%%%%%%%%%%%%%%%%%%%%%%%%%%%%%%%%%%%%%%%%%%%%%%%%%%%
\theorem\label{thmintofuniversalset}
$\Xt=\opair{\X}{\topology{}}$
is taken as a topological-space.
\begin{equation}
\func{\Int{\Xt}}{\X}=\X.
\end{equation}
\prooff
Considering that $\topology{}\subseteq\CSs{\X}$,
it is clear that,
\begin{equation}
\[\topology{}\cap\CSs{\X}\]=\topology{}.
\end{equation}
So, according to \refdef{defintofset}, and \refdef{deftopologicalspace},
\begin{align}
\func{\Int{\Xt}}{\X}&=\union{\[\topology{}\cap\CSs{\X}\]}\cr
&=\union{\topology{}}\cr
&=\X.
\end{align}
\endthm
%%%%%%%%%%%%%%%%%%%%%%%%%%%%%%%%%%%%%%%%%%%%%%%%%%%%%%%%%%%%%%%%%%%%%%%%%%%%%%%%%%%%%
\theorem\label{thmintofsetissetofintpoints}
$\Xt=\opair{\X}{\topology{}}$
is taken as a topological-space.
\begin{equation}
\Foreach{\asubset}{\CSs{\X}}
\func{\Int{\Xt}}{\asubset}=\Sint{\Xt}{\asubset}.
\end{equation}
\prooff
$\asubset$
is taken as an arbitrary subset of $\X$.
\begin{itemize}
\item[${\textbf{\textsf{p1}}}$]
$\x$
is taken as an arbitrary element of $\func{\Int{\Xt}}{\asubset}$.\\
According to \refcor{corintofset0},
\begin{equation}
\func{\Int{\Xt}}{\asubset}\in\[\topology{}\cap\CSs{\asubset}\].
\end{equation}
So, according to \refdef{definteriorpoint},
\begin{equation}
\x\in\Sint{\Xt}{\asubset}.
\end{equation}
\endp
\item[${\textbf{\textsf{p2}}}$]
$\p{\x}$
is taken as an arbitrary element of $\Sint{\Xt}{\asubset}$.
According to \refdef{definteriorpoint},
\begin{equation}
\Exists{\U}{\[\topology{}\cap\CSs{\asubset}\]}\p{\x}\in\U.
\end{equation}
$\U$
is taken as such an element of $\[\topology{}\cap\CSs{\asubset}\]$ that,
\begin{equation}
\p{\x}\in\U.
\end{equation}
According to \refcor{corintofset0},
\begin{equation}
\U\subseteq\func{\Int{\Xt}}{\asubset}.
\end{equation}
Therefore, it is clear that,
\begin{equation}
\p{\x}\in\func{\Int{\Xt}}{\asubset}.
\end{equation}
\endp
\end{itemize}
\endthm
%%%%%%%%%%%%%%%%%%%%%%%%%%%%%%%%%%%%%%%%%%%%%%%%%%%%%%%%%%%%%%%%%%%%%%%%%%%%%%%%
\theorem\label{thmintofopenset}
$\Xt=\opair{\X}{\topology{}}$
is taken as a topological-space.
For every
$\asubset$
in
$\CSs{\X}$,
$\asubset$
is an open set of
$\Xt$
if-and-only-if
$\func{\Int{\Xt}}{\asubset}$
equals
$\asubset$.
\begin{equation}
\Foreach{\asubset}{\CSs{\X}}
\[\(\asubset\in\topology{}\)\thenn\(\func{\Int{\Xt}}{\asubset}=\asubset\)\].
\end{equation}
\prooff
$\asubset$
is taken as an arbitrary element of $\X$.
\begin{itemize}
\item[${\textbf{\textsf{p1}}}$]
It is assumed that $\func{\Int{\Xt}}{\asubset}=\asubset$.
Then according to \refcor{corintofset0},
it is clear that,
\begin{equation}
\asubset\in\topology{}.
\end{equation}
\endp
\item[${\textbf{\textsf{p2}}}$]
It is assumed that,
$\asubset\in\topology{}$.
Then considering that $\asubset\subseteq\asubset$,
it is clear that,
\begin{equation}
\asubset\in\[\topology{}\cap\CSs{\asubset}\],
\end{equation}
and hence according to \refcor{corintofset0},
\begin{equation}
\asubset\subseteq\func{\Int{\Xt}}{\asubset}.
\end{equation}
In addition, according to \refcor{corintofset0},
\begin{equation}
\func{\Int{\Xt}}{\asubset}\subseteq\asubset.
\end{equation}
Therefore,
\begin{equation}
\func{\Int{\Xt}}{\asubset}=\asubset.
\end{equation}
\endp
\end{itemize}
\endthm
%%%%%%%%%%%%%%%%%%%%%%%%%%%%%%%%%%%%%%%%%%%%%%%%%%%%%%%%%%%%%%%%%%%%%%%%%%%%%%%%
\theorem\label{thmintofintofset}
$\Xt=\opair{\X}{\topology{}}$
is taken as a topological-space.
\begin{equation}\label{thmintofintofseteq1}
\Foreach{\asubset}{\CSs{\X}}
\func{\Int{\Xt}}{\func{\Int{\Xt}}{\asubset}}=\func{\Int{\Xt}}{\asubset}.
\end{equation}
\prooff
$\asubset$
is taken as an arbitrary subset of
$\X$.
According to \refcor{corintofset0},
\begin{equation}
\func{\Int{\Xt}}{\asubset}\in\topology{}.
\end{equation}
So, according to \refthm{thmintofopenset},
\Ref{thmintofintofseteq1}
is obtained.
\endthm
%%%%%%%%%%%%%%%%%%%%%%%%%%%%%%%%%%%%%%%%%%%%%%%%%%%%%%%%%%%%%%%%%%%%%%%%%%%%%%%%%%%%%%%%%%%%%
\theorem\label{thmintofasubsetofset}
$\Xt=\opair{\X}{\topology{}}$
is taken as a topological-space.
\begin{equation}
\Foreach{\opair{\asubset}{\p{\asubset}}}{\CSs{\X}\times\CSs{\X}}
\left\{\(\p{\asubset}\subseteq\asubset\)\then\[\func{\Int{\Xt}}{\p{\asubset}}\subseteq\func{\Int{\Xt}}{\asubset}\]\right\}.
\end{equation}
\prooff
$\asubset$
and
$\p{\asubset}$
are taken as such subsets of
$\X$
that
$\p{\asubset}\subseteq\asubset$.
Then,
\begin{equation}
\CSs{\p{\asubset}}\subseteq\CSs{\asubset},
\end{equation}
and accordingly,
\begin{equation}
\[\topology{}\cap\CSs{\p{\asubset}}\]\subseteq\[\topology{}\cap\CSs{\asubset}\].
\end{equation}
Thus clearly,
\begin{equation}
\union{\[\topology{}\cap\CSs{\p{\asubset}}\]}\subseteq\union{\[\topology{}\cap\CSs{\asubset}\]}.
\end{equation}
Hence according to \refdef{defintofset},
\begin{equation}
\func{\Int{\Xt}}{\p{\asubset}}\subseteq\func{\Int{\Xt}}{\asubset}.
\end{equation}
\endthm
%%%%%%%%%%%%%%%%%%%%%%%%%%%%%%%%%%%%%%%%%%%%%%%%%%%%%%%%%%%%%%%%%%%%%%%%%%%%%%%%%%%%%%%%%%%%%%%%%%%%%%
%%%%%%%%%%%%%%%%%%%%%%%%%%%%%%%%%%%%%%%%%%%%%%%%%%%%%%%%%%%%%%%%%%%%%%%%%%%%%%%%%%%%%%%%%%%%%%%%%%%%%%
\theorem\label{thmintofintersectionofsets0}
$\Xt=\opair{\X}{\topology{}}$
is taken as a topological-space.
\begin{align}
&\Foreach{\sCi}{\[\compl{\CSs{\CSs{\X}}}{\seta{\empty}}\]}\cr
&\func{\Int{\Xt}}{\intersection{\sCi}}\subseteq
\intersection{\defset{\U}{\CSs{\X}}{\[\Exists{\asubset}{\sCi}\U=\func{\Int{\Xt}}{\asubset}\]}}.
\end{align}
\prooff
$\sCi$
is taken as a non-empty collection of subsets of
$\X$.
Considering that,
\begin{equation}
\Foreach{\asubset}{\sCi}
\(\intersection{\sCi}\)\subseteq\asubset,
\end{equation}
According to \refthm{thmintofasubsetofset},
\begin{equation}
\Foreach{\asubset}{\sCi}
\func{\Int{\Xt}}{\intersection{\sCi}}\subseteq\func{\Int{\Xt}}{\asubset}.
\end{equation}
Hence,
\begin{equation}
\func{\Int{\Xt}}{\intersection{\sCi}}\subseteq
\bigcap_{\asubset\in\sCi}\func{\Int{\Xt}}{\asubset}.
\end{equation}
\endthm
%%%%%%%%%%%%%%%%%%%%%%%%%%%%%%%%%%%%%%%%%%%%%%%%%%%%%%%%%%%%%%%%%%%%%%%%%%%%%%%%%%%%%%%%%%%%%%%%%%%%%%
%%%%%%%%%%%%%%%%%%%%%%%%%%%%%%%%%%%%%%%%%%%%%%%%%%%%%%%%%%%%%%%%%%%%%%%%%%%%%%%%%%%%%%%%%%%%%%%%%%%%%%
\theorem\label{thmintofintersectionofsets}
$\Xt=\opair{\X}{\topology{}}$
is taken as a topological-space.
\begin{equation}
\Foreach{\opair{\asubset}{\p{\asubset}}}{\CSs{\X}\times\CSs{\X}}
\[\func{\Int{\Xt}}{\asubset\cap\p{\asubset}}=
\func{\Int{\Xt}}{\asubset}\cap\func{\Int{\Xt}}{\p{\asubset}}\].
\end{equation}
\prooff
$\asubset$
and
$\p{\asubset}$
are taken as arbitrary subsets of
$\X$.
\begin{itemize}
\item[${\textbf{\textsf{p1}}}$]
According to \refcor{corintofset0},
\begin{align}
\func{\Int{\Xt}}{\asubset}&\subseteq\asubset,\label{thmintofintersectionofsetsp1eq1}\\
\func{\Int{\Xt}}{\p{\asubset}}&\subseteq\p{\asubset}.\label{thmintofintersectionofsetsp1eq2}
\end{align}
Therefore,
\begin{equation}\label{thmintofintersectionofsetsp1eq3}
\[\func{\Int{\Xt}}{\asubset}\cap\func{\Int{\Xt}}{\p{\asubset}}\]
\in\CSs{\asubset\cap\p{\asubset}}.
\end{equation}
In addition, according to, \refcor{corintofset0},
\begin{align}
\func{\Int{\Xt}}{\asubset}&\in\topology{},\label{thmintofintersectionofsetsp1eq4}\\
\func{\Int{\Xt}}{\p{\asubset}}&\in\topology{}.\label{thmintofintersectionofsetsp1eq5}
\end{align}
So according to \refdef{deftopologicalspace},
\begin{equation}\label{thmintofintersectionofsetsp1eq6}
\[\func{\Int{\Xt}}{\asubset}\cap\func{\Int{\Xt}}{\p{\asubset}}\]\in\topology{}.
\end{equation}
\Ref{thmintofintersectionofsetsp1eq3}
and
\Ref{thmintofintersectionofsetsp1eq6}
imply,
\begin{equation}
\[\func{\Int{\Xt}}{\asubset}\cap\func{\Int{\Xt}}{\p{\asubset}}\]
\in\[\topology{}\cap\CSs{\asubset\cap\p{\asubset}}\].
\end{equation}
Thus according to \refcor{corintofset0},
\begin{equation}
\[\func{\Int{\Xt}}{\asubset}\cap\func{\Int{\Xt}}{\p{\asubset}}\]\subseteq
\func{\Int{\Xt}}{\asubset\cap\p{\asubset}}.
\end{equation}
\endp
\item[${\textbf{\textsf{p2}}}$]
According to
\refthm{thmintofintersectionofsets0},
\begin{equation}
\func{\Int{\Xt}}{\asubset\cap\p{\asubset}}\subseteq
\[\func{\Int{\Xt}}{\asubset}\cap\func{\Int{\Xt}}{\p{\asubset}}\].
\end{equation}
\endp
\end{itemize}
Therefore,
\begin{equation}
\func{\Int{\Xt}}{\asubset\cap\p{\asubset}}=
\[\func{\Int{\Xt}}{\asubset}\cap\func{\Int{\Xt}}{\p{\asubset}}\].
\end{equation}
\endthm
%%%%%%%%%%%%%%%%%%%%%%%%%%%%%%%%%%%%%%%%%%%%%%%%%%%%%%%%%%%%%%%%%%%%%%%%%%%%%%%%%%%%
\theorem\label{thmintofunionofsets}
$\Xt=\opair{\X}{\topology{}}$
is taken as a topological-space.
\begin{align}
&\Foreachs{\sCi}{\CSs{\X}}\cr
&\func{\Int{\Xt}}{\union{\sCi}}\supseteq
\union{\defset{\U}{\CSs{\X}}{\[\Exists{\asubset}{\sCi}\U=\func{\Int{\Xt}}{\asubset}\]}}.
\end{align}
\prooff
$\sCi$
is taken as an arbitrary collection of subsets of
$\X$.
Considering that,
\begin{equation}
\Foreach{\asubset}{\sCi}\func{\Int{\Xt}}{\asubset}\subseteq\asubset,
\end{equation}
it is clear that,
\begin{equation}
\union{\defset{\U}{\CSs{\X}}{\[\Exists{\asubset}{\sCi}\U=\func{\Int{\Xt}}{\asubset}\]}}
\subseteq\union{\sCi}.
\end{equation}
In addition,
$\defset{\U}{\CSs{\X}}{\[\Exists{\asubset}{\sCi}\U=\func{\Int{\Xt}}{\asubset}\]}$
is a collection of open sets of $\Xt$. That is,
\begin{equation}
\Foreach{\asubset}{\sCi}\func{\Int{\Xt}}{\asubset}\in\topology{}.
\end{equation}
So according to \refdef{deftopologicalspace},
\begin{equation}
\union{\defset{\U}{\CSs{\X}}{\[\Exists{\asubset}{\sCi}\U=\func{\Int{\Xt}}{\asubset}\]}}\in\topology{}.
\end{equation}
Therefore,
\begin{equation}
\union{\defset{\U}{\CSs{\X}}{\[\Exists{\asubset}{\sCi}\U=\func{\Int{\Xt}}{\asubset}\]}}\in
\[\topology{}\cap\CSs{\union{\sCi}}\].
\end{equation}
Hence according to \refcor{corintofset0},
\begin{equation}
\union{\defset{\U}{\CSs{\X}}{\[\Exists{\asubset}{\sCi}\U=\func{\Int{\Xt}}{\asubset}\]}}\subseteq
\func{\Int{\Xt}}{\union{\sCi}}.
\end{equation}
\endthm
%%%%%%%%%%%%%%%%%%%%%%%%%%%%%%%%%%%%%%%%%%%%%%%%%%%%%%%%%%%%%%%%%%%%%%%%%%%%%%%%%%%%%
\theorem
$\Xt=\opair{\X}{\topology{}}$
is taken as a topological-space, and
$\S$
as a subset of
$\X$.
\begin{equation}
\Foreach{\asubset}{\topology{}}
\[\asubset\subseteq\S\thenn\asubset\subseteq\func{\Int{\Xt}}{\S}\].
\end{equation}
\prooff
$\asubset$
is taken as an arbitrary element of
$\topology{}$
(an open set of $\Xt$).
\begin{itemize}
\item[${\textbf{\textsf{p1}}}$]
It is assumed that
$\asubset\subseteq\S$.
Then,
\begin{equation}
\asubset\in\[\topology{}\cap\CSs{\S}\],
\end{equation}
and according to \refcor{corintofset0},
\begin{equation}
\asubset\subseteq\func{\Int{\Xt}}{\S}.
\end{equation}
\endp
\item[${\textbf{\textsf{p2}}}$]
It is assumed that
$\asubset\subseteq\func{\Int{\Xt}}{\S}$.
Then considering that
$\func{\Int{\Xt}}{\S}\subseteq\S$,
it is clear that,
\begin{equation}
\asubset\subseteq\S.
\end{equation}
\endp
\end{itemize}
\endthm
%%%%%%%%%%%%%%%%%%%%%%%%%%%%%%%%%%%%%%%%%%%%%%%%%%%%%%%%%%%%%%%%%%%%%%%%%%%
%%%%%%%%%%%%%%%%%%%%%%%%%%%%%%%%%%%%%%%%%%%%%%%%%%%%%%%%%%%%%%%%%%%%%%%%%%%
%%%%%%%%%%%%%%%%%%%%%%%%%%%%%%%%%%%%%%%%%%%%%%%%%%%%%%%%%%%%%%%%%%%%%%%%%%%
%%%%%%%%%%%%%%%%%%%%%%%%%%%%%%%%%%%%%%%%%%%%%%%%%%%%%%%%%%%%%%%%%%%%%%%%%%%%%%%%%%%%%%%%%%%%%%%%%%%%%%%%%%%%%%%%%%%%%%%%%%%%%%%%%%%%%%%%%%%%%%%%%%%%%%
%%%%%%%%%%%%%%%%%%%%%%%%%%%%%%%%%%%%%%%%%%%%%%%%%%%%%%%%%%%%%%%%%%%%%%%%%%%
%%%%%%%%%%%%%%%%%%%%%%%%%%%%%%%%%%%%%%%%%%%%%%%%%%%%%%%%%%%%%%%%%%%%%%%%%%%
%%%%%%%%%%%%%%%%%%%%%%%%%%%%%%%%%%%%%%%%%%%%%%%%%%%%%%%%%%%%%%%%%%%%%%%%%%%
\theorem\label{thmsubspaceinterior}
$\opair{\X}{\topology{}}$
is taken as a topological-space, and $\Y$
as a subset of $\X$.
\begin{equation}
\Foreach{\asubset}{\CSs{\Y}}
\[\func{\Int{\Xt}}{\asubset}\subseteq
\func{\Int{\opair{\Y}{\stopology{\topology{}}{\Y}}}}{\asubset}\].
\end{equation}
\prooff
$\asubset$
is taken as an arbitrary subset of $\Y$.
According to \refdef{defsubspacetopology1},
\begin{align}\label{thmsubspaceinteriorp1}
\stopology{\topology{}}{\Y}=
\defset{\V}{\CSs{\Y}}{\[\Exists{\U}{\topology{}}\V=\Y\cap\U\]}.
\end{align}
\begin{itemize}
\item[${\textbf{\textsf{p}}}$]
$\U$
is taken as an element of $\[\topology{}\cap\CSs{\asubset}\]$.
Hence,
\begin{align}
\U&\in\topology{},\label{thmsubspaceinteriorpp1}\\
\U&\subseteq\asubset.\label{thmsubspaceinteriorpp2}
\end{align}
According to \Ref{thmsubspaceinteriorpp2},
and considering that $\asubset\subseteq\Y$,
it is clear that,
\begin{equation}\label{thmsubspaceinteriorpp3}
\U\subseteq\Y,
\end{equation}
and accordingly,
\begin{equation}\label{thmsubspaceinteriorpp4}
\Y\cap\U=\U
\end{equation}
According to
\Ref{thmsubspaceinteriorp1},
\Ref{thmsubspaceinteriorpp1},
and
\Ref{thmsubspaceinteriorpp4},
\begin{equation}\label{thmsubspaceinteriorpp5}
\U\in\stopology{\topology{}}{\Y}.
\end{equation}
\Ref{thmsubspaceinteriorpp2}
and
\Ref{thmsubspaceinteriorpp5}
imply,
\begin{equation}
\U\in\[\stopology{\topology{}}{\Y}\cap\CSs{\asubset}\].
\end{equation}
\endp
\end{itemize}
Therefore,
\begin{equation}\label{thmsubspaceinteriorp2}
\[\topology{}\cap\CSs{\asubset}\]\subseteq
\[\stopology{\topology{}}{\Y}\cap\CSs{\asubset}\].
\end{equation}
Thus according to \refdef{defintofset},
\begin{align}
\func{\Int{\Xt}}{\asubset}&=
\union{\[\topology{}\cap\CSs{\asubset}\]}\cr
&\subseteq\union{\[\stopology{\topology{}}{\Y}\cap\CSs{\asubset}\]}\cr
&=\func{\Int{\opair{\Y}{\stopology{\topology{}}{\Y}}}}{\asubset}.
\end{align}
\endthm
%%%%%%%%%%%%%%%%%%%%%%%%%%%%%%%%%%%%%%%%%%%%%%%%%%%%%%%%%%%%%%%%%%%%%%%%%%%%%%%%%%%%%%%%%%
%%%%%%%%%%%%%%%%%%%%%%%%%%%%%%%%%%%%%%%%%%%%%%%%%%%%%%%%%%%%%%%%%%%%%%%%%%%
\theorem\label{thmintofsetwithdifferenttopologies}
Each $\Xt_{1}=\opair{\X}{\topology{1}}$ and $\Xt_{2}=\opair{\X}{\topology{2}}$
is taken as a topological-space.
\begin{equation}
\(\topology{2}\subseteq\topology{1}\)\then\[\Foreach{\asubset}{\CSs{\X}}
\func{\Int{\Xt_{2}}}{\asubset}\subseteq\func{\Int{\Xt_{1}}}{\asubset}\].
\end{equation}
\prooff
It is assumed that,
\begin{equation}
\topology{2}\subseteq\topology{1},
\end{equation}
and
$\asubset$ is taken as a subset of $\X$.
According to
\refdef{defintofset},
\begin{align}
\func{\Int{\Xt_{1}}}{\asubset}&=\union{\[\topology{1}\cap\CSs{\asubset}\]},\\
\func{\Int{\Xt_{2}}}{\asubset}&=\union{\[\topology{2}\cap\CSs{\asubset}\]}.
\end{align}
Therefore,
\begin{equation}
\func{\Int{\Xt_{2}}}{\asubset}\subseteq\func{\Int{\Xt_{1}}}{\asubset}.
\end{equation}
\endthm
%%%%%%%%%%%%%%%%%%%%%%%%%%%%%%%%%%%%%%%%%%%%%%%%%%%%%%%%%%%%%%%%%%%%%%%%%%%
%%%%%%%%%%%%%%%%%%%%%%%%%%%%%%%%%%%%%%%%%%%%%%%%%%%%%%%%%%%%%%%%%%%%%%%%%%%
%%%%%%%%%%%%%%%%%%%%%%%%%%%%%%%%%%%%%%%%%%%%%%%%%%%%%%%%%%%%%%%%%%%%%%%%%%%
%%%%%%%%%%%%%%%%%%%%%%%%%%%%%%%%%%%%%%%%%%%%%%%%%%%%%%%%%%%%%%%%%%%%%%%%%%%
\subsection{
Exterior of a Subset of a Topological Space}
\definition\label{defextofset}
$\Xt=\opair{\X}{\topology{}}$
is taken as a topological-space. The mapping $\Ext{\Xt}$
is defined as,
\begin{itemize}
\item[${\textbf{\textsf{ext1}}}$]
\hfill
$\Ext{\Xt}\indef\Func{\CSs{\X}}{\CSs{\X}}.$
\item[${\textbf{\textsf{ext2}}}$]
\hfill
$\Foreach{\asubset}{\CSs{\X}}
\func{\Ext{\Xt}}{\asubset}\eqdef\func{\Int{\Xt}}{\compl{\X}{\asubset}}.$
\end{itemize}
For every
$\asubset$
in
$\CSs{\X}$,
$\func{\Ext{\Xt}}{\asubset}$
is referred to as the $\quotl$exterior of $\asubset$ in the topological-space $\Xt$$\quotr$.
\endef
%%%%%%%%%%%%%%%%%%%%%%%%%%%%%%%%%%%%%%%%%%%%%%%%%%%%%%%%%%%%%%%%%%%%%%%%%%%%%%%%%%%%%%%%%
\theorem\label{thmextofset0}
$\Xt=\opair{\X}{\topology{}}$
is taken as a topological-space.
\begin{align}
\Foreach{\asubset}{\CSs{\X}}\left\{
\begin{aligned}
&\func{\Ext{\Xt}}{\asubset}=\union{\[\topology{}\cap\CSs{\compl{\X}{\asubset}}\]},\\
&\func{\Ext{\Xt}}{\asubset}\in\topology{},\\
&\func{\Ext{\Xt}}{\asubset}\cap\asubset=\empty,\\
&\Foreach{\U}{\[\topology{}\cap\CSs{\compl{\X}{\asubset}}\]}\U\subseteq\func{\Ext{\Xt}}{\asubset}.
\end{aligned}
\right.
\end{align}
\prooff
According to
\refdef{defextofset}
and
\refcor{corintofset0},
it is clear.
\endthm
%%%%%%%%%%%%%%%%%%%%%%%%%%%%%%%%%%%%%%%%%%%%%%%%%%%%%%%%%%%%%%%%%%%%%%%%%%%%%%%%%%%%%%%%%%%%%%%%%%%%%%%%
\theorem\label{thmextofemptyset}
$\Xt=\opair{\X}{\topology{}}$
is taken as a topological-space.
\begin{equation}
\func{\Ext{\Xt}}{\empty}=\X.
\end{equation}
\prooff
According to
\refdef{defextofset}
and
\refthm{thmintofuniversalset},
\begin{align}
\func{\Ext{\Xt}}{\empty}&=\func{\Int{\Xt}}{\compl{\X}{\empty}}\cr
&=\func{\Int{\Xt}}{\X}\cr
&=\X.
\end{align}
\endthm
%%%%%%%%%%%%%%%%%%%%%%%%%%%%%%%%%%%%%%%%%%%%%%%%%%%%%%%%%%%%%%%%%%%%%%%%%%%%%%%%%%%%%%%%%%%%%%%%%%%%%%%%
%%%%%%%%%%%%%%%%%%%%%%%%%%%%%%%%%%%%%%%%%%%%%%%%%%%%%%%%%%%%%%%%%%%%%%%%%%%%%%%%%%%%%%%%%%%%%%%%%%%%%%%%
\theorem\label{thmextofuniversalset}
$\Xt=\opair{\X}{\topology{}}$
is taken as a topological-space.
\begin{equation}
\func{\Ext{\Xt}}{\X}=\empty.
\end{equation}
\prooff
According to
\refdef{defextofset}
and
\refthm{thmintofemptyset},
\begin{align}
\func{\Ext{\Xt}}{\X}&=\func{\Int{\Xt}}{\compl{\X}{\X}}\cr
&=\func{\Int{\Xt}}{\empty}\cr
&=\empty.
\end{align}
\endthm
%%%%%%%%%%%%%%%%%%%%%%%%%%%%%%%%%%%%%%%%%%%%%%%%%%%%%%%%%%%%%%%%%%%%%%%%%%%%%%%%%%%%%%%%%%%%%%%%%%%%%%%%
\theorem\label{thmextofsetissetofextpoints}
$\Xt=\opair{\X}{\topology{}}$
is taken as a topological-space.
\begin{equation}
\Foreach{\asubset}{\CSs{\X}}\func{\Ext{\Xt}}{\asubset}=\Sext{\Xt}{\asubset}.
\end{equation}
\prooff
According to
\refcor{corSintSext},
\refthm{thmintofsetissetofintpoints},
and
\refdef{defextofset},
it is clear.
\endthm
%%%%%%%%%%%%%%%%%%%%%%%%%%%%%%%%%%%%%%%%%%%%%%%%%%%%%%%%%%%%%%%%%%%%%%%%%%%%%%%%%%%%%%%%%%%%%%%%%%%%%%%%
\theorem
$\Xt=\opair{\X}{\topology{}}$
is taken as a topological-space.
\begin{equation}
\Foreach{\asubset}{\CSs{\X}}
\[\asubset\in\Fclosed{\X}{\topology{}}
\thenn\func{\Ext{\Xt}}{\asubset}=\(\compl{\X}{\asubset}\)\].
\end{equation}
\prooff
$\asubset$
is taken as an arbitrary subset of $\X$.
According to \refthm{thmintofopenset},
and
\refdef{defextofset},
\begin{align}
\asubset\in\Fclosed{\X}{\topology{}}&\thenn
\(\compl{\X}{\asubset}\)\in\topology{},\\
\(\compl{\X}{\asubset}\)\in\topology{}&\thenn
\func{\Int{\Xt}}{\compl{\X}{\asubset}}=\(\compl{\X}{\asubset}\),\\
\func{\Int{\Xt}}{\compl{\X}{\asubset}}=\(\compl{\X}{\asubset}\)&\thenn
\func{\Ext{\Xt}}{\asubset}=\(\compl{\X}{\asubset}\).
\end{align}
Therefore it is clear that,
\begin{equation}
\asubset\in\Fclosed{\X}{\topology{}}\thenn
\func{\Ext{\Xt}}{\asubset}=\(\compl{\X}{\asubset}\).
\end{equation}
\endthm
%%%%%%%%%%%%%%%%%%%%%%%%%%%%%%%%%%%%%%%%%%%%%%%%%%%%%%%%%%%%%%%%%%%%%%%%%%%%%%%%%%%%%%%%%%%%%%%%%%%%%%%%
%%%%%%%%%%%%%%%%%%%%%%%%%%%%%%%%%%%%%%%%%%%%%%%%%%%%%%%%%%%%%%%%%%%%%%%%%%%%%%%%%%%%%%%%%%%%%%%%%%%%%%%%
%%%%%%%%%%%%%%%%%%%%%%%%%%%%%%%%%%%%%%%%%%%%%%%%%%%%%%%%%%%%%%%%%%%%%%%%%%%%%%%%%%%%%%%%%%%%%%%%%%%%%%%%
\theorem\label{thmextofasubsetofset}
$\Xt=\opair{\X}{\topology{}}$
is taken as a topological-space.
\begin{equation}
\Foreach{\opair{\asubset}{\p{\asubset}}}{\CSs{\X}\times\CSs{\X}}
\left\{\(\p{\asubset}\subseteq\asubset\)\then\[\func{\Ext{\Xt}}{\p{\asubset}}\supseteq\func{\Ext{\Xt}}{\asubset}\]\right\}.
\end{equation}
\prooff
$\asubset$
and
$\p{\asubset}$
are taken as such subsets of $\X$ that
$\p{\asubset}\subseteq\asubset$.
Then,
\begin{equation}
\(\compl{\X}{\asubset}\)\subseteq\(\compl{\X}{\p{\asubset}}\),
\end{equation}
and hence according to \refthm{thmintofasubsetofset},
\begin{equation}
\func{\Int{\Xt}}{\compl{\X}{\asubset}}\subseteq
\func{\Int{\Xt}}{\compl{\X}{\p{\asubset}}}.
\end{equation}
Thus according to
\refdef{defextofset},
\begin{equation}
\func{\Ext{\Xt}}{\p{\asubset}}\supseteq\func{\Ext{\Xt}}{\asubset}.
\end{equation}
\endthm
%%%%%%%%%%%%%%%%%%%%%%%%%%%%%%%%%%%%%%%%%%%%%%%%%%%%%%%%%%%%%%%%%%%%%%%%%%%%%%%%%
\theorem
$\Xt=\opair{\X}{\topology{}}$
is taken as a topological-space.
\begin{equation}
\Foreach{\asubset}{\CSs{\X}}
\func{\Ext{\Xt}}{\func{\Ext{\Xt}}{\asubset}}\supseteq
\func{\Int{\Xt}}{\asubset}.
\end{equation}
\prooff
$\asubset$
is taken as a subset of $\X$.
According to
\refthm{thmextofset0},
\begin{equation}
\func{\Ext{\Xt}}{\asubset}\subseteq\(\compl{\X}{\asubset}\).
\end{equation}
Thus according to
\refthm{thmextofasubsetofset},
\begin{equation}
\func{\Ext{\Xt}}{\func{\Ext{\Xt}}{\asubset}}\supseteq\func{\Ext{\Xt}}{\compl{\X}{\asubset}}.
\end{equation}
In addition, according to
\refdef{defextofset},
\begin{align}
\func{\Ext{\Xt}}{\compl{\X}{\asubset}}=\func{\Int{\Xt}}{\asubset}.
\end{align}
Therefore,
\begin{equation}
\func{\Ext{\Xt}}{\func{\Ext{\Xt}}{\asubset}}\supseteq\func{\Int{\Xt}}{\asubset}.
\end{equation}
\endthm
%%%%%%%%%%%%%%%%%%%%%%%%%%%%%%%%%%%%%%%%%%%%%%%%%%%%%%%%%%%%%%%%%%%%%%%%%%%%%%%%%%%%%%%%%%%%%%%%%%
\theorem\label{thmextofunionofsets0}
$\Xt=\opair{\X}{\topology{}}$
is taken as a topological-space.
\begin{align}
\Foreach{\sCi}{\[\compl{\CSs{\CSs{\X}}}{\seta{\empty}}\]}
\[\func{\Ext{\Xt}}{\union{\sCi}}\subseteq
\bigcap_{\asubset\in\sCi}\func{\Ext{\Xt}}{\asubset}\].
\end{align}
\prooff
$\sCi$
is taken as a non-empty collection of subsets of $\X$.
According to
\refdef{defextofset}
and
\refthm{thmintofintersectionofsets0},
\begin{align}
\func{\Ext{\Xt}}{\union{\sCi}}&=
\func{\Int{\Xt}}{\compl{\X}{\union{\sCi}}}\cr
&=\func{\Int{\Xt}}{\bigcap_{\asubset\in\sCi}\(\compl{\X}{\asubset}\)}\cr
&\subseteq\bigcap_{\asubset\in\sCi}\func{\Int{\Xt}}{\compl{\X}{\asubset}}\cr
&=\bigcap_{\asubset\in\sCi}\func{\Ext{\Xt}}{\asubset}.
\end{align}
\endthm
%%%%%%%%%%%%%%%%%%%%%%%%%%%%%%%%%%%%%%%%%%%%%%%%%%%%%%%%%%%%%%%%%%%%%%%%%%%%%%%%%%%%%%%%%%%%%%%%%%
\theorem\label{thmextofunionofsets}
$\Xt=\opair{\X}{\topology{}}$
is taken as a topological-space.
\begin{equation}
\Foreach{\opair{\asubset}{\p{\asubset}}}{\CSs{\X}\times\CSs{\X}}
\[\func{\Ext{\Xt}}{\asubset\cup\p{\asubset}}=
\func{\Ext{\Xt}}{\asubset}\cap\func{\Ext{\Xt}}{\p{\asubset}}\].
\end{equation}
\prooff
Each
$\asubset$
and
$\p{\asubset}$
is taken as an arbitrary subset of $\X$.
According to
\refthm{thmintofintersectionofsets}
and
\refdef{defextofset},
\begin{align}
\func{\Ext{\Xt}}{\asubset\cup\p{\asubset}}&=
\func{\Int{\Xt}}{\compl{\X}{\(\asubset}\cup\p{\asubset}\)}\cr
&=\func{\Int{\Xt}}{\(\compl{\X}{\asubset}\)\cap\(\compl{\X}{\p{\asubset}}\)}\cr
&=\func{\Int{\Xt}}{\compl{\X}{\asubset}}\cap\func{\Int{\Xt}}{\compl{\X}{\p{\asubset}}}\cr
&=\func{\Ext{\Xt}}{\asubset}\cap\func{\Ext{\Xt}}{\p{\asubset}}.
\end{align}
\endthm
%%%%%%%%%%%%%%%%%%%%%%%%%%%%%%%%%%%%%%%%%%%%%%%%%%%%%%%%%%%%%%%%%%%%%%%%%%%%%%%%%%%%%%%%%%%%%%%%%%%
\theorem\label{thmextofintersectionofsets}
$\Xt=\opair{\X}{\topology{}}$
is taken as a topological-space.
\begin{align}
&\Foreach{\sCi}{\[\compl{\CSs{\CSs{\X}}}{\seta{\empty}}\]}\cr
&\[\func{\Ext{\Xt}}{\intersection{\sCi}}\supseteq
\union{\defset{\U}{\CSs{\X}}{\Exists{\asubset}{\sCi}\U=\func{\Ext{\Xt}}{\asubset}}}\].
\end{align}
\prooff
$\sCi$
is taken as a non-empty collection of subsets of $\X$.
According to
\refthm{thmintofunionofsets}
and
\refdef{defextofset},
\begin{align}
\func{\Ext{\Xt}}{\intersection{\sCi}}&=
\func{\Int{\Xt}}{\compl{\X}{\(\intersection{\sCi}\)}}\cr
&=\func{\Int{\Xt}}{\bigcup_{\asubset\in\sCi}\[\compl{\X}{\asubset}\]}\cr
&\supseteq\bigcup_{\asubset\in\sCi}\func{\Int{\Xt}}{\compl{\X}{\asubset}}\cr
&=\bigcup_{\asubset\in\sCi}\func{\Ext{\Xt}}{\asubset}.
\end{align}
\endthm
%%%%%%%%%%%%%%%%%%%%%%%%%%%%%%%%%%%%%%%%%%%%%%%%%%%%%%%%%%%%%%%%%%%%%%%%%%%%%%%%%%%%%%%%%%%%%%%%%%%
\theorem\label{thmsubspaceexterior}
$\Xt=\opair{\X}{\topology{}}$
is taken as a topological-space, and
$\Y$ as a subset of $\X$.
\begin{equation}
\Foreach{\asubset}{\CSs{\Y}}
\left\{
\[\func{\Int{\Xt}}{\Y}\cap\func{\Ext{\Xt}}{\asubset}\]\subseteq
\func{\Ext{\opair{\Y}{\stopology{\topology{}}{\Y}}}}{\asubset}
\right\}.
\end{equation}
\prooff
$\asubset$
is taken as an arbitrary subset of $\Y$.
It is clear that,
\begin{equation}
\(\compl{\Y}{\asubset}\)=\Y\cap\(\compl{\X}{\asubset}\).
\end{equation}
Thus according to
\refthm{thmintofintersectionofsets},
\refthm{thmsubspaceinterior}
and
\refdef{defextofset},
\begin{align}
\func{\Ext{\opair{\Y}{\stopology{\topology{}}{\Y}}}}{\asubset}&=
\func{\Int{\opair{\Y}{\stopology{\topology{}}{\Y}}}}{\compl{\Y}{\asubset}}\cr
&\supseteq\func{\Int{\Xt}}{\compl{\Y}{\asubset}}\cr
&=\func{\Int{\Xt}}{\Y\cap\[\compl{\X}{\asubset}\]}\cr
&=\[\func{\Int{\Xt}}{\Y}\cap\func{\Int{\Xt}}{\compl{\X}{\asubset}}\]\cr
&=\func{\Int{\Xt}}{\Y}\cap\func{\Ext{\Xt}}{\asubset}.
\end{align}
\endthm
%%%%%%%%%%%%%%%%%%%%%%%%%%%%%%%%%%%%%%%%%%%%%%%%%%%%%%%%%%%%%%%%%%%%%%%%%%%%%%%%%%%%%%%%%%%%%%%%%%%%
%%%%%%%%%%%%%%%%%%%%%%%%%%%%%%%%%%%%%%%%%%%%%%%%%%%%%%%%%%%%%%%%%%%%%%%%%%%%%%%%%%%%%%%%%%%%%%%%%%%%
%%%%%%%%%%%%%%%%%%%%%%%%%%%%%%%%%%%%%%%%%%%%%%%%%%%%%%%%%%%%%%%%%%%%%%%%%%%%%%%%%%%%%%%%%%%%%%%%%%%%
%%%%%%%%%%%%%%%%%%%%%%%%%%%%%%%%%%%%%%%%%%%%%%%%%%%%%%%%%%%%%%%%%%%%%%%%%%%%%%%%%%%%%%%%%%%%%%%%%%%%
%%%%%%%%%%%%%%%%%%%%%%%%%%%%%%%%%%%%%%%%%%%%%%%%%%%%%%%%%%%%%%%%%%%%%%%%%%%%%%%%%%%%%%%%%%%%%%%%%%%%
%%%%%%%%%%%%%%%%%%%%%%%%%%%%%%%%%%%%%%%%%%%%%%%%%%%%%%%%%%%%%%%%%%%%%%%%%%%%%%%%%%%%%%%%%%%%%%%%%%%%
%%%%%%%%%%%%%%%%%%%%%%%%%%%%%%%%%%%%%%%%%%%%%%%%%%%%%%%%%%%%%%%%%%%%%%%%%%%%%%%%%%%%%%%%%%%%%%%%%%%%
%%%%%%%%%%%%%%%%%%%%%%%%%%%%%%%%%%%%%%%%%%%%%%%%%%%%%%%%%%%%%%%%%%%%%%%%%%%%%%%%%%%%%%%%%%%%%%%%%%%%
%%%%%%%%%%%%%%%%%%%%%%%%%%%%%%%%%%%%%%%%%%%%%%%%%%%%%%%%%%%%%%%%%%%%%%%%%%%%%%%%%%%%%%%%%%%%%%%%%%%%
\theorem\label{thmextofsetwithdifferenttopologies}
Each
$\Xt_{1}=\opair{\X}{\topology{1}}$
and
$\Xt_{2}=\opair{\X}{\topology{2}}$
is taken as a topological-space.
\begin{equation}
\(\topology{2}\subseteq\topology{1}\)\then\[\Foreach{\asubset}{\CSs{\X}}
\func{\Ext{\Xt_{2}}}{\asubset}\subseteq\func{\Ext{\Xt_{1}}}{\asubset}\].
\end{equation}
\prooff
According to
\refthm{thmintofsetwithdifferenttopologies}
and
\refdef{defextofset},
it is clear.
\endthm
%%%%%%%%%%%%%%%%%%%%%%%%%%%%%%%%%%%%%%%%%%%%%%%%%%%%%%%%%%%%%%%%%%%%%%%%%%%%%%%%%%%%%%%%%%
%%%%%%%%%%%%%%%%%%%%%%%%%%%%%%%%%%%%%%%%%%%%%%%%%%%%%%%%%%%%%%%%%%%%%%%%%%%%%%%%%%%%%%%%%%
%%%%%%%%%%%%%%%%%%%%%%%%%%%%%%%%%%%%%%%%%%%%%%%%%%%%%%%%%%%%%%%%%%%%%%%%%%%%%%%%%%%%%%%%%%
%%%%%%%%%%%%%%%%%%%%%%%%%%%%%%%%%%%%%%%%%%%%%%%%%%%%%%%%%%%%%%%%%%%%%%%%%%%%%%%%%%%%%%%%%%
%%%%%%%%%%%%%%%%%%%%%%%%%%%%%%%%%%%%%%%%%%%%%%%%%%%%%%%%%%%%%%%%%%%%%%%%%%%%%%%%%%%%%%%%%%
\subsection{
Closure of a Subset of a Topological Space}
\definition\label{defCinc}
$\X$
is taken as a set.
The mapping $\Cinc{\X}$
is defined as,
\begin{itemize}
\item[${\textbf{\textsf{Inc1}}}$]
\hfill
$\Cinc{\X}\indef\Func{\CSs{\X}}{\CSs{\CSs{\X}}}.$
\item[${\textbf{\textsf{Inc2}}}$]
\hfill
$\Foreach{\asubset}{\X}\func{\Cinc{\X}}{\asubset}\eqdef
\defset{\U}{\CSs{\X}}{\asubset\subseteq\U}.$
\end{itemize}
\endef
%%%%%%%%%%%%%%%%%%%%%%%%%%%%%%%%%%%%%%%%%%%%%%%%%%%%%%%%%%%%%%%%%%%%%%%%%%%%%%%%%%%%%%
\theorem\label{thmpreclosure0}
$\Xt=\opair{\X}{\topology{}}$
is taken as a topological-space, and
$\asubset$
as a subset of
$\X$.
There exists only one closed set of $\Xt$ including $\asubset$
that is a subset of every closed set of $\Xt$ including $\asubset$.
That is,
\begin{equation}
\Existsu{\U}{\[\Fclosed{\X}{\topology{}}\cap\func{\Cinc{\X}}{\asubset}\]}
\left\{\Foreach{\V}{\[\Fclosed{\X}{\topology{}}\cap\func{\Cinc{\X}}{\asubset}\]}\U\subseteq\V\right\}.
\end{equation}
More precisely,
\begin{gather}
\left\{\U\in\[\Fclosed{\X}{\topology{}}\cap\func{\Cinc{\X}}{\asubset}\],~
\Foreach{\V}{\[\Fclosed{\X}{\topology{}}\cap\func{\Cinc{\X}}{\asubset}\]}\V\subseteq\U\right\}\cr
\vthenn\cr
\left\{\U=\intersection{\[\Fclosed{\X}{\topology{}}\cap\func{\Cinc{\X}}{\asubset}\]}\right\}.
\end{gather}
\prooff
The set
$\U$
is defined as
\begin{equation}\label{thmpreint0peq1}
\U:=\intersection{\[\Fclosed{\X}{\topology{}}\cap\func{\Cinc{\X}}{\asubset}\]}.
\end{equation}
\begin{itemize}
\item[${\textbf{\textsf{p1}}}$]
Considering that,
\begin{equation}\label{thmpreclosure0p1eq1}
\[\Fclosed{\X}{\topology{}}\cap
\func{\Cinc{\X}}{\asubset}\]\subseteq\Fclosed{\X}{\topology{}},
\end{equation}
According to
\refthm{thmclosedsets},
it is clear that,
\begin{equation}\label{thmpreclosure0p1eq2}
\intersection{\[\Fclosed{\X}{\topology{}}
\cap\func{\Cinc{\X}}{\asubset}\]}\in\Fclosed{\X}{\topology{}}.
\end{equation}
In addition, considering that,
\begin{equation}\label{thmpreclosure0p1eq3}
\[\Fclosed{\X}{\topology{}}\cap
\func{\Cinc{\X}}{\asubset}\]\subseteq\func{\Cinc{\X}}{\asubset},
\end{equation}
it is clear that,
\begin{equation}\label{thmpreclosure0p1eq4}
\intersection{\[\Fclosed{\X}{\topology{}}
\cap\func{\Cinc{\X}}{\asubset}\]}\in\func{\Cinc{\X}}{\asubset},
\end{equation}
That is, the intersection of all elements of a collection of subset of $\X$
includinf
$\asubset$
must be a subset of
$\X$
including
$\asubset$.
Therefore,
\begin{equation}\label{thmpreclosure0p1eq5}
\intersection{\[\Fclosed{\X}{\topology{}}
\cap\func{\Cinc{\X}}{\asubset}\]}\in
\[\Fclosed{\X}{\topology{}}
\cap\func{\Cinc{\X}}{\asubset}\].
\end{equation}
It is clear that,
\begin{equation}\label{thmpreclosure0p1eq6}
\Foreach{\V}{\[\Fclosed{\X}{\topology{}}
\cap\func{\Cinc{\X}}{\asubset}\]}
\intersection{\[\Fclosed{\X}{\topology{}}
\cap\func{\Cinc{\X}}{\asubset}\]}\subseteq\V.
\end{equation}
Therefore,
\begin{align}\label{thmpreclosure0p1eq7}
\Existsis{\U}{\[\Fclosed{\X}{\topology{}}\cap\func{\Cinc{\X}}{\asubset}\]}
\left\{\Foreach{\V}{\[\Fclosed{\X}{\topology{}}\cap\func{\Cinc{\X}}{\asubset}\]}\U\subseteq\V\right\}.
\end{align}
\endp
\item[${\textbf{\textsf{p2}}}$]
$\U^{\prime}$
is taken as such an element of
$\Fclosed{\X}{\topology{}}\cap\func{\Cinc{\X}}{\asubset}$
that,
\begin{equation}
\Foreach{\V}{\[\Fclosed{\X}{\topology{}}\cap\func{\Cinc{\X}}{\asubset}\]}\p{\U}\subseteq\V.
\end{equation}
Then according to \Ref{thmpreclosure0p1eq7},
it is clear that,
\begin{equation}
\left\{
\begin{aligned}
\U^{\prime}\subseteq\U,\\
\U\subseteq\U^{\prime},
\end{aligned}
\right.
\end{equation}
which means,
\begin{equation}
\U^{\prime}=\U.
\end{equation}
\endp
\end{itemize}
\endthm
%%%%%%%%%%%%%%%%%%%%%%%%%%%%%%%%%%%%%%%%%%%%%%%%%%%%%%%%%%%%%%%%%%%%%%%%%%%%%%%%%%%%
\definition\label{defclosureofset}
$\Xt=\opair{\X}{\topology{}}$
is taken as a topological-space.
The mapping
$\Cl{\Xt}$
is defined as,
\begin{itemize}
\item[${\textbf{\textsf{Cl1}}}$]
\hfill
$\Cl{\Xt}\indef\Func{\CSs{\X}}{\CSs{\X}}.$
\item[${\textbf{\textsf{Cl2}}}$]
\hfill
$\Foreach{\asubset}{\CSs{\X}}
\func{\Cl{\Xt}}{\asubset}\eqdef
\bigcap\[\Fclosed{\X}{\topology{}}\cap\func{\Cinc{\X}}{\asubset}\].$
\end{itemize}
For every
$\asubset$
in
$\CSs{\X}$,
$\func{\Cl{\Xt}}{\asubset}$
is referred to as the $\quotl$closure of $\asubset$ in the topological-space $\Xt$$\quotr$.
\endef
%%%%%%%%%%%%%%%%%%%%%%%%%%%%%%%%%%%%%%%%%%%%%%%%%%%%%%%%%%%%%%%%%%%%%%%%%%%%%%%%%%%%
\corollary\label{corclosureofset0}
$\Xt=\opair{\X}{\topology{}}$
is taken as a topological-space.
\begin{align}
\Foreach{\asubset}{\CSs{\X}}\left\{
\begin{aligned}
&\func{\Cl{\Xt}}{\asubset}\in\Fclosed{\X}{\topology{}},\\
&\func{\Cl{\Xt}}{\asubset}\supseteq\asubset,\\
&\Foreach{\U}{\[\Fclosed{\X}{\topology{}}\cap\func{\Cinc{\X}}{\asubset}\]}
\U\supseteq\func{\Cl{\Xt}}{\asubset}.
\end{aligned}
\right.
\end{align}
In other words, for every
$\asubset$
in
$\CSs{\X}$,
$\func{\Cl{\Xt}}{\asubset}$
is a closed set of $\Xt$ including
$\asubset$, and a subset of every closed set of $\Xt$ including $\asubset$.
In addition, for every
$\asubset$
in
$\X$,
such subset of
$\X$
is unique.
\endcor
%%%%%%%%%%%%%%%%%%%%%%%%%%%%%%%%%%%%%%%%%%%%%%%%%%%%%%%%%%%%%%%%%%%%%%%%%%%%%%%%%%%%%%%%%%%
\theorem\label{thmclosureofsetissetofadhpoints}
$\Xt=\opair{\X}{\topology{}}$
is taken as a topological-space.
\begin{equation}
\Foreach{\asubset}{\CSs{\X}}
\func{\Cl{\Xt}}{\asubset}=\Sadh{\Xt}{\asubset}.
\end{equation}
\prooff
$\asubset$
is taken as an arbitrary subset of
$\X$.
\begin{itemize}
\item[${\textbf{\textsf{p1}}}$]
$\x$
is taken as an arbitrary element of
$\compl{\X}{\func{\Cl{\Xt}}{\asubset}}$.\\
According to
\refcor{corclosureofset0},
$\func{\Cl{\Xt}}{\asubset}$
is a closed set of $\Xt$,
and hence
\begin{equation}
\[\compl{\X}{\func{\Cl{\Xt}}{\asubset}}\]\in\topology{},
\end{equation}
and thus according to \refdef{defnbdclassofsets},
\begin{equation}
\[\compl{\X}{\func{\Cl{\Xt}}{\asubset}}\]\in
\func{\nei{\Xt}}{\left\{\x\right\}}.
\end{equation}
In addition, according to \refcor{corclosureofset0},
$\func{\Cl{\Xt}}{\asubset}$
includes $\asubset$, and hence
\begin{equation}
\[\compl{\X}{\func{\Cl{\Xt}}{\asubset}}\]\cap\asubset=\empty.
\end{equation}
Therefore,
\begin{equation}
\Exists{\U}{\func{\nei{\Xt}}{\left\{\x\right\}}}\U\cap\asubset=\empty.
\end{equation}
Thus according to
\refdef{defadherentpoint},
\begin{equation}
\x\in\(\compl{\X}{\Sadh{\Xt}{\asubset}}\).
\end{equation}
\endp
\end{itemize}
Therefore,
\begin{equation}
\[\compl{\X}{\func{\Cl{\Xt}}{\asubset}}\]\subseteq
\(\compl{\X}{\Sadh{\Xt}{\asubset}}\),
\end{equation}
which means,
\begin{equation}
\Sadh{\Xt}{\asubset}\subseteq
\func{\Cl{\Xt}}{\asubset}.
\end{equation}
\begin{itemize}
\item[${\textbf{\textsf{p2}}}$]
$\p{\x}$
is taken as an arbitrary element of
$\compl{\X}{\Sadh{\Xt}{\asubset}}$.
Considering \refdef{defadherentpoint},
$\p{\U}$
is taken as an element of
$\func{\nei{\Xt}}{\left\{\p{\x}\right\}}$
that,
\begin{equation}
\p{\U}\cap\asubset=\empty.
\end{equation}
Thus considering that $\p{\U}$
is an open set of $\Xt$,
\begin{align}
\(\compl{\X}{\p{\U}}\)&\in\Fclosed{\X}{\topology{}},\\
\(\compl{\X}{\p{\U}}\)&\supseteq\asubset,
\end{align}
and hence
\begin{equation}
\(\compl{\X}{\p{\U}}\)\in\[\Fclosed{\X}{\topology{}}\cap\func{\Cinc{\X}}{\asubset}\].
\end{equation}
Thus according to \refcor{corclosureofset0},
\begin{equation}
\func{\Cl{\Xt}}{\asubset}\subseteq\(\compl{\X}{\p{\U}}\),
\end{equation}
which means,
\begin{equation}
\p{\U}\subseteq\[\compl{\X}{\func{\Cl{\Xt}}{\asubset}}\].
\end{equation}
Thus according to $\p{\x}\in\p{\U}$,
it is clear that,
\begin{equation}
\p{\x}\in\[\compl{\X}{\func{\Cl{\Xt}}{\asubset}}\].
\end{equation}
\endp
\end{itemize}
Therefore,
\begin{equation}
\(\compl{\X}{\Sadh{\Xt}{\asubset}}\)\subseteq
\[\compl{\X}{\func{\Cl{\Xt}}{\asubset}}\],
\end{equation}
which means,
\begin{equation}
\[\func{\Cl{\Xt}}{\asubset}\]\subseteq\(\Sadh{\Xt}{\asubset}\).
\end{equation}
\endthm
%%%%%%%%%%%%%%%%%%%%%%%%%%%%%%%%%%%%%%%%%%%%%%%%%%%%%%%%%%%%%%%%%%%%%%%%%%%%%%%%%%%%%%%%%%%%%%%%%%%%%%%%%%%%%%%%%%%%%%%%%%%%%%%%%%%%%%%%%%%%%%%%%%%%%%%%%%%%%%%%
%%%%%%%%%%%%%%%%%%%%%%%%%%%%%%%%%%%%%%%%%%%%%%%%%%%%%%%%%%%%%%%%%%%%%%%%%%%%%%%%
%%%%%%%%%%%%%%%%%%%%%%%%%%%%%%%%%%%%%%%%%%%%%%%%%%%%%%%%%%%%%%%%%%%%%%%%%%%%%%%%
\theorem\label{thmclosureofsetistheunionofsetandsetoflimitpoints}
$\Xt=\opair{\X}{\topology{}}$
is taken as a topological-space.
\begin{equation}
\Foreach{\asubset}{\CSs{\X}}
\[\func{\Cl{\Xt}}{\asubset}=\asubset\cup\Slim{\Xt}{\asubset}\].
\end{equation}
\prooff
$\asubset$
is taken as an arbitrary subset of $\X$.
\begin{itemize}
\item[${\textbf{\textsf{p1}}}$]
Considering that
$\Slim{\Xt}{\asubset}\subseteq\Sadh{\Xt}{\asubset}$
)according to
\refcor{correlationofpointsetpositions}),
and
$\func{\Cl{\Xt}}{\asubset}=\Sadh{\Xt}{\asubset}$
(according to
\refthm{thmclosureofsetissetofadhpoints}),
it is clear that,
\begin{equation}\label{thmclosureofsetistheunionofsetandsetoflimitpointsp1eq1}
\Slim{\Xt}{\asubset}\subseteq\func{\Cl{\Xt}}{\asubset}.
\end{equation}
Thus considering that
$\asubset\subseteq\func{\Cl{\Xt}}{\asubset}$
(according to
\refcor{corclosureofset0}),
\begin{equation}\label{thmclosureofsetistheunionofsetandsetoflimitpointsp1eq2}
\func{\Cl{\Xt}}{\asubset}\supseteq\asubset\cup\Slim{\Xt}{\asubset}.
\end{equation}
\endp
\item[${\textbf{\textsf{p2}}}$]
$~$
\begin{itemize}
\item[${\textbf{\textsf{p2-1}}}$]
$\point$
is taken as an arbitrary element of
$\func{\Cl{\Xt}}{\asubset}$
Then according to
\refthm{thmclosureofsetissetofadhpoints}
and
\refdef{defadherentpoint},
\begin{equation}\label{thmclosureofsetistheunionofsetandsetoflimitpointsp2-1eq1},
\Foreach{\U}{\func{\nei{\Xt}}{\seta{\point}}}
\U\cap\asubset\neq\empty.
\end{equation}
It is trivial that either $\point$ is an element of $\asubset$ or not an element of $\asubset$. That is,
\begin{equation}\label{thmclosureofsetistheunionofsetandsetoflimitpointsp2-1eq2}
\OR{\(\point\in\asubset\)}{\(\point\notin\asubset\)}.
\end{equation}
In addition, it is clear that,
\begin{equation}\label{thmclosureofsetistheunionofsetandsetoflimitpointsp2-1eq3}
\(\point\in\asubset\)\then\[\point\in\(\asubset\cup\Slim{\Xt}{\asubset}\)\].
\end{equation}
\begin{itemize}
\item[${\textbf{\textsf{p2-1-1}}}$]
It is assumed that,
\begin{equation}
\point\notin\asubset.
\end{equation}
Then according to
\Ref{thmclosureofsetistheunionofsetandsetoflimitpointsp2-1eq1},
\begin{equation}
\Foreach{\U}{\func{\nei{\Xt}}{\seta{\point}}}
\[\U\cap\(\compl{\asubset}{\seta{\point}}\)\]\neq\empty,
\end{equation}
and hence according to
\refdef{deflimitpoint},
\begin{equation}
\point\in\Slim{\Xt}{\asubset}.
\end{equation}
Thus it is clear that,
\begin{equation}
\point\in\[\asubset\cup\Slim{\Xt}{\asubset}\].
\end{equation}
\endp
\end{itemize}
Therefore,
\begin{equation}\label{thmclosureofsetistheunionofsetandsetoflimitpointsp2-1eq4}
\(\point\notin\asubset\)\then\[\point\in\(\asubset\cup\Slim{\Xt}{\asubset}\)\].
\end{equation}
\Ref{thmclosureofsetistheunionofsetandsetoflimitpointsp2-1eq2},
\Ref{thmclosureofsetistheunionofsetandsetoflimitpointsp2-1eq3},
and
\Ref{thmclosureofsetistheunionofsetandsetoflimitpointsp2-1eq4}
imply,
\begin{equation}
\point\in\(\asubset\cup\Slim{\Xt}{\asubset}\).
\end{equation}
\endp
\end{itemize}
Therefore,
\begin{equation}
\Foreach{\point}{\func{\Cl{\Xt}}{\asubset}}
\[\point\in\(\asubset\cup\Slim{\Xt}{\asubset}\)\],
\end{equation}
which means,
\begin{equation}\label{thmclosureofsetistheunionofsetandsetoflimitpointsp2eq1}
\func{\Cl{\Xt}}{\asubset}\subseteq
\(\asubset\cup\Slim{\Xt}{\asubset}\).
\end{equation}
\endp
\end{itemize}
\Ref{thmclosureofsetistheunionofsetandsetoflimitpointsp1eq2}
and
\Ref{thmclosureofsetistheunionofsetandsetoflimitpointsp2eq1}
clearly imply that,
\begin{equation}
\func{\Cl{\Xt}}{\asubset}=
\(\asubset\cup\Slim{\Xt}{\asubset}\).
\end{equation}
\endthm
%%%%%%%%%%%%%%%%%%%%%%%%%%%%%%%%%%%%%%%%%%%%%%%%%%%%%%%%%%%%%%%%%%%%%%%%%%%%%%%%
%%%%%%%%%%%%%%%%%%%%%%%%%%%%%%%%%%%%%%%%%%%%%%%%%%%%%%%%%%%%%%%%%%%%%%%%%%%%%%%%
%%%%%%%%%%%%%%%%%%%%%%%%%%%%%%%%%%%%%%%%%%%%%%%%%%%%%%%%%%%%%%%%%%%%%%%%%%%%%%%%%%%%%%%%%%%%%%%%%%%%%%%%%%%%%%%%%%%%%%%%%%%%%%%%%%%%%%%%%%%%%%%%%%%%%%%%%%%%%%%%
\theorem\label{thmopensetsintersectingclosure}
$\Xt=\opair{\X}{\topology{}}$
is taken as a topological-space, and
$\asubset$
as a subset of
$\X$.
Every open set of
$\Xt$
that intersects
$\func{\Cl{\Xt}}{\asubset}$, also intersects $\asubset$. That is,
\begin{equation}
\Foreach{\U}{\topology{}}
\[\U\cap\func{\Cl{\Xt}}{\asubset}\neq\empty\then\U\cap\asubset\neq\empty\].
\end{equation}
\prooff
$\U$
is taken as such an element of
$\topology{}$ that,
\begin{equation}
\U\cap\func{\Cl{\Xt}}{\asubset}\neq\empty.
\end{equation}
Then,
\begin{equation}
\Exists{\x}{\X}\[\x\in\U\cap\func{\Cl{\Xt}}{\asubset}\].
\end{equation}
$\point$
is taken as an element of
$\U\cap\func{\Cl{\Xt}}{\asubset}$.
Hence,
$\U$
is a neighbourhood of
$\seta{\point}$
in
$\Xt$,
and
$\point$
is an element of
$\func{\Cl{\Xt}}{\asubset}$.
Thus according to
\refthm{thmclosureofsetissetofadhpoints},
\begin{equation}
\U\cap\asubset\neq\empty.
\end{equation}
\endthm
%%%%%%%%%%%%%%%%%%%%%%%%%%%%%%%%%%%%%%%%%%%%%%%%%%%%%%%%%%%%%%%%%%%%%%%%%%%%%%%%%
\theorem\label{thmclosureofclosedset}
$\Xt=\opair{\X}{\topology{}}$
is taken as a topological-space. For every
$\asubset$
in
$\CSs{\X}$,
$\asubset$
is a closed set of $\Xt$ if-and-only-if
$\func{\Cl{\Xt}}{\asubset}$
equals
$\asubset$.
That is,
\begin{equation}
\Foreach{\asubset}{\CSs{\X}}
\[\asubset\in\Fclosed{\X}{\topology{}}\thenn\func{\Cl{\Xt}}{\asubset}=\asubset\].
\end{equation}
\prooff
$\asubset$
is taken as an arbitrary subset of $\X$.
\begin{itemize}
\item[${\textbf{\textsf{p1}}}$]
It is assumed that
$\func{\Cl{\Xt}}{\asubset}=\asubset$.
Then according to
\refcor{corclosureofset0}
it is clear that,
\begin{equation}
\asubset\in\Fclosed{\X}{\topology{}}.
\end{equation}
\endp
\item[${\textbf{\textsf{p2}}}$]
It is assumed that
$\asubset\in\Fclosed{\X}{\topology{}}$.
Then considering that
$\asubset\supseteq\asubset$,
it is clear that
\begin{equation}
\asubset\in\[\Fclosed{\X}{\topology{}}\cap\func{\Cinc{\X}}{\asubset}\],
\end{equation}
and hence according to
\refcor{corclosureofset0},
\begin{equation}
\asubset\supseteq\func{\Cl{\Xt}}{\asubset}.
\end{equation}
In addition, according to
\refcor{corclosureofset0},
\begin{equation}
\func{\Cl{\Xt}}{\asubset}\supseteq\asubset.
\end{equation}
Therefore,
\begin{equation}
\func{\Cl{\Xt}}{\asubset}=\asubset.
\end{equation}
\endp
\end{itemize}
\endthm
%%%%%%%%%%%%%%%%%%%%%%%%%%%%%%%%%%%%%%%%%%%%%%%%%%%%%%%%%%%%%%%%%%%%%%%%%%%%%%%%
\theorem\label{thmsetisclosediffSlimisasubsetofset}
$\Xt=\opair{\X}{\topology{}}$
is taken as a topological-space. For every
$\asubset$
in
$\CSs{\X}$,
$\asubset$
is a closed set of
$\Xt$
if-and-only-if
$\asubset$
contains all limit-points of $\asubset$ in $\Xt$. That is,
\begin{equation}
\Foreach{\asubset}{\CSs{\X}}
\left\{\[\asubset\in\Fclosed{\X}{\topology{}}\]\thenn\(\Slim{\Xt}{\asubset}\subseteq\asubset\)\right\}.
\end{equation}
\prooff
$\asubset$
is taken as an arbitrary element of $\X$.
According to \refthm{thmclosureofclosedset},
\begin{equation}
\[\asubset\in\Fclosed{\X}{\topology{}}\]\thenn
\[\func{\Cl{\Xt}}{\asubset}=\asubset\].
\end{equation}
In addition, according to
\refthm{thmclosureofsetistheunionofsetandsetoflimitpoints},
\begin{equation}
\[\func{\Cl{\Xt}}{\asubset}=\asubset\]
\thenn
\(\asubset=\asubset\cup\Slim{\Xt}{\asubset}\).
\end{equation}
Additionally, it is clear that,
\begin{equation}
\(\asubset=\asubset\cup\Slim{\Xt}{\asubset}\)\thenn
\(\Slim{\Xt}{\asubset}\subseteq\asubset\).
\end{equation}
Therefore,
\begin{equation}
\[\asubset\in\Fclosed{\X}{\topology{}}\]\thenn\(\Slim{\Xt}{\asubset}\subseteq\asubset\).
\end{equation}
\endthm
%%%%%%%%%%%%%%%%%%%%%%%%%%%%%%%%%%%%%%%%%%%%%%%%%%%%%%%%%%%%%%%%%%%%%%%%%%%%%%%%
\theorem\label{thmclosureandextrelation}
$\Xt=\opair{\X}{\topology{}}$
is taken as a topological-space.
\begin{equation}
\Foreach{\asubset}{\CSs{\X}}
\func{\Cl{\Xt}}{\asubset}=\compl{\X}{\[\func{\Ext{\Xt}}{\asubset}\]}.
\end{equation}
\prooff
According to
\refcor{corSextSadh},
\refthm{thmextofsetissetofextpoints},
and
\refthm{thmclosureofsetissetofadhpoints},
it is clear.
\endthm
%%%%%%%%%%%%%%%%%%%%%%%%%%%%%%%%%%%%%%%%%%%%%%%%%%%%%%%%%%%%%%%%%%%%%%%%%%%%%%%%%%%%%%%%%%%%%%%%%%%%%%%%
\theorem\label{thmclosureofemptyset}
$\Xt=\opair{\X}{\topology{}}$
is taken as a topological-space.
\begin{equation}
\func{\Cl{\Xt}}{\empty}=\empty.
\end{equation}
\prooff
According to
\refthm{thmclosureandextrelation}
and
\refthm{thmextofemptyset},
\begin{align}
\func{\Cl{\Xt}}{\empty}&=\compl{\X}{\func{\Ext{\Xt}}{\empty}}\cr
&=\compl{\X}{\X}\cr
&=\empty.
\end{align}
\endthm
%%%%%%%%%%%%%%%%%%%%%%%%%%%%%%%%%%%%%%%%%%%%%%%%%%%%%%%%%%%%%%%%%%%%%%%%%%%%%%%%%%%%%%%%%%%%%%%%%%%%%%%%
\theorem\label{thmclosureofuniversalset}
$\Xt=\opair{\X}{\topology{}}$
is taken as a topological-space.
\begin{equation}
\func{\Cl{\Xt}}{\X}=\X.
\end{equation}
\prooff
According to
\refthm{thmclosureandextrelation}
and
\refthm{thmextofuniversalset},
\begin{align}
\func{\Cl{\Xt}}{\X}&=\compl{\X}{\func{\Ext{\Xt}}{\X}}\cr
&=\compl{\X}{\empty}\cr
&=\X.
\end{align}
\endthm
%%%%%%%%%%%%%%%%%%%%%%%%%%%%%%%%%%%%%%%%%%%%%%%%%%%%%%%%%%%%%%%%%%%%%%%%%%%%%%%%%%%%%%%%%%%%%%%%%%%%%%%%
%%%%%%%%%%%%%%%%%%%%%%%%%%%%%%%%%%%%%%%%%%%%%%%%%%%%%%%%%%%%%%%%%%%%%%%%%%%%%%%%%%%
\theorem\label{thmclosureofasubsetofset}
$\Xt=\opair{\X}{\topology{}}$
is taken as a topological-space.
\begin{equation}
\Foreach{\opair{\asubset}{\p{\asubset}}}{\CSs{\X}\times\CSs{\X}}
\left\{\(\p{\asubset}\subseteq\asubset\)\then
\[\func{\Cl{\Xt}}{\p{\asubset}}\subseteq\func{\Cl{\Xt}}{\asubset}\]\right\}.
\end{equation}
\prooff
According to
\refthm{thmextofasubsetofset}
and
\refthm{thmclosureandextrelation},
it is clear.
\endthm
%%%%%%%%%%%%%%%%%%%%%%%%%%%%%%%%%%%%%%%%%%%%%%%%%%%%%%%%%%%%%%%%%%%%%%%%%%%%%%%%%%%%%%%
\theorem\label{thmclosureofunionofsets0}
$\Xt=\opair{\X}{\topology{}}$
is taken as a topological-space.
\begin{align}
\Foreach{\sCi}{\[\compl{\CSs{\CSs{\X}}}{\seta{\empty}}\]}
\[\func{\Cl{\Xt}}{\union{\sCi}}\supseteq
\bigcup_{\asubset\in\sCi}\func{\Cl{\Xt}}{\asubset}\].
\end{align}
\prooff
$\sCi$
is taken as a non-empty collection of subsets of $\X$.
According to
\refthm{thmclosureandextrelation}
and
\refthm{thmextofunionofsets0},
\begin{align}
\func{\Cl{\Xt}}{\union{\sCi}}&=
\compl{\X}{\func{\Ext{\Xt}}{\union{\sCi}}}\cr
&\supseteq\compl{\X}{\[\bigcap_{\asubset\in\sCi}\func{\Ext{\Xt}}{\asubset}\]}\cr
&=\bigcup_{\asubset\in\sCi}\[\compl{\X}{\func{\Ext{\Xt}}{\asubset}}\]\cr
&=\bigcup_{\asubset\in\sCi}\func{\Cl{\Xt}}{\asubset}.
\end{align}
\endthm
%%%%%%%%%%%%%%%%%%%%%%%%%%%%%%%%%%%%%%%%%%%%%%%%%%%%%%%%%%%%%%%%%%%%%%%%%%%%%%%%%%%%%%%
\theorem\label{thmclosureofunionofsets}
$\Xt=\opair{\X}{\topology{}}$
is taken as a topological-space.
\begin{equation}
\Foreach{\opair{\asubset}{\p{\asubset}}}{\CSs{\X}\times\CSs{\X}}
\[\func{\Cl{\Xt}}{\asubset\cup\p{\asubset}}=
\func{\Cl{\Xt}}{\asubset}\cup\func{\Cl{\Xt}}{\p{\asubset}}\].
\end{equation}
\prooff
$\asubset$
and
$\p{\asubset}$
are taken as arbitrary subsets of $\X$.
According to \refthm{thmextofunionofsets}
and
\refthm{thmclosureandextrelation},
\begin{align}
\func{\Cl{\Xt}}{\asubset\cup\p{\asubset}}&=
\compl{\X}{\[\func{\Ext{\Xt}}{\asubset\cup\p{\asubset}}\]}\cr
&=\compl{\X}{\[\func{\Ext{\Xt}}{\asubset}\cap\func{\Ext{\Xt}}{\p{\asubset}}\]}\cr
&=\[\compl{\X}{\func{\Ext{\Xt}}{\asubset}}\]\cup
\[\compl{\X}{\func{\Ext{\Xt}}{\p{\asubset}}}\]\cr
&=\func{\Cl{\Xt}}{\asubset}\cup\func{\Cl{\Xt}}{\p{\asubset}}.
\end{align}
\endthm
%%%%%%%%%%%%%%%%%%%%%%%%%%%%%%%%%%%%%%%%%%%%%%%%%%%%%%%%%%%%%%%%%%%%%%%%%%%%%%%%%%%%%%%%%%
\theorem\label{thmclosureofintersectionofsets}
$\Xt=\opair{\X}{\topology{}}$
is taken as a topological-space.
\begin{equation}
\Foreach{\sCi}{\[\compl{\CSs{\CSs{\X}}}{\seta{\empty}}\]}
\[\func{\Cl{\Xt}}{\intersection{\sCi}}\subseteq
\bigcap_{\asubset\in\sCi}\func{\Cl{\Xt}}{\asubset}\].
\end{equation}
\prooff
$\sCi$
is taken as a non-empty collection of subsets of $\X$.
According to \refthm{thmextofintersectionofsets}
and
\refthm{thmclosureandextrelation},
\begin{align}
\compl{\X}{\func{\Cl{\Xt}}{\intersection{\sCi}}}&=
\func{\Ext{\Xt}}{\intersection{\sCi}}\cr
&\supseteq\bigcup_{\asubset\in\sCi}\func{\Ext{\Xt}}{\asubset}\cr
&=\bigcup_{\asubset\in\sCi}\[\compl{\X}{\func{\Cl{\Xt}}{\asubset}}\].
\end{align}
Thus,
\begin{align}
\func{\Cl{\Xt}}{\intersection{\sCi}}&\subseteq
\compl{\X}{\bigcup_{\asubset\in\sCi}\[\compl{\X}{\func{\Cl{\Xt}}{\asubset}}\]}\cr
&=\bigcap_{\asubset\in\sCi}\func{\Cl{\Xt}}{\asubset}.
\end{align}
\endthm
%%%%%%%%%%%%%%%%%%%%%%%%%%%%%%%%%%%%%%%%%%%%%%%%%%%%%%%%%%%%%%%%%%%%%%%%%%%%%%%%%%%%%%%
\theorem\label{thmclosureofclosureofset}
$\Xt=\opair{\X}{\topology{}}$
is taken as a topological-space.
\begin{equation}
\Foreach{\asubset}{\CSs{\X}}
\[\func{\Cl{\Xt}}{\func{\Cl{\Xt}}{\asubset}}=
\func{\Cl{\Xt}}{\asubset}\].
\end{equation}
\prooff
According to
\refthm{thmintofintofset},
\refdef{defextofset},
and
\refthm{thmclosureandextrelation},
\begin{align}
\Foreach{\asubset}{\CSs{\X}}
\func{\Cl{\Xt}}{\func{\Cl{\Xt}}{\asubset}}&=
\func{\Cl{\Xt}}{\compl{\X}{\func{\Ext{\Xt}}{\asubset}}}\cr
&=\compl{\X}{\[\func{\Ext{\Xt}}{\compl{\X}{\func{\Ext{\Xt}}{\asubset}}}\]}\cr
&=\compl{\X}{\[\func{\Int{\Xt}}{\func{\Ext{\Xt}}{\asubset}}\]}\cr
&=\compl{\X}{\[\func{\Int{\Xt}}{\func{\Int{\Xt}}{\compl{\X}{\asubset}}}\]}\cr
&=\compl{\X}{\[\func{\Int{\Xt}}{\compl{\X}{\asubset}}\]}\cr
&=\compl{\X}{\[\func{\Ext{\Xt}}{\asubset}\]}\cr
&=\func{\Cl{\Xt}}{\asubset}.
\end{align}
\endthm
%%%%%%%%%%%%%%%%%%%%%%%%%%%%%%%%%%%%%%%%%%%%%%%%%%%%%%%%%%%%%%%%%%%%%%%%%%%%%%%%%%%%%%%%%%
\theorem\label{thmsubspaceclosure}
$\opair{\X}{\topology{}}$
is taken as a topological-space, and
$\Y$
as a subset of $\X$.
\begin{equation}
\Foreach{\asubset}{\CSs{\Y}}
\[\func{\Cl{\opair{\Y}{\stopology{\topology{}}{\Y}}}}{\asubset}=
\Y\cap\func{\Cl{\Xt}}{\asubset}\].
\end{equation}
\prooff
$\asubset$
is taken as an arbitrary subset of $\Y$.
According to
\refthm{thmsubspaceclosedsets},
\begin{align}\label{thmsubspaceclosurep3}
\Fclosed{\Y}{\stopology{\topology{}}{\Y}}=
\defset{\V}{\CSs{\Y}}{\[\Exists{\U}{\Fclosed{\X}{\topology{}}}\V=\Y\cap\U\]}.
\end{align}
Thus considering that,
\begin{equation}\label{thmsubspaceclosurep2}
\func{\Cinc{\Y}}{\asubset}=\[\func{\Cinc{\X}}{\asubset}\cap\CSs{\Y}\].
%\defset{\U}{\func{\Cinc{\X}}{\asubset}}{\U\subseteq\Y},
%\defset{\V}{\CSs{\Y}}{\Exists{\U}{\func{\Cinc{\X}}{\asubset}}\V=\Y\cap\U},
\end{equation}
it can be easily seen that,
\begin{align}
\[\Fclosed{\Y}{\stopology{\topology{}}{\Y}}\cap\func{\Cinc{\Y}}{\asubset}\]=
\defset{\V}{\CSs{\Y}}{\[\Exists{\U}{\Fclosed{\X}{\topology{}}\cap\func{\Cinc{\X}}{\asubset}}\V=\Y\cap\U\]}.
\end{align}
Thus according to \refdef{defclosureofset},
\begin{align}
\func{\Cl{\opair{\Y}{\stopology{\topology{}}{\Y}}}}{\asubset}&=
\intersection{\[\Fclosed{\Y}{\stopology{\topology{}}{\Y}}\cap\func{\Cinc{\Y}}{\asubset}\]}\cr
&=\Y\cap\left\{\intersection{\[\Fclosed{\X}{\topology{}}\cap\func{\Cinc{\X}}{\asubset}\]}\right\}\cr
&=\Y\cap\func{\Cl{\Xt}}{\asubset}.
\end{align}
\endthm
%%%%%%%%%%%%%%%%%%%%%%%%%%%%%%%%%%%%%%%%%%%%%%%%%%%%%%%%%%%%%%%%%%%%%%%%%%%%%%%%%%%%%%%%%%
\theorem\label{thmclosureofsetwithdifferenttopologies}
Each
$\Xt_{1}=\opair{\X}{\topology{1}}$
and
$\Xt_{2}=\opair{\X}{\topology{2}}$
is taken as a topological-space.
\begin{equation}
\(\topology{2}\subseteq\topology{1}\)\then\[\Foreach{\asubset}{\CSs{\X}}
\func{\Cl{\Xt_{2}}}{\asubset}\supseteq\func{\Cl{\Xt_{1}}}{\asubset}\].
\end{equation}
\prooff
According to
\refthm{thmextofsetwithdifferenttopologies}
and
\refthm{thmclosureandextrelation},
it is clear.
\endthm
%%%%%%%%%%%%%%%%%%%%%%%%%%%%%%%%%%%%%%%%%%%%%%%%%%%%%%%%%%%%%%%%%%%%%%%%%%%%%%%%%%%%%%%%%%
\definition\label{defperfectset}
$\Xt=\opair{\X}{\topology{}}$
is taken as a topological-space. For every
$\asubset$
in
$\CSs{\X}$,
$\asubset$
is referred to as a $\quotl$perfect set of $\Xt$$\quotr$ iff
$\asubset$
is a closed and dense-in-itself set of $\Xt$.
\endef
%%%%%%%%%%%%%%%%%%%%%%%%%%%%%%%%%%%%%%%%%%%%%%%%%%%%%%%%%%%%%%%%%%%%%%%%%%%%%%%%%%%%%%%%%%
\corollary\label{corperfectset0}
$\Xt=\opair{\X}{\topology{}}$
is taken as a topological-space. For every $\asubset$
in
$\CSs{\X}$,
$\asubset$
is a perfect set of $\Xt$ if-and-only-if,
\begin{equation}
\AND{\[\asubset\in\Fclosed{\X}{\topology{}}\]}{\[\Siso{\Xt}{\asubset}=\empty\]}.
\end{equation}
\endcor
%%%%%%%%%%%%%%%%%%%%%%%%%%%%%%%%%%%%%%%%%%%%%%%%%%%%%%%%%%%%%%%%%%%%%%%%%%%%%%%%%%%%%%%%%%
\theorem\label{thmperfectsetisthesetofalllimitpointsofset}
$\Xt=\opair{\X}{\topology{}}$
is taken as a topological-space. For every
$\asubset$
in
$\CSs{\X}$,
$\asubset$
is a perfect set of $\Xt$ if-and-only-if
\begin{equation}
\asubset=\Slim{\Xt}{\asubset}.
\end{equation}
\prooff
$\asubset$
is taken as an arbitrary subset of $\X$.
\begin{itemize}
\item[${\textbf{\textsf{p1}}}$]
It is assumed that $\asubset$ is a perfect set of $\Xt$.
Then according to \refcor{corperfectset0},
\begin{align}
&\asubset\in\Fclosed{\Xt}{\topology{}},\label{thmperfectsetisthesetofalllimitpointsofsetp1eq1}\\
&\Siso{\Xt}{\asubset}=\empty.\label{thmperfectsetisthesetofalllimitpointsofsetp1eq2}
\end{align}
\Ref{thmperfectsetisthesetofalllimitpointsofsetp1eq1}
and
\refthm{thmsetisclosediffSlimisasubsetofset}
imply,
\begin{equation}
\Slim{\Xt}{\asubset}\subseteq\asubset,
\end{equation}
and hence,
\begin{equation}
\asubset\cap\Slim{\Xt}{\asubset}=\Slim{\Xt}{\asubset},
\end{equation}
and hence according to \refcor{corlimitpointsandisolatedpointsaredistinct},
\begin{align}
\asubset&=\(\asubset\cap\Slim{\Xt}{\asubset}\)\sqcup\Siso{\Xt}{\asubset}\cr
&=\Slim{\Xt}{\asubset}\sqcup\empty\cr
&=\Slim{\Xt}{\asubset}.
\end{align}
\endp
\item[${\textbf{\textsf{p2}}}$]
It is assumed that $\asubset=\Slim{\Xt}{\asubset}$.
Then according to
\refthm{corlimitpointsandisolatedpointsaredistinct},
\begin{align}\label{thmperfectsetisthesetofalllimitpointsofsetp2eq1}
\Siso{\Xt}{\asubset}&=\compl{\asubset}{\Slim{\Xt}{\asubset}}\cr
&=\compl{\asubset}{\asubset}\cr
&=\empty.
\end{align}
In addition, according to \refthm{thmsetisclosediffSlimisasubsetofset}
(considering that $\Slim{\Xt}{\asubset}\subseteq\asubset$),
\begin{equation}
\asubset\in\Fclosed{\X}{\topology{}}.
\end{equation}
Thus according to \refcor{corperfectset0},
$\asubset$
is a perfect set of $\Xt$.
\endp
\end{itemize}
\endthm
%%%%%%%%%%%%%%%%%%%%%%%%%%%%%%%%%%%%%%%%%%%%%%%%%%%%%%%%%%%%%%%%%%%%%%%%%%%%%%%%%%%%%%%%%%
\definition\label{defsparatedsets}
$\Xt=\opair{\X}{\topology{}}$
is taken as a topological-space, and each
$\asubset$
and
$\bsubset$
as a subset of $\X$.
\begin{itemize}
\item
It is said that $\quotl$$\asubset$ and $\bsubset$ are glued together
in the topological space $\Xt$$\quotr$ iff
\begin{equation*}
\asubset\cap\func{\Cl{\Xt}}{\bsubset}\neq\empty,~
\bsubset\cap\func{\Cl{\Xt}}{\asubset}\neq\empty.
\end{equation*}
\item
It is said that $\quotl$$\asubset$ and $\bsubset$ are free-of-eachother in the topological-space
$\Xt$$\quotr$ iff
\begin{equation*}
\asubset\cap\func{\Cl{\Xt}}{\bsubset}=
\bsubset\cap\func{\Cl{\Xt}}{\asubset}=\empty.
\end{equation*}
\end{itemize}
\endef
%%%%%%%%%%%%%%%%%%%%%%%%%%%%%%%%%%%%%%%%%%%%%%%%%%%%%%%%%%%%%%%%%%%%%%%%%%%%%%%%%%%%%%%%%%
%%%%%%%%%%%%%%%%%%%%%%%%%%%%%%%%%%%%%%%%%%%%%%%%%%%%%%%%%%%%%%%%%%%%%%%%%%%%%%%%%%%%%%%%%%
%%%%%%%%%%%%%%%%%%%%%%%%%%%%%%%%%%%%%%%%%%%%%%%%%%%%%%%%%%%%%%%%%%%%%%%%%%%%%%%%%%%%%%%%%%
%%%%%%%%%%%%%%%%%%%%%%%%%%%%%%%%%%%%%%%%%%%%%%%%%%%%%%%%%%%%%%%%%%%%%%%%%%%%%%%%%%%%%%%%%%
%%%%%%%%%%%%%%%%%%%%%%%%%%%%%%%%%%%%%%%%%%%%%%%%%%%%%%%%%%%%%%%%%%%%%%%%%%%%%%%%%%%%%%%%%%
%%%%%%%%%%%%%%%%%%%%%%%%%%%%%%%%%%%%%%%%%%%%%%%%%%%%%%%%%%%%%%%%%%%%%%%%%%%%%%%%%%%%%%%%%%
\subsection{
Boundary of a Subset of a Topological Space}
\definition\label{deffrontierofset}
$\Xt=\opair{\X}{\topology{}}$
is taken as a topological-space.
The mapping $\Fr{\Xt}$ is defined as,
\begin{itemize}
\item[${\textbf{\textsf{Fr1}}}$]
\hfill
$\Fr{\Xt}\indef\Func{\CSs{\X}}{\CSs{\X}}.$
\item[${\textbf{\textsf{Fr2}}}$]
\hfill
$\Foreach{\asubset}{\CSs{\X}}
\func{\Fr{\Xt}}{\asubset}\eqdef
\[\compl{\func{\Cl{\Xt}}{\asubset}}{\func{\Int{\Xt}}{\asubset}}\].$
\end{itemize}
For every
$\asubset$
in
$\CSs{\X}$,
$\func{\Fr{\Xt}}{\asubset}$
is referred to as the $\quotl$boundary of $\asubset$ in the topological-space $\Xt$$\quotr$.
\endef
%%%%%%%%%%%%%%%%%%%%%%%%%%%%%%%%%%%%%%%%%%%%%%%%%%%%%%%%%%%%%%%%%%%%%%%%%%%%%%%%
\theorem\label{thmfronierofsetissetofboundarypoints}
$\Xt=\opair{\X}{\topology{}}$
is taken as a topological-space.
\begin{equation}
\Foreach{\asubset}{\CSs{\X}}
\func{\Fr{\Xt}}{\asubset}=\Sbound{\Xt}{\asubset}.
\end{equation}
\prooff
According to
\refdef{definteriorpoint},
\refdef{defboundarypoint},
\refdef{defadherentpoint},
\refthm{thmintofsetissetofintpoints},
\refthm{thmclosureofsetissetofadhpoints},
and
\refdef{deffrontierofset},
it is clear.
\endthm
%%%%%%%%%%%%%%%%%%%%%%%%%%%%%%%%%%%%%%%%%%%%%%%%%%%%%%%%%%%%%%%%%%%%%%%%%%%%%%%%%
%%%%%%%%%%%%%%%%%%%%%%%%%%%%%%%%%%%%%%%%%%%%%%%%%%%%%%%%%%%%%%%%%%%%%%%%%%%%%%%%%%%%%%%%%%%%%%%%%%%%%%%%
\theorem\label{thmclosurefrontierrelation}
$\Xt=\opair{\X}{\topology{}}$
is taken as a topological-space.
\begin{equation}
\Foreach{\asubset}{\CSs{\X}}
\[\func{\Cl{\Xt}}{\asubset}=
\func{\Int{\Xt}}{\asubset}\cup\func{\Fr{\Xt}}{\asubset}=
\asubset\cup\func{\Fr{\Xt}}{\asubset}\].
\end{equation}
\prooff
$\asubset$
is taken as an arbitrary subset of $\X$.
Considering that
$\func{\Int{\Xt}}{\asubset}\subseteq\func{\Cl{\Xt}}{\asubset}$,
according to
\refdef{deffrontierofset},
\begin{align}
\func{\Cl{\Xt}}{\asubset}&=
\[\compl{\func{\Cl{\Xt}}{\asubset}}
{\func{\Int{\Xt}}{\asubset}}\]\cup\func{\Int{\Xt}}{\asubset}\cr
&=\func{\Int{\Xt}}{\asubset}\cup\func{\Fr{\Xt}}{\asubset}.
\end{align}
Thus according to
$\func{\Int{\Xt}}{\asubset}\subseteq\asubset\subseteq\func{\Cl{\Xt}}{\asubset}$,
it is clear that,
\begin{equation}
\func{\Cl{\Xt}}{\asubset}=\asubset\cup\func{\Fr{\Xt}}{\asubset}
\end{equation}
\endthm
%%%%%%%%%%%%%%%%%%%%%%%%%%%%%%%%%%%%%%%%%%%%%%%%%%%%%%%%%%%%%%%%%%%%%%%%%%%%%%%%%%%%%%%%%%%%%%%%%%%%%%%%
\theorem\label{thmfrontierofemptyset}
$\Xt=\opair{\X}{\topology{}}$
is taken as a topological-space.
\begin{equation}
\func{\Fr{\Xt}}{\empty}=\empty.
\end{equation}
\prooff
According to
\refdef{deffrontierofset},
\refthm{thmintofemptyset},
and
\refthm{thmclosureofemptyset},
\begin{align}
\func{\Fr{\Xt}}{\empty}&=\compl{\func{\Cl{\Xt}}{\empty}}{\func{\Int{\Xt}}{\empty}}\cr
&=\compl{\empty}{\empty}\cr
&=\empty.
\end{align}
\endthm
%%%%%%%%%%%%%%%%%%%%%%%%%%%%%%%%%%%%%%%%%%%%%%%%%%%%%%%%%%%%%%%%%%%%%%%%%%%%%%%%%
\theorem\label{thmfrontierofuniversalset}
$\Xt=\opair{\X}{\topology{}}$
is taken as a topological-space.
\begin{equation}
\func{\Fr{\Xt}}{\X}=\empty.
\end{equation}
\prooff
According to
\refdef{deffrontierofset},
\refthm{thmintofuniversalset},
and
\refthm{thmclosureofuniversalset},
\begin{align}
\func{\Fr{\Xt}}{\X}&=\compl{\func{\Cl{\Xt}}{\X}}{\func{\Int{\Xt}}{\X}}\cr
&=\compl{\X}{\X}\cr
&=\empty.
\end{align}
\endthm
%%%%%%%%%%%%%%%%%%%%%%%%%%%%%%%%%%%%%%%%%%%%%%%%%%%%%%%%%%%%%%%%%%%%%%%%%%%%%%%%%
\theorem\label{thmfrontierofclosedset}
$\Xt=\opair{\X}{\topology{}}$
is taken as a topological-space. For every
$\asubset$
in
$\CSs{\X}$,
$\asubset$
is a closed set of
$\Xt$
if-and-only-if
$\func{\Fr{\Xt}}{\asubset}$
is a subset of $\asubset$. That is,
\begin{equation}
\Foreach{\asubset}{\CSs{\X}}
\[\asubset\in\Fclosed{\X}{\topology{}}\thenn\func{\Fr{\Xt}}{\asubset}\subseteq\asubset\].
\end{equation}
\prooff
$\asubset$
is taken as an arbitrary subset of $\X$.
\begin{itemize}
\item[${\textbf{\textsf{p1}}}$]
It is assumed that
$\asubset\in\Fclosed{\X}{\topology{}}$.
Then according to
\refthm{thmclosureofclosedset},
it is clear that,
\begin{equation}
\func{\Cl{\Xt}}{\asubset}=\asubset.
\end{equation}
Thus according to
\refdef{deffrontierofset},
\begin{align}
\func{\Fr{\Xt}}{\asubset}&=\[\compl{\asubset}{\func{\Int{\Xt}}{\asubset}}\]\cr
&\subseteq\asubset.
\end{align}
\endp
\item[${\textbf{\textsf{p2}}}$]
It is assumed that
$\func{\Fr{\Xt}}{\asubset}\subseteq\asubset$.
Then according to
\refdef{deffrontierofset},
\begin{equation}
\[\compl{\func{\Cl{\Xt}}{\asubset}}{\func{\Int{\Xt}}{\asubset}}\]\subseteq\asubset.
\end{equation}
Thus considering that
$\func{\Int{\Xt}}{\asubset}\subseteq\asubset$
it is clear that,
\begin{equation}
\func{\Cl{\Xt}}{\asubset}\subseteq\asubset,
\end{equation}
and hence considering that
$\asubset\subseteq\func{\Cl{\Xt}}{\asubset}$,
\begin{equation}
\func{\Cl{\Xt}}{\asubset}=\asubset.
\end{equation}
Thus according to
\refthm{thmclosureofclosedset},
\begin{equation}
\asubset\in\Fclosed{\X}{\topology{}}.
\end{equation}
\endp
\end{itemize}
\endthm
%%%%%%%%%%%%%%%%%%%%%%%%%%%%%%%%%%%%%%%%%%%%%%%%%%%%%%%%%%%%%%%%%%%%%%%%%%%%%%%%%%%%%%%%%%%%%%%%
\theorem\label{thmfrontier1}
$\Xt=\opair{\X}{\topology{}}$
is taken as a topological-space.
\begin{equation}
\Foreach{\asubset}{\CSs{\X}}
\[\func{\Fr{\Xt}}{\asubset}=
\func{\Cl{\Xt}}{\asubset}\cap\func{\Cl{\Xt}}{\compl{\X}{\asubset}}\].
\end{equation}
\prooff
$\asubset$
is taken as an arbitrary subset of
$\X$.
According to \refdef{defextofset}
and
\refthm{thmclosureandextrelation},
\begin{align}
\func{\Int{\Xt}}{\asubset}&=\func{\Ext{\Xt}}{\compl{\X}{\asubset}}\cr
&=\compl{\X}{\func{\Cl{\Xt}}{\compl{\X}{\asubset}}}.
\end{align}
Thus according to
\refdef{deffrontierofset},
\begin{align}
\func{\Fr{\Xt}}{\asubset}&=
\compl{\func{\Cl{\Xt}}{\asubset}}{\func{\Int{\Xt}}{\asubset}}\cr
&=\compl{\func{\Cl{\Xt}}{\asubset}}{\[\compl{\X}{\func{\Cl{\Xt}}{\compl{\X}{\asubset}}}\]}\cr
&=\func{\Cl{\Xt}}{\asubset}\cap\func{\Cl{\Xt}}{\compl{\X}{\asubset}}.
\end{align}
\endthm
%%%%%%%%%%%%%%%%%%%%%%%%%%%%%%%%%%%%%%%%%%%%%%%%%%%%%%%%%%%%%%%%%%%%%%%%%%%%%%%%%%%%%%%%%%%%%%%%%%
\theorem\label{thmfrontierofcomplementofset}
$\Xt=\opair{\X}{\topology{}}$
is taken as a topological-space.
\begin{equation}
\Foreach{\asubset}{\CSs{\X}}
\[\func{\Fr{\Xt}}{\compl{\X}{\asubset}}=\func{\Fr{\Xt}}{\asubset}\].
\end{equation}
\prooff
$\asubset$
is taken as an arbitrary subset of $\X$.
According to \refthm{thmfrontier1},
\begin{align}
\func{\Fr{\Xt}}{\compl{\X}{\asubset}}&=
\func{\Cl{\Xt}}{\compl{\X}{\asubset}}\cap
\func{\Cl{\Xt}}{\compl{\X}{\[\compl{\X}{\asubset}\]}}\cr
&=\func{\Cl{\Xt}}{\compl{\X}{\asubset}}\cap
\func{\Cl{\Xt}}{\asubset}\cr
&=\func{\Cl{\Xt}}{\asubset}\cap\func{\Cl{\Xt}}{\compl{\X}{\asubset}}\cr
&=\func{\Fr{\Xt}}{\asubset}.
\end{align}
\endthm
%%%%%%%%%%%%%%%%%%%%%%%%%%%%%%%%%%%%%%%%%%%%%%%%%%%%%%%%%%%%%%%%%%%%%%%%%%%%%%%%%%%%%%%%%%%%%%%%%%%%%%%%%
\theorem\label{thmpartitionofuniversalset}
$\Xt=\opair{\X}{\topology{}}$
is taken as a topological-space.
\begin{align}
\Foreach{\asubset}{\CSs{\X}}
\left\{
\begin{aligned}
&\func{\Int{\Xt}}{\asubset}\cap\func{\Ext{\Xt}}{\asubset}=\empty,\\
&\func{\Int{\Xt}}{\asubset}\cap\func{\Fr{\Xt}}{\asubset}=\empty,\\
&\func{\Ext{\Xt}}{\asubset}\cap\func{\Fr{\Xt}}{\asubset}=\empty,\\
&\X=\func{\Int{\Xt}}{\asubset}\cup\func{\Fr{\Xt}}{\asubset}\cup\func{\Ext{\Xt}}{\asubset}.
\end{aligned}
\right.
\end{align}
\prooff
$\asubset$
is taken as an arbitrary subset of
$\X$.
Considering that
$\func{\Int{\Xt}}{\asubset}\subseteq\asubset$,
and
$\func{\Int{\Xt}}{\compl{\X}{\asubset}}\subseteq\(\compl{\X}{\asubset}\)$,
according to \refdef{defextofset},
\begin{align}
\func{\Int{\Xt}}{\asubset}\cap\func{\Ext{\Xt}}{\asubset}&=
\func{\Int{\Xt}}{\asubset}\cap\func{\Int{\Xt}}{\compl{\X}{\asubset}}\cr
&=\empty.
\end{align}
According to
\refdef{deffrontierofset},
\begin{align}
\func{\Int{\Xt}}{\asubset}\cap\func{\Fr{\Xt}}{\asubset}&=
\func{\Int{\Xt}}{\asubset}\cap\(\compl{\func{\Cl{\Xt}}{\asubset}}{\func{\Int{\Xt}}{\asubset}}\)\cr
&=\empty.
\end{align}
According to
\refdef{defextofset},
\refdef{deffrontierofset},
and
\refthm{thmfrontierofcomplementofset},
\begin{align}
\func{\Ext{\Xt}}{\asubset}\cap\func{\Fr{\Xt}}{\asubset}&=
\func{\Int{\Xt}}{\compl{\X}{\asubset}}\cap\func{\Fr{\Xt}}{\compl{\X}{\asubset}}\cr
&=\func{\Int{\Xt}}{\compl{\X}{\asubset}}\cap
\(\compl{\func{\Cl{\Xt}}{\compl{\X}{\asubset}}}{\func{\Int{\Xt}}{\compl{\X}{\asubset}}}\)\cr
&=\empty.
\end{align}
According to
\refthm{thmclosurefrontierrelation}
and
\refthm{thmclosureandextrelation},
\begin{align}
\[\func{\Int{\Xt}}{\asubset}\cup\func{\Fr{\Xt}}{\asubset}\]\cup\func{\Ext{\Xt}}{\asubset}&=
\func{\Cl{\Xt}}{\asubset}\cup\[\compl{\X}{\func{\Cl{\Xt}}{\asubset}}\]\cr
&=\X.
\end{align}
\endthm
%%%%%%%%%%%%%%%%%%%%%%%%%%%%%%%%%%%%%%%%%%%%%%%%%%%%%%%%%%%%%%%%%%%%%%%%%%%%%%%%%%%%%%%%%%%%%
\theorem\label{thmfrontierofsetisclosed}
$\Xt=\opair{\X}{\topology{}}$
is taken as a topological-space.
The boundary of every subset of $\X$
in the topological-space $\Xt$
is a closed set of $\Xt$. That is,
\begin{equation}
\Foreach{\asubset}{\CSs{\X}}
\[\func{\Fr{\Xt}}{\asubset}\in\Fclosed{\X}{\asubset}\].
\end{equation}
\prooff
$\asubset$
is taken as an arbitrary subset of $\X$.
According to
\refcor{corclosureofset0},
\begin{align}
\func{\Cl{\Xt}}{\asubset}\in\Fclosed{\X}{\topology{}},\\
\func{\Cl{\Xt}}{\compl{\X}{\asubset}}\in\Fclosed{\X}{\topology{}}.
\end{align}
Thus according to
\refthm{thmfrontier1}
and
\refthm{thmclosedsets},
\begin{align}
\func{\Fr{\Xt}}{\asubset}&=\func{\Cl{\Xt}}{\asubset}\cap
\func{\Cl{\Xt}}{\compl{\X}{\asubset}}\cr
&\in\Fclosed{\X}{\topology{}}.
\end{align}
\endthm
%%%%%%%%%%%%%%%%%%%%%%%%%%%%%%%%%%%%%%%%%%%%%%%%%%%%%%%%%%%%%%%%%%%%%%%%%%%%%%%%%
\theorem\label{thmfrontierofopenset}
$\Xt=\opair{\X}{\topology{}}$
is taken as a topological-space.
For every
$\asubset$
in
$\CSs{\X}$,
$\asubset$
is an open set of $\Xt$ if-and-only if
$\func{\Fr{\Xt}}{\asubset}$
and
$\asubset$
do not intersect each other. That is,
\begin{equation}
\Foreach{\asubset}{\CSs{\X}}
\[\asubset\in\topology{}\thenn
\func{\Fr{\Xt}}{\asubset}\cap\asubset=\empty\].
\end{equation}
\prooff
$\asubset$
is taken as an arbitrary subset of $\X$.
According to \refdef{deffamilyofclosedsets},
\begin{equation}
\(\asubset\in\topology{}\)\thenn\[\(\compl{\X}{\asubset}\)\in\Fclosed{\X}{\topology{}}\].
\end{equation}
According to
\refthm{thmfrontierofclosedset},
\begin{equation}
\[\(\compl{\X}{\asubset}\)\in\Fclosed{\X}{\topology{}}\]\thenn
\[\func{\Fr{\Xt}}{\compl{\X}{\asubset}}\subseteq\(\compl{\X}{\asubset}\)\].
\end{equation}
According to
\refthm{thmfrontierofcomplementofset},
\begin{equation}
\[\func{\Fr{\Xt}}{\compl{\X}{\asubset}}\subseteq\(\compl{\X}{\asubset}\)\]\thenn
\[\func{\Fr{\Xt}}{\asubset}\subseteq\(\compl{\X}{\asubset}\)\].
\end{equation}
It is clear that,
\begin{equation}
\[\func{\Fr{\Xt}}{\asubset}\subseteq\(\compl{\X}{\asubset}\)\]\thenn
\[\func{\Fr{\Xt}}{\asubset}\cap\asubset=\empty\].
\end{equation}
Therefore,
\begin{equation}
\(\asubset\in\topology{}\)\thenn
\[\func{\Fr{\Xt}}{\asubset}\cap\asubset=\empty\].
\end{equation}
\endthm
%%%%%%%%%%%%%%%%%%%%%%%%%%%%%%%%%%%%%%%%%%%%%%%%%%%%%%%%%%%%%%%%%%%%%%%%%%%%%%%%%%
\theorem\label{thmfrontierofsetwithnullint}
$\Xt=\opair{\X}{\topology{}}$
is taken as a topological-space.
For every
$\asubset$
in
$\CSs{\X}$,
the boundary of
$\asubset$
in the topological-space $\Xt$
equals
$\asubset$, if-and-only-if
$\asubset$
is a closed set of $\Xt$ and the interior of
$\asubset$
in
$\Xt$
is empty. That is,
\begin{equation}
\Foreach{\asubset}{\CSs{\X}}
\left\{\[\asubset\in\Fclosed{\X}{\topology{}},~
\func{\Int{\Xt}}{\asubset}=\empty\]\thenn
\func{\Fr{\Xt}}{\asubset}=\asubset\right\}.
\end{equation}
\prooff
$\asubset$
is taken as an arbitrary subset of
$\X$.
\begin{itemize}
\item[${\textbf{\textsf{p1}}}$]
It is assumed that,
\begin{align}
&\asubset\in\Fclosed{\X}{\topology{}},\label{thmfrontierofsetwithnullintp1eq1}\\
&\func{\Int{\Xt}}{\asubset}=\empty.\label{thmfrontierofsetwithnullintp1eq2}
\end{align}
According to
\refthm{thmclosureofclosedset},
\Ref{thmfrontierofsetwithnullintp1eq1}
implies,
\begin{equation}
\func{\Cl{\Xt}}{\asubset}=\asubset.
\end{equation}
Thus according to
\refdef{deffrontierofset}
and
\Ref{thmfrontierofsetwithnullintp1eq2},
\begin{align}
\func{\Fr{\Xt}}{\asubset}&=
\compl{\func{\Cl{\Xt}}{\asubset}}{\func{\Int{\Xt}}{\asubset}}\cr
&=\compl{\asubset}{\empty}\cr
&=\asubset.
\end{align}
\endp
\item[${\textbf{\textsf{p2}}}$]
It is assumed that,
\begin{equation}
\func{\Fr{\Xt}}{\asubset}=\asubset.
\end{equation}
Then, according to
\refthm{thmfrontierofsetisclosed},
$\asubset$
is a closed set of $\Xt$. That is,
$\asubset\in\Fclosed{\X}{\topology{}}$,
and hence according to
\refthm{thmclosureofclosedset}
and
\refdef{deffrontierofset},
\begin{align}
\asubset&=\func{\Fr{\Xt}}{\asubset}\cr
&=\compl{\func{\Cl{\Xt}}{\asubset}}{\func{\Int{\Xt}}{\asubset}}\cr
&=\compl{\asubset}{\func{\Int{\Xt}}{\asubset}}.
\end{align}
Thus considering that
$\func{\Int{\Xt}}{\asubset}\subseteq\asubset$,
clearly,
\begin{equation}
\func{\Int{\Xt}}{\asubset}=\empty.
\end{equation}
\endp
\end{itemize}
\endthm
%%%%%%%%%%%%%%%%%%%%%%%%%%%%%%%%%%%%%%%%%%%%%%%%%%%%%%%%%%%%%%%%%%%%%%%%%%%%%%%%%%
\theorem\label{thmfrontierofclopenset}
$\Xt=\opair{\X}{\topology{}}$
is taken as a topological-space.
For every
$\asubset$
in
$\CSs{\X}$,
the boundary of
$\asubset$
in the topological-space $\Xt$
is empty, if-and-only-if
$\asubset$
is a clopen set of $\Xt$. That is,
\begin{equation}
\Foreach{\asubset}{\CSs{\X}}
\[\asubset\in\topology{}\cap\Fclosed{\X}{\topology{}}\thenn
\func{\Fr{\Xt}}{\asubset}=\empty\].
\end{equation}
\prooff
$\asubset$
is taken as an arbitrary subset of $\X$.
According to
\refthm{thmintofopenset}
and
\refthm{thmclosureofclosedset},
\begin{equation}
\asubset\in\[\topology{}\cap\Fclosed{\X}{\topology{}}\]\thenn
\func{\Int{\Xt}}{\asubset}=\func{\Cl{\Xt}}{\asubset}=\asubset.
\end{equation}
In addition, considering that,
$\func{\Int{\Xt}}{\asubset}\subseteq\func{\Cl{\Xt}}{\asubset}$,
\begin{equation}
\func{\Cl{\Xt}}{\asubset}=\func{\Int{\Xt}}{\asubset}\thenn
\compl{\func{\Cl{\Xt}}{\asubset}}{\func{\Int{\Xt}}{\asubset}}=\empty.
\end{equation}
Thus, according to
\refdef{deffrontierofset},
\begin{equation}
\func{\Cl{\Xt}}{\asubset}=\func{\Int{\Xt}}{\asubset}\thenn\func{\Fr{\Xt}}{\asubset}=\empty.
\end{equation}
Therefore,
\begin{equation}
\asubset\in\topology{}\cap\Fclosed{\X}{\topology{}}\thenn
\func{\Fr{\Xt}}{\asubset}=\empty.
\end{equation}
\endthm
%%%%%%%%%%%%%%%%%%%%%%%%%%%%%%%%%%%%%%%%%%%%%%%%%%%%%%%%%%%%%%%%%%%%%%%%%%%%%%%%%%%%%%%%%%%%%%%%%%%%%%%%
\theorem\label{thmfrontierofunionofsets}
$\Xt=\opair{\X}{\topology{}}$
is taken as a topological-space.
\begin{equation}
\Foreach{\opair{\asubset}{\p{\asubset}}}{\CSs{\X}\times\CSs{\X}}
\[\func{\Fr{\Xt}}{\asubset\cup\p{\asubset}}\subseteq
\func{\Fr{\Xt}}{\asubset}\cup\func{\Fr{\Xt}}{\p{\asubset}}\].
\end{equation}
\prooff
$\asubset$
and
$\p{\asubset}$
are taken as arbitrary subsets of $\X$.\\
According to
\refthm{thmintofunionofsets}
($\func{\Int{\Xt}}{\asubset\cup\p{\asubset}}\supseteq\func{\Int{\Xt}}{\asubset}\cup\func{\Int{\Xt}}{\p{\asubset}}$),
it is clear that
\begin{equation}\label{thmfrontierofunionofsetsp1}
\compl{\func{\Cl{\Xt}}{\asubset\cup\p{\asubset}}}
{\func{\Int{\Xt}}{\asubset\cup\p{\asubset}}}\subseteq
\compl{\func{\Cl{\Xt}}{\asubset\cup\p{\asubset}}}
{\[\func{\Int{\Xt}}{\asubset}\cup\func{\Int{\Xt}}{\p{\asubset}}\]}.
\end{equation}
Thus according to
\refthm{thmclosureofunionofsets},
($\func{\Cl{\Xt}}{\asubset\cup\p{\asubset}}=\func{\Cl{\Xt}}{\asubset}\cup\func{\Cl{\Xt}}{\p{\asubset}}$),
\begin{equation}
\compl{\func{\Cl{\Xt}}{\asubset\cup\p{\asubset}}}
{\func{\Int{\Xt}}{\asubset\cup\p{\asubset}}}\subseteq
\compl{\[\func{\Cl{\Xt}}{\asubset}\cup\func{\Cl{\Xt}}{\p{\asubset}}\]}
{\[\func{\Int{\Xt}}{\asubset}\cup\func{\Int{\Xt}}{\p{\asubset}}\]}.
\end{equation}
In addition, considering that,
\begin{align}
\left\{
\begin{aligned}
&\func{\Int{\Xt}}{\asubset}\subseteq\func{\Cl{\Xt}}{\asubset},\\
&\func{\Int{\Xt}}{\p{\asubset}}\subseteq\func{\Cl{\Xt}}{\p{\asubset}},
\end{aligned}
\right.
\end{align}
it is clear that,
\begin{gather}
\compl{\[\func{\Cl{\Xt}}{\asubset}\cup\func{\Cl{\Xt}}{\p{\asubset}}\]}
{\[\func{\Int{\Xt}}{\asubset}\cup\func{\Int{\Xt}}{\p{\asubset}}\]}\cr
\subseteq\cr
\left\{\[\compl{\func{\Cl{\Xt}}{\asubset}}{\func{\Int{\Xt}}{\asubset}}\]\cup
\[\compl{\func{\Cl{\Xt}}{\p{\asubset}}}{\func{\Int{\Xt}}{\p{\asubset}}}\]\right\}.
\end{gather}
Hence according to
\Ref{thmfrontierofunionofsetsp1}
and
\refdef{deffrontierofset},
\begin{equation}
\func{\Fr{\Xt}}{\asubset\cup\p{\asubset}}\subseteq
\[\func{\Fr{\Xt}}{\asubset}\cup\func{\Fr{\Xt}}{\p{\asubset}}\].
\end{equation}
\endthm
%%%%%%%%%%%%%%%%%%%%%%%%%%%%%%%%%%%%%%%%%%%%%%%%%%%%%%%%%%%%%%%%%%%%%%%%%%%%%%%%%%%%%%%%
\theorem\label{thmfrontieroffrontierofset}
$\Xt=\opair{\X}{\topology{}}$
is taken as a topological-space.
\begin{equation}
\Foreach{\asubset}{\CSs{\X}}
\[\func{\Fr{\Xt}}{\func{\Fr{\Xt}}{\asubset}}\subseteq\func{\Fr{\Xt}}{\asubset}\].
\end{equation}
\prooff
$\asubset$
is taken as an arbitrary subset of $\X$.
According to \refthm{thmfrontierofsetisclosed},
$\func{\Fr{\Xt}}{\asubset}$
is a closed set of $\Xt$
($\func{\Fr{\Xt}}{\asubset}\in\Fclosed{\X}{\topology{}}$).
Thus according to
\refthm{thmfrontierofclosedset},
$\func{\Fr{\Xt}}{\func{\Fr{\Xt}}{\asubset}}\subseteq\func{\Fr{\Xt}}{\asubset}$.
\endthm
%%%%%%%%%%%%%%%%%%%%%%%%%%%%%%%%%%%%%%%%%%%%%%%%%%%%%%%%%%%%%%%%%%%%%%%%%%%%%%%%%%
\theorem\label{thmintoffrontierofclosedset}
$\Xt=\opair{\X}{\topology{}}$
is taken as a topological-space.
The interior of boundary of every closed set of $\Xt$ is empty. That is,
\begin{equation}
\Foreach{\asubset}{\Fclosed{\X}{\topology{}}}
\[\func{\Int{\Xt}}{\func{\Fr{\Xt}}{\asubset}}=\empty\].
\end{equation}
\prooff
$\asubset$
is taken as a closed set of $\Xt$.
Then according to \refthm{thmfrontierofclosedset},
\begin{equation}
\func{\Fr{\Xt}}{\asubset}\subseteq\asubset,
\end{equation}
and hence according to
\refdef{thmintofasubsetofset},
\begin{equation}
\func{\Int{\Xt}}{\func{\Fr{\Xt}}{\asubset}}\subseteq\func{\Int{\Xt}}{\asubset}.
\end{equation}
In addition, according to \refcor{corintofset0},
\begin{equation}
\func{\Int{\Xt}}{\func{\Fr{\Xt}}{\asubset}}\subseteq\func{\Fr{\Xt}}{\asubset}.
\end{equation}
Therefore, according to \refthm{thmpartitionofuniversalset},
\begin{align}
\func{\Int{\Xt}}{\func{\Fr{\Xt}}{\asubset}}&\subseteq\[\func{\Fr{\Xt}}{\asubset}
\cap\func{\Int{\Xt}}{\asubset}\]\cr
&=\empty
\end{align}
and hence,
\begin{equation}
\func{\Int{\Xt}}{\func{\Fr{\Xt}}{\asubset}}=\empty.
\end{equation}
\endthm
%%%%%%%%%%%%%%%%%%%%%%%%%%%%%%%%%%%%%%%%%%%%%%%%%%%%%%%%%%%%%%%%%%%%%%%%%%%%%%%%%%%%%
\theorem\label{thmintoffrontierofopenset}
$\Xt=\opair{\X}{\topology{}}$
is taken as a topological-space.
The interior of boundary of every open set of $\Xt$ is empty. That is,
\begin{equation}
\Foreach{\asubset}{\topology{}}
\[\func{\Int{\Xt}}{\func{\Fr{\Xt}}{\asubset}}=\empty\].
\end{equation}
\prooff
$\asubset$
is taken as an arbitrary open set of $\Xt$.
Then,
$\(\compl{\X}{\asubset}\)\in\Fclosed{\X}{\topology{}}$,
and hence according to
\refthm{thmintoffrontierofclosedset},
\begin{equation}
\func{\Int{\Xt}}{\func{\Fr{\Xt}}{\compl{\X}{\asubset}}}=\empty.
\end{equation}
Thus according to
\refthm{thmfrontierofcomplementofset}
($\func{\Fr{\Xt}}{\compl{\X}{\asubset}}=\func{\Fr{\Xt}}{\asubset}$),
$\func{\Int{\Xt}}{\func{\Fr{\Xt}}{\asubset}}=\empty$.
\endthm
%%%%%%%%%%%%%%%%%%%%%%%%%%%%%%%%%%%%%%%%%%%%%%%%%%%%%%%%%%%%%%%%%%%%%%%%%%%%%%%%%%%%%
\theorem\label{thmintoffrontierofclosureofset}
$\Xt=\opair{\X}{\topology{}}$
is taken as a topological-space.
The interior of boundary of closure of every subset of $\X$ is empty. That is,
\begin{equation}
\Foreach{\asubset}{\CSs{\X}}
\[\func{\Int{\Xt}}{\func{\Fr{\Xt}}{\func{\Cl{\Xt}}{\asubset}}}=\empty\].
\end{equation}
\prooff
According to
\refcor{corclosureofset0}
and
\refthm{thmintoffrontierofclosedset}
it is clear.
\endthm
%%%%%%%%%%%%%%%%%%%%%%%%%%%%%%%%%%%%%%%%%%%%%%%%%%%%%%%%%%%%%%%%%%%%%%%%%%%%%%%%%%%%%
\theorem\label{thmintoffrontierofinteriorofset}
$\Xt=\opair{\X}{\topology{}}$
is taken as a topological-space.
The interior of boundary of interior of every subset of $\X$ is empty. That is,
\begin{equation}
\Foreach{\asubset}{\CSs{\X}}
\[\func{\Int{\Xt}}{\func{\Fr{\Xt}}{\func{\Int{\Xt}}{\asubset}}}=\empty\].
\end{equation}
\prooff
According to \refcor{corintofset0} and \refthm{thmintoffrontierofopenset},
it is clear.
\endthm
%%%%%%%%%%%%%%%%%%%%%%%%%%%%%%%%%%%%%%%%%%%%%%%%%%%%%%%%%%%%%%%%%%%%%%%%%%%%%%%%%%%%%%
\theorem\label{thmpowersoffrontierofset}
$\Xt=\opair{\X}{\topology{}}$
is taken as a topological-space.
\begin{equation}\label{thmpowersoffrontierofseteq1}
\Foreach{n}{\(\compl{\Zp}{\seta{1}}\)}
\power{\Fr{\Xt}}{n}=\power{\Fr{\Xt}}{2}.
\end{equation}
\prooff
\begin{itemize}
\item[${\textbf{\textsf{p}}}$]
$n$
is taken as an arbitrary element of
$\(\compl{\Zp}{\seta{1}}\)$.
\begin{itemize}
\item[${\textbf{\textsf{p1}}}$]
$\asubset$
is taken as an arbitrary subset of $\X$.\\
According to
\refthm{thmfrontierofsetisclosed},
each
$\pfunc{\Fr{\Xt}}{\asubset}{n-1}$
and
$\pfunc{\Fr{\Xt}}{\func{\Fr{\Xt}}{\asubset}}{n}$
is a closed set of $\Xt$. That is,
\begin{align}
\pfunc{\Fr{\Xt}}{\asubset}{n-1}&\in\Fclosed{\X}{\topology{}},\label{thmpowersoffrontierofsetp1}\\
\pfunc{\Fr{\Xt}}{\asubset}{n}&\in\Fclosed{\X}{\topology{}}.\label{thmpowersoffrontierofsetp2}
\end{align}
\Ref{thmpowersoffrontierofsetp1}
and
\refthm{thmintoffrontierofclosedset}
imply,
\begin{align}\label{thmpowersoffrontierofsetp3}
\func{\Int{\Xt}}{\pfunc{\Fr{\Xt}}{\asubset}{n}}&=
\func{\Int{\Xt}}{\func{\Fr{\Xt}}{\pfunc{\Fr{\Xt}}{\asubset}{n-1}}}\cr
&=\empty.
\end{align}
\Ref{thmpowersoffrontierofsetp2},
\Ref{thmpowersoffrontierofsetp3},
and
\refthm{thmfrontierofsetwithnullint}
imply,
\begin{align}
\pfunc{\Fr{\Xt}}{\asubset}{n+1}&=
\func{\Fr{\Xt}}{\pfunc{\Fr{\Xt}}{\asubset}{n}}\cr
&=\pfunc{\Fr{\Xt}}{\asubset}{n}.
\end{align}
\endp
\end{itemize}
Therefore,
\begin{equation}
\Foreach{\asubset}{\CSs{\X}}
\pfunc{\Fr{\Xt}}{\asubset}{n+1}=
\pfunc{\Fr{\Xt}}{\asubset}{n}.
\end{equation}
which means,
\begin{equation}
\power{\Fr{\Xt}}{n+1}=\power{\Fr{\Xt}}{n}.
\end{equation}
\endp
\end{itemize}
Therefore,
\begin{equation}
\Foreach{n}{\(\compl{\Zp}{\seta{1}}\)}
\power{\Fr{\Xt}}{n+1}=\power{\Fr{\Xt}}{n}.
\end{equation}
It is equivalent to \Ref{thmpowersoffrontierofseteq1}.
\endthm
%%%%%%%%%%%%%%%%%%%%%%%%%%%%%%%%%%%%%%%%%%%%%%%%%%%%%%%%%%%%%%%%%%%%%%%%%%%%%%%%%%%%%%%%%%
\theorem\label{thmsubspacefrontier}
$\opair{\X}{\topology{}}$
is taken as a topological-space, and
$\Y$
as a subset of
$\X$.
\begin{equation}
\Foreach{\asubset}{\CSs{\Y}}
\[\func{\Fr{\opair{\Y}{\stopology{\topology{}}{\Y}}}}{\asubset}\subseteq
\func{\Fr{\Xt}}{\asubset}\].
\end{equation}
\prooff
$\asubset$
is taken as an arbitrary subset of $\Y$.
(It is obvious that
$\asubset$
is also a subset of
$\X$.)
According to
\refthm{thmsubspaceinterior},
\begin{equation}\label{thmsubspacefrontierp1}
\func{\Int{\Xt}}{\asubset}\subseteq
\func{\Int{\opair{\Y}{\stopology{\topology{}}{\Y}}}}{\asubset},
\end{equation}
and hence according to
\begin{align}\label{thmsubspacefrontierp2},
\[\compl{\func{\Cl{\opair{\Y}{\stopology{\topology{}}{\Y}}}}{\asubset}}
{\func{\Int{\opair{\Y}{\stopology{\topology{}}{\Y}}}}{\asubset}}\]\subseteq
\[\compl{\func{\Cl{\opair{\Y}{\stopology{\topology{}}{\Y}}}}{\asubset}}
{\func{\Int{\Xt}}{\asubset}}\].
\end{align}
In addition, according to
\refthm{thmsubspaceclosure},
\begin{equation}\label{thmsubspacefrontierp3}
\func{\Cl{\opair{\Y}{\stopology{\topology{}}{\Y}}}}{\asubset}
\subseteq\func{\Cl{\Xt}}{\asubset},
\end{equation}
and hence,
\begin{equation}\label{thmsubspacefrontierp4}
\[\compl{\func{\Cl{\opair{\Y}{\stopology{\topology{}}{\Y}}}}{\asubset}}
{\func{\Int{\Xt}}{\asubset}}\]\subseteq
\[\compl{\func{\Cl{\Xt}}{\asubset}}
{\func{\Int{\Xt}}{\asubset}}\].
\end{equation}
According to
\Ref{thmsubspacefrontierp2},
\Ref{thmsubspacefrontierp4},
and
\refdef{deffrontierofset},
\begin{align}
\func{\Fr{\opair{\Y}{\stopology{\topology{}}{\Y}}}}{\asubset}&=
\[\compl{\func{\Cl{\opair{\Y}{\stopology{\topology{}}{\Y}}}}{\asubset}}
{\func{\Int{\opair{\Y}{\stopology{\topology{}}{\Y}}}}{\asubset}}\]\cr
&\subseteq\[\compl{\func{\Cl{\Xt}}{\asubset}}
{\func{\Int{\Xt}}{\asubset}}\]\cr
&=\func{\Fr{\Xt}}{\asubset}.
\end{align}
\endthm
%%%%%%%%%%%%%%%%%%%%%%%%%%%%%%%%%%%%%%%%%%%%%%%%%%%%%%%%%%%%%%%%%%%%%%%%%%%%%%%%%%%%%%%%%%%%
\theorem\label{thmfrontierofsetwithdifferenttopologies}
Each
$\Xt_{1}=\opair{\X}{\topology{1}}$
and
$\Xt_{2}=\opair{\X}{\topology{2}}$
is taken as a topological-space.
\begin{equation}
\(\topology{2}\subseteq\topology{1}\)\then\[\Foreach{\asubset}{\CSs{\X}}
\func{\Fr{\Xt_{2}}}{\asubset}\supseteq\func{\Fr{\Xt_{1}}}{\asubset}\].
\end{equation}
\prooff
It is assumed that,
\begin{equation}\label{thmfrontierofsetwithdifferenttopologiesp1}
\topology{2}\subseteq\topology{1},
\end{equation}
and
$\asubset$
is taken as an arbitrary subset of $\X$.
According to \refthm{thmintofsetwithdifferenttopologies},
\begin{equation}\label{thmfrontierofsetwithdifferenttopologiesp2}
\func{\Int{\Xt_{2}}}{\asubset}\subseteq\func{\Int{\Xt_{1}}}{\asubset}.
\end{equation}
Thus,
\begin{equation}\label{thmfrontierofsetwithdifferenttopologiesp3}
\[\compl{\func{\Cl{\Xt_1}}{\asubset}}{\func{\Int{\Xt_1}}{\asubset}}\]\subseteq
\[\compl{\func{\Cl{\Xt_1}}{\asubset}}{\func{\Int{\Xt_2}}{\asubset}}\].
\end{equation}
In addition, according to \refthm{thmclosureofsetwithdifferenttopologies},
\begin{equation}\label{thmfrontierofsetwithdifferenttopologiesp4}
\func{\Cl{\Xt_{1}}}{\asubset}\subseteq\func{\Cl{\Xt_{2}}}{\asubset}.
\end{equation}
Thus,
\begin{equation}\label{thmfrontierofsetwithdifferenttopologiesp5}
\[\compl{\func{\Cl{\Xt_1}}{\asubset}}{\func{\Int{\Xt_2}}{\asubset}}\]\subseteq
\[\compl{\func{\Cl{\Xt_2}}{\asubset}}{\func{\Int{\Xt_2}}{\asubset}}\].
\end{equation}
\Ref{thmfrontierofsetwithdifferenttopologiesp3},
\Ref{thmfrontierofsetwithdifferenttopologiesp5},
and
\refdef{deffrontierofset} imply,
\begin{align}\label{thmfrontierofsetwithdifferenttopologiesp6}
\func{\Fr{\Xt_{1}}}{\asubset}&=
\[\compl{\func{\Cl{\Xt_1}}{\asubset}}{\func{\Int{\Xt_1}}{\asubset}}\]\cr
&\subseteq\[\compl{\func{\Cl{\Xt_2}}{\asubset}}{\func{\Int{\Xt_2}}{\asubset}}\]\cr
&=\func{\Fr{\Xt_{2}}}{\asubset}.
\end{align}
\endthm
%%%%%%%%%%%%%%%%%%%%%%%%%%%%%%%%%%%%%%%%%%%%%%%%%%%%%%%%%%%%%%%%%%%%%%%%%%%%%%%%%%%%%%%%%%
%%%%%%%%%%%%%%%%%%%%%%%%%%%%%%%%%%%%%%%%%%%%%%%%%%%%%%%%%%%%%%%%%%%%%%%%%%%%%%%%%%%%%%%%%%
%%%%%%%%%%%%%%%%%%%%%%%%%%%%%%%%%%%%%%%%%%%%%%%%%%%%%%%%%%%%%%%%%%%%%%%%%%%%%%%%%%%%%%%%%%
%%%%%%%%%%%%%%%%%%%%%%%%%%%%%%%%%%%%%%%%%%%%%%%%%%%%%%%%%%%%%%%%%%%%%%%%%%%%%%%%%%%%%%%%%%
%%%%%%%%%%%%%%%%%%%%%%%%%%%%%%%%%%%%%%%%%%%%%%%%%%%%%%%%%%%%%%%%%%%%%%%%%%%%%%%%%%%%%%%%%%
%%%%%%%%%%%%%%%%%%%%%%%%%%%%%%%%%%%%%%%%%%%%%%%%%%%%%%%%%%%%%%%%%%%%%%%%%%%%%%%%%%%%%%%%%%
%%%%%%%%%%%%%%%%%%%%%%%%%%%%%%%%%%%%%%%%%%%%%%%%%%%%%%%%%%%%%%%%%%%%%%%%%%%%%%%%%%%%%%%%%%
%%%%%%%%%%%%%%%%%%%%%%%%%%%%%%%%%%%%%%%%%%%%%%%%%%%%%%%%%%%%%%%%%%%%%%%%%%%%%%%%%%%%%%%%%%
%%%%%%%%%%%%%%%%%%%%%%%%%%%%%%%%%%%%%%%%%%%%%%%%%%%%%%%%%%%%%%%%%%%%%%%%%%%%%%%%%%%%%%%%%%
%%%%%%%%%%%%%%%%%%%%%%%%%%%%%%%%%%%%%%%%%%%%%%%%%%%%%%%%%%%%%%%%%%%%%%%%%%%%%%%%%%%%%%%%%%
%%%%%%%%%%%%%%%%%%%%%%%%%%%%%%%%%%%%%%%%%%%%%%%%%%%%%%%%%%%%%%%%%%%%%%%%%%%%%%%%%%%%%%%%%%
\subsection{
Dense Sets of a Topological Space}
\definition\label{defdenseset}
$\opair{\X}{\topology{}}$
is taken as a topological-space, and
$\asubset$
as a subset of
$\X$.
It is said that
$\quotl$$\asubset$ is dense in the topological-space $\Xt$$\quotr$,
and
$\asubset$
is referred to as a $\quotl$dense set of $\Xt$$\quotr$ iff
\begin{equation}
\func{\Cl{\Xt}}{\asubset}=\X.
\end{equation}
\endef
%%%%%%%%%%%%%%%%%%%%%%%%%%%%%%%%%%%%%%%%%%%%%%%%%%%%%%%%%%%%%%%%%%%%%%%%%%%%%%%%%%%%%%%%%%%%%
\definition\label{defclassofdensesets}
$\opair{\X}{\topology{}}$
is taken as a topological-space.
The set of all dense sets of $\Xt$
is denoted by $\Fdense{\X}{\topology{}}$. That is,
\begin{equation}
\Fdense{\X}{\topology{}}:=\defset{\asubset}{\CSs{\X}}
{\[\func{\Cl{\Xt}}{\asubset}=\X\]}.
\end{equation}
\endef
%%%%%%%%%%%%%%%%%%%%%%%%%%%%%%%%%%%%%%%%%%%%%%%%%%%%%%%%%%%%%%%%%%%%%%%%%%%%%%%%%%%%%%%%%%%%%
\theorem\label{thmifintofsetisdensethensetisdense}
$\opair{\X}{\topology{}}$
is taken as a topological-space.
\begin{equation}
\Foreach{\asubset}{\CSs{\X}}
\[\func{\Int{\Xt}}{\asubset}\in\Fdense{\X}{\topology{}}\then
\asubset\in\Fdense{\X}{\topology{}}\].
\end{equation}
\prooff
$\asubset$
is taken as an arbitrary subset of
$\X$,
and it is assumed that,
\begin{equation}
\func{\Int{\Xt}}{\asubset}\in\Fdense{\X}{\topology{}}.
\end{equation}
Then according to
\refdef{defclassofdensesets},
\begin{equation}
\func{\Cl{\Xt}}{\func{\Int{\Xt}}{\asubset}}=\X.
\end{equation}
Considering that
$\func{\Int{\Xt}}{\asubset}\subseteq\asubset$,
according to
\refthm{thmclosureofasubsetofset},
\begin{equation}
\func{\Cl{\Xt}}{\asubset}\supseteq
\func{\Cl{\Xt}}{\func{\Int{\Xt}}{\asubset}}.
\end{equation}
Therefore,
\begin{equation}
\func{\Cl{\Xt}}{\asubset}=\X,
\end{equation}
and hence according to
\refdef{defclassofdensesets},
\begin{equation}
\asubset\in\Fdense{\X}{\topology{}}.
\end{equation}
\endthm
%%%%%%%%%%%%%%%%%%%%%%%%%%%%%%%%%%%%%%%%%%%%%%%%%%%%%%%%%%%%%%%%%%%%%%%%%%%%%%%%%%%%%%%%%%%%%
\theorem
$\opair{\X}{\topology{}}$
is taken as a topological-space.
$\X$
is dense in
$\Xt$. That is,
\begin{equation}
\X\in\Fdense{\X}{\topology{}}.
\end{equation}
\prooff
Considering that
$\func{\Cl{\Xt}}{\X}=\X$,
and according to
\refdef{defdenseset},
it is clear.
\endthm
%%%%%%%%%%%%%%%%%%%%%%%%%%%%%%%%%%%%%%%%%%%%%%%%%%%%%%%%%%%%%%%%%%%%%%%%%%%%%%%%%%%%%%%%%%%%%%
\theorem\label{thmdenseNandSconditions0}
$\Xt=\opair{\X}{\topology{}}$
is taken as a topological-space.
For every
$\asubset$
in
$\CSs{\X}$,
$\asubset$
is a dense set of $\Xt$ if-and-only-if the exterior of $\asubset$
in the topological-space $\Xt$ is empty. That is,
\begin{equation}
\Foreach{\asubset}{\CSs{\X}}
\left\{\[\asubset\in\Fdense{\X}{\topology{}}\]\thenn
\[\func{\Ext{\Xt}}{\asubset}=\empty\]\right\}.
\end{equation}
\prooff
$\asubset$
is taken as an arbitrary subset of $\X$.
According to
\refdef{defclassofdensesets},
\begin{equation}
\[\asubset\in\Fdense{\X}{\topology{}}\]\thenn\[\func{\Cl{\Xt}}{\asubset}=\X\].
\end{equation}
In addition, according to
\refthm{thmclosureandextrelation},
\begin{equation}
\[\func{\Cl{\Xt}}{\asubset}=\X\]\thenn\[\func{\Ext{\Xt}}{\asubset}=\empty\].
\end{equation}
Therefore,
\begin{equation}
\[\asubset\in\Fdense{\X}{\topology{}}\]\thenn
\[\func{\Ext{\Xt}}{\asubset}=\empty\].
\end{equation}
\endthm
%%%%%%%%%%%%%%%%%%%%%%%%%%%%%%%%%%%%%%%%%%%%%%%%%%%%%%%%%%%%%%%%%%%%%%%%%%%%%%%%%%%%%%%%%%%%%%
\theorem\label{thmdenseNandSconditions}
$\Xt=\opair{\X}{\topology{}}$
is taken as a topological-space.
For every
$\asubset$
in
$\CSs{\X}$,
$\asubset$
is a dense set of $\Xt$ if-and-only-if
$\asubset$
intersects every non-empty open set of $\Xt$.That is,
\begin{equation}
\Foreach{\asubset}{\CSs{\X}}
\left\{\[\asubset\in\Fdense{\X}{\topology{}}\]
\thenn\[\Foreach{\U}{\(\compl{\topology{}}{\seta{\empty}}\)}\U\cap\asubset\neq\empty\]\right\}.
\end{equation}
\prooff
$\asubset$
is taken as an arbitrary subset of
$\X$
\begin{itemize}
\item[${\textbf{\textsf{p1}}}$]
It is assumed that
$\asubset$
is a dense set of $\Xt$
($\asubset\in\Fdense{\X}{\topology{}}$).
Then according to
\refdef{defclassofdensesets},
\begin{equation}\label{thmdenseNandSconditionsp1eq1}
\func{\Cl{\Xt}}{\asubset}=\X.
\end{equation}
\begin{itemize}
\item[${\textbf{\textsf{p1-1}}}$]
$\U$
is taken as an arbitrary element of
$\compl{\topology{}}{\seta{\empty}}$.
Then according to
\Ref{thmdenseNandSconditionsp1eq1},
\begin{align}
\U\cap\func{\Cl{\Xt}}{\asubset}&=\U\cap\X\cr
&=\U\cr
&\neq\empty,
\end{align}
and hence according to
\refthm{thmopensetsintersectingclosure},
\begin{equation}
\U\cap\asubset\neq\empty.
\end{equation}
\endp
\item[${\textbf{\textsf{p2}}}$]
It is assumed that,
\begin{equation}
\Foreach{\U}{\(\compl{\topology{}}{\seta{\empty}}\)}
\U\cap\asubset\neq\empty.
\end{equation}
Then
$\empty$
is the only open set of
$\Xt$
and a subset of
$\X$. That is.
\begin{equation}
\topology{}\cap\CSs{\compl{\X}{\asubset}}=\seta{\empty}.
\end{equation}
Then according to
\refdef{defintofset}
and
\refdef{defextofset},
\begin{align}
\func{\Ext{\Xt}}{\asubset}&=\func{\Int{\Xt}}{\compl{\X}{\asubset}}\cr
&=\union{\topology{}\cap\CSs{\compl{\X}{\asubset}}}\cr
&=\union{\seta{\empty}}\cr
&=\empty,
\end{align}
and hence according to
\refthm{thmclosureandextrelation},
\begin{align}
\func{\Cl{\Xt}}{\asubset}&=\compl{\X}{\func{\Ext{\Xt}}{\compl{\X}{\asubset}}}\cr
&=\compl{\X}{\empty}\cr
&=\X.
\end{align}
This means,
\begin{equation}
\asubset\in\Fdense{\X}{\topology{}}.
\end{equation}
\endp
\end{itemize}
\end{itemize}
\endthm
%%%%%%%%%%%%%%%%%%%%%%%%%%%%%%%%%%%%%%%%%%%%%%%%%%%%%%%%%%%%%%%%%%%%%%%%%%%%%%%%%%%%%%%%%%%%%
\theorem\label{thmunionofdensesetsisdense}
$\Xt=\opair{\X}{\topology{}}$
is taken as a topological-space.
The union of any non-empty collection of dense-sets of
$\Xt$ is a dense set of $\Xt$. That is,
\begin{align}
\Foreach{\sCi}{\[\compl{\CSs{\Fdense{\X}{\topology{}}}}{\seta{\empty}}\]}
\(\union{\sCi}\)\in\Fdense{\X}{\topology{}}.
%\[\compl{\CSs{\CSs{\X}}}{\seta{\empty}}\]
\end{align}
\prooff
$\sCi$
is taken as an arbitrary non-empty subset of.
$\Fdense{\X}{\topology{}}$
Then according to
\refdef{defclassofdensesets},
\begin{equation}
\Foreach{\asubset}{\sCi}
\func{\Cl{\Xt}}{\asubset}=\X.
\end{equation}
Thus,
\begin{equation}
\[\bigcup_{\asubset\in\sCi}\func{\Cl{\Xt}}{\asubset}\]=\X.
\end{equation}
In addition, according to
\refthm{thmclosureofunionofsets0},
\begin{equation}
\func{\Cl{\Xt}}{\union{\sCi}}\supseteq
\[\bigcup_{\asubset\in\sCi}\func{\Cl{\Xt}}{\asubset}\].
\end{equation}
Therefore,
\begin{equation}
\func{\Cl{\Xt}}{\union{\sCi}}\supseteq\X,
\end{equation}
which means,
\begin{equation}
\func{\Cl{\Xt}}{\union{\sCi}}=\X.
\end{equation}
Hence according to
\refdef{defclassofdensesets},
\begin{equation}
\(\union{\sCi}\)\in\Fdense{\X}{\topology{}}.
\end{equation}
\endthm
%%%%%%%%%%%%%%%%%%%%%%%%%%%%%%%%%%%%%%%%%%%%%%%%%%%%%%%%%%%%%%%%%%%%%%%%%%%%%%%%%%%%%%%%%%%%%%%
%%%%%%%%%%%%%%%%%%%%%%%%%%%%%%%%%%%%%%%%%%%%%%%%%%%%%%%%%%%%%%%%%%%%%%%%%%%%%%%%%%%%%%%%%%%%%%%
%%%%%%%%%%%%%%%%%%%%%%%%%%%%%%%%%%%%%%%%%%%%%%%%%%%%%%%%%%%%%%%%%%%%%%%%%%%%%%%%%%%%%%%%%%%%%%%
\theorem\label{thmintersectionofopendensesets}
$\Xt=\opair{\X}{\topology{}}$
is taken as a topological-space.
The intersection of any pair of dense sets of $\Xt$
that at least one of them is an open set of $\Xt$,
is a dense set of $\Xt$. That is,
\begin{equation}\label{thmintersectionofopendensesetseq1}
\Foreach{\opair{\asubset}{\p{\asubset}}}
{\left\{\[\topology{}\cap\Fdense{\X}{\topology{}}\]
\times\Fdense{\X}{\topology{}}\right\}}
\(\asubset\cap\p{\asubset}\)\in\Fdense{\X}{\topology{}}.
\end{equation}
\prooff
$\asubset$
is taken as an arbitrary element of
$\topology{}\cap\Fdense{\X}{\topology{}}$,
and
$\p{\asubset}$
as an arbitrary element of
$\Fdense{\X}{\topology{}}$.
According to
\refthm{thmdenseNandSconditions},
\begin{align}
&\Foreach{\U}{\compl{\topology{}}{\seta{\empty}}}
\U\cap\asubset\neq\empty,\label{thmintersectionofopendensesetsp1}\\
&\Foreach{\U}{\compl{\topology{}}{\seta{\empty}}}
\U\cap\p{\asubset}\neq\empty.\label{thmintersectionofopendensesetsp2}
\end{align}
\begin{itemize}
\item[${\textbf{\textsf{p1}}}$]
$\U$
is taken as an arbitrary element of
$\compl{\topology{}}{\seta{\empty}}$.
According to
\Ref{thmintersectionofopendensesetsp1},
\begin{equation}
\(\U\cap\asubset\)\neq\empty.
\end{equation}
In addition, considering that
$\asubset\in\topology{}$,
and
$\U\in\topology{}$,
according to
\refdef{deftopologicalspace},
\begin{equation}
\(\U\cap\asubset\)\in\topology{}.
\end{equation}
Therefore,
\begin{equation}
\(\U\cap\asubset\)\in\compl{\topology{}}{\seta{\empty}},
\end{equation}
and hence according to
\Ref{thmintersectionofopendensesetsp2},
\begin{align}
\U\cap\(\asubset\cap\p{\asubset}\)&=
\(\U\cap\asubset\)\cap\p{\asubset}\cr
&\neq\empty.
\end{align}
\endp
\end{itemize}
Therefore,
\begin{equation}
\Foreach{\U}{\compl{\topology{}}{\seta{\empty}}}
\U\cap\(\asubset\cap\p{\asubset}\)\neq\empty,
\end{equation}
and hence according to
\refthm{thmdenseNandSconditions},
\begin{equation}
\(\asubset\cap\p{\asubset}\)\in\Fdense{\X}{\topology{}}.
\end{equation}
\endthm
%%%%%%%%%%%%%%%%%%%%%%%%%%%%%%%%%%%%%%%%%%%%%%%%%%%%%%%%%%%%%%%%%%%%%%%%%%%%%%%%%%%%%%%%%%%%%%%
%%%%%%%%%%%%%%%%%%%%%%%%%%%%%%%%%%%%%%%%%%%%%%%%%%%%%%%%%%%%%%%%%%%%%%%%%%%%%%%%%%%%%%%%%%%%%%%
\theorem
$\X$
is taken as a set.
\begin{equation}
\defset{\topology{}}{\Ctops{\X}}{\[\Fdense{\X}{\topology{}}=\seta{\X}\]}=
\seta{\CSs{\X}}.
\end{equation}
\prooff
\begin{itemize}
\item[${\textbf{\textsf{p1}}}$]
It is known that
$\CSs{\X}$
is a topology on $\X$.
\begin{equation}
\CSs{\X}\in\Ctops{\X}.
\end{equation}
In addition, according to
\refthm{thmdiscretetopologyclopensets},
every subset of $\X$
is a closed set of $\opair{\X}{\CSs{\X}}$. That is,
\begin{equation}
\Foreach{\asubset}{\CSs{\X}}
\asubset\in\Fclosed{\X}{\CSs{\X}}.
\end{equation}
Hence according to
\refthm{thmclosureofclosedset},
\begin{equation}
\Foreach{\asubset}{\CSs{\X}}
\func{\Cl{\opair{\X}{\CSs{\X}}}}{\asubset}=\asubset.
\end{equation}
Thus,
\begin{equation}
\Foreach{\asubset}{\CSs{\X}}
\left\{\[\func{\Cl{\opair{\X}{\CSs{\X}}}}{\asubset}=\X\]\then\[\asubset=\X\]\right\}.
\end{equation}
Hence according to
\refdef{defdenseset},
\begin{equation}
\Fdense{\X}{\CSs{\X}}=\seta{\X}.
\end{equation}
\endp
\item[${\textbf{\textsf{p2}}}$]
$\topology{}$
is taken as such an element of
$\Ctops{\X}$
that
\begin{equation}
\Fdense{\X}{\topology{}}=\seta{\X}.
\end{equation}
Then according to
\refdef{defclassofdensesets},
\begin{equation}
\Foreach{\asubset}{\CSs{\X}}
\left\{\[\func{\Cl{\opair{\X}{\topology{}}}}{\asubset}=\X\]\then\[\asubset=\X\]\right\}.
\end{equation}
This clearly implies,
\begin{equation}
\Foreach{\x}{\X}\func{\Cl{\opair{\X}{\topology{}}}}{\compl{\X}{\seta{\x}}}\neq\X.
\end{equation}
Hence considering that,
\begin{equation}
\Foreach{\x}{\X}
\(\compl{\X}{\seta{\x}}\)\subseteq
\func{\Cl{\opair{\X}{\topology{}}}}{\compl{\X}{\seta{\x}}}\subseteq\X,
\end{equation}
it is clear that,
\begin{equation}
\Foreach{\x}{\X}
\func{\Cl{\opair{\X}{\topology{}}}}{\compl{\X}{\seta{\x}}}=\(\compl{\X}{\seta{\x}}\).
\end{equation}
According to
\refthm{thmclosureofclosedset},
this means,
\begin{equation}
\Foreach{\x}{\X}
\(\compl{\X}{\seta{\x}}\)\in\Fclosed{\X}{\topology{}},
\end{equation}
and hence according to
\refdef{deffamilyofclosedsets},
\begin{equation}
\Foreach{\x}{\X}
\seta{\x}\in\topology{}.
\end{equation}
According to
\refthm{thmNandSdiscretetopology},
this implies,
\begin{equation}
\topology{}=\CSs{\X}.
\end{equation}
\end{itemize}
\endthm
%%%%%%%%%%%%%%%%%%%%%%%%%%%%%%%%%%%%%%%%%%%%%%%%%%%%%%%%%%%%%%%%%%%%%%%%%%%%%%%%%%%%%%%%%%%%%%%%%%%%%%%%%%%%%%%%%%%
%%%%%%%%%%%%%%%%%%%%%%%%%%%%%%%%%%%%%%%%%%%%%%%%%%%%%%%%%%%%%%%%%%%%%%%%%%%%%%%%%%%%%%%%%%%%%%%%%%%%%%%%%%%%%%%%%%%
\definition\label{defsubsetsdense}
$\Xt=\opair{\X}{\topology{}}$
is taken as a topological-space, and each
$\asubset$
and
$\p{\asubset}$
as a subset of$\X$.
It is said that, $\quotl$$\p{\asubset}$ is dense in $\asubset$ within the topological-space $\Xt$$\quotr$, iff
\begin{equation}
\func{\Cl{\Xt}}{\p{\asubset}}\supseteq\asubset.
\end{equation}
\endef
%%%%%%%%%%%%%%%%%%%%%%%%%%%%%%%%%%%%%%%%%%%%%%%%%%%%%%%%%%%%%%%%%%%%%%%%%%%%%%%%%%%%%%%%%%%%%%%
%%%%%%%%%%%%%%%%%%%%%%%%%%%%%%%%%%%%%%%%%%%%%%%%%%%%%%%%%%%%%%%%%%%%%%%%%%%%%%%%%%%%%%%%%%%%%%%
\theorem\label{thmsubspacedense}
$\Xt=\opair{\X}{\topology{}}$
is taken as a topological-space, and
$\Y$
and
$\asubset$
are taken as such subsets of $\X$ that
$\asubset\subseteq\Y$.
$\asubset$ is dense in $\Y$ within $\Xt$
if-and-only-if $\asubset$
is a dense set of
$\opair{\Y}{\stopology{\topology{}}{\Y}}$.
That is.
\begin{equation}\label{thmsubspacedenseeq1}
\[\func{\Cl{\Xt}}{\asubset}\supseteq\Y\]\thenn
\[\func{\Cl{\opair{\Y}{\stopology{\topology{}}{\Y}}}}{\asubset}=\Y\].
\end{equation}
In other words,
\begin{equation}\label{thmsubspacedenseeq1}
\[\func{\Cl{\Xt}}{\asubset}\supseteq\Y\]\thenn
\[\asubset\in\Fdense{\Y}{\stopology{\topology{}}{\Y}}\].
\end{equation}
\prooff
According to
\refthm{thmsubspaceclosure},
\begin{equation}
\func{\Cl{\opair{\Y}{\stopology{\topology{}}{\Y}}}}{\asubset}=
\Y\cap\func{\Cl{\Xt}}{\asubset}.
\end{equation}
This clearly implies \Ref{thmsubspacedenseeq1}.
\endthm
%%%%%%%%%%%%%%%%%%%%%%%%%%%%%%%%%%%%%%%%%%%%%%%%%%%%%%%%%%%%%%%%%%%%%%%%%%%%%%%%%%%%%%%%%%%%%%%%%%%
\theorem\label{thmdensetransitivity}
$\Xt=\opair{\X}{\topology{}}$
is taken as a topological-space, and each
$\asubset$,
$\bsubset$,
and
$\csubset$
as a subset of $\X$.
If $\asubset$ is dense in $\bsubset$ within $\Xt$, and
$\bsubset$
is dense in $\csubset$ within $\Xt$, then
$\asubset$ is dense in $\csubset$ within $\Xt$. That is,
\begin{equation}
\[\func{\Cl{\Xt}}{\asubset}\supseteq\bsubset,
~\func{\Cl{\Xt}}{\bsubset}\supseteq\csubset\]\then
\[\func{\Cl{\Xt}}{\asubset}\supseteq\csubset\].
\end{equation}
\prooff
It is assumed that,
\begin{align}
\func{\Cl{\Xt}}{\asubset}&\supseteq\bsubset,\\
\func{\Cl{\Xt}}{\bsubset}&\supseteq\csubset.
\end{align}
Then according to
\refthm{thmclosureofasubsetofset}
and
\refthm{thmclosureofclosureofset},
\begin{align}
\func{\Cl{\Xt}}{\bsubset}&\subseteq\func{\Cl{\Xt}}{\func{\Cl{\Xt}}{\asubset}}\cr
&=\func{\Cl{\Xt}}{\asubset},
\end{align}
and
\begin{align}
\func{\Cl{\Xt}}{\csubset}&\subseteq\func{\Cl{\Xt}}{\func{\Cl{\Xt}}{\bsubset}}\cr
&=\func{\Cl{\Xt}}{\bsubset}.
\end{align}
Therefore clearly,
\begin{equation}
\func{\Cl{\Xt}}{\csubset}\subseteq\func{\Cl{\Xt}}{\asubset}.
\end{equation}
Thus considering that
$\csubset\subseteq\func{\Cl{\Xt}}{\csubset}$,
clearly,
\begin{equation}
\func{\Cl{\Xt}}{\asubset}\supseteq\csubset.
\end{equation}
\endthm
%%%%%%%%%%%%%%%%%%%%%%%%%%%%%%%%%%%%%%%%%%%%%%%%%%%%%%%%%%%%%%%%%%%%%%%%%%%%%%%%%%%%%%%%%%%%%
\theorem\label{thmdensetransitivity1}
$\Xt=\opair{\X}{\topology{}}$
is taken as a topological-space, and each
$\Y$
and
$\asubset$
as such a subset of $\X$ that
$\asubset\subseteq\Y$.
If $\asubset$ is a dense set of $\opair{\Y}{\stopology{\topology{}}{\Y}}$, and
$\Y$
is a dense set of $\Xt$, then
$\asubset$
is a dense set of $\Xt$. That is,
\begin{equation}
\[\func{\Cl{\opair{\Y}{\stopology{\topology{}}{\Y}}}}{\asubset}=\Y,~
\func{\Cl{\Xt}}{\Y}=\X\]\then\[\func{\Cl{\Xt}}{\asubset}=\X\].
\end{equation}
\prooff
It is assumed that,
\begin{align}
\func{\Cl{\opair{\Y}{\stopology{\topology{}}{\Y}}}}{\asubset}&=\Y,\label{thmdensetransitivity1p1}\\
\func{\Cl{\Xt}}{\Y}&=\X.\label{thmdensetransitivity1p2}
\end{align}
According to
\refthm{thmsubspacedense},
\Ref{thmdensetransitivity1p1}
implies,
\begin{equation}\label{thmdensetransitivity1p3}
\func{\Cl{\Xt}}{\asubset}\supseteq\Y.
\end{equation}
According to
\refthm{thmdensetransitivity}, clearly
\Ref{thmdensetransitivity1p2},
and
\Ref{thmdensetransitivity1p3}
imply,
\begin{equation}
\func{\Cl{\Xt}}{\asubset}\supseteq\X,
\end{equation}
which means,
\begin{equation}
\func{\Cl{\Xt}}{\asubset}=\X.
\end{equation}
\endthm
%%%%%%%%%%%%%%%%%%%%%%%%%%%%%%%%%%%%%%%%%%%%%%%%%%%%%%%%%%%%%%%%%%%%%%%%%%%%%%%%%%%%%%%%%%%%%%%%%%%%%%%%%%%%%%%%%%%
\subsection{
Nowhere-Dense Sets of a Topological Space}
\definition\label{defnowheredenseset}
$\Xt=\opair{\X}{\topology{}}$
is taken as a topological-space,
$\asubset$
as a subset of $\X$.
It is said that $\quotl$$\asubset$ is nowhere dense in the topological-space $\Xt$$\quotr$,
and $\asubset$ is referred to as a $\quotl$nowhere-dense set of $\Xt$$\quotr$, iff
\begin{equation}
\func{\Ext{\Xt}}{\asubset}\in\Fdense{\X}{\topology{}}.
\end{equation}
\endef
%%%%%%%%%%%%%%%%%%%%%%%%%%%%%%%%%%%%%%%%%%%%%%%%%%%%%%%%%%%%%%%%%%%%%%%%%%%%%%%%%%%%%%%%%%%%%%%%%%%%%%%%%%%%%%%%%%%
\definition\label{defclassofnowheredensesets}
$\opair{\X}{\topology{}}$
is taken as a topological-space.
The set of all nowhere-dense sets o $\Xt$
is denoted by $\Fnwdense{\X}{\topology{}}$. That is,
\begin{equation}
\Fnwdense{\X}{\topology{}}:=\defset{\asubset}{\CSs{\X}}
{\[\func{\Ext{\Xt}}{\asubset}\in\Fdense{\X}{\topology{}}\]}.
\end{equation}
\endef
%%%%%%%%%%%%%%%%%%%%%%%%%%%%%%%%%%%%%%%%%%%%%%%%%%%%%%%%%%%%%%%%%%%%%%%%%%%%%%%%%%%%%%%%%%%%%%%%%%%%%%%%%%%%%%%%%%%
\corollary\label{corclassofnowheredensesets0}
$\opair{\X}{\topology{}}$
is taken as a topological-space.
\begin{equation}
\Fnwdense{\X}{\topology{}}:=\defset{\asubset}{\CSs{\X}}
{\[\func{\Cl{\Xt}}{\func{\Ext{\Xt}}{\asubset}}=\X\]}.
\end{equation}
\endcor
%%%%%%%%%%%%%%%%%%%%%%%%%%%%%%%%%%%%%%%%%%%%%%%%%%%%%%%%%%%%%%%%%%%%%%%%%%%%%%%%%%%%%%%%%%%%%%%%%%%%%%%%%%%%%%%%%%%
\theorem\label{thmnwdensesetsproperty1}
$\opair{\X}{\topology{}}$
is taken as a topological-space.
\begin{equation}\label{thmnwdensesetsproperty1eq1}
\Fnwdense{\X}{\topology{}}=\defset{\asubset}{\CSs{\X}}
{\[\func{\Int{\Xt}}{\func{\Cl{\Xt}}{\asubset}}=\empty\]}.
\end{equation}
\prooff
According to
\refthm{thmclosureandextrelation}
and
\refdef{defextofset},
\begin{align}
\Foreach{\asubset}{\CSs{\X}}
\compl{\X}{\[\func{\Cl{\Xt}}{\func{\Ext{\Xt}}{\asubset}}\]}
&=\func{\Ext{\Xt}}{\func{\Ext{\Xt}}{\asubset}}\cr
&=\func{\Int{\Xt}}{\compl{\X}{\func{\Ext{\Xt}}{\asubset}}}\cr
&=\func{\Int{\Xt}}{\func{\Cl{\Xt}}{\asubset}}.
\end{align}
Hence,
\begin{equation}
\Foreach{\asubset}{\CSs{\X}}
\left\{\[\func{\Cl{\Xt}}{\func{\Ext{\Xt}}{\asubset}}=\X\]\thenn
\[\func{\Int{\Xt}}{\func{\Cl{\Xt}}{\asubset}}=\empty\]\right\}.
\end{equation}
Thus according to
\refcor{corclassofnowheredensesets0},
\begin{equation}
\Foreach{\asubset}{\CSs{\X}}
\left\{\[\asubset\in\Fnwdense{\X}{\topology{}}\]\thenn
\[\func{\Int{\Xt}}{\func{\Cl{\Xt}}{\asubset}}=\empty\]\right\}.
\end{equation}
\endthm
%%%%%%%%%%%%%%%%%%%%%%%%%%%%%%%%%%%%%%%%%%%%%%%%%%%%%%%%%%%%%%%%%%%%%%%%%%%%%%%%%%%%%%%%%%%%%%%%%%%%%%%%%%%%%%%%%%%%%
\theorem\label{thmemptysetisanwdenseset}
$\opair{\X}{\topology{}}$
is taken as a topological-space.
\begin{equation}
\empty\in\Fnwdense{\X}{\topology{}}.
\end{equation}
\prooff
According to
\refthm{thmintofemptyset}
and
\refthm{thmclosureofemptyset},
\begin{align}
\func{\Int{\Xt}}{\func{\Cl{\Xt}}{\empty}}&=
\func{\Int{\Xt}}{\empty}\cr
&=\empty,
\end{align}
and hence according to
\refthm{thmnwdensesetsproperty1},
it is clear that
$\empty$
is a nowhere-dense set of $\Xt$.
\endthm
%%%%%%%%%%%%%%%%%%%%%%%%%%%%%%%%%%%%%%%%%%%%%%%%%%%%%%%%%%%%%%%%%%%%%%%%%%%%%%%%%%%%%%%%%%%%%%%%%%%%%%%%%%%%%%%%%%%%%
\theorem\label{thmclosureofanwdensetsetisanwdenseset}
$\opair{\X}{\topology{}}$
is taken as a topological-space.
For every $\asubset$ in $\CSs{\X}$,
$\asubset$
is a nowhere-dense set of $\Xt$ if-and-only-if the closure of
$\asubset$ (in $\Xt$)
is a nowhere-dense set of $\Xt$. That is,
\begin{equation}
\Foreach{\asubset}{\CSs{\X}}
\[\asubset\in\Fnwdense{\X}{\topology{}}\thenn
\func{\Cl{\Xt}}{\asubset}\in\Fnwdense{\X}{\topology{}}\].
\end{equation}
\prooff
$\asubset$
is taken as an arbitrary subset of $\X$.
According to
\refthm{thmclosureofclosureofset},
\begin{equation}
\func{\Int{\Xt}}{\func{\Cl{\Xt}}{\func{\Cl{\Xt}}{\asubset}}}=
\func{\Int{\Xt}}{\func{\Cl{\Xt}}{\asubset}},
\end{equation}
and hence
\begin{equation}
\func{\Int{\Xt}}{\func{\Cl{\Xt}}{\func{\Cl{\Xt}}{\asubset}}}=\empty\thenn
\func{\Int{\Xt}}{\func{\Cl{\Xt}}{\asubset}}=\empty,
\end{equation}
and thus according to,
\refthm{thmnwdensesetsproperty1},
This means,
\begin{equation}
\asubset\in\Fnwdense{\X}{\topology{}}\thenn
\func{\Cl{\Xt}}{\asubset}\in\Fnwdense{\X}{\topology{}}.
\end{equation}
\endthm
%%%%%%%%%%%%%%%%%%%%%%%%%%%%%%%%%%%%%%%%%%%%%%%%%%%%%%%%%%%%%%%%%%%%%%%%%%%%%%%%%%%%%%%%%%%%%%%%%%%%%%%%%%%%%%%%%%%%%
\theorem\label{thmnwdensesetsproperty2}
$\opair{\X}{\topology{}}$
is taken as a topological-space.
\begin{equation}
\Foreach{\asubset}{\Fnwdense{\X}{\topology{}}}
\func{\Int{\Xt}}{\asubset}=\empty.
\end{equation}
\prooff
$\asubset$
is taken as an arbitrary subset of $\X$.
Considering that $\asubset\subseteq\func{\Cl{\Xt}}{\asubset}$,
According to \refthm{thmintofasubsetofset},
\begin{equation}
\func{\Int{\Xt}}{\asubset}\subseteq
\func{\Int{\Xt}}{\func{\Cl{\Xt}}{\asubset}},
\end{equation}
and hence according to \refthm{thmnwdensesetsproperty1},
\begin{equation}
\func{\Int{\Xt}}{\asubset}=\empty.
\end{equation}
\endthm
%%%%%%%%%%%%%%%%%%%%%%%%%%%%%%%%%%%%%%%%%%%%%%%%%%%%%%%%%%%%%%%%%%%%%%%%%%%%%%%%%%%%%%%%%%%%%%%%%%%%%%%%%%%%%%%%%%%%%
\theorem\label{thmnwdensesetsproperty3}
$\opair{\X}{\topology{}}$
is taken as a topological-space, and
$\asubset$
as a subset of $\X$.
\begin{gather}
\[\asubset\in\Fnwdense{\X}{\topology{}}\]\cr
\vthenn\cr
\left\{\Foreach{\U}{\(\compl{\topology{}}{\seta{\empty}}\)}\[\Exists{\V}
{\[\(\compl{\topology{}}{\seta{\empty}}\)\cap\CSs{\U}\]}\V\cap\asubset=\empty\]\right\}.
\end{gather}
\prooff
\begin{itemize}
\item[${\textbf{\textsf{p1}}}$]
$\asubset$
is taken as an arbitrary element of $\Fnwdense{\X}{\topology{}}$.
Then according to
\refthm{thmnwdensesetsproperty1},
\begin{equation}\label{thmnwdensesetsproperty3p1eq1}
\func{\Int{\Xt}}{\func{\Cl{\Xt}}{\asubset}}=\empty.
\end{equation}
\begin{itemize}
\item[${\textbf{\textsf{p1-1}}}$]
$\U$
is taken as an arbitrary element of $\(\compl{\topology{}}{\seta{\empty}}\)$
(a non-empty open set of $\X$).
Then according to \Ref{thmnwdensesetsproperty3p1eq1},
\refthm{thmintofsetissetofintpoints},
and
\refdef{definteriorpoint},
\begin{equation}
\U\nsubseteq\func{\Cl{\Xt}}{\asubset},
\end{equation}
which means,
\begin{equation}
\U\cap\(\compl{\X}{\func{\Cl{\Xt}}{\asubset}}\)\neq\empty,
\end{equation}
and hence according to
\refthm{thmclosureandextrelation},
\begin{equation}
\U\cap\func{\Ext{\Xt}}{\asubset}\neq\empty.
\end{equation}
Thus considering that $\func{\Ext{\Xt}}{\asubset}\in\topology{}$ and
$\U\in\topology{}$,
based on \refdef{deftopologicalspace},
it is clear that,
\begin{equation}
\[\U\cap\func{\Ext{\Xt}}{\asubset}\]\in\(\compl{\topology{}}{\seta{\empty}}\).
\end{equation}
In addition, it is clear that,
\begin{equation}
\[\U\cap\func{\Ext{\Xt}}{\asubset}\]\in\CSs{\U}.
\end{equation}
Additionally, considering that,
\begin{equation}
\[\U\cap\func{\Ext{\Xt}}{\asubset}\]\subseteq\func{\Ext{\Xt}}{\asubset},
\end{equation}
and
\begin{equation}
\func{\Ext{\Xt}}{\asubset}\cap\asubset=\empty,
\end{equation}
it is clear that,
\begin{equation}
\[\U\cap\func{\Ext{\Xt}}{\asubset}\]\cap\asubset=\empty.
\end{equation}
Therefore,
\begin{equation}
\Exists{\V}
{\[\(\compl{\topology{}}{\seta{\empty}}\)\cap\CSs{\U}\]}\V\cap\asubset=\empty.
\end{equation}
\endp
\end{itemize}
\endp
\item[${\textbf{\textsf{p2}}}$]
It is assumed that,
\begin{equation}\label{thmnwdensesetsproperty3p2eq1}
\Foreach{\U}{\(\compl{\topology{}}{\seta{\empty}}\)}\[\Exists{\V}
{\[\(\compl{\topology{}}{\seta{\empty}}\)\cap\CSs{\U}\]}\V\cap\asubset=\empty\].
\end{equation}
\begin{itemize}
\item[${\textbf{\textsf{p2-1}}}$]
$\point$
is taken as an arbitrary element of $\X$.
\begin{itemize}
\item[${\textbf{\textsf{p2-1-1}}}$]
$\U$
is taken as an arbitrary neighbourhood of $\seta{\point}$ in $\Xt$
(an arbitrary element of $\func{\nei{\Xt}}{\seta{\point}}$).
It is clear that $\U$
is a non-empty open set of $\Xt$ (because it contains $\point$). That is,
\begin{equation}
\U\in\(\compl{\topology{}}{\seta{\empty}}\).
\end{equation}
Thus according to \Ref{thmnwdensesetsproperty3p2eq1},
\begin{equation}
\Existsis{\V}
{\[\(\compl{\topology{}}{\seta{\empty}}\)\cap\CSs{\U}\]}\V\subseteq\(\compl{\X}{\asubset}\).
\end{equation}
So it is clear that,
\begin{equation}
\V\in\[\topology{}\cap\CSs{\compl{\X}{\asubset}}\],
\end{equation}
and hence according to
\refcor{corintofset0},
\begin{equation}
\V\subseteq\func{\Int{\Xt}}{\compl{\X}{\asubset}},
\end{equation}
and thus according to \refdef{defextofset},
\begin{equation}
\V\subseteq\func{\Ext{\Xt}}{\asubset}.
\end{equation}
Thus considering that,
$\V\neq\empty$
and
$\V\subseteq\U$,
it is clear that,
\begin{equation}
\U\cap\func{\Ext{\Xt}}{\asubset}\neq\empty.
\end{equation}
\endp
\end{itemize}
Therefore,
\begin{equation}
\Foreach{\U}{\func{\nei{\Xt}}{\seta{\point}}}
\U\cap\func{\Ext{\Xt}}{\asubset}\neq\empty.
\end{equation}
Thus according to
\refthm{thmclosureofsetissetofadhpoints}
and
\refdef{defadherentpoint},
\begin{equation}
\point\in\func{\Cl{\Xt}}{\func{\Ext{\Xt}}{\asubset}},
\end{equation}
\endp
\end{itemize}
Therefore,
\begin{equation}
\Foreach{\point}{\X}
\point\in\func{\Cl{\Xt}}{\func{\Ext{\Xt}}{\asubset}},
\end{equation}
which means,
\begin{equation}
\func{\Cl{\Xt}}{\func{\Ext{\Xt}}{\asubset}}=\X.
\end{equation}
Hence according to
\refcor{corclassofnowheredensesets0},
\begin{equation}
\asubset\in\Fnwdense{\X}{\topology{}}.
\end{equation}
\endp
\end{itemize}
\endthm
%%%%%%%%%%%%%%%%%%%%%%%%%%%%%%%%%%%%%%%%%%%%%%%%%%%%%%%%%%%%%%%%%%%%%%%%%%%%%%%%%%%%%%%%%%%%%%%%%%%%%%%%%%%%%%%%%%%%%
\theorem\label{thmsubsetofnwdenseset}
$\opair{\X}{\topology{}}$
is taken as a topological-space.
Every subset of every nowhere-dense set of $\Xt$
is a nowhere-dense set of $\Xt$.
\begin{equation}
\Foreach{\asubset}{\Fnwdense{\X}{\topology{}}}
\[\Foreach{\bsubset}{\CSs{\asubset}}\bsubset\in\Fnwdense{\X}{\topology{}}\].
\end{equation}
\prooff
\begin{itemize}
\item[${\textbf{\textsf{p1}}}$]
$\asubset$
is taken as an arbitrary element of $\Fnwdense{\X}{\topology{}}$.
Then according to \refthm{thmnwdensesetsproperty1},
\begin{equation}\label{thmsubsetofnwdensesetp1eq1}
\func{\Int{\Xt}}{\func{\Cl{\Xt}}{\asubset}}=\empty.
\end{equation}
\begin{itemize}
\item[${\textbf{\textsf{p1-1}}}$]
$\bsubset$
is taken as an arbitrary subset of $\asubset$.
Then according to \refthm{thmclosureofasubsetofset},
\begin{equation}
\func{\Cl{\Xt}}{\bsubset}\subseteq
\func{\Cl{\Xt}}{\asubset},
\end{equation}
and hence according to \refthm{thmintofasubsetofset},
\begin{equation}
\func{\Int{\Xt}}{\func{\Cl{\Xt}}{\bsubset}}\subseteq
\func{\Int{\Xt}}{\func{\Cl{\Xt}}{\asubset}}.
\end{equation}
This and \Ref{thmsubsetofnwdensesetp1eq1}
imply,
\begin{equation}
\func{\Int{\Xt}}{\func{\Cl{\Xt}}{\bsubset}}=\empty,
\end{equation}
which according to \refthm{thmnwdensesetsproperty1}, means,
\begin{equation}
\bsubset\in\Fnwdense{\X}{\topology{}}.
\end{equation}
\endp
\end{itemize}
\endp
\end{itemize}
\endthm
%%%%%%%%%%%%%%%%%%%%%%%%%%%%%%%%%%%%%%%%%%%%%%%%%%%%%%%%%%%%%%%%%%%%%%%%%%%%%%%%%%%%%%%%%%%%%%%%%%%%%%%%%%%%%%%%%%%%%
\theorem\label{thmunionofnwdensesets}
$\opair{\X}{\topology{}}$
is taken as a topological-space.
The union of any pair of nowhere-dense sets of $\Xt$
is a nowhere-dense set of $\Xt$. That is,
\begin{equation}
\Foreach{\opair{\asubset}{\p{\asubset}}}
{\[\Fnwdense{\X}{\topology{}}\times\Fnwdense{\X}{\topology{}}\]}
\(\asubset\cup\p{\asubset}\)\in\Fnwdense{\X}{\topology{}}.
\end{equation}
\prooff
$\asubset$
and
$\p{\asubset}$
are taken as arbitrary element of $\Fnwdense{\X}{\topology{}}$.
Then according to
\refthm{thmnwdensesetsproperty3},
\begin{align}
&\Foreach{\U}{\(\compl{\topology{}}{\seta{\empty}}\)}\[\Exists{\V}
{\[\(\compl{\topology{}}{\seta{\empty}}\)\cap\CSs{\U}\]}\V\cap\asubset=\empty\],\label{thmunionofnwdensesetspeq1}\\
&\Foreach{\U}{\(\compl{\topology{}}{\seta{\empty}}\)}\[\Exists{\V}
{\[\(\compl{\topology{}}{\seta{\empty}}\)\cap\CSs{\U}\]}\V\cap\p{\asubset}=\empty\].\label{thmunionofnwdensesetspeq2}
\end{align}
\begin{itemize}
\item[${\textbf{\textsf{p1}}}$]
$\U$
is taken as an arbitrary element of $\(\compl{\topology{}}{\seta{\empty}}\)$.
According to \Ref{thmunionofnwdensesetspeq1},
\begin{equation}\label{thmunionofnwdensesetsp1eq1}
\Existsis{\V}
{\[\(\compl{\topology{}}{\seta{\empty}}\)\cap\CSs{\U}\]}\V\cap\asubset=\empty.
\end{equation}
Thus considering that $\V\in\(\compl{\topology{}}{\seta{\empty}}\)$,
according to \Ref{thmunionofnwdensesetspeq2},
\begin{equation}\label{thmunionofnwdensesetsp1eq2}
\Existsis{\p{\V}}
{\[\(\compl{\topology{}}{\seta{\empty}}\)\cap\CSs{\V}\]}\p{\V}\cap\p{\asubset}=\empty.
\end{equation}
Considering that $\p{\V}\subseteq\V$
and
$\V\subseteq\U$,
it is clear that,
\begin{equation}\label{thmunionofnwdensesetsp1eq3}
\p{\V}\subseteq\U.
\end{equation}
In addition, considering that $\p{\V}\subseteq\V$ and $\V\cap\asubset=\empty$,
it is clear that,
\begin{equation}\label{thmunionofnwdensesetsp1eq4}
\p{\V}\cap\asubset=\empty.
\end{equation}
This and \Ref{thmunionofnwdensesetsp1eq2} imply that,
\begin{align}\label{thmunionofnwdensesetsp1eq5}
\p{\V}\cap\(\asubset\cup\p{\asubset}\)&=
\(\p{\V}\cap\asubset\)\cup\(\p{\V}\cap\p{\asubset}\)\cr
&=\empty\cup\empty\cr
&=\empty.
\end{align}
\Ref{thmunionofnwdensesetsp1eq2},
\Ref{thmunionofnwdensesetsp1eq3},
and
\Ref{thmunionofnwdensesetsp1eq5}
imply,
\begin{equation}
\Existsis{\p{\V}}
{\[\(\compl{\topology{}}{\seta{\empty}}\)\cap\CSs{\U}\]}
\p{\V}\cap\(\asubset\cup\p{\asubset}\)=\empty.
\end{equation}
\endp
\end{itemize}
Therefore,
\begin{equation}
\Foreach{\U}{\(\compl{\topology{}}{\seta{\empty}}\)}\[\Exists{\V}
{\[\(\compl{\topology{}}{\seta{\empty}}\)\cap\CSs{\U}\]}\V\cap\(\asubset\cup\p{\asubset}\)=\empty\].
\end{equation}
According to \refthm{thmnwdensesetsproperty3},
This implies,
\begin{equation}
\(\asubset\cup\p{\asubset}\)\in\Fnwdense{\X}{\topology{}}.
\end{equation}
\endthm
%%%%%%%%%%%%%%%%%%%%%%%%%%%%%%%%%%%%%%%%%%%%%%%%%%%%%%%%%%%%%%%%%%%%%%%%%%%%%%%%%%%%%%%%%%%%%%%%%%%%%%%%%%%%%%%%%%%%%
\theorem\label{thmfrontierofclosedset}
$\opair{\X}{\topology{}}$
is taken as a topological-space.
The boundary of every closed set of $\Xt$
is a nowhere-dense set of $\Xt$. That is,
\begin{equation}
\Foreach{\asubset}{\Fclosed{\X}{\topology{}}}
\[\func{\Fr{\Xt}}{\asubset}\in\Fnwdense{\X}{\topology{}}\].
\end{equation}
\prooff
$\asubset$
is taken as an arbitrary element of $\Fclosed{\X}{\topology{}}$.
According to \refthm{thmfrontierofsetisclosed},
the boundary of every subset of $\X$
in the topological-space $\Xt$, is a closed set of $\Xt$. Thus,
\begin{equation}\label{thmfrontierofclosedsetpeq1}
\func{\Fr{\Xt}}{\asubset}\in\Fclosed{\X}{\topology{}},
\end{equation}
and hence according to
\refthm{thmclosureofclosedset},
\begin{equation}\label{thmfrontierofclosedsetpeq2}
\func{\Cl{\Xt}}{\func{\Fr{\Xt}}{\asubset}}=
\func{\Fr{\Xt}}{\asubset}.
\end{equation}
In addition, considering that
$\asubset\in\Fclosed{\X}{\topology{}}$,
according to
\refthm{thmintoffrontierofclosedset},
\begin{equation}\label{thmfrontierofclosedsetpeq3}
\func{\Int{\Xt}}{\func{\Fr{\Xt}}{\asubset}}=\empty.
\end{equation}
\Ref{thmfrontierofclosedsetpeq2}
and
\Ref{thmfrontierofclosedsetpeq3}
imply,
\begin{equation}
\func{\Int{\Xt}}{\func{\Cl{\Xt}}{\func{\Fr{\Xt}}{\asubset}}}=\empty.
\end{equation}
According to
\refthm{thmnwdensesetsproperty1},
this implies,
\begin{equation}
\func{\Fr{\Xt}}{\asubset}\in\Fnwdense{\X}{\topology{}}.
\end{equation}
\endthm
%%%%%%%%%%%%%%%%%%%%%%%%%%%%%%%%%%%%%%%%%%%%%%%%%%%%%%%%%%%%%%%%%%%%%%%%%%%%%%%%%%%%%%%%%%%%%%%%%%%%%%%%%%%%%%%%%%%%%
\theorem\label{thmfrontierofopenset}
$\opair{\X}{\topology{}}$
is taken as a topological-space.
The boundary of every open set of $\Xt$
is a nowhere-dense set of $\Xt$. That is,
\begin{equation}
\Foreach{\asubset}{\topology{}}
\[\func{\Fr{\Xt}}{\asubset}\in\Fnwdense{\X}{\topology{}}\].
\end{equation}
\prooff
$\asubset$
is taken as an arbitrary element of $\topology{}$. Then according to \refdef{deffamilyofclosedsets},
\begin{equation}
\(\compl{\X}{\asubset}\)\in\Fclosed{\X}{\topology{}},
\end{equation}
and according to \refthm{thmfrontierofclosedset},
\begin{equation}
\func{\Fr{\Xt}}{\compl{\X}{\asubset}}\in\Fnwdense{\X}{\topology{}}.
\end{equation}
In addition, according to
\refthm{thmfrontierofcomplementofset},
\begin{equation}
\func{\Fr{\Xt}}{\compl{\X}{\asubset}}=\func{\Fr{\Xt}}{\asubset}.
\end{equation}
Therefore,
\begin{equation}
\func{\Fr{\Xt}}{\asubset}\in\Fnwdense{\X}{\topology{}}.
\end{equation}
\endthm
%%%%%%%%%%%%%%%%%%%%%%%%%%%%%%%%%%%%%%%%%%%%%%%%%%%%%%%%%%%%%%%%%%%%%%%%%%%%%%%%%%%%%%%%%%%%%%%%%%%%%%%%%%%%%%%%%%%%%
\theorem\label{thmcomplementofclosednwdenseset}
$\opair{\X}{\topology{}}$
is taken as a topological-space.
For every $\asubset$ in $\CSs{\X}$, $\asubset$
is a closed and nowhere-dense set of $\Xt$
if-and-only-if the complement of $\asubset$
is an open and dense set of $\Xt$. That is,
\begin{equation}
\Foreach{\asubset}{\CSs{\X}}
\left\{\asubset\in\Fclosed{\X}{\topology{}}\cap\Fnwdense{\X}{\topology{}}\thenn
\(\compl{\X}{\asubset}\)\in\[\topology{}\cap\Fdense{\X}{\topology{}}\]\right\}.
\end{equation}
\prooff
$\asubset$
is taken as an arbitrary element of $\[\Fclosed{\X}{\topology{}}\cap\Fnwdense{\X}{\topology{}}\]$.
\begin{itemize}
\item[${\textbf{\textsf{p1}}}$]
It is assumed that,
\begin{equation}
\asubset\in\[\Fclosed{\X}{\topology{}}\cap\Fnwdense{\X}{\topology{}}\].
\end{equation}
Then according to
\refthm{thmclosureofclosedset},
\begin{equation}
\func{\Cl{\Xt}}{\asubset}=\asubset.
\end{equation}
In addition, according to
\refthm{thmnwdensesetsproperty1},
\begin{equation}
\func{\Int{\Xt}}{\func{\Cl{\Xt}}{\asubset}}=\empty.
\end{equation}
These imply,
\begin{equation}
\func{\Int{\Xt}}{\asubset}=\empty,
\end{equation}
and hence according to
\refthm{thmclosureandextrelation}
and
\refdef{defextofset},
\begin{align}
\func{\Cl{\Xt}}{\compl{\X}{\asubset}}&=
\compl{\X}{\func{\Ext{\Xt}}{\compl{\X}{\asubset}}}\cr
&=\compl{\X}{\func{\Int{\Xt}}{\asubset}}\cr
&=\compl{\X}{\empty}\cr
&=\X.
\end{align}
According to
\refdef{defclassofdensesets},
this implies,
\begin{equation}
\(\compl{\X}{\asubset}\)\in\Fdense{\X}{\topology{}}.
\end{equation}
In addition, considering that $\asubset$ is a closed set of $\Xt$, it is clear that,
\begin{equation}
\asubset\in\topology{}.
\end{equation}
\endp
\item[${\textbf{\textsf{p2}}}$]
It is assumed that,
\begin{equation}
\asubset\in\[\topology{}\cap\Fdense{\X}{\topology{}}\].
\end{equation}
Then,
\begin{equation}
\(\compl{\X}{\asubset}\)\in\Fclosed{\X}{\topology{}},
\end{equation}
and hence according to
\refthm{thmclosureofclosedset},
\begin{equation}
\func{\Cl{\Xt}}{\compl{\X}{\asubset}}=\(\compl{\X}{\asubset}\),
\end{equation}
and hence according to \refthm{thmclosureandextrelation} and \refdef{defextofset},
\begin{align}
\func{\Int{\Xt}}{\func{\Cl{\Xt}}{\compl{\X}{\asubset}}}&=
\func{\Int{\Xt}}{\compl{\X}{\asubset}}\cr
&=\func{\Ext{\Xt}}{\asubset}.
\end{align}
In addition, considering that
$\asubset\in\Fdense{\X}{\topology{}}$,
according to \refthm{thmdenseNandSconditions0},
\begin{equation}
\func{\Ext{\Xt}}{\asubset}=\empty.
\end{equation}
Therefore,
\begin{equation}
\func{\Int{\Xt}}{\func{\Cl{\Xt}}{\compl{\X}{\asubset}}}=\empty,
\end{equation}
and hence according to \refthm{thmnwdensesetsproperty1},
\begin{equation}
\(\compl{\X}{\asubset}\)\in\Fnwdense{\X}{\topology{}}.
\end{equation}
\endp
\end{itemize}
\endthm
%%%%%%%%%%%%%%%%%%%%%%%%%%%%%%%%%%%%%%%%%%%%%%%%%%%%%%%%%%%%%%%%%%%%%%%%%%%%%%%%%%%%%%%%%%%%%%%%%%%%%%%%%%%%%%%%%%%%%
\theorem\label{thmintofcomplementofanwdenseset}
$\opair{\X}{\topology{}}$
is taken as a topological-space.
The interior of complement of every nowhere-dense set of $\Xt$
(in the topological-space $\Xt$)
is a dense set of $\Xt$. That is,
\begin{equation}
\Foreach{\asubset}{\Fnwdense{\X}{\topology{}}}
\func{\Int{\Xt}}{\compl{\X}{\asubset}}\in\Fdense{\X}{\topology{}}.
\end{equation}
\prooff
$\asubset$
is taken as an arbitrary element of $\Fnwdense{\X}{\topology{}}$. Then according to
\refthm{thmclosureofanwdensetsetisanwdenseset},
\begin{equation}
\func{\Cl{\Xt}}{\asubset}\in\Fnwdense{\X}{\topology{}}.
\end{equation}
Considering that
$\func{\Cl{\Xt}}{\asubset}\in\Fclosed{\X}{\topology{}}$,
this and
\refthm{thmcomplementofclosednwdenseset}
imply that,
\begin{equation}
\(\compl{\X}{\func{\Cl{\Xt}}{\asubset}}\)\in\Fdense{\X}{\topology{}}.
\end{equation}
Hence according to \refthm{thmclosureandextrelation},
\begin{equation}
\func{\Ext{\Xt}}{\asubset}\in\Fdense{\X}{\topology{}}.
\end{equation}
According to
\refdef{defclassofdensesets},
this means,
\begin{equation}
\func{\Cl{\Xt}}{\func{\Ext{\Xt}}{\asubset}}=\X,
\end{equation}
and hence according to \refdef{defextofset},
\begin{equation}
\func{\Cl{\Xt}}{\func{\Int{\Xt}}{\compl{\X}{\asubset}}}=\X.
\end{equation}
According to \refdef{defclassofdensesets},
this means,
\begin{equation}
\func{\Int{\Xt}}{\compl{\X}{\asubset}}\in\Fdense{\X}{\topology{}}.
\end{equation}
\endthm
%%%%%%%%%%%%%%%%%%%%%%%%%%%%%%%%%%%%%%%%%%%%%%%%%%%%%%%%%%%%%%%%%%%%%%%%%%%%%%%%%%%%%%%%%%%%%%%%%%%%%%%%%%%%%%%%%%%%%
%%%%%%%%%%%%%%%%%%%%%%%%%%%%%%%%%%%%%%%%%%%%%%%%%%%%%%%%%%%%%%%%%%%%%%%%%%%%%%%%%%%%%%%%%%%%%%%%%%%%%%%%%%%%%%%%%%%%%
%%%%%%%%%%%%%%%%%%%%%%%%%%%%%%%%%%%%%%%%%%%%%%%%%%%%%%%%%%%%%%%%%%%%%%%%%%%%%%%%%%%%%%%%%%%%%%%%%%%%%%%%%%%%%%%%%%%%%
%%%%%%%%%%%%%%%%%%%%%%%%%%%%%%%%%%%%%%%%%%%%%%%%%%%%%%%%%%%%%%%%%%%%%%%%%%%%%%%%%%%%%%%%%%%%%%%%%%%%%%%%%%%%%%%%%%%%%
%%%%%%%%%%%%%%%%%%%%%%%%%%%%%%%%%%%%%%%%%%%%%%%%%%%%%%%%%%%%%%%%%%%%%%%%%%%%%%%%%%%%%%%%%%%%%%%%%%%%%%%%%%%%%%%%%%%%%
%%%%%%%%%%%%%%%%%%%%%%%%%%%%%%%%%%%%%%%%%%%%%%%%%%%%%%%%%%%%%%%%%%%%%%%%%%%%%%%%%%%%%%%%%%%%%%%%%%%%%%%%%%%%%%%%%%%%%
\theorem
$\opair{\X}{\topology{}}$
is taken as a topological-space.
The complement of every nowhere-dense set of $\Xt$
is a dense set of $\Xt$. That is.
\begin{equation}
\Foreach{\asubset}{\Fnwdense{\X}{\topology{}}}
\(\compl{\X}{\asubset}\)\in\Fdense{\X}{\topology{}}.
\end{equation}
\prooff
According to
\refthm{thmifintofsetisdensethensetisdense}
and
\refthm{thmintofcomplementofanwdenseset},
it is trivial.
\endthm
%%%%%%%%%%%%%%%%%%%%%%%%%%%%%%%%%%%%%%%%%%%%%%%%%%%%%%%%%%%%%%%%%%%%%%%%%%%%%%%%%%%%%%%%%%%%%%%%%%%%%%%%%%%%%%%%%%%%%
%%%%%%%%%%%%%%%%%%%%%%%%%%%%%%%%%%%%%%%%%%%%%%%%%%%%%%%%%%%%%%%%%%%%%%%%%%%%%%%%%%%%%%%%%%%%%%%%%%%%%%%%%%%%%%%%%%%%%
%%%%%%%%%%%%%%%%%%%%%%%%%%%%%%%%%%%%%%%%%%%%%%%%%%%%%%%%%%%%%%%%%%%%%%%%%%%%%%%%%%%%%%%%%%%%%%%%%%%%%%%%%%%%%%%%%%%%%
%%%%%%%%%%%%%%%%%%%%%%%%%%%%%%%%%%%%%%%%%%%%%%%%%%%%%%%%%%%%%%%%%%%%%%%%%%%%%%%%%%%%%%%%%%%%%%%%%%%%%%%%%%%%%%%%%%%%%
%%%%%%%%%%%%%%%%%%%%%%%%%%%%%%%%%%%%%%%%%%%%%%%%%%%%%%%%%%%%%%%%%%%%%%%%%%%%%%%%%%%%%%%%%%%%%%%%%%%%%%%%%%%%%%%%%%%%%
%%%%%%%%%%%%%%%%%%%%%%%%%%%%%%%%%%%%%%%%%%%%%%%%%%%%%%%%%%%%%%%%%%%%%%%%%%%%%%%%%%%%%%%%%%%%%%%%%%%%%%%%%%%%%%%%%%%%%
\section{Product Topology}
\theorem
The set $\defSet{\Xt_{\alpha}=\opair{\X_\alpha}{\topology{\alpha}}}{\alpha\in\index}$ is taken as a collection of
topological spaces parametrized by a set $\index$. Denoting by $\Proj{\alpha}{}$ the projection map of
$\Cproduct{\defSet{\X_\alpha}{\alpha\in\index}}$ onto $\X_{\alpha}$, for every $\alpha\in\index$,
the set $\defSet{\func{\pimage{{\Proj{\alpha}{}}}}{\U_{\alpha}}}{\alpha\in\index,~\U_{\alpha}\in\topology{\alpha}}$
is a sub-base for a topology on $\Cproduct{\defSet{\X_\alpha}{\alpha\in\index}}$. Consequently,
the set $\displaystyle\defSet{\CProduct{\alpha\in\index}{\U_{\alpha}}}{\[\Foreach{\alpha}{\index}\U_{\alpha}\in\topology{\alpha}\]}$
is a base for this topology.
\proof
It is left as an exercise.
\endthm
%%%%%%%%%%%%%%%%%%%%%%%%%%%%%%%%%%%%%%%%%%%%%%%%%%%%%%%%%%%%%%%%%%%%%%%%%%%%%%%%%%%%%%%%%%%%%%%%%%%%%%%%%%%%%%%%%%%%%
\definition
The set $\defSet{\Xt_{\alpha}=\opair{\X_\alpha}{\topology{\alpha}}}{\alpha\in\index}$ is taken as a collection of
topological spaces parametrized by a set $\index$. The topology on
$\Cproduct{\defSet{\X_\alpha}{\alpha\in\index}}$ generated by\\
$\displaystyle\defSet{\CProduct{\alpha\in\index}{\U_{\alpha}}}{\[\Foreach{\alpha}{\index}\U_{\alpha}\in\topology{\alpha}\]}$
will be denoted by $\producttop{\defSet{\Xt_{\alpha}}{\alpha\in\index}}$ and called the
$\quotl$product topology of the collection $\defSet{\Xt_{\alpha}=\opair{\X_\alpha}{\topology{\alpha}}}{\alpha\in\index}$
of topological spaces$\quotr$.
\begin{align}
\producttop{\defSet{\Xt_{\alpha}}{\alpha\in\index}}=
\topgen{\Cproduct{\defSet{\X_\alpha}{\alpha\in\index}}}
{\defSet{\CProduct{\alpha\in\index}{\U_{\alpha}}}{\[\Foreach{\alpha}{\index}\U_{\alpha}\in\topology{\alpha}\]}}.
\end{align}
The topological space
$\opair{\Cproduct{\defSet{\X_\alpha}{\alpha\in\index}}}{\producttop{\defSet{\Xt_{\alpha}}{\alpha\in\index}}}$
will be denoted by $\CProduct{\alpha\in\index}{\Xt_{\alpha}}$.
Specifically, when $\index=\seta{\binary{1}{2}}$, the topological product of $\seta{\binary{\Xt_1}{\Xt_2}}$
will be denoted alternatively by $\topprod{\Xt_1}{\Xt_2}$.
\endef
%%%%%%%%%%%%%%%%%%%%%%%%%%%%%%%%%%%%%%%%%%%%%%%%%%%%%%%%%%%%%%%%%%%%%%%%%%%%%%%%%%%%%%%%%%%%%%%%%%%%%%%%%%%%%%%%%%%%%
\theorem
Each
$\Xt_{1}=\opair{\X_1}{\topology{1}}$ and $\Xt_{2}=\opair{\X_2}{\topology{2}}$ is taken as a topological space.
For every $\asubset_1\subseteq\X_1$ and $\asubset_2\subseteq\X_2$,
the product topology of
$\seta{\binary{\opair{\asubset_1}{\stopology{\topology{1}}{\asubset_1}}}{\opair{\asubset_2}{\stopology{\topology{2}}{\asubset_2}}}}$
is the same as a the topology on $\Cprod{\asubset_1}{\asubset_2}$ induced from
$\topprod{\Xt_1}{\Xt_2}$. That is,
\begin{align}
\producttop{\seta{\binary{\opair{\asubset_1}{\stopology{\topology{1}}{\asubset_1}}}{\opair{\asubset_2}{\stopology{\topology{2}}{\asubset_2}}}}}=
\stopology{\producttop{\seta{\binary{\Xt_1}{\Xt_2}}}}{\Cprod{\asubset_1}{\asubset_2}}.
\end{align}
\proof
It is left as an exercise.
\endthm
%%%%%%%%%%%%%%%%%%%%%%%%%%%%%%%%%%%%%%%%%%%%%%%%%%%%%%%%%%%%%%%%%%%%%%%%%%%%%%%%%%%%%%%%%%%%%%%%%%%%%%%%%%%%%%%%%%%%%
%%%%%%%%%%%%%%%%%%%%%%%%%%%%%%%%%%%%%%%%%%%%%%%%%%%%%%%%%%%%%%%%%%%%%%%%%%%%%%%%%%%%%%%%%%%%%%%%%%%%%%%%%%%%%%%%%%%%%
%%%%%%%%%%%%%%%%%%%%%%%%%%%%%%%%%%%%%%%%%%%%%%%%%%%%%%%%%%%%%%%%%%%%%%%%%%%%%%%%%%%%%%%%%%%%%%%%%%%%%%%%%%%%%%%%%%%%%
%%%%%%%%%%%%%%%%%%%%%%%%%%%%%%%%%%%%%%%%%%%%%%%%%%%%%%%%%%%%%%%%%%%%%%%%%%%%%%%%%%%%%%%%%%%%%%%%%%%%%%%%%%%%%%%%%%%%%
%%%%%%%%%%%%%%%%%%%%%%%%%%%%%%%%%%%%%%%%%%%%%%%%%%%%%%%%%%%%%%%%%%%%%%%%%%%%%%%%%%%%%%%%%%%%%%%%%%%%%%%%%%%%%%%%%%%%%
%%%%%%%%%%%%%%%%%%%%%%%%%%%%%%%%%%%%%%%%%%%%%%%%%%%%%%%%%%%%%%%%%%%%%%%%%%%%%%%%%%%%%%%%%%%%%%%%%%%%%%%%%%%%%%%%%%%%%
\section{Quotient Topology}
\definition
$\Xt=\opair{\X}{\topology{\X}}$ is taken as a topological space, and $\eqrel{}$ as an equivalence relation on $\X$.
\begin{equation}
\quotienttop{\Xt}{\eqrel{}}:=\defsets{\U}{\EqClass{\X}{\eqrel{}}}{\(\Union{u}{\U}{u}\)\in\topology{}}.
\end{equation}
\endef
%%%%%%%%%%%%%%%%%%%%%%%%%%%%%%%%%%%%%%%%%%%%%%%%%%%%%%%%%%%%%%%%%%%%%%%%%%%%%%%%%%%%%%%%%%%%%%%%%%%%%%%%%%%%%%%%%%%%%
\theorem
$\Xt=\opair{\X}{\topology{\X}}$ is taken as a topological space, and $\eqrel{}$ as an equivalence relation on $\X$.
$\quotienttop{\Xt}{\eqrel{}}$ is a topology on the set of all equivalence classes of $\eqrel{}$, that is
$\EqClass{\X}{\eqrel{}}$.
\begin{equation}
\quotienttop{\Xt}{\eqrel{}}\in\Ctops{\EqClass{\X}{\eqrel{}}}.
\end{equation}
\proof
It is left as an exercise.
\endthm
%%%%%%%%%%%%%%%%%%%%%%%%%%%%%%%%%%%%%%%%%%%%%%%%%%%%%%%%%%%%%%%%%%%%%%%%%%%%%%%%%%%%%%%%%%%%%%%%%%%%%%%%%%%%%%%%%%%%%
\definition
$\Xt=\opair{\X}{\topology{\X}}$ is taken as a topological space, and $\eqrel{}$ as an equivalence relation on $\X$.
$\quotienttop{\Xt}{\eqrel{}}$ is referred to as the $\quotl$quotient topology of the topological space $\Xt$
relative to the equivalence relation $\eqrel{}$$\quotr$. Furthermore, the topological space
$\opair{\EqClass{\X}{\eqrel{}}}{\quotienttop{\Xt}{\eqrel{}}}$ will be denoted by $\topq{\Xt}{\eqrel{}}$
and called $\quotl$quotient topological space of $\Xt$ relative to $\eqrel{}$$\quotr$.
\endef
%%%%%%%%%%%%%%%%%%%%%%%%%%%%%%%%%%%%%%%%%%%%%%%%%%%%%%%%%%%%%%%%%%%%%%%%%%%%%%%%%%%%%%%%%%%%%%%%%%%%%%%%%%%%%%%%%%%%%
%%%%%%%%%%%%%%%%%%%%%%%%%%%%%%%%%%%%%%%%%%%%%%%%%%%%%%%%%%%%%%%%%%%%%%%%%%%%%%%%%%%%%%%%%%%%%%%%%%%%%%%%%%%%%%%%%%%%%
%%%%%%%%%%%%%%%%%%%%%%%%%%%%%%%%%%%%%%%%%%%%%%%%%%%%%%%%%%%%%%%%%%%%%%%%%%%%%%%%%%%%%%%%%%%%%%%%%%%%%%%%%%%%%%%%%%%%%
%%%%%%%%%%%%%%%%%%%%%%%%%%%%%%%%%%%%%%%%%%%%%%%%%%%%%%%%%%%%%%%%%%%%%%%%%%%%%%%%%%%%%%%%%%%%%%%%%%%%%%%%%%%%%%%%%%%%%
%%%%%%%%%%%%%%%%%%%%%%%%%%%%%%%%%%%%%%%%%%%%%%%%%%%%%%%%%%%%%%%%%%%%%%%%%%%%%%%%%%%%%%%%%%%%%%%%%%%%%%%%%%%%%%%%%%%%%
%%%%%%%%%%%%%%%%%%%%%%%%%%%%%%%%%%%%%%%%%%%%%%%%%%%%%%%%%%%%%%%%%%%%%%%%%%%%%%%%%%%%%%%%%%%%%%%%%%%%%%%%%%%%%%%%%%%%%
\section{Metric Topology and Metrizable Spaces}
\definition
$\metricspace{}=\opair{\M}{d}$ is taken as a metric space. For every $\point\in\M$ and every positive real number $r$,
\begin{equation}
\metricball{\metricspace{}}{\point}{r}:=\defset{\x}{\M}{\func{d}{\binary{\point}{\x}}<r},
\end{equation}
and is called the $\quotl$open ball of the metric space $\Xt$ centered at $\point$ with radius $r$$\quotr$.
\endef
%%%%%%%%%%%%%%%%%%%%%%%%%%%%%%%%%%%%%%%%%%%%%%%%%%%%%%%%%%%%%%%%%%%%%%%%%%%%%%%%%%%%%%%%%%%%%%%%%%%%%%%%%%%%%%%%%%%%%
\theorem
$\metricspace{}=\opair{\M}{d}$ is taken as a metric space. The set of all open balls of $\metricspace{}$,
that is $\defSet{\metricball{\metricspace{}}{\point}{r}}{\point\in\M,~r\in\Rp}$, is a base for a topology on $\M$.
\proof
It is left as an exercise.
\endthm
%%%%%%%%%%%%%%%%%%%%%%%%%%%%%%%%%%%%%%%%%%%%%%%%%%%%%%%%%%%%%%%%%%%%%%%%%%%%%%%%%%%%%%%%%%%%%%%%%%%%%%%%%%%%%%%%%%%%%
\definition
$\metricspace{}=\opair{\M}{d}$ is taken as a metric space. The topology on $\M$ generated by
$\defSet{\metricball{\metricspace{}}{\point}{r}}{\point\in\M,~r\in\Rp}$ will be denoted by
$\Metrictop{\metricspace{}}$ and called $\quotl$topology on $\M$ induced by the metric $d$$\quotr$.
\endef
%%%%%%%%%%%%%%%%%%%%%%%%%%%%%%%%%%%%%%%%%%%%%%%%%%%%%%%%%%%%%%%%%%%%%%%%%%%%%%%%%%%%%%%%%%%%%%%%%%%%%%%%%%%%%%%%%%%%%
\definition
$\Xt=\opair{\X}{\topology{\X}}$ is taken as a topological space. $\Xt$ is called a
$\quotl$metrizable topological space$\quotr$ iff there exists a metric on $\X$ so that
$\topology{}$ is induced by it.
\endef
\chapteR{
Continuous Maps}
\thispagestyle{fancy}
\section{
Definition and Basic Properties of Continuous Maps}
\definition\label{defcontinuousfunction}
Each
$\Xt=\opair{\X}{\topology{\X}}$,
and
$\Yt=\opair{\Y}{\topology{\Y}}$
is taken as a topological-space.
For every $\cf$ in $\Func{\X}{\Y}$
(Every function $\cf$ from $\X$ to $\Y$),
$\cf$
is referred to as a $\quotl$continuous mapping from the topological-space $\Xt$
to the topological-space $\Yt$$\quotr$ iff these properties are hold.
\begin{itemize}
\item[${\textbf{\textsf{CF1}}}$]
%\hfill
$\Foreach{\U}{\topology{\Y}}
\func{\pimage{\cf}}{\U}\in\topology{\X}.$
\end{itemize}
In addition, the set of all continuous mappings from $\Xt$ to $\Yt$
is denoted by $\CF{\Xt}{\Yt}$. That is,
\begin{equation}
\CF{\Xt}{\Yt}:=\defset{\cf}{\Func{\X}{\Y}}
{\(\Foreach{\U}{\topology{\Y}}\func{\pimage{\cf}}{\U}\in\topology{\X}\)}.
\end{equation}
\endef
%%%%%%%%%%%%%%%%%%%%%%%%%%%%%%%%%%%%%%%%%%%%%%%%%%%%%%%%%%%%%%%%%%%%%%%%%%%%%%%%%%%%%%%%%%%%%%%%%%
\theorem\label{thmcontiniuityequiv1}
Each
$\Xt=\opair{\X}{\topology{\X}}$
and
$\Yt=\opair{\Y}{\topology{\Y}}$
is taken as a topological-space. For every
$\cf$
in
$\Func{\X}{\Y}$,
$\cf$
is a continuous map from $\Xt$ to $\Yt$ if-and-only-if
the inverse-image of every closed set of $\Yt$ under $\cf$
is a closed set of $\Xt$. That is,
\begin{align}\label{thmcontiniuityequiv1eq1}
\Foreach{\cf}{\Func{\X}{\Y}}
\left\{\[\cf\in\CF{\Xt}{\Yt}\]\thenn
\[\Foreach{\U}{\Fclosed{\Y}{\topology{\Y}}}
\func{\pimage{\cf}}{\U}\in\Fclosed{\X}{\topology{\X}}\]\right\}.
\end{align}
\prooff
$\cf$
is taken as an arbitrary element of $\Func{\X}{\Y}$. According to \refdef{defcontinuousfunction},
\begin{equation}\label{thmcontiniuityequiv1peq1}
\[\cf\in\CF{\Xt}{\Yt}\]\thenn
\(\Foreach{\U}{\topology{\Y}}\func{\pimage{\cf}}{\U}\in\topology{\X}\).
\end{equation}
\begin{itemize}
\item[${\textbf{\textsf{p1}}}$]
It is assumed that,
\begin{equation}\label{thmcontiniuityequiv1p1eq1}
\Foreach{\U}{\topology{\Y}}
\func{\pimage{\cf}}{\U}\in\topology{\X}.
\end{equation}
\begin{itemize}
\item[${\textbf{\textsf{p1-1}}}$]
$\U$
is taken as an arbitrary element of $\Fclosed{\Y}{\topology{\Y}}$.
Then according to
\refdef{deffamilyofclosedsets},
\begin{equation}\label{thmcontiniuityequiv1p1-1eq1}
\(\compl{\Y}{\U}\)\in\topology{\Y},
\end{equation}
and hence according to
\Ref{thmcontiniuityequiv1p1eq1},
\begin{equation}\label{thmcontiniuityequiv1p1-1eq2}
\func{\pimage{\cf}}{\compl{\Y}{\U}}\in\topology{\X}.
\end{equation}
In addition, it is known that,
\begin{equation}\label{thmcontiniuityequiv1p1-1eq3}
\func{\pimage{\cf}}{\compl{\Y}{\U}}=
\[\compl{\X}{\func{\pimage{\cf}}{\U}}\].
\end{equation}
\Ref{thmcontiniuityequiv1p1-1eq2}
and
\Ref{thmcontiniuityequiv1p1-1eq3}
imply that,
\begin{equation}
\[\compl{\X}{\func{\pimage{\cf}}{\U}}\]\in\topology{\X},
\end{equation}
and hence according to
\refdef{deffamilyofclosedsets},
\begin{equation}
\func{\pimage{\cf}}{\U}\in\Fclosed{\X}{\topology{\X}}.
\end{equation}
\endp
\end{itemize}
Therefore,
\begin{equation}
\[\Foreach{\U}{\Fclosed{\Y}{\topology{\Y}}}
\func{\pimage{\cf}}{\U}\in\Fclosed{\X}{\topology{\X}}\].
\end{equation}
\endp
\end{itemize}
Thus,
\begin{equation}
\[\Foreach{\U}{\topology{\Y}}
\func{\pimage{\cf}}{\U}\in\topology{\X}\]\then
\[\Foreach{\U}{\Fclosed{\Y}{\topology{\Y}}}
\func{\pimage{\cf}}{\U}\in\Fclosed{\X}{\topology{\X}}\].
\end{equation}
Similarly, it can be seen that,
\begin{equation}
\[\Foreach{\U}{\Fclosed{\Y}{\topology{\Y}}}
\func{\pimage{\cf}}{\U}\in\Fclosed{\X}{\topology{\X}}\]\then
\[\Foreach{\U}{\topology{\Y}}
\func{\pimage{\cf}}{\U}\in\topology{\X}\].
\end{equation}
These and
\Ref{thmcontiniuityequiv1peq1}
imply,
\begin{equation}
\[\cf\in\CF{\Xt}{\Yt}\]\thenn
\[\Foreach{\U}{\Fclosed{\Y}{\topology{\Y}}}
\func{\pimage{\cf}}{\U}\in\Fclosed{\X}{\topology{\X}}\].
\end{equation}
\endthm
%%%%%%%%%%%%%%%%%%%%%%%%%%%%%%%%%%%%%%%%%%%%%%%%%%%%%%%%%%%%%%%%%%%%%%%%%%%%%%%%%%%%%%%%%%%%%%%%%%
\theorem\label{thmcontiniuityequiv4}
Each
$\Xt=\opair{\X}{\topology{\X}}$,
and
$\Yt=\opair{\Y}{\topology{\Y}}$
is taken as a topological-space.
If
$\cf$
is a function from
$\X$
to
$\Y$, and
$\baseof{\Yt}$ is base for the topological-space $\Yt$, then
$\cf$
is a continuous map from $\Xt$ to $\Yt$ if-and-only-if the inverse-image of every element of $\baseof{\Yt}$ under $\cf$
is an open set of $\Xt$. That is,
\begin{align}
&\Foreach{\cf}{\Func{\X}{\Y}}\cr
&\(\Foreach{\baseof{\Yt}}{\Cbase{\Yt}}
\left\{\[\cf\in\CF{\Xt}{\Yt}\]\thenn
\[\Foreach{\B}{\baseof{\Yt}}
\func{\pimage{\cf}}{\B}\in\topology{\X}\]\right\}\).
\end{align}
\prooff
$\cf$
is taken as an arbitrary element of $\Func{\X}{\Y}$, and $\baseof{\Yt}$
as an arbitrary element of $\Cbase{\Yt}$
(a base for the topological-space $\Yt$).
\begin{itemize}
\item[${\textbf{\textsf{p1}}}$]
It is assumed that $\cf$
is a continuous map from $\Xt$ to $\Yt$. That is,
\begin{equation}
\cf\in\CF{\Xt}{\Yt}.
\end{equation}
Then considering that,
\begin{equation}
\baseof{\Yt}\subseteq\topology{\Y},
\end{equation}
and according to \refdef{defcontinuousfunction},
\begin{equation*}
\Foreach{\B}{\baseof{\Yt}}
\func{\pimage{\cf}}{\B}\in\topology{\X}.
\end{equation*}
\endp
\end{itemize}
\begin{itemize}
\item[${\textbf{\textsf{p2}}}$]
It is assumed that,
\begin{equation}\label{thmcontiniuityequiv4p2eq1}
\Foreach{\B}{\baseof{\Yt}}
\func{\pimage{\cf}}{\B}\in\topology{\X}.
\end{equation}
\begin{itemize}
\item[${\textbf{\textsf{p2-1}}}$]
$\U$
is taken as an arbitrary element of $\topology{\Y}$.
Then according to \refdef{defbase},
\begin{equation}\label{thmcontiniuityequiv4p2-1eq1}
\Existsis{\sC}{\CSs{\baseof{\Yt}}}
\U=\union{\sC}.
\end{equation}
Hence,
\begin{align}\label{thmcontiniuityequiv4p2-1eq2}
\func{\pimage{\cf}}{\U}&=
\func{\pimage{\cf}}{\union{\sC}}\cr
&=\Union{\B}{\sC}{\func{\pimage{\cf}}{\B}}.
\end{align}
In addition,
\Ref{thmcontiniuityequiv4p2eq1}
and
\refdef{deftopologicalspace}
imply,
\begin{equation}\label{thmcontiniuityequiv4p2-1eq3}
\[\Union{\B}{\sC}{\func{\pimage{\cf}}{\B}}\]\in\topology{\X}.
\end{equation}
\Ref{thmcontiniuityequiv4p2-1eq2}
and
\Ref{thmcontiniuityequiv4p2-1eq3}
imply,
\begin{equation}
\func{\pimage{\cf}}{\U}\in\topology{\X}.
\end{equation}
\endp
\end{itemize}
Therefore,
\begin{equation}
\Foreach{\U}{\topology{\Y}}
\func{\pimage{\cf}}{\U}\in\topology{\X},
\end{equation}
which according to \refdef{defcontinuousfunction}, means,
\begin{equation*}
\cf\in\CF{\Xt}{\Yt}.
\end{equation*}
\endp
\end{itemize}
\endthm
%%%%%%%%%%%%%%%%%%%%%%%%%%%%%%%%%%%%%%%%%%%%%%%%%%%%%%%%%%%%%%%%%%%%%%%%%%%%%%%%%%%%%%%%%%%%%%%%%%
\theorem\label{thmcontiniuityequiv2}
Each $\Xt=\opair{\X}{\topology{\X}}$
and
$\Yt=\opair{\Y}{\topology{\Y}}$
is taken as a topological-space.
\begin{align}
&\Foreach{\cf}{\Func{\X}{\Y}}\cr
&\[\cf\in\CF{\Xt}{\Yt}\]\thenn
\left\{\Foreach{\Asubset{\X}}{\CSs{\X}}
\[\func{\resd{\cf}}{\Asubset{\X}}\in
\CF{\opair{\Asubset{\X}}{\stopology{\topology{\X}}{\Asubset{\X}}}}{\Yt}\]\right\}.\cr
&{}
\end{align}
\prooff
$\cf$
is taken as an arbitrary element of $\Func{\X}{\Y}$
(an arbitrary function from $\X$ to $\Y$).
\begin{itemize}
\item[${\textbf{\textsf{p1}}}$]
It is assumed that $\cf$ is a continuous map from $\Xt$ to $\Yt$.
That is,
$\cf\in\CF{\Xt}{\Yt}$.
\begin{itemize}
\item[${\textbf{\textsf{p1-1}}}$]
$\Asubset{\X}$
is taken as an arbitrary element of $\CSs{\X}$.
Then according to \refdef{defsubspacetopology1},
\begin{equation}\label{thmcontiniuityequiv2p1-1eq1}
\Foreach{\U}{\topology{\X}}
\(\Asubset{\X}\cap\U\)\in\stopology{\topology{\X}}{\Asubset{\X}}.
\end{equation}
In addition, it is known that,
\begin{equation}\label{thmcontiniuityequiv2p1-1eq2}
\Foreach{\U}{\CSs{\Y}}
\func{\pimage{\[\func{\resd{\cf}}{\Asubset{\X}}\]}}{\U}=
\[\Asubset{\X}\cap\func{\pimage{\cf}}{\U}\].
\end{equation}
\begin{itemize}
\item[${\textbf{\textsf{p1-1-1}}}$]
$\U$
is taken as an arbitrary element of $\topology{\Y}$.
Then considering that $\cf$ is a continuous map from $\Xt$ to $\Yt$,
according to \refdef{defcontinuousfunction} it is clear that,
\begin{equation}\label{thmcontiniuityequiv2p1-1-1eq1}
\func{\pimage{\cf}}{\U}\in\topology{\X},
\end{equation}
and hence according to \Ref{thmcontiniuityequiv2p1-1eq1},
\begin{equation}\label{thmcontiniuityequiv2p1-1-1eq2}
\[\Asubset{\X}\cap\func{\pimage{\cf}}{\U}\]\in
\stopology{\topology{\X}}{\Asubset{\X}}.
\end{equation}
In addition, according to \Ref{thmcontiniuityequiv2p1-1eq2},
\begin{equation}\label{thmcontiniuityequiv2p1-1-1eq3}
\func{\pimage{\[\func{\resd{\cf}}{\Asubset{\X}}\]}}{\U}=
\[\Asubset{\X}\cap\func{\pimage{\cf}}{\U}\].
\end{equation}
\Ref{thmcontiniuityequiv2p1-1-1eq2}
and
\Ref{thmcontiniuityequiv2p1-1-1eq3}
imply,
\begin{equation}
\func{\pimage{\[\func{\resd{\cf}}{\Asubset{\X}}\]}}{\U}\in
\stopology{\topology{\X}}{\Asubset{\X}}.
\end{equation}
\endp
\end{itemize}
Therefore,
\begin{equation}
\Foreach{\U}{\topology{\Y}}
\left\{
\func{\pimage{\[\func{\resd{\cf}}{\Asubset{\X}}\]}}{\U}\in
\stopology{\topology{\X}}{\Asubset{\X}}
\right\}.
\end{equation}
According to \refdef{defcontinuousfunction},
this imply,
\begin{equation}
\func{\resd{\cf}}{\Asubset{\X}}\in
\CF{\opair{\Asubset{\X}}{\stopology{\topology{\X}}{\Asubset{\X}}}}{\Yt}.
\end{equation}
\endp
\end{itemize}
\endp
\end{itemize}
\begin{itemize}
\item[${\textbf{\textsf{p2}}}$]
It is assumed that,
\begin{equation}
\Foreach{\Asubset{\X}}{\CSs{\X}}
\[\func{\resd{\cf}}{\Asubset{\X}}\in
\CF{\opair{\Asubset{\X}}{\stopology{\topology{\X}}{\Asubset{\X}}}}{\Yt}\].
\end{equation}
Then considering that,
\begin{align}
\X&\in\CSs{\X},\\
\func{\resd{\cf}}{\X}&=\cf,\\
\stopology{\topology{\X}}{\X}&=\topology{\X},
\end{align}
it is clear that,
\begin{equation}
\cf\in\CF{\Xt}{\Yt}.
\end{equation}
\endp
\end{itemize}
\endthm
%%%%%%%%%%%%%%%%%%%%%%%%%%%%%%%%%%%%%%%%%%%%%%%%%%%%%%%%%%%%%%%%%%%%%%%%%%%%%%%%%%%%%%%%%%%%%%%%%%
\theorem\label{thmcontiniuityequiv3}
Each
$\Xt=\opair{\X}{\topology{\X}}$
and
$\Yt=\opair{\Y}{\topology{\Y}}$
is taken as a topological-space.
\begin{align}
&\Foreach{\cf}{\Func{\X}{\Y}}\cr
&\[\cf\in\CF{\Xt}{\Yt}\]\thenn\left\{
\Foreach{\Asubset{\Y}}{\func{\Cinc{\Y}}{\func{\image{\cf}}{\X}}}
\[\func{\rescd{\cf}}{\Asubset{\Y}}\in\CF{\Xt}{\opair{\Asubset{\Y}}{\stopology{\topology{\Y}}{\Asubset{\Y}}}}\]\right\}.\cr
&{}
\end{align}
\prooff
$\cf$
is taken as an arbitrary element of $\Func{\X}{\Y}$
(an arbitrary function from $\X$ to $\Y$).
\begin{itemize}
\item[${\textbf{\textsf{p1}}}$]
It is assumed that $\cf$
is a continuous map from $\Xt$ to $\Yt$. That is,
$\cf\in\CF{\Xt}{\Yt}$.
\begin{itemize}
\item[${\textbf{\textsf{p1-1}}}$]
$\Asubset{\Y}$
is taken as an arbitrary element of $\func{\Cinc{\Y}}{\func{\image{\cf}}{\X}}$.
Then according to \refdef{defsubspacetopology1},
\begin{equation}\label{thmcontiniuityequiv3p1-1eq1}
\Foreach{\V}{\stopology{\topology{\Y}}{\Asubset{\Y}}}
\[\Exists{\U}{\topology{Y}}\V=\Asubset{\Y}\cap\U\].
\end{equation}
In addition, it is known that,
\begin{equation}\label{thmcontiniuityequiv3p1-1eq2}
\Foreach{\U}{\CSs{\Y}}
\func{\pimage{\[\func{\rescd{\cf}}{\Asubset{\Y}}\]}}{\Asubset{\Y}\cap\U}=
\func{\pimage{\cf}}{\U}.
\end{equation}
\begin{itemize}
\item[${\textbf{\textsf{p1-1-1}}}$]
$\V$
is taken as an arbitrary element of $\stopology{\topology{\Y}}{\Asubset{\Y}}$.
Then according to \Ref{thmcontiniuityequiv3p1-1eq1},
\begin{equation}\label{thmcontiniuityequiv3p1-1-1eq1}
\Existsis{\U}{\topology{Y}}\V=\Asubset{\Y}\cap\U.
\end{equation}
\Ref{thmcontiniuityequiv3p1-1eq2}
and
\Ref{thmcontiniuityequiv3p1-1-1eq1}
imply,
\begin{align}\label{thmcontiniuityequiv3p1-1-1eq2}
\func{\pimage{\[\func{\rescd{\cf}}{\Asubset{\Y}}\]}}{\V}&=
\func{\pimage{\[\func{\rescd{\cf}}{\Asubset{\Y}}\]}}{\Asubset{\Y}\cap\U}\cr
&=\func{\pimage{\cf}}{\U}.
\end{align}
Considering that $\cf$
is a continuous map from $\Xt$ to $\Yt$, and
$\U\in\topology{\Y}$,
according to \refdef{defcontinuousfunction},
it is clear that,
\begin{equation}\label{thmcontiniuityequiv3p1-1-1eq3}
\func{\pimage{\cf}}{\U}\in\topology{\X},
\end{equation}
\Ref{thmcontiniuityequiv3p1-1-1eq2}
and
\Ref{thmcontiniuityequiv3p1-1-1eq3}
imply,
\begin{equation}
\func{\pimage{\[\func{\rescd{\cf}}{\Asubset{\Y}}\]}}{\V}\in\topology{\X}.
\end{equation}
\endp
\end{itemize}
Therefore,
\begin{equation}
\Foreach{\V}{\stopology{\topology{\Y}}{\Asubset{\Y}}}
\left\{
\func{\pimage{\[\func{\rescd{\cf}}{\Asubset{\Y}}\]}}{\V}\in\topology{\X}
\right\},
\end{equation}
which according to \refdef{defcontinuousfunction}, means,
\begin{equation}
\func{\rescd{\cf}}{\Asubset{\Y}}\in
\CF{\Xt}{\opair{\Asubset{\Y}}{\stopology{\topology{\Y}}{\Asubset{\Y}}}}.
\end{equation}
\endp
\end{itemize}
\endp
\end{itemize}
\begin{itemize}
\item[${\textbf{\textsf{p2}}}$]
It is assumed that,
\begin{equation}
\Foreach{\Asubset{\Y}}{\func{\Cinc{\Y}}{\func{\image{\cf}}{\X}}}
\[\func{\rescd{\cf}}{\Asubset{\Y}}\in
\CF{\Xt}{\opair{\Asubset{\Y}}{\stopology{\topology{\Y}}{\Asubset{\Y}}}}\].
\end{equation}
Then considering that,
\begin{align}
\Y&\in\func{\Cinc{\Y}}{\func{\image{\cf}}{\X}},\\
\func{\rescd{\cf}}{\Y}&=\cf,\\
\stopology{\topology{\Y}}{\Y}&=\topology{\Y},
\end{align}
it is clear that,
\begin{equation}
\cf\in\CF{\Xt}{\Yt}.
\end{equation}
\endp
\end{itemize}
\endthm
%%%%%%%%%%%%%%%%%%%%%%%%%%%%%%%%%%%%%%%%%%%%%%%%%%%%%%%%%%%%%%%%%%%%%%%%%%%%%%%%%%%%%%%%%%%%%%%%%%
\corollary\label{correstrictionofcontinuousfunction}
Each
$\Xt=\opair{\X}{\topology{\X}}$
and
$\Yt=\opair{\Y}{\topology{\Y}}$
is taken as a topological-space.
\begin{align}
&\Foreach{\hf}{\CF{\Xt}{\Yt}}\cr
&\bigg[
\Foreach{\Asubset{\X}}{\CSs{\X}}\cr
&~~~\func{\rescd{\func{\resd{\hf}}{\Asubset{\X}}}}{\func{\image{\hf}}{\Asubset{\X}}}
\in\CF{\opair{\Asubset{\X}}{\stopology{\topology{\X}}{\Asubset{\X}}}}
{\opair{\func{\image{\hf}}{\Asubset{\X}}}{\stopology{\topology{\Y}}{\func{\image{\hf}}{\Asubset{\X}}}}}
\bigg].\cr
&{}
\end{align}
\endcor
%%%%%%%%%%%%%%%%%%%%%%%%%%%%%%%%%%%%%%%%%%%%%%%%%%%%%%%%%%%%%%%%%%%%%%%%%%%%%%%%%%%%%%%%%%%%%%%%%%
\theorem\label{thmcontiniuityandclosure}
Each
$\Xt=\opair{\X}{\topology{\X}}$
and
$\Yt=\opair{\Y}{\topology{\Y}}$
is taken as a topological-space.
\begin{align}
&\Foreach{\cf}{\Func{\X}{\Y}}\cr
&\[\cf\in\CF{\Xt}{\Yt}\]\thenn
\left\{
\Foreach{\Asubset{\Y}}{\CSs{\Y}}
\[\func{\Cl{\Xt}}{\func{\pimage{\cf}}{\Asubset{\Y}}}\subseteq
\func{\pimage{\cf}}{\func{\Cl{\Yt}}{\Asubset{\Y}}}\]\right\}.\cr
&{}
\end{align}
\prooff
$\cf$
is taken as an arbitrary element of
$\Func{\X}{\Y}$
(an arbitrary function from $\X$ to $\Y$).
\begin{itemize}
\item[${\textbf{\textsf{p1}}}$]
It is assumed that $\cf$
is a continuous map from $\Xt$ to $\Yt$. That is,
\begin{equation}\label{thmcontiniuityandclosurep1eq1}
\cf\in\CF{\Xt}{\Yt}.
\end{equation}
Then according to
\refthm{thmcontiniuityequiv1},
\begin{equation}\label{thmcontiniuityandclosurep1eq2}
\Foreach{\U}{\Fclosed{\Y}{\topology{\Y}}}
\func{\pimage{\cf}}{\U}\in\Fclosed{\X}{\topology{\X}}.
\end{equation}
\begin{itemize}
\item[${\textbf{\textsf{p1-1}}}$]
$\Asubset{\Y}$
is taken as an arbitrary subset of $\Y$. According to \refcor{corclosureofset0},
\begin{equation}\label{thmcontiniuityandclosurep1-1eq1}
\func{\Cl{\Yt}}{\Asubset{\Y}}\in\Fclosed{\Y}{\topology{\Y}}.
\end{equation}
\Ref{thmcontiniuityandclosurep1eq2}
and
\Ref{thmcontiniuityandclosurep1-1eq1}
imply,
\begin{equation}\label{thmcontiniuityandclosurep1-1eq2}
\func{\pimage{\cf}}{\func{\Cl{\Yt}}{\Asubset{\Y}}}\in\Fclosed{\X}{\topology{\X}}.
\end{equation}
Additionally, considering that
$\func{\Cl{\Yt}}{\Asubset{\Y}}\supseteq\Asubset{\Y}$
i is clear that,
\begin{equation}\label{thmcontiniuityandclosurep1-1eq3}
\func{\pimage{\cf}}{\func{\Cl{\Yt}}{\Asubset{\Y}}}
\supseteq
\func{\pimage{\cf}}{\Asubset{\Y}}.
\end{equation}
According to \refcor{corclosureofset0},
\Ref{thmcontiniuityandclosurep1-1eq2}
and
\Ref{thmcontiniuityandclosurep1-1eq3}
imply,
\begin{equation}
\func{\Cl{\Xt}}{\func{\pimage{\cf}}{\Asubset{\Y}}}
\subseteq
\func{\pimage{\cf}}{\func{\Cl{\Yt}}{\Asubset{\Y}}}.
\end{equation}
\endp
\end{itemize}
Thus,
\begin{equation*}
\Foreach{\Asubset{\Y}}{\CSs{\Y}}
\[\func{\Cl{\Xt}}{\func{\pimage{\cf}}{\Asubset{\Y}}}
\subseteq
\func{\pimage{\cf}}{\func{\Cl{\Yt}}{\Asubset{\Y}}}\].
\end{equation*}
\endp
\end{itemize}
\begin{itemize}
\item[${\textbf{\textsf{p2}}}$]
It is assumed that,
\begin{equation}\label{thmcontiniuityandclosurep2eq1}
\Foreach{\Asubset{\Y}}{\CSs{\Y}}
\[\func{\Cl{\Xt}}{\func{\pimage{\cf}}{\Asubset{\Y}}}
\subseteq
\func{\pimage{\cf}}{\func{\Cl{\Yt}}{\Asubset{\Y}}}\].
\end{equation}
\begin{itemize}
\item[${\textbf{\textsf{p2-1}}}$]
$\U$
is taken as an arbitrary element of $\Fclosed{\Y}{\topology{\Y}}$.
Then according to \refthm{thmclosureofclosedset},
\begin{equation}\label{thmcontiniuityandclosurep2-1eq1}
\func{\Cl{\Yt}}{\U}=\U.
\end{equation}
\Ref{thmcontiniuityandclosurep2eq1}
and
\Ref{thmcontiniuityandclosurep2-1eq1}
imply,
\begin{equation}
\func{\Cl{\Xt}}{\func{\pimage{\cf}}{\U}}
\subseteq
\func{\pimage{\cf}}{\U}.
\end{equation}
Hence considering that,
\begin{equation}
\func{\Cl{\Xt}}{\func{\pimage{\cf}}{\U}}
\supseteq
\func{\pimage{\cf}}{\U},
\end{equation}
it is evident that,
\begin{equation}
\func{\Cl{\Xt}}{\func{\pimage{\cf}}{\U}}=
\func{\pimage{\cf}}{\U},
\end{equation}
which according to \refthm{thmclosureofclosedset}, means,
\begin{equation}
\func{\pimage{\cf}}{\U}\in\Fclosed{\X}{\topology{\X}}.
\end{equation}
\endp
\end{itemize}
Therefore,
\begin{equation}
\Foreach{\U}{\Fclosed{\Y}{\topology{\Y}}}
\[\func{\pimage{\cf}}{\U}\in\Fclosed{\X}{\topology{\X}}\],
\end{equation}
which according to \refthm{thmcontiniuityequiv1}, means,
\begin{equation*}
\cf\in\CF{\Xt}{\Yt}.
\end{equation*}
\endp
\end{itemize}
\endthm
%%%%%%%%%%%%%%%%%%%%%%%%%%%%%%%%%%%%%%%%%%%%%%%%%%%%%%%%%%%%%%%%%%%%%%%%%%%%%%%%%%%%%%%%%%%%%%%%%%
\theorem\label{thmcontiniuityandinterior}
Each
$\Xt=\opair{\X}{\topology{\X}}$
and
$\Yt=\opair{\Y}{\topology{\Y}}$
is taken as a topological-space.
\begin{align}
&\Foreach{\cf}{\Func{\X}{\Y}}\cr
&\[\cf\in\CF{\Xt}{\Yt}\]\thenn
\left\{
\Foreach{\Asubset{\Y}}{\CSs{\Y}}
\[\func{\pimage{\cf}}{\func{\Int{\Yt}}{\Asubset{\Y}}}\subseteq
\func{\Int{\Xt}}{\func{\pimage{\cf}}{\Asubset{\Y}}}\]\right\}.\cr
&{}
\end{align}
\prooff
$\cf$
is taken as an arbitrary element of $\Func{\X}{\Y}$
(an arbitrary function from $\X$ to $\Y$).
\begin{itemize}
\item[${\textbf{\textsf{p1}}}$]
It is assumed that $\cf$
is a continuous map from $\Xt$ to $\Yt$.
That is,
\begin{equation}\label{thmcontiniuityandinteriorp1eq1}
\cf\in\CF{\Xt}{\Yt}.
\end{equation}
Then according to \refdef{defcontinuousfunction},
\begin{equation}\label{thmcontiniuityandinteriorp1eq2}
\Foreach{\U}{\topology{\Y}}
\func{\pimage{\cf}}{\U}\in\topology{\X}.
\end{equation}
\begin{itemize}
\item[${\textbf{\textsf{p1-1}}}$]
$\Asubset{\Y}$
is taken as an arbitrary subset of $\Y$. According to \refcor{corintofset0},
\begin{equation}\label{thmcontiniuityandinteriorp1-1eq1}
\func{\Int{\Yt}}{\Asubset{\Y}}\in\topology{\Y}.
\end{equation}
\Ref{thmcontiniuityandinteriorp1eq2}
and
\Ref{thmcontiniuityandinteriorp1-1eq1}
imply,
\begin{equation}\label{thmcontiniuityandinteriorp1-1eq2}
\func{\pimage{\cf}}{\func{\Int{\Yt}}{\Asubset{\Y}}}\in\topology{\X}.
\end{equation}
Additionally, considering that
$\func{\Int{\Yt}}{\Asubset{\Y}}\subseteq\Asubset{\Y}$,
clearly,
\begin{equation}\label{thmcontiniuityandinteriorp1-1eq3}
\func{\pimage{\cf}}{\func{\Int{\Yt}}{\Asubset{\Y}}}
\subseteq
\func{\pimage{\cf}}{\Asubset{\Y}}.
\end{equation}
As a result of \refcor{corintofset0},
\Ref{thmcontiniuityandinteriorp1-1eq2}
and
\Ref{thmcontiniuityandinteriorp1-1eq3}
imply,
\begin{equation}
\func{\pimage{\cf}}{\func{\Int{\Yt}}{\Asubset{\Y}}}
\subseteq
\func{\Int{\Xt}}{\func{\pimage{\cf}}{\Asubset{\Y}}}.
\end{equation}
\endp
\end{itemize}
Thus,
\begin{equation*}
\Foreach{\Asubset{\Y}}{\CSs{\Y}}
\[\func{\pimage{\cf}}{\func{\Int{\Yt}}{\Asubset{\Y}}}
\subseteq
\func{\Int{\Xt}}{\func{\pimage{\cf}}{\Asubset{\Y}}}\].
\end{equation*}
\endp
\end{itemize}
\begin{itemize}
\item[${\textbf{\textsf{p2}}}$]
It is assumed that,
\begin{equation}\label{thmcontiniuityandinteriorp2eq1}
\Foreach{\Asubset{\Y}}{\CSs{\Y}}
\[\func{\pimage{\cf}}{\func{\Int{\Yt}}{\Asubset{\Y}}}
\subseteq
\func{\Int{\Xt}}{\func{\pimage{\cf}}{\Asubset{\Y}}}\].
\end{equation}
\begin{itemize}
\item[${\textbf{\textsf{p2-1}}}$]
$\U$
is taken as an arbitrary element of $\topology{\Y}$.
Then according to \refthm{thmintofopenset},
\begin{equation}\label{thmcontiniuityandinteriorp2-1eq1}
\func{\Int{\Yt}}{\U}=\U.
\end{equation}
\Ref{thmcontiniuityandinteriorp2eq1}
and
\Ref{thmcontiniuityandinteriorp2-1eq1}
imply,
\begin{equation}
\func{\pimage{\cf}}{\U}\subseteq
\func{\Int{\Xt}}{\func{\pimage{\cf}}{\U}}.
\end{equation}
Hence considering that,
\begin{equation}
\func{\Int{\Xt}}{\func{\pimage{\cf}}{\U}}
\subseteq
\func{\pimage{\cf}}{\U},
\end{equation}
it is evident that,
\begin{equation}
\func{\Int{\Xt}}{\func{\pimage{\cf}}{\U}}=
\func{\pimage{\cf}}{\U},
\end{equation}
which according to \refthm{thmintofopenset}, means,
\begin{equation}
\func{\pimage{\cf}}{\U}\in\topology{\X}.
\end{equation}
\endp
\end{itemize}
Therefore,
\begin{equation}
\Foreach{\U}{\topology{\Y}}
\[\func{\pimage{\cf}}{\U}\in\topology{\X}\],
\end{equation}
which according to \refdef{defcontinuousfunction}\,
means,
\begin{equation*}
\cf\in\CF{\Xt}{\Yt}.
\end{equation*}
\endp
\end{itemize}
\endthm
%%%%%%%%%%%%%%%%%%%%%%%%%%%%%%%%%%%%%%%%%%%%%%%%%%%%%%%%%%%%%%%%%%%%%%%%%%%%%%%%%%%%%%%%%%%%%%%%%%
\theorem\label{thmcontiniuityandclosure1}
Each
$\Xt=\opair{\X}{\topology{\X}}$
and
$\Yt=\opair{\Y}{\topology{\Y}}$
is taken as a topological-space.
\begin{align}
&\Foreach{\cf}{\Func{\X}{\Y}}\cr
&\[\cf\in\CF{\Xt}{\Yt}\]\thenn
\left\{
\Foreach{\Asubset{\X}}{\CSs{\X}}
\[\func{\image{\cf}}{\func{\Cl{\Xt}}{\Asubset{\X}}}\subseteq
\func{\Cl{\Yt}}{\func{\image{\cf}}{\Asubset{\X}}}\]\right\}.\cr
&{}
\end{align}
\prooff
$\cf$
is taken as an arbitrary element of $\Func{\X}{\Y}$
(an arbitrary function from $\X$ to $\Y$).
As a result of
\refthm{thmcontiniuityandclosure},
\begin{align}
\[\cf\in\CF{\Xt}{\Yt}\]\thenn
\left\{
\Foreach{\Asubset{\Y}}{\CSs{\Y}}
\[\func{\Cl{\Xt}}{\func{\pimage{\cf}}{\Asubset{\Y}}}\subseteq
\func{\pimage{\cf}}{\func{\Cl{\Yt}}{\Asubset{\Y}}}\]\right\}.
\end{align}
Additionally, considering that,
\begin{align*}
&\Foreach{\Asubset{\X}}{\CSs{\X}}
\func{\pimage{\cf}}{\func{\image{\cf}}{\Asubset{\X}}}\supseteq\Asubset{\X},\\
&\Foreach{\Asubset{\Y}}{\CSs{\Y}}
\func{\image{\cf}}{\func{\pimage{\cf}}{\Asubset{\Y}}}\subseteq\Asubset{\Y},
\end{align*}
According to
\refthm{thmclosureofasubsetofset},
\begin{gather}
\left\{
\Foreach{\Asubset{\Y}}{\CSs{\Y}}
\[\func{\Cl{\Xt}}{\func{\pimage{\cf}}{\Asubset{\Y}}}\subseteq
\func{\pimage{\cf}}{\func{\Cl{\Yt}}{\Asubset{\Y}}}\]\right\}\cr
\vthenn\cr
\left\{
\Foreach{\Asubset{\X}}{\CSs{\X}}
\[\func{\Cl{\Xt}}{\Asubset{\X}}\subseteq
\func{\pimage{\cf}}{\func{\Cl{\Yt}}{\func{\image{\cf}}{\Asubset{\X}}}}\]\right\},
\end{gather}
and also,
\begin{gather}
\left\{
\Foreach{\Asubset{\X}}{\CSs{\X}}
\[\func{\Cl{\Xt}}{\Asubset{\X}}\subseteq
\func{\pimage{\cf}}{\func{\Cl{\Yt}}{\func{\image{\cf}}{\Asubset{\X}}}}\]\right\}\cr
\vthenn\cr
\left\{\Foreach{\Asubset{\X}}{\CSs{\X}}
\[\func{\image{\cf}}{\func{\Cl{\Xt}}{\Asubset{\X}}}\subseteq
\func{\Cl{\Yt}}{\func{\image{\cf}}{\Asubset{\X}}}\]\right\}.
\end{gather}
Therefore,
\begin{align}
\[\cf\in\CF{\Xt}{\Yt}\]\thenn
\left\{\Foreach{\Asubset{\X}}{\CSs{\X}}
\[\func{\image{\cf}}{\func{\Cl{\Xt}}{\Asubset{\X}}}\subseteq
\func{\Cl{\Yt}}{\func{\image{\cf}}{\Asubset{\X}}}\]\right\}.
\end{align}
\endthm
%%%%%%%%%%%%%%%%%%%%%%%%%%%%%%%%%%%%%%%%%%%%%%%%%%%%%%%%%%%%%%%%%%%%%%%%%%%%%%%%%%%%%%%%%%%%%%%%%%
%%%%%%%%%%%%%%%%%%%%%%%%%%%%%%%%%%%%%%%%%%%%%%%%%%%%%%%%%%%%%%%%%%%%%%%%%%%%%%%%%%%%%%%%%%%%%%%%%%
\definition\label{deflocalcontinuity}
Each
$\Xt=\opair{\X}{\topology{\X}}$
and
$\Yt=\opair{\Y}{\topology{\Y}}$
is taken as a topological-space, and $\cf$
as an element of $\Func{\X}{\Y}$. For every
$\point$ in $\X$,
it is said that
$\quotl$$\cf$ is locally-continuous at $\point$ with respect to $\opair{\topology{\X}}{\topology{\Y}}$$\quotr$
iff
\begin{equation}
\Foreach{\U}{\func{\nei{\Yt}}{\seta{\func{\cf}{\point}}}}
\[\Exists{\V}{\func{\nei{\Xt}}{\seta{\point}}}
\func{\image{\cf}}{\V}\subseteq\U\].
\end{equation}
\endef
%%%%%%%%%%%%%%%%%%%%%%%%%%%%%%%%%%%%%%%%%%%%%%%%%%%%%%%%%%%%%%%%%%%%%%%%%%%%%%%%%%%%%%%%%%%%%%%%%%
\theorem\label{thmcontinuityandlocalcontinuity}
Each
$\Xt=\opair{\X}{\topology{\X}}$
and
$\Yt=\opair{\Y}{\topology{\Y}}$
is taken as a topological-space.
For every $\cf$ in $\Func{\X}{\Y}$,
$\cf$
is a continuous map from $\Xt$ to $\Yt$
if-and-only-if for every $\point$ in $\X$,
$\cf$ is locally-continuous at $\point$
with respect to
$\opair{\topology{\X}}{\topology{\Y}}$
That is,
\begin{gather}
\Foreach{\cf}{\Func{\X}{\Y}}\cr
\[\cf\in\CF{\Xt}{\Yt}\right.\cr
\vthenn\cr
\Foreach{\point}{\X}
\left.\left\{\Foreach{\U}{\func{\nei{\Yt}}{\seta{\func{\cf}{\point}}}}
\[\Exists{\V}{\func{\nei{\Xt}}{\seta{\point}}}
\func{\image{\cf}}{\V}\subseteq\U\]\right\}\].
\end{gather}
\prooff
$\cf$
is taken as an arbitrary element of $\Func{\X}{\Y}$.
\begin{itemize}
\item[${\textbf{\textsf{p1}}}$]
It is assumed that,
\begin{equation}
\cf\in\CF{\Xt}{\Yt}.
\end{equation}
\begin{itemize}
\item[${\textbf{\textsf{p1-1}}}$]
$\point$
is taken as an arbitrary element of $\X$.
\begin{itemize}
\item[${\textbf{\textsf{p1-1-1}}}$]
$\U$
is taken as an arbitrary element of $\func{\nei{\Yt}}{\seta{\func{\cf}{\point}}}$
(a neighbourhood of $\seta{\point}$ in $\Yt$). Then according to \refdef{defnbdclassofsets},
\begin{align}
\U&\in\topology{\Y},\label{thmcontinuityandlocalcontinuityp1-1-1eq1}\\
\func{\cf}{\point}&\in\U.\label{thmcontinuityandlocalcontinuityp1-1-1eq2}
\end{align}
As a result of \refdef{defcontinuousfunction},
\Ref{thmcontinuityandlocalcontinuityp1-1-1eq1}
implies,
\begin{equation}\label{thmcontinuityandlocalcontinuityp1-1-1eq3}
\func{\pimage{\cf}}{\U}\in\topology{\X}.
\end{equation}
\Ref{thmcontinuityandlocalcontinuityp1-1-1eq2}
implies,
\begin{equation}\label{thmcontinuityandlocalcontinuityp1-1-1eq4}
\point\in\func{\pimage{\cf}}{\U}.
\end{equation}
According to \refdef{defnbdclassofsets},
\Ref{thmcontinuityandlocalcontinuityp1-1-1eq3}
and
\Ref{thmcontinuityandlocalcontinuityp1-1-1eq4}
imply,
\begin{equation}
\func{\pimage{\cf}}{\U}\in\func{\nei{\Xt}}{\seta{\point}}.
\end{equation}
Additionally, considering that,
\begin{equation*}
\Foreach{\Asubset{\Y}}{\CSs{\Y}}
\func{\image{\cf}}{\func{\pimage{\cf}}{\Asubset{\Y}}}\subseteq\Asubset{\Y},
\end{equation*}
it is evident that,
\begin{equation}
\func{\image{\cf}}{\func{\pimage{\cf}}{\U}}\subseteq\U.
\end{equation}
\endp
\end{itemize}
Therefore it is seen that,
\begin{equation}
\Foreach{\U}{\func{\nei{\Yt}}{\seta{\func{\cf}{\point}}}}
\[\Exists{\V}{\func{\nei{\Xt}}{\seta{\point}}}
\func{\image{\cf}}{\V}\subseteq\U\].
\end{equation}
\endp
\end{itemize}
Therefore,
\begin{equation*}
\Foreach{\point}{\X}
\left\{\Foreach{\U}{\func{\nei{\Yt}}{\seta{\func{\cf}{\point}}}}
\[\Exists{\V}{\func{\nei{\Xt}}{\seta{\point}}}
\func{\image{\cf}}{\V}\subseteq\U\]\right\}.
\end{equation*}
\endp
\end{itemize}
\begin{itemize}
\item[${\textbf{\textsf{p2}}}$]
It is assumed that,
\begin{align}\label{thmcontinuityandlocalcontinuityp2eq1}
\Foreach{\point}{\X}
\left\{\Foreach{\U}{\func{\nei{\Yt}}{\seta{\func{\cf}{\point}}}}
\[\Exists{\V}{\func{\nei{\Xt}}{\seta{\point}}}
\func{\image{\cf}}{\V}\subseteq\U\]\right\}.
\end{align}
\begin{itemize}
\item[${\textbf{\textsf{p2-1}}}$]
$\U$
is taken as an arbitrary element of $\topology{\Y}$
Then,
\begin{equation}\label{thmcontinuityandlocalcontinuityp2-1eq1}
\Foreach{\point}{\func{\pimage{\cf}}{\U}}
\func{\cf}{\point}\in\U,
\end{equation}
and
\begin{equation}\label{thmcontinuityandlocalcontinuityp2-1eq2}
\Foreach{\x}{\U}\U\in\func{\nei{\Yt}}{\seta{\x}}.
\end{equation}
Thus,
\begin{equation}\label{thmcontinuityandlocalcontinuityp2-1eq3}
\Foreach{\point}{\func{\pimage{\cf}}{\U}}
\U\in\func{\nei{\Yt}}{\seta{\func{\cf}{\point}}}.
\end{equation}
\Ref{thmcontinuityandlocalcontinuityp2eq1}
and
\Ref{thmcontinuityandlocalcontinuityp2-1eq3}
imply,
\begin{equation}
\Foreach{\point}{\func{\pimage{\cf}}{\U}}
\[\Exists{\V}{\func{\nei{\Xt}}{\seta{\point}}}
\func{\image{\cf}}{\V}\subseteq\U\].
\end{equation}
As a result of \refthm{thmopensetpointstopology},
this implies,
\begin{equation}
\func{\pimage{\cf}}{\U}\in\topology{\X}.
\end{equation}
\endp
\end{itemize}
Therefore,
\begin{equation}
\Foreach{\U}{\topology{\Y}}
\func{\pimage{\cf}}{\U}\in\topology{\X},
\end{equation}
which according to \refdef{defcontinuousfunction}, means,
\begin{equation*}
\cf\in\CF{\Xt}{\Yt}.
\end{equation*}
\endp
\end{itemize}
\endthm
%%%%%%%%%%%%%%%%%%%%%%%%%%%%%%%%%%%%%%%%%%%%%%%%%%%%%%%%%%%%%%%%%%%%%%%%%%%%%%%%%%%%%%%%%%%%%%%%%%
\theorem\label{thmcompositionofcontinuousfunctions}
Each
$\Xt_{1}=\opair{\X_1}{\topology{\X_1}}$,
$\Xt_{2}=\opair{\X_2}{\topology{\X_2}}$,
and
$\Xt_{3}=\opair{\X_3}{\topology{\X_3}}$
is taken as a topological-space.
\begin{equation}
\Foreach{\opair{\cf}{\cg}}{\[\CF{\Xt_1}{\Xt_2}\times\CF{\Xt_2}{\Xt_3}\]}
\[\cmp{\cg}{\cf}\in\CF{\Xt_1}{\Xt_3}\].
\end{equation}
\prooff
$\cf$
is taken as a continuous map from
$\Xt_{1}$ to $\Xt_{2}$
(an arbitrary element of $\CF{\Xt_{1}}{\Xt_{2}}$), and
$\cg$
as a continuous map from $\Xt_{2}$ to $\Xt_{3}$
(an arbitrary element of $\CF{\Xt_{2}}{\Xt_{3}}$).
Then according to
\refdef{defcontinuousfunction},
\begin{align}
&\Foreach{\U}{\topology{\X_2}}
\func{\pimage{\cf}}{\U}\in\topology{\X_1},\label{thmcompositionofcontinuousfunctionspeq1}\\
&\Foreach{\U}{\topology{\X_3}}
\func{\pimage{\cg}}{\U}\in\topology{\X_2}.\label{thmcompositionofcontinuousfunctionspeq2}
\end{align}
In addition,
\begin{equation}\label{thmcompositionofcontinuousfunctionspeq3}
\(\cmp{\cg}{\cf}\)\in\Func{\X_1}{\X_3},
\end{equation}
and
\begin{equation}\label{thmcompositionofcontinuousfunctionspeq4}
\Foreach{\asubset}{\CSs{\X_3}}
\func{\pimage{\(\cmp{\cg}{\cf}\)}}{\asubset}=
\func{\pimage{\cf}}{\func{\pimage{\cg}}{\asubset}}.
\end{equation}
Therefore considering that
$\topology{\X_3}\subseteq\CSs{\X_3}$,
\begin{equation}\label{thmcompositionofcontinuousfunctionspeq5}
\Foreach{\U}{\topology{\X_3}}
\func{\pimage{\(\cmp{\cg}{\cf}\)}}{\U}
\in\topology{\X_1}.
\end{equation}
As a result of \refdef{defcontinuousfunction},
\Ref{thmcompositionofcontinuousfunctionspeq3}
and
\Ref{thmcompositionofcontinuousfunctionspeq5}
imply,
\begin{equation}
\(\cmp{\cg}{\cf}\)\in\CF{\Xt_{1}}{\Xt_{3}}.
\end{equation}
\endthm
%%%%%%%%%%%%%%%%%%%%%%%%%%%%%%%%%%%%%%%%%%%%%%%%%%%%%%%%%%%%%%%%%%%%%%%%%%%%%%%%%%%%%%%%%%%%%%%%%%
\theorem\label{thmconstantfunctioniscontinuous}
Each
$\Xt=\opair{\X}{\topology{\X}}$
and
$\Yt=\opair{\Y}{\topology{\Y}}$
is taken as a topological-space.
For every $\cf$ in $\Func{\X}{\Y}$, if
$\cf$
is a constant function from $\X$ to $\Y$, then $\cf$
is a continuous map from $\Xt$ to $\Yt$. That is,
\begin{equation}
\Foreach{\cf}{\Func{\X}{\Y}}
\(\left\{\Existsu{\const}{\Y}
\[\Foreach{\point}{\X}\func{\cf}{\point}=\const\]\right\}\then
\[\cf\in\CF{\Xt}{\Yt}\]\).
\end{equation}
\prooff
$\cf$
is taken as an arbitrary element of $\Func{\X}{\Y}$, and it is assumed that,
\begin{equation}
\Existsuis{\const}{\Y}
\[\Foreach{\point}{\X}\func{\cf}{\point}=\const\].
\end{equation}
\begin{itemize}
\item[${\textbf{\textsf{p1}}}$]
$\U$
is taken as an arbitrary element of $\topology{\Y}$.
It is clear that,
\begin{equation}
\OR{\(\const\in\U\)}{\(\const\notin\U\)},
\end{equation}
and
\begin{align}
&\(\const\in\U\)\then\func{\pimage{\cf}}{\U}=\X,\\
&\(\const\notin\U\)\then\func{\pimage{\cf}}{\U}=\empty.
\end{align}
Therefore considering that,
\begin{equation}
\opair{\X}{\empty}\in\(\topology{\X}\times\topology{\X}\),
\end{equation}
it is evident that,
\begin{equation}
\func{\pimage{\cf}}{\U}\in\topology{\X}.
\end{equation}
\endp
\end{itemize}
Hence,
\begin{equation}
\Foreach{\U}{\topology{\Y}}\func{\pimage{\cf}}{\U}\in\topology{\X},
\end{equation}
which according to\refdef{defcontinuousfunction}, means,
\begin{equation}
\cf\in\CF{\Xt}{\Yt}.
\end{equation}
\endthm
%%%%%%%%%%%%%%%%%%%%%%%%%%%%%%%%%%%%%%%%%%%%%%%%%%%%%%%%%%%%%%%%%%%%%%%%%%%%%%%%%%%%%%%%%%%%%%%%%%
\theorem\label{thmcontinuityofidentityfunction}
$\X$
is taken as a set, and each
$\topology{1}$
and
$\topology{2}$
as a topology on $\X$
(an element of $\Ctops{\X}$).
The identity-function on $\X$
is a continuous map from $\opair{\X}{\topology{1}}$
to
$\opair{\X}{\topology{2}}$ if-and-only-if
$\topology{1}$ is finer than $\topology{2}$. That is,
\begin{equation}\label{thmcontinuityofidentityfunctioneq1}
\[\idf{\X}\in\CF{\opair{\X}{\topology{1}}}{\opair{\X}{\topology{2}}}\]
\thenn
\(\topology{2}\subseteq\topology{1}\).
\end{equation}
\prooff
According to \refdef{defcontinuousfunction},
\begin{equation}
\[\idf{\X}\in\CF{\opair{\X}{\topology{1}}}{\opair{\X}{\topology{2}}}\]
\thenn
\[\Foreach{\U}{\topology{2}}\func{\pimage{\idf{\X}}}{\U}\in\topology{1}\].
\end{equation}
Additionally, considering that,
\begin{equation}
\Foreach{\asubset}{\CSs{\X}}
\func{\pimage{\idf{\X}}}{\asubset}=\asubset,
\end{equation}
it is evident that,
\begin{equation}
\[\Foreach{\U}{\topology{2}}\func{\pimage{\idf{\X}}}{\U}\in\topology{1}\]
\thenn
\(\topology{2}\subseteq\topology{1}\).
\end{equation}
Therefore
\Ref{thmcontinuityofidentityfunctioneq1}
is clearly obtained.
\endthm
%%%%%%%%%%%%%%%%%%%%%%%%%%%%%%%%%%%%%%%%%%%%%%%%%%%%%%%%%%%%%%%%%%%%%%%%%%%%%%%%%%%%%%%%%%%%%%%%%%
\theorem\label{thmtrivialcontinuousmapofaspace}
$\Xt=\opair{\X}{\topology{}}$
is taken as a topological-space.
The identity-function on $\X$
is a continuous map from $\Xt$ to $\Xt$. That is,
\begin{equation}
\idf{\X}\in\CF{\Xt}{\Xt}.
\end{equation}
\prooff
According to
\refthm{thmcontinuityofidentityfunction},
it is trivial.
\endthm
%%%%%%%%%%%%%%%%%%%%%%%%%%%%%%%%%%%%%%%%%%%%%%%%%%%%%%%%%%%%%%%%%%%%%%%%%%%%%%%%%%%%%%%%%%%%%%%%%%
\theorem\label{thmcontinuityofinclusionfunction}
$\Xt=\opair{\X}{\topology{}}$
is taken as a topological-space.
For every $\asubset$ in $\CSs{\X}$,
the inclusion-function of $\asubset$ into $\X$
is a continuous map from $\opair{\asubset}{\stopology{\topology{}}{\asubset}}$ to $\Xt$. That is,
\begin{equation}
\Foreach{\asubset}{\CSs{\X}}
\(\incf{\asubset}{\X}\)\in\CF{\opair{\asubset}{\stopology{\topology{}}{\asubset}}}{\Xt}.
\end{equation}
\prooff
$\asubset$
is taken as an arbitrary subset of $\X$. Considering that,
\begin{equation}
\(\incf{\asubset}{\X}\)=\func{\resd{\idf{\X}}}{\asubset},
\end{equation}
According to \refthm{thmcontiniuityequiv2} and
\refthm{thmtrivialcontinuousmapofaspace},
it is evident that,
\begin{equation}
\(\incf{\asubset}{\X}\)\in
\CF{\opair{\asubset}{\stopology{\topology{}}{\asubset}}}{\Xt}.
\end{equation}
\endthm
%%%%%%%%%%%%%%%%%%%%%%%%%%%%%%%%%%%%%%%%%%%%%%%%%%%%%%%%%%%%%%%%%%%%%%%%%%%%%%%%%%%%%%%%%%%%%%%%%%
\theorem\label{thmcoarsesttopologyofcontinuousinclusion}
$\Xt=\opair{\X}{\topology{}}$
is taken as a topological-space, and
$\asubset$
as a subset of $\X$.
\begin{align}
&\stopology{\topology{}}{\asubset}\in
\defset{\p{\atopology}}{\Ctops{\asubset}}
{\(\incf{\asubset}{\X}\)\in\CF{\opair{\asubset}{\p{\atopology}}}{\Xt}},
\label{thmcoarsesttopologyofcontinuousinclusioneq1}\\
&\Foreach{\atopology}{\defset{\p{\atopology}}{\Ctops{\asubset}}
{\(\incf{\asubset}{\X}\)\in\CF{\opair{\asubset}{\p{\atopology}}}{\Xt}}}
\(\stopology{\topology{}}{\asubset}\subseteq\atopology\).
\label{thmcoarsesttopologyofcontinuousinclusioneq2}
\end{align}
\prooff
According to
\refthm{thmcontinuityofinclusionfunction},
\Ref{thmcoarsesttopologyofcontinuousinclusioneq1}
is clear.
\begin{itemize}
\item[${\textbf{\textsf{p1}}}$]
$\atopology$
is taken as an arbitrary element of
$\defset{\p{\atopology}}{\Ctops{\asubset}}
{\(\incf{\asubset}{\X}\)\in\CF{\opair{\asubset}{\p{\atopology}}}{\Xt}}$. This means,
$\opair{\asubset}{\atopology}$
is a topological-space, and,
\begin{equation}
\(\incf{\asubset}{\X}\)\in\CF{\opair{\asubset}{\atopology}}{\Xt},
\end{equation}
and hence according to
\refdef{defcontinuousfunction},
\begin{equation}
\Foreach{\U}{\topology{}}
\func{\pimage{\(\incf{\asubset}{\X}\)}}{\U}\in\atopology.
\end{equation}
Hence considering that,
\begin{equation}
\Foreach{\U}{\CSs{\X}}
\func{\pimage{\(\incf{\asubset}{\X}\)}}{\U}=\asubset\cap\U,
\end{equation}
it is evident that,
\begin{equation}
\Foreach{\U}{\topology{}}
\(\asubset\cap\U\)\in\atopology.
\end{equation}
Hence considering that,
\begin{equation}
\stopology{\topology{}}{\asubset}=
\defset{\V}{\CSs{\asubset}}{\[\Exists{\U}{\topology{}}\V=\asubset\cap\U\]},
\end{equation}
clearly,
\begin{equation}
\stopology{\topology{}}{\asubset}\subseteq\atopology.
\end{equation}
\endp
\end{itemize}
\endthm
%%%%%%%%%%%%%%%%%%%%%%%%%%%%%%%%%%%%%%%%%%%%%%%%%%%%%%%%%%%%%%%%%%%%%%%%%%%%%%%%%%%%%%%%%%%%%%%%%%
\theorem\label{thmcontinuitywrtfinerandcoarsertopologies}
Each
$\X$
and
$\Y$
is taken as a set, each
$\topology{\X}$
and
$\p{\topology{\X}}$
as a topology on $\X$, and each $\topology{\Y}$ and $\p{\topology{\Y}}$
as a topology on $\Y$.
If
$\p{\topology{\X}}$
is finer than $\topology{\X}$, and $\p{\topology{\Y}}$
is coarser than $\topology{\Y}$, the every continuous map from $\opair{\X}{\topology{\X}}$ to
$\opair{\X}{\topology{\Y}}$
is a continuous map from $\opair{\X}{\p{\topology{\X}}}$ to $\opair{\Y}{\p{\topology{\Y}}}$.
That is,
\begin{align}
\[\AND{\(\p{\topology{\X}}\supseteq\topology{\X}\)}{\(\p{\topology{\Y}}\subseteq\topology{\Y}\)}\]
\then
\[\CF{\opair{\X}{\topology{\X}}}{\opair{\X}{\topology{\Y}}}\subseteq
\CF{\opair{\X}{\p{\topology{\X}}}}{\opair{\Y}{\p{\topology{\Y}}}}\].
\end{align}
\prooff
It is assumed that,
\begin{align}
\p{\topology{\X}}&\supseteq\topology{\X}
\label{thmcontinuitywrtfinerandcoarsertopologiespeq1}\\
\p{\topology{\Y}}&\subseteq\topology{\Y}.
\label{thmcontinuitywrtfinerandcoarsertopologiespeq2}
\end{align}
\begin{itemize}
\item[${\textbf{\textsf{p1}}}$]
$\cf$
is taken as an arbitrary element of $\CF{\opair{\X}{\topology{\X}}}{\opair{\Y}{\topology{\Y}}}$.
Then according to \refdef{defcontinuousfunction},
\begin{equation}\label{thmcontinuitywrtfinerandcoarsertopologiesp1eq1}
\Foreach{\U}{\topology{\Y}}
\func{\pimage{\cf}}{\U}\in\topology{\X}.
\end{equation}
\Ref{thmcontinuitywrtfinerandcoarsertopologiespeq1},
\Ref{thmcontinuitywrtfinerandcoarsertopologiespeq2},
and
\Ref{thmcontinuitywrtfinerandcoarsertopologiesp1eq1}
imply,
\begin{equation}
\Foreach{\U}{\p{\topology{\Y}}}
\func{\pimage{\cf}}{\U}\in\p{\topology{\X}},
\end{equation}
which according to \refdef{defcontinuousfunction},
means,
\begin{equation}
\cf\in\CF{\opair{\X}{\p{\topology{\X}}}}{\opair{\Y}{\p{\topology{\Y}}}}.
\end{equation}
\endp
\end{itemize}
\endthm
%%%%%%%%%%%%%%%%%%%%%%%%%%%%%%%%%%%%%%%%%%%%%%%%%%%%%%%%%%%%%%%%%%%%%%%%%%%%%%%%%%%%%%%%%%%%%%%%%%
\theorem\label{thmcontinuousimageofdensesetisdense}
Each
$\Xt=\opair{\X}{\topology{\X}}$,
and
$\Yt=\opair{\Y}{\topology{\Y}}$
is taken as a topological-space.
For every
$\cf$
in
$\Func{\X}{\Y}$,
if
$\cf$
is a surjective continuous map from $\Xt$ to $\Yt$, then the image of
every dense set of $\Xt$ under $\cf$ is a dense set of $\Yt$. That is,
\begin{align}
&\Foreach{\cf}{\Func{\X}{\Y}}\cr
&\left\{
\[\AND{\cf\in\CF{\Xt}{\Yt}}{\func{\image{\cf}}{\X}=\Y}\]\then
\[\Foreach{\Asubset{\X}}{\Fdense{\X}{\topology{\X}}}
\func{\image{\cf}}{\Asubset{\X}}\in\Fdense{\Y}{\topology{\Y}}\]\right\}.\cr
&{}
\end{align}
\prooff
$\cf$
is taken as a function from $\X$ to $\Y$. It is assumed that,
\begin{align}
&\cf\in\CF{\Xt}{\Yt},\label{thmcontinuousimageofdensesetisdensepeq1}\\
&\func{\image{\cf}}{\X}=\Y.\label{thmcontinuousimageofdensesetisdensepeq2}
\end{align}
Then according to \refthm{thmcontiniuityandclosure1},
\begin{equation}\label{thmcontinuousimageofdensesetisdensepeq3}
\Foreach{\asubset}{\CSs{\X}}
\[\func{\image{\cf}}{\func{\Cl{\Xt}}{\asubset}}\subseteq
\func{\Cl{\Yt}}{\func{\image{\cf}}{\asubset}}\].
\end{equation}
\begin{itemize}
\item[${\textbf{\textsf{p1}}}$]
$\Asubset{\X}$
is taken as an arbitrary element of $\Fdense{\X}{\topology{\X}}$.
Then according to \refdef{defclassofdensesets},
\begin{equation}
\func{\Cl{\Xt}}{\Asubset{\X}}=\X.
\end{equation}
\Ref{thmcontinuousimageofdensesetisdensepeq1},
\Ref{thmcontinuousimageofdensesetisdensepeq2},
and
\Ref{thmcontinuousimageofdensesetisdensepeq3}
imply,
\begin{equation}
\Y\subseteq\func{\Cl{\Yt}}{\func{\image{\cf}}{\Asubset{\X}}},
\end{equation}
which means,
\begin{equation}
\func{\Cl{\Yt}}{\func{\image{\cf}}{\Asubset{\X}}}=\Y,
\end{equation}
and hence according to
\refdef{defclassofdensesets},
\begin{equation}
\func{\image{\cf}}{\Asubset{\X}}\in\Fdense{\Y}{\topology{\Y}}.
\end{equation}
\endp
\end{itemize}
\endthm
%%%%%%%%%%%%%%%%%%%%%%%%%%%%%%%%%%%%%%%%%%%%%%%%%%%%%%%%%%%%%%%%%%%%%%%%%%%%%%%%%%%%%%%%%%%%%%%%%%
\theorem\label{thmfunctionwithdiscretedomainiscontinuous}
$\X$
is taken as a set, and
$\Yt=\opair{\Y}{\topology{\Y}}$
as a topological-space.
Every function from $\X$ to $\Y$
is a continuous map from the discrete topological-space $\opair{\X}{\CSs{\X}}$ to $\Yt$. That is,
\begin{equation}
\Func{\X}{\Y}\subseteq
\CF{\opair{\X}{\CSs{\X}}}{\Yt}.
\end{equation}
\prooff
$\cf$
is taken as an arbitrary element of $\Func{\X}{\Y}$. It is evident that,
\begin{equation}
\Foreach{\U}{\CSs{\Y}}
\func{\pimage{\cf}}{\U}\in\CSs{\X},
\end{equation}
and hence more specifically,
\begin{equation}
\Foreach{\U}{\topology{\Y}}
\func{\pimage{\cf}}{\U}\in\CSs{\X}.
\end{equation}
Thus according to \refdef{defcontinuousfunction},
\begin{equation}
\cf\in\CF{\opair{\X}{\CSs{\X}}}{\Yt}.
\end{equation}
\endthm
%%%%%%%%%%%%%%%%%%%%%%%%%%%%%%%%%%%%%%%%%%%%%%%%%%%%%%%%%%%%%%%%%%%%%%%%%%%%%%%%%%%%%%%%%%%%%%%%%%
\theorem\label{thmfunctionwithindescretecodomainiscontinuos}
$\Xt=\opair{\X}{\topology{\X}}$
is taken as a topological-space, and $\Y$ as a set.
Every function from $\X$ to $\Y$
is a continuous map from $\Xt$
to the indiscrete topological-space $\opair{\Y}{\seta{\binary{\empty}{\Y}}}$. That is,
\begin{equation}
\Func{\X}{\Y}\subseteq
\CF{\Xt}{\opair{\Y}{\seta{\binary{\empty}{\Y}}}}.
\end{equation}
\prooff
$\cf$
is taken as an arbitrary element of $\Func{\X}{\Y}$. Considering that,
\begin{align}
\func{\pimage{\cf}}{\Y}&=\X\cr
&\in\topology{\X},
\end{align}
and
\begin{align}
\func{\pimage{\cf}}{\empty}&=\empty
\cr\in\topology{\X},
\end{align}
it is evident that,
\begin{equation}
\Foreach{\U}{\seta{\binary{\empty}{\Y}}}
\func{\pimage{\cf}}{\U}\in\topology{\X},
\end{equation}
and hence according to
\refdef{defcontinuousfunction},
\begin{equation}
\cf\in\CF{\Xt}{\opair{\Y}{\seta{\binary{\empty}{\Y}}}}.
\end{equation}
\endthm
%%%%%%%%%%%%%%%%%%%%%%%%%%%%%%%%%%%%%%%%%%%%%%%%%%%%%%%%%%%%%%%%%%%%%%%%%%%%%%%%%%%%%%%%%%%%%%%%%%%%%%%%%%%%%%%%%%%%%
\theorem
Each
$\Xt_{1}=\opair{\X_1}{\topology{1}}$, $\Xt_{2}=\opair{\X_2}{\topology{2}}$, and
$\Xt_{3}=\opair{\X_3}{\topology{3}}$ is taken as a topological space.
\begin{itemize}
\item
The projection maps $\function{\Proj{1}{}}{\Cprod{\X_1}{\X_2}}{\X_1}$ and
$\function{\Proj{2}{}}{\Cprod{\X_1}{\X_2}}{\X_2}$ are continuous maps from
$\topprod{\Xt_1}{\Xt_2}$ to $\Xt_1$ and $\Xt_2$, respectibely, and
the topology $\producttop{\seta{\binary{\Xt_1}{\Xt_2}}}$ is the coarsest topology on
$\Cprod{\X_1}{\X_2}$ that make these projection maps continuous.
\item
For every function $\function{\cf}{\X_3}{\Cprod{\X_1}{\X_2}}$, $\cf$ is a continuous map
from $\Xt_3$ to $\topprod{\Xt_1}{\Xt_2}$ if and only if
$\cmp{\Proj{1}{}}{f}\in\CF{\Xt_3}{\Xt_1}$ and $\cmp{\Proj{2}{}}{f}\in\CF{\Xt_3}{\Xt_2}$.
\end{itemize}
\proof
It is left as an exercise.
\endthm
%%%%%%%%%%%%%%%%%%%%%%%%%%%%%%%%%%%%%%%%%%%%%%%%%%%%%%%%%%%%%%%%%%%%%%%%%%%%%%%%%%%%%%%%%%%%%%%%%%%%%%%%%%%%%%%%%%%%%
\theorem
$\Xt=\opair{\X}{\topology{\X}}$ is taken as a topological space, and $\eqrel{}$ as an equivalence relation on $\X$.
The canonical projection $\function{\pi}{\X}{\EqClass{\X}{\eqrel{}}}$ is a continuous map from $\Xt$
to the quotient space $\topq{\Xt}{\eqrel{}}$, and $\quotienttop{\Xt}{\eqrel{}}$ is the finest topology on
$\EqClass{\X}{\eqrel{}}$ that makes this map continuous.
\proof
It is left as an exercise.
\endthm
%%%%%%%%%%%%%%%%%%%%%%%%%%%%%%%%%%%%%%%%%%%%%%%%%%%%%%%%%%%%%%%%%%%%%%%%%%%%%%%%%%%%%%%%%%%%%%%%%%
%%%%%%%%%%%%%%%%%%%%%%%%%%%%%%%%%%%%%%%%%%%%%%%%%%%%%%%%%%%%%%%%%%%%%%%%%%%%%%%%%%%%%%%%%%%%%%%%%%
%%%%%%%%%%%%%%%%%%%%%%%%%%%%%%%%%%%%%%%%%%%%%%%%%%%%%%%%%%%%%%%%%%%%%%%%%%%%%%%%%%%%%%%%%%%%%%%%%%
%%%%%%%%%%%%%%%%%%%%%%%%%%%%%%%%%%%%%%%%%%%%%%%%%%%%%%%%%%%%%%%%%%%%%%%%%%%%%%%%%%%%%%%%%%%%%%%%%%
%%%%%%%%%%%%%%%%%%%%%%%%%%%%%%%%%%%%%%%%%%%%%%%%%%%%%%%%%%%%%%%%%%%%%%%%%%%%%%%%%%%%%%%%%%%%%%%%%%
%%%%%%%%%%%%%%%%%%%%%%%%%%%%%%%%%%%%%%%%%%%%%%%%%%%%%%%%%%%%%%%%%%%%%%%%%%%%%%%%%%%%%%%%%%%%%%%%%%
%%%%%%%%%%%%%%%%%%%%%%%%%%%%%%%%%%%%%%%%%%%%%%%%%%%%%%%%%%%%%%%%%%%%%%%%%%%%%%%%%%%%%%%%%%%%%%%%%%
%%%%%%%%%%%%%%%%%%%%%%%%%%%%%%%%%%%%%%%%%%%%%%%%%%%%%%%%%%%%%%%%%%%%%%%%%%%%%%%%%%%%%%%%%%%%%%%%%%
%%%%%%%%%%%%%%%%%%%%%%%%%%%%%%%%%%%%%%%%%%%%%%%%%%%%%%%%%%%%%%%%%%%%%%%%%%%%%%%%%%%%%%%%%%%%%%%%%%
%%%%%%%%%%%%%%%%%%%%%%%%%%%%%%%%%%%%%%%%%%%%%%%%%%%%%%%%%%%%%%%%%%%%%%%%%%%%%%%%%%%%%%%%%%%%%%%%%%
%%%%%%%%%%%%%%%%%%%%%%%%%%%%%%%%%%%%%%%%%%%%%%%%%%%%%%%%%%%%%%%%%%%%%%%%%%%%%%%%%%%%%%%%%%%%%%%%%%
\section{
Coverings of a Topological Space
}
\definition\label{defcovermapofset}
$\X$
is taken as a set.
The mapping $\fcover{\X}$ is defined as the following.
\begin{itemize}
\item[${\textbf{\textsf{COVM1}}}$]
%\hfill
$\fcover{\X}\indef\Func{\CSs{\X}}{\CSs{\CSs{\CSs{\X}}}}$.
\item[${\textbf{\textsf{COVM2}}}$]
%\hfill
$\Foreach{\Y}{\CSs{\X}}
\func{\fcover{\X}}{\Y}\eqdef\defset{\acover}{\CSs{\CSs{\X}}}{\(\union{\acover}\)\supseteq\Y}$.
\end{itemize}
\begin{itemize}
\item
$\fcover{\X}$
is referred to as the $\quotl$covering transformation of $\X$$\quotr$.
\item
For every subset $\Y$ of $\X$,
every element of $\func{\fcover{\X}}{\Y}$ is referred to as a
$\quotl$covering of $\Y$ in $\X$$\quotr$.
\item
$\func{\fcover{\X}}{\X}$
is also denoted by $\covers{\X}$.
Additionally, for every element $\acover$ of $\covers{\X}$,
$\opair{\X}{\acover}$
is also referred to as a $\quotl$covering of $\X$$\quotr$.
\item
Every finite element of $\covers{\X}$
is referred to as a $\quotl$finite covering of $\X$$\quotr$.
\end{itemize}
\endef
%%%%%%%%%%%%%%%%%%%%%%%%%%%%%%%%%%%%%%%%%%%%%%%%%%%%%%%%%%%%%%%%%%%%%%%%%%%%%%%%%%%%%%%%%%%%%%%%%%
\definition\label{defsubcover}
$\X$
is taken as a set,
$\acover$
an element of
$\covers{\X}$
(a covering of $\X$), and
$\p{\acover}$
an element of $\CSs{\CSs{\X}}$
(a cllection of subsets of $\X$).
$\opair{\X}{\p{\acover}}$
is referred to as a $\quotl$subcover of $\opair{\X}{\acover}$$\quotr$ iff,
\begin{align}
&\p{\acover}\in\covers{\X},\\
&\p{\acover}\subseteq\acover.
\end{align}
\endef
%%%%%%%%%%%%%%%%%%%%%%%%%%%%%%%%%%%%%%%%%%%%%%%%%%%%%%%%%%%%%%%%%%%%%%%%%%%%%%%%%%%%%%%%%%%%%%%%%%
\definition\label{defrefinedcover}
$\X$
is taken as a set,
$\acover$
as an element of
$\covers{\X}$
(a cover of $\X$), and
$\p{\acover}$
an element of $\CSs{\CSs{\X}}$
(a collection of subsets of $\X$).
$\opair{\X}{\p{\acover}}$
is referred to as a $\quotl$refinement of $\opair{\X}{\acover}$$\quotr$ iff
\begin{align}
&\p{\acover}\in\covers{\X},\\
&\Foreach{\p{\asubset}}{\p{\acover}}
\[\Exists{\asubset}{\acover}\p{\asubset}\subseteq\asubset\].
\end{align}
\endef
%%%%%%%%%%%%%%%%%%%%%%%%%%%%%%%%%%%%%%%%%%%%%%%%%%%%%%%%%%%%%%%%%%%%%%%%%%%%%%%%%%%%%%%%%%%%%%%%%%
%%%%%%%%%%%%%%%%%%%%%%%%%%%%%%%%%%%%%%%%%%%%%%%%%%%%%%%%%%%%%%%%%%%%%%%%%%%%%%%%%%%%%%%%%%%%%%%%%%
\subsection{
Fundamental Covers of a Topological Space}
\definition\label{deffundamentalcover}
$\Xt=\opair{\X}{\topology{}}$
is taken as a topological-space.
For every $\acover$ in $\CSs{\CSs{\X}}$,
$\acover$
is referred to as a $\quotl$fundamental cover of $\Xt$$\quotr$ iff these properties are hold.
\begin{itemize}
\item[${\textbf{\textsf{FCOV1}}}$]
$\acover$
is a cover of $\X$.
%\hfill
$\acover\in\covers{\X}$.
\item[${\textbf{\textsf{FCOV2}}}$]
%\hfill
$\defset{\U}{\CSs{\X}}{\[\Foreach{\covelm}{\acover}\covelm\cap\U\in\stopology{\topology{}}{\covelm}\]}
\subseteq\topology{}$.
\end{itemize}
The set of all fundamental covers of $\Xt$ is denoted by $\Fcov{\Xt}$. That is,
\begin{align}
&\Fcov{\Xt}:=\cr
&\defset{\acover}{\covers{\X}}
{\(\defset{\U}{\CSs{\X}}{\[\Foreach{\covelm}{\acover}\covelm\cap\U\in\stopology{\topology{}}{\covelm}\]}
\subseteq\topology{}\)}.
\end{align}
\endef
%%%%%%%%%%%%%%%%%%%%%%%%%%%%%%%%%%%%%%%%%%%%%%%%%%%%%%%%%%%%%%%%%%%%%%%%%%%%%%%%%%%%%%%%%%%%%%%%%%
\theorem\label{thmfundamentalcoverequiv1}
$\Xt=\opair{\X}{\topology{}}$
is taken as a topological-space.
\begin{align}
&\Foreach{\acover}{\covers{\X}}\cr
&\[\(\acover\in\Fcov{\Xt}\)\thenn
\(\defset{\U}{\CSs{\X}}{\[\Foreach{\covelm}{\acover}\covelm\cap\U\in\stopology{\topology{}}{\covelm}\]}
=\topology{}\)\].\cr
&{}
\end{align}
\prooff
According to \refdef{defsubspacetopology1}
and
\refdef{deffundamentalcover},
it is clear.
\endthm
%%%%%%%%%%%%%%%%%%%%%%%%%%%%%%%%%%%%%%%%%%%%%%%%%%%%%%%%%%%%%%%%%%%%%%%%%%%%%%%%%%%%%%%%%%%%%%%%%%
\theorem\label{thmfundamentalcoverequivclosedsets}
$\Xt=\opair{\X}{\topology{}}$
is taken as a topological-space.
\begin{align}
&\Foreach{\acover}{\covers{\X}}\cr
&\[\(\acover\in\Fcov{\Xt}\)\thenn
\(\defset{\U}{\CSs{\X}}
{\[\Foreach{\covelm}{\acover}\covelm\cap\U\in\Fclosed{\covelm}{\stopology{\topology{}}{\covelm}}\]}
\subseteq\Fclosed{\X}{\topology{}}\)\].\cr
&{}
\end{align}
\prooff
$\acover$
is taken as an arbitrary element of $\covers{\X}$.
\begin{itemize}
\item[${\textbf{\textsf{p1}}}$]
It is assumed that,
\begin{equation}
\acover\in\Fcov{\Xt}.
\end{equation}
Then according to \refdef{deffundamentalcover},
\begin{equation}\label{thmfundamentalcoverequivclosedsetsp1eq1}
\defset{\U}{\CSs{\X}}{\[\Foreach{\covelm}{\acover}\covelm\cap\U\in\stopology{\topology{}}{\covelm}\]}
\subseteq\topology{}.
\end{equation}
\begin{itemize}
\item[${\textbf{\textsf{p1-1}}}$]
$\U$
is taken as such an arbitrary element of $\CSs{\X}$ that,
\begin{equation}\label{thmfundamentalcoverequivclosedsetsp1-1eq1}
\Foreach{\covelm}{\acover}\covelm\cap\U
\in\Fclosed{\covelm}{\stopology{\topology{}}{\covelm}}.
\end{equation}
Then according to
\refdef{deffamilyofclosedsets},
\begin{equation}\label{thmfundamentalcoverequivclosedsetsp1-1eq2}
\Foreach{\covelm}{\acover}\[\compl{\covelm}{\(\covelm\cap\U\)}\]
\in\stopology{\topology{}}{\covelm}.
\end{equation}
Hence considering that,
\begin{equation}
\Foreach{\covelm}{\acover}
\[\compl{\covelm}{\(\covelm\cap\U\)}\]=\[\covelm\cap\(\compl{\X}{\U}\)\],
\end{equation}
it is evident that,
\begin{equation}
\Foreach{\covelm}{\acover}\[\covelm\cap\(\compl{\X}{\U}\)\]
\in\stopology{\topology{}}{\covelm}.
\end{equation}
This and
\Ref{thmfundamentalcoverequivclosedsetsp1eq1}
imply,
\begin{equation}
\(\compl{\X}{\U}\)\in\topology{},
\end{equation}
which according to \refdef{deffamilyofclosedsets},
means,
\begin{equation}
\U\in\Fclosed{\X}{\topology{}}.
\end{equation}
\endp
\end{itemize}
Thus,
\begin{equation*}
\defset{\U}{\CSs{\X}}
{\[\Foreach{\covelm}{\acover}\covelm\cap\U\in\Fclosed{\covelm}{\stopology{\topology{}}{\covelm}}\]}
\subseteq\Fclosed{\X}{\topology{}}.
\end{equation*}
\endp
\end{itemize}
\begin{itemize}
\item[${\textbf{\textsf{p2}}}$]
It is assumed that,
\begin{equation}\label{thmfundamentalcoverequivclosedsetsp2eq1}
\defset{\U}{\CSs{\X}}
{\[\Foreach{\covelm}{\acover}\covelm\cap\U\in\Fclosed{\covelm}{\stopology{\topology{}}{\covelm}}\]}
\subseteq\Fclosed{\X}{\topology{}}.
\end{equation}
\begin{itemize}
\item[${\textbf{\textsf{p2-1}}}$]
$\U$
is taken as such an arbitrary element of $\CSs{\X}$ that,
\begin{equation}\label{thmfundamentalcoverequivclosedsetsp2-1eq1}
\Foreach{\covelm}{\acover}\covelm\cap\U
\in\stopology{\topology{}}{\covelm}.
\end{equation}
Then according to \refdef{deffamilyofclosedsets},
\begin{equation}\label{thmfundamentalcoverequivclosedsetsp2-1eq2}
\Foreach{\covelm}{\acover}\[\compl{\covelm}{\(\covelm\cap\U\)}\]
\in\Fclosed{\covelm}{\stopology{\topology{}}{\covelm}}.
\end{equation}
Thus considering that,
\begin{equation*}
\Foreach{\covelm}{\acover}
\[\compl{\covelm}{\(\covelm\cap\U\)}\]=\[\covelm\cap\(\compl{\X}{\U}\)\],
\end{equation*}
it is clear that,
\begin{equation}
\Foreach{\covelm}{\acover}\[\covelm\cap\(\compl{\X}{\U}\)\]
\in\Fclosed{\covelm}{\stopology{\topology{}}{\covelm}}.
\end{equation}
This and \Ref{thmfundamentalcoverequivclosedsetsp2eq1}
imply,
\begin{equation}
\(\compl{\X}{\U}\)\in\Fclosed{\X}{\topology{}},
\end{equation}
which according to \refdef{deffamilyofclosedsets}, means,
\begin{equation}
\U\in\topology{}.
\end{equation}
\endp
\end{itemize}
Thus,
\begin{equation}
\defset{\U}{\CSs{\X}}
{\[\Foreach{\covelm}{\acover}\covelm\cap\U\in\stopology{\topology{}}{\covelm}\]}
\subseteq\topology{},
\end{equation}
and hence according to \refdef{deffundamentalcover},
\begin{equation*}
\acover\in\Fcov{\Xt}.
\end{equation*}
\endp
\end{itemize}
\endthm
%%%%%%%%%%%%%%%%%%%%%%%%%%%%%%%%%%%%%%%%%%%%%%%%%%%%%%%%%%%%%%%%%%%%%%%%%%%%%%%%%%%%%%%%%%%%%%%%%%
\theorem\label{thmfundamentalcoverequivclosedsets1}
$\Xt=\opair{\X}{\topology{}}$
is taken as a topological-space.
\begin{align}
&\Foreach{\acover}{\covers{\X}}\cr
&\[\(\acover\in\Fcov{\Xt}\)\thenn
\(\defset{\U}{\CSs{\X}}
{\[\Foreach{\covelm}{\acover}\covelm\cap\U\in\Fclosed{\covelm}{\stopology{\topology{}}{\covelm}}\]}
=\Fclosed{\X}{\topology{}}\)\].\cr
&{}
\end{align}
\prooff
According to \refthm{thmsubspaceclosedsets}
and
\refthm{thmfundamentalcoverequivclosedsets},
it is clear.
\endthm
%%%%%%%%%%%%%%%%%%%%%%%%%%%%%%%%%%%%%%%%%%%%%%%%%%%%%%%%%%%%%%%%%%%%%%%%%%%%%%%%%%%%%%%%%%%%%%%%%%
\theorem\label{thmrefinementfundamentalcover}
$\Xt=\opair{\X}{\topology{}}$
is taken as a topological-space,
$\opair{\X}{\acover}$
as a cover of $\X$, and
$\opair{\X}{\p{\acover}}$
as a refinement of $\opair{\X}{\acover}$. If $\p{\acover}$
is a fundamental cover of $\Xt$, then $\acover$
is also a fundamental cover of $\Xt$. That is,
\begin{equation}
\p{\acover}\in\Fcov{\Xt}\then\acover\in\Fcov{\Xt}.
\end{equation}
\prooff
According to \refdef{defrefinedcover},
\begin{equation}\label{thmrefinementfundamentalcoverp1}
\Foreach{\p{\covelm}}{\p{\acover}}
\[\Exists{\covelm}{\acover}\p{\covelm}\subseteq\covelm\].
\end{equation}
Thus there exists a function from $\p{\acover}$ to $\acover$, like $\eta$, such that,
\begin{equation}\label{thmrefinementfundamentalcoverp2}
\Foreach{\p{\covelm}}{\p{\acover}}
\p{\covelm}\subseteq\func{\eta}{\p{\covelm}}.
\end{equation}
(Why?)
\begin{itemize}
\item[${\textbf{\textsf{p1}}}$]
It is assumed that $\p{\acover}$
is a fundamental cover of $\Xt$. Then,
\begin{equation}\label{thmrefinementfundamentalcoverp1eq1}
\defset{\U}{\CSs{\X}}
{\[\Foreach{\p{\covelm}}{\p{\acover}}
\(\p{\covelm}\cap\U\)\in\stopology{\topology{}}{\p{\covelm}}}\]\subseteq\topology{}.
\end{equation}
\begin{itemize}
\item[${\textbf{\textsf{p1-1}}}$]
$\U$
is taken as such an arbitrary element of $\CSs{\X}$ that,
\begin{equation}\label{thmrefinementfundamentalcoverp1-1eq1}
\Foreach{\covelm}{\acover}\(\covelm\cap\U\)\in\stopology{\topology{}}{\covelm}.
\end{equation}
Hence considering that the codomain of $\eta$ is $\acover$, it is clear that,
\begin{equation}
\Foreach{\p{\covelm}}{\p{\acover}}\label{thmrefinementfundamentalcoverp1-1eq2}
\[\func{\eta}{\p{\covelm}}\cap\U\]\in\stopology{\topology{}}{\func{\eta}{\p{\covelm}}}.
\end{equation}
\Ref{thmrefinementfundamentalcoverp2}
and
\refdef{defsubspacetopology1}
imply,
\begin{equation}\label{thmrefinementfundamentalcoverp1-1eq3}
\Foreach{\p{\covelm}}{\p{\acover}}
\p{\covelm}\cap\[\func{\eta}{\p{\covelm}}\cap\U\]\in
\stopology{\stopology{\topology{}}{\func{\eta}{\p{\covelm}}}}{\p{\covelm}}.
\end{equation}
\Ref{thmrefinementfundamentalcoverp2}
and
\refthm{thmsubspacetopologytransitivity}
imply,
\begin{equation}\label{thmrefinementfundamentalcoverp1-1eq4}
\Foreach{\p{\covelm}}{\p{\acover}}
\stopology{\stopology{\topology{}}{\func{\eta}{\p{\covelm}}}}{\p{\covelm}}=
\stopology{\topology{}}{\p{\covelm}}.
\end{equation}
Additionally, according to \Ref{thmrefinementfundamentalcoverp2},
it is evident that,
\begin{equation}\label{thmrefinementfundamentalcoverp1-1eq5}
\Foreach{\p{\covelm}}{\p{\acover}}
\p{\covelm}\cap\[\func{\eta}{\p{\covelm}}\cap\U\]=
\p{\covelm}\cap\U.
\end{equation}
\Ref{thmrefinementfundamentalcoverp1-1eq3},
\Ref{thmrefinementfundamentalcoverp1-1eq4},
and
\Ref{thmrefinementfundamentalcoverp1-1eq5}
imply,
\begin{equation}
\Foreach{\p{\covelm}}{\p{\acover}}
\(\p{\covelm}\cap\U\)\in\stopology{\topology{}}{\p{\covelm}}.
\end{equation}
\endp
\end{itemize}
Therefore,
\begin{gather}
\defset{\U}{\CSs{\X}}
{\[\Foreach{\covelm}{\acover}
\(\covelm\cap\U\)\in\stopology{\topology{}}{\covelm}}\]\cr
\subseteq\cr
\defset{\U}{\CSs{\X}}
{\[\Foreach{\p{\covelm}}{\p{\acover}}
\(\p{\covelm}\cap\U\)\in\stopology{\topology{}}{\p{\covelm}}}\].
\end{gather}
This and
\Ref{thmrefinementfundamentalcoverp1eq1}
imply,
\begin{equation}
\defset{\U}{\CSs{\X}}
{\[\Foreach{\covelm}{\acover}
\(\covelm\cap\U\)\in\stopology{\topology{}}{\covelm}}\]\subseteq\topology{},
\end{equation}
which means
$\acover$
is a fundamental cover of .$\Xt$
\endp
\end{itemize}
\endthm
%%%%%%%%%%%%%%%%%%%%%%%%%%%%%%%%%%%%%%%%%%%%%%%%%%%%%%%%%%%%%%%%%%%%%%%%%%%%%%%%%%%%%%%%%%%%%%%%%%
\theorem\label{thmrefinementfundamentalcover}
$\Xt=\opair{\X}{\topology{}}$
is taken as a topological-space.
For every cover of $\X$, like $\acover$,
if every point $\point$ of $\Xt$
possesses a neighbourhood $\V$ in $\Xt$
such that the set $\left\{\V\cap\U~|~\U\in\acover\right\}$
is a fundamental cover of $\Xt$, then $\acover$
is a fundamental cover of $\Xt$. That is,
\begin{align}
&\Foreach{\acover}{\covers{\X}}\bigg(\bigg[\Foreach{\point}{\X}\cr
&\left.\bigg(
\Exists{\V}{\func{\nei{\Xt}}{\seta{\point}}}\bigg[\defset{\asubset}{\CSs{\X}}
{\(\Exists{\U}{\acover}\asubset=\V\cap\U\)}\in\Fcov{\Xt}\bigg]\bigg)\]\cr
&\then\acover\in\Fcov{\Xt}\bigg).
\end{align}
\prooff

\endthm
%%%%%%%%%%%%%%%%%%%%%%%%%%%%%%%%%%%%%%%%%%%%%%%%%%%%%%%%%%%%%%%%%%%%%%%%%%%%%%%%%%%%%%%%%%%%%%%%%%
\theorem\label{thmcontinuityandfundamentalcover}
Each
$\Xt=\opair{\X}{\topology{\X}}$
and
$\Xt=\opair{\Y}{\topology{\Y}}$
is taken as a topological-space, and $\acover$ an element of $\Fcov{\Xt}$
(a fundamental cover of $\Xt$).
For every $\cf$ in $\Func{\X}{\Y}$,
if for every $\covelm$ in $\acover$,
the domain-restriction of $\cf$ to $\covelm$
is a continuous map from $\opair{\covelm}{\stopology{\topology{\X}}{\covelm}}$ to $\Yt$, then
$\cf$
is a continuous map from $\Xt$ to $\Yt$. That is,
\begin{align}
\defset{\cf}{\Func{\X}{\Y}}
{\[\Foreach{\covelm}{\acover}
\func{\resd{\cf}}{\covelm}\in\CF{\opair{\covelm}{\stopology{\topology{\X}}{\covelm}}}{\Yt}\]}
\subseteq\CF{\Xt}{\Yt}.
\end{align}
\prooff
Considering that $\acover$ is a fundamental cover of $\Xt$, according to \refdef{deffundamentalcover},
\begin{equation}\label{thmcontinuityandfundamentalcoverpeq1}
\defset{\U}{\CSs{\X}}
{\[\Foreach{\covelm}{\acover}
\(\covelm\cap\U\)\in\stopology{\topology{\X}}{\covelm}}\]\subseteq\topology{\X}.
\end{equation}
\begin{itemize}
\item[${\textbf{\textsf{p1}}}$]
$\cf$
is taken as an arbitrary element of $\Func{\X}{\Y}$.
As an evident property of the domain-restriction-map of a function,
\begin{equation}\label{thmcontinuityandfundamentalcoverp1eq2}
\Foreach{\covelm}{\acover}
\bigg(\Foreach{\V}{\CSs{\Y}}
\func{\pimage{\[\func{\resd{\cf}}{\covelm}\]}}{\V}=
\[\covelm\cap\func{\pimage{\cf}}{\V}\]\bigg).
\end{equation}
\begin{itemize}
\item[${\textbf{\textsf{p1-1}}}$]
It is assumed that for every $\covelm$ in $\acover$, the restriction of $\cf$ to $\covelm$
is a continuous map from the topological-space $\opair{\covelm}{\stopology{\topology{\X}}{\covelm}}$ to $\Yt$.
That is,
\begin{equation}\label{thmcontinuityandfundamentalcoverp1-1eq1}
\Foreach{\covelm}{\acover}
\[\func{\resd{\cf}}{\covelm}\in
\CF{\opair{\covelm}{\stopology{\topology{\X}}{\covelm}}}{\Yt}\].
\end{equation}
Then according to \refdef{defcontinuousfunction},
\begin{equation}\label{thmcontinuityandfundamentalcoverp1-1eq2}
\Foreach{\covelm}{\acover}
\bigg(\Foreach{\V}{\topology{\Y}}
\func{\[\pimage{\func{\resd{\cf}}{\covelm}\]}}{\V}
\in\stopology{\topology{\X}}{\covelm}\bigg).
\end{equation}
This and
\Ref{thmcontinuityandfundamentalcoverp1eq2}
imply,
\begin{equation}
\Foreach{\covelm}{\acover}
\bigg(\Foreach{\V}{\topology{\Y}}
\[\covelm\cap\func{\pimage{\cf}}{\V}\]
\in\stopology{\topology{\X}}{\covelm}\bigg),
\end{equation}
or equivalently,
\begin{equation}
\Foreach{\V}{\topology{\Y}}
\bigg(\Foreach{\covelm}{\acover}
\[\covelm\cap\func{\pimage{\cf}}{\V}\]
\in\stopology{\topology{\X}}{\covelm}\bigg).
\end{equation}
This and \Ref{thmcontinuityandfundamentalcoverpeq1} imply,
\begin{equation}
\Foreach{\V}{\topology{\Y}}
\func{\pimage{\cf}}{\V}\in\topology{\X},
\end{equation}
which according to \refdef{defcontinuousfunction}, means,
\begin{equation}
\cf\in\CF{\Xt}{\Yt}.
\end{equation}
\endp
\end{itemize}
\endp
\end{itemize}
\endthm
%%%%%%%%%%%%%%%%%%%%%%%%%%%%%%%%%%%%%%%%%%%%%%%%%%%%%%%%%%%%%%%%%%%%%%%%%%%%%%%%%%%%%%%%%%%%%%%%%%
\subsection{
Open Covers of a Topological Space}
\definition\label{defopencover}
$\Xt=\opair{\X}{\topology{}}$
is taken as a topological-space. For every $\acover$ in $\CSs{\CSs{\X}}$, $\acover$
is referred to as an $\quotl$open cover of the topological-space $\Xt$$\quotr$
iff these properties are hold.
\begin{itemize}
\item[${\textbf{\textsf{OCOV1}}}$]
$\acover$ is a cover of $\X$, that is
$\acover\in\covers{\X}$.
\item[${\textbf{\textsf{OCOV2}}}$]
$\acover\subseteq\topology{}.$
\end{itemize}
The set of all open covers of $\Xt$ is denoted by $\Ocov{\Xt}$.
\begin{equation}
\Ocov{\Xt}:=\defset{\acover}{\covers{\X}}{\acover\subseteq\topology{}}.
\end{equation}
\endef
%%%%%%%%%%%%%%%%%%%%%%%%%%%%%%%%%%%%%%%%%%%%%%%%%%%%%%%%%%%%%%%%%%%%%%%%%%%%%%%%%%%%%%%%%%%%%%%%%%
\theorem\label{thmeveryopencoverisafundamentalcover}
$\Xt=\opair{\X}{\topology{}}$
is taken as a topological-space.
Every open cover of $\Xt$ is a fundamental cover of $\Xt$. That is,
\begin{equation}
\Ocov{\Xt}\subseteq\Fcov{\Xt}.
\end{equation}
\prooff
$\acover$
is taken as an arbitrary element of $\Ocov{\Xt}$. Then according to \refdef{defopencover},
\begin{align}
&\acover\in\covers{\X},\label{thmeveryopencoverisafundamentalcoverpeq1}\\
&\Foreach{\covelm}{\acover}\covelm\in\topology{}.\label{thmeveryopencoverisafundamentalcoverpeq2}
\end{align}
\begin{itemize}
\item[${\textbf{\textsf{p1}}}$]
$\U$
is taken as such an element of $\CSs{\X}$ that,
\begin{equation}\label{thmeveryopencoverisafundamentalcoverp1eq1}
\Foreach{\covelm}{\acover}
\covelm\cap\U\in\stopology{\topology{}}{\covelm}.
\end{equation}
Consequently, based on \refthm{thmsubsubopen} and according to \Ref{thmeveryopencoverisafundamentalcoverpeq2},
\begin{equation}\label{thmeveryopencoverisafundamentalcoverp1eq2}
\Foreach{\covelm}{\acover}\(\covelm\cap\U\)\in\topology{},
\end{equation}
and hence according to \refdef{deftopologicalspace},
\begin{equation}\label{thmeveryopencoverisafundamentalcoverp1eq3}
\[\Union{\covelm}{\acover}{\(\covelm\cap\U\)}\]\in\topology{}.
\end{equation}
In addition, considering that,
\begin{equation}\label{thmeveryopencoverisafundamentalcoverp1eq4}
\(\union{\acover}\)=\X,
\end{equation}
it is clear that,
\begin{align}\label{thmeveryopencoverisafundamentalcoverp1eq5}
\[\Union{\covelm}{\acover}{\(\covelm\cap\U\)}\]&=
\(\Union{\covelm}{\acover}{\covelm}\)\cap\U\cr
&=\X\cap\U\cr
&=\U.
\end{align}
\Ref{thmeveryopencoverisafundamentalcoverp1eq3}
and
\Ref{thmeveryopencoverisafundamentalcoverp1eq5}
imply,
\begin{equation}
\U\in\topology{}.
\end{equation}
\endp
\end{itemize}
Therefore,
\begin{equation}
\defset{\U}{\CSs{\X}}
{\[\Foreach{\covelm}{\acover}\covelm\cap\U\in\stopology{\topology{}}{\covelm}\]}
\subseteq\topology{}.
\end{equation}
Consequently, based on \refdef{deffundamentalcover}, and according to \Ref{thmeveryopencoverisafundamentalcoverpeq1},
\begin{equation}
\acover\in\Fcov{\Xt}.
\end{equation}
\endthm
%%%%%%%%%%%%%%%%%%%%%%%%%%%%%%%%%%%%%%%%%%%%%%%%%%%%%%%%%%%%%%%%%%%%%%%%%%%%%%%%%%%%%%%%%%%%%%%%%%
\theorem\label{thmcontinuityandopencover}
Each
$\Xt=\opair{\X}{\topology{\X}}$
and
$\Xt=\opair{\Y}{\topology{\Y}}$
is taken as a topological-space, and $\acover$ as an element of $\Ocov{\Xt}$
(an open cover of $\Xt$).
For every $\cf$ in $\Func{\X}{\Y}$, if for every $\covelm$ in $\acover$,
the domain-restriction of $\cf$ to $\covelm$
is a continuous map from the topological-space $\opair{\covelm}{\stopology{\topology{\X}}{\covelm}}$ to $\Yt$, then
$\cf$ is a continuous map from $\Xt$ to $\Yt$. That is,
\begin{align*}
\defset{\cf}{\Func{\X}{\Y}}
{\[\Foreach{\covelm}{\acover}
\func{\resd{\cf}}{\covelm}\in\CF{\opair{\covelm}{\stopology{\topology{\X}}{\covelm}}}{\Yt}\]}
\subseteq\CF{\Xt}{\Yt}.
\end{align*}
\prooff
According to \refthm{thmcontinuityandfundamentalcover}
and
\refthm{thmeveryopencoverisafundamentalcover}, it is clear.
\endthm
%%%%%%%%%%%%%%%%%%%%%%%%%%%%%%%%%%%%%%%%%%%%%%%%%%%%%%%%%%%%%%%%%%%%%%%%%%%%%%%%%%%%%%%%%%%%%%%%%%
\subsection{Closed Covers of a Topological Space}
\definition\label{defclosedcover}
$\Xt=\opair{\X}{\topology{}}$
is taken as a topological-space.
For every $\acover$ in $\CSs{\CSs{\X}}$,
$\acover$
is referred to as a $\quotl$closed cover of the topological-space $\Xt$$\quotr$ iff
these properties are hold.
\begin{itemize}
\item[${\textbf{\textsf{CCOV1}}}$]
$\acover$
is a cover of $\X$, that is
$\acover\in\covers{\X}$.
\item[${\textbf{\textsf{CCOV2}}}$]
$\acover\subseteq\Fclosed{\X}{\topology{}}.$
\end{itemize}
The set of all closed covers of $\Xt$ is denoted by $\Ccov{\Xt}$.
\begin{equation}
\Ccov{\Xt}:=\defset{\acover}{\covers{\X}}{\acover\subseteq\Fclosed{\X}{\topology{}}}.
\end{equation}
\begin{itemize}
\item
Every finite element of $\Ccov{\Xt}$
is referred to as a $\quotl$finite closed cover of $\Xt$$\quotr$.
\end{itemize}
\endef
%%%%%%%%%%%%%%%%%%%%%%%%%%%%%%%%%%%%%%%%%%%%%%%%%%%%%%%%%%%%%%%%%%%%%%%%%%%%%%%%%%%%%%%%%%%%%%%%%%
\theorem\label{thmeveryfiniteclosedcoverisafundamentalcover}
$\Xt=\opair{\X}{\topology{}}$
is taken as a topological-space.
Every non-empty finite closed cover of
\footnote{
The non-emptiness constraint is imposed in order to discard  the case where
$\empty$ is a closed cover of the topological-space
$\opair{\empty}{\seta{\empty}}$, when $\Xt$ is taken to be this topological-space.
}
$\Xt$
is a fundamental cover of $\Xt$. That is
\begin{equation}
\defset{\acover}{\Ccov{\Xt}}{\CarD{\acover}\in\Zp}\subseteq\Fcov{\Xt}.
\end{equation}
\prooff
$\acover$
is taken as a finite closed cover of $\Xt$. Then according to \refdef{defclosedcover},
\begin{align}
&\acover\in\covers{\X},
\label{thmeveryfiniteclosedcoverisafundamentalcoverpeq1}\\
&\Foreach{\covelm}{\acover}\covelm\in\Fclosed{\X}{\topology{}},
\label{thmeveryfiniteclosedcoverisafundamentalcoverpeq2}\\
&\CarD{\acover}\in\Zp.
\label{thmeveryfiniteclosedcoverisafundamentalcoverpeq3}
\end{align}
\begin{itemize}
\item[${\textbf{\textsf{p1}}}$]
$\U$
is taken as such an arbitrary element of $\CSs{\X}$ that,
\begin{equation}\label{thmeveryfiniteclosedcoverisafundamentalcoverp1eq1}
\Foreach{\covelm}{\acover}
\covelm\cap\U\in\Fclosed{\covelm}{\stopology{\topology{}}{\covelm}}.
\end{equation}
Base on \refthm{thmsubsubclosed}, this and \Ref{thmeveryfiniteclosedcoverisafundamentalcoverpeq2}
imply,
\begin{equation}\label{thmeveryfiniteclosedcoverisafundamentalcoverp1eq2}
\Foreach{\covelm}{\acover}\(\covelm\cap\U\)\in\Fclosed{\X}{\topology{}},
\end{equation}
Based on \refthm{thmclosedsets},
this and \Ref{thmeveryfiniteclosedcoverisafundamentalcoverpeq3}
imply,
\begin{equation}\label{thmeveryfiniteclosedcoverisafundamentalcoverp1eq3}
\[\Union{\covelm}{\acover}{\(\covelm\cap\U\)}\]\in\Fclosed{\X}{\topology{}}.
\end{equation}
In addition, considering that,
\begin{equation}\label{thmeveryfiniteclosedcoverisafundamentalcoverp1eq4}
\(\union{\acover}\)=\X,
\end{equation}
it is evident that,
\begin{align}\label{thmeveryfiniteclosedcoverisafundamentalcoverp1eq5}
\[\Union{\covelm}{\acover}{\(\covelm\cap\U\)}\]&=
\(\Union{\covelm}{\acover}{\covelm}\)\cap\U\cr
&=\X\cap\U\cr
&=\U.
\end{align}
\Ref{thmeveryfiniteclosedcoverisafundamentalcoverp1eq3}
and
\Ref{thmeveryfiniteclosedcoverisafundamentalcoverp1eq5}
imply,
\begin{equation}
\U\in\Fclosed{\X}{\topology{}}.
\end{equation}
\endp
\end{itemize}
Therefore,
\begin{equation}
\defset{\U}{\CSs{\X}}
{\[\Foreach{\covelm}{\acover}\covelm\cap\U\in\Fclosed{\covelm}{\stopology{\topology{}}{\covelm}}\]}
\subseteq\Fclosed{\X}{\topology{}}.
\end{equation}
Based on \refthm{thmfundamentalcoverequivclosedsets},
this and
\Ref{thmeveryfiniteclosedcoverisafundamentalcoverpeq1}
imply,
\begin{equation}
\acover\in\Fcov{\Xt}.
\end{equation}
\endthm
%%%%%%%%%%%%%%%%%%%%%%%%%%%%%%%%%%%%%%%%%%%%%%%%%%%%%%%%%%%%%%%%%%%%%%%%%%%%%%%%%%%%%%%%%%%%%%%%%%
\theorem\label{thmcontinuityandfiniteclosecover}
Each
$\Xt=\opair{\X}{\topology{\X}}$
and
$\Xt=\opair{\Y}{\topology{\Y}}$
is taken as a topological-space, and
$\acover$
an element of $\defset{\acover}{\Ccov{\Xt}}{\CarD{\acover}\in\Zp}$
(a finite closed cover of $\Xt$).
For every $\cf$ in $\Func{\X}{\Y}$,
if for every $\covelm$ in $\acover$,
the domain-restriction of $\cf$ to $\covelm$
is a continuous map from the topological-space $\opair{\covelm}{\stopology{\topology{\X}}{\covelm}}$ to $\Yt$, then
$\cf$
is a continuous map from $\Xt$ to . That is,$\Yt$
\begin{align*}
\defset{\cf}{\Func{\X}{\Y}}
{\[\Foreach{\covelm}{\acover}
\func{\resd{\cf}}{\covelm}\in\CF{\opair{\covelm}{\stopology{\topology{\X}}{\covelm}}}{\Yt}\]}
\subseteq\CF{\Xt}{\Yt}.
\end{align*}
\prooff
According to
\refthm{thmcontinuityandfundamentalcover}
and
\refthm{thmeveryfiniteclosedcoverisafundamentalcover},
it is clear.
\endthm
%%%%%%%%%%%%%%%%%%%%%%%%%%%%%%%%%%%%%%%%%%%%%%%%%%%%%%%%%%%%%%%%%%%%%%%%%%%%%%%%%%%%%%%%%%%%%%%%%%
\subsection{Locally-Finite Covers}
\definition\label{deflocallyfinitecover}
$\Xt=\opair{\X}{\topology{}}$
is taken as a topological-space. For every $\acover$ in $\CSs{\CSs{\X}}$,
$\acover$
is referred to as a $\quotl$locally-finite cover of the topological-space $\Xt$$\quotr$ iff
these properties are hold.
\begin{itemize}
\item[${\textbf{\textsf{LfCOV1}}}$]
$\acover$
is a cover of $\X$, that is
$\acover\in\covers{\X}$.
\item[${\textbf{\textsf{LfCOV2}}}$]
$\Foreach{\point}{\X}\[\Exists{\U}{\func{\nei{\Xt}}{\seta{\point}}}
\CarD{\defset{\covelm}{\acover}{\U\cap\covelm\neq\empty}}\in\Zpz\].$
\end{itemize}
The set of all locally-finite covers of $\Xt$ is denoted by $\Lfcov{\Xt}$.
\begin{itemize}
\item
Each element of $\Ccov{\Xt}\cap\Lfcov{\Xt}$
is referred to as a $\quotl$closed locally-finite cover of $\Xt$$\quotr$.
\end{itemize}
\endef
%%%%%%%%%%%%%%%%%%%%%%%%%%%%%%%%%%%%%%%%%%%%%%%%%%%%%%%%%%%%%%%%%%%%%%%%%%%%%%%%%%%%%%%%%%%%%%%%%%
\theorem\label{thmeverylocallyfiniteclosedcoverisafundamentalcover}
$\Xt=\opair{\X}{\topology{}}$
is taken as a topological-space.
Every closed locally-finite cover of $\Xt$ is a fundamental cover of $\Xt$. That is,
\begin{equation}
\(\Ccov{\Xt}\cap\Lfcov{\Xt}\)\subseteq\Fcov{\Xt}.
\end{equation}
\prooff
$\acover$
is taken as an arbitrary element of $\Ccov{\Xt}\cap\Lfcov{\Xt}$. Then according to \refdef{defclosedcover}
and
\refdef{deflocallyfinitecover},
\begin{align}
&\acover\in\covers{\X},
\label{thmeverylocallyfiniteclosedcoverisafundamentalcoverpeq1}\\
&\Foreach{\covelm}{\acover}\covelm\in\Fclosed{\X}{\topology{}},
\label{thmeverylocallyfiniteclosedcoverisafundamentalcoverpeq2}\\
&\Foreach{\point}{\X}\[\Exists{\U}{\func{\nei{\Xt}}{\seta{\point}}}
\CarD{\defset{\covelm}{\acover}{\U\cap\covelm\neq\empty}}\in\Zpz\].
\label{thmeverylocallyfiniteclosedcoverisafundamentalcoverpeq3}
\end{align}
\begin{itemize}
\item[${\textbf{\textsf{p1}}}$]
$\asubset$
is taken as such an element of $\CSs{\X}$ that,
\begin{equation}
\Foreach{\covelm}{\acover}\(\covelm\cap\asubset\)\in
\Fclosed{\covelm}{\stopology{\topology{}}{\covelm}}.
\end{equation}
Based on \refthm{thmsubsubclosed},
this and
\Ref{thmeverylocallyfiniteclosedcoverisafundamentalcoverpeq2}
imply,
\begin{equation}\label{thmeveryfiniteclosedcoverisafundamentalcoverp1eq2}
\Foreach{\covelm}{\acover}\(\covelm\cap\asubset\)\in\Fclosed{\X}{\topology{}},
\end{equation}
\begin{itemize}
\item[${\textbf{\textsf{p1-1}}}$]
$\point$
is taken as an arbitrary element of $\(\compl{\X}{\asubset}\)$.

\end{itemize}
\end{itemize}
\endthm
%%%%%%%%%%%%%%%%%%%%%%%%%%%%%%%%%%%%%%%%%%%%%%%%%%%%%%%%%%%%%%%%%%%%%%%%%%%%%%%%%%%%%%%%%%%%%%%%%%
\theorem\label{thmcontinuityandlocallyfiniteclosecover}
Each
$\Xt=\opair{\X}{\topology{\X}}$
and
$\Xt=\opair{\Y}{\topology{\Y}}$
is taken as a topological-space, and $\acover$ an element of $\Ccov{\Xt}\cap\Lfcov{\Xt}$
(a closed locally-finite cover of $\Xt$). For every $\cf$ in $\Func{\X}{\Y}$,
if for every $\covelm$ in $\acover$ the domain restriction of $\cf$ to $\covelm$
is a continuous map from $\opair{\covelm}{\stopology{\topology{\X}}{\covelm}}$ to $\Yt$, then
$\cf$ is a continuous map from $\Xt$ to $\Yt$. That is,
\begin{align*}
\defset{\cf}{\Func{\X}{\Y}}
{\[\Foreach{\covelm}{\acover}
\func{\resd{\cf}}{\covelm}\in\CF{\opair{\covelm}{\stopology{\topology{\X}}{\covelm}}}{\Yt}\]}
\subseteq\CF{\Xt}{\Yt}.
\end{align*}
\prooff
According to \refthm{thmcontinuityandfundamentalcover} and
\refthm{thmeverylocallyfiniteclosedcoverisafundamentalcover},
it is clear.
\endthm
%%%%%%%%%%%%%%%%%%%%%%%%%%%%%%%%%%%%%%%%%%%%%%%%%%%%%%%%%%%%%%%%%%%%%%%%%%%%%%%%%%%%%%%%%%%%%%%%%%%%%%%%%%%%%%%%%%%%%%%%%%%%%
%%%%%%%%%%%%%%%%%%%%%%%%%%%%%%%%%%%%%%%%%%%%%%%%%%%%%%%%%%%%%%%%%%%%%%%%%%%%%%%%%%%%%%%%%%%%%%%%%%%%%%%%%%%%%%%%%%%%%%%%%%%%%
%%%%%%%%%%%%%%%%%%%%%%%%%%%%%%%%%%%%%%%%%%%%%%%%%%%%%%%%%%%%%%%%%%%%%%%%%%%%%%%%%%%%%%%%%%%%%%%%%%%%%%%%%%%%%%%%%%%%%%%%%%%%%
%%%%%%%%%%%%%%%%%%%%%%%%%%%%%%%%%%%%%%%%%%%%%%%%%%%%%%%%%%%%%%%%%%%%%%%%%%%%%%%%%%%%%%%%%%%%%%%%%%%%%%%%%%%%%%%%%%%%%%%%%%%%%
%%%%%%%%%%%%%%%%%%%%%%%%%%%%%%%%%%%%%%%%%%%%%%%%%%%%%%%%%%%%%%%%%%%%%%%%%%%%%%%%%%%%%%%%%%%%%%%%%%%%%%%%%%%%%%%%%%%%%%%%%%%%%
%%%%%%%%%%%%%%%%%%%%%%%%%%%%%%%%%%%%%%%%%%%%%%%%%%%%%%%%%%%%%%%%%%%%%%%%%%%%%%%%%%%%%%%%%%%%%%%%%%%%%%%%%%%%%%%%%%%%%%%%%%%%%
%%%%%%%%%%%%%%%%%%%%%%%%%%%%%%%%%%%%%%%%%%%%%%%%%%%%%%%%%%%%%%%%%%%%%%%%%%%%%%%%%%%%%%%%%%%%%%%%%%%%%%%%%%%%%%%%%%%%%%%%%%%%%
%%%%%%%%%%%%%%%%%%%%%%%%%%%%%%%%%%%%%%%%%%%%%%%%%%%%%%%%%%%%%%%%%%%%%%%%%%%%%%%%%%%%%%%%%%%%%%%%%%%%%%%%%%%%%%%%%%%%%%%%%%%%%
%%%%%%%%%%%%%%%%%%%%%%%%%%%%%%%%%%%%%%%%%%%%%%%%%%%%%%%%%%%%%%%%%%%%%%%%%%%%%%%%%%%%%%%%%%%%%%%%%%%%%%%%%%%%%%%%%%%%%%%%%%%%%
%%%%%%%%%%%%%%%%%%%%%%%%%%%%%%%%%%%%%%%%%%%%%%%%%%%%%%%%%%%%%%%%%%%%%%%%%%%%%%%%%%%%%%%%%%%%%%%%%%%%%%%%%%%%%%%%%%%%%%%%%%%%%
%%%%%%%%%%%%%%%%%%%%%%%%%%%%%%%%%%%%%%%%%%%%%%%%%%%%%%%%%%%%%%%%%%%%%%%%%%%%%%%%%%%%%%%%%%%%%%%%%%%%%%%%%%%%%%%%%%%%%%%%%%%%%
%%%%%%%%%%%%%%%%%%%%%%%%%%%%%%%%%%%%%%%%%%%%%%%%%%%%%%%%%%%%%%%%%%%%%%%%%%%%%%%%%%%%%%%%%%%%%%%%%%%%%%%%%%%%%%%%%%%%%%%%%%%%%
\section{
Open Maps and Closed Maps
}
\definition\label{defopenmap}
Each
$\Xt=\opair{\X}{\topology{\X}}$
and
$\Yt=\opair{\Y}{\topology{\Y}}$
is taken as a topological-space.
For every function $\cf$ in $\Func{\X}{\Y}$,
$\cf$
is reffered to as a $\quotl$open map from $\Xt$ to $\Yt$$\quotr$ iff,
\begin{equation}
\Foreach{\U}{\topology{\X}}
\[\func{\image{\cf}}{\U}\in\topology{\Y}\].
\end{equation}
The set of all open maps from $\Xt$ to $\Yt$ is denoted by $\OM{\Xt}{\Yt}$. That is,
\begin{equation}
\OM{\Xt}{\Yt}:=\defset{\om}{\Func{\X}{\Y}}
{\big(\Foreach{\U}{\topology{\X}}
\[\func{\image{\cf}}{\U}\in\topology{\Y}\]\big)}.
\end{equation}
\endef
%%%%%%%%%%%%%%%%%%%%%%%%%%%%%%%%%%%%%%%%%%%%%%%%%%%%%%%%%%%%%%%%%%%%%%%%%%%%%%%%%%%%%%%
\definition\label{defclosedmap}
Each
$\Xt=\opair{\X}{\topology{\X}}$
and
$\Yt=\opair{\Y}{\topology{\Y}}$
is taken as a topological-space.
For every function $\cf$ in $\Func{\X}{\Y}$,
$\cf$
is referred to as a $\quotl$closed map from $\Xt$ to $\Yt$$\quotr$ iff,
\begin{equation}
\Foreach{\U}{\Fclosed{\X}{\topology{\X}}}
\[\func{\image{\cf}}{\U}\in\Fclosed{\Y}{\topology{\Y}}\].
\end{equation}
The set of all open maps from $\Xt$ to $\Yt$ is denoted by $\CM{\Xt}{\Yt}$. That is,
\begin{align}
\CM{\Xt}{\Yt}:=\defset{\om}{\Func{\X}{\Y}}
{\big(\Foreach{\U}{\Fclosed{\X}{\topology{\X}}}
\[\func{\image{\cf}}{\U}\in\Fclosed{\Y}{\topology{\Y}}\]\big)}.
\end{align}
\endef
%%%%%%%%%%%%%%%%%%%%%%%%%%%%%%%%%%%%%%%%%%%%%%%%%%%%%%%%%%%%%%%%%%%%%%%%%%%%%%%%%%%%%%%
\theorem\label{thmclosedmapproperty1}
Each
$\Xt=\opair{\X}{\topology{\X}}$
and
$\Yt=\opair{\Y}{\topology{\Y}}$
is taken as a topological-space, and $\om$ as an element of $\CM{\Xt}{\Yt}$
(a closed map from $\Xt$ to $\Yt$).
\begin{align}
&\Foreach{\Asubset{\Y}}{\CSs{\Y}}\cr
&\[\Foreach{\U}{\defset{\U}{\topology{\X}}{\U\supseteq\func{\pimage{\om}}{\Asubset{\Y}}}}
\bigg(\Exists{\V}{\defset{\V}{\topology{\Y}}{\V\supseteq\Asubset{\Y}}}
\func{\pimage{\om}}{\V}\subseteq\U\bigg)\].\cr
&{}
\end{align}
\prooff
According to \refdef{defclosedmap},
\begin{equation}\label{thmclosedmapproperty1peq1}
\Foreach{\U}{\Fclosed{\X}{\topology{\X}}}
\func{\image{\om}}{\U}\in\Fclosed{\Y}{\topology{\Y}}.
\end{equation}
\begin{itemize}
\item[${\textbf{\textsf{p1}}}$]
$\Asubset{\Y}$
is taken as an element of $\CSs{\Y}$.
\begin{itemize}
\item[${\textbf{\textsf{p1-1}}}$]
$\U$
is taken as such an element of $\topology{\X}$ that,
\begin{equation}\label{thmclosedmapproperty1p1-1eq1}
\U\supseteq\func{\pimage{\om}}{\Asubset{\Y}}.
\end{equation}
$\V$
is defined as,
\begin{equation}\label{thmclosedmapproperty1p1-1eq2}
\V:=\compl{\Y}{\func{\image{\om}}{\compl{\X}{\U}}}.
\end{equation}
Considering that,
$\U\in\topology{\X}$,
According to,
\refdef{deffamilyofclosedsets},
\begin{equation}\label{thmclosedmapproperty1p1-1eq3}
\(\compl{\X}{\U}\)\in\Fclosed{\X}{\topology{\X}},
\end{equation}
and hence according to \Ref{thmclosedmapproperty1peq1},
\begin{equation}\label{thmclosedmapproperty1p1-1eq4}
\func{\image{\om}}{\compl{\X}{\U}}\in\Fclosed{\Y}{\topology{\Y}},
\end{equation}
and thus according to \refdef{deffamilyofclosedsets},
\begin{equation}\label{thmclosedmapproperty1p1-1eq5}
\[\compl{\Y}{\func{\image{\om}}{\compl{\X}{\U}}}\]\in\topology{\Y}.
\end{equation}
According to \Ref{thmclosedmapproperty1p1-1eq1},
\begin{equation}\label{thmclosedmapproperty1p1-1eq6}
\(\compl{\X}{\U}\)\subseteq\[\compl{\X}{\func{\pimage{\om}}{\Asubset{\Y}}}\],
\end{equation}
and thus,
\begin{equation}\label{thmclosedmapproperty1p1-1eq7}
\func{\image{\om}}{\compl{\X}{\U}}\subseteq
\func{\image{\om}}{\compl{\X}{\func{\pimage{\om}}{\Asubset{\Y}}}},
\end{equation}
and thus,
\begin{equation}\label{thmclosedmapproperty1p1-1eq8}
\[\compl{\Y}{\func{\image{\om}}{\compl{\X}{\func{\pimage{\om}}{\Asubset{\Y}}}}}\]
\subseteq
\[\compl{\Y}{\func{\image{\om}}{\compl{\X}{\U}}}\].
\end{equation}
Additionally, it can be easily seen that,
\begin{equation}\label{thmclosedmapproperty1p1-1eq9}
\Asubset{\Y}\subseteq
\[\compl{\Y}{\func{\image{\om}}{\compl{\X}{\func{\pimage{\om}}{\Asubset{\Y}}}}}\].
\end{equation}
Thus according to \Ref{thmclosedmapproperty1p1-1eq8},
and
\Ref{thmclosedmapproperty1p1-1eq9},
\begin{equation}\label{thmclosedmapproperty1p1-1eq10}
\Asubset{\Y}\subseteq
\[\compl{\Y}{\func{\image{\om}}{\compl{\X}{\U}}}\].
\end{equation}
In addition,
\begin{equation}\label{thmclosedmapproperty1p1-1eq11}
\func{\pimage{\om}}{\compl{\Y}{\func{\image{\om}}{\compl{\X}{\U}}}}=
\compl{\X}{\func{\pimage{\om}}{\func{\image{\om}}{\compl{\X}{\U}}}}.
\end{equation}
Thus considering that,
\begin{equation}\label{thmclosedmapproperty1p1-1eq12}
\func{\pimage{\om}}{\func{\image{\om}}{\compl{\X}{\U}}}
\supseteq\(\compl{\X}{\U}\),
\end{equation}
it is evident that,
\begin{align}\label{thmclosedmapproperty1p1-1eq13}
\U&=\[\compl{\X}{\(\compl{\X}{\U}\)}\]\cr
&\supseteq\compl{\X}{\func{\pimage{\om}}{\func{\image{\om}}{\compl{\X}{\U}}}}\cr
&=\func{\pimage{\om}}{\compl{\Y}{\func{\image{\om}}{\compl{\X}{\U}}}}.
\end{align}
\Ref{thmclosedmapproperty1p1-1eq2},
\Ref{thmclosedmapproperty1p1-1eq5},
\Ref{thmclosedmapproperty1p1-1eq10},
and
\Ref{thmclosedmapproperty1p1-1eq13}
imply,
\begin{equation}
\Exists{\V}{\defset{\V}{\topology{\Y}}{\V\supseteq\Asubset{\Y}}}
\func{\pimage{\om}}{\V}\subseteq\U.
\end{equation}
\endp
\end{itemize}
\endp
\end{itemize}
\endthm
%%%%%%%%%%%%%%%%%%%%%%%%%%%%%%%%%%%%%%%%%%%%%%%%%%%%%%%%%%%%%%%%%%%%%%%%%%%%%%%%%%%%%%%
\theorem\label{thmopenmapproperty1}
Each
$\Xt=\opair{\X}{\topology{\X}}$
and
$\Yt=\opair{\Y}{\topology{\Y}}$
is taken as a topological-space, and $\om$ as an element of $\OM{\Xt}{\Yt}$
(an open map from $\Xt$ to $\Yt$).
\begin{align}
&\Foreach{\Asubset{\Y}}{\CSs{\Y}}\cr
&\bigg[\Foreach{\U}{\defset{\U}{\Fclosed{\X}{\topology{\X}}}
{\U\supseteq\func{\pimage{\om}}{\Asubset{\Y}}}}\cr
&~~\bigg(\Exists{\V}{\defset{\V}{\Fclosed{\Y}{\topology{\Y}}}{\V\supseteq\Asubset{\Y}}}
\func{\pimage{\om}}{\V}\subseteq\U\bigg)\bigg].\cr
&{}
\end{align}
\prooff
Based on \refdef{defopenmap},
\begin{equation}\label{thmopenmapproperty1peq1}
\Foreach{\U}{\topology{\X}}
\func{\image{\om}}{\U}\in\topology{\Y}.
\end{equation}
\begin{itemize}
\item[${\textbf{\textsf{p1}}}$]
$\Asubset{\Y}$
is taken as an arbitrary element of $\CSs{\Y}$.
\begin{itemize}
\item[${\textbf{\textsf{p1-1}}}$]
$\U$
is taken as such an element of $\Fclosed{\X}{\topology{\X}}$ that,
\begin{equation}\label{thmopenmapproperty1p1-1eq1}
\U\supseteq\func{\pimage{\om}}{\Asubset{\Y}}.
\end{equation}
$\V$
is defined as,
\begin{equation}\label{thmopenmapproperty1p1-1eq2}
\V:=\compl{\Y}{\func{\image{\om}}{\compl{\X}{\U}}}.
\end{equation}
Considering that $\U\in\Fclosed{\X}{\topology{\X}}$, and based on
\refdef{deffamilyofclosedsets},
\begin{equation}\label{thmopenmapproperty1p1-1eq3}
\(\compl{\X}{\U}\)\in\topology{\X},
\end{equation}
and hence accordng to \Ref{thmopenmapproperty1peq1},
\begin{equation}\label{thmopenmapproperty1p1-1eq4}
\func{\image{\om}}{\compl{\X}{\U}}\in\topology{\Y},
\end{equation}
and thus according to \refdef{deffamilyofclosedsets},
\begin{equation}\label{thmopenmapproperty1p1-1eq5}
\[\compl{\Y}{\func{\image{\om}}{\compl{\X}{\U}}}\]\in\Fclosed{\Y}{\topology{\Y}}.
\end{equation}
According to \Ref{thmopenmapproperty1p1-1eq1},
\begin{equation}\label{thmopenmapproperty1p1-1eq6}
\(\compl{\X}{\U}\)\subseteq\[\compl{\X}{\func{\pimage{\om}}{\Asubset{\Y}}}\],
\end{equation}
and hence,
\begin{equation}\label{thmopenmapproperty1p1-1eq7}
\func{\image{\om}}{\compl{\X}{\U}}\subseteq
\func{\image{\om}}{\compl{\X}{\func{\pimage{\om}}{\Asubset{\Y}}}},
\end{equation}
and thus,
\begin{equation}\label{thmopenmapproperty1p1-1eq8}
\[\compl{\Y}{\func{\image{\om}}{\compl{\X}{\func{\pimage{\om}}{\Asubset{\Y}}}}}\]
\subseteq
\[\compl{\Y}{\func{\image{\om}}{\compl{\X}{\U}}}\].
\end{equation}
In addition, it can be easily seen that,
\begin{equation}\label{thmopenmapproperty1p1-1eq9}
\Asubset{\Y}\subseteq
\[\compl{\Y}{\func{\image{\om}}{\compl{\X}{\func{\pimage{\om}}{\Asubset{\Y}}}}}\].
\end{equation}
Thus according to \Ref{thmopenmapproperty1p1-1eq8} and
\Ref{thmopenmapproperty1p1-1eq9},
\begin{equation}\label{thmopenmapproperty1p1-1eq10}
\Asubset{\Y}\subseteq
\[\compl{\Y}{\func{\image{\om}}{\compl{\X}{\U}}}\].
\end{equation}
In addition,
\begin{equation}\label{thmopenmapproperty1p1-1eq11}
\func{\pimage{\om}}{\compl{\Y}{\func{\image{\om}}{\compl{\X}{\U}}}}=
\compl{\X}{\func{\pimage{\om}}{\func{\image{\om}}{\compl{\X}{\U}}}}.
\end{equation}
Thus considering that,
\begin{equation}\label{thmopenmapproperty1p1-1eq12}
\func{\pimage{\om}}{\func{\image{\om}}{\compl{\X}{\U}}}
\supseteq\(\compl{\X}{\U}\),
\end{equation}
it is evident that,
\begin{align}\label{thmopenmapproperty1p1-1eq13}
\U&=\[\compl{\X}{\(\compl{\X}{\U}\)}\]\cr
&\supseteq\compl{\X}{\func{\pimage{\om}}{\func{\image{\om}}{\compl{\X}{\U}}}}\cr
&=\func{\pimage{\om}}{\compl{\Y}{\func{\image{\om}}{\compl{\X}{\U}}}}.
\end{align}
\Ref{thmopenmapproperty1p1-1eq2},
\Ref{thmopenmapproperty1p1-1eq5},
\Ref{thmopenmapproperty1p1-1eq10},
and
\Ref{thmopenmapproperty1p1-1eq13}
imply,
\begin{equation}
\Exists{\V}{\defset{\V}{\Fclosed{\Y}{\topology{\Y}}}{\V\supseteq\Asubset{\Y}}}
\func{\pimage{\om}}{\V}\subseteq\U.
\end{equation}
\endp
\end{itemize}
\endp
\end{itemize}
\endthm
%%%%%%%%%%%%%%%%%%%%%%%%%%%%%%%%%%%%%%%%%%%%%%%%%%%%%%%%%%%%%%%%%%%%%%%%%%%%%%%%%%%%%%%
\theorem\label{thmopenmapandinterior}
Each
$\Xt=\opair{\X}{\topology{\X}}$
and
$\Yt=\opair{\Y}{\topology{\Y}}$
is taken as a topological-space.
\begin{align}
&\Foreach{\om}{\Func{\X}{\Y}}\cr
&\bigg(\om\in\OM{\Xt}{\Yt}\thenn
\bigg[\Foreach{\Asubset{\X}}{\CSs{\X}}
\func{\image{\om}}{\func{\Int{\Xt}}{\Asubset{\X}}}\subseteq
\func{\Int{\Yt}}{\func{\image{\om}}{\Asubset{\X}}}\bigg]\bigg).\cr
&{}
\end{align}
\prooff
$\om$
is taken as an arbitrary element of $\Func{\X}{\Y}$.
\begin{itemize}
\item[${\textbf{\textsf{p1}}}$]
It is assumed that,
\begin{equation}\label{thmopenmapandinteriorp1eq1}
\om\in\OM{\Xt}{\Yt}.
\end{equation}
Then based on \refdef{defopenmap},
and considering that(\refcor{corintofset0}),
\begin{equation}\label{thmopenmapandinteriorp1eq2}
\Foreach{\Asubset{\X}}{\CSs{\X}}
\func{\Int{\Xt}}{\Asubset{\X}}\in\topology{\X},
\end{equation}
it is clear that,
\begin{equation}\label{thmopenmapandinteriorp1eq3}
\Foreach{\Asubset{\X}}{\CSs{\X}}
\func{\image{\om}}{\func{\Int{\Xt}}{\Asubset{\X}}}\in\topology{\Y}.
\end{equation}
In addition, according to \refcor{corintofset0},
\begin{equation}\label{thmopenmapandinteriorp1eq4}
\Foreach{\Asubset{\X}}{\CSs{\X}}
\func{\Int{\Xt}}{\Asubset{\X}}\subseteq\Asubset{\X},
\end{equation}
and hence,
\begin{equation}\label{thmopenmapandinteriorp1eq5}
\Foreach{\Asubset{\X}}{\CSs{\X}}
\func{\image{\om}}{\func{\Int{\Xt}}{\Asubset{\X}}}\subseteq
\func{\image{\om}}{\Asubset{\X}}.
\end{equation}
Based on \refcor{corintofset0},
\Ref{thmopenmapandinteriorp1eq3}
and
\Ref{thmopenmapandinteriorp1eq5}
imply that,
\begin{equation*}
\Foreach{\Asubset{\X}}{\CSs{\X}}
\func{\image{\om}}{\func{\Int{\Xt}}{\Asubset{\X}}}\subseteq
\func{\Int{\Yt}}{\func{\image{\om}}{\Asubset{\X}}}.
\end{equation*}
\endp
\end{itemize}
\begin{itemize}
\item[${\textbf{\textsf{p2}}}$]
It is assumed that,
\begin{equation}\label{thmopenmapandinteriorp2eq1}
\Foreach{\Asubset{\X}}{\CSs{\X}}
\func{\image{\om}}{\func{\Int{\Xt}}{\Asubset{\X}}}\subseteq
\func{\Int{\Yt}}{\func{\image{\om}}{\Asubset{\X}}}.
\end{equation}
Based on \refthm{thmintofopenset},
\begin{equation}\label{thmopenmapandinteriorp2eq2}
\Foreach{\U}{\topology{\X}}
\func{\Int{\Xt}}{\U}=\U.
\end{equation}
\Ref{thmopenmapandinteriorp2eq1}
and
\Ref{thmopenmapandinteriorp2eq2}
imply,
\begin{equation}
\Foreach{\U}{\topology{\X}}
\func{\image{\om}}{\U}\subseteq
\func{\Int{\Yt}}{\func{\image{\om}}{\U}}.
\end{equation}
Based on \refcor{corintofset0}, this means,
\begin{equation}
\Foreach{\U}{\topology{\X}}
\func{\image{\om}}{\U}=
\func{\Int{\Yt}}{\func{\image{\om}}{\U}},
\end{equation}
and hence according to \refthm{thmintofopenset},
\begin{equation}
\Foreach{\U}{\topology{\X}}
\func{\image{\om}}{\U}\in\topology{\Y}.
\end{equation}
According to \refdef{defopenmap},
this means,
\begin{equation*}
\om\in\OM{\Xt}{\Yt}.
\end{equation*}
\endp
\end{itemize}
\endthm
%%%%%%%%%%%%%%%%%%%%%%%%%%%%%%%%%%%%%%%%%%%%%%%%%%%%%%%%%%%%%%%%%%%%%%%%%%%%%%%%%%%%%%%
\theorem\label{thmopenmapandbasis}
Each
$\Xt=\opair{\X}{\topology{\X}}$
and
$\Yt=\opair{\Y}{\topology{\Y}}$
is taken as a topological-space.
For a base of $\Xt$, the set of all functions from $\X$ to $\Y$ that send
every element of the base to an open set of $\Yt$, equals the set of aa open maps from
$\Xt$ to $\Yt$. That is,
\begin{align}
\Foreach{\baseof{\Xt}}{\Cbase{\Xt}}
\defset{\om}{\Func{\X}{\Y}}
{\bigg[\Foreach{\B}{\baseof{\Xt}}
\func{\image{\om}}{\B}\in\topology{\Y}\bigg]}=\OM{\X}{\Y}.
\end{align}
\marginpar{\rotatebox{90}{\fbox{$
\defset{\baseof{\Xt}}{\Cbase{\Xt}}
{\defset{\om}{\Func{\X}{\Y}}
{\bigg[\defset{\B}{\baseof{\Xt}}
{\bigg[\func{\image{\om}}{\B}\in\topology{\Y}\bigg]}=\baseof{\Xt}\bigg]}=\OM{\Xt}{\Yt}}=\baseof{\Xt}.
$}}}
\prooff
$\baseof{\Xt}$
is taken as an arbitrary element of $\Cbase{\Xt}$
(a base for the topological-space $\Xt$).
\begin{itemize}
\item[${\textbf{\textsf{p1}}}$]
$\om$
is taken as such an element of $\Func{\X}{\Y}$ that,
\begin{equation}\label{thmopenmapandbasisp1eq1}
\Foreach{\B}{\baseof{\Xt}}
\func{\image{\om}}{\B}\in\topology{\Y}.
\end{equation}
\begin{itemize}
\item[${\textbf{\textsf{p1-1}}}$]
$\U$
is taken as an arbitrary element of $\topology{\X}$. Then according to \refdef{defbase},
\begin{equation}\label{thmopenmapandbasisp1-1eq1}
\Existsis{\sC}{\CSs{\baseof{\Yt}}}
\U=\union{\sC}.
\end{equation}
\Ref{thmopenmapandbasisp1eq1}
implies,
\begin{equation}\label{thmopenmapandbasisp1-1eq2}
\defset{\Asubset{\Y}}{\CSs{\Y}}
{\big[\Exists{\V}{\sC}\Asubset{\Y}=\func{\image{\om}}{\V}\big]}\subseteq\topology{\Y},
\end{equation}
and hence according to \refdef{deftopologicalspace},
\begin{equation}\label{thmopenmapandbasisp1-1eq3}
\union{\defset{\Asubset{\Y}}{\CSs{\Y}}
{\big[\Exists{\V}{\sC}\Asubset{\Y}=\func{\image{\om}}{\V}\big]}}
\in\topology{\Y}.
\end{equation}
In addition, according to \Ref{thmopenmapandbasisp1-1eq1},
\begin{align}\label{thmopenmapandbasisp1-1eq4}
\func{\image{\om}}{\U}&=
\func{\image{\om}}{\union{\sC}}\cr
&=\union{\defset{\Asubset{\Y}}{\CSs{\Y}}
{\big[\Exists{\V}{\sC}\Asubset{\Y}=\func{\image{\om}}{\V}\big]}}.\cr
&{}
\end{align}
\Ref{thmopenmapandbasisp1-1eq3} and \Ref{thmopenmapandbasisp1-1eq4} imply,
\begin{equation}
\func{\image{\om}}{\U}\in\topology{\Y}.
\end{equation}
\endp
\end{itemize}
Thus,
\begin{equation}
\Foreach{\U}{\topology{\X}}
\func{\image{\om}}{\U}\in\topology{\Y},
\end{equation}
which according to \refdef{defopenmap}, means,
\begin{equation}
\om\in\OM{\Xt}{\Yt}.
\end{equation}
\endp
\end{itemize}
\begin{itemize}
\item[${\textbf{\textsf{p2}}}$]
$\om$
is taken as an arbitrary element of $\OM{\Xt}{\Yt}$. Then according to \refdef{defopenmap},
and considering that $\baseof{\Xt}\subseteq\topology{\X}$, $\om$ sends every open set of $\Xt$,
and hence every element of $\baseof{\Xt}$, to an open set of $\Yt$.
\endp
\end{itemize}
\endthm
%%%%%%%%%%%%%%%%%%%%%%%%%%%%%%%%%%%%%%%%%%%%%%%%%%%%%%%%%%%%%%%%%%%%%%%%%%%%%%%%%%%%%%%
\theorem\label{thmopenmapequiv0}
Each
$\Xt=\opair{\X}{\topology{\X}}$
and
$\Yt=\opair{\Y}{\topology{\Y}}$
is taken as a topological-space, and $\om$ as an element of $\Func{\X}{\Y}$.
\begin{gather}
\om\in\OM{\Xt}{\Yt}\cr
\vthenn\cr
\Foreach{\point}{\X}
\[\Foreach{\U}{\func{\nei{\Xt}}{\seta{\point}}}
\bigg(\Exists{\V}{\func{\nei{\Yt}}{\seta{\func{\om}{\point}}}}
\V\subseteq\func{\image{\om}}{\U}\bigg)\].
\end{gather}
\prooff
\begin{itemize}
\item[${\textbf{\textsf{p1}}}$]
It is assumed that,
\begin{equation}
\om\in\OM{\Xt}{\Yt}.
\end{equation}
That is,
\begin{equation}
\Foreach{\U}{\topology{\X}}
\func{\image{\om}}{\U}\in\topology{\Y}.
\end{equation}
This and \refdef{defnbdclassofsets} imply,
\begin{equation}
\Foreach{\point}{\X}
\[\Foreach{\U}{\func{\nei{\Xt}}{\seta{\point}}}
\func{\image{\om}}{\U}\in\func{\nei{\Yt}}{\seta{\func{\om}{\point}}}\].
\end{equation}
This directly implies,
\begin{equation*}
\Foreach{\point}{\X}
\[\Foreach{\U}{\func{\nei{\Xt}}{\seta{\point}}}
\bigg(\Exists{\V}{\func{\nei{\Yt}}{\seta{\func{\om}{\point}}}}
\V\subseteq\func{\image{\om}}{\U}\bigg)\].
\end{equation*}
\endp
\end{itemize}
\begin{itemize}
\item[${\textbf{\textsf{p2}}}$]
It is assuemd that,
\begin{align}\label{thmopenmapequiv0p2eq1}
\Foreach{\point}{\X}
\[\Foreach{\U}{\func{\nei{\Xt}}{\seta{\point}}}
\bigg(\Exists{\V}{\func{\nei{\Yt}}{\seta{\func{\om}{\point}}}}
\V\subseteq\func{\image{\om}}{\U}\bigg)\].
\end{align}
\begin{itemize}
\item[${\textbf{\textsf{p2-1}}}$]
$\U$
is taken as an arbitrary element of $\topology{\X}$. According to \refdef{defnbdclassofsets},
\begin{equation}\label{thmopenmapequiv0p2-1eq1}
\Foreach{\point}{\U}
\U\in\func{\nei{\Xt}}{\seta{\point}}.
\end{equation}
\Ref{thmopenmapequiv0p2eq1}
and
\Ref{thmopenmapequiv0p2-1eq1}
imply,
\begin{equation}
\Foreach{\point}{\U}
\bigg(\Exists{\V}{\func{\nei{\Yt}}{\seta{\func{\om}{\point}}}}
\V\subseteq\func{\image{\om}}{\U}\bigg).
\end{equation}
$\eta$
is defined as,
\begin{equation}
\left\{
\begin{split}
&\eta\in\Func{\U}{\topology{\Y}},\cr
&\Foreach{\point}{\U}\func{\om}{\point}\in
\func{\eta}{\point}\subseteq\func{\image{\om}}{\U}.
\end{split}
\right.
\end{equation}
This and \refdef{deftopologicalspace} imply,
\begin{equation}
\[\Union{\point}{\U}{\func{\eta}{\point}}\]\in\topology{\Y}.
\end{equation}
In addition, it can be easily seen that,
\begin{equation}
\func{\image{\om}}{\U}=
\[\Union{\point}{\U}{\func{\eta}{\point}}\].
\end{equation}
Thus,
\begin{equation}
\func{\image{\om}}{\U}\in\topology{\Y}.
\end{equation}
\endp
\end{itemize}
Therefore,
\begin{equation}
\Foreach{\U}{\topology{\X}}
\func{\image{\om}}{\U}\in\topology{\Y},
\end{equation}
which based on \refdef{defopenmap} means,
\begin{equation*}
\om\in\OM{\Xt}{\Yt}.
\end{equation*}
\endp
\end{itemize}
\endthm
%%%%%%%%%%%%%%%%%%%%%%%%%%%%%%%%%%%%%%%%%%%%%%%%%%%%%%%%%%%%%%%%%%%%%%%%%%%%%%%%%%%%%%%
\theorem\label{thmclosedmapandclosure}
Each
$\Xt=\opair{\X}{\topology{\X}}$
and
$\Yt=\opair{\Y}{\topology{\Y}}$
is taken as a topological-space, and $\om$ as an element of $\Func{\X}{\Y}$.
\begin{align}
\[\om\in\CM{\Xt}{\Yt}\]\thenn
\bigg[\Foreach{\Asubset{\X}}{\CSs{\X}}
\func{\Cl{\Yt}}{\func{\image{\om}}{\Asubset{\X}}}\subseteq
\func{\image{\om}}{\func{\Cl{\Xt}}{\Asubset{\X}}}\bigg].
\end{align}
\prooff
\begin{itemize}
\item[${\textbf{\textsf{p1}}}$]
It is assumed that,
\begin{equation}\label{thmclosedmapandclosurep1eq1}
\om\in\CM{\Xt}{\Yt}.
\end{equation}
That is,
\begin{equation}\label{thmclosedmapandclosurep1eq2}
\Foreach{\U}{\Fclosed{\X}{\topology{\X}}}
\func{\image{\om}}{\U}\in\Fclosed{\Y}{\topology{\Y}}.
\end{equation}
According to \refcor{corclosureofset0},
\begin{equation}\label{thmclosedmapandclosurep1eq3}
\Foreach{\Asubset{\X}}{\CSs{\X}}
\bigg[\func{\Cl{\Xt}}{\Asubset{\X}}\in\Fclosed{\X}{\topology{\X}}\bigg]
\end{equation}
According to \Ref{thmclosedmapandclosurep1eq2} and \Ref{thmclosedmapandclosurep1eq3},
\begin{equation}\label{thmclosedmapandclosurep1eq4}
\Foreach{\Asubset{\X}}{\CSs{\X}}
\func{\image{\om}}{\func{\Cl{\Xt}}{\Asubset{\X}}}\in\Fclosed{\Y}{\topology{\Y}}.
\end{equation}
In addition, considering that
$\Foreach{\Asubset{\X}}{\CSs{\X}}\Asubset{\X}\subseteq\func{\Cl{\Xt}}{\Asubset{\X}}$,
it is evident that,
\begin{equation}\label{thmclosedmapandclosurep1eq5}
\Foreach{\Asubset{\X}}{\CSs{\X}}
\func{\image{\om}}{\Asubset{\X}}
\subseteq\func{\image{\om}}{\func{\Cl{\Xt}}{\Asubset{\X}}}.
\end{equation}
Based on \refcor{corclosureofset0},
\Ref{thmclosedmapandclosurep1eq4} and \Ref{thmclosedmapandclosurep1eq5} imply,
\begin{equation*}
\Foreach{\Asubset{\X}}{\CSs{\X}}
\func{\Cl{\Yt}}{\func{\image{\om}}{\Asubset{\X}}}\subseteq
\func{\image{\om}}{\func{\Cl{\Xt}}{\Asubset{\X}}}.
\end{equation*}
\endp
\end{itemize}
\begin{itemize}
\item[${\textbf{\textsf{p2}}}$]
It is assumed that,
\begin{equation}\label{thmclosedmapandclosurep2eq1}
\Foreach{\Asubset{\X}}{\CSs{\X}}
\func{\Cl{\Yt}}{\func{\image{\om}}{\Asubset{\X}}}\subseteq
\func{\image{\om}}{\func{\Cl{\Xt}}{\Asubset{\X}}}.
\end{equation}
According to \refthm{thmclosureofclosedset},
\begin{equation}\label{thmclosedmapandclosurep2eq2}
\Foreach{\Asubset{\X}}{\Fclosed{\X}{\topology{\X}}}
\func{\Cl{\Xt}}{\Asubset{\X}}=\Asubset{\X},
\end{equation}
and hence according to \Ref{thmclosedmapandclosurep2eq1},
\begin{equation}
\Foreach{\Asubset{\X}}{\Fclosed{\X}{\topology{\X}}}
\func{\Cl{\Yt}}{\func{\image{\om}}{\Asubset{\X}}}\subseteq
\func{\image{\om}}{\Asubset{\X}},
\end{equation}
and thus according to \refcor{corclosureofset0},
it is clear that,
\begin{equation}
\Foreach{\Asubset{\X}}{\Fclosed{\X}{\topology{\X}}}
\func{\Cl{\Yt}}{\func{\image{\om}}{\Asubset{\X}}}=
\func{\image{\om}}{\Asubset{\X}}.
\end{equation}
On the basis of \refthm{thmclosureofclosedset}, this means,
\begin{equation}
\Foreach{\Asubset{\X}}{\Fclosed{\X}{\topology{\X}}}
\func{\image{\om}}{\Asubset{\X}}\in\Fclosed{\Y}{\topology{\Y}},
\end{equation}
which according to \refdef{defclosedmap}, means,
\begin{equation*}
\om\in\CM{\Xt}{\Yt}.
\end{equation*}
\endp
\end{itemize}
\endthm
%%%%%%%%%%%%%%%%%%%%%%%%%%%%%%%%%%%%%%%%%%%%%%%%%%%%%%%%%%%%%%%%%%%%%%%%%%%%%%%%%%%%%%%%%%%%%%%%%%%%%%%%%%%%%%%%%%%%%%%%
%%%%%%%%%%%%%%%%%%%%%%%%%%%%%%%%%%%%%%%%%%%%%%%%%%%%%%%%%%%%%%%%%%%%%%%%%%%%%%%%%%%%%%%%%%%%%%%%%%%%%%%%%%%%%%%%%%%%%%%%
%%%%%%%%%%%%%%%%%%%%%%%%%%%%%%%%%%%%%%%%%%%%%%%%%%%%%%%%%%%%%%%%%%%%%%%%%%%%%%%%%%%%%%%%%%%%%%%%%%%%%%%%%%%%%%%%%%%%%%%%
%%%%%%%%%%%%%%%%%%%%%%%%%%%%%%%%%%%%%%%%%%%%%%%%%%%%%%%%%%%%%%%%%%%%%%%%%%%%%%%%%%%%%%%%%%%%%%%%%%%%%%%%%%%%%%%%%%%%%%%%
%%%%%%%%%%%%%%%%%%%%%%%%%%%%%%%%%%%%%%%%%%%%%%%%%%%%%%%%%%%%%%%%%%%%%%%%%%%%%%%%%%%%%%%%%%%%%%%%%%%%%%%%%%%%%%%%%%%%%%%%
%%%%%%%%%%%%%%%%%%%%%%%%%%%%%%%%%%%%%%%%%%%%%%%%%%%%%%%%%%%%%%%%%%%%%%%%%%%%%%%%%%%%%%%%%%%%%%%%%%%%%%%%%%%%%%%%%%%%%%%%
%%%%%%%%%%%%%%%%%%%%%%%%%%%%%%%%%%%%%%%%%%%%%%%%%%%%%%%%%%%%%%%%%%%%%%%%%%%%%%%%%%%%%%%%%%%%%%%%%%%%%%%%%%%%%%%%%%%%%%%%
\section{Homeomorphisms}
\definition\label{defhomeomorphism}
Each
$\Xt=\opair{\X}{\topology{\X}}$
and
$\Yt=\opair{\Y}{\topology{\Y}}$
is taken as a topological-space.
For every $\hf$ in $\Func{\X}{\Y}$,
$\hf$
is reffered to as a $\quotl$homeomorphism from (the topological-space) $\Xt$
to (the topological-space) $\Yt$$\quotr$ iff these properties are hold.
\begin{itemize}
\item[${\textbf{\textsf{HF1}}}$]
$\hf$
is a bijection from $\X$ to $\Y$.
\hfill
$\hf\in\IF{\X}{\Y}.$
\item[${\textbf{\textsf{HF2}}}$]
$\hf$
is a continuous map from $\Xt$ to $\Yt$.
\hfill
$\hf\in\CF{\Xt}{\Yt}.$
\item[${\textbf{\textsf{HF3}}}$]
$\finv{\hf}$
(the inverse of $\hf$)
is a continuous map from $\Yt$ to $\Xt$.
\hfill
$\finv{\hf}\in\CF{\Yt}{\Xt}.$
\end{itemize}
The set of all homeomorphisms from $\Xt$ to $\Yt$ is denoted by $\HOF{\Xt}{\Yt}$. That is,
\begin{align}
\HOF{\Xt}{\Yt}:=\defset{\hf}{\IF{\X}{\Y}}
{\[\AND{\bigg(\hf\in\CF{\Xt}{\Yt}\bigg)}{\bigg(\finv{\hf}\in\CF{\Yt}{\Xt}\bigg)}\]}.
\end{align}
\endef
%%%%%%%%%%%%%%%%%%%%%%%%%%%%%%%%%%%%%%%%%%%%%%%%%%%%%%%%%%%%%%%%%%%%%%%%%%%%%%%%%%%%%%%
\theorem\label{thmtrivialhomeomorphismofaspace}
$\Xt=\opair{\X}{\topology{}}$
is taken as a topological-space.
The identity-function on $\X$ is a homeomorphism from $\Xt$ to $\Xt$. That is,
\begin{equation}\label{thmtrivialhomeomorphismofaspaceeq1}
\idf{\X}\in\HOF{\Xt}{\Xt}.
\end{equation}
\prooff
It is known that,
\begin{equation}
\idf{\X}\in\IF{\X}{\X},
\end{equation}
and,
\begin{equation}
\finv{\idf{\X}}=\idf{\X}.
\end{equation}
So based on \refthm{thmtrivialcontinuousmapofaspace},
\begin{equation}
\opair{\idf{\X}}{\finv{\idf{\X}}}\in
\[\CF{\Xt}{\Xt}\times\CF{\Xt}{\Xt}\].
\end{equation}
Thus based on \refdef{defhomeomorphism},
\Ref{thmtrivialhomeomorphismofaspaceeq1}
is clear.
\endthm
%%%%%%%%%%%%%%%%%%%%%%%%%%%%%%%%%%%%%%%%%%%%%%%%%%%%%%%%%%%%%%%%%%%%%%%%%%%%%%%%%%%%%%%
\theorem\label{thminverseofhomeomorphism}
Each
$\Xt=\opair{\X}{\topology{\X}}$
and
$\Yt=\opair{\Y}{\topology{\Y}}$
is taken as a topological-space.
The inverse of every homeomorphism from $\Xt$ to $\Yt$
is a homeomorphism from $\Yt$ to $\Xt$. That is,
\begin{equation}\label{thminverseofhomeomorphismpeq1}
\Foreach{\hf}{\HOF{\Xt}{\Yt}}
\finv{\hf}\in\HOF{\Yt}{\Xt}.
\end{equation}
\prooff
$\hf$
is taken as an arbitrary element of $\HOF{\Xt}{\Yt}$. Then according to \refdef{defhomeomorphism},
\begin{align}
&\hf\in\IF{\X}{\Y},\label{thminverseofhomeomorphismpeq2}\\
&\opair{\hf}{\finv{\hf}}\in
\[\CF{\Xt}{\Yt}\times\CF{\Yt}{\Xt}\].\label{thminverseofhomeomorphismpeq3}
\end{align}
In addition, it is known that,
\begin{align}
&\finv{\hf}\in\IF{\Y}{\X},\\
&\finv{\(\finv{\hf}\)}=\hf.
\end{align}
Thus according to \Ref{thminverseofhomeomorphismpeq3},
\begin{equation}
\opair{\finv{\hf}}{\finv{\(\finv{\hf}\)}}\in
\[\CF{\Yt}{\Xt}\times\CF{\Xt}{\Yt}\],
\end{equation}
and hence according to \refdef{defhomeomorphism},
\begin{equation}
\finv{\hf}\in\HOF{\Yt}{\Xt}.
\end{equation}
\endthm
%%%%%%%%%%%%%%%%%%%%%%%%%%%%%%%%%%%%%%%%%%%%%%%%%%%%%%%%%%%%%%%%%%%%%%%%%%%%%%%%%%%%%%%
\theorem\label{thmcompositionofhomeomorphisms}
Each
$\Xt_{1}=\opair{\X_{1}}{\topology{1}}$,
$\Xt_2=\opair{\X_2}{\topology{2}}$,
and
$\Xt_3=\opair{\X_3}{\topology{3}}$
is taken as a topological-space
The composition of every homeomorphism from
$\Xt_2$
to
$\Xt_3$
with every homeomorphism from
$\Xt_1$
to
$\Xt_2$
is a homeomorphism from
$\Xt_1$
to
$\Xt_3$. That is,
\begin{equation}
\Foreach{\opair{\hf}{\p{\hf}}}
{\[\HOF{\Xt_{1}}{\Xt_{2}}\times\HOF{\Xt_{2}}{\Xt_{3}}\]}
\(\cmp{\p{\hf}}{\hf}\)\in\HOF{\Xt_{1}}{\Xt_{3}}.
\end{equation}
\prooff
$\hf$
is taken as an arbitrary element of
$\HOF{\Xt_1}{\Xt_2}$,
and
$\p{\hf}$
an arbitrary element of
$\HOF{\Xt_2}{\Xt_3}$. Then according to \refdef{defhomeomorphism},
\begin{align}
\hf&\in\CF{\Xt_1}{\Xt_2}\label{thmcompositionofhomeomorphismspeq1}\\
\finv{\hf}&\in\CF{\Xt_2}{\Xt_1}\label{thmcompositionofhomeomorphismspeq2},
\end{align}
and,
\begin{align}
\p{\hf}&\in\CF{\Xt_2}{\Xt_3}\label{thmcompositionofhomeomorphismspeq3}\\
\finv{\p{\hf}}&\in\CF{\Xt_3}{\Xt_2}.\label{thmcompositionofhomeomorphismspeq4}
\end{align}
Therefore, according to \refthm{thmcompositionofcontinuousfunctions},
\begin{align}
\cmp{\p{\hf}}{\hf}&\in\CF{\Xt_1}{\Xt_3},\label{thmcompositionofhomeomorphismspeq5}\\
\cmp{\finv{\hf}}{\finv{\p{\hf}}}&\in\CF{\Xt_3}{\Xt_1}.\label{thmcompositionofhomeomorphismspeq6}
\end{align}
In addition, it is known that,
\begin{align}
&\(\cmp{\p{\hf}}{\hf}\)\in\IF{\X_1}{\X_3},\label{thmcompositionofhomeomorphismspeq7}\\
&\finv{\(\cmp{\p{\hf}}{\hf}\)}=\cmp{\finv{\hf}}{\finv{\p{\hf}}}.\label{thmcompositionofhomeomorphismspeq8}
\end{align}
According to \Ref{thmcompositionofhomeomorphismspeq6}
and
\Ref{thmcompositionofhomeomorphismspeq8},
it is clear that,
\begin{equation}\label{thmcompositionofhomeomorphismspeq9}
\finv{\(\cmp{\p{\hf}}{\hf}\)}\in\CF{\Xt_3}{\Xt_1}.
\end{equation}
According to \refdef{defhomeomorphism},
\Ref{thmcompositionofhomeomorphismspeq5}
and
\Ref{thmcompositionofhomeomorphismspeq9}
imply that,
\begin{equation}
\(\cmp{\p{\hf}}{\hf}\)\in\HOF{\Xt_1}{\Xt_3}.
\end{equation}
\endthm
%%%%%%%%%%%%%%%%%%%%%%%%%%%%%%%%%%%%%%%%%%%%%%%%%%%%%%%%%%%%%%%%%%%%%%%%%%%%%%%%%%%%%%%
\definition\label{defhomeomorphic}
Each
$\Xt=\opair{\X}{\topology{\X}}$
and
$\Yt=\opair{\Y}{\topology{\Y}}$
is taken as a topological-space.
It is said that $\quotl$$\Xt$ is homeomorphic to $\Yt$$\quotr$,
denoted by
$\quotl$$\homeomorphic{\Xt}{\Yt}$$\quotr$
iff there exists a homeomorphism from
$\Xt$
to
$\Yt$.
In other words,
\begin{equation}
\(\homeomorphic{\Xt}{\Yt}\)\iffdef\[\HOF{\Xt}{\Yt}\neq\empty\].
\end{equation}
\endef
%%%%%%%%%%%%%%%%%%%%%%%%%%%%%%%%%%%%%%%%%%%%%%%%%%%%%%%%%%%%%%%%%%%%%%%%%%%%%%%%%%%%%%%
\corollary\label{corhomeomorphismequivalencerelation}
Each
$\Xt_{1}=\opair{\X_{1}}{\topology{1}}$,
$\Xt_2=\opair{\X_2}{\topology{2}}$,
and
$\Xt_3=\opair{\X_3}{\topology{3}}$
is taken as a topological-space.
\begin{equation}
\left\{
\begin{split}
&\homeomorphic{\Xt_1}{\Xt_1}.\cr
&\(\homeomorphic{\Xt_1}{\Xt_2}\)
\then\(\homeomorphic{\Xt_2}{\Xt_1}\).\cr
&\[\(\homeomorphic{\Xt_1}{\Xt_2}\),~\(\homeomorphic{\Xt_2}{\Xt_3}\)\]\then
\(\homeomorphic{\Xt_1}{\Xt_3}\).
\end{split}
\right.
\end{equation}
\endcor
%%%%%%%%%%%%%%%%%%%%%%%%%%%%%%%%%%%%%%%%%%%%%%%%%%%%%%%%%%%%%%%%%%%%%%%%%%%%%%%%%%%%%%%
\theorem\label{thmhomeomorphismsandopensets}
Each
$\Xt=\opair{\X}{\topology{\X}}$
and
$\Yt=\opair{\Y}{\topology{\Y}}$
is taken as a topological-space.
\begin{align}
\HOF{\Xt}{\Yt}=
\defset{\hf}{\IF{\X}{\Y}}{\bigg[
\defset{\U}{\CSs{\X}}{\func{\image{\hf}}{\U}\in\topology{\Y}}
=\topology{\X}
\bigg]}.
\end{align}
\prooff
\begin{itemize}
\item[${\textbf{\textsf{p1}}}$]
$\hf$
is taken as an arbitrary element of $\HOF{\Xt}{\Yt}$.
Then according to \refdef{defhomeomorphism},
\begin{align}
\hf&\in\IF{\X}{\Y},\label{thmhomeomorphismsandopensetsp1eq1}\\
\hf&\in\CF{\Xt}{\Yt},\label{thmhomeomorphismsandopensetsp1eq2}\\
\finv{\hf}&\in\CF{\Yt}{\Xt},\label{thmhomeomorphismsandopensetsp1eq3}
\end{align}
\begin{itemize}
\item[${\textbf{\textsf{p1-1}}}$]
$\U$
is taken as an arbitrary element of $\topology{\X}$. Then according to \Ref{thmhomeomorphismsandopensetsp1eq3}
and \refdef{defcontinuousfunction},
\begin{equation}
\func{\pimage{\(\finv{\hf}\)}}{\U}\in\topology{\Y}.
\end{equation}
It is known that,
\begin{equation}
\Foreach{\V}{\CSs{\X}}
\func{\pimage{\(\finv{\hf}\)}}{\V}=\func{\image{\hf}}{\V}.
\end{equation}
Thus,
\begin{equation}
\func{\image{\hf}}{\U}\in\topology{\Y}.
\end{equation}
\endp
\end{itemize}
\begin{itemize}
\item[${\textbf{\textsf{p2-1}}}$]
$\U$
is taken as such an arbitrary element of
$\CSs{\X}$ that,
\begin{equation}
\func{\image{\hf}}{\U}\in\topology{\Y}.
\end{equation}
Then according to \Ref{thmhomeomorphismsandopensetsp1eq2}
and
\refdef{defcontinuousfunction},
\begin{equation}
\func{\pimage{\hf}}{\func{\image{\hf}}{\U}}\in\topology{\X}.
\end{equation}
Additionally, according to \Ref{thmhomeomorphismsandopensetsp1eq1},
\begin{equation}
\Foreach{\V}{\CSs{\X}}
\func{\pimage{\hf}}{\func{\image{\hf}}{\V}}=\V.
\end{equation}
Therefore,
\begin{equation}
\U\in\topology{\X}.
\end{equation}
\endp
\end{itemize}
\endp
\end{itemize}
\begin{itemize}
\item[${\textbf{\textsf{p2}}}$]
$\hf$
is taken as such an arbitrary element of $\IF{\X}{\Y}$ such that,
\begin{equation}\label{thmhomeomorphismsandopensetsp2eq1}
\defset{\U}{\CSs{\X}}{\func{\image{\hf}}{\U}\in\topology{\Y}}
=\topology{\X}.
\end{equation}
Considering that $\hf$ is bijective,
\begin{equation}\label{thmhomeomorphismsandopensetsp2eq2}
\Foreach{\Asubset{\Y}}{\CSs{\Y}}
\func{\image{\hf}}{\func{\pimage{\hf}}{\Asubset{\Y}}}=\Asubset{\Y},
\end{equation}
and,
\begin{equation}\label{thmhomeomorphismsandopensetsp2eq3}
\Foreach{\Asubset{\X}}{\CSs{\X}}
\func{\pimage{\[\finv{\hf}\]}}{\Asubset{\X}}=\func{\image{\hf}}{\Asubset{\X}}.
\end{equation}
\begin{itemize}
\item[${\textbf{\textsf{p2-1}}}$]
$\U$
is taken as an arbitrary element of $\topology{\Y}$. Then according to \Ref{thmhomeomorphismsandopensetsp2eq2},
\begin{equation}
\func{\image{\hf}}{\func{\pimage{\hf}}{\U}}\in\topology{\Y},
\end{equation}
and hence according to \Ref{thmhomeomorphismsandopensetsp2eq1},
\begin{equation}
\func{\pimage{\hf}}{\U}\in\topology{\X}.
\end{equation}
\endp
\end{itemize}
Therefore,
\begin{equation}\label{thmhomeomorphismsandopensetsp2eq3}
\Foreach{\U}{\topology{\Y}}
\func{\pimage{\hf}}{\U}\in\topology{\X},
\end{equation}
which according to \refdef{defcontinuousfunction}, means,
\begin{equation}\label{thmhomeomorphismsandopensetsp2eq4}
\hf\in\CF{\Xt}{\Yt}.
\end{equation}
\begin{itemize}
\item[${\textbf{\textsf{p2-2}}}$]
$\U$
is taken as an arbitrary element of $\topology{\X}$. Then according to,
\Ref{thmhomeomorphismsandopensetsp2eq1}
and
\Ref{thmhomeomorphismsandopensetsp2eq3},
\begin{align}
\func{\pimage{\[\finv{\hf}\]}}{\U}&=\func{\image{\hf}}{\U}\cr
&\in\topology{\Y}.
\end{align}
\endp
\end{itemize}
Therefore,
\begin{equation}\label{thmhomeomorphismsandopensetsp2eq5}
\Foreach{\U}{\topology{\X}}
\func{\pimage{\[\finv{\hf}\]}}{\U}\in\topology{\Y},
\end{equation}
which according to \refdef{defcontinuousfunction}, means,
\begin{equation}\label{thmhomeomorphismsandopensetsp2eq6}
\finv{\hf}\in\CF{\Yt}{\Xt}.
\end{equation}
Based on \refdef{defhomeomorphism},
\Ref{thmhomeomorphismsandopensetsp2eq4}
and
\Ref{thmhomeomorphismsandopensetsp2eq6}
imply that,
\begin{equation}
\hf\in\HOF{\Xt}{\Yt}.
\end{equation}
\endp
\end{itemize}
\endthm
%%%%%%%%%%%%%%%%%%%%%%%%%%%%%%%%%%%%%%%%%%%%%%%%%%%%%%%%%%%%%%%%%%%%%%%%%%%%%%%%%%%%%%%
\theorem\label{thmhomeomorphismsisopenandcontinuousmap}
Each
$\Xt=\opair{\X}{\topology{\X}}$
and
$\Yt=\opair{\Y}{\topology{\Y}}$
is taken as a topological-space.
The set of all homeomorphisms from $\Xt$ to $\Yt$
equals the set of all bijective maps from $\X$ to $\Y$
that are simultaneously an open map from $\Xt$ to $\Yt$ and a continuous map from $\Xt$ to $\Yt$.
That is,
\begin{equation}
\HOF{\Xt}{\Yt}=\bigg[\IF{\X}{\Y}\cap\OM{\Xt}{\Yt}\cap\CF{\Xt}{\Yt}\bigg].
\end{equation}
\prooff
\begin{itemize}
\item[${\textbf{\textsf{p1}}}$]

According to \refdef{defopenmap},
\refdef{defhomeomorphism},
and
\refthm{thmhomeomorphismsandopensets},
it is evident that,
\begin{equation}
\HOF{\Xt}{\Yt}\subseteq\bigg[\IF{\X}{\Y}\cap\OM{\Xt}{\Yt}\cap\CF{\Xt}{\Yt}\bigg].
\end{equation}
\endp
\end{itemize}
\begin{itemize}
\item[${\textbf{\textsf{p2}}}$]
$\hf$
is taken as an arbitrary element of
$\bigg[\IF{\X}{\Y}\cap\OM{\Xt}{\Yt}\cap\CF{\Xt}{\Yt}\bigg]$. Then,
\begin{equation}
\Foreach{\Asubset{\X}}{\CSs{\X}}
\func{\pimage{\[\finv{\hf}\]}}{\Asubset{\X}}=\func{\image{\hf}}{\Asubset{\X}},
\end{equation}
and according to \refdef{defopenmap},
\begin{equation}
\Foreach{\U}{\topology{\X}}
\func{\image{\hf}}{\U}\in\topology{\Y}.
\end{equation}
Therefore,
\begin{equation}
\Foreach{\U}{\topology{\X}}
\func{\pimage{\[\finv{\hf}\]}}{\U}\in\topology{\Y},
\end{equation}
which according to \refdef{defcontinuousfunction}, means,
\begin{equation}
\finv{\hf}\in\CF{\Yt}{\Xt}.
\end{equation}
Therefore, according to \refdef{defhomeomorphism} it is clear that,
\begin{equation}
\hf\in\HOF{\Xt}{\Yt}.
\end{equation}
\endp
\end{itemize}
\endthm
%%%%%%%%%%%%%%%%%%%%%%%%%%%%%%%%%%%%%%%%%%%%%%%%%%%%%%%%%%%%%%%%%%%%%%%%%%%%%%%%%%%%%%%
\theorem\label{thmhomeomorphismsandclosedsets}
Each
$\Xt=\opair{\X}{\topology{\X}}$
and
$\Yt=\opair{\Y}{\topology{\Y}}$
is taken as a topological-space.
\begin{align}
\HOF{\Xt}{\Yt}=
\defset{\hf}{\IF{\X}{\Y}}{\bigg[
\defset{\U}{\CSs{\X}}{\func{\image{\hf}}{\U}\in\Fclosed{\Y}{\topology{\Y}}}
=\Fclosed{\X}{\topology{\X}}\bigg]}.
\end{align}
\proof
\begin{itemize}
\item[${\textbf{\textsf{p1}}}$]
$\hf$
is taken as an arbitrary element of $\HOF{\Xt}{\Yt}$. Then according to \refdef{defhomeomorphism},
\begin{align}
\hf&\in\IF{\X}{\Y},\label{thmhomeomorphismsandclosedsetsp1eq1}\\
\hf&\in\CF{\Xt}{\Yt},\label{thmhomeomorphismsandclosedsetsp1eq2}\\
\finv{\hf}&\in\CF{\Yt}{\Xt},\label{thmhomeomorphismsandclosedsetsp1eq3}
\end{align}
In addition, according to \refdef{deffamilyofclosedsets},
\begin{equation}\label{thmhomeomorphismsandclosedsetsp1eq4}
\Foreach{\U}{\CSs{\X}}
\bigg[\func{\image{\hf}}{\U}\in\Fclosed{\Y}{\topology{\Y}}
\thenn
\(\compl{\Y}{\func{\image{\hf}}{\U}}\)\in\topology{\Y}\bigg].
\end{equation}
According to \Ref{thmhomeomorphismsandclosedsetsp1eq1},
\begin{equation}\label{thmhomeomorphismsandclosedsetsp1eq5}
\Foreach{\U}{\CSs{\X}}
\bigg[\(\compl{\Y}{\func{\image{\hf}}{\U}}\)=\func{\image{\hf}}{\compl{\X}{\U}}\bigg].
\end{equation}
\Ref{thmhomeomorphismsandclosedsetsp1eq4}
and
\Ref{thmhomeomorphismsandclosedsetsp1eq5}
imply,
\begin{equation}\label{thmhomeomorphismsandclosedsetsp1eq6}
\Foreach{\U}{\CSs{\X}}
\bigg[\func{\image{\hf}}{\U}\in\Fclosed{\Y}{\topology{\Y}}
\thenn
\func{\image{\hf}}{\compl{\X}{\U}}\in\topology{\Y}\bigg].
\end{equation}
According to
\refthm{thmhomeomorphismsandopensets},
\begin{equation}\label{thmhomeomorphismsandclosedsetsp1eq7}
\Foreach{\U}{\CSs{\X}}
\bigg[\func{\image{\hf}}{\compl{\X}{\U}}\in\topology{\Y}
\thenn
\(\compl{\X}{\U}\)\in\topology{\X}\bigg].
\end{equation}
According to
\refdef{deffamilyofclosedsets},
\begin{equation}\label{thmhomeomorphismsandclosedsetsp1eq8}
\Foreach{\U}{\CSs{\X}}
\bigg[\(\compl{\X}{\U}\)\in\topology{\X}\thenn
\U\in\Fclosed{\X}{\topology{\X}}\bigg].
\end{equation}
\Ref{thmhomeomorphismsandclosedsetsp1eq4},
\Ref{thmhomeomorphismsandclosedsetsp1eq6},
\Ref{thmhomeomorphismsandclosedsetsp1eq7},
and
\Ref{thmhomeomorphismsandclosedsetsp1eq8}
imply that,
\begin{equation}
\Foreach{\U}{\CSs{\X}}
\bigg(\U\in\Fclosed{\X}{\topology{\X}}
\thenn
\func{\image{\hf}}{\U}\in\Fclosed{\Y}{\topology{\Y}}\bigg)
\end{equation}
\endp
\end{itemize}
\begin{itemize}
\item[${\textbf{\textsf{p2}}}$]
$\hf$
is taken as such an arbitrary element of $\IF{\X}{\Y}$ that,
\begin{equation}\label{thmhomeomorphismsandclosedsetsp2eq1}
\defset{\U}{\CSs{\X}}{\func{\image{\hf}}{\U}\in\Fclosed{\Y}{\topology{\Y}}}
=\Fclosed{\X}{\topology{\X}}.
\end{equation}
\begin{itemize}
\item[${\textbf{\textsf{p2-1}}}$]
$\U$
is taken as an arbitrary element of
$\CSs{\X}$. According to \refdef{deffamilyofclosedsets},
\begin{equation}\label{thmhomeomorphismsandclosedsetsp2-1eq1}
\U\in\topology{\X}\thenn
\(\compl{\X}{\U}\)\in\Fclosed{\X}{\topology{\X}}.
\end{equation}
According to \Ref{thmhomeomorphismsandclosedsetsp2eq1},
\begin{equation}\label{thmhomeomorphismsandclosedsetsp2-1eq2}
\(\compl{\X}{\U}\)\in\Fclosed{\X}{\topology{\X}}\thenn
\func{\image{\hf}}{\compl{\X}{\U}}\in\Fclosed{\Y}{\topology{\Y}}.
\end{equation}
Considering the bijectivity of $\hf$,
\begin{equation}\label{thmhomeomorphismsandclosedsetsp2-1eq3}
\Foreach{\Asubset{\X}}{\CSs{\X}}
\func{\image{\hf}}{\compl{\X}{\Asubset{\X}}}=
\[\compl{\Y}{\func{\image{\hf}}{\Asubset{\X}}}\].
\end{equation}
According to \refdef{deffamilyofclosedsets}
and
\Ref{thmhomeomorphismsandclosedsetsp2-1eq3}
imply that,
\begin{equation}\label{thmhomeomorphismsandclosedsetsp2-1eq4}
\func{\image{\hf}}{\compl{\X}{\U}}\in\Fclosed{\Y}{\topology{\Y}}\thenn
\func{\image{\hf}}{\U}\in\topology{\Y}.
\end{equation}
\Ref{thmhomeomorphismsandclosedsetsp2-1eq1},
\Ref{thmhomeomorphismsandclosedsetsp2-1eq2}
and
\Ref{thmhomeomorphismsandclosedsetsp2-1eq4}
imply that,
\begin{equation}
\U\in\topology{\X}\thenn
\func{\image{\hf}}{\U}\in\topology{\Y}.
\end{equation}
\endp
\end{itemize}
Therefore,
\begin{equation}
\defset{\U}{\CSs{\X}}{\func{\image{\hf}}{\U}\in\topology{\Y}}
=\topology{\X},
\end{equation}
and hence according to,
\refthm{thmhomeomorphismsandopensets},
\begin{equation}
\hf\in\HOF{\Xt}{\Yt}.
\end{equation}
\endp
\end{itemize}
\endthm
%%%%%%%%%%%%%%%%%%%%%%%%%%%%%%%%%%%%%%%%%%%%%%%%%%%%%%%%%%%%%%%%%%%%%%%%%%%%%%%%%%%%%%%
\theorem\label{thmhomeomorphismsisclosedandcontinuousmap}
Each
$\Xt=\opair{\X}{\topology{\X}}$
and
$\Yt=\opair{\Y}{\topology{\Y}}$
is taken as a topological-space.
The set of all homeomorphisms from
$\Xt$
to
$\Yt$
equals the set of all bijective maps from $\X$
to $\Y$
that are simultaneously a closed map from $\Xt$
to
$\Yt$
and a continuous map from $\Xt$ to $\Yt$. That is,
\begin{equation}
\HOF{\Xt}{\Yt}=\bigg[\IF{\X}{\Y}\cap\CM{\Xt}{\Yt}\cap\CF{\Xt}{\Yt}\bigg].
\end{equation}
\prooff
\begin{itemize}
\item[${\textbf{\textsf{p1}}}$]
According to
\refdef{defclosedmap},
\refdef{defhomeomorphism},
and
\refthm{thmhomeomorphismsandclosedsets},
it is evident that,
\begin{equation}
\HOF{\Xt}{\Yt}\subseteq\bigg[\IF{\X}{\Y}\cap\CM{\Xt}{\Yt}\cap\CF{\Xt}{\Yt}\bigg].
\end{equation}
\endp
\end{itemize}
\begin{itemize}
\item[${\textbf{\textsf{p2}}}$]
$\hf$
is taken as an arbitrary element of
$\bigg[\IF{\X}{\Y}\cap\CM{\Xt}{\Yt}\cap\CF{\Xt}{\Yt}\bigg]$. Then,
\begin{equation}
\Foreach{\Asubset{\X}}{\CSs{\X}}
\func{\pimage{\[\finv{\hf}\]}}{\Asubset{\X}}=\func{\image{\hf}}{\Asubset{\X}},
\end{equation}
and according to \refdef{defclosedmap},
\begin{equation}
\Foreach{\U}{\Fclosed{\X}{\topology{\X}}}
\func{\image{\hf}}{\U}\in\Fclosed{\Y}{\topology{\Y}}.
\end{equation}
Therefore,
\begin{equation}
\Foreach{\U}{\Fclosed{\X}{\topology{\X}}}
\func{\pimage{\[\finv{\hf}\]}}{\U}\in\Fclosed{\Y}{\topology{\Y}},
\end{equation}
which according to \refthm{thmcontiniuityequiv1}, means,
\begin{equation}
\finv{\hf}\in\CF{\Yt}{\Xt}.
\end{equation}
Therefore, according to \refdef{defhomeomorphism}, it is evident that,
\begin{equation}
\hf\in\HOF{\Xt}{\Yt}.
\end{equation}
\endp
\end{itemize}
\endthm
%%%%%%%%%%%%%%%%%%%%%%%%%%%%%%%%%%%%%%%%%%%%%%%%%%%%%%%%%%%%%%%%%%%%%%%%%%%%%%%%%%%%%%%
\corollary\label{corhomeomorphismisopenclosedcontinuousmap}
Each
$\Xt=\opair{\X}{\topology{\X}}$
and
$\Yt=\opair{\Y}{\topology{\Y}}$
is taken as a topological-space.
\begin{align}
\HOF{\Xt}{\Yt}&=\bigg[\IF{\X}{\Y}\cap\OM{\Xt}{\Yt}\cap\CF{\Xt}{\Yt}\bigg]\cr
&=\bigg[\IF{\X}{\Y}\cap\CM{\Xt}{\Yt}\cap\CF{\Xt}{\Yt}\bigg].
\end{align}
\endcor
%%%%%%%%%%%%%%%%%%%%%%%%%%%%%%%%%%%%%%%%%%%%%%%%%%%%%%%%%%%%%%%%%%%%%%%%%%%%%%%%%%%%%%%
\theorem\label{thmhomeomorphismstructurepreservation}
Each
$\Xt=\opair{\X}{\topology{\X}}$
and
$\Yt=\opair{\Y}{\topology{\Y}}$
is taken as a topological-space.
For every map
$\hf$
from
$\X$
to
$\Y$,
$\hf$
is a homeomorphism from
$\Xt$
to
$\Yt$ if and only if $\hf$
is a bijective map from $\X$ to $\Y$,
and the restriction of domain and codomain of
$\image{\hf}$
to
$\topology{\X}$
and
$\topology{\Y}$, respectively,
is a bijection from
$\topology{\X}$
to
$\topology{\Y}$. That is,
\begin{align}
&\Foreach{\hf}{\Func{\X}{\Y}}\\
&\[\hf\in\HOF{\Xt}{\Yt}\thenn\bigg(
\hf\in\IF{\X}{\Y},~\func{\rescd{\func{\resd{\image{\hf}}}
{\topology{\X}}}}{\topology{\Y}}\in\IF{\topology{\X}}{\topology{\Y}}\bigg)\].
\end{align}
This means any homeomorphism from
$\Xt$
to
$\Yt$
forms a one-to-one correspondence both between the set of points of
$\Xt$
and the that of,
$\Yt$
and between the set of open sets of
$\Xt$
and that of
$\Yt$.\\
An interpretation of this assertion is that a pair of homeomorphic topological-spaces have the same
topological structure.
\proof
$\hf$
is taken as an arbitrary map from
$\X$
to
$\Y$.
\begin{itemize}
\item[${\textbf{\textsf{p1}}}$]
It is assumed that,
\begin{equation}\label{thmhomeomorphismstructurepreservationp1eq0}
\hf\in\HOF{\Xt}{\Yt}.
\end{equation}
Then according to
\refdef{defhomeomorphism}
and
\refthm{thminverseofhomeomorphism},
\begin{align}
\hf&\in\IF{\X}{\Y},\label{thmhomeomorphismstructurepreservationp1eq1}\\
\finv{\hf}&\in\HOF{\Yt}{\Xt},\label{thmhomeomorphismstructurepreservationp1eq2}
\end{align}
According to
\Ref{thmhomeomorphismstructurepreservationp1eq1},
it is clear that,
\begin{align}
\image{\hf}&\in\IF{\CSs{\X}}{\CSs{\Y}},\label{thmhomeomorphismstructurepreservationp1eq3}\\
\pimage{\hf}&\in\IF{\CSs{\Y}}{\CSs{\X}}.\label{thmhomeomorphismstructurepreservationp1eq4}
\end{align}
\Ref{thmhomeomorphismstructurepreservationp1eq3}
implies that,
$\func{\resd{\image{\hf}}}{{\topology{\X}}}$
is an injective map from
$\topology{\X}$
to
$\CSs{\Y}$. That is,
\begin{align}\label{thmhomeomorphismstructurepreservationp1eq5}
&\Foreach{\opair{\U_1}{\U_2}}{\[\topology{\X}\times\topology{\X}\]}\cr
&\bigg(\func{\[\func{\resd{\image{\hf}}}{{\topology{\X}}}\]}{\U_1}=
\func{\[\func{\resd{\image{\hf}}}{{\topology{\X}}}\]}{\U_2}\bigg)\then
\(\U_1=\U_2\).
\end{align}
Additionally, according to
\refthm{thmhomeomorphismsandopensets},
\begin{equation}\label{thmhomeomorphismstructurepreservationp1eq6}
\func{\image{\[\func{\resd{\image{\hf}}}{{\topology{\X}}}\]}}{\topology{\X}}
\subseteq\topology{\Y}.
\end{equation}
Hence the function
$\func{\rescd{\func{\resd{\image{\hf}}}
{\topology{\X}}}}{\topology{\Y}}$
is well-defined, and an injective map from
$\topology{\X}$
to
$\topology{\Y}$.
\Ref{thmhomeomorphismstructurepreservationp1eq0},
\Ref{thmhomeomorphismstructurepreservationp1eq2},
and
\refthm{thmhomeomorphismsandopensets}
imply that,
\begin{align}
&\Foreach{\U}{\topology{\X}}
\func{\image{\hf}}{\U}\in\topology{\Y},\label{thmhomeomorphismstructurepreservationp1eq7}\\
&\Foreach{\U}{\topology{\Y}}
\func{\image{\[\finv{\hf}\]}}{\U}\in\topology{\X}.\label{thmhomeomorphismstructurepreservationp1eq8}
\end{align}
In addition, it is known that,
\begin{equation}\label{thmhomeomorphismstructurepreservationp1eq9}
\Foreach{\Asubset{\Y}}{\CSs{\Y}}
\func{\image{\hf}}{\func{\image{\[\finv{\hf}\]}}{\Asubset{\Y}}}=\Asubset{\Y}.
\end{equation}
\Ref{thmhomeomorphismstructurepreservationp1eq7},
\Ref{thmhomeomorphismstructurepreservationp1eq8},
and
\Ref{thmhomeomorphismstructurepreservationp1eq9}
imply that,
\begin{equation}
\Foreach{\U}{\topology{\Y}}
\Existsis{\func{\image{\[\finv{\hf}\]}}{\U}}{\topology{\X}}
\func{\image{\hf}}{\func{\image{\[\finv{\hf}\]}}{\U}}=\U.
\end{equation}
This means,
$\func{\rescd{\func{\resd{\image{\hf}}}
{\topology{\X}}}}{\topology{\Y}}$
is a surjective map from
$\topology{\X}$
to
$\topology{\Y}$.\\
Therefore,
\begin{equation}
\func{\rescd{\func{\resd{\image{\hf}}}
{\topology{\X}}}}{\topology{\Y}}\in\IF{\topology{\X}}{\topology{\Y}}.
\end{equation}
Therefore,
\begin{equation}
\IF{\topology{\X}}{\topology{\Y}}\neq\empty,
\end{equation}
which means,
\begin{equation}
\Card{\topology{\X}}\cardeq\Card{\topology{\Y}}.
\end{equation}
Briefly, it becomes clear that
$\X$
and
$\Y$
have the same cardinality. Moreover,
$\topology{\X}$
and
$\topology{\Y}$
also have the same cardinality.
\endp
\end{itemize}
\begin{itemize}
\item[${\textbf{\textsf{p2}}}$]
It is assumed that,
\begin{align}
&\hf\in\IF{\X}{\Y},\\
&\func{\rescd{\func{\resd{\image{\hf}}}
{\topology{\X}}}}{\topology{\Y}}\in\IF{\topology{\X}}{\topology{\Y}}.
\end{align}
Therefore, clearly,
\begin{align}
&\Foreach{\U}{\topology{\Y}}
\func{\pimage{\hf}}{\U}\in\topology{\X},\\
&\Foreach{\U}{\topology{\X}}
\func{\pimage{\[\finv{\hf}\]}}{\U}\in\topology{\Y},
\end{align}
which according to
\refdef{defcontinuousfunction},
means,
\begin{align}
\hf&\in\CF{\Xt}{\Yt},\\
\finv{\hf}&\in\CF{\Yt}{\Xt}.
\end{align}
Thus according to \refdef{defhomeomorphism},
it is clear that,
\begin{equation}
\hf\in\HOF{\Xt}{\Yt}.
\end{equation}
\endp
\end{itemize}
\endthm
%%%%%%%%%%%%%%%%%%%%%%%%%%%%%%%%%%%%%%%%%%%%%%%%%%%%%%%%%%%%%%%%%%%%%%%%%%%%%%%%%%%%%%%
\theorem\label{thmhomeomorphismandclosureofsets}
Each
$\Xt=\opair{\X}{\topology{\X}}$
and
$\Yt=\opair{\Y}{\topology{\Y}}$
is taken as a topological-space.
\begin{align}
&\HOF{\Xt}{\Yt}=\cr
&\defset{\hf}{\IF{\X}{\Y}}
{\bigg[\Foreach{\Asubset{\X}}{\CSs{\X}}
\func{\image{\hf}}{\func{\Cl{\Xt}}{\Asubset{\X}}}=
\func{\Cl{\Yt}}{\func{\image{\hf}}{\Asubset{\X}}}\bigg]}.\cr
&{}
\end{align}
\proof
$\hf$
is taken as a map from $\X$ to $\Y$.
\begin{itemize}
\item[${\textbf{\textsf{p1}}}$]
It is assumed that,
\begin{equation}\label{thmhomeomorphismandclosureofsetsp1eq0}
\hf\in\HOF{\Xt}{\Yt}.
\end{equation}
Then according to \refcor{corhomeomorphismisopenclosedcontinuousmap},
\begin{align}
\hf&\in\IF{\X}{\Y},\label{thmhomeomorphismandclosureofsetsp1eq1}\\
%\hf&\in\CF{\Xt}{\Yt},\label{thmhomeomorphismandclosureofsetsp1eq2}\\
\hf&\in\CM{\Xt}{\Yt}.\label{thmhomeomorphismandclosureofsetsp1eq2}
%\finv{\hf}&\in\CF{\Yt}{\Xt}.
\end{align}
According to \refthm{thmclosedmapandclosure},
\Ref{thmhomeomorphismandclosureofsetsp1eq2}
implies,
\begin{equation}\label{thmhomeomorphismandclosureofsetsp1eq3}
\Foreach{\Asubset{\X}}{\CSs{\X}}
\func{\image{\hf}}{\func{\Cl{\Xt}}{\Asubset{\X}}}\supseteq
\func{\Cl{\Yt}}{\func{\image{\hf}}{\Asubset{\X}}}.
\end{equation}
In addition, according to \refthm{thminverseofhomeomorphism}, it is clear that,
\begin{equation}\label{thmhomeomorphismandclosureofsetsp1eq4}
\finv{\hf}\in\HOF{\Yt}{\Xt},
\end{equation}
and hence according to \refcor{corhomeomorphismisopenclosedcontinuousmap},
\begin{equation}\label{thmhomeomorphismandclosureofsetsp1eq5}
\finv{\hf}\in\CM{\Yt}{\Xt},
\end{equation}
and thus according to \refthm{thmclosedmapandclosure},
\begin{equation}\label{thmhomeomorphismandclosureofsetsp1eq6}
\Foreach{\Asubset{\Y}}{\CSs{\Y}}
\func{\image{\[\finv{\hf}\]}}{\func{\Cl{\Yt}}{\Asubset{\Y}}}\supseteq
\func{\Cl{\Xt}}{\func{\image{\[\finv{\hf}\]}}{\Asubset{\Y}}},
\end{equation}
Thus,
Considering the bijectivity of $\hf$,
\begin{align}\label{thmhomeomorphismandclosureofsetsp1eq7}
\Foreach{\Asubset{\X}}{\CSs{\X}}
\func{\image{\[\finv{\hf}\]}}{\func{\Cl{\Yt}}{\func{\image{\hf}}{\Asubset{\X}}}}&\supseteq
\func{\Cl{\Xt}}{\func{\image{\[\finv{\hf}\]}}{\func{\image{\hf}}{\Asubset{\X}}}}\cr
&=\func{\Cl{\Xt}}{\Asubset{\X}},
\end{align}
and hence considering again the bijectivity of $\hf$,
\begin{align}\label{thmhomeomorphismandclosureofsetsp1eq8}
\Foreach{\Asubset{\X}}{\CSs{\X}}
\func{\Cl{\Yt}}{\func{\image{\hf}}{\Asubset{\X}}}&=
\func{\image{\hf}}
{\func{\image{\[\finv{\hf}\]}}{\func{\Cl{\Yt}}{\func{\image{\hf}}{\Asubset{\X}}}}}\cr
&\supseteq\func{\image{\hf}}{\func{\Cl{\Xt}}{\Asubset{\X}}}.
\end{align}
\Ref{thmhomeomorphismandclosureofsetsp1eq3}
and
\Ref{thmhomeomorphismandclosureofsetsp1eq8}
imply,
\begin{equation*}
\Foreach{\Asubset{\X}}{\CSs{\X}}
\func{\image{\hf}}{\func{\Cl{\Xt}}{\Asubset{\X}}}=
\func{\Cl{\Yt}}{\func{\image{\hf}}{\Asubset{\X}}}.
\end{equation*}
\endp
\end{itemize}
\begin{itemize}
\item[${\textbf{\textsf{p2}}}$]
$\hf$
is taken as such an arbitrary element of $\IF{\X}{\Y}$ such that,
\begin{equation}\label{thmhomeomorphismandclosureofsetsp2eq1}
\Foreach{\Asubset{\X}}{\CSs{\X}}
\func{\image{\hf}}{\func{\Cl{\Xt}}{\Asubset{\X}}}=
\func{\Cl{\Yt}}{\func{\image{\hf}}{\Asubset{\X}}}.
\end{equation}
Then,
\begin{align}
&\Foreach{\Asubset{\X}}{\CSs{\X}}
\func{\image{\hf}}{\func{\Cl{\Xt}}{\Asubset{\X}}}\supseteq
\func{\Cl{\Yt}}{\func{\image{\hf}}{\Asubset{\X}}},\label{thmhomeomorphismandclosureofsetsp2eq2}\\
&\Foreach{\Asubset{\X}}{\CSs{\X}}
\func{\image{\hf}}{\func{\Cl{\Xt}}{\Asubset{\X}}}\subseteq
\func{\Cl{\Yt}}{\func{\image{\hf}}{\Asubset{\X}}}.\label{thmhomeomorphismandclosureofsetsp2eq3}
\end{align}
According to
\refthm{thmclosedmapandclosure},
\Ref{thmhomeomorphismandclosureofsetsp2eq2} means,
\begin{equation}\label{thmhomeomorphismandclosureofsetsp2eq4}
\hf\in\CM{\Xt}{\Yt}.
\end{equation}
The bijectivity of $\hf$ and
\Ref{thmhomeomorphismandclosureofsetsp2eq3}
imply that,
\begin{align}\label{thmhomeomorphismandclosureofsetsp2eq5}
\Foreach{\Asubset{\Y}}{\CSs{\Y}}
\func{\image{\hf}}{\func{\Cl{\Xt}}{\func{\pimage{\hf}}{\Asubset{\Y}}}}&\subseteq
\func{\Cl{\Yt}}{\func{\image{\hf}}{\func{\pimage{\hf}}{\Asubset{\Y}}}}\cr
&=\func{\Cl{\Yt}}{\Asubset{\Y}},
\end{align}
and hence considering the bijectivity of $\hf$,
\begin{align}\label{thmhomeomorphismandclosureofsetsp2eq6}
\Foreach{\Asubset{\Y}}{\CSs{\Y}}
\func{\Cl{\Xt}}{\func{\pimage{\hf}}{\Asubset{\Y}}}&=
\func{\pimage{\hf}}
{\func{\image{\hf}}{\func{\Cl{\Xt}}{\func{\pimage{\hf}}{\Asubset{\Y}}}}}\cr
&\subseteq\func{\pimage{\hf}}{\func{\Cl{\Yt}}{\Asubset{\Y}}}.
\end{align}
Thus, according to
\refthm{thmcontiniuityandclosure},
\begin{equation}\label{thmhomeomorphismandclosureofsetsp2eq7}
\hf\in\CF{\Xt}{\Yt}.
\end{equation}
Therefore, based on \refcor{corhomeomorphismisopenclosedcontinuousmap},
and considering that
$\hf\in\IF{\X}{\Y}$,
\Ref{thmhomeomorphismandclosureofsetsp2eq4}
and
\Ref{thmhomeomorphismandclosureofsetsp2eq7},
imply that,
\begin{equation}
\hf\in\HOF{\Xt}{\Yt}.
\end{equation}
\endp
\end{itemize}
\endthm
%%%%%%%%%%%%%%%%%%%%%%%%%%%%%%%%%%%%%%%%%%%%%%%%%%%%%%%%%%%%%%%%%%%%%%%%%%%%%%%%%%%%%%%
\theorem\label{thmhomeomorphismandinteriorofsets}
Each
$\Xt=\opair{\X}{\topology{\X}}$
and
$\Yt=\opair{\Y}{\topology{\Y}}$
is taken as a topological-space.
\begin{align}
&\HOF{\Xt}{\Yt}=\cr
&\defset{\hf}{\IF{\X}{\Y}}
{\bigg[\Foreach{\Asubset{\X}}{\CSs{\X}}
\func{\image{\hf}}{\func{\Int{\Xt}}{\Asubset{\X}}}=
\func{\Int{\Yt}}{\func{\image{\hf}}{\Asubset{\X}}}\bigg]}.\cr
&{}
\end{align}
\proof
$\hf$
is taken as a map from
$\X$
to
$\Y$.
\begin{itemize}
\item[${\textbf{\textsf{p1}}}$]
It is assumed that,
\begin{equation}\label{thmhomeomorphismandinteriorofsetsp1eq0}
\hf\in\HOF{\Xt}{\Yt}.
\end{equation}
Then according to
\refcor{corhomeomorphismisopenclosedcontinuousmap},
\begin{align}
\hf&\in\IF{\X}{\Y},\label{thmhomeomorphismandinteriorofsetsp1eq1}\\
\hf&\in\CF{\Xt}{\Yt},\label{thmhomeomorphismandinteriorofsetsp1eq2}\\
\hf&\in\OM{\Xt}{\Yt}.\label{thmhomeomorphismandinteriorofsetsp1eq3}
%\finv{\hf}&\in\CF{\Yt}{\Xt}.
\end{align}
According to
\refthm{thmopenmapandinterior},
\Ref{thmhomeomorphismandinteriorofsetsp1eq3}
implies that,
\begin{equation}\label{thmhomeomorphismandinteriorofsetsp1eq4}
\Foreach{\Asubset{\X}}{\CSs{\X}}
\func{\image{\hf}}{\func{\Int{\Xt}}{\Asubset{\X}}}\subseteq
\func{\Int{\Yt}}{\func{\image{\hf}}{\Asubset{\X}}}.
\end{equation}
In addition, according to \refthm{thmcontiniuityandinterior} it is clear that,
\begin{equation}\label{thmhomeomorphismandinteriorofsetsp1eq5}
\Foreach{\Asubset{\Y}}{\CSs{\Y}}
\func{\pimage{\hf}}{\func{\Int{\Yt}}{\Asubset{\Y}}}\subseteq
\func{\Int{\Xt}}{\func{\pimage{\hf}}{\Asubset{\Y}}},
\end{equation}
and hence considering the bijectivity of $\hf$,
\begin{align}\label{thmhomeomorphismandinteriorofsetsp1eq6}
\Foreach{\Asubset{\X}}{\CSs{\X}}
\func{\pimage{\hf}}{\func{\Int{\Yt}}{\func{\image{\hf}}{\Asubset{\X}}}}&\subseteq
\func{\Int{\Xt}}{\func{\pimage{\hf}}{\func{\image{\hf}}{\Asubset{\X}}}}\cr
&=\func{\Int{\Xt}}{\Asubset{\X}},
\end{align}
and thus by considering again the bijectivity of $\hf$,
\begin{align}\label{thmhomeomorphismandinteriorofsetsp1eq7}
\Foreach{\Asubset{\X}}{\CSs{\X}}
\func{\Int{\Yt}}{\func{\image{\hf}}{\Asubset{\X}}}&=
\func{\image{\hf}}
{\func{\pimage{\hf}}{\func{\Int{\Yt}}{\func{\image{\hf}}{\Asubset{\X}}}}}\cr
&\subseteq
\func{\image{\hf}}{\func{\Int{\Xt}}{\Asubset{\X}}}.
\end{align}
\Ref{thmhomeomorphismandinteriorofsetsp1eq4}
and
\Ref{thmhomeomorphismandinteriorofsetsp1eq7}
imply,
\begin{equation*}
\Foreach{\Asubset{\X}}{\CSs{\X}}
\func{\image{\hf}}{\func{\Int{\Xt}}{\Asubset{\X}}}=
\func{\Int{\Yt}}{\func{\image{\hf}}{\Asubset{\X}}}.
\end{equation*}
\endp
\end{itemize}
\begin{itemize}
\item[${\textbf{\textsf{p2}}}$]
$\hf$
is taken as such an arbitrary element of
$\IF{\X}{\Y}$ such that,
\begin{equation}\label{thmhomeomorphismandinteriorofsetsp2eq1}
\Foreach{\Asubset{\X}}{\CSs{\X}}
\func{\image{\hf}}{\func{\Int{\Xt}}{\Asubset{\X}}}=
\func{\Int{\Yt}}{\func{\image{\hf}}{\Asubset{\X}}}.
\end{equation}
Then,
\begin{align}
&\Foreach{\Asubset{\X}}{\CSs{\X}}
\func{\image{\hf}}{\func{\Int{\Xt}}{\Asubset{\X}}}\subseteq
\func{\Int{\Yt}}{\func{\image{\hf}}{\Asubset{\X}}},\label{thmhomeomorphismandinteriorofsetsp2eq2}\\
&\Foreach{\Asubset{\X}}{\CSs{\X}}
\func{\image{\hf}}{\func{\Int{\Xt}}{\Asubset{\X}}}\supseteq
\func{\Int{\Yt}}{\func{\image{\hf}}{\Asubset{\X}}}.\label{thmhomeomorphismandinteriorofsetsp2eq3}
\end{align}
Based on
\refthm{thmopenmapandinterior},
\Ref{thmhomeomorphismandinteriorofsetsp2eq2}
means,
\begin{equation}\label{thmhomeomorphismandinteriorofsetsp2eq4}
\hf\in\OM{\Xt}{\Yt}.
\end{equation}
Considering the bijectivity of
$\hf$,
\Ref{thmhomeomorphismandinteriorofsetsp2eq3}
implies,
\begin{align}\label{thmhomeomorphismandinteriorofsetsp2eq5}
\Foreach{\Asubset{\Y}}{\CSs{\Y}}
\func{\image{\hf}}{\func{\Int{\Xt}}{\func{\pimage{\hf}}{\Asubset{\Y}}}}&\supseteq
\func{\Int{\Yt}}{\func{\image{\hf}}{\func{\pimage{\hf}}{\Asubset{\Y}}}}\cr
&=\func{\Int{\Yt}}{\Asubset{\Y}},
\end{align}
and hence by considering the bijectivity of $\hf$,
\begin{align}\label{thmhomeomorphismandinteriorofsetsp2eq6}
\Foreach{\Asubset{\Y}}{\CSs{\Y}}
\func{\Int{\Xt}}{\func{\pimage{\hf}}{\Asubset{\Y}}}&=
\func{\pimage{\hf}}
{\func{\image{\hf}}{\func{\Int{\Xt}}{\func{\pimage{\hf}}{\Asubset{\Y}}}}}\cr
&\supseteq\func{\pimage{\hf}}{\func{\Int{\Yt}}{\Asubset{\Y}}}.
\end{align}
Thus, according to
\refthm{thmcontiniuityandinterior},
\begin{equation}\label{thmhomeomorphismandinteriorofsetsp2eq7}
\hf\in\CF{\Xt}{\Yt}.
\end{equation}
Therefore, based on
\refcor{corhomeomorphismisopenclosedcontinuousmap},
and considering that
$\hf\in\IF{\X}{\Y}$,
\Ref{thmhomeomorphismandinteriorofsetsp2eq4}
and
\Ref{thmhomeomorphismandinteriorofsetsp2eq7},
imply,
\begin{equation}
\hf\in\HOF{\Xt}{\Yt}.
\end{equation}
\endp
\end{itemize}
\endthm
%%%%%%%%%%%%%%%%%%%%%%%%%%%%%%%%%%%%%%%%%%%%%%%%%%%%%%%%%%%%%%%%%%%%%%%%%%%%%%%%%%%%%%%
\theorem\label{thmhomeomorphismandfrontierofsets}
Each
$\Xt=\opair{\X}{\topology{\X}}$
and
$\Yt=\opair{\Y}{\topology{\Y}}$
is taken as a topological-space, and
$\hf$
a homeomorphism from
$\Xt$
to
$\Yt$.
\begin{equation}
\Foreach{\Asubset{\X}}{\CSs{\X}}
\func{\image{\hf}}{\func{\Fr{\Xt}}{\Asubset{\X}}}=
\func{\Fr{\Yt}}{\func{\image{\hf}}{\Asubset{\X}}}.
\end{equation}
\proof
Considering the bijectivity of $\hf$,
and according to \refdef{deffrontierofset},
\refthm{thmhomeomorphismandclosureofsets},
and
\refthm{thmhomeomorphismandinteriorofsets},
it is evident that,
\begin{align}
\Foreach{\Asubset{\X}}{\CSs{\X}}
\func{\image{\hf}}{\func{\Fr{\Xt}}{\Asubset{\X}}}&=
\func{\image{\hf}}{\compl{\func{\Cl{\Xt}}{\Asubset{\X}}}{\func{\Int{\Xt}}{\Asubset{\X}}}}\cr
&=\compl{\func{\image{\hf}}{\func{\Cl{\Xt}}{\Asubset{\X}}}}
{\func{\image{\hf}}{\func{\Int{\Xt}}{\Asubset{\X}}}}\cr
&=\compl{\func{\Cl{\Yt}}{\func{\image{\hf}}{\Asubset{\X}}}}
{\func{\Int{\Yt}}{\func{\image{\hf}}{\Asubset{\X}}}}\cr
&=\func{\Fr{\Xt}}{\func{\image{\hf}}{\Asubset{\X}}}.
\end{align}
\endthm
%%%%%%%%%%%%%%%%%%%%%%%%%%%%%%%%%%%%%%%%%%%%%%%%%%%%%%%%%%%%%%%%%%%%%%%%%%%%%%%%%%%%%%%
\theorem\label{thmhomeomorphismandfrontierofsets}
Each
$\Xt=\opair{\X}{\topology{\X}}$
and
$\Yt=\opair{\Y}{\topology{\Y}}$
is taken as a topological-space, and
$\hf$
a homeomorphism from
$\Xt$
to
$\Yt$.
\begin{equation}
\Foreach{\point}{\X}
\defset{\U}{\CSs{\X}}
{\func{\image{\hf}}{\U}\in\func{\nei{\Yt}}{\seta{\func{\hf}{\point}}}}=
\func{\nei{\Xt}}{\seta{\point}}.
\end{equation}
\prooff
According to
\refdef{defnbdclassofsets}
and
\refthm{thmhomeomorphismsisopenandcontinuousmap},
it is obvious.
\endthm
%%%%%%%%%%%%%%%%%%%%%%%%%%%%%%%%%%%%%%%%%%%%%%%%%%%%%%%%%%%%%%%%%%%%%%%%%%%%%%%%%%%%%%%
\theorem\label{thmrestrictionofhomeomorphism}
Each
$\Xt=\opair{\X}{\topology{\X}}$
and
$\Yt=\opair{\Y}{\topology{\Y}}$
is taken as a topological-space.
\begin{align}
&\Foreach{\hf}{\HOF{\Xt}{\Yt}}\cr
&\bigg[\Foreach{\Asubset{\X}}{\CSs{\X}}\cr
&~~~\func{\rescd{\func{\resd{\hf}}{\Asubset{\X}}}}{\func{\image{\hf}}{\Asubset{\X}}}
\in\HOF{\opair{\Asubset{\X}}{\stopology{\topology{\X}}{\Asubset{\X}}}}
{\opair{\func{\image{\hf}}{\Asubset{\X}}}{\stopology{\topology{\Y}}{\func{\image{\hf}}{\Asubset{\X}}}}}
\bigg].\cr
&{}
\end{align}
\prooff
$\hf$
is taken as an arbitrary element of
$\HOF{\Xt}{\Yt}$.
Then according to
\refcor{corhomeomorphismisopenclosedcontinuousmap},
\begin{align}
\hf&\in\IF{\X}{\Y},\label{thmrestrictionofhomeomorphismpeq1}\\
\hf&\in\CF{\Xt}{\Yt},\label{thmrestrictionofhomeomorphismpeq2}\\
\finv{\hf}&\in\CF{\Yt}{\Xt}.\label{thmrestrictionofhomeomorphismpeq3}
\end{align}
\Ref{thmrestrictionofhomeomorphismpeq1}
clearly implies,
\begin{equation}\label{thmrestrictionofhomeomorphismpeq4}
\Foreach{\Asubset{\X}}{\CSs{\X}}
\func{\rescd{\func{\resd{\hf}}{\Asubset{\X}}}}{\func{\image{\hf}}{\Asubset{\X}}}
\in\IF{\Asubset{\X}}{\func{\image{\hf}}{\Asubset{\X}}}.
\end{equation}
According to
\refcor{correstrictionofcontinuousfunction} and
\Ref{thmrestrictionofhomeomorphismpeq2},
\begin{align}\label{thmrestrictionofhomeomorphismpeq5}
&\Foreach{\Asubset{\X}}{\CSs{\X}}\cr
&\func{\rescd{\func{\resd{\hf}}{\Asubset{\X}}}}{\func{\image{\hf}}{\Asubset{\X}}}
\in\CF{\opair{\Asubset{\X}}{\stopology{\topology{\X}}{\Asubset{\X}}}}
{\opair{\func{\image{\hf}}{\Asubset{\X}}}{\stopology{\topology{\Y}}{\func{\image{\hf}}{\Asubset{\X}}}}}.\cr
&{}
\end{align}
According to
\refcor{correstrictionofcontinuousfunction}
and
\Ref{thmrestrictionofhomeomorphismpeq3},
\begin{align}\label{thmrestrictionofhomeomorphismpeq6}
&\Foreach{\Asubset{\Y}}{\CSs{\Y}}\cr
&\bigg[\func{\rescd{\func{\resd{\finv{\hf}}}{\Asubset{\Y}}}}{\func{\image{\[\finv{\hf}\]}}{\Asubset{\Y}}}\cr
&~~\in\CF{\opair{\Asubset{\Y}}{\stopology{\topology{\Y}}{\Asubset{\Y}}}}
{\opair{\func{\image{\[\finv{\hf}\]}}{\Asubset{\Y}}}{\stopology{\topology{\X}}{\func{\image{\[\finv{\hf}\]}}{\Asubset{\Y}}}}}\bigg],\cr
&{}
\end{align}
and hence,
\begin{align}\label{thmrestrictionofhomeomorphismpeq7}
&\Foreach{\Asubset{\X}}{\CSs{\X}}\cr
&\bigg[\func{\rescd{\func{\resd{\finv{\hf}}}{\func{\image{\hf}}{\Asubset{\X}}}}}{\func{\image{\[\finv{\hf}\]}}{\func{\image{\hf}}{\Asubset{\X}}}}\cr
&~~\in\CF{\opair{\func{\image{\hf}}{\Asubset{\X}}}{\stopology{\topology{\Y}}{\func{\image{\hf}}{\Asubset{\X}}}}}
{\opair{\func{\image{\[\finv{\hf}\]}}{\func{\image{\hf}}{\Asubset{\X}}}}{\stopology{\topology{\X}}{\func{\image{\[\finv{\hf}\]}}{\func{\image{\hf}}{\Asubset{\X}}}}}}\bigg].\cr
&{}
\end{align}
According to
\Ref{thmrestrictionofhomeomorphismpeq1},
\begin{align}\label{thmrestrictionofhomeomorphismpeq8}
&\Foreach{\Asubset{\X}}{\CSs{\X}}\cr
&\bigg[\func{\rescd{\func{\resd{\finv{\hf}}}{\func{\image{\hf}}{\Asubset{\X}}}}}{\func{\image{\[\finv{\hf}\]}}{\func{\image{\hf}}{\Asubset{\X}}}}=
\finv{\[\func{\rescd{\func{\resd{\hf}}{\Asubset{\X}}}}{\func{\image{\hf}}{\Asubset{\X}}}\]}\bigg],\cr
&{}
\end{align}
and
\begin{align}\label{thmrestrictionofhomeomorphismpeq9}
\Foreach{\Asubset{\X}}{\CSs{\X}}
\func{\image{\[\finv{\hf}\]}}{\func{\image{\hf}}{\Asubset{\X}}}=\Asubset{\X}.
\end{align}
\Ref{thmrestrictionofhomeomorphismpeq7},
\Ref{thmrestrictionofhomeomorphismpeq8},
and
\Ref{thmrestrictionofhomeomorphismpeq9}
imply,
\begin{align}\label{thmrestrictionofhomeomorphismpeq10}
&\Foreach{\Asubset{\X}}{\CSs{\X}}\cr
&\bigg[\finv{\[\func{\rescd{\func{\resd{\hf}}{\Asubset{\X}}}}{\func{\image{\hf}}{\Asubset{\X}}}\]}\in
\CF{\opair{\func{\image{\hf}}{\Asubset{\X}}}{\stopology{\topology{\Y}}{\func{\image{\hf}}{\Asubset{\X}}}}}
{\opair{\Asubset{\X}}{\stopology{\topology{\X}}{\Asubset{\X}}}}\bigg].\cr
&{}
\end{align}
Based on
\refdef{defhomeomorphism},
\Ref{thmrestrictionofhomeomorphismpeq4},
\Ref{thmrestrictionofhomeomorphismpeq5},
and
\Ref{thmrestrictionofhomeomorphismpeq10}
imply that,
\begin{align}
&\Foreach{\Asubset{\X}}{\CSs{\X}}\cr
&\bigg[\func{\rescd{\func{\resd{\hf}}{\Asubset{\X}}}}{\func{\image{\hf}}{\Asubset{\X}}}
\in\HOF{\opair{\Asubset{\X}}{\stopology{\topology{\X}}{\Asubset{\X}}}}
{\opair{\func{\image{\hf}}{\Asubset{\X}}}{\stopology{\topology{\Y}}{\func{\image{\hf}}{\Asubset{\X}}}}}
\bigg].\cr
&{}
\end{align}
\endthm
%%%%%%%%%%%%%%%%%%%%%%%%%%%%%%%%%%%%%%%%%%%%%%%%%%%%%%%%%%%%%%%%%%%%%%%%%%%%%%%%%%%%%%%%%%%%%%%%%%%%%%%%%%%%%%%%%%%%%%
%%%%%%%%%%%%%%%%%%%%%%%%%%%%%%%%%%%%%%%%%%%%%%%%%%%%%%%%%%%%%%%%%%%%%%%%%%%%%%%%%%%%%%%%%%%%%%%%%%%%%%%%%%%%%%%%%%%%%%
%%%%%%%%%%%%%%%%%%%%%%%%%%%%%%%%%%%%%%%%%%%%%%%%%%%%%%%%%%%%%%%%%%%%%%%%%%%%%%%%%%%%%%%%%%%%%%%%%%%%%%%%%%%%%%%%%%%%%%
%%%%%%%%%%%%%%%%%%%%%%%%%%%%%%%%%%%%%%%%%%%%%%%%%%%%%%%%%%%%%%%%%%%%%%%%%%%%%%%%%%%%%%%%%%%%%%%%%%%%%%%%%%%%%%%%%%%%%%
%%%%%%%%%%%%%%%%%%%%%%%%%%%%%%%%%%%%%%%%%%%%%%%%%%%%%%%%%%%%%%%%%%%%%%%%%%%%%%%%%%%%%%%%%%%%%%%%%%%%%%%%%%%%%%%%%%%%%%
%%%%%%%%%%%%%%%%%%%%%%%%%%%%%%%%%%%%%%%%%%%%%%%%%%%%%%%%%%%%%%%%%%%%%%%%%%%%%%%%%%%%%%%%%%%%%%%%%%%%%%%%%%%%%%%%%%%%%%
%%%%%%%%%%%%%%%%%%%%%%%%%%%%%%%%%%%%%%%%%%%%%%%%%%%%%%%%%%%%%%%%%%%%%%%%%%%%%%%%%%%%%%%%%%%%%%%%%%%%%%%%%%%%%%%%%%%%%%
%%%%%%%%%%%%%%%%%%%%%%%%%%%%%%%%%%%%%%%%%%%%%%%%%%%%%%%%%%%%%%%%%%%%%%%%%%%%%%%%%%%%%%%%%%%%%%%%%%%%%%%%%%%%%%%%%%%%%%
%%%%%%%%%%%%%%%%%%%%%%%%%%%%%%%%%%%%%%%%%%%%%%%%%%%%%%%%%%%%%%%%%%%%%%%%%%%%%%%%%%%%%%%%%%%%%%%%%%%%%%%%%%%%%%%%%%%%%%
%%%%%%%%%%%%%%%%%%%%%%%%%%%%%%%%%%%%%%%%%%%%%%%%%%%%%%%%%%%%%%%%%%%%%%%%%%%%%%%%%%%%%%%%%%%%%%%%%%%%%%%%%%%%%%%%%%%%%%
\section{
Embeddings
}
\definition\label{defembedding}
Each
$\Xt=\opair{\X}{\topology{\X}}$
and
$\Yt=\opair{\Y}{\topology{\Y}}$
is taken as a topological-space.
For every
$\em$
in
$\Func{\X}{\Y}$,
$\em$
is referred to as a $\quotl$embedding of the topological-space $\Xt$ in
the topological-space $\Yt$$\quotr$
iff
\begin{itemize}
\item[${\textbf{\textsf{EF1}}}$]
$\func{\rescd{\em}}{\func{\image{\em}}{\X}}\in
\HOF{\Xt}{\opair{\func{\image{\em}}{\X}}{\stopology{\topology{\Y}}{\func{\image{\em}}{\X}}}}.$
\end{itemize}
The set of all embeddings of
$\Xt$
in
$\Yt$
is denoted by
$\EM{\Xt}{\Yt}$. That is,
\begin{align}
\EM{\Xt}{\Yt}=\defset{\em}{\Func{\X}{\Y}}
{\func{\rescd{\em}}{\func{\image{\em}}{\X}}\in\HOF{\Xt}
{\opair{\func{\image{\em}}{\X}}{\stopology{\topology{\Y}}{\func{\image{\em}}{\X}}}}}.
\end{align}
\endef
%%%%%%%%%%%%%%%%%%%%%%%%%%%%%%%%%%%%%%%%%%%%%%%%%%%%%%%%%%%%%%%%%%%%%%%%%%%%%%%%%%%%%%%%%%%%%%%%%%%%%%%%%%%%%%%%%%%%%%
\definition\label{defembeddable}
Each
$\Xt=\opair{\X}{\topology{\X}}$
and
$\Yt=\opair{\Y}{\topology{\Y}}$
is taken as a topological-space.
it is said that $\quotl$the topological-space $\Xt$ can be embedded in
$\Yt$$\quotr$,
which is denoted by $\quotl$$\embedded{\Xt}{\Yt}$$\quotr$ iff
there exists at least one embedding of $\Xt$ in $\Yt$. That is,
\begin{equation}
\(\embedded{\Xt}{\Yt}\)\iffdef
\[\EM{\Xt}{\Yt}\neq\empty\].
\end{equation}
\endef
%%%%%%%%%%%%%%%%%%%%%%%%%%%%%%%%%%%%%%%%%%%%%%%%%%%%%%%%%%%%%%%%%%%%%%%%%%%%%%%%%%%%%%%%%%%%%%%%%%%%%%%%%%%%%%%%%%%%%%
\theorem\label{thminclusionmapisanembedding}
$\Xt=\opair{\X}{\topology{}}$
is taken as a topological-space.
For every
$\asubset$
in
$\CSs{\X}$,
the injection of
$\asubset$
into
$\X$
is an embedding of
$\opair{\asubset}{\stopology{\topology{}}{\asubset}}$
in
$\Xt$. That is,
\begin{equation}\label{thminclusionmapisanembeddingeq1}
\(\incf{\asubset}{\X}\)\in
\EM{\opair{\asubset}{\stopology{\topology{}}{\asubset}}}{\Xt}.
\end{equation}
\proof
$\asubset$
is taken as an arbitrary element of
$\CSs{\X}$. It is clear that,
\begin{align}
&\(\incf{\asubset}{\X}\)\in\Func{\X}{\Y},\\
&\func{\image{\(\incf{\asubset}{\X}\)}}{\asubset}=\asubset,
\end{align}
and,
\begin{equation}
\func{\rescd{\(\incf{\asubset}{\X}\)}}{\func{\image{\(\incf{\asubset}{\X}\)}}{\asubset}}
=\idf{\asubset}.
\end{equation}
Thus according to
\refthm{thmtrivialhomeomorphismofaspace},
\begin{align}
&\func{\rescd{\(\incf{\asubset}{\X}\)}}{\func{\image{\(\incf{\asubset}{\X}\)}}{\asubset}}\cr
\in
&\HOF{\opair{\asubset}{\stopology{\topology{}}{\asubset}}}{\opair{\asubset}{\stopology{\topology{}}{\asubset}}}\cr
=&\HOF{\opair{\asubset}{\stopology{\topology{}}{\asubset}}}
{\opair{\func{\image{\(\incf{\asubset}{\X}\)}}{\asubset}}
{\stopology{\topology{}}{\func{\image{\(\incf{\asubset}{\X}\)}}{\asubset}}}},
\end{align}
and hence according to
\refdef{defembedding},
\Ref{thminclusionmapisanembeddingeq1}
is obtained.
\endthm
%%%%%%%%%%%%%%%%%%%%%%%%%%%%%%%%%%%%%%%%%%%%%%%%%%%%%%%%%%%%%%%%%%%%%%%%%%%%%%%%%%%%%%%%%%%%%%%%%%%%%%%%%%%%%%%%%%%%%%
\theorem\label{thmcompositionofembeddings}
Each
$\Xt_1=\opair{\X_1}{\topology{1}}$,
$\Xt_2=\opair{\X_2}{\topology{2}}$,
and
$\Xt_3=\opair{\X_3}{\topology{3}}$
is taken as a topological-space.
The composition of an embedding of
$\Xt_2$
in
$\Xt_3$
with an embedding of
$\Xt_1$
in
$\Xt_2$,
is an embedding of
$\Xt_1$
in
$\Xt_3$. That is,
\begin{equation}
\Foreach{\opair{\em}{\p{\em}}}{\[\EM{\Xt_1}{\Xt_2}\times\EM{\Xt_2}{\Xt_3}\]}
\(\cmp{\p{\em}}{\em}\)\in\EM{\Xt_1}{\Xt_3}.
\end{equation}
\prooff
$\em$
is taken as an embedding of
$\Xt_1$
in
$\Xt_2$, and $\p{\em}$ as an embedding of $\Xt_2$ in $\Xt_3$.
Then according to \refdef{defembedding},
\begin{align}
&\func{\rescd{\em}}{\func{\image{\em}}{\X_1}}\in
\HOF{\Xt_1}{\opair{\func{\image{\em}}{\X_1}}{\stopology{\topology{2}}{\func{\image{\em}}{\X_1}}}},
\label{thmcompositionofembeddingspeq1}\\
&\func{\rescd{\p{\em}}}{\func{\image{\p{\em}}}{\X_2}}\in
\HOF{\Xt_2}{\opair{\func{\image{\p{\em}}}{\X_2}}{\stopology{\topology{3}}{\func{\image{\p{\em}}}{\X_2}}}}.
\label{thmcompositionofembeddingspeq2}
\end{align}
According to
\refthm{thmrestrictionofhomeomorphism},
and by defining,
\begin{equation}\label{thmcompositionofembeddingspeq3}
\xi=\func{\rescd{\p{\em}}}{\func{\image{\p{\em}}}{\X_2}},
\end{equation}
\Ref{thmcompositionofembeddingspeq2}
implies that,
\begin{align}\label{thmcompositionofembeddingspeq4}
&\func{\rescd{\func{\resd{\xi}}{\func{\image{\em}}{\X_1}}}}{\func{\image{\xi}}{\func{\image{\em}}{\X_1}}}\cr
\in&\HOF{\opair{\func{\image{\em}}{\X_1}}{\stopology{\topology{2}}{\func{\image{\em}}{\X_1}}}}
{\opair{\func{\image{\xi}}{\func{\image{\em}}{\X_1}}}
{\stopology{\stopology{\topology{3}}{\func{\image{\p{\em}}}{\X_2}}}{\func{\image{\xi}}{\func{\image{\em}}{\X_1}}}}}.\cr
&{}
\end{align}
Therefore, considering that,
\begin{align}
&\func{\image{\xi}}{\func{\image{\em}}{\X_1}}=
\func{\image{\(\cmp{\p{\em}}{\em}\)}}{\X_1},\label{thmcompositionofembeddingspeq5}\\
&\stopology{\stopology{\topology{3}}{\func{\image{\p{\em}}}{\X_2}}}{\func{\image{\xi}}{\func{\image{\em}}{\X_1}}}=
\stopology{\topology{3}}{\func{\image{\xi}}{\func{\image{\em}}{\X_1}}},\label{thmcompositionofembeddingspeq6}
\end{align}
it is evident that,
\begin{align}\label{thmcompositionofembeddingspeq7}
&\cmp{\[\func{\rescd{\func{\resd{\xi}}{\func{\image{\em}}{\X_1}}}}{\func{\image{\xi}}{\func{\image{\em}}{\X_1}}}\]}
{\[\func{\rescd{\em}}{\func{\image{\em}}{\X_1}}\]}\cr
\in&\HOF{\Xt_1}
{\opair{\func{\image{\(\cmp{\p{\em}}{\em}\)}}{\X_1}}
{\stopology{\topology{3}}{\func{\image{\(\cmp{\p{\em}}{\em}\)}}{\X_1}}}}.
\end{align}
In addition, it can be easily seen that,
\begin{align}\label{thmcompositionofembeddingspeq8}
&\cmp{\[\func{\rescd{\func{\resd{\xi}}{\func{\image{\em}}{\X_1}}}}{\func{\image{\xi}}{\func{\image{\em}}{\X_1}}}\]}
{\[\func{\rescd{\em}}{\func{\image{\em}}{\X_1}}\]}\cr
=
&\func{\rescd{\(\cmp{\p{\em}}{\em}\)}}{\func{\image{\(\cmp{\p{\em}}{\em}\)}}{\X_1}}.
\end{align}
\Ref{thmcompositionofembeddingspeq7}
and
\Ref{thmcompositionofembeddingspeq8}
imply that,
\begin{equation}
\func{\rescd{\(\cmp{\p{\em}}{\em}\)}}{\func{\image{\(\cmp{\p{\em}}{\em}\)}}{\X_1}}
\in
\HOF{\Xt_1}
{\opair{\func{\image{\(\cmp{\p{\em}}{\em}\)}}{\X_1}}
{\stopology{\topology{3}}{\func{\image{\(\cmp{\p{\em}}{\em}\)}}{\X_1}}}},
\end{equation}
and thus according to \refdef{defembedding},
\begin{equation}
\(\cmp{\p{\em}}{\em}\)\in\EM{\Xt_1}{\Xt_3}.
\end{equation}
\endthm
%%%%%%%%%%%%%%%%%%%%%%%%%%%%%%%%%%%%%%%%%%%%%%%%%%%%%%%%%%%%%%%%%%%%%%%%%%%%%%%%%%%%%%%%%%%%%%%%%%%%%%%%%%%%%%%%%%%%%%
\definition\label{defequivalentembeddings}
Each
$\Xt=\opair{\X}{\topology{\X}}$
and
$\Yt=\opair{\Y}{\topology{\Y}}$
is taken as a topological-space,
and each
$\em_1$
and
$\em_2$
as an embedding of
$\Xt$
in
$\Yt$. It is said that
$\quotl$$\em_1$ is equivalent to $\em_2$$\quotr$, which is denoted by $\quotl$$\embeq{\em_1}{\em_2}{\Xt}{\Yt}$$\quotr$, iff
\begin{equation*}
\Exists{\opair{\hf_{1}}{\hf_2}}{\HOF{\Xt}{\Xt}\times\HOF{\Yt}{\Yt}}
\(\cmp{\em_1}{\hf_1}=\cmp{\hf_2}{\em_2}\).
\end{equation*}
\endef
\chapteR{
Connectedness
}
\thispagestyle{fancy}
\section{
Partitions of a Topological Space
}
\subsection{
Partition of a Set
}
\definition\label{defpartitionofset}
$\X$
is taken as a set.
For every
$\apartition$
in
$\CSs{\CSs{\X}}$,
$\apartition$
is referred to as a $\quotl$partition of $\X$$\quotr$ iff these properties are satisfied.
\begin{itemize}
\item[${\textbf{\textsf{PS1}}}$]
$\opair{\X}{\apartition}$
is a cover of $\X$.
\hfill
$\apartition\in\covers{\X}.$
\item[${\textbf{\textsf{PS2}}}$]
The empty-set is not an element of $\apartition$.
\hfill
$\seta{\empty}\cap\apartition=\empty.$
\item[${\textbf{\textsf{PS3}}}$]
The element of
$\apartition$\\
are pairwise distinct.
\hfill
$\Foreach{\covelm_1}{\apartition}
\bigg(\Foreach{\covelm_2}{\compl{\apartition}{\seta{\covelm_1}}}
\covelm_1\cap\covelm_2=\empty\bigg).$
\end{itemize}
\begin{itemize}
\item
The set of all partitions of $\X$ is denoted by $\Cpart{\X}$.
\item
Each element of
$\compl{\Cpart{\X}}{\seta{\seta{\X}}}$
is referred to as a $\quotl$non-trivial partition of $\X$$\quotr$.
\item
Each finite element of $\Cpart{\X}$
is referred to as a $\quotl$finite partition of $\X$$\quotr$.
\end{itemize}
\endef
%%%%%%%%%%%%%%%%%%%%%%%%%%%%%%%%%%%%%%%%%%%%%%%%%%%%%%%%%%%%%%%%%%%%%%%%%%%%%%%%%%%%%%%%%%%%%%%%%%%%%%%%%%%
\theorem\label{thmpartitionofsetequiv0}
$\X$
is taken as a set.
\begin{equation}
\Cpart{\X}=\defset{\apartition}{\bigg(\compl{\CSs{\compl{\CSs{\X}}{\seta{\empty}}}}{\seta{\empty}}\bigg)}
{\X=\disunion{\apartition}}.
\end{equation}
\proof
According to
\refdef{defcovermapofset}
and
\refdef{defpartitionofset},
it is clear.
\endthm
%%%%%%%%%%%%%%%%%%%%%%%%%%%%%%%%%%%%%%%%%%%%%%%%%%%%%%%%%%%%%%%%%%%%%%%%%%%%%%%%%%%%%%%%%%%%%%%%%%%%%%%%%%%
\corollary\label{coremptysetpartition}
$\empty$
is the only partition of the empty-set. That is,
\begin{equation}
\Cpart{\empty}=\seta{\empty}.
\end{equation}
\endcor
%%%%%%%%%%%%%%%%%%%%%%%%%%%%%%%%%%%%%%%%%%%%%%%%%%%%%%%%%%%%%%%%%%%%%%%%%%%%%%%%%%%%%%%%%%%%%%%%%%%%%%%%%%%
%%%%%%%%%%%%%%%%%%%%%%%%%%%%%%%%%%%%%%%%%%%%%%%%%%%%%%%%%%%%%%%%%%%%%%%%%%%%%%%%%%%%%%%%%%%%%%%%%%%%%%%%%%%
%%%%%%%%%%%%%%%%%%%%%%%%%%%%%%%%%%%%%%%%%%%%%%%%%%%%%%%%%%%%%%%%%%%%%%%%%%%%%%%%%%%%%%%%%%%%%%%%%%%%%%%%%%%
\subsection{
Open and Closed Partitions of a Topological Space
}
\definition\label{defopenpartition}
$\Xt=\opair{\X}{\topology{}}$
is taken as a topological-space.
For every
$\apartition$
in
$\CSs{\CSs{\X}}$,
$\apartition$
is referred to as an $\quotl$open partition of the topological-space $\Xt$$\quotr$ iff these properties are satisfied.
\begin{itemize}
\item[${\textbf{\textsf{OP1}}}$]
$\apartition$
is a partition of $\X$.
\hfill
$\apartition\in\Cpart{\X}.$
\item[${\textbf{\textsf{OP2}}}$]
Each element of $\apartition$ is an open set of $\Xt$.
\hfill
$\apartition\subseteq\topology{}.$
\end{itemize}
\begin{itemize}
\item
The set of all open partitions of $\Xt$ is denoted by $\Opart{\Xt}$. That is,
\begin{equation}
\Opart{\Xt}:=\defset{\apartition}{\Cpart{\X}}
{\apartition\subseteq\topology{}}.
\end{equation}
\item
Each element of
$\compl{\Opart{\X}}{\seta{\seta{\X}}}$
is referred to as a $\quotl$non-trivial open partition of $\Xt$$\quotr$.
\end{itemize}
\endef
%%%%%%%%%%%%%%%%%%%%%%%%%%%%%%%%%%%%%%%%%%%%%%%%%%%%%%%%%%%%%%%%%%%%%%%%%%%%%%%%%%%%%%%%%%%%%%%%%%%%%%%%%%%%%%%%%%%%%%%%%%%%%%%%
\definition\label{defclosedpartition}
$\Xt=\opair{\X}{\topology{}}$
is taken as a topological-space.
For every $\apartition$ in
$\CSs{\CSs{\X}}$,
$\apartition$
is referred to as a $\quotl$closed partition of the topological-space $\Xt$$\quotr$
iff these properties are satisfied.
\begin{itemize}
\item[${\textbf{\textsf{CP1}}}$]
$\apartition$
is a partition of $\X$.
\hfill
$\apartition\in\Cpart{\X}.$
\item[${\textbf{\textsf{CP2}}}$]
Each element of $\apartition$
is a closd set of $\Xt$.
\hfill
$\apartition\subseteq\Fclosed{\X}{\topology{}}.$
\end{itemize}
\begin{itemize}
\item
The set of all closed partitions of $\Xt$ is denoted by $\Clpart{\Xt}$. That is,
\begin{equation}
\Clpart{\Xt}:=\defset{\apartition}{\Cpart{\X}}
{\apartition\subseteq\Fclosed{\X}{\topology{}}}.
\end{equation}
\item
Each element of
$\compl{\Clpart{\X}}{\seta{\seta{\X}}}$
is referred to as a $\quotl$non-trivial closed partition of $\Xt$$\quotr$.
\end{itemize}
\endef
%%%%%%%%%%%%%%%%%%%%%%%%%%%%%%%%%%%%%%%%%%%%%%%%%%%%%%%%%%%%%%%%%%%%%%%%%%%%%%%%%%%%%%%%%%%%%%%%%%%%%%%%%%%%%%%%%%%%%%%%%%%%%%%%
%%%%%%%%%%%%%%%%%%%%%%%%%%%%%%%%%%%%%%%%%%%%%%%%%%%%%%%%%%%%%%%%%%%%%%%%%%%%%%%%%%%%%%%%%%%%%%%%%%%%%%%%%%%%%%%%%%%%%%%%%%%%%%%%
%%%%%%%%%%%%%%%%%%%%%%%%%%%%%%%%%%%%%%%%%%%%%%%%%%%%%%%%%%%%%%%%%%%%%%%%%%%%%%%%%%%%%%%%%%%%%%%%%%%%%%%%%%%%%%%%%%%%%%%%%%%%%%%%
\subsection{
Non-Attached Partitions of a Topological Space
}
\definition\label{defnonclingingpartition}
$\Xt=\opair{\X}{\topology{}}$
is taken as a topological-space.
For every
$\apartition$
in
$\CSs{\CSs{\X}}$,
$\apartition$
is referred to as a $\quotl$non-attached partition of $\Xt$$\quotr$ iff these properties are satisfied.
\begin{itemize}
\item[${\textbf{\textsf{NcP1}}}$]
$\apartition$
is a partition of $\X$.
\hfill
$\apartition\in\Cpart{\X}.$
\item[${\textbf{\textsf{NcP2}}}$]
\begin{equation}
\Foreach{\covelm_1}{\apartition}
\bigg(\Foreach{\covelm_2}{\compl{\apartition}{\seta{\covelm_1}}}
\covelm_1\cap\func{\Cl{\Xt}}{\covelm_2}=\empty\bigg).
\end{equation}
\end{itemize}
\begin{itemize}
\item
The set of all non-attached partitions of $\Xt$ is denoted by $\Ncpart{\Xt}$.
\item
Each element of $\compl{\Ncpart{\X}}{\seta{\seta{\X}}}$
is referred to as a $\quotl$non-trivial non-attached partition of $\Xt$$\quotr$.
\end{itemize}
\endef
%%%%%%%%%%%%%%%%%%%%%%%%%%%%%%%%%%%%%%%%%%%%%%%%%%%%%%%%%%%%%%%%%%%%%%%%%%%%%%%%%%%%%%%%%%%%%%%%%%%%%%%%%%%%%%%%%%%%%%%%%%%%%%%%
%%%%%%%%%%%%%%%%%%%%%%%%%%%%%%%%%%%%%%%%%%%%%%%%%%%%%%%%%%%%%%%%%%%%%%%%%%%%%%%%%%%%%%%%%%%%%%%%%%%%%%%%%%%%%%%%%%%%%%%%%%%%%%%%
%%%%%%%%%%%%%%%%%%%%%%%%%%%%%%%%%%%%%%%%%%%%%%%%%%%%%%%%%%%%%%%%%%%%%%%%%%%%%%%%%%%%%%%%%%%%%%%%%%%%%%%%%%%%%%%%%%%%%%%%%%%%%%%%
%%%%%%%%%%%%%%%%%%%%%%%%%%%%%%%%%%%%%%%%%%%%%%%%%%%%%%%%%%%%%%%%%%%%%%%%%%%%%%%%%%%%%%%%%%%%%%%%%%%%%%%%%%%%%%%%%%%%%%%%%%%%%%%%
%%%%%%%%%%%%%%%%%%%%%%%%%%%%%%%%%%%%%%%%%%%%%%%%%%%%%%%%%%%%%%%%%%%%%%%%%%%%%%%%%%%%%%%%%%%%%%%%%%%%%%%%%%%%%%%%%%%%%%%%%%%%%%%%
%%%%%%%%%%%%%%%%%%%%%%%%%%%%%%%%%%%%%%%%%%%%%%%%%%%%%%%%%%%%%%%%%%%%%%%%%%%%%%%%%%%%%%%%%%%%%%%%%%%%%%%%%%%%%%%%%%%%%%%%%%%%%%%%
%%%%%%%%%%%%%%%%%%%%%%%%%%%%%%%%%%%%%%%%%%%%%%%%%%%%%%%%%%%%%%%%%%%%%%%%%%%%%%%%%%%%%%%%%%%%%%%%%%%%%%%%%%%%%%%%%%%%%%%%%%%%%%%%
%%%%%%%%%%%%%%%%%%%%%%%%%%%%%%%%%%%%%%%%%%%%%%%%%%%%%%%%%%%%%%%%%%%%%%%%%%%%%%%%%%%%%%%%%%%%%%%%%%%%%%%%%%%%%%%%%%%%%%%%%%%%%%%%
%%%%%%%%%%%%%%%%%%%%%%%%%%%%%%%%%%%%%%%%%%%%%%%%%%%%%%%%%%%%%%%%%%%%%%%%%%%%%%%%%%%%%%%%%%%%%%%%%%%%%%%%%%%%%%%%%%%%%%%%%%%%%%%%
\section{
Connected Topological Spaces
}
\definition\label{defconnectedness}
$\Xt=\opair{\X}{\topology{}}$
is taken as a topological-space.
$\Xt$
is referred to as a $\quotl$connected topological-space$\quotr$ iff
$\empty$
and
$\X$
are the only open-and-closed (clopen) sets of $\Xt$, that is,
\begin{equation*}
\topology{}\cap\Fclosed{\X}{\topology{}}=\seta{\binary{\empty}{\X}}.
\end{equation*}
\begin{itemize}
\item
It is said that $\quotl$$\X$ is connected under the topology $\topology{}$$\quotr$ iff
$\opair{\X}{\topology{}}$ is a connected topological-space.
\item
The set of all subsets of $\X$ that are connected under the topology induced from $\topology{}$ on them,
is denoted by $\connecteds{\Xt}$. That is,
\begin{align}
\connecteds{\Xt}:=
\defset{\asubset}{\CSs{\X}}
{\bigg(\stopology{\topology{}}{\asubset}\cap\Fclosed{\asubset}
{\stopology{\topology{}}{\asubset}}=\seta{\binary{\empty}{\asubset}}\bigg)}.
\end{align}
\item
For every
$\asubset$
in
$\CSs{\X}$,
$\asubset$ is called a $\quotl$connected set of the topological-space $\Xt$$\quotr$ iff,
$\opair{\asubset}{\stopology{\topology{}}{\asubset}}$
is a connected topological-space.
\item
$\Xt$
is referred to as a $\quotl$disconnected topological-space$\quotr$ iff
$\Xt$ is not a connected topological-space, that is,
$\Xt$ possesses at least one clopen set other than
$\empty$
and
$\X$, which means,
\begin{equation*}
\compl{\(\topology{}\cap\Fclosed{\X}{\topology{}}\)}{\seta{\binary{\empty}{\X}}}\neq\empty.
\end{equation*}
\end{itemize}
\endef
%%%%%%%%%%%%%%%%%%%%%%%%%%%%%%%%%%%%%%%%%%%%%%%%%%%%%%%%%%%%%%%%%%%%%%%%%%%%%%%%%%%%%%%%%%%%%%%%%%%%%%%%%%%%%%%%%%%%%%%%%%%%%%%%
\theorem\label{thmemptyspaceisconnected}
$\opair{\empty}{\seta{\empty}}$
is a connected topological-space. Moreover,
$\empty$
is the only connected set of $\opair{\empty}{\seta{\empty}}$. That is,
\begin{equation}
\connecteds{\opair{\empty}{\seta{\empty}}}=\seta{\empty}.
\end{equation}
\proof
According to \refdef{defconnectedness},
it is obvious.
\endthm
%%%%%%%%%%%%%%%%%%%%%%%%%%%%%%%%%%%%%%%%%%%%%%%%%%%%%%%%%%%%%%%%%%%%%%%%%%%%%%%%%%%%%%%%%%%%%%%%%%%%%%%%%%%%%%%%%%%%%%%%%%%%%%%%
\theorem\label{thmemptysetisaconnectedsetofeverytopologicalspace}
$\Xt=\opair{\X}{\topology{}}$
is taken as a topological-space.
$\empty$
is a connected set of
$\Xt$. That is,
\begin{equation}
\empty\in\connecteds{\Xt}.
\end{equation}
\proof
According to
\refdef{defsubspacetopology1},
\refdef{defconnectedness},
and
\refthm{thmemptyspaceisconnected},
it is obvious.
\endthm
%%%%%%%%%%%%%%%%%%%%%%%%%%%%%%%%%%%%%%%%%%%%%%%%%%%%%%%%%%%%%%%%%%%%%%%%%%%%%%%%%%%%%%%%%%%%%%%%%%%%%%%%%%%%%%%%%%%%%%%%%%%%%%%%
\theorem\label{thmindiscretespaceisconnected}
$\X$
is taken as a set.
The indiscrete topological-space
$\opair{\X}{\seta{\binary{\empty}{\X}}}$
is connected.
\begin{equation}
\X\in\connecteds{\opair{\X}{\seta{\binary{\empty}{\X}}}}.
\end{equation}
\proof
According to \refdef{defconnectedness},
it is clear.
\endthm
%%%%%%%%%%%%%%%%%%%%%%%%%%%%%%%%%%%%%%%%%%%%%%%%%%%%%%%%%%%%%%%%%%%%%%%%%%%%%%%%%%%%%%%%%%%%%%%%%%%%%%%%%%%%%%%%%%%%%%%%%%%%%%%%
\theorem\label{thmeverysubsetofindiscretespaceisaconnectedset}
$\X$
is taken as a set.
Every subset of $\X$
is a connected set of the indiscrete topological-space $\opair{\X}{\seta{\binary{\empty}{\X}}}$. That is,
\begin{equation}
\connecteds{\opair{\X}{\seta{\binary{\empty}{\X}}}}=\CSs{\X}.
\end{equation}
\proof
$\asubset$
is taken as an arbitrary subset of $\X$. According to \refdef{defsubspacetopology1},
\begin{align}
\stopology{\seta{\binary{\empty}{\X}}}{\asubset}&=
\seta{\binary{\empty\cap\asubset}{\X\cap\asubset}}\cr
&=\seta{\binary{\empty}{\asubset}}.
\end{align}
Thus, according to \refthm{thmindiscretespaceisconnected},
$\opair{\asubset}{\stopology{\seta{\binary{\empty}{\X}}}{\asubset}}$
is a connected topological space and hence according to \refdef{defconnectedness},
\begin{equation}
\asubset\in\connecteds{\opair{\X}{\seta{\binary{\empty}{\X}}}}.
\end{equation}
\endthm
%%%%%%%%%%%%%%%%%%%%%%%%%%%%%%%%%%%%%%%%%%%%%%%%%%%%%%%%%%%%%%%%%%%%%%%%%%%%%%%%%%%%%%%%%%%%%%%%%%%%%%%%%%%%%%%%%%%%%%%%%%%%%%%%
\theorem\label{thmsingletonspaceisconnected}
$\x$
is taken as a set.
The singleton topological-space
$\singletonTS{\x}=\opair{\seta{\x}}{\CSs{\seta{\x}}}$ is connected. That is,
\begin{equation}
\seta{\x}\in\connecteds{\singletonTS{\seta{\x}}}.
\end{equation}
\proof
According to
\refthm{thmsingletontopology0},
$\singletonTS{\x}$
is an indiscrete topological-space and hence according to \refthm{thmindiscretespaceisconnected},
it is connected.
\endthm
%%%%%%%%%%%%%%%%%%%%%%%%%%%%%%%%%%%%%%%%%%%%%%%%%%%%%%%%%%%%%%%%%%%%%%%%%%%%%%%%%%%%%%%%%%%%%%%%%%%%%%%%%%%%%%%%%%%%%%%%%%%%%%%%
\theorem\label{thmsingletonsubspaceisconnected}
$\Xt=\opair{\X}{\topology{}}$
is taken as a topological-space.
For every
$\x$
in
$\X$,
the singleton
$\seta{\x}$
is a connected set of
$\Xt$. That is,
\begin{equation}
\Foreach{\x}{\X}
\bigg[\seta{\x}\in\connecteds{\Xt}\bigg].
\end{equation}
\prooff
$\x$
is taken as an arbitrary element of
$\X$. According to
\refthm{thmsingletonsubspacetopology},
\begin{align}
\opair{\seta{\x}}{\stopology{\topology{}}{\seta{\x}}}&=
\opair{\seta{\x}}{\CSs{\seta{\x}}}\cr
&=\singletonTS{\x},
\end{align}
and hence according to
\refthm{thmsingletonspaceisconnected},
$\opair{\seta{\x}}{\stopology{\topology{}}{\seta{\x}}}$
is a connected topological-space. According to
\refdef{defconnectedness},
this means,
\begin{equation}
\seta{\x}\in\connecteds{\Xt}.
\end{equation}
\endthm
%%%%%%%%%%%%%%%%%%%%%%%%%%%%%%%%%%%%%%%%%%%%%%%%%%%%%%%%%%%%%%%%%%%%%%%%%%%%%%%%%%%%%%%%%%%%%%%%%%%%%%%%%%%%%%%%%%%%%%%%%%%%%%%%
\theorem\label{thmconnectednessandopenpartitions}
$\Xt=\opair{\X}{\topology{}}$
is taken as a topological-space.
$\Xt$
is a connected topological-space if and only if
$\Xt$
does not have any (non-trivial) open partition consisting of two elements. That is,
\begin{equation}
\bigg(\X\in\connecteds{\Xt}\bigg)\thenn
\bigg(\defset{\apartition}{\Opart{\Xt}}
{\CarD{\apartition}=2}=\empty\bigg).
\end{equation}
\prooff
\begin{itemize}
\item[${\textbf{\textsf{p1}}}$]
It is assumed that,
\begin{equation}\label{thmconnectednessandopenpartitionsp1eq1}
\defset{\apartition}{\Opart{\Xt}}
{\CarD{\apartition}=2}\neq\empty.
\end{equation}
Then there exists a partition of $\Xt$ such that,
\begin{align}
&\apartition=\seta{\binary{\U_1}{\U_2}},\label{thmconnectednessandopenpartitionsp1eq2}\\
&\U_1\notin\seta{\binary{\empty}{\X}},\label{thmconnectednessandopenpartitionsp1eq3}\\
&\U_1=\compl{\X}{\U_2},\label{thmconnectednessandopenpartitionsp1eq4}\\
&\opair{\U_1}{\U_2}\in\topology{}\times\topology{}.\label{thmconnectednessandopenpartitionsp1eq5}
\end{align}
According to \refdef{deffamilyofclosedsets},
\Ref{thmconnectednessandopenpartitionsp1eq4}
and
\Ref{thmconnectednessandopenpartitionsp1eq5}
imply that,
\begin{equation}\label{thmconnectednessandopenpartitionsp1eq6}
\U_1\in\bigg(\topology{}\cap\Fclosed{\X}{\topology{}}\bigg).
\end{equation}
\Ref{thmconnectednessandopenpartitionsp1eq3}
and
\Ref{thmconnectednessandopenpartitionsp1eq6}
imply,
\begin{equation}
\bigg(\topology{}\cap\Fclosed{\X}{\topology{}}\bigg)\neq
\seta{\binary{\empty}{\X}},
\end{equation}
which according to \refdef{defconnectedness}, means,
$\Xt$
is disconnected.
\endp
\end{itemize}
\begin{itemize}
\item[${\textbf{\textsf{p2}}}$]
It is assumed that
$\Xt$
is not connected. Then according to, \refdef{defconnectedness},
\begin{equation}
\Existsis{\U}{\bigg(\topology{}\cap\Fclosed{\X}{\topology{}}\bigg)}
\U\notin\seta{\binary{\empty}{\X}}.
\end{equation}
Then according to \refdef{defpartitionofset},
\begin{equation}
\seta{\binary{\U}{\compl{\X}{\U}}}\in\Cpart{\X},
\end{equation}
and according to \refdef{deffamilyofclosedsets}
(considering that $\U\in\Fclosed{\X}{\topology{}}$),
\begin{equation}
\compl{\X}{\U}\in\topology{}.
\end{equation}
Therefore, based on \refdef{defopenpartition},
\begin{equation}
\seta{\binary{\U}{\compl{\X}{\U}}}\in
\defset{\apartition}{\Opart{\Xt}}
{\CarD{\apartition}=2},
\end{equation}
which means,
\begin{equation}
\defset{\apartition}{\Opart{\Xt}}
{\CarD{\apartition}=2}\neq\empty.
\end{equation}
\endp
\end{itemize}
\endthm
%%%%%%%%%%%%%%%%%%%%%%%%%%%%%%%%%%%%%%%%%%%%%%%%%%%%%%%%%%%%%%%%%%%%%%%%%%%%%%%%%%%%%%%%%%%%%%%%%%%%%%%%%%%%%%%%%%%%%%%%%%%%%%%%
\theorem\label{thmconnectednessandclosedpartitions}
$\Xt=\opair{\X}{\topology{}}$
is taken as a topological-space.
$\Xt$
is a connected topological-space if and only if
$\Xt$
does not have any (non-trivial) closed partition consisting of two elements. That is,
\begin{equation}\label{thmconnectednessandclosedpartitionseq1}
\bigg(\X\in\connecteds{\Xt}\bigg)\thenn
\bigg(\defset{\apartition}{\Clpart{\Xt}}
{\CarD{\apartition}=2}=\empty\bigg).
\end{equation}
\prooff
Every open partition of $\Xt$ with two elements is a closed partition of $\Xt$ with two elements, or vice versa.
This can be verified easily according to
\refdef{deffamilyofclosedsets},
\refdef{defopenpartition},
and
\refdef{defclosedpartition}.
Thus,
\begin{equation}
\defset{\apartition}{\Opart{\Xt}}
{\CarD{\apartition}=2}
=
\defset{\apartition}{\Clpart{\Xt}}
{\CarD{\apartition}=2}.
\end{equation}
Hence, according to \refthm{thmconnectednessandopenpartitions},
\Ref{thmconnectednessandclosedpartitionseq1}
is obtained.
\endthm
%%%%%%%%%%%%%%%%%%%%%%%%%%%%%%%%%%%%%%%%%%%%%%%%%%%%%%%%%%%%%%%%%%%%%%%%%%%%%%%%%%%%%%%%%%%%%%%%%%%%%%%%%%%%%%%%%%%%%%%%%%%%%%%%
\theorem\label{thmconnectednessandnonclingingpartitions}
$\Xt=\opair{\X}{\topology{}}$
is taken as a topological-space.
$\Xt$
is a connected topological-space if and only if
$\Xt$
does not have any (non-trivial) non-attached partition consisting of two elements. That is,
\begin{equation}
\bigg(\X\in\connecteds{\Xt}\bigg)\thenn
\bigg(\defset{\apartition}{\Ncpart{\Xt}}
{\CarD{\apartition}=2}=\empty\bigg).
\end{equation}
\prooff
\begin{itemize}
\item[${\textbf{\textsf{p1}}}$]
It is assumed that
$\Xt$
is disconnected:
\begin{equation}
\X\notin\connecteds{\Xt}.
\end{equation}
Then according to \refdef{thmconnectednessandopenpartitions},
there exists at least one open partition of $\Xt$ with two elements, like
$\seta{\binary{\U}{\(\compl{\X}{\U}\)}}$.
It is clear that,
$\U$
is also a closed set of $\Xt$ (because $\compl{\X}{\U}$ is an open set of $\Xt$) Therefore, according to \refthm{thmclosureofclosedset},
\begin{equation}
\U=\func{\Cl{\Xt}}{\U},
\end{equation}
and hence,
\begin{align}
\(\compl{\X}{\U}\)\cap\func{\Cl{\Xt}}{\U}&=
\(\compl{\X}{\U}\)\cap\U\cr
&=\empty.
\end{align}
Thus, according to \refthm{defnonclingingpartition},
$\seta{\binary{\U}{\(\compl{\X}{\U}\)}}$
is a non-attached partition of $\Xt$ with two elements. Therefore,
\begin{equation*}
\defset{\apartition}{\Ncpart{\Xt}}
{\CarD{\apartition}=2}\neq\empty.
\end{equation*}
\endp
\end{itemize}
\begin{itemize}
\item[${\textbf{\textsf{p2}}}$]
It is assumed that,
\begin{equation}
\defset{\apartition}{\Ncpart{\Xt}}
{\CarD{\apartition}=2}\neq\empty.
\end{equation}
Then, there exists at least one non-attached partition of $\Xt$ with two elements,
like $\seta{\binary{\U}{\(\compl{\X}{\U}\)}}$.
According to \refdef{defnonclingingpartition},
\begin{align}
\U\cap\func{\Cl{\Xt}}{\compl{\X}{\U}}&=\empty,\\
\(\compl{\X}{\U}\)\cap\func{\Cl{\Xt}}{\U}&=\empty.
\end{align}
Therefore,
\begin{align}
\func{\Cl{\Xt}}{\U}&=\U,\\
\func{\Cl{\Xt}}{\compl{\X}{\U}}&=\(\compl{\X}{\U}\),
\end{align}
which according to \refthm{thmclosureofclosedset}, means,
\begin{equation}
\opair{\U}{\compl{\X}{\U}}\in\Fclosed{\X}{\topology{}}\times\Fclosed{\X}{\topology{}},
\end{equation}
and hence according to \refdef{defclosedpartition},
$\seta{\binary{\Y}{\compl{\X}{\U}}}$
is a closed partition of $\Xt$ with two elements. That is,
\begin{equation}
\seta{\binary{\U}{\compl{\X}{\U}}}\in
\defset{\apartition}{\Clpart{\Xt}}
{\CarD{\apartition}=2}.
\end{equation}
Thus,
\begin{equation}
\defset{\apartition}{\Clpart{\Xt}}
{\CarD{\apartition}=2}\neq\empty,
\end{equation}
and hence according to \refthm{thmconnectednessandclosedpartitions},
$\Xt$
is disconnected. That is,
\begin{equation*}
\X\notin\connecteds{\Xt}.
\end{equation*}
\end{itemize}
\endthm
%%%%%%%%%%%%%%%%%%%%%%%%%%%%%%%%%%%%%%%%%%%%%%%%%%%%%%%%%%%%%%%%%%%%%%%%%%%%%%%%%%%%%%%%%%%%%%%%%%%%%%%%%%%%%%%%%%%%%%%%%%%%%%%%
\theorem\label{thmconnectedsubsetsofdisconnectedspace}
$\Xt=\opair{\X}{\topology{}}$
is taken as a connected topological-space, and
$\seta{\binary{\covelm}{\p{\covelm}}}$
as an open partition of $\Xt$ with two elements
(an element os
$\defset{\apartition}{\Opart{\Xt}}{\CarD{\apartition}=2}$).
Every connected set of $\Xt$
is a subset of one of
$\covelm$
or
$\p{\covelm}$. That is,
\begin{equation}
\Foreach{\asubset}{\connecteds{\Xt}}
\bigg[\OR{\asubset\subseteq\covelm}{\asubset\subseteq\p{\covelm}}\bigg].
\end{equation}
\proof
$\asubset$
is taken as such a subset of $\X$
that is neither a subset of $\covelm$ nor a subset of $\p{\covelm}$, that is,
\begin{equation}
\AND{\asubset\nsubseteq\covelm}{\asubset\nsubseteq\p{\covelm}}.
\end{equation}
Then considering that
$\seta{\binary{\covelm}{\p{\covelm}}}$
is a cover of $\X$, it is clear that,
\begin{align}
\asubset\cap\covelm&\neq\empty,\\
\asubset\cap\p{\covelm}&\neq\empty.
\end{align}
Additionally, considering that
$\covelm$
and
$\p{\covelm}$
do not intersect each other, it is clear that,
\begin{equation}
\(\asubset\cap\covelm\)\cap\(\asubset\cap\p{\covelm}\)=\empty.
\end{equation}
In addition, it is clear that,
\begin{align}
\(\asubset\cap\covelm\)\cup\(\asubset\cap\p{\covelm}\)&=
\asubset\cap\(\covelm\cup\p{\covelm}\)\cr
&=\asubset\cap\X\cr
&=\asubset.
\end{align}
Therefore, according to \refdef{defpartitionofset}, it is evident that,
$\seta{\binary{\asubset\cap\covelm}{\asubset\cap\p{\covelm}}}$
is a partition of $\X$ with two elements, that is,
\begin{equation}
\seta{\binary{\asubset\cap\covelm}{\asubset\cap\p{\covelm}}}\in
\Cpart{\X}.
\end{equation}
In addition, considering that
$\covelm$
and
$\p{\covelm}$
are open sets of $\Xt$, according to \refdef{defsubspacetopology1},
\begin{align}
\(\asubset\cap\covelm\)&\in\stopology{\topology{}}{\asubset},\\
\(\asubset\cap\p{\covelm}\)&\in\stopology{\topology{}}{\asubset}.
\end{align}
Therefore, according to \refdef{defopenpartition},
it is evident that
$\seta{\binary{\asubset\cap\covelm}{\asubset\cap\p{\covelm}}}$
is an open partition of $\Xt$ with two elements, that is,
\begin{equation}
\seta{\binary{\asubset\cap\covelm}{\asubset\cap\p{\covelm}}}\in
\defset{\apartition}{\Opart{\opair{\asubset}{\stopology{\topology{}}{\asubset}}}}{\CarD{\apartition}=2},
\end{equation}
and therefore,
\begin{equation}
\defset{\apartition}{\Opart{\opair{\asubset}
{\stopology{\topology{}}{\asubset}}}}{\CarD{\apartition}=2}\neq\empty.
\end{equation}
According to \refthm{thmconnectednessandopenpartitions}, this means,
$\asubset$
is a disconnected set of $\Xt$, that is,
\begin{equation}
\asubset\notin\connecteds{\Xt}.
\end{equation}
\endthm
%%%%%%%%%%%%%%%%%%%%%%%%%%%%%%%%%%%%%%%%%%%%%%%%%%%%%%%%%%%%%%%%%%%%%%%%%%%%%%%%%%%%%%%%%%%%%%%%%%%%%%%%%%%%%%%%%%%%%%%%%%%%%%%%
\theorem\label{thmclopenandconnectedsetsofspace}
$\Xt=\opair{\X}{\topology{}}$
is taken as a topological-space.
For every pair of subsets $\opair{\csubset}{\asubset}$ of $\X$,if
$\csubset$
is a connected set of $\Xt$, and $\asubset$
is am open-and-closed set of $\Xt$, then $\asubset$ includes either $\csubset$ or $\compl{\X}{\asubset}$. That is,
\begin{equation}
\Foreach{\opair{\csubset}{\asubset}}
{\bigg(\connecteds{\Xt}\times\[\topology{}\cap\Fclosed{\X}{\topology{}}\]\bigg)}
\bigg[\OR{\csubset\subseteq\asubset}{\csubset\subseteq\(\compl{\X}{\asubset}\)}\bigg].
\end{equation}
\proof
$\csubset$
is taken as an arbitrary element of $\connecteds{\Xt}$
(a connected set of $\Xt$),
and $\asubset$ as an arbitrary element of $\[\topology{}\cap\Fclosed{\X}{\topology{}}\]$
(an open-and-closed set of $\Xt$). Then according to \refdef{defconnectedness},
\begin{equation}
\[\stopology{\topology{}}{\csubset}\cap
\Fclosed{\X}{\stopology{\topology{}}{\csubset}}\]=\seta{\binary{\empty}{\csubset}},
\end{equation}
and according to \refdef{defsubspacetopology1}
and
\refthm{thmsubspaceclosedsets},
\begin{equation}
\(\asubset\cap\csubset\)\in
\[\stopology{\topology{}}{\csubset}\cap
\Fclosed{\X}{\stopology{\topology{}}{\csubset}}\].
\end{equation}
Therefore,
\begin{equation}
\(\asubset\cap\csubset\)\in\seta{\binary{\empty}{\csubset}},
\end{equation}
which means,
\begin{equation}
\OR{\csubset\subseteq\asubset}{\csubset\subseteq\(\compl{\X}{\asubset}\)}.
\end{equation}
\endthm
%%%%%%%%%%%%%%%%%%%%%%%%%%%%%%%%%%%%%%%%%%%%%%%%%%%%%%%%%%%%%%%%%%%%%%%%%%%%%%%%%%%%%%%%%%%%%%%%%%%%%%%%%%%%%%%%%%%%%%%%%%%%%%%%
\theorem\label{thmclosureofconnectedset}
$\Xt=\opair{\X}{\topology{}}$
is taken as a topological-space.
For every
$\asubset$
in
$\CSs{\X}$,
if
$\asubset$
is a connected set of
$\Xt$, then every set including $\asubset$
and included in $\func{\Cl{\Xt}}{\asubset}$
(the closure of $\asubset$ in
$\Xt$)
is also a connected set of $\Xt$. That is,
\begin{equation}
\Foreach{\asubset}{\connecteds{\Xt}}
\defset{\bsubset}{\CSs{\X}}{\[\asubset\subseteq\bsubset\subseteq\func{\Cl{\Xt}}{\asubset}\]}
\subseteq\connecteds{\Xt}.
\end{equation}
\proof
$\asubset$
is taken as an arbitrary element of $\connecteds{\Xt}$
(a connected set of $\Xt$). Then according to \refdef{defconnectedness},
\begin{equation}\label{thmclosureofconnectedsetpeq1}
\stopology{\topology{}}{\asubset}\cap\Fclosed{\asubset}{\stopology{\topology{}}{\asubset}}=
\seta{\binary{\empty}{\asubset}}.
\end{equation}
$\bsubset$
is taken as such a subset of $\X$ that,
\begin{equation}\label{thmclosureofconnectedsetpeq2}
\asubset\subseteq\bsubset\subseteq\func{\Cl{\Xt}}{\asubset}.
\end{equation}
\begin{itemize}
\item[${\textbf{\textsf{p1}}}$]
$\U$
is taken as an arbitrary element of
$\compl{\[\stopology{\topology{}}{\bsubset}\cap\Fclosed{\bsubset}{\stopology{\topology{}}{\bsubset}}\]}{\seta{\empty}}$
(a non-empty open-and-closed set of the topological-space
$\opair{\bsubset}{\stopology{\topology{}}{\bsubset}}$). Then according to \refdef{defsubspacetopology1},
\begin{equation}\label{thmclosureofconnectedsetp1eq1}
\Existsis{\V}{\(\compl{\topology{}}{\seta{\empty}}\)}
\U=\B\cap\V.
\end{equation}
Therefore, considering that
$\U\subseteq\func{\Cl{\Xt}}{\asubset}$, and $\U\neq\empty$,
it is clear that,
\begin{equation}\label{thmclosureofconnectedsetp1eq2}
\func{\Cl{\Xt}}{\asubset}\cap\V\neq\empty.
\end{equation}
Hence, according to
\refthm{thmopensetsintersectingclosure},
and considering that $\V\in\topology{}$,
\begin{equation}\label{thmclosureofconnectedsetp1eq3}
\asubset\cap\V\neq\empty.
\end{equation}
In addition, according to \Ref{thmclosureofconnectedsetp1eq1},
and considering that
$\asubset\subseteq\bsubset$,
\begin{align}\label{thmclosureofconnectedsetp1eq4}
\asubset\cap\U&=\asubset\cap\(\bsubset\cap\V\)\cr
&=\(\asubset\cap\bsubset\)\cap\V\cr
&=\asubset\cap\V.
\end{align}
\Ref{thmclosureofconnectedsetp1eq3}
and
\Ref{thmclosureofconnectedsetp1eq4}
imply,
\begin{equation}\label{thmclosureofconnectedsetp1eq5}
\asubset\cap\U\neq\empty.
\end{equation}
Moreover, considering that $\asubset\subseteq\bsubset$,
and ($\U$
is an open-and-closed set of $\opair{\bsubset}{\stopology{\topology{}}{\bsubset}}$)
$\U\in\[\stopology{\topology{}}{\bsubset}\cap\Fclosed{\bsubset}{\stopology{\topology{}}{\bsubset}}\]$,
according to
\refdef{defsubspacetopology1}
and
\refthm{thmsubspaceclosedsets},
it is clear that,
\begin{equation}\label{thmclosureofconnectedsetp1eq6}
\(\asubset\cap\U\)\in\[\stopology{\stopology{\topology{}}{\bsubset}}{\asubset}\cap
\Fclosed{\asubset}{\stopology{\stopology{\topology{}}{\bsubset}}{\asubset}}\],
\end{equation}
and hence according to \refthm{thmsubspacetopologytransitivity},
($\stopology{\stopology{\topology{}}{\bsubset}}{\asubset}=\stopology{\topology{}}{\asubset}$),
\begin{equation}\label{thmclosureofconnectedsetp1eq7}
\(\asubset\cap\U\)\in\[\stopology{\topology{}}{\asubset}
\cap\Fclosed{\asubset}{\stopology{\topology{}}{\asubset}}\].
\end{equation}
\Ref{thmclosureofconnectedsetpeq1},
\Ref{thmclosureofconnectedsetp1eq5},
and
\Ref{thmclosureofconnectedsetp1eq7}
imply that,
\begin{equation}\label{thmclosureofconnectedsetp1eq8}
\(\asubset\cap\U\)=\asubset,
\end{equation}
which means,
\begin{equation}\label{thmclosureofconnectedsetp1eq9}
\asubset\subseteq\U,
\end{equation}
According to
\refcor{corclosureofset0},
and considering that (
$\U$
is a closed set of $\opair{\bsubset}{\stopology{\topology{}}{\bsubset}}$)
$\U\in\Fclosed{\bsubset}{\stopology{\topology{}}{\bsubset}}$,
\Ref{thmclosureofconnectedsetp1eq9}
implies that the closure of $\asubset$ in the topological-space
$\opair{\bsubset}{\stopology{\topology{}}{\bsubset}}$
is a subset of $\U$. That is,
\begin{equation}\label{thmclosureofconnectedsetp1eq10}
\func{\Cl{\opair{\bsubset}{\stopology{\topology{}}{\bsubset}}}}{\asubset}\subseteq\U.
\end{equation}
Moreover, according to \refthm{thmsubspaceclosure},
and considering that $\asubset\subseteq\bsubset$,
\begin{equation}\label{thmclosureofconnectedsetp1eq11}
\func{\Cl{\opair{\bsubset}{\stopology{\topology{}}{\bsubset}}}}{\asubset}=
\bsubset\cap\func{\Cl{\Xt}}{\asubset},
\end{equation}
and thus considering that
$\bsubset\subseteq\func{\Cl{\Xt}}{\asubset}$,
\begin{equation}\label{thmclosureofconnectedsetp1eq12}
\func{\Cl{\opair{\bsubset}{\stopology{\topology{}}{\bsubset}}}}{\asubset}=\bsubset.
\end{equation}
\Ref{thmclosureofconnectedsetp1eq10}
and
\Ref{thmclosureofconnectedsetp1eq12}
yield,
\begin{equation}
\bsubset\subseteq\U,
\end{equation}
and hence considering that $\U\subseteq\bsubset$
(because $\U$ is an open set of $\opair{\bsubset}{\stopology{\topology{}}{\bsubset}}$), it is evident that,
\begin{equation}
\U=\bsubset.
\end{equation}
\endp
\end{itemize}
Therefore,
\begin{equation}
\[\stopology{\topology{}}{\bsubset}\cap\Fclosed{\bsubset}
{\stopology{\topology{}}{\bsubset}}\]=\seta{\binary{\empty}{\bsubset}},
\end{equation}
which according to \refdef{defconnectedness}, means,
$\bsubset$
is a connected set of $\Xt$, that is,
\begin{equation}
\bsubset\in\connecteds{\Xt}.
\end{equation}
\endthm
%%%%%%%%%%%%%%%%%%%%%%%%%%%%%%%%%%%%%%%%%%%%%%%%%%%%%%%%%%%%%%%%%%%%%%%%%%%%%%%%%%%%%%%%%%%%%%%%%%%%%%%%%%%%%%%%%%%%%%%%%%%%%%%%
\theorem\label{thmclosureofconnectedsetisconnected}
$\Xt=\opair{\X}{\topology{}}$
is taken as a topological-space.
For every $\asubset$ in $\CSs{\X}$, if $\asubset$ is a connected set of $\Xt$, then
$\func{\Cl{\Xt}}{\asubset}$
(the closure of $\asubset$ in $\Xt$)
is also a connected set of $\Xt$. That is,
\begin{equation}
\Foreach{\asubset}{\connecteds{\Xt}}
\func{\Cl{\Xt}}{\asubset}\in\connecteds{\Xt}.
\end{equation}
\proof
According to \refcor{corclosureofset0}
and
\refthm{thmclosureofconnectedset},
it is obvious.
\endthm
%%%%%%%%%%%%%%%%%%%%%%%%%%%%%%%%%%%%%%%%%%%%%%%%%%%%%%%%%%%%%%%%%%%%%%%%%%%%%%%%%%%%%%%%%%%%%%%%%%%%%%%%%%%%%%%%%%%%%%%%%%%%%%%%
\theorem\label{thmunionofafamilyofconnectedsets}
$\Xt=\opair{\X}{\topology{}}$
is taken as a topological-space.
The union of elements of any collection of connected sets of $\Xt$ with pairwise intersecting elements
is a connected set of $\Xt$. That is,
\begin{align}
&\defsets{\sCi}{\connecteds{\Xt}}
{\Foreach{\opair{\Asubset{1}}{\Asubset{2}}}{\sCi\times\sCi}\Asubset{1}\cap\Asubset{2}\neq\empty}\cr
\subseteq
&\defsets{\sCi}{\connecteds{\Xt}}{\(\union{\sCi}\)\in\connecteds{\Xt}}.
\end{align}
\proof
$\sCi$
is taken as such a subset of $\connecteds{\Xt}$ that,
\begin{equation}\label{thmunionofafamilyofconnectedsetspeq1}
\Foreach{\opair{\Asubset{1}}{\Asubset{2}}}{\sCi\times\sCi}
\Asubset{1}\cap\Asubset{2}\neq\empty.
\end{equation}
Considering that $\sCi$
is a subset of $\connecteds{\Xt}$, according to \refdef{defconnectedness},
\begin{equation}\label{thmunionofafamilyofconnectedsetspeq2}
\Foreach{\asubset}{\sCi}
\stopology{\topology{}}{\asubset}\cap\Fclosed{\asubset}{\stopology{\topology{}}{\asubset}}=\seta{\binary{\empty}{\asubset}}.
\end{equation}
\begin{itemize}
\item[${\textbf{\textsf{p1}}}$]
$\U$
is taken  as an arbitrary element of
$\compl{\bigg[\stopology{\topology{}}{\union{\sCi}}\cap
\Fclosed{\union{\sCi}}{\stopology{\topology{}}
{\union{\sCi}}}\bigg]}{\seta{\empty}}$.
Then considering that,
\begin{equation}\label{thmunionofafamilyofconnectedsetsp1eq1}
\Foreach{\asubset}{\sCi}
\asubset\subseteq\(\union{\sCi}\),
\end{equation}
according to \refdef{defsubspacetopology1} and \refthm{thmsubspaceclosedsets},
it is clear that,
\begin{equation}\label{thmunionofafamilyofconnectedsetsp1eq2}
\Foreach{\asubset}{\sCi}
\(\asubset\cap\U\)\in\stopology{\stopology{\topology{}}{\union{\sCi}}}{\asubset}
\cap
\Fclosed{\asubset}{\stopology{\stopology{\topology{}}{\union{\sCi}}}{\asubset}}.
\end{equation}
Moreover, according to \refthm{thmsubspacetopologytransitivity},
\begin{equation}\label{thmunionofafamilyofconnectedsetsp1eq3}
\Foreach{\asubset}{\sCi}
\stopology{\stopology{\topology{}}{\union{\sCi}}}{\asubset}=\stopology{\topology{}}{\asubset}.
\end{equation}
Therefore,
\begin{equation}\label{thmunionofafamilyofconnectedsetsp1eq4}
\Foreach{\asubset}{\sCi}
\(\asubset\cap\U\)\in\stopology{\topology{}}{\asubset}\cap
\Fclosed{\asubset}{\stopology{\topology{}}{\asubset}}.
\end{equation}
Thus, according to \Ref{thmunionofafamilyofconnectedsetspeq2},
\begin{equation}\label{thmunionofafamilyofconnectedsetsp1eq5}
\Foreach{\asubset}{\sCi}
\(\asubset\cap\U\)\in\seta{\binary{\empty}{\asubset}}.
\end{equation}
Furthermore, considering that,
\begin{align}\label{thmunionofafamilyofconnectedsetsp1eq6}
\Union{\asubset}{\sCi}{\(\asubset\cap\U\)}&=
\(\union{\sCi}\)\cap\U\cr
&=\U,
\end{align}
and that
$\U$ is a non-empty set, it is evident that,
\begin{equation}\label{thmunionofafamilyofconnectedsetsp1eq7}
\Union{\asubset}{\sCi}{\(\asubset\cap\U\)}\neq\empty,
\end{equation}
and hence according to \Ref{thmunionofafamilyofconnectedsetsp1eq5},
\begin{equation}\label{thmunionofafamilyofconnectedsetsp1eq8}
\Existsis{\Asubset{0}}{\sCi}
\(\Asubset{0}\cap\U\)=\Asubset{0}.
\end{equation}
According to \Ref{thmunionofafamilyofconnectedsetspeq1},
\begin{equation}\label{thmunionofafamilyofconnectedsetsp1eq9}
\Foreach{\asubset}{\sCi}
\(\asubset\cap\Asubset{0}\)\neq\empty.
\end{equation}
Thus, according to \Ref{thmunionofafamilyofconnectedsetsp1eq8},
\begin{align}\label{thmunionofafamilyofconnectedsetsp1eq9}
\Foreach{\asubset}{\sCi}
\empty&\neq\asubset\cap\Asubset{0}\cr
&=\asubset\cap\(\Asubset{0}\cap\U\)\cr
&=\(\asubset\cap\U\)\cap\Asubset{0},
\end{align}
and hence,
\begin{equation}\label{thmunionofafamilyofconnectedsetsp1eq10}
\Foreach{\asubset}{\sCi}
\(\asubset\cap\U\)\neq\empty.
\end{equation}
Hence, according to \Ref{thmunionofafamilyofconnectedsetsp1eq5},
\begin{equation}\label{thmunionofafamilyofconnectedsetsp1eq11}
\Foreach{\asubset}{\sCi}
\(\asubset\cap\U\)=\asubset.
\end{equation}
Therefore,
\begin{equation}\label{thmunionofafamilyofconnectedsetsp1eq12}
\Union{\asubset}{\sCi}{\(\asubset\cap\U\)}=\(\union{\sCi}\).
\end{equation}
Thus, according to \Ref{thmunionofafamilyofconnectedsetsp1eq6},
\begin{equation}
\U=\(\union{\sCi}\).
\end{equation}
\endp
\end{itemize}
Therefore,
\begin{equation}
\bigg[\stopology{\topology{}}{\union{\sCi}}\cap
\Fclosed{\union{\sCi}}{\stopology{\topology{}}
{\union{\sCi}}}\bigg]=\seta{\binary{\empty}{\union{\sCi}}},
\end{equation}
which means
$\opair{\union{\sCi}}{\stopology{\topology{}}{\union{\sCi}}}$
is connected, that is,
\begin{equation}
\(\union{\sCi}\)\in\connecteds{\Xt}.
\end{equation}
\endthm
%%%%%%%%%%%%%%%%%%%%%%%%%%%%%%%%%%%%%%%%%%%%%%%%%%%%%%%%%%%%%%%%%%%%%%%%%%%%%%%%%%%%%%%%%%%%%%%%%%%%%%%%%%%%%%%%%%%%%%%%%%%%%%%%
\theorem\label{thmunionofafamilyofconnectedsets1}
$\Xt=\opair{\X}{\topology{}}$
is taken as a topological-space.
The union of elements of any collection of connected sets of $\Xt$ with non-empty intersection of all elements,
is a connected set of $\Xt$. That is,
\begin{align}\label{thmunionofafamilyofconnectedsets1eq1}
&\defsets{\sCi}{\connecteds{\Xt}}
{\(\intersection{\sCi}\)\neq\empty}\cr
\subseteq
&\defsets{\sCi}{\connecteds{\Xt}}{\(\union{\sCi}\)\in\connecteds{\Xt}}.
\end{align}
$\caution$
$\(\intersection{\sCi}\)$
is not defined for $\sCi=\empty$. This means,
$\empty$
is not an element of
$\defsets{\sCi}{\connecteds{\Xt}}
{\(\intersection{\sCi}\)\neq\empty}$.
\proof
It is clear that for every subset $\sCi$ of $\connecteds{\Xt}$,
if the intersection of elements of $\sCi$ is non-empty, then the elements of $\sCi$
are pairwise intersecting, that is,
\begin{align}
&\defsets{\sCi}{\connecteds{\Xt}}
{\(\intersection{\sCi}\)\neq\empty}\cr
\subseteq
&\defsets{\sCi}{\connecteds{\Xt}}
{\Foreach{\opair{\Asubset{1}}{\Asubset{2}}}{\sCi\times\sCi}\Asubset{1}\cap\Asubset{2}\neq\empty}.
\end{align}
Therefore, base on \refthm{thmunionofafamilyofconnectedsets},
\Ref{thmunionofafamilyofconnectedsets1eq1}
is obtained
\endthm
%%%%%%%%%%%%%%%%%%%%%%%%%%%%%%%%%%%%%%%%%%%%%%%%%%%%%%%%%%%%%%%%%%%%%%%%%%%%%%%%%%%%%%%%%%%%%%%%%%%%%%%%%%%%%%%%%%%%%%%%%%%%%%%%
\theorem\label{thmunionofafamilyofconnectedsets2}
$\Xt=\opair{\X}{\topology{}}$
is taken as a topological-space.
The union of elements of any countable collection
$\Zfamily{\Asubset{k}}{k\in\Zp}$ of connected sets of $\Xt$ indexed by $\Zp$
such that $\Asubset{k}\cap\Asubset{k+1}\neq\empty$ for any $k$,
is a connected set of $\Xt$. That is,
\begin{align}
&\defsets{\sCi}{\connecteds{\Xt}}
{\bigg[\Exists{\indexf}{\surFunc{\Zp}{\sCi}}\bigg(
\Foreach{k}{\Zp}\func{\indexf}{k}\cap\func{\indexf}{k+1}\neq\empty\bigg)\bigg]}\cr
\subseteq
&\defsets{\sCi}{\connecteds{\Xt}}{\(\union{\sCi}\)\in\connecteds{\Xt}}.
\end{align}
\proof
$\sCi$
is taken as a countable subset of $\connecteds{\Xt}$ and it is assumed that there is an index function $\indexf$ in
$\surFunc{\Zp}{\sCi}$
(the set of all surjective functions from $\Zp$ to $\sCi$)
such that,
\begin{equation}\label{thmunionofafamilyofconnectedsets2peq1}
\Foreach{k}{\Zp}
\func{\indexf}{k}\cap\func{\indexf}{k+1}\neq\empty.
\end{equation}
Then it is clear that,
\begin{align}
&\Foreach{k}{\Zp}\func{\indexf}{k}\neq\empty,\label{thmunionofafamilyofconnectedsets2peq2}\\
&\Foreach{k}{\Zp}\func{\indexf}{k}\in\connecteds{\Xt}.\label{thmunionofafamilyofconnectedsets2peq3}
\end{align}
The mapping
$\mu$
is defined as,
\begin{gather}\label{thmunionofafamilyofconnectedsets2peq4}
\begin{aligned}
&\mu\indef\Func{\Zp}{\CSs{\X}},\\
&\Foreach{n}{\Zp}\func{\mu}{n}\eqdef
\[\fUnion{k}{1}{n}{\func{\indexf}{k}}\]
\end{aligned}
\end{gather}
It is clear that,
\begin{equation}\label{thmunionofafamilyofconnectedsets2peq5}
\func{\mu}{1}=\func{\indexf}{1},
\end{equation}
and hence according to,
\Ref{thmunionofafamilyofconnectedsets2peq3},
\begin{equation}\label{thmunionofafamilyofconnectedsets6}
\func{\mu}{1}\in\connecteds{\Xt}.
\end{equation}
\begin{itemize}
\item[${\textbf{\textsf{p1}}}$]
$n$
is taken as an arbitrary element of $\Zp$, and it is assumed that,
$\func{\mu}{n}$
is a connected set of $\Xt$:
\begin{equation}\label{thmunionofafamilyofconnectedsets2p1eq1}
\func{\mu}{n}\in\connecteds{\Xt}.
\end{equation}
It is clear that,
\begin{align}\label{thmunionofafamilyofconnectedsets2p1eq2}
\func{\mu}{n+1}&=\[\fUnion{k}{1}{n+1}{\func{\indexf}{k}}\]\cr
&=\func{\indexf}{n+1}\cup\[\fUnion{k}{1}{n}{\func{\indexf}{k}}\]\cr
&=\func{\indexf}{n+1}\cup\func{\mu}{n}.
\end{align}
In addition, according to
\Ref{thmunionofafamilyofconnectedsets2peq1},
\begin{align}\label{thmunionofafamilyofconnectedsets2p1eq3}
\func{\indexf}{n+1}\cap\func{\mu}{n}&=
\func{\indexf}{n+1}\cap\[\fUnion{k}{1}{n}{\func{\indexf}{k}}\]\cr
&=\fUnion{k}{1}{n}{\bigg[\func{\indexf}{n+1}\cap\func{\indexf}{k}\bigg]}\cr
&\neq\empty.
\end{align}
Therefore according to
\Ref{thmunionofafamilyofconnectedsets2peq3},
\Ref{thmunionofafamilyofconnectedsets2p1eq1},
and
\Ref{thmunionofafamilyofconnectedsets2p1eq3},
it is evident that,
$\func{\indexf}{n+1}$
and
$\func{\mu}{n}$
are two intersecting connected sets.
Thus according to
\refthm{thmunionofafamilyofconnectedsets1},
$\[\func{\indexf}{n+1}\cup\func{\mu}{n}\]$
is a connected set of $\Xt$.
\begin{equation}
\[\func{\indexf}{n+1}\cup\func{\mu}{n}\]\in\connecteds{\Xt}.
\end{equation}
Therefore according to \Ref{thmunionofafamilyofconnectedsets2p1eq2}
it is clear that,
$\func{\mu}{n+1}$
is a connected set of $\Xt$.
\begin{equation}
\func{\mu}{n+1}\in\connecteds{\Xt}.
\end{equation}
\endp
\end{itemize}
Therefore,
\begin{equation}\label{thmunionofafamilyofconnectedsets2peq7}
\Foreach{n}{\Zp}
\bigg[\func{\mu}{n}\in\connecteds{\Xt}\then
\func{\mu}{n+1}\in\connecteds{\Xt}\bigg].
\end{equation}
Based on induction,
\Ref{thmunionofafamilyofconnectedsets2peq5}
and
\Ref{thmunionofafamilyofconnectedsets2peq7}
imply that every
$\func{\mu}{n}$
is a connected set of
$\Xt$.
\begin{equation}\label{thmunionofafamilyofconnectedsets2peq8}
\Foreach{n}{\Zp}\func{\mu}{n}\in\connecteds{\Xt}.
\end{equation}
According to
\Ref{thmunionofafamilyofconnectedsets2peq4},
\begin{equation}\label{thmunionofafamilyofconnectedsets2peq9}
\Foreach{n}{\Zp}
\func{\indexf}{1}\subseteq\func{\mu}{n},
\end{equation}
and thus,
\begin{align}\label{thmunionofafamilyofconnectedsets2peq10}
\func{\indexf}{1}\subseteq
\[\Intersection{n}{\Zp}{\func{\mu}{n}}\],
\end{align}
and hence according to \Ref{thmunionofafamilyofconnectedsets2peq2},
\begin{equation}\label{thmunionofafamilyofconnectedsets2peq11}
\[\Intersection{n}{\Zp}{\func{\mu}{n}}\]\neq\empty.
\end{equation}
\Ref{thmunionofafamilyofconnectedsets2peq8}
and
\Ref{thmunionofafamilyofconnectedsets2peq11}
imply that,
$\defSet{\func{\mu}{n}}{n\in\Zp}$
is a collection of connected sets of $\Xt$ whose intersection is non-empty.
\begin{equation}\label{thmunionofafamilyofconnectedsets2peq12}
\defSet{\func{\mu}{n}}{n\in\Zp}\in
\defsets{\sCi}{\connecteds{\Xt}}
{\(\intersection{\sCi}\)\neq\empty}.
\end{equation}
Thus according to \refthm{thmunionofafamilyofconnectedsets1},
\begin{equation}\label{thmunionofafamilyofconnectedsets2peq13}
\[\Union{n}{\Zp}{\func{\mu}{n}}\]\in\connecteds{\Xt}.
\end{equation}
In addition, according to \Ref{thmunionofafamilyofconnectedsets2peq4},
it is clear that,
\begin{equation}\label{thmunionofafamilyofconnectedsets2peq14}
\[\Union{n}{\Zp}{\func{\mu}{n}}\]=
\[\Union{n}{\Zp}\func{\indexf}{n}\].
\end{equation}
\Ref{thmunionofafamilyofconnectedsets2peq13}
and
\Ref{thmunionofafamilyofconnectedsets2peq14}
imply that,
\begin{equation}\label{thmunionofafamilyofconnectedsets2peq15}
\[\Union{n}{\Zp}\func{\indexf}{n}\]\in\connecteds{\Xt}.
\end{equation}
Considering that $\indexf$ is surjective,
it is clear that,
\begin{equation}\label{thmunionofafamilyofconnectedsets2peq16}
\[\Union{n}{\Zp}\func{\indexf}{n}\]=\union{\sCi}.
\end{equation}
Therefore,
\begin{equation}
\(\union{\sCi}\)\in\connecteds{\Xt}.
\end{equation}
\endthm
%%%%%%%%%%%%%%%%%%%%%%%%%%%%%%%%%%%%%%%%%%%%%%%%%%%%%%%%%%%%%%%%%%%%%%%%%%%%%%%%%%%%%%%%%%%%%%%%%%%%%%%%%%%%%%%%%%%%%%%%%%%%%%%%
\theorem\label{thmconnectednessandnonclingingpartitions}
$\Xt=\opair{\X}{\topology{}}$
is taken as a topological-space.
$\Xt$
is a connected topological-space if and only if
$\Xt$
has no non-clinging partitions consisting of two elements. That is,
\begin{equation}
\bigg(\X\in\connecteds{\Xt}\bigg)\thenn
\bigg(\defset{\apartition}{\Ncpart{\Xt}}
{\CarD{\apartition}=2}=\empty\bigg).
\end{equation}
\prooff
\begin{itemize}
\item[${\textbf{\textsf{p1}}}$]
It is assumed that $\Xt$ is disconnected:
\begin{equation}
\X\notin\connecteds{\Xt}.
\end{equation}
Then according to
\refdef{thmconnectednessandopenpartitions},
there exists at least one open partition $\seta{\binary{\U}{\(\compl{\X}{\U}\)}}$ of $\Xt$ consisting of two elements.
It is clear that $\U$ is a closed set of $\Xt$
(because $\compl{\X}{\U}$ is an open set of $\Xt$).
Thus according to \refthm{thmclosureofclosedset},
\begin{equation}
\U=\func{\Cl{\Xt}}{\U},
\end{equation}
and hence,
\begin{align}
\(\compl{\X}{\U}\)\cap\func{\Cl{\Xt}}{\U}&=
\(\compl{\X}{\U}\)\cap\U\cr
&=\empty.
\end{align}
Thus according to \refthm{defnonclingingpartition},
$\seta{\binary{\U}{\(\compl{\X}{\U}\)}}$
is a non-clinging partition of $\Xt$ with two elements. Therefore,
\begin{equation*}
\defset{\apartition}{\Ncpart{\Xt}}
{\CarD{\apartition}=2}\neq\empty.
\end{equation*}
\endp
\end{itemize}
\begin{itemize}
\item[${\textbf{\textsf{p2}}}$]
It is assumed that,
\begin{equation}
\defset{\apartition}{\Ncpart{\Xt}}
{\CarD{\apartition}=2}\neq\empty.
\end{equation}
Then there exists at least one non-clinging partition $\seta{\binary{\U}{\(\compl{\X}{\U}\)}}$ of $\Xt$.
According to \refdef{defnonclingingpartition},
\begin{align}
\U\cap\func{\Cl{\Xt}}{\compl{\X}{\U}}&=\empty,\\
\(\compl{\X}{\U}\)\cap\func{\Cl{\Xt}}{\U}&=\empty.
\end{align}
Therefore,
\begin{align}
\func{\Cl{\Xt}}{\U}&=\U,\\
\func{\Cl{\Xt}}{\compl{\X}{\U}}&=\(\compl{\X}{\U}\),
\end{align}
which according to \refthm{thmclosureofclosedset}
means,
\begin{equation}
\opair{\U}{\compl{\X}{\U}}\in\Fclosed{\X}{\topology{}}\times\Fclosed{\X}{\topology{}},
\end{equation}
and hence according to
\refdef{defclosedpartition},
$\seta{\binary{\Y}{\compl{\X}{\U}}}$
is a closed partition of $\Xt$ with two elements.
\begin{equation}
\seta{\binary{\U}{\compl{\X}{\U}}}\in
\defset{\apartition}{\Clpart{\Xt}}
{\CarD{\apartition}=2}.
\end{equation}
Thus,
\begin{equation}
\defset{\apartition}{\Clpart{\Xt}}
{\CarD{\apartition}=2}\neq\empty,
\end{equation}
and hence according to \refthm{thmconnectednessandclosedpartitions},
$\Xt$
is disconnected.
\begin{equation*}
\X\notin\connecteds{\Xt}.
\end{equation*}
\end{itemize}
\endthm
%%%%%%%%%%%%%%%%%%%%%%%%%%%%%%%%%%%%%%%%%%%%%%%%%%%%%%%%%%%%%%%%%%%%%%%%%%%%%%%%%%%%%%%%%%%%%%%%%%%%%%%%%%%%%%%%%%%%%%%%%%%%%%%%
\theorem\label{thmunionoftwoconnectedsets1}
$\Xt=\opair{\X}{\topology{}}$
is taken as a topological-space, and $\asubset$ and $\bsubset$
as subsets of $\X$. If each $\asubset$ and $\bsubset$
is a connected set of $\Xt$, and $\asubset$
intersects the closure of $\bsubset$ in $\Xt$, then the union of $\asubset$ and $\bsubset$
is a connected set of $\Xt$.
\begin{align}
&\[\bigg(\opair{\asubset}{\bsubset}\in\connecteds{\Xt}\times\connecteds{\Xt}\bigg),~
\asubset\cap\func{\Cl{\Xt}}{\bsubset}\neq\empty\]\cr
\then
&\bigg(\(\asubset\cup\bsubset\)\in\connecteds{\Xt}\bigg).
\end{align}
\prooff
$\asubset$
and
$\bsubset$
are taken as such elements of $\connecteds{\Xt}$ that,
\begin{equation}\label{thmunionoftwoconnectedsets1peq1}
\asubset\cap\func{\Cl{\Xt}}{\bsubset}\neq\empty.
\end{equation}
Therefore, $\asubset$ and $\bsubset$ are non-empty. Considering that $\asubset$ and $\bsubset$ are connected in $\Xt$,
\refdef{defconnectedness} implies that,
\begin{align}
\stopology{\topology{}}{\asubset}\cap
\Fclosed{\asubset}{\stopology{\topology{}}{\asubset}}&=
\seta{\binary{\empty}{\asubset}},
\label{thmunionoftwoconnectedsets1peq2}\\
\stopology{\topology{}}{\bsubset}\cap
\Fclosed{\bsubset}{\stopology{\topology{}}{\bsubset}}&=
\seta{\binary{\empty}{\bsubset}}.\label{thmunionoftwoconnectedsets1peq3}
\end{align}
In addition, according to \refthm{thmsubspaceclosure},
\begin{align}\label{thmunionoftwoconnectedsets1peq4}
\func{\Cl{\opair{\asubset\cup\bsubset}{\stopology{\topology{}}{\asubset\cup\bsubset}}}}{\bsubset}
&=\(\asubset\cup\bsubset\)\cap\func{\Cl{\Xt}}{\bsubset}\cr
&=\(\asubset\cap\func{\Cl{\Xt}}{\bsubset}\)\cup
\(\bsubset\cap\func{\Cl{\Xt}}{\bsubset}\)\cr
&=\(\asubset\cap\func{\Cl{\Xt}}{\bsubset}\)\cup\bsubset.
\end{align}
\begin{itemize}
\item[${\textbf{\textsf{p1}}}$]
$\U$
is taken as an element of
$\compl{\[\stopology{\topology{}}{\asubset\cup\bsubset}
\cap\Fclosed{\asubset\cup\bsubset}
{\stopology{\topology{}}{\asubset\cup\bsubset}}\]}{\seta{\empty}}$
(a non-empty open-and-closed set of $\opair{\asubset\cup\bsubset}
{\stopology{\topology{}}{\asubset\cup\bsubset}}$).
Then considering that each $\asubset$ and $\bsubset$
 is a subset of $\asubset\cup\bsubset$, according to \refdef{defsubspacetopology1} and \refthm{thmsubspacetopologytransitivity},
\begin{align}
\(\asubset\cap\U\)&\in\bigg(\stopology{\topology{}}{\asubset}\cap
\Fclosed{\asubset}{\stopology{\topology{}}{\asubset}}\bigg),
\label{thmunionoftwoconnectedsets1p1eq1}\\
\(\bsubset\cap\U\)&\in\bigg(\stopology{\topology{}}{\bsubset}\cap
\Fclosed{\bsubset}{\stopology{\topology{}}{\bsubset}}\bigg).
\label{thmunionoftwoconnectedsets1p1eq2}
\end{align}
Therefore according to \Ref{thmunionoftwoconnectedsets1peq2} and \Ref{thmunionoftwoconnectedsets1peq3},
\begin{align}
\(\asubset\cap\U\)&\in\seta{\binary{\empty}{\asubset}},
\label{thmunionoftwoconnectedsets1p1eq3}\\
\(\bsubset\cap\U\)&\in\seta{\binary{\empty}{\bsubset}}.
\label{thmunionoftwoconnectedsets1p1eq4}
\end{align}
Considering that $\U$ is an open subset of $\asubset\cup\bsubset$,
\begin{align}\label{thmunionoftwoconnectedsets1p1eq5}
\(\asubset\cap\U\)\cup\(\bsubset\cap\U\)&=
\(\asubset\cup\bsubset\)\cap\U\cr
&=\U,
\end{align}
and hence at least one of $\asubset\cap\U$ and $\bsubset\cap\U$ is non-empty.
\begin{equation}\label{thmunionoftwoconnectedsets1p1eq6}
\OR{\(\asubset\cap\U\neq\empty\)}{\(\bsubset\cap\U\neq\empty\)}.
\end{equation}
Thus according to
\Ref{thmunionoftwoconnectedsets1p1eq3} and \Ref{thmunionoftwoconnectedsets1p1eq4},
\begin{equation}\label{thmunionoftwoconnectedsets1p1eq7}
\OR{\(\asubset\cap\U=\asubset\)}{\(\bsubset\cap\U=\bsubset\)},
\end{equation}
which means,
\begin{equation}\label{thmunionoftwoconnectedsets1p1eq8}
\OR{\(\asubset\subseteq\U\)}{\(\bsubset\subseteq\U\)}.
\end{equation}
\begin{itemize}
\item[${\textbf{\textsf{p1-1}}}$]
It is assumed that,
\begin{equation}\label{thmunionoftwoconnectedsets1p1-1eq1}
\asubset\subseteq\U.
\end{equation}
\Ref{thmunionoftwoconnectedsets1peq4}

\Ref{thmunionoftwoconnectedsets1p1-1eq1}
imply that,
\begin{equation}
\U\cap\func{\Cl{\opair{\asubset\cup\bsubset}{\stopology{\topology{}}{\asubset\cup\bsubset}}}}{\bsubset}=
\(\asubset\cap\func{\Cl{\Xt}}{\bsubset}\)\cup\(\U\cap\bsubset\),
\end{equation}
and hence according to \Ref{thmunionoftwoconnectedsets1peq1},
\begin{equation}
\U\cap\func{\Cl{\opair{\asubset\cup\bsubset}{\stopology{\topology{}}{\asubset\cup\bsubset}}}}{\bsubset}\neq\empty.
\end{equation}
Therefore considering that,
\begin{equation}
\U\in\stopology{\topology{}}{\asubset\cup\bsubset},
\end{equation}
and based on \refthm{thmopensetsintersectingclosure},
\begin{equation}
\bsubset\cap\U\neq\empty.
\end{equation}
Therefore according to \Ref{thmunionoftwoconnectedsets1p1eq4},
\begin{equation}
\bsubset\cap\U=\bsubset,
\end{equation}
which means,
\begin{equation}
\bsubset\subseteq\U,
\end{equation}
and hence according to \Ref{thmunionoftwoconnectedsets1p1-1eq1},
\begin{equation}
\U=\asubset\cup\bsubset.
\end{equation}
\endp
\end{itemize}
Therefore,
\begin{equation}\label{thmunionoftwoconnectedsets1p1eq9}
\(\asubset\subseteq\U\)\then\(\U=\asubset\cup\bsubset\).
\end{equation}
\begin{itemize}
\item[${\textbf{\textsf{p1-2}}}$]
It is assumed that,
\begin{equation}\label{thmunionoftwoconnectedsets1p1-2eq1}
\bsubset\subseteq\U.
\end{equation}
Considering that $\U$ is a closed set of $\opair{\asubset\cup\bsubset}{\stopology{\topology{}}{\asubset\cup\bsubset}}$
and according to \refcor{corclosureofset0},
\begin{equation}
\U\supseteq\func{\Cl{\opair{\asubset\cup\bsubset}{\stopology{\topology{}}{\asubset\cup\bsubset}}}}{\bsubset}.
\end{equation}
Therefore according to \Ref{thmunionoftwoconnectedsets1peq4},
\begin{equation}
\U\supseteq\[\(\asubset\cap\func{\Cl{\Xt}}{\bsubset}\)\cup\bsubset\],
\end{equation}
and hence according to \Ref{thmunionoftwoconnectedsets1peq1},
\begin{equation}
\asubset\cap\U\neq\empty.
\end{equation}
Therefore according to \Ref{thmunionoftwoconnectedsets1p1eq3},
\begin{equation}
\asubset\cap\U=\asubset,
\end{equation}
which means,
\begin{equation}
\asubset\subseteq\U,
\end{equation}
and hence according to \Ref{thmunionoftwoconnectedsets1p1-2eq1},
\begin{equation}
\U=\asubset\cup\bsubset.
\end{equation}
\endp
\end{itemize}
Therefore,
\begin{equation}\label{thmunionoftwoconnectedsets1p1eq10}
\(\bsubset\subseteq\U\)\then
\(\U=\asubset\cup\bsubset\).
\end{equation}
\Ref{thmunionoftwoconnectedsets1p1eq8},
\Ref{thmunionoftwoconnectedsets1p1eq9},
and
\Ref{thmunionoftwoconnectedsets1p1eq10}
imply that,
\begin{equation}
\U=\asubset\cup\bsubset.
\end{equation}
\endp
\end{itemize}
Therefore,
\begin{equation}
\[\stopology{\topology{}}{\asubset\cup\bsubset}
\cap\Fclosed{\asubset\cup\bsubset}
{\stopology{\topology{}}{\asubset\cup\bsubset}}\]=
\seta{\binary{\empty}{\asubset\cup\bsubset}},
\end{equation}
which according to \refdef{defconnectedness}, means
$\asubset\cup\bsubset$
is a connected set of $\Xt$.
\begin{equation}
\(\asubset\cup\bsubset\)\in\connecteds{\Xt}.
\end{equation}
\endthm
%%%%%%%%%%%%%%%%%%%%%%%%%%%%%%%%%%%%%%%%%%%%%%%%%%%%%%%%%%%%%%%%%%%%%%%%%%%%%%%%%%%%%%%%%%%%%%%%%%%%%%%%%%%%%%%%%%%%%%%%%%%%%%%%
\theorem\label{thmsubspaceconnectedness}
$\Xt=\opair{\X}{\topology{}}$
is taken as a topological-space, and $\Y$ as a subset of $\X$. For every $\asubset$ in
$\CSs{\Y}$, $\asubset$
is a connected set of $\Xt$ if and only if
$\asubset$ is a connected set of $\opair{\Y}{\stopology{\topology{}}{\Y}}$. That is,
\begin{equation}
\Foreach{\asubset}{\CSs{\Y}}
\bigg[
\asubset\in\connecteds{\Xt}\thenn
\asubset\in\connecteds{\opair{\Y}{\stopology{\topology{}}{\Y}}}
\bigg].
\end{equation}
\proof
$\asubset$
is taken as an arbitrary subset of $\Y$. According to \refdef{defconnectedness},
\begin{align}
&\bigg[\asubset\in\connecteds{\Xt}\bigg]
\thenn
\bigg[\stopology{\topology{}}{\asubset}\cap
\Fclosed{\asubset}{\stopology{\topology{}}{\asubset}}=\seta{\binary{\empty}{\asubset}}\bigg],\\
&\bigg[\asubset\in\connecteds{\opair{\Y}{\stopology{\topology{}}{\Y}}}\bigg]
\thenn
\bigg[\stopology{\stopology{\topology{}}{\Y}}{\asubset}\cap
\Fclosed{\asubset}{\stopology{\stopology{\topology{}}{\Y}}{\asubset}}=
\seta{\binary{\empty}{\asubset}}\bigg].
\end{align}
In addition, according to \refthm{thmsubspacetopologytransitivity},
\begin{equation}
\stopology{\stopology{\topology{}}{\Y}}{\asubset}=
\stopology{\topology{}}{\asubset}.
\end{equation}
Therefore,
\begin{equation}
\asubset\in\connecteds{\Xt}\thenn
\asubset\in\connecteds{\opair{\Y}{\stopology{\topology{}}{\Y}}}.
\end{equation}
\endthm
%%%%%%%%%%%%%%%%%%%%%%%%%%%%%%%%%%%%%%%%%%%%%%%%%%%%%%%%%%%%%%%%%%%%%%%%%%%%%%%%%%%%%%%%%%%%%%%%%%%%%%%%%%%%%%%%%%%%%%%%%%%%%%%%
\theorem\label{thmconnectednessandfrontier}
$\Xt=\opair{\X}{\topology{}}$
is taken as a topological-space.
$\Xt$
is a connected topological-space if and only if $\empty$ and $\X$
are the only subsets of $\X$ with empty boundary in $\Xt$. That is,
\begin{equation}\label{thmconnectednessandfrontiereq1}
\bigg(\X\in\connecteds{\Xt}\bigg)\thenn
\bigg(\defset{\asubset}{\CSs{\X}}{\[\func{\Fr{\Xt}}{\asubset}=\empty\]}=
\seta{\binary{\empty}{\X}}\bigg).
\end{equation}
\proof
According to \refthm{thmfrontierofclopenset},
\begin{equation}
\defset{\asubset}{\CSs{\X}}{\[\func{\Fr{\Xt}}{\asubset}=\empty\]}=
\[\topology{}\cap\Fclosed{\X}{\topology{}}\].
\end{equation}
In addition, according to \refdef{defconnectedness},
\begin{equation}
\bigg(\X\in\connecteds{\Xt}\bigg)\thenn
\bigg(\[\topology{}\cap\Fclosed{\X}{\topology{}}\]=\seta{\binary{\empty}{\X}}\bigg).
\end{equation}
Therefore
\Ref{thmconnectednessandfrontiereq1} is clearly obtained.
\endthm
%%%%%%%%%%%%%%%%%%%%%%%%%%%%%%%%%%%%%%%%%%%%%%%%%%%%%%%%%%%%%%%%%%%%%%%%%%%%%%%%%%%%%%%%%%%%%%%%%%%%%%%%%%%%%%%%%%%%%%%%%%%%%%%%
\theorem\label{thmintersectionofaconnectedsetwithfrontierofanotherset}
$\Xt=\opair{\X}{\topology{}}$
is taken as a topological-space.
For every connected set $\aconnectedset$ of $\Xt$,
and for every subset $\asubset$ of $\X$,
If $\aconnectedset$ intersects both $\asubset$ and $\(\compl{\X}{\asubset}\)$, then
$\aconnectedset$ intersects the boundary of $\asubset$ in $\Xt$ (that is, $\func{\Fr{\Xt}}{\asubset}$).
That is,
\begin{align}
&\Foreach{\opair{\aconnectedset}{\asubset}}{\bigg(\connecteds{\Xt}\times{\CSs{\X}}\bigg)}\cr
&\bigg[
\bigg(\AND{\aconnectedset\cap\asubset\neq\empty}{\aconnectedset\cap\(\compl{\X}{\asubset}\)\neq\empty}\bigg)
\then\bigg(\aconnectedset\cap\func{\Fr{\Xt}}{\asubset}\neq\empty\bigg)
\bigg].
\end{align}
\prooff
$\aconnectedset$
is taken as an arbitrary element of $\connecteds{\Xt}$, and $\asubset$ as an arbitrary element of $\CSs{\X}$,
and it is assumed that,
\begin{align}
\aconnectedset\cap\asubset&\neq\empty,
\label{thmintersectionofaconnectedsetwithfrontierofanothersetpeq1}\\
\aconnectedset\cap\(\compl{\X}{\asubset}\)&\neq\empty.
\label{thmintersectionofaconnectedsetwithfrontierofanothersetpeq2}
\end{align}
Then considering that every subset of $\X$ is a subset of its closure, it is clear that,
\begin{align}
\aconnectedset\cap\func{\Cl{\Xt}}{\asubset}&\neq\empty,
\label{thmintersectionofaconnectedsetwithfrontierofanothersetpeq3}\\
\aconnectedset\cap\func{\Cl{\Xt}}{\(\compl{\X}{\asubset}\)}&\neq\empty.
\label{thmintersectionofaconnectedsetwithfrontierofanothersetpeq4}
\end{align}
According to \refcor{corclosureofset0},
Each $\func{\Cl{\Xt}}{\asubset}$ and $\func{\Cl{\Xt}}{\compl{\X}{\asubset}}$
is a closed set of $\Xt$.That is,
\begin{align}
\func{\Cl{\Xt}}{\asubset}&\in\Fclosed{\X}{\topology{}},
\label{thmintersectionofaconnectedsetwithfrontierofanothersetpeq5}\\
\func{\Cl{\Xt}}{\compl{\X}{\asubset}}&\in\Fclosed{\X}{\topology{}}.
\label{thmintersectionofaconnectedsetwithfrontierofanothersetpeq6}
\end{align}
Based on \refthm{thmsubspaceclosedsets},
\Ref{thmintersectionofaconnectedsetwithfrontierofanothersetpeq5}
and
\Ref{thmintersectionofaconnectedsetwithfrontierofanothersetpeq6} imply that,
\begin{align}
\aconnectedset\cap\func{\Cl{\Xt}}{\asubset}
&\in\Fclosed{\aconnectedset}{\stopology{\topology{}}{\aconnectedset}},
\label{thmintersectionofaconnectedsetwithfrontierofanothersetpeq7}\\
\aconnectedset\cap\func{\Cl{\Xt}}{\compl{\X}{\asubset}}
&\in\Fclosed{\aconnectedset}{\stopology{\topology{}}{\aconnectedset}}.
\label{thmintersectionofaconnectedsetwithfrontierofanothersetpeq8}
\end{align}
In addition, it can be easily seen that,
\begin{align}\label{thmintersectionofaconnectedsetwithfrontierofanothersetpeq9}
\bigg(\aconnectedset\cap\func{\Cl{\Xt}}{\asubset}\bigg)\cup
\bigg(\aconnectedset\cap\func{\Cl{\Xt}}{\compl{\X}{\asubset}}\bigg)&=
\aconnectedset\cap\bigg(\func{\Cl{\Xt}}{\asubset}\cup\func{\Cl{\Xt}}{\compl{\X}{\asubset}}\bigg)\cr
&=\aconnectedset\cap\X\cr
&=\aconnectedset.
\end{align}
According to \refdef{defclosedpartition},
\Ref{thmintersectionofaconnectedsetwithfrontierofanothersetpeq3},
\Ref{thmintersectionofaconnectedsetwithfrontierofanothersetpeq4},
\Ref{thmintersectionofaconnectedsetwithfrontierofanothersetpeq7},
\Ref{thmintersectionofaconnectedsetwithfrontierofanothersetpeq8},
and
\Ref{thmintersectionofaconnectedsetwithfrontierofanothersetpeq9}
imply that if
$\bigg(\aconnectedset\cap\func{\Cl{\Xt}}{\asubset}\bigg)$
does not intersect
$\bigg(\aconnectedset\cap\func{\Cl{\Xt}}{\compl{\X}{\asubset}}\bigg)$, then
$\seta{\binary{\bigg(\aconnectedset\cap\func{\Cl{\Xt}}{\asubset}\bigg)}
{\bigg(\aconnectedset\cap\func{\Cl{\Xt}}{\compl{\X}{\asubset}}\bigg)}}$
is a closed partition of $\Xt$ with two elements, and hence
$\opair{\aconnectedset}{\stopology{\topology{}}{\aconnectedset}}$
the set of all closed partitions of $\opair{\aconnectedset}{\stopology{\topology{}}{\aconnectedset}}$ with two elements in not empty. That is,
\begin{gather}
\[\bigg(\aconnectedset\cap\func{\Cl{\Xt}}{\asubset}\bigg)\cap
\bigg(\aconnectedset\cap\func{\Cl{\Xt}}{\compl{\X}{\asubset}}\bigg)=\empty\]\cr
\vthen\cr
\bigg(\defset{\apartition}{\Clpart{\opair{\aconnectedset}{\stopology{\topology{}}{\aconnectedset}}}}
{\CarD{\apartition}=2}\neq\empty\bigg),
\end{gather}
and hence according to \refthm{thmconnectednessandclosedpartitions},
\begin{gather}
\[\bigg(\aconnectedset\cap\func{\Cl{\Xt}}{\asubset}\bigg)\cap
\bigg(\aconnectedset\cap\func{\Cl{\Xt}}{\compl{\X}{\asubset}}\bigg)=\empty\]\cr
\vthen\cr
\aconnectedset\notin\connecteds{\Xt}.
\end{gather}
Therefore considering that $\aconnectedset\in\connecteds{\Xt}$,
\begin{equation}
\bigg(\aconnectedset\cap\func{\Cl{\Xt}}{\asubset}\bigg)\cap
\bigg(\aconnectedset\cap\func{\Cl{\Xt}}{\compl{\X}{\asubset}}\bigg)\neq\empty,
\end{equation}
or equivalently,
\begin{equation}
\aconnectedset\cap\bigg(
\func{\Cl{\Xt}}{\asubset}\cap\func{\Cl{\Xt}}{\compl{\X}{\asubset}}\bigg)\neq\empty,
\end{equation}
which according to \refthm{thmfrontier1} means,
\begin{equation}
\aconnectedset\cap\func{\Fr{\Xt}}{\asubset}\neq\empty.
\end{equation}
\endthm
%%%%%%%%%%%%%%%%%%%%%%%%%%%%%%%%%%%%%%%%%%%%%%%%%%%%%%%%%%%%%%%%%%%%%%%%%%%%%%%%%%%%%%%%%%%%%%%%%%%%%%%%%%%%%%%%%%%%%%%%%%%%%%%%
\theorem\label{thmconnectednessincoarsertopology}
$\X$
is taken as a set. Each $\topology{}$ and $\p{\topology{}}$ is taken as a topology in $\X$ such that
$\topology{}$ is coarser than $\p{\topology{}}$.
($\topology{}\subseteq\p{\topology{}}$).
If the topological space $\opair{\X}{\p{\topology{}}}$
is connected then the topological space $\opair{\X}{\topology{}}$ is also connected.
\proof
Suppose that $\opair{\X}{\p{\topology{}}}$ is connected, that is,
\begin{equation}\label{thmconnectednessincoarsertopologypeq1}
\p{\topology{}}\cap\Fclosed{\X}{\p{\topology{}}}=\seta{\binary{\empty}{\X}}.
\end{equation}
Considering that
$\topology{}\subseteq\p{\topology{}}$,
\refdef{deffinersoarsertopology},
and
\refthm{thmclosedsetsofcoarsertopology} yield,
\begin{equation}\label{thmconnectednessincoarsertopologypeq2}
\topology{}\cap\Fclosed{\X}{\topology{}}\subseteq
\p{\topology{}}\cap\Fclosed{\X}{\p{\topology{}}},
\end{equation}
and hence according to \Ref{thmconnectednessincoarsertopologypeq1},
and considering that $\topology{}$ is a topology in $\X$,
\begin{equation}
\topology{}\cap\Fclosed{\X}{\topology{}}=\seta{\binary{\empty}{\X}}.
\end{equation}
According to \refdef{defconnectedness},
this means $\opair{\X}{\topology{}}$ is connected.
\endthm
%%%%%%%%%%%%%%%%%%%%%%%%%%%%%%%%%%%%%%%%%%%%%%%%%%%%%%%%%%%%%%%%%%%%%%%%%%%%%%%%%%%%%%%%%%%%%%%%%%%%%%%%%%%%%%%%%%%%%%%%%%%%%%%%
\theorem
Each $\Xt_1$ and $\Xt_2$ is taken as a topological space.
If $\Xt_1$ and $\Xt_2$ are connected, then so is $\topprod{\Xt_1}{\Xt_2}$.
\proof
It is left as an exercise.
\endthm
%%%%%%%%%%%%%%%%%%%%%%%%%%%%%%%%%%%%%%%%%%%%%%%%%%%%%%%%%%%%%%%%%%%%%%%%%%%%%%%%%%%%%%%%%%%%%%%%%%%%%%%%%%%%%%%%%%%%%%%%%%%%%%%%
\theorem
$\Xt=\opair{\X}{\topology{\X}}$ is taken as a topological space, and $\eqrel{}$ as an equivalence relation on $\X$.
If $\Xt$ is connected, then so is $\topq{\Xt}{\eqrel{}}$.
\proof
It is left as an exercise.
\endthm
%%%%%%%%%%%%%%%%%%%%%%%%%%%%%%%%%%%%%%%%%%%%%%%%%%%%%%%%%%%%%%%%%%%%%%%%%%%%%%%%%%%%%%%%%%%%%%%%%%%%%%%%%%%%%%%%%%%%%%%%%%%%%%%%
\subsection{
Connectedness As a Topological Property
}
\theorem\label{thmconnectednessandcontinuousfunctions}
Each $\Xt=\opair{\X}{\topology{\X}}$ and $\Yt=\opair{\Y}{\topology{\Y}}$
is taken as a topological space. If $\Xt$ is connected, then for every $\cf$ in $\CF{\Xt}{\Yt}$
(every continuous map from $\cf$ $\Xt$ to $\Yt$),
$\opair{\func{\image{\cf}}{\X}}{\stopology{\topology{\Y}}{\func{\image{\cf}}{\X}}}$
is also connected. That is,
\begin{align}
\bigg(\X\in\connecteds{\Xt}\bigg)\then
\[\Foreach{\cf}{\CF{\Xt}{\Yt}}
\bigg(\func{\image{\cf}}{\X}\in
\connecteds{
\Yt
}
\bigg)\].
\end{align}
\prooff
It is assumed that $\Xt$ is connected. Then according to \refdef{defconnectedness},
\begin{equation}\label{thmconnectednessandcontinuousfunctionspeq1}
\topology{\X}\cap\Fclosed{\X}{\topology{\X}}=\seta{\binary{\empty}{\X}}.
\end{equation}
$\cf$
is taken as an arbitrary element of $\CF{\Xt}{\Yt}$. Then according to \refcor{correstrictionofcontinuousfunction},
$\func{\rescd{\cf}}{\func{\image{\cf}}{\X}}$
is a continuous map from $\Xt$ to the topological space
$\opair{\func{\image{\cf}}{\X}}{\stopology{\topology{\Y}}{\func{\image{\cf}}{\X}}}$:
\begin{equation}\label{thmconnectednessandcontinuousfunctionspeq2}
\func{\rescd{\cf}}{\func{\image{\cf}}{\X}}\in
\CF{\Xt}{\opair{\func{\image{\cf}}{\X}}{\stopology{\topology{\Y}}{\func{\image{\cf}}{\X}}}}.
\end{equation}
Thus, according to \refdef{defcontinuousfunction} and \refthm{thmcontiniuityequiv1},
\begin{align}
&\Foreach{\U}{\stopology{\topology{\Y}}{\func{\image{\cf}}{\X}}}
\func{\pimage{\[\func{\rescd{\cf}}{\func{\image{\cf}}{\X}}\]}}{\U}\in\topology{\X},
\label{thmconnectednessandcontinuousfunctionspeq3}\\
&\Foreach{\U}{\Fclosed{\Y}{\stopology{\topology{\Y}}{\func{\image{\cf}}{\X}}}}
\func{\pimage{\[\func{\rescd{\cf}}{\func{\image{\cf}}{\X}}\]}}{\U}\in\Fclosed{\X}{\topology{\X}}.
\label{thmconnectednessandcontinuousfunctionspeq4}
\end{align}
In addition, considering that $\X$ is the domain of $\cf$, it is clear that
$\func{\rescd{\cf}}{\func{\image{\cf}}{\X}}$
is a surjective map from $\X$ to $\func{\image{\cf}}{\X}$:
\begin{equation}\label{thmconnectednessandcontinuousfunctionspeq5}
\func{\rescd{\cf}}{\func{\image{\cf}}{\X}}\in\surFunc{\X}{\func{\image{\cf}}{\X}}.
\end{equation}
\begin{itemize}
\item[${\textbf{\textsf{p1}}}$]
$\U$
is taken as an arbitrary element of $\[\stopology{\topology{\Y}}{\func{\image{\cf}}{\X}}
\cap\Fclosed{\Y}{\stopology{\topology{\Y}}{\func{\image{\cf}}{\X}}}\]$
(a clopen set of the topological space $\opair{\func{\image{\cf}}{\X}}{\stopology{\topology{\Y}}{\func{\image{\cf}}{\X}}}$)
Then acccording to \Ref{thmconnectednessandcontinuousfunctionspeq3} and
\Ref{thmconnectednessandcontinuousfunctionspeq4},
it is clear that,
\begin{equation}\label{thmconnectednessandcontinuousfunctionsp1eq1}
\func{\pimage{\[\func{\rescd{\cf}}{\func{\image{\cf}}{\X}}\]}}{\U}
\in\[\topology{\X}\cap\Fclosed{\X}{\topology{\X}}\],
\end{equation}
and hus according to \Ref{thmconnectednessandcontinuousfunctionspeq1},
\begin{equation}\label{thmconnectednessandcontinuousfunctionsp1eq2}
\func{\pimage{\[\func{\rescd{\cf}}{\func{\image{\cf}}{\X}}\]}}{\U}
\in\seta{\binary{\empty}{\X}}.
\end{equation}
Moreover, according to \Ref{thmconnectednessandcontinuousfunctionspeq5}
(that is the surjectivity of $\func{\rescd{\cf}}{\func{\image{\cf}}{\X}}$),
it is obvious that,
\begin{align}
&\bigg(\func{\pimage{\[\func{\rescd{\cf}}{\func{\image{\cf}}{\X}}\]}}{\U}=\empty\bigg)\then
\bigg(\U=\empty\bigg),\label{thmconnectednessandcontinuousfunctionsp1eq3}\\
&\bigg(\func{\pimage{\[\func{\rescd{\cf}}{\func{\image{\cf}}{\X}}\]}}{\U}=\X\bigg)\then
\bigg(\U=\func{\image{\cf}}{\X}\bigg).\label{thmconnectednessandcontinuousfunctionsp1eq4}
\end{align}
\Ref{thmconnectednessandcontinuousfunctionsp1eq2},
\Ref{thmconnectednessandcontinuousfunctionsp1eq3},
and
\Ref{thmconnectednessandcontinuousfunctionsp1eq4}
yield,
\begin{equation}
\U\in\seta{\binary{\empty}{\func{\image{\cf}}{\X}}}.
\end{equation}
\endp
\end{itemize}
Therefore,
\begin{equation}
\[\stopology{\topology{\Y}}{\func{\image{\cf}}{\X}}
\cap\Fclosed{\Y}{\stopology{\topology{\Y}}{\func{\image{\cf}}{\X}}}\]\subseteq\seta{\binary{\empty}{\func{\image{\cf}}{\X}}}.
\end{equation}
Moreover, it is clear that,
\begin{equation}
\seta{\binary{\empty}{\func{\image{\cf}}{\X}}}\subseteq\[\stopology{\topology{\Y}}{\func{\image{\cf}}{\X}}
\cap\Fclosed{\Y}{\stopology{\topology{\Y}}{\func{\image{\cf}}{\X}}}\].
\end{equation}
Thus,
\begin{equation}
\[\stopology{\topology{\Y}}{\func{\image{\cf}}{\X}}
\cap\Fclosed{\Y}{\stopology{\topology{\Y}}{\func{\image{\cf}}{\X}}}\]=\seta{\binary{\empty}{\func{\image{\cf}}{\X}}}.
\end{equation}
According to \refdef{defconnectedness},
this means,
$\func{\image{\cf}}{\X}$
is a connected set of the topological space $\Yt$:
\begin{equation}
\func{\image{\cf}}{\X}\in\connecteds{\Yt},
\end{equation}
or in other words,
$\opair{\func{\image{\cf}}{\X}}{\stopology{\topology{\Y}}{\func{\image{\cf}}{\X}}}$
is a connected topological space:
\begin{equation}
\func{\image{\cf}}{\X}\in
\connecteds{\opair{\func{\image{\cf}}{\X}}{\stopology{\topology{\Y}}{\func{\image{\cf}}{\X}}}}.
\end{equation}
\endthm
%%%%%%%%%%%%%%%%%%%%%%%%%%%%%%%%%%%%%%%%%%%%%%%%%%%%%%%%%%%%%%%%%%%%%%%%%%%%%%%%%%%%%%%%%%%%%%%%%%%%%%%%%%%%%%%%%%%%%%%%%%%%%%%%
\theorem\label{thmconnectedsetsandcontinuousfunctions}
Each
$\Xt=\opair{\X}{\topology{\X}}$
and
$\Yt=\opair{\Y}{\topology{\Y}}$
is taken as a topological space. For every continuous map $\cf$ from $\Xt$ to $\Yt$,
the image under $\cf$ of every connected set of $\Xt$ is a connected set of $\Yt$:
\begin{equation}
\Foreach{\cf}{\CF{\Xt}{\Yt}}
\[\Foreach{\aconnectedset}{\connecteds{\Xt}}
\bigg(\func{\image{\cf}}{\aconnectedset}\in\connecteds{\Yt}\bigg)\].
\end{equation}
\proof
Let $\cf$ be an arbitrary element of $\CF{\Xt}{\Yt}$, and let
$\aconnectedset$ be an arbitrary element of $\connecteds{\Xt}$. According to \refthm{thmcontiniuityequiv2},
$\func{\resd{\cf}}{\aconnectedset}$
is a continuous map from the topological space
$\opair{\aconnectedset}{\stopology{\topology{\X}}{\aconnectedset}}$ to $\Yt$:
\begin{equation}
\func{\resd{\cf}}{\aconnectedset}
\in\CF{\opair{\aconnectedset}{\stopology{\topology{}}{\aconnectedset}}}{\Yt}.
\end{equation}
Moreover, according to
\refdef{defconnectedness},
\begin{equation}
\aconnectedset\in\connecteds{\opair{\aconnectedset}{\stopology{\topology{\X}}{\aconnectedset}}}.
\end{equation}
Thus, according to \refthm{thmconnectednessandcontinuousfunctions},
\begin{align}
\func{\image{\[\func{\resd{\cf}}{\aconnectedset}\]}}{\aconnectedset}\in
\connecteds{
\Yt
}.
\end{align}
Therefore, considering that
\begin{equation}
\func{\image{\[\func{\resd{\cf}}{\aconnectedset}\]}}{\aconnectedset}=
\func{\image{\cf}}{\aconnectedset},
\end{equation}
it is clear that,
\begin{equation}
\func{\image{\cf}}{\aconnectedset}\in\connecteds{\Yt}.
\end{equation}
\endthm
%%%%%%%%%%%%%%%%%%%%%%%%%%%%%%%%%%%%%%%%%%%%%%%%%%%%%%%%%%%%%%%%%%%%%%%%%%%%%%%%%%%%%%%%%%%%%%%%%%%%%%%%%%%%%%%%%%%%%%%%%%%%%%%%
\theorem\label{thmconnectednessisatopologicalproperty}
Each $\Xt=\opair{\X}{\topology{\X}}$ and $\Yt=\opair{\Y}{\topology{\Y}}$
is taken as a topological space. If $\Xt$ and $\Yt$
are homeomorphic, then $\Xt$ is connected if and only if $\Yt$ is connected. That is,
\begin{equation}
\bigg(\homeomorphic{\Xt}{\Yt}\bigg)\then
\bigg(\X\in\connecteds{\Xt}\thenn\Y\in\connecteds{\Yt}\bigg).
\end{equation}
$\caution$
In other words, connectedness is a topological property.
\proof
Suppose that $\Xt$ is homeomorphic to $\Yt$, that is,
\begin{equation}\label{thmconnectednessisatopologicalpropertypeq1}
\homeomorphic{\Xt}{\Yt}.
\end{equation}
According to \refdef{defhomeomorphic},
this means there exists at least one homeomorphism from $\Xt$ to $\Yt$:
\begin{equation}\label{thmconnectednessisatopologicalpropertypeq2}
\HOF{\Xt}{\Yt}\neq\empty.
\end{equation}
$\hf$ is taken as an element of $\HOF{\Xt}{\Yt}$
($\HOF{\Xt}{\Yt}=\seta{\binary{\hf}{\ldots}}$).
According to \refdef{defhomeomorphism},
\begin{align}
\hf&\in\IF{\X}{\Y},
\label{thmconnectednessisatopologicalpropertypeq3}\\
\hf&\in\CF{\Xt}{\Yt},
\label{thmconnectednessisatopologicalpropertypeq4}\\
\finv{\hf}&\in\CF{\Yt}{\Xt}.
\label{thmconnectednessisatopologicalpropertypeq5}
\end{align}
According to \Ref{thmconnectednessisatopologicalpropertypeq3}
(bijectivity of $\hf$),
it is obvious that,
\begin{align}
\func{\image{\hf}}{\X}&=\Y,
\label{thmconnectednessisatopologicalpropertypeq6}\\
\func{\image{\[\finv{\hf}\]}}{\Y}&=\X.
\label{thmconnectednessisatopologicalpropertypeq7}
\end{align}
based on
\refthm{thmconnectednessandcontinuousfunctions},
\Ref{thmconnectednessisatopologicalpropertypeq4},
\Ref{thmconnectednessisatopologicalpropertypeq5},
\Ref{thmconnectednessisatopologicalpropertypeq6},
and
\Ref{thmconnectednessisatopologicalpropertypeq7}
yield,
\begin{equation}
\bigg(\X\in\connecteds{\Xt}\thenn\Y\in\connecteds{\Yt}\bigg).
\end{equation}
\endthm
%%%%%%%%%%%%%%%%%%%%%%%%%%%%%%%%%%%%%%%%%%%%%%%%%%%%%%%%%%%%%%%%%%%%%%%%%%%%%%%%%%%%%%%%%%%%%%%%%%%%%%%%%%%%%%%%%%%%%%%%%%%%%%%%
%%%%%%%%%%%%%%%%%%%%%%%%%%%%%%%%%%%%%%%%%%%%%%%%%%%%%%%%%%%%%%%%%%%%%%%%%%%%%%%%%%%%%%%%%%%%%%%%%%%%%%%%%%%%%%%%%%%%%%%%%%%%%%%%
%%%%%%%%%%%%%%%%%%%%%%%%%%%%%%%%%%%%%%%%%%%%%%%%%%%%%%%%%%%%%%%%%%%%%%%%%%%%%%%%%%%%%%%%%%%%%%%%%%%%%%%%%%%%%%%%%%%%%%%%%%%%%%%%
%%%%%%%%%%%%%%%%%%%%%%%%%%%%%%%%%%%%%%%%%%%%%%%%%%%%%%%%%%%%%%%%%%%%%%%%%%%%%%%%%%%%%%%%%%%%%%%%%%%%%%%%%%%%%%%%%%%%%%%%%%%%%%%%
%%%%%%%%%%%%%%%%%%%%%%%%%%%%%%%%%%%%%%%%%%%%%%%%%%%%%%%%%%%%%%%%%%%%%%%%%%%%%%%%%%%%%%%%%%%%%%%%%%%%%%%%%%%%%%%%%%%%%%%%%%%%%%%%
%%%%%%%%%%%%%%%%%%%%%%%%%%%%%%%%%%%%%%%%%%%%%%%%%%%%%%%%%%%%%%%%%%%%%%%%%%%%%%%%%%%%%%%%%%%%%%%%%%%%%%%%%%%%%%%%%%%%%%%%%%%%%%%%
%%%%%%%%%%%%%%%%%%%%%%%%%%%%%%%%%%%%%%%%%%%%%%%%%%%%%%%%%%%%%%%%%%%%%%%%%%%%%%%%%%%%%%%%%%%%%%%%%%%%%%%%%%%%%%%%%%%%%%%%%%%%%%%%
\section{
Maximally-Connected Sets of a Topological Space
}
\definition\label{defconnectedcomponent}
$\Xt=\opair{\X}{\topology{}}$
is taken as a topological space.
For every $\asubset$ in $\CSs{\X}$, $\asubset$
is referred to as a $\quotl$maximally-connected set of the topological space $\Xt$$\quotr$
iff the following properties are hold.
\begin{itemize}
\item[${\textbf{\textsf{MC1}}}$]
$\asubset$
is a connected set of $\Xt$.
\hfill
$\asubset\in\connecteds{\Xt}.$
\item[${\textbf{\textsf{MC2}}}$]
No connected set of $\Xt$ strictly includes $\asubset$.
\hfill
$\Foreach{\aconnectedset}{\[\compl{\connecteds{\Xt}}{\seta{\asubset}}\]}
\asubset\nsubseteq\aconnectedset.$
\item[${\textbf{\textsf{MC3}}}$]
$\asubset$ is non-empty.
\footnote{
One reason for assuming the non-emptiness condition for maximally-connected sets
is to make the collection of all maximally-connected sets of the empty topological space
a partition of it.}
\hfill
$\asubset\neq\empty.$
\end{itemize}
The collection of all maximally connected sets of $\Xt$ will be denoted by $\maxcon{\Xt}$:
\begin{align}
\maxcon{\Xt}:=\defset{\asubset}{\[\compl{\connecteds{\Xt}}{\seta{\empty}}\]}
{\[\Foreach{\aconnectedset}{\[\compl{\connecteds{\Xt}}{\seta{\asubset}}\]}
\asubset\nsubseteq\aconnectedset\]}.
\end{align}
$\caution$
Each element of $\maxcon{\Xt}$
ia alternatively called a $\quotl$connected component of the topological space $\Xt$$\quotr$.
\endef
%%%%%%%%%%%%%%%%%%%%%%%%%%%%%%%%%%%%%%%%%%%%%%%%%%%%%%%%%%%%%%%%%%%%%%%%%%%%%%%%%%%%%%%%%%%%%%%%%%%%%%%%%%%%%%%%%%%%%%%%%%%%%%%%
\theorem\label{thmconnectedcomponentsofemptyspace}
The topological space $\opair{\empty}{\seta{\empty}}$
does not have any maximally-connected set. That is,
\begin{equation}
\maxcon{\opair{\empty}{\seta{\empty}}}=\empty.
\end{equation}

\proof
According to \refdef{defconnectedcomponent} and \refthm{thmemptyspaceisconnected},
it is clear.
\endthm
%%%%%%%%%%%%%%%%%%%%%%%%%%%%%%%%%%%%%%%%%%%%%%%%%%%%%%%%%%%%%%%%%%%%%%%%%%%%%%%%%%%%%%%%%%%%%%%%%%%%%%%%%%%%%%%%%%%%%%%%%%%%%%%%
\theorem\label{thmconnectedspacehasoneconnectedcomponent}
$\Xt=\opair{\X}{\topology{}}$
is taken as a non-empty topological space.
($\X\neq\empty$).
$\Xt$ is connected if and only if $\X$ is the only maximally-connected set of $\Xt$. That is,
\begin{equation}
\bigg(\X\in\connecteds{\Xt}\bigg)\thenn\bigg(\maxcon{\Xt}=\seta{\X}\bigg).
\end{equation}
\proof
\begin{itemize}
\item[${\textbf{\textsf{p1}}}$]
It is assumed that,
$\Xt$ is connected:
\begin{equation}
\X\in\connecteds{\Xt}.
\end{equation}
According to \refdef{defconnectedcomponent}, it is clear that,
\begin{equation}
\X\in\maxcon{\Xt}.
\end{equation}
Moreover,
\begin{equation}
\Foreach{\asubset}{\compl{\CSs{\X}}{\seta{\X}}}
\[\Existsis{\X}{\compl{\connecteds{\Xt}}{\seta{\asubset}}}\asubset\subseteq\X\].
\end{equation}
According to \refdef{defconnectedcomponent}, this means
no other subset of $\X$ can be a maximally-connected set of $\Xt$:
\begin{equation}
\Foreach{\asubset}{\compl{\CSs{\X}}{\seta{\X}}}
\asubset\notin\maxcon{\Xt}.
\end{equation}
Therefore,
\begin{equation*}
\maxcon{\Xt}=\seta{\X}.
\end{equation*}
\endp
\end{itemize}
\begin{itemize}
\item[${\textbf{\textsf{p2}}}$]
It is assumed that,
$\X$ is the only maximally-connected set of $\Xt$:
\begin{equation}
\maxcon{\Xt}=\seta{\X}.
\end{equation}
Then according to \refdef{defconnectedcomponent}, it is clear that,
\begin{equation*}
\X\in\connecteds{\Xt}.
\end{equation*}
\endp
\end{itemize}
\endthm
%%%%%%%%%%%%%%%%%%%%%%%%%%%%%%%%%%%%%%%%%%%%%%%%%%%%%%%%%%%%%%%%%%%%%%%%%%%%%%%%%%%%%%%%%%%%%%%%%%%%%%%%%%%%%%%%%%%%%%%%%%%%%%%%
\definition\label{defconnectedclassofaset}
$\Xt=\opair{\X}{\topology{}}$
is taken as a topological space.
The mapping
$\mcp{\Xt}$ is defined as,
\begin{align}
&\mcp{\Xt}\indef\Func{\CSs{\X}}{\CSs{\X}}\\
&\Foreach{\asubset}{\CSs{\X}}\func{\mcp{\Xt}}{\asubset}\eqdef
\union{\defset{\csubset}{\connecteds{\Xt}}{\asubset\subseteq\csubset}}.
\end{align}
This means, for every $\asubset$ in $\CSs{\X}$,
$\func{\mcp{\Xt}}{\point}$
is the union of all connected sets of $\Xt$ containing $\asubset$.
\endef
%%%%%%%%%%%%%%%%%%%%%%%%%%%%%%%%%%%%%%%%%%%%%%%%%%%%%%%%%%%%%%%%%%%%%%%%%%%%%%%%%%%%%%%%%%%%%%%%%%%%%%%%%%%%%%%%%%%%%%%%%%%%%%%%
\theorem\label{thmconnectedclassofapoint}
$\Xt=\opair{\X}{\topology{}}$
is taken as a topological space.
For every point $\point$ in $\X$,
$\func{\mcp{\Xt}}{\seta{\point}}$
is a maximally-connected set of $\Xt$:
\begin{equation}
\Foreach{\point}{\X}
\func{\mcp{\Xt}}{\seta{\point}}\in\maxcon{\Xt}.
\end{equation}
\prooff
$\point$
is taken as an arbitrary element of $\X$. According to \refdef{defconnectedclassofaset},
$\func{\mcp{\Xt}}{\seta{\point}}$ is the union of the elements of a collection of intersectiong connected sets of $\Xt$:
\begin{equation}\label{thmconnectedclassofapointpeq1}
\func{\mcp{\Xt}}{\seta{\point}}=
\union{\defset{\csubset}{\connecteds{\Xt}}{\seta{\point}\subseteq\csubset}}.
\end{equation}
Thus, according to \refthm{thmunionofafamilyofconnectedsets1},
$\func{\mcp{\Xt}}{\seta{\point}}$
is a connected set of $\Xt$ containing $\point$:
\begin{equation}\label{thmconnectedclassofapointpeq2}
\func{\mcp{\Xt}}{\seta{\point}}\in\connecteds{\Xt}.
\end{equation}
Moreover,
$\func{\mcp{\Xt}}{\seta{\point}}$ can not be a subset of any connected set other than itself, which can be verified in the
following way.
\begin{itemize}
\item[${\textbf{\textsf{p1}}}$]
$\bsubset$
is taken as an arbitrary element of $\connecteds{\Xt}$ such that
\begin{equation}\label{thmconnectedclassofapointp1eq1}
\func{\mcp{\Xt}}{\seta{\point}}\subseteq\bsubset.
\end{equation}
Then considering that $\func{\mcp{\Xt}}{\seta{\point}}$ contains $\point$,
$\bsubset$ also contains $\point$, and hence,
\begin{equation}\label{thmconnectedclassofapointp1eq2}
\bsubset\in\defset{\csubset}{\connecteds{\Xt}}{\seta{\point}\subseteq\csubset},
\end{equation}
and therefore it is clear that,
\begin{align}\label{thmconnectedclassofapointp1eq3}
\bsubset&\subseteq
\union{\defset{\csubset}{\connecteds{\Xt}}{\seta{\point}\subseteq\csubset}}\cr
&=\func{\mcp{\Xt}}{\seta{\point}}.
\end{align}
\Ref{thmconnectedclassofapointp1eq1}
and
\Ref{thmconnectedclassofapointp1eq3}
imply that,
\begin{equation}
\bsubset=\func{\mcp{\Xt}}{\seta{\point}}.
\end{equation}
\endp
\end{itemize}
Therefore,
\begin{equation}\label{thmconnectedclassofapointpeq3}
\Foreach{\aconnectedset}{\[\compl{\connecteds{\Xt}}{\seta{\func{\mcp{\Xt}}{\seta{\point}}}}\]}
\func{\mcp{\Xt}}{\seta{\point}}\nsubseteq\aconnectedset.
\end{equation}
Based on \refdef{defconnectedcomponent},
\Ref{thmconnectedclassofapointpeq2} and
\Ref{thmconnectedclassofapointpeq3}
yield,
\begin{equation}
\func{\mcp{\Xt}}{\seta{\point}}\in\maxcon{\Xt}.
\end{equation}
\endthm
%%%%%%%%%%%%%%%%%%%%%%%%%%%%%%%%%%%%%%%%%%%%%%%%%%%%%%%%%%%%%%%%%%%%%%%%%%%%%%%%%%%%%%%%%%%%%%%%%%%%%%%%%%%%%%%%%%%%%%%%%%%%%%%%
\corollary\label{corconnectedclassofapointisconnected}
$\Xt=\opair{\X}{\topology{}}$
is taken as a topological space.
For every point $\point$ in $\X$,
$\func{\mcp{\Xt}}{\seta{\point}}$ is a connected set of $\Xt$:
\begin{equation}
\Foreach{\point}{\X}
\func{\mcp{\Xt}}{\seta{\point}}\in\connecteds{\Xt}.
\end{equation}
\endcor
%%%%%%%%%%%%%%%%%%%%%%%%%%%%%%%%%%%%%%%%%%%%%%%%%%%%%%%%%%%%%%%%%%%%%%%%%%%%%%%%%%%%%%%%%%%%%%%%%%%%%%%%%%%%%%%%%%%%%%%%%%%%%%%%
\theorem\label{thmconnectedclassofapointcontainsit}
$\Xt=\opair{\X}{\topology{}}$
is taken as a topological space.
For every point $\point$ in $\X$,
$\func{\mcp{\Xt}}{\seta{\point}}$
contains $\point$:
\begin{equation}
\Foreach{\point}{\X}
\point\in\func{\mcp{\Xt}}{\seta{\point}}.
\end{equation}
\proof
$\point$
is taken as an arbitrary element of $\X$. According to \refthm{thmsingletonsubspaceisconnected},
\begin{equation}
\Foreach{\point}{\X}
\seta{\point}\in\connecteds{\Xt},
\end{equation}
and hence,
\begin{equation}
\Foreach{\point}{\X}
\seta{\point}\subseteq
\defset{\csubset}{\connecteds{\Xt}}{\seta{\point}\subseteq\csubset},
\end{equation}
and hence acccording to \refdef{defconnectedclassofaset},
\begin{align}
\Foreach{\point}{\X}
\point&\in\union{\defset{\csubset}{\connecteds{\Xt}}{\seta{\point}\subseteq\csubset}}\cr
&=\func{\mcp{\Xt}}{\seta{\point}}.
\end{align}
\endthm
%%%%%%%%%%%%%%%%%%%%%%%%%%%%%%%%%%%%%%%%%%%%%%%%%%%%%%%%%%%%%%%%%%%%%%%%%%%%%%%%%%%%%%%%%%%%%%%%%%%%%%%%%%%%%%%%%%%%%%%%%%%%%%%%
\theorem\label{thmconnectedcomponentanditspoints}
$\Xt=\opair{\X}{\topology{}}$
is taken as a topological space.
\begin{equation}\label{thmconnectedcomponentanditspointseq1}
\Foreach{\csubset}{\maxcon{\Xt}}
\bigg[\Foreach{\point}{\csubset}
\csubset=\func{\mcp{\Xt}}{\seta{\point}}\bigg].
\end{equation}
\proof
$\csubset$
is taken as an arbitrary element of $\maxcon{\Xt}$
(a maximally-connected set of $\Xt$). The according to \refdef{defconnectedcomponent},
$\csubset$
is a connected set of $\Xt$:
\begin{equation}\label{thmconnectedcomponentanditspointspeq1}
\csubset\in\connecteds{\X}.
\end{equation}
\begin{itemize}
\item[${\textbf{\textsf{p1}}}$]
$\point$
is taken as an arbitrary element of $\csubset$. Then according to \Ref{thmconnectedcomponentanditspointspeq1},
\begin{equation}
\csubset\in\defset{\csubset}{\connecteds{\Xt}}{\seta{\point}\subseteq\csubset},
\end{equation}
and hence according to \refdef{defconnectedclassofaset},
\begin{align}
\csubset&\subseteq\union{\defset{\csubset}{\connecteds{\Xt}}{\seta{\point}\subseteq\csubset}}\cr
&=\func{\mcp{\Xt}}{\seta{\x}}.
\end{align}
Therefore considering that $\csubset$ is amaximally-connected set of $\Xt$, and
$\func{\mcp{\Xt}}{\seta{\point}}$ is a connected set of $\Xt$
(\refcor{corconnectedclassofapointisconnected}),
based on \refdef{defconnectedcomponent}, it ibecomes clear that,
\begin{equation}
\csubset=\func{\mcp{\Xt}}{\seta{\point}}.
\end{equation}
\endp
\end{itemize}
\endthm
%%%%%%%%%%%%%%%%%%%%%%%%%%%%%%%%%%%%%%%%%%%%%%%%%%%%%%%%%%%%%%%%%%%%%%%%%%%%%%%%%%%%%%%%%%%%%%%%%%%%%%%%%%%%%%%%%%%%%%%%%%%%%%%%
\theorem\label{thmmaximallyconnectedsetsandconnectedclassofpoints}
$\Xt=\opair{\X}{\topology{}}$
is taken as a topological space.
\begin{equation}\label{thmmaximallyconnectedsetsandconnectedclassofpointseq1}
\maxcon{\Xt}=\defset{\csubset}{\CSs{\X}}
{\[\Exists{\point}{\X}\csubset=\func{\mcp{\Xt}}{\seta{\point}}\]}.
\end{equation}
\proof
According to \refthm{thmconnectedclassofapoint},
\begin{equation}\label{thmmaximallyconnectedsetsandconnectedclassofpointspeq1}
\defset{\csubset}{\CSs{\X}}
{\[\Exists{\point}{\X}\csubset=\func{\mcp{\Xt}}{\seta{\point}}\]}
\subseteq\maxcon{\Xt}.
\end{equation}
\begin{itemize}
\item[${\textbf{\textsf{p1}}}$]
$\csubset$
is taken as an arbitrary element of $\maxcon{\Xt}$
(a maximally-connected set of $\Xt$). According to \refdef{defconnectedcomponent},
$\csubset$ is non-empty, and hence possesses an element like $\x$:
\begin{equation}
\Existsis{\x}{\X}\x\in\csubset.
\end{equation}
Then according to \refdef{thmconnectedcomponentanditspoints},
\begin{equation}
\csubset=\func{\mcp{\Xt}}{\seta{\x}},
\end{equation}
and hence,
\begin{equation}
\csubset\in
\defset{\csubset}{\CSs{\X}}
{\[\Exists{\point}{\X}\csubset=\func{\mcp{\Xt}}{\seta{\point}}\]}.
\end{equation}
\endp
\end{itemize}
Therefore,
\begin{equation}\label{thmmaximallyconnectedsetsandconnectedclassofpointspeq2}
\defset{\csubset}{\CSs{\X}}
{\[\Exists{\point}{\X}\csubset=\func{\mcp{\Xt}}{\seta{\point}}\]}\subseteq
\maxcon{\Xt}.
\end{equation}
\Ref{thmmaximallyconnectedsetsandconnectedclassofpointspeq1} and
\Ref{thmmaximallyconnectedsetsandconnectedclassofpointspeq2} clearly imply
\Ref{thmmaximallyconnectedsetsandconnectedclassofpointseq1}.
\endthm
%%%%%%%%%%%%%%%%%%%%%%%%%%%%%%%%%%%%%%%%%%%%%%%%%%%%%%%%%%%%%%%%%%%%%%%%%%%%%%%%%%%%%%%%%%%%%%%%%%%%%%%%%%%%%%%%%%%%%%%%%%%%%%%%
\theorem\label{thmconnectedcomponentscoverspace}
$\Xt=\opair{\X}{\topology{}}$
is taken as a topological space.
$\opair{\X}{\maxcon{\Xt}}$ is a cover of $\X$:
\begin{equation}
\maxcon{\Xt}\in\covers{\X}.
\end{equation}
\proof
According to \refthm{thmconnectedclassofapointcontainsit},
\begin{align}
\union{\defset{\csubset}{\CSs{\X}}
{\[\Exists{\point}{\X}\csubset=\func{\mcp{\Xt}}{\seta{\point}}\]}}&=
\Union{\point}{\X}{\func{\mcp{\Xt}}{\seta{\point}}}\cr
&=\X,
\end{align}
and hence according to \refthm{thmmaximallyconnectedsetsandconnectedclassofpoints},
\begin{equation}
\union{\[\maxcon{\Xt}\]}=\X.
\end{equation}
According to \refdef{defcovermapofset}, this means $\opair{\X}{\maxcon{\Xt}}$
is a cover of $\X$.
\endthm
%%%%%%%%%%%%%%%%%%%%%%%%%%%%%%%%%%%%%%%%%%%%%%%%%%%%%%%%%%%%%%%%%%%%%%%%%%%%%%%%%%%%%%%%%%%%%%%%%%%%%%%%%%%%%%%%%%%%%%%%%%%%%%%%
\theorem\label{thmconnectedcomponentsaredistinct}
$\Xt=\opair{\X}{\topology{}}$
is taken as a topological space.
Every pair of different maximally-connected sets of $\Xt$ are disjoint. That is,
\begin{equation}
\Foreach{\opair{\asubset}{\bsubset}}{\[\maxcon{\Xt}\times\maxcon{\Xt}\]}
\bigg(\asubset\cap\bsubset\neq\empty\then\asubset=\bsubset\bigg).
\end{equation}
\proof
Each
$\asubset$
and
$\bsubset$
is taken as an arbitrary element of $\maxcon{\Xt}$
(a maximally-connected set of $\Xt$), and it is assumed that,
\begin{equation}
\asubset\cap\bsubset\neq\empty.
\end{equation}
Therefore there exists a point like $\x$ in $\asubset\cap\bsubset$:
\begin{equation}
\Existsis{\x}{\X}\x\in\(\asubset\cap\bsubset\).
\end{equation}
Then according to \refthm{thmconnectedcomponentanditspoints},
\begin{align}
\asubset&=\func{\mcp{\Xt}}{\seta{\x}},\\
\bsubset&=\func{\mcp{\Xt}}{\seta{\x}},
\end{align}
and hence,
\begin{equation}
\asubset=\bsubset.
\end{equation}
\endthm
%%%%%%%%%%%%%%%%%%%%%%%%%%%%%%%%%%%%%%%%%%%%%%%%%%%%%%%%%%%%%%%%%%%%%%%%%%%%%%%%%%%%%%%%%%%%%%%%%%%%%%%%%%%%%%%%%%%%%%%%%%%%%%%%
\theorem\label{thmconnectedcomponentsformapartitionofspace}
$\Xt=\opair{\X}{\topology{}}$
is taken as a topological space.
The collection of all maximally-connected set of $\Xt$ ($\maxcon{\Xt}$) is a partition of $\X$:
\begin{equation}
\maxcon{\Xt}\in\Cpart{\X}.
\end{equation}
\proof
According to
\refdef{defconnectedcomponent},
\refthm{thmconnectedcomponentscoverspace},
and
\refthm{thmconnectedcomponentsaredistinct},
$\maxcon{\Xt}$
is such an element of $\CSs{\CSs{\X}}$ that covers $\X$, and its elements are non-empty and pair-wise disjoint.
Thus, according to \refdef{defpartitionofset}, this set is a partition of $\X$
\footnote{When $\X$ equals $\empty$, the collection of all maximally-connected sets of $\Xt$
must equal $\empty$, and hence is a partition of $\X$.
}
\endthm
%%%%%%%%%%%%%%%%%%%%%%%%%%%%%%%%%%%%%%%%%%%%%%%%%%%%%%%%%%%%%%%%%%%%%%%%%%%%%%%%%%%%%%%%%%%%%%%%%%%%%%%%%%%%%%%%%%%%%%%%%%%%%%%%
\definition\label{defmaximallyconnectedrelation}
$\Xt=\opair{\X}{\topology{}}$
is taken as a topological space.
The relation $\mcR{\X}{\topology{}}$ in $\X$ is defined as,
\begin{align}
\mcR{\X}{\topology{}}:=\defset{\opair{\point_1}{\point_2}}{\(\X\times\X\)}
{\[\Exists{\csubset}{\maxcon{\Xt}}\seta{\binary{\point_1}{\point_2}}\subseteq\csubset\]}.
\end{align}
This means, two points $\point_1$ and $\point_2$ of $\Xt$
are related by
$\mcR{\X}{\topology{}}$
iff they belong to the same maximally-connected set of $\Xt$.
\endef
%%%%%%%%%%%%%%%%%%%%%%%%%%%%%%%%%%%%%%%%%%%%%%%%%%%%%%%%%%%%%%%%%%%%%%%%%%%%%%%%%%%%%%%%%%%%%%%%%%%%%%%%%%%%%%%%%%%%%%%%%%%%%%%%
\theorem\label{thmmaximallyconnectedrelationisequivalence}
$\Xt=\opair{\X}{\topology{}}$
is taken as a topological space.
The relation $\mcR{\X}{\topology{}}$ in $\X$ is an equivalence relation.
\begin{equation}
\mcR{\X}{\topology{}}\in\EqR{\X}.
\end{equation}
\prooff
\begin{itemize}
\item[${\textbf{\textsf{p1}}}$]
According to \refthm{thmconnectedclassofapoint} and
\refthm{thmconnectedclassofapointcontainsit},
\begin{align}
\Foreach{\point}{\X}
\Existsis{\func{\mcp{\Xt}}{\seta{\point}}}{\maxcon{\Xt}}
\bigg(\seta{\binary{\point}{\point}}&=\seta{\point}\cr
&\subseteq\func{\mcp{\Xt}}{\seta{\point}}\bigg).\cr
&{}
\end{align}
Thus according to \refdef{defmaximallyconnectedrelation},
\begin{equation}
\Foreach{\point}{\X}
\opair{\point}{\point}\in\mcR{\X}{\topology{}}.
\end{equation}
\endp
\end{itemize}
\begin{itemize}
\item[${\textbf{\textsf{p2}}}$]
$\opair{\point_1}{\point_2}$
is taken as an arbitrary element of $\mcR{\X}{\topology{}}$. Then according to \refdef{defmaximallyconnectedrelation},
\begin{align}
\Existsis{\csubset}{\maxcon{\Xt}}
\seta{\binary{\point_1}{\point_2}}\subseteq\csubset,
\end{align}
and hence according to \refdef{defmaximallyconnectedrelation}, clearly,
\begin{equation}
\opair{\point_2}{\point_1}\in\mcR{\X}{\topology{}}.
\end{equation}
\endp
\end{itemize}
\begin{itemize}
\item[${\textbf{\textsf{p3}}}$]
$\point_1$,
$\point_2$,
and
$\point_3$
are taken as such arbitrary elements of $\X$ such that,
\begin{align}
\opair{\point_1}{\point_2}&\in\mcR{\X}{\topology{}},
\label{thmmaximallyconnectedrelationisequivalencep3eq1}\\
\opair{\point_2}{\point_3}&\in\mcR{\X}{\topology{}}.
\label{thmmaximallyconnectedrelationisequivalencep3eq2}
\end{align}
Then based on \refdef{defmaximallyconnectedrelation},
\begin{align}
&\Existsis{\csubset}{\maxcon{\Xt}}\seta{\binary{\point_1}{\point_2}}\subseteq\csubset,
\label{thmmaximallyconnectedrelationisequivalencep3eq3}\\
&\Existsis{\p{\csubset}}{\maxcon{\Xt}}\seta{\binary{\point_2}{\point_3}}\subseteq\p{\csubset},
\label{thmmaximallyconnectedrelationisequivalencep3eq4}
\end{align}
and hence according to \refthm{thmconnectedcomponentanditspoints},
\begin{align}
\csubset&=\func{\mcp{\Xt}}{\seta{\point_2}},
\label{thmmaximallyconnectedrelationisequivalencep3eq5}\\
\p{\csubset}&=\func{\mcp{\Xt}}{\seta{\point_2}},
\label{thmmaximallyconnectedrelationisequivalencep3eq6}
\end{align}
and hence,
\begin{equation}\label{thmmaximallyconnectedrelationisequivalencep3eq7}
\csubset=\p{\csubset}.
\end{equation}
\Ref{thmmaximallyconnectedrelationisequivalencep3eq3},
\Ref{thmmaximallyconnectedrelationisequivalencep3eq4},
and
\Ref{thmmaximallyconnectedrelationisequivalencep3eq7}
imply that,
\begin{equation}
\Existsis{\csubset}{\maxcon{\Xt}}
\seta{\binary{\point_1}{\point_2}}\subseteq\csubset,
\end{equation}
which according to \refdef{defmaximallyconnectedrelation}, means,
\begin{equation}
\opair{\point_1}{\point_2}\in\mcR{\X}{\topology{}}.
\end{equation}
\endp
\end{itemize}
Therefore,
\begin{align}
&\Foreach{\point}{\X}\opair{\point}{\point}\in\mcR{\X}{\topology{}},\\
&\Foreach{\opair{\point_1}{\point_2}}{\(\X\times\X\)}
\bigg[\opair{\point_1}{\point_2}\in\mcR{\X}{\topology{}}\then\opair{\point_2}{\point_1}\in\mcR{\X}{\topology{}}\bigg],\\
&\Foreach{\opair{\binary{\point_1}{\point_2}}{\point_3}}{\(\X\times\X\times\X\)}\cr
&\[\bigg(\opair{\point_1}{\point_2}\in\mcR{\X}{\topology{}},~
\opair{\point_2}{\point_3}\in\mcR{\X}{\topology{}}\bigg)
\then\opair{\point_1}{\point_3}\in\mcR{\X}{\topology{}}\].\cr
&{}
\end{align}
Thus $\mcR{\X}{\topology{}}$ is an equivalence relation in $\X$.
\endthm
%%%%%%%%%%%%%%%%%%%%%%%%%%%%%%%%%%%%%%%%%%%%%%%%%%%%%%%%%%%%%%%%%%%%%%%%%%%%%%%%%%%%%%%%%%%%%%%%%%%%%%%%%%%%%%%%%%%%%%%%%%%%%%%%
\theorem\label{thmmaximallyconnectedrelationequivalenceclasses}
$\Xt=\opair{\X}{\topology{}}$
is taken as a topological space.
The collection of all equivalence classes of the equivalence relation $\mcR{\X}{\topology{}}$
equals the collection of all maximally-connected sets of $\Xt$:
\begin{equation}\label{thmmaximallyconnectedrelationequivalenceclasseseq1}
\EqClass{\X}{\mcR{\X}{\topology{}}}=\maxcon{\Xt}.
\end{equation}
\prooff
It is known that,
\begin{equation}\label{thmmaximallyconnectedrelationequivalenceclassespeq1}
\EqClass{\X}{\mcR{\X}{\topology{}}}=\defSet{\pEqclass{\point}{\mcR{\X}{\topology{}}}}{\point\in\X},
\end{equation}
where,
\begin{equation}\label{thmmaximallyconnectedrelationequivalenceclassespeq2}
\Foreach{\point}{\X}
\pEqclass{\point}{\mcR{\X}{\topology{}}}=
\defset{\x}{\X}
{\opair{\point}{\x}\in\mcR{\X}{\topology{}}}.
\end{equation}
$\pEqclass{\point}{\mcR{\X}{\topology{}}}$
is the equivalence class of the point $\point$ in the equivalence relation $\mcR{\X}{\topology{}}$
(which can be the same for different points).
\begin{itemize}
\item[${\textbf{\textsf{p1}}}$]
$\point$
is taken as an arbitrary point of $\Xt$.
\begin{itemize}
\item[${\textbf{\textsf{p1-1}}}$]
$\x$
is taken as an arbitrary point of $\pEqclass{\point}{\mcR{\X}{\topology{}}}$.
Then according to \Ref{thmmaximallyconnectedrelationequivalenceclassespeq2},
\begin{equation}\label{thmmaximallyconnectedrelationequivalenceclassesp1-1eq1}
\opair{\point}{\x}\in\mcR{\X}{\topology{}},
\end{equation}
and hence according to \refdef{defmaximallyconnectedrelation},
\begin{equation}\label{thmmaximallyconnectedrelationequivalenceclassesp1-1eq2}
\Existsis{\csubset}{\maxcon{\Xt}}
\seta{\binary{\point}{\x}}\subseteq\csubset,
\end{equation}
and thus according to \refthm{thmconnectedcomponentanditspoints},
\begin{equation}\label{thmmaximallyconnectedrelationequivalenceclassesp1-1eq3}
\csubset=\func{\mcp{\Xt}}{\seta{\point}}.
\end{equation}
\Ref{thmmaximallyconnectedrelationequivalenceclassesp1-1eq2} and
\Ref{thmmaximallyconnectedrelationequivalenceclassesp1-1eq3}
imply that,
\begin{equation}
\x\in\func{\mcp{\Xt}}{\seta{\point}}.
\end{equation}
\endp
\end{itemize}
\begin{itemize}
\item[${\textbf{\textsf{p1-2}}}$]
$\x$
is taken as an arbitrary point of $\func{\mcp{\Xt}}{\seta{\point}}$. Then according to \refthm{thmconnectedclassofapoint} and
\refthm{thmconnectedclassofapointcontainsit},
\begin{equation}
\Existsis{\func{\mcp{\Xt}}{\seta{\point}}}{\maxcon{\Xt}}
\seta{\binary{\point}{\x}}\subseteq\func{\mcp{\Xt}}{\seta{\point}},
\end{equation}
and hence according to \refdef{defmaximallyconnectedrelation},
\begin{equation}
\opair{\point}{\x}\in\mcR{\X}{\topology{}},
\end{equation}
and thus according to \Ref{thmmaximallyconnectedrelationequivalenceclassespeq2},
\begin{equation}
\x\in\pEqclass{\point}{\mcR{\X}{\topology{}}}.
\end{equation}
\endp
\end{itemize}
\endp
\end{itemize}
Therefore,
\begin{equation}
\Foreach{\point}{\X}
\pEqclass{\point}{\mcR{\X}{\topology{}}}=
\func{\mcp{\Xt}}{\seta{\point}},
\end{equation}
and thus according to \Ref{thmmaximallyconnectedrelationequivalenceclassespeq1},
\begin{equation}
\EqClass{\X}{\mcR{\X}{\topology{}}}=
\defSet{\func{\mcp{\Xt}}{\seta{\point}}}{\point\in\X}.
\end{equation}
According to \refthm{thmmaximallyconnectedsetsandconnectedclassofpoints}, this yields
\Ref{thmmaximallyconnectedrelationequivalenceclasseseq1}.
\endthm
%%%%%%%%%%%%%%%%%%%%%%%%%%%%%%%%%%%%%%%%%%%%%%%%%%%%%%%%%%%%%%%%%%%%%%%%%%%%%%%%%%%%%%%%%%%%%%%%%%%%%%%%%%%%%%%%%%%%%%%%%%%%%%%%
\theorem\label{thmconnectedcomponentsareclosed}
$\Xt=\opair{\X}{\topology{}}$
is taken as a topological space.
Every maximally-connected set of $\Xt$, is a closed set of $\Xt$. That is,
\begin{equation}
\maxcon{\Xt}\subseteq\Fclosed{\X}{\topology{}}.
\end{equation}
\proof
$\csubset$
is taken as an arbitrary element of $\maxcon{\Xt}$. Then according to \refcor{corconnectedclassofapointisconnected},
$\csubset$ is a connected set of $\Xt$:
\begin{equation}\label{thmconnectedcomponentsareclosedpeq1}
\csubset\in\connecteds{\Xt},
\end{equation}
and hence according to \refthm{thmclosureofconnectedsetisconnected},
$\func{\Cl{\Xt}}{\csubset}$ is also a connected set of $\Xt$:
\begin{equation}\label{thmconnectedcomponentsareclosedpeq2}
\func{\Cl{\Xt}}{\csubset}\in\connecteds{\Xt}.
\end{equation}
In addition, according to \refcor{corclosureofset0}, it is clear that
\begin{equation}\label{thmconnectedcomponentsareclosedpeq3}
\csubset\subseteq\func{\Cl{\Xt}}{\csubset}.
\end{equation}
According to \refdef{defconnectedcomponent},
and considering that $\csubset$ is a maximally-connected set of $\Xt$,
\Ref{thmconnectedcomponentsareclosedpeq2}
and
\Ref{thmconnectedcomponentsareclosedpeq3}
imply that,
\begin{equation}
\func{\Cl{\Xt}}{\csubset}=\csubset,
\end{equation}
and hence according to \refthm{thmclosureofclosedset},
\begin{equation}
\csubset\in\Fclosed{\X}{\topology{}}.
\end{equation}
\endthm
%%%%%%%%%%%%%%%%%%%%%%%%%%%%%%%%%%%%%%%%%%%%%%%%%%%%%%%%%%%%%%%%%%%%%%%%%%%%%%%%%%%%%%%%%%%%%%%%%%%%%%%%%%%%%%%%%%%%%%%%%%%%%%%%
\theorem\label{thmconnectedsetaresubsetsofconnectedcomponents}
$\Xt=\opair{\X}{\topology{}}$
is taken as a non-empty topological space.
($\X\neq\empty$).
Every connected set of $\Xt$ is a subset of a maximally-connected set of $\Xt$. That is,
\begin{equation}
\Foreach{\csubset}{\connecteds{\Xt}}
\Exists{\p{\csubset}}{\maxcon{\Xt}}
\csubset\subseteq\p{\csubset}.
\end{equation}
\proof
Considering the non-emptiness of $\X$,
\begin{equation}
\maxcon{\Xt}\neq\empty.
\end{equation}
So, because
$\empty$
is a subset of every maximally-connected set of $\Xt$, it is enough to show that every non-empty connected set of $\Xt$
is a subset of a maximally-connected set of $\Xt$.\\
$\csubset$
is taken as an arbitrary element of $\compl{\connecteds{\Xt}}{\seta{\empty}}$.
Then,
\begin{equation}
\Existsis{\point}{\X}
\point\in\csubset.
\end{equation}
Thus according to \refdef{defconnectedclassofaset},
\begin{equation}
\csubset\subseteq\func{\mcp{\Xt}}{\seta{\point}},
\end{equation}
and hence according to \refthm{thmconnectedclassofapoint},
\begin{equation}
\Existsis{\func{\mcp{\Xt}}{\seta{\point}}}{}{\maxcon{\Xt}}
\csubset\subseteq\func{\mcp{\Xt}}{\seta{\point}}.
\end{equation}
\endthm
%%%%%%%%%%%%%%%%%%%%%%%%%%%%%%%%%%%%%%%%%%%%%%%%%%%%%%%%%%%%%%%%%%%%%%%%%%%%%%%%%%%%%%%%%%%%%%%%%%%%%%%%%%%%%%%%%%%%%%%%%%%%%%%%
\theorem\label{thmcontinuousimageofconnectedcomponents}
Each
$\Xt=\opair{\X}{\topology{\X}}$
and
$\Yt=\opair{\Y}{\topology{\Y}}$
is taken as a topological space.
For every continuous map $\cf$ from $\Xt$ to $\Yt$, the image under $\cf$ of every maximally-connected set of $\Xt$
is a subset of a maximally-connected set of $\Yt$:
\begin{align}\label{thmcontinuousimageofconnectedcomponentseq1}
\Foreach{\cf}{\CF{\Xt}{\Yt}}
\[\Foreach{\Csubset{\X}}{\maxcon{\Xt}}
\bigg(\Exists{\Csubset{\Y}}{\maxcon{\Yt}}
\Csubset{\X}\subseteq\Csubset{\Y}\bigg)\].
\end{align}
\proof
Considering that every maximally-connected set of $\Xt$ is a connected set of $\Xt$
(\refcor{corconnectedclassofapointisconnected}),
and for every continuous map from $\Xt$ to $\Yt$, the image under $\cf$ of every connected set of $\Xt$
is a connected set of $\Yt$ (\refthm{thmconnectedsetsandcontinuousfunctions}),
and every connected set of $\Yt$ is a subset of a maximally-connected set of $\Yt$
\footnote{when there exists a map from $\X$ to $\Y$, $\Y$ is non-empty, and hence possesses at least one
maximally-connected set.
}
(\refthm{thmconnectedsetaresubsetsofconnectedcomponents}),
\Ref{thmcontinuousimageofconnectedcomponentseq1} is clearly seen to be true.
\endthm
%%%%%%%%%%%%%%%%%%%%%%%%%%%%%%%%%%%%%%%%%%%%%%%%%%%%%%%%%%%%%%%%%%%%%%%%%%%%%%%%%%%%%%%%%%%%%%%%%%%%%%%%%%%%%%%%%%%%%%%%%%%%%%%%
\theorem\label{thmcardinalityofconnectedcomponentsisatopp}
Each
$\Xt=\opair{\X}{\topology{\X}}$
and
$\Yt=\opair{\Y}{\topology{\Y}}$
is taken as a topological space.
For every homeomorphism $\hf$ from $\Xt$ to $\Yt$,
the image under $\hf$ of every maximally-connected set of $\Xt$
is a maximally-connected set of $\Yt$:
\begin{equation}
\Foreach{\hf}{\HOF{\Xt}{\Yt}}
\bigg[\Foreach{\Csubset{\X}}{\maxcon{\Xt}}
\func{\image{\hf}}{\Csubset{\X}}\in\maxcon{\Yt}\bigg].
\end{equation}
\proof
$\hf$
is taken as an arbitrary element of $\HOF{\Xt}{\Yt}$
Then according to \refdef{defhomeomorphism},
\begin{align}
\hf&\in\IF{\X}{\Y},
\label{thmcardinalityofconnectedcomponentsisatopppeq1}\\
\hf&\in\CF{\Xt}{\Yt},
\label{thmcardinalityofconnectedcomponentsisatopppeq2}\\
\finv{\hf}&\in\CF{\Yt}{\Xt}.
\label{thmcardinalityofconnectedcomponentsisatopppeq3}
\end{align}
\begin{itemize}
\item[${\textbf{\textsf{p1}}}$]
$\Csubset{\X}$
is taken as an arbitrary element of $\maxcon{\Xt}$. Then according to \Ref{thmcardinalityofconnectedcomponentsisatopppeq2}, and
\refthm{thmcontinuousimageofconnectedcomponents},
there exists a maximally-connected set in $\Yt$ that includes
$\func{\image{\hf}}{\Csubset{\X}}$:
\begin{equation}\label{thmcardinalityofconnectedcomponentsisatoppp1eq1}
\Existsis{\Csubset{\Y}}{\maxcon{\Yt}}
\func{\image{\hf}}{\Csubset{\X}}\subseteq\Csubset{\Y}.
\end{equation}
Considering \Ref{thmcardinalityofconnectedcomponentsisatopppeq3},
$\finv{\hf}$
is a continuous map from $\Yt$ to $\Xt$, and hence according to \refthm{thmcontinuousimageofconnectedcomponents},
considering that $\Csubset{\Y}$ is a maximally-connected set of $\Yt$, the image under $\finv{\hf}$ of $\Csubset{\Y}$
is a subset of a maximally-connected set of $\Xt$:
\begin{equation}\label{thmcardinalityofconnectedcomponentsisatoppp1eq2}
\Existsis{\p{\Csubset{\X}}}{\maxcon{\Xt}}
\func{\image{\[\finv{\hf}\]}}{\Csubset{\Y}}\subseteq\p{\Csubset{\X}}.
\end{equation}
According to \Ref{thmcardinalityofconnectedcomponentsisatoppp1eq1},
\begin{equation}\label{thmcardinalityofconnectedcomponentsisatoppp1eq3}
\func{\image{\[\finv{\hf}\]}}{\func{\image{\hf}}{\Csubset{\X}}}
\subseteq\func{\image{\[\finv{\hf}\]}}{\Csubset{\Y}}.
\end{equation}
Moreover, it is clear that,
\begin{equation}\label{thmcardinalityofconnectedcomponentsisatoppp1eq4}
\func{\image{\[\finv{\hf}\]}}{\func{\image{\hf}}{\Csubset{\X}}}=
\Csubset{\X}.
\end{equation}
\Ref{thmcardinalityofconnectedcomponentsisatoppp1eq2},
\Ref{thmcardinalityofconnectedcomponentsisatoppp1eq3},
and
\Ref{thmcardinalityofconnectedcomponentsisatoppp1eq4}
imply that,
\begin{equation}\label{thmcardinalityofconnectedcomponentsisatoppp1eq5}
\Csubset{\X}\subseteq\p{\Csubset{\X}}.
\end{equation}
\Ref{thmcardinalityofconnectedcomponentsisatopppeq1},
\Ref{thmcardinalityofconnectedcomponentsisatoppp1eq2},
and
\Ref{thmcardinalityofconnectedcomponentsisatoppp1eq5}
imply that,
\begin{align}\label{thmcardinalityofconnectedcomponentsisatoppp1eq6}
\Csubset{\Y}&=\func{\image{\hf}}{\func{\pimage{\hf}}{\Csubset{\Y}}}\cr
&=\func{\image{\hf}}{\func{\image{\[\finv{\h}\]}{\Csubset{\Y}}}}\cr
&\subseteq\func{\image{\hf}}{\p{\Csubset{\X}}}\cr
&=\func{\image{\hf}}{\Csubset{\X}}.
\end{align}
\Ref{thmcardinalityofconnectedcomponentsisatoppp1eq1}
and
\Ref{thmcardinalityofconnectedcomponentsisatoppp1eq6}
imply that,
\begin{equation}
\func{\image{\hf}}{\Csubset{\X}}=\Csubset{\Y},
\end{equation}
and hence according to \Ref{thmcardinalityofconnectedcomponentsisatoppp1eq1},
\begin{equation}
\func{\image{\hf}}{\Csubset{\X}}\in\maxcon{\Yt}.
\end{equation}
\endp
\end{itemize}
\endthm
%%%%%%%%%%%%%%%%%%%%%%%%%%%%%%%%%%%%%%%%%%%%%%%%%%%%%%%%%%%%%%%%%%%%%%%%%%%%%%%%%%%%%%%%%%%%%%%%%%%%%%%%%%%%%%%%%%%%%%%%%%%%%%%%
\theorem\label{thmcardinalityofconnectedcomponentsisatopp1}
Each $\Xt=\opair{\X}{\topology{\X}}$ and
$\Yt=\opair{\Y}{\topology{\Y}}$
is taken as a non-empty topological space.
($\X\neq\empty,~\Y\neq\empty$).
If $\Xt$ and $\Yt$ are homeomorphic, then there exists a bijection between the collections of all
maximally-connected sets of $\Xt$ and $\Yt$. That is,
\begin{equation}
\bigg(\homeomorphic{\Xt}{\Yt}\bigg)\then
\bigg(\Card{\maxcon{\Xt}}\cardeq\Card{\maxcon{\Yt}}\bigg).
\end{equation}
\proof
It is assumed that,
$\homeomorphic{\Xt}{\Yt}$.
This means there exists a homeomorphism like $\hf$ from $\Xt$ to $\Yt$. According to \refthm{thmcardinalityofconnectedcomponentsisatopp},
\begin{equation}
\codomain{\func{\resd{\image{\hf}}}{\maxcon{\Xt}}}\subseteq\maxcon{\Yt},
\end{equation}
and hence
$\func{\rescd{\func{\resd{\image{\hf}}}{\maxcon{\Xt}}}}{\maxcon{\Yt}}$ is well-defined. Considering that $\hf$
is an element of $\HOF{\Xt}{\Yt}$, according to \refdef{defhomeomorphism},
\begin{equation}
\hf\in\IF{\X}{\Y},
\end{equation}
and hence
$\func{\rescd{\func{\resd{\image{\hf}}}{\maxcon{\Xt}}}}{\maxcon{\Yt}}$
is injective.
\begin{align}
&\Foreach{\opair{\asubset}{\bsubset}}{\maxcon{\Xt}\times\maxcon{\Xt}}\cr
&~~~\(\func{\[\func{\rescd{\func{\resd{\image{\hf}}}{\maxcon{\Xt}}}}{\maxcon{\Yt}}\]}{\asubset}=
\func{\func{\rescd{\func{\resd{\image{\hf}}}{\maxcon{\Xt}}}}{\maxcon{\Yt}}}{\bsubset}\)\cr
&\then\func{\image{\hf}}{\asubset}=\func{\image{\hf}}{\bsubset}\cr
&\then\asubset=\bsubset.
\end{align}
\begin{itemize}
\item[${\textbf{\textsf{p1}}}$]
$\Csubset{\Y}$ is taken as an arbitrary element of $\maxcon{\Yt}$. Considering that $\hf$ is a homeomorphism from $\Xt$
to $\Yt$, according to \refdef{defhomeomorphism},
$\finv{\hf}$ is a continuous map from $\Yt$ to $\Xt$, and hence according to \refthm{thmcardinalityofconnectedcomponentsisatopp},
\begin{equation}
\func{\image{\[\finv{\hf}\]}}{\Csubset{\Y}}\in\connecteds{\Xt}.
\end{equation}
Moreover, it is clear that,
\begin{align}
\func{\func{\rescd{\func{\resd{\image{\hf}}}{\maxcon{\Xt}}}}{\maxcon{\Yt}}}
{\func{\image{\[\finv{\hf}\]}}{\Csubset{\Y}}}&=
\func{\image{\hf}}{\func{\image{\[\finv{\hf}\]}}{\Csubset{\Y}}}\cr
&=\Csubset{\Y}.
\end{align}
\endp
\end{itemize}
Therefore,
\begin{align}
\Foreach{\Csubset{\Y}}{\maxcon{\Yt}}
\Existsis{\func{\image{\[\finv{\hf}\]}}{\Csubset{\Y}}}{\maxcon{\Xt}}
\func{\image{\hf}}{\func{\image{\[\finv{\hf}\]}}{\Csubset{\Y}}}=\Csubset{\Y}.
\end{align}
This means,
$\hf$ is surjective.
\begin{equation}
\func{\rescd{\func{\resd{\image{\hf}}}{\maxcon{\Xt}}}}{\maxcon{\Yt}}
\in\surFunc{\maxcon{\Xt}}{\maxcon{\Yt}}.
\end{equation}
Thus,
$\func{\rescd{\func{\resd{\image{\hf}}}{\maxcon{\Xt}}}}{\maxcon{\Yt}}$
is a bijective map from $\maxcon{\Xt}$ to $\maxcon{\Yt}$. Thus,
\begin{align}
\Existsis{\phi}{\Func{\maxcon{\Xt}}{\maxcon{\Yt}}}
\phi\in\IF{\maxcon{\Xt}}{\maxcon{\Yt}},
\end{align}
which means the cardinality of $\maxcon{\Xt}$ is the same as the cardinality of $\maxcon{\Yt}$:
\begin{equation}
\Card{\maxcon{\Xt}}\cardeq\Card{\maxcon{\Yt}}.
\end{equation}
\endthm
%%%%%%%%%%%%%%%%%%%%%%%%%%%%%%%%%%%%%%%%%%%%%%%%%%%%%%%%%%%%%%%%%%%%%%%%%%%%%%%%%%%%%%%%%%%%%%%%%%%%%%%%%%%%%%%%%%%%%%%%%%%%%%%%
%%%%%%%%%%%%%%%%%%%%%%%%%%%%%%%%%%%%%%%%%%%%%%%%%%%%%%%%%%%%%%%%%%%%%%%%%%%%%%%%%%%%%%%%%%%%%%%%%%%%%%%%%%%%%%%%%%%%%%%%%%%%%%%%
%%%%%%%%%%%%%%%%%%%%%%%%%%%%%%%%%%%%%%%%%%%%%%%%%%%%%%%%%%%%%%%%%%%%%%%%%%%%%%%%%%%%%%%%%%%%%%%%%%%%%%%%%%%%%%%%%%%%%%%%%%%%%%%%
\section{
Totally-Disconnected Topological Spaces
}
\definition\label{deftotallydisconnectedspace}
$\Xt=\opair{\X}{\topology{}}$
is taken as a topological space.
$\Xt$ is referred to as a $\quotl$totally-disconnected topological space$\quotr$ iff
\begin{equation}
\connecteds{\Xt}=\seta{\empty}\cup
\defset{\asubset}{\X}
{\[\Exists{\point}{\X}\asubset=\seta{\point}\]}.
\end{equation}
In other words,
$\Xt$ is called a totally-disconnected topological space iff
the collection of all connected sets of $\Xt$
($\connecteds{\Xt}$) does not have any element other than $\empty$
and the singleton subsets of $\X$.
$\caution$
According to \refthm{thmemptysetisaconnectedsetofeverytopologicalspace} and \refthm{thmsingletonsubspaceisconnected},
it is clear that in any case,
\begin{equation}
\seta{\empty}\cup\defset{\asubset}{\X}
{\[\Exists{\point}{\X}\asubset=\seta{\point}\]}
\subseteq\connecteds{\Xt}.
\end{equation}
\endef
%%%%%%%%%%%%%%%%%%%%%%%%%%%%%%%%%%%%%%%%%%%%%%%%%%%%%%%%%%%%%%%%%%%%%%%%%%%%%%%%%%%%%%%%%%%%%%%%%%%%%%%%%%%%%%%%%%%%%%%%%%%%%%%%
\theorem\label{thmdiscretespaceistotallydisconnected}
$\X$
is taken as a set.
The discrete topological space $\opair{\X}{\CSs{\X}}$ is a totally-disconnected topological space.
\proof
\begin{itemize}
\item[${\textbf{\textsf{p1}}}$]
$\asubset$
is taken as an arbitrary subset of $\X$  other that $\empty$ and singletons. That is, $\asubset$
is taken as an arbitrary element of $\compl{\CSs{\X}}
{\(\seta{\empty}\cup\defset{\asubset}{\X}
{\[\Exists{\point}{\X}\asubset=\seta{\point}\]}\)}$
Then,
$\CarD{\asubset}>1$, and hence $\X$ has an element like $\x$, and
\begin{equation}
\seta{\x}\neq\X.
\end{equation}
Thus according to \refthm{thmdiscretetopologyclopensets},
$\seta{\x}$
is a clopen set of $\opair{\asubset}{\CSs{\asubset}}$ that neither equals $\empty$ nor $\asubset$. That is,
\begin{equation}
\Existsis{\seta{\x}}{\CSs{\asubset}}
\seta{\x}\in\[\CSs{\asubset}\cap\Fclosed{\asubset}{\CSs{\asubset}}\].
\end{equation}
Moreover, according to \refthm{thmsubspacesofdiscretespace},
\begin{equation}
\CSs{\asubset}=\stopology{\CSs{\X}}{\asubset}.
\end{equation}
Therefore,
\begin{equation}
\Existsis{\seta{\x}}{\CSs{\asubset}}
\seta{\x}\in\compl{\(\bigg[\stopology{\CSs{\X}}{\asubset}\cap\Fclosed{\asubset}{\stopology{\CSs{\X}}{\asubset}}\bigg]}
\seta{\binary{\empty}{\asubset}}\).
\end{equation}
This means there exists a non-empty proper subset of $\asubset$ like $\seta{\x}$ that is a clopen set of
$\opair{\asubset}{\stopology{\CSs{\X}}{\asubset}}$.
Thus according to \refdef{defconnectedness},
$\asubset$ is not a connected set of $\Xt$:
\begin{equation}
\asubset\neq\connecteds{\Xt}.
\end{equation}
\endp
\end{itemize}
Therefore,
\begin{equation}
\connecteds{\Xt}=\seta{\empty}\cup
\defset{\asubset}{\X}
{\[\Exists{\point}{\X}\asubset=\seta{\point}\]},
\end{equation}
which according to \refdef{deftotallydisconnectedspace}, means $\Xt$ is a totally-disconnected topological space.
\endthm
%%%%%%%%%%%%%%%%%%%%%%%%%%%%%%%%%%%%%%%%%%%%%%%%%%%%%%%%%%%%%%%%%%%%%%%%%%%%%%%%%%%%%%%%%%%%%%%%%%%%%%%%%%%%%%%%%%%%%%%%%%%%%%%%
\theorem\label{thmemptysetistotallydisconnected}
$\opair{\empty}{\seta{\empty}}$
is a totally-disconnected topological space.
\proof
According to \refthm{thmemptyspaceisconnected} and
\refdef{deftotallydisconnectedspace}, it is clear.
Alternatively, considering that $\opair{\empty}{\seta{\empty}}$
is a discrete topological space, this proposition is an immediate consequence of the \refthm{thmdiscretespaceistotallydisconnected}.
\endthm
%%%%%%%%%%%%%%%%%%%%%%%%%%%%%%%%%%%%%%%%%%%%%%%%%%%%%%%%%%%%%%%%%%%%%%%%%%%%%%%%%%%%%%%%%%%%%%%%%%%%%%%%%%%%%%%%%%%%%%%%%%%%%%%%
\theorem\label{thmsingletonspaceistotallydisconnected}
$\x$ is taken as a set.
The singleton topological space $\opair{\seta{\x}}{\CSs{\X}}$ is a totally-disconnected topological space.
\proof
According to \refthm{thmsingletontopology0},
$\opair{\seta{\x}}{\CSs{\X}}$
is a discrete topological space, and hence according to \refthm{thmdiscretespaceistotallydisconnected},
$\opair{\seta{\x}}{\CSs{\X}}$ is a totally-disconnected topological space.
\endthm
%%%%%%%%%%%%%%%%%%%%%%%%%%%%%%%%%%%%%%%%%%%%%%%%%%%%%%%%%%%%%%%%%%%%%%%%%%%%%%%%%%%%%%%%%%%%%%%%%%%%%%%%%%%%%%%%%%%%%%%%%%%%%%%%
%%%%%%%%%%%%%%%%%%%%%%%%%%%%%%%%%%%%%%%%%%%%%%%%%%%%%%%%%%%%%%%%%%%%%%%%%%%%%%%%%%%%%%%%%%%%%%%%%%%%%%%%%%%%%%%%%%%%%%%%%%%%%%%%
%%%%%%%%%%%%%%%%%%%%%%%%%%%%%%%%%%%%%%%%%%%%%%%%%%%%%%%%%%%%%%%%%%%%%%%%%%%%%%%%%%%%%%%%%%%%%%%%%%%%%%%%%%%%%%%%%%%%%%%%%%%%%%%%
\section{
Local Connectedness
}
\definition\label{deflocallyconnectedness}
$\Xt=\opair{\X}{\topology{}}$
is taken as a topological space.
$\Xt$
is called a $\quotl$locally-connected topological space$\quotr$ iff there exists a base of the topological space $\Xt$,
each element of which is a connected set of $\Xt$. In other words, $\Xt$
 is called a locally-connected topological space iff
\begin{equation}
\Exists{\base}{\Cbase{\Xt}}
\base\subseteq\connecteds{\Xt}.
\end{equation}
\endef
%%%%%%%%%%%%%%%%%%%%%%%%%%%%%%%%%%%%%%%%%%%%%%%%%%%%%%%%%%%%%%%%%%%%%%%%%%%%%%%%%%%%%%%%%%%%%%%%%%%%%%%%%%%%%%%%%%%%%%%%%%%%%%%%
\theorem\label{thmlocallyconnectednessequiv1}
$\Xt=\opair{\X}{\topology{}}$
is taken as a topological space.
$\Xt$
is a locally-connected topological space if and only if for every open set $\U$ of $\Xt$,
every maximally-connected set of the topological space $\opair{\U}{\stopology{\topology{}}{\U}}$
is an open set of $\Xt$. That is,
$\Xt$
is a locally-connected topological space if and only if
\begin{equation}
\Foreach{\U}{\topology{}}
\maxcon{\opair{\U}{\stopology{\topology{}}{\U}}}
\subseteq\topology{}.
\end{equation}
\proof
\begin{itemize}
\item[${\textbf{\textsf{p1}}}$]
It is assumed that
$\Xt$
is locally-connected, that is,
$\Xt$
has a base such as $\base$, each element of which is a connected set of $\Xt$:
\begin{equation}\label{thmlocallyconnectednessequiv1p1eq1}
\Existsis{\base}{\Cbase{\Xt}}
\base\subseteq\connecteds{\Xt}.
\end{equation}
\begin{itemize}
\item[${\textbf{\textsf{p1-1}}}$]
$\U$
is taken as an arbitrary element of $\topology{}$.
$\base_{\U}$ is defined as,
\begin{equation}\label{thmlocallyconnectednessequiv1p1-1eq1}
\base_{\U}:=\defset{\B}{\base}{\B\subseteq\U}.
\end{equation}
According to \refthm{thmbaseofopensubspace},
$\base_{\U}$ is a base for the topological space $\opair{\U}{\stopology{\topology{}}{\U}}$:
\begin{equation}\label{thmlocallyconnectednessequiv1p1-1eq2}
\base_{\U}\in\Cbase{\opair{\U}{\stopology{\topology{}}{\U}}}.
\end{equation}
According to \Ref{thmlocallyconnectednessequiv1p1eq1}
and
\Ref{thmlocallyconnectednessequiv1p1-1eq1},
it is clear that every element of $\base_{\U}$ is a connected set of $\Xt$ and a subset of $\U$:
\begin{equation}\label{thmlocallyconnectednessequiv1p1-1eq3}
\Foreach{\B}{\base_{\U}}\bigg(\B\subseteq\U,~\B\in\connecteds{\Xt}\bigg),
\end{equation}
and hence according to \refthm{thmsubspaceconnectedness},
every element of $\base_{\U}$ is a connected set of the topological space $\opair{\U}{\stopology{\topology{}}{\U}}$:
\begin{equation}\label{thmlocallyconnectednessequiv1p1-1eq4}
\Foreach{\B}{\base_{\U}}
\B\in\connecteds{\opair{\U}{\stopology{\topology{}}{\U}}}.
\end{equation}
\begin{itemize}
\item[${\textbf{\textsf{p1-1-1}}}$]
$\csubset$
is taken as an arbitrary element of $\maxcon{\opair{\U}{\stopology{\topology{}}{\U}}}$
(a maximally-connected set of $\opair{\U}{\stopology{\topology{}}{\U}}$)
Knowing that a base of a topological space covers it
(\refthm{thmbaseNcondition1}),
\begin{equation}\label{thmlocallyconnectednessequiv1p1-1-1eq1}
\Foreach{\point}{\U}
\bigg(\Exists{\B}{\base_{\U}}\point\in\B\bigg),
\end{equation}
and hence considering that
$\csubset\subseteq\U$,
\begin{equation}\label{thmlocallyconnectednessequiv1p1-1-1eq2}
\Foreach{\point}{\csubset}
\bigg(\Exists{\B}{\base_{\U}}\point\in\B\bigg).
\end{equation}
Accordingly, the function $\eta$ is defined as,
\begin{align}
&\eta\in\Func{\csubset}{\CSs{\U}},
\label{thmlocallyconnectednessequiv1p1-1-1eq3}\\
&\Foreach{\point}{\csubset}\point\in\func{\eta}{\point}\in\base_{\U}.
\label{thmlocallyconnectednessequiv1p1-1-1eq4}
\end{align}
Therefore,
\begin{equation}\label{thmlocallyconnectednessequiv1p1-1-1eq5}
\Foreach{\point}{\csubset}
\seta{\point}\subseteq
\[\func{\eta}{\point}\cap\csubset\],
\end{equation}
and moreover, according to \Ref{thmlocallyconnectednessequiv1p1-1eq4},
\begin{equation}\label{thmlocallyconnectednessequiv1p1-1-1eq6}
\Foreach{\point}{\csubset}
\func{\eta}{\point}\in\connecteds{\opair{\U}{\stopology{\topology{}}{\U}}}.
\end{equation}
Moreover, according to \refdef{defconnectedcomponent}, and
Considering that $\csubset$ is a maximally-connected set of $\opair{\U}{\stopology{\topology{}}{\U}}$,
it is a connected set of this topological space too.
\begin{equation}\label{thmlocallyconnectednessequiv1p1-1-1eq7}
\csubset\in\connecteds{\opair{\U}{\stopology{\topology{}}{\U}}}.
\end{equation}
According to \refthm{thmunionofafamilyofconnectedsets1},
\Ref{thmlocallyconnectednessequiv1p1-1-1eq5},
\Ref{thmlocallyconnectednessequiv1p1-1-1eq5},
and
\Ref{thmlocallyconnectednessequiv1p1-1-1eq7}
imply that for every $\point$ in $\csubset$,
$\[\func{\eta}{\point}\cup\csubset\]$ is a connected set of $\opair{\U}{\stopology{\topology{}}{\U}}$,
because the union of a pair of intersecting connected sets is again connected:
\begin{equation}\label{thmlocallyconnectednessequiv1p1-1-1eq8}
\Foreach{\point}{\csubset}
\[\func{\eta}{\point}\cup\csubset\]\in\connecteds{\opair{\U}{\stopology{\topology{}}{\U}}}.
\end{equation}
According to \refdef{defconnectedcomponent},
and considering that $\csubset$ is an element of $\maxcon{\opair{\U}{\stopology{\topology{}}{\U}}}$,
and for every $\point$ in $\csubset$,
$\csubset\subseteq\[\func{\eta}{\point}\cup\csubset\]$,
\Ref{thmlocallyconnectednessequiv1p1-1-1eq7}
implies that,
\begin{equation}\label{thmlocallyconnectednessequiv1p1-1-1eq9}
\Foreach{\point}{\csubset}
\[\func{\eta}{\point}\cup\csubset\]=\csubset.
\end{equation}
This means that,
\begin{equation}\label{thmlocallyconnectednessequiv1p1-1-1eq10}
\Foreach{\point}{\csubset}\func{\eta}{\point}\subseteq\csubset,
\end{equation}
and hence,
\begin{equation}\label{thmlocallyconnectednessequiv1p1-1-1eq11}
\csubset=\Union{\point}{\csubset}{\func{\eta}{\point}}.
\end{equation}
According to \refthm{thmSandNconditionsofbase0},
\Ref{thmlocallyconnectednessequiv1p1-1-1eq4}
and
\Ref{thmlocallyconnectednessequiv1p1-1-1eq11}
imply that,
\begin{equation}
\csubset\in\stopology{\topology{}}{\U},
\end{equation}
and hence considering that $\U\in\topology{}$, and according to \refthm{thmsubsubopen},
\begin{equation}
\csubset\in\topology{}.
\end{equation}
\endp
\end{itemize}
\endp
\end{itemize}
Thus it is seen that,
\begin{equation*}
\Foreach{\U}{\topology{}}
\maxcon{\opair{\U}{\stopology{\topology{}}{\U}}}\subseteq\topology{}.
\end{equation*}
\endp
\end{itemize}
\begin{itemize}
\item[${\textbf{\textsf{p2}}}$]
It is assumed that,
\begin{equation}\label{thmlocallyconnectednessequiv1p2eq1}
\Foreach{\U}{\topology{}}
\maxcon{\opair{\U}{\stopology{\topology{}}{\U}}}\subseteq\topology{}.
\end{equation}
The collection $\base_{\rm{MC}}$ of the subsets of $\X$ is defined as,
\begin{equation}\label{thmlocallyconnectednessequiv1p2eq2}
\base_{\rm{MC}}:=\defset{\B}{\CSs{\X}}
{\[\Exists{\U}{\topology{}}\B\in\maxcon{\opair{\U}{\stopology{\topology{}}{\U}}}\]}.
\end{equation}
This means,
$\base_{\rm{MC}}$
consists of all connected components of open subspaces of $\Xt$.
According to \Ref{thmlocallyconnectednessequiv1p2eq1} and
\Ref{thmlocallyconnectednessequiv1p2eq2},
\begin{align}
&\base_{\rm{MC}}\subseteq\topology{},
\label{thmlocallyconnectednessequiv1p2eq3}\\
&\Foreach{\U}{\topology{}}\maxcon{\opair{\U}{\stopology{\topology{}}{\U}}}
\subseteq\base_{\rm{MC}}.
\label{thmlocallyconnectednessequiv1p2eq4}
\end{align}
it is clear that,
\begin{equation}\label{thmlocallyconnectednessequiv1p2eq5}
\(\union{\base_{\rm{MC}}}\)=
\Union{\U}{\topology{}}{\[\union{\maxcon{\opair{\U}{\stopology{\topology{}}{\U}}}}\]}.
\end{equation}
According to \refthm{thmconnectedcomponentscoverspace} and \refdef{defcovermapofset},
every open set of $\Xt$
equals the union of all of its connected components.
\begin{equation}\label{thmlocallyconnectednessequiv1p2eq6}
\Foreach{\U}{\topology{}}
\U=\[\union{\maxcon{\opair{\U}{\stopology{\topology{}}{\U}}}}\].
\end{equation}
\Ref{thmlocallyconnectednessequiv1p2eq5},
\Ref{thmlocallyconnectednessequiv1p2eq6},
and
\refdef{deftopologicalspace}
imply that,
\begin{align}\label{thmlocallyconnectednessequiv1p2eq7}
\(\union{\base_{\rm{MC}}}\)&=\Union{\U}{\topology{}}{\U}\cr
&=\X.
\end{align}
According to \refdef{defbase},
\Ref{thmlocallyconnectednessequiv1p2eq3},
\Ref{thmlocallyconnectednessequiv1p2eq4},
\Ref{thmlocallyconnectednessequiv1p2eq6},
and
\Ref{thmlocallyconnectednessequiv1p2eq7},
$\base_{\rm{MC}}$
is a base for $\Xt$:
\begin{equation}\label{thmlocallyconnectednessequiv1p2eq8}
\base_{\rm{MC}}\in\Cbase{\Xt}.
\end{equation}
According to \refdef{defconnectedcomponent},
for every $\U$ in $\topology{}$,
every maximally-connected set of the topological space $\opair{\U}{\stopology{\topology{}}{\U}}$
is a connected set of $\opair{\U}{\stopology{\topology{}}{\U}}$:
\begin{equation}\label{thmlocallyconnectednessequiv1p2eq9}
\Foreach{\U}{\topology{}}
\maxcon{\opair{\U}{\stopology{\topology{}}{\U}}}\subseteq
\connecteds{\opair{\U}{\stopology{\topology{}}{\U}}}.
\end{equation}
Moreover, according to \refthm{thmsubspaceconnectedness},
for every $\U$ in $\topology{}$,
every connected set of the topological space $\opair{\U}{\stopology{\topology{}}{\U}}$
is a connected set of $\Xt$:
\begin{equation}\label{thmlocallyconnectednessequiv1p2eq10}
\Foreach{\U}{\topology{}}
\connecteds{\opair{\U}{\stopology{\topology{}}{\U}}}\subseteq\connecteds{\Xt}.
\end{equation}
\Ref{thmlocallyconnectednessequiv1p2eq9}
and
\Ref{thmlocallyconnectednessequiv1p2eq10}
imply that every maximally-connected set of the topological space $\opair{\U}{\stopology{\topology{}}{\U}}$
is a connected set of $\Xt$:
\begin{equation}\label{thmlocallyconnectednessequiv1p2eq11}
\Foreach{\U}{\topology{}}
\maxcon{\opair{\U}{\stopology{\topology{}}{\U}}}\subseteq
\connecteds{\Xt}.
\end{equation}
\Ref{thmlocallyconnectednessequiv1p2eq2}
and
\Ref{thmlocallyconnectednessequiv1p2eq11}
imply that every element of $\base_{\rm{MC}}$ is a connected set of $\Xt$:
\begin{equation}\label{thmlocallyconnectednessequiv1p2eq12}
\base_{\rm{MC}}\subseteq\connecteds{\Xt}.
\end{equation}
Therefore according to \Ref{thmlocallyconnectednessequiv1p2eq8} and \Ref{thmlocallyconnectednessequiv1p2eq11},
it is clear that $\Xt$ has a base the elements of which are all connected.
\begin{equation}
\Existsis{\base_{\rm{MC}}}{\Cbase{\Xt}}
\base_{\rm{MC}}\subseteq\connecteds{\Xt}.
\end{equation}
According to \refdef{deflocallyconnectedness}, this means
$\Xt$ is a locally-connected topological space.
\endp
\end{itemize}
\endthm
%%%%%%%%%%%%%%%%%%%%%%%%%%%%%%%%%%%%%%%%%%%%%%%%%%%%%%%%%%%%%%%%%%%%%%%%%%%%%%%%%%%%%%%%%%%%%%%%%%%%%%%%%%%%%%%%%%%%%%%%%%%%%%%%
\definition\label{deflocallyconnectednessatapoint}
$\Xt=\opair{\X}{\topology{}}$
is taken as a topological space, and $\point$ as an element of $\X$.
It is said that $\quotl$the topological space $\Xt$ is locally-connected at the point $\point$ iff
for every neighborhood $\U$ of $\point$ in $\Xt$, there exists a neighborhood $\V$ of $\point$ in $\Xt$
which is a subset of $\U$ and a connectes set of $\Xt$.
In other words, it is said that $\Xt$ is locally-connected at $\point$ iff
\begin{equation}
\Foreach{\U}{\func{\nei{\Xt}}{\seta{\point}}}
\bigg(\Exists{\V}{\[\func{\nei{\Xt}}{\seta{\point}}\cap\connecteds{\Xt}\]}
\V\subseteq\U\bigg).
\end{equation}
\endef
%%%%%%%%%%%%%%%%%%%%%%%%%%%%%%%%%%%%%%%%%%%%%%%%%%%%%%%%%%%%%%%%%%%%%%%%%%%%%%%%%%%%%%%%%%%%%%%%%%%%%%%%%%%%%%%%%%%%%%%%%%%%%%%%
\theorem\label{thmlocallyconnectednessequiv2}
$\Xt=\opair{\X}{\topology{}}$
is taken as a topological space.
$\Xt$
is a locally-connected topological space if and only if $\Xt$ is locally-connected at every point $\point$ of $\Xt$.
\proof
\begin{itemize}
\item[${\textbf{\textsf{p1}}}$]
It is assumed that $\Xt$
is a locally-connected topological space. Then according to \refdef{deflocallyconnectedness},
\begin{equation}\label{thmlocallyconnectednessequiv2p1eq1}
\Existsis{\base}{\Cbase{\Xt}}
\base\subseteq\connecteds{\Xt}.
\end{equation}
\begin{itemize}
\item[${\textbf{\textsf{p1-1}}}$]
$\point$
is taken as an element of $\X$, and $\U$ as an element of $\func{\nei{\Xt}}{\point}$. Then according to \refdef{defnbdclassofsets},
\begin{align}
\U&\in\topology{},
\label{thmlocallyconnectednessequiv2p1-1eq1}\\
\point&\in\U,
\label{thmlocallyconnectednessequiv2p1-1eq12}
\end{align}
and hence according to \refthm{thmopensetpointsbase},
\begin{equation}\label{thmlocallyconnectednessequiv2p1-1eq3}
\Existsis{\B}{\base}\[\point\in\B,~\B\subseteq\U\].
\end{equation}
According to \Ref{thmlocllyconnectednessequiv2p1eq1}, \Ref{thmlocallyconnectednessequiv2p1-1eq3}, and
\refdef{defnbdclassofsets},
and considering that $\base\subseteq\topology{}$, it is clear that,
\begin{equation}
\Existsis{\B}{\[\func{\nei{\Xt}}{\seta{\point}}\cap\connecteds{\Xt}\]}
\B\subseteq\U.
\end{equation}
\endp
\end{itemize}
Therefore,
\begin{align*}
&\Foreach{\point}{\X}\cr
&\[\Foreach{\U}{\func{\nei{\Xt}}{\seta{\point}}}
\bigg(\Exists{\V}{\[\func{\nei{\Xt}}{\seta{\point}}\cap\connecteds{\Xt}\]}
\V\subseteq\U\bigg)\].\cr
&{}
\end{align*}
According to \refdef{deflocallyconnectednessatapoint},
this means $\Xt$ is locally-connected at every point.
\endp
\end{itemize}
\begin{itemize}
\item[${\textbf{\textsf{p2}}}$]
It is assumed that $\Xt$ is locally-connected at every point. That is,
\begin{align}
&\Foreach{\point}{\X}\cr
&\[\Foreach{\U}{\func{\nei{\Xt}}{\seta{\point}}}
\bigg(\Exists{\V}{\[\func{\nei{\Xt}}{\seta{\point}}\cap\connecteds{\Xt}\]}
\V\subseteq\U\bigg)\].\cr
&{}
\end{align}
Then according to \refthm{thmopensetpointsbase}, it is clear that,
\begin{equation*}
\Existsis{\base}{\Cbase{\Xt}}
\base\subseteq\connecteds{\Xt},
\end{equation*}
which means $\Xt$ is a locally-connected topological space.
\end{itemize}
\endthm
\chapteR{
Separation Axioms
}
\thispagestyle{fancy}
\section{
Separation Relations in a Topological Space
}
\definition\label{defpseudoCartesianproduct}
Each $\X$ and $\Y$ is taken as a set.
\begin{align}
\psCprod{\X}{\Y}&:=\defset{\opair{\x}{\y}}{\Cprod{\X}{\Y}}{\x\neq\y},\\
\iCprod{\X}{\Y}&:=\defset{\opair{\x}{\y}}{\Cprod{\X}{\Y}}{\x=\y}.
\end{align}
\endef
%%%%%%%%%%%%%%%%%%%%%%%%%%%%%%%%%%%%%%%%%%%%%%%%%%%%%%%%%%%%%%%%%%%%%%%%%%%%%%%%%%%%%%%%%%%%%%%%%%%%%%%%%%%%%%%%%%%
\theorem\label{thmpseudoCartesianproductscompl}
Each $\X$ and $\Y$ is taken as a set. The relations $\psCprod{\X}{\Y}$ and $\iCprod{\X}{\Y}$ in $\X$
are complements of each other:
\begin{equation}
\compl{\(\Cprod{\X}{\X}\)}{\(\iCprod{\X}{\Y}\)}=
\psCprod{\X}{\Y}.
\end{equation}
\prooff
It is trivial.
\endthm
%%%%%%%%%%%%%%%%%%%%%%%%%%%%%%%%%%%%%%%%%%%%%%%%%%%%%%%%%%%%%%%%%%%%%%%%%%%%%%%%%%%%%%%%%%%%%%%%%%%%%%%%%%%%%%%%%%%
\definition\label{defindistinguishablepoints}
$\Xt=\opair{\X}{\topology{}}$
is taken as a topological space, and $\opair{\point}{\x}$ as an element of $\Cprod{\X}{\X}$.
$\opair{\point}{\x}$ is referred to as a $\quotl$ pair of indistinguishable points of the topological space $\Xt$$\quotr$
iff the collection of all neighborhoods of $\seta{\point}$ in $\Xt$
equals the collection of all neighborhoods of $\seta{\x}$ in $\Xt$.\\
The set of all pairs of indistinguishable points of $\Xt$ will be denoted by $\InDist{\Xt}$:
\begin{equation}
\InDist{\Xt}:=\defset{\opair{\point}{\x}}{\Cprod{\X}{\X}}
{\func{\nei{\Xt}}{\seta{\point}}=\func{\nei{\Xt}}{\seta{\x}}}.
\end{equation}
\endef
%%%%%%%%%%%%%%%%%%%%%%%%%%%%%%%%%%%%%%%%%%%%%%%%%%%%%%%%%%%%%%%%%%%%%%%%%%%%%%%%%%%%%%%%%%%%%%%%%%%%%%%%%%%%%%%%%%%
\definition\label{defdistinguishablepoints}
$\Xt=\opair{\X}{\topology{}}$
is taken as a topological space, and $\opair{\point}{\x}$ as an element of $\Cprod{\X}{\X}$.
$\opair{\point}{\x}$
is referred to as a $\quotl$pair of partially-distinguishable points of the topological space $\Xt$$\quotr$
iff the collection of all neighborhoods of $\seta{\point}$ in $\Xt$
and the collection of all neighborhoods of $\seta{\x}$ in $\Xt$ are distinct.\\
The set of all pairs of partially-distinguishable points of $\Xt$ will be denoted by $\Dist{\Xt}$:
\begin{equation}
\Dist{\Xt}:=\defset{\opair{\point}{\x}}{\Cprod{\X}{\X}}
{\func{\nei{\Xt}}{\seta{\point}}\neq\func{\nei{\Xt}}{\seta{\x}}}.
\end{equation}
\endef
%%%%%%%%%%%%%%%%%%%%%%%%%%%%%%%%%%%%%%%%%%%%%%%%%%%%%%%%%%%%%%%%%%%%%%%%%%%%%%%%%%%%%%%%%%%%%%%%%%%%%%%%%%%%%%%%%%%
\definition\label{defseparatedpoints}
$\Xt=\opair{\X}{\topology{}}$
is taken as a topological space, and $\opair{\point}{\x}$
as an element of $\Cprod{\X}{\X}$.
$\opair{\point}{\x}$ is referred to as a $\quotl$pair of distinguishable points of the topological space $\Xt$$\quotr$
iff the collection of all neighborhoods of $\seta{\point}$ in $\Xt$
and the collection of all neighborhoods of $\seta{\x}$ in $\Xt$
do not have an inclusion relation.
The set of all distinguishable pairs of $\Xt$ will be denoted by $\Sep{\Xt}$:
\begin{align}
&\Sep{\Xt}:=\cr
&\defset{\opair{\point}{\x}}{\Cprod{\X}{\X}}
{\[\AND{\bigg(\func{\nei{\Xt}}{\seta{\point}}\nsubseteq\func{\nei{\Xt}}{\seta{\x}}\bigg)}
{\bigg(\func{\nei{\Xt}}{\seta{\x}}\nsubseteq\func{\nei{\Xt}}{\seta{\point}}\bigg)}\]}.\cr
&{}
\end{align}
\endef
%%%%%%%%%%%%%%%%%%%%%%%%%%%%%%%%%%%%%%%%%%%%%%%%%%%%%%%%%%%%%%%%%%%%%%%%%%%%%%%%%%%%%%%%%%%%%%%%%%%%%%%%%%%%%%%%%%%
\definition\label{defliberatedpoints}
$\Xt=\opair{\X}{\topology{}}$
is taken as a topological space, and $\opair{\point}{\x}$ an an element of $\Cprod{\X}{\X}$.
$\opair{\point}{\x}$
is referred to as a $\quotl$pair of separated points of the topological space $\Xt$$\quotr$
iff there exists at least one neighborhood $\U$ of $\seta{\point}$ in $\Xt$,
and a neighborhood $\V$ of $\seta{\x}$ in $\Xt$, that do not interset each other.\\
The set of all pairs of separated points of $\Xt$
will be denoted by $\Lib{\Xt}$:
\begin{align}
&\Lib{\Xt}:=\cr
&\defset{\opair{\point}{\x}}{\Cprod{\X}{\X}}
{\bigg[\Exists{\opair{\U}{\V}}{\func{\nei{\Xt}}{\seta{\point}}\times\func{\nei{\Xt}}{\seta{\x}}}
\(\U\cap\V=\empty\)\bigg]}.\cr
&{}
\end{align}
\endef
%%%%%%%%%%%%%%%%%%%%%%%%%%%%%%%%%%%%%%%%%%%%%%%%%%%%%%%%%%%%%%%%%%%%%%%%%%%%%%%%%%%%%%%%%%%%%%%%%%%%%%%%%%%%%%%%%%%
\theorem\label{thm(in)distinguishablepoints}
$\Xt=\opair{\X}{\topology{}}$
is taken as a topological space. For every $\opair{\point}{\x}$ in $\Cprod{\X}{\X}$,
$\opair{\point}{\x}$ is a pair of partially-distinguishable points of $\Xt$ iff $\opair{\point}{\x}$
is not a pair of indistinguishable points of $\Xt$. Inb other words, the relations $\InDist{\Xt}$
and $\Dist{\Xt}$ in $\X$ are complements of each other:
\begin{align}
\Dist{\Xt}=\compl{\(\Cprod{\X}{\X}\)}{\InDist{\Xt}}.
\end{align}
\prooff
According to \refdef{defindistinguishablepoints} and \refdef{defdistinguishablepoints},
it is clear.
\endthm
%%%%%%%%%%%%%%%%%%%%%%%%%%%%%%%%%%%%%%%%%%%%%%%%%%%%%%%%%%%%%%%%%%%%%%%%%%%%%%%%%%%%%%%%%%%%%%%%%%%%%%%%%%%%%%%%%%%
\theorem\label{thmidenticalpointsareindistinguishable}
$\Xt=\opair{\X}{\topology{}}$
is taken as a topological space. For every $\point$ in $\X$,
$\opair{\point}{\point}$
is a pair of indistinguishable points of $\Xt$. That is,
\begin{equation}
\iCprod{\X}{\X}\subseteq\InDist{\Xt}.
\end{equation}
\proof
According to \refdef{defindistinguishablepoints},
it is clear.
\endthm
%%%%%%%%%%%%%%%%%%%%%%%%%%%%%%%%%%%%%%%%%%%%%%%%%%%%%%%%%%%%%%%%%%%%%%%%%%%%%%%%%%%%%%%%%%%%%%%%%%%%%%%%%%%%%%%%%%%
\theorem\label{thmdistinguishablepointsaredifferent}
$\Xt=\opair{\X}{\topology{}}$
is taken as a topological space.
For every pair $\opair{\point}{\x}$ of partially-distinguished points of the topological space $\Xt$,
$\point$ and $\x$ are distinct.
\begin{equation}
\Dist{\Xt}\subseteq\(\psCprod{\X}{\X}\).
\end{equation}
\proof
According to \refdef{defdistinguishablepoints}, it is clear.
\endthm
%%%%%%%%%%%%%%%%%%%%%%%%%%%%%%%%%%%%%%%%%%%%%%%%%%%%%%%%%%%%%%%%%%%%%%%%%%%%%%%%%%%%%%%%%%%%%%%%%%%%%%%%%%%%%%%%%%%
\theorem\label{thmlibsepdis}
$\Xt=\opair{\X}{\topology{}}$
is taken as a topological space.
\begin{align}
\Lib{\Xt}\subseteq\Sep{\Xt}\subseteq\Dist{\Xt}.
\end{align}
\proof
According to
\refdef{defdistinguishablepoints},
\refdef{defseparatedpoints},
and
\refdef{defliberatedpoints},
it is clear.
\endthm
%%%%%%%%%%%%%%%%%%%%%%%%%%%%%%%%%%%%%%%%%%%%%%%%%%%%%%%%%%%%%%%%%%%%%%%%%%%%%%%%%%%%%%%%%%%%%%%%%%%%%%%%%%%%%%%%%%%
\theorem\label{thmseparatedpointsaredifferent}
$\Xt=\opair{\X}{\topology{}}$
is taken as a topological space.
For every pair $\opair{\point}{\x}$ of distinguishable points of $\Xt$,
$\point$ and $\x$ are distinct.
Moreover, for every pair $\opair{\point}{\x}$ of separated points of $\Xt$,
$\point$ and $\x$ are distinct. That is,
\begin{align}
\Sep{\Xt}&\subseteq\(\psCprod{\X}{\X}\),\\
\Lib{\Xt}&\subseteq\(\psCprod{\X}{\X}\).
\end{align}
\prooff
According to \refthm{thmdistinguishablepointsaredifferent} and \refthm{thmlibsepdis},
it is clear.
\endthm
%%%%%%%%%%%%%%%%%%%%%%%%%%%%%%%%%%%%%%%%%%%%%%%%%%%%%%%%%%%%%%%%%%%%%%%%%%%%%%%%%%%%%%%%%%%%%%%%%%%%%%%%%%%%%%%%%%%
\theorem\label{thmindistinguishablityequiv0}
$\Xt=\opair{\X}{\topology{}}$
is taken as a topological space.
For every $\opair{\point}{\x}$ in $\Cprod{\X}{\X}$,
$\opair{\point}{\x}$
is a pair of indistinguishable points of $\Xt$
if and only if every open set of $\Xt$ either contains $\point$ and $\x$, or contains non of them. That is,
\begin{align}
&\InDist{\Xt}=\cr
&\defset{\opair{\point}{\x}}{\Cprod{\X}{\X}}
{\bigg[\topology{}=\defset{\U}{\topology{}}{\seta{\binary{\point}{\x}}\subseteq\U}\cup
\defset{\U}{\topology{}}{\seta{\binary{\point}{\x}}\subseteq\(\compl{\X}{\U}\)}\bigg]}.\cr
&{}
\end{align}
\proof
According to \refdef{defnbdclassofsets},
\begin{align}\label{thmindistinguishablityequiv0peq1}
\Foreach{\opair{\point}{\x}}{\Cprod{\X}{\X}}\qquad
&\func{\nei{\Xt}}{\seta{\point}}\cap\func{\nei{\Xt}}{\seta{\x}}\cr
=&\defset{\U}{\topology{}}{\point\in\U}\cap
\defset{\U}{\topology{}}{\x\in\U}\cr
=&\defset{\U}{\topology{}}{\(\point\in\U,~\x\in\U\)}\cr
=&\defset{\U}{\topology{}}{\seta{\binary{\point}{\x}}\subseteq\U},
\end{align}
and
\begin{align}\label{thmindistinguishablityequiv0peq2}
\Foreach{\opair{\point}{\x}}{\Cprod{\X}{\X}}\qquad
&\compl{\topology{}}{\bigg(\func{\nei{\Xt}}{\seta{\point}}\cup\func{\nei{\Xt}}{\seta{\x}}\bigg)}\cr
=&\bigg(\compl{\topology{}}{\func{\nei{\Xt}}{\seta{\point}}}\bigg)
\cap\bigg(\compl{\topology{}}{\func{\nei{\Xt}}{\seta{\x}}}\bigg)\cr
=&\defset{\U}{\topology{}}{\point\notin\U}\cap
\defset{\U}{\topology{}}{\x\notin\U}\cr
=&\defset{\U}{\topology{}}{\point\in\(\compl{\X}{\U}\)}\cap
\defset{\U}{\topology{}}{\x\in\(\compl{\X}{\U}\)}\cr
=&\defset{\U}{\topology{}}{\[\point\in\(\compl{\X}{\U}\),~\x\in\(\compl{\X}{\U}\)\]}\cr
=&\defset{\U}{\topology{}}{\seta{\binary{\point}{\x}}\subseteq\(\compl{\X}{\U}\)}.
\end{align}
\begin{itemize}
\item[${\textbf{\textsf{p1}}}$]
$\opair{\point}{\x}$
is taken as an arbitrary element of $\InDist{\Xt}$. Then according to \refdef{defindistinguishablepoints},
\begin{equation}\label{thmindistinguishablityequiv0p1eq1}
\func{\nei{\Xt}}{\point}=\func{\nei{\Xt}}{\x}.
\end{equation}
and hence,
\begin{align}
\func{\nei{\Xt}}{\seta{\point}}&=
\func{\nei{\Xt}}{\seta{\point}}\cap\func{\nei{\Xt}}{\seta{\x}},
\label{thmindistinguishablityequiv0p1eq2}\\
\compl{\topology{}}{\func{\nei{\Xt}}{\seta{\point}}}&=
\compl{\topology{}}{\bigg(\func{\nei{\Xt}}{\seta{\point}}\cup\func{\nei{\Xt}}{\seta{\x}}\bigg)}.
\label{thmindistinguishablityequiv0p1eq3}
\end{align}
Moreover, it is clear that every open set of $\Xt$
is either a neighborhood of $\seta{\point}$, or not a neighborhood of $\seta{\point}$:
\begin{equation}\label{thmindistinguishablityequiv0p1eq4}
\topology{}=\func{\nei{\Xt}}{\seta{\point}}\cup
\(\compl{\topology{}}{\func{\nei{\Xt}}{\seta{\point}}}\).
\end{equation}
\Ref{thmindistinguishablityequiv0peq1},
\Ref{thmindistinguishablityequiv0peq2},
\Ref{thmindistinguishablityequiv0p1eq2},
\Ref{thmindistinguishablityequiv0p1eq3}
and
\Ref{thmindistinguishablityequiv0p1eq4}
imply that,
\begin{equation*}
\topology{}=\defset{\U}{\topology{}}{\seta{\binary{\point}{\x}}\subseteq\U}\cup
\defset{\U}{\topology{}}{\seta{\binary{\point}{\x}}\subseteq\(\compl{\X}{\U}\)}.
\end{equation*}
\endp
\end{itemize}
\begin{itemize}
\item[${\textbf{\textsf{p2}}}$]
$\opair{\point}{\x}$
is taken as such an arbitrary element of $\Cprod{\X}{\X}$ that,
\begin{equation}\label{thmindistinguishablityequiv0p2eq1}
\topology{}=\defset{\U}{\topology{}}{\seta{\binary{\point}{\x}}\subseteq\U}\cup
\defset{\U}{\topology{}}{\seta{\binary{\point}{\x}}\subseteq\(\compl{\X}{\U}\)}.
\end{equation}
Thus, according to \Ref{thmindistinguishablityequiv0peq1} and \Ref{thmindistinguishablityequiv0peq2},
\begin{align}
\topology{}=&\bigg[\func{\nei{\Xt}}{\seta{\point}}\cap\func{\nei{\Xt}}{\seta{\x}}\bigg]
\cup
\bigg[\compl{\topology{}}{\bigg(\func{\nei{\Xt}}{\seta{\point}}\cup\func{\nei{\Xt}}{\seta{\x}}\bigg)}\bigg]\cr
=&~~~\(\func{\nei{\Xt}}{\seta{\point}}\cup
\bigg[\compl{\topology{}}{\bigg(\func{\nei{\Xt}}{\seta{\point}}\cup\func{\nei{\Xt}}{\seta{\x}}\bigg)}\bigg]\)\cr
&\cap\(\func{\nei{\Xt}}{\seta{\x}}\cup
\bigg[\compl{\topology{}}{\bigg(\func{\nei{\Xt}}{\seta{\point}}\cup\func{\nei{\Xt}}{\seta{\x}}\bigg)}\bigg]\),
\end{align}
and thus,
\begin{align}
\(\func{\nei{\Xt}}{\seta{\point}}\cup
\bigg[\compl{\topology{}}{\bigg(\func{\nei{\Xt}}
{\seta{\point}}\cup\func{\nei{\Xt}}{\seta{\x}}\bigg)}\bigg]\)=\topology{},\\
\(\func{\nei{\Xt}}{\seta{\x}}\cup
\bigg[\compl{\topology{}}{\bigg(\func{\nei{\Xt}}
{\seta{\point}}\cup\func{\nei{\Xt}}{\seta{\x}}\bigg)}\bigg]\)=\topology{}.
\end{align}
Therefore,
\begin{align}
\func{\nei{\Xt}}{\seta{\x}}&\subseteq\func{\nei{\Xt}}{\seta{\point}},\\
\func{\nei{\Xt}}{\seta{\point}}&\subseteq\func{\nei{\Xt}}{\seta{\x}},
\end{align}
which means,
\begin{equation}
\func{\nei{\Xt}}{\seta{\point}}=\func{\nei{\Xt}}{\seta{\x}},
\end{equation}
and hence according to \refdef{defindistinguishablepoints},
\begin{equation}
\opair{\point}{\x}\in\InDist{\Xt}.
\end{equation}
\endp
\end{itemize}
\endthm
%%%%%%%%%%%%%%%%%%%%%%%%%%%%%%%%%%%%%%%%%%%%%%%%%%%%%%%%%%%%%%%%%%%%%%%%%%%%%%%%%%%%%%%%%%%%%%%%%%%%%%%%%%%%%%%%%%%
\theorem\label{thmindistinguishablityequiv1}
$\Xt=\opair{\X}{\topology{}}$
is taken as a topological space.
For every $\opair{\point}{\x}$ in $\Cprod{\X}{\X}$,
$\opair{\point}{\x}$ is a pair of indistinguishable points of $\Xt$
if and only if the collection of all closed sets of $\Xt$ containing $\point$
equals the collection of all closed sets of $\Xt$ containing $\x$:
\begin{align}
\InDist{\Xt}=\defset{\opair{\point}{\x}}{\Cprod{\X}{\X}}
{\[\func{\cnei{\Xt}}{\seta{\point}}=\func{\cnei{\Xt}}{\seta{\x}}\]}.
\end{align}
\prooff
According to \refdef{defcnbdclassofsets},
\begin{align}\label{thmindistinguishablityequiv1peq1}
\Foreach{\opair{\point}{\x}}{\Cprod{\X}{\X}}\qquad
&\func{\cnei{\Xt}}{\seta{\point}}\cap
\func{\cnei{\Xt}}{\seta{\x}}\cr
=&\defset{\V}{\Fclosed{\X}{\topology{}}}{\point\in\V}\cap
\defset{\V}{\Fclosed{\X}{\topology{}}}{\x\in\V}\cr
=&\defset{\V}{\Fclosed{\X}{\topology{}}}{\(\point\in\V,~\x\in\V\)}\cr
=&\defset{\V}{\Fclosed{\X}{\topology{}}}{\seta{\binary{\point}{\x}}\subseteq\V}
\end{align}
and
\begin{align}\label{thmindistinguishablityequiv1peq2}
\Foreach{\opair{\point}{\x}}{\Cprod{\X}{\X}}\qquad
&\compl{\Fclosed{\X}{\topology{}}}
{\bigg(\func{\cnei{\Xt}}{\seta{\point}}\cup\func{\cnei{\Xt}}{\seta{\x}}\bigg)}\cr
=&\bigg(\compl{\Fclosed{\X}{\topology{}}}
{\func{\cnei{\Xt}}{\seta{\point}}}\bigg)\cap
\bigg(\compl{\Fclosed{\X}{\topology{}}}
{\func{\cnei{\Xt}}{\seta{\x}}}\bigg)\cr
=&\defset{\V}{\Fclosed{\X}{\topology{}}}{\point\notin\V}\cap
\defset{\V}{\Fclosed{\X}{\topology{}}}{\x\notin\V}\cr
=&\defset{\V}{\Fclosed{\X}{\topology{}}}{\(\point\notin\V,~\x\notin\V\)}\cr
=&\defset{\V}{\Fclosed{\X}{\topology{}}}{\opair{\point}{\x}\subseteq\(\compl{\X}{\V}\)}.
\end{align}
\begin{itemize}
\item[${\textbf{\textsf{p1}}}$]
$\opair{\point}{\x}$
is taken as an arbitrary element of $\InDist{\Xt}$. Then according to \refdef{defindistinguishablepoints},
\begin{equation}\label{thmindistinguishablityequiv1p1eq1}
\topology{}=\defset{\U}{\topology{}}{\seta{\binary{\point}{\x}}\subseteq\U}\cup
\defset{\U}{\topology{}}{\seta{\binary{\point}{\x}}\subseteq\(\compl{\X}{\U}\)}.
\end{equation}
Thus, according to \refdef{deffamilyofclosedsets},
\begin{align}\label{thmindistinguishablityequiv1p1eq2}
\Fclosed{\X}{\topology{}}=
\defset{\V}{\Fclosed{\X}{\topology{}}}{\opair{\point}{\x}\subseteq\(\compl{\X}{\V}\)}\cup
\defset{\V}{\Fclosed{\X}{\topology{}}}{\opair{\point}{\x}\subseteq\V}.
\end{align}
\Ref{thmindistinguishablityequiv1peq1},
\Ref{thmindistinguishablityequiv1peq2},
and
\Ref{thmindistinguishablityequiv1p1eq2}
yield,
\begin{align}
\Fclosed{\X}{\topology{}}=&
\bigg[\func{\cnei{\Xt}}{\seta{\point}}\cap
\func{\cnei{\Xt}}{\seta{\x}}\bigg]\cup
\[\compl{\Fclosed{\X}{\topology{}}}
{\bigg(\func{\cnei{\Xt}}{\seta{\point}}\cup\func{\cnei{\Xt}}{\seta{\x}}\bigg)}\]\cr
=&~~~\(\func{\cnei{\Xt}}{\seta{\point}}\cup
\[\compl{\Fclosed{\X}{\topology{}}}
{\bigg(\func{\cnei{\Xt}}{\seta{\point}}\cup\func{\cnei{\Xt}}{\seta{\x}}\bigg)}\]\)\cr
&\cap\(\func{\cnei{\Xt}}{\seta{\x}}\cup
\[\compl{\Fclosed{\X}{\topology{}}}
{\bigg(\func{\cnei{\Xt}}{\seta{\point}}\cup\func{\cnei{\Xt}}{\seta{\x}}\bigg)}\]\),
\end{align}
and hence,
\begin{align}
\(\func{\cnei{\Xt}}{\seta{\point}}\cup
\bigg[\compl{\Fclosed{\X}{\topology{}}}{\bigg(\func{\cnei{\Xt}}
{\seta{\point}}\cup\func{\cnei{\Xt}}{\seta{\x}}\bigg)}\bigg]\)=\Fclosed{\X}{\topology{}},\\
\(\func{\nei{\Xt}}{\seta{\x}}\cup
\bigg[\compl{\Fclosed{\X}{\topology{}}}{\bigg(\func{\cnei{\Xt}}
{\seta{\point}}\cup\func{\cnei{\Xt}}{\seta{\x}}\bigg)}\bigg]\)=\Fclosed{\X}{\topology{}}.
\end{align}
Therefore,
\begin{align}
\func{\cnei{\Xt}}{\seta{\x}}&\subseteq\func{\cnei{\Xt}}{\seta{\point}},\\
\func{\cnei{\Xt}}{\seta{\point}}&\subseteq\func{\cnei{\Xt}}{\seta{\x}},
\end{align}
which means,
\begin{equation*}
\func{\cnei{\Xt}}{\seta{\point}}=\func{\cnei{\Xt}}{\seta{\x}},
\end{equation*}
\endp
\end{itemize}
\begin{itemize}
\item[${\textbf{\textsf{p2}}}$]
$\opair{\point}{\x}$
is taken as such an element of $\Cprod{\X}{\X}$ that
\begin{equation}\label{thmindistinguishablityequiv1p2eq1}
\func{\cnei{\Xt}}{\seta{\point}}=\func{\cnei{\Xt}}{\seta{\x}},
\end{equation}
Then,
\begin{align}
\func{\cnei{\Xt}}{\seta{\point}}&=
\func{\cnei{\Xt}}{\seta{\point}}\cap
\func{\cnei{\Xt}}{\seta{\x}}
\label{thmindistinguishablityequiv1p2eq2}\\
\compl{\Fclosed{\X}{\topology{}}}
{\func{\cnei{\Xt}}{\seta{\point}}}&=
\compl{\Fclosed{\X}{\topology{}}}
{\bigg(\func{\cnei{\Xt}}{\seta{\point}}\cup
\func{\cnei{\Xt}}{\seta{\x}}\bigg)},
\label{thmindistinguishablityequiv1p2eq3}
\end{align}
and hence according to \Ref{thmindistinguishablityequiv1peq1} and \Ref{thmindistinguishablityequiv1peq2},
\begin{align}
\func{\cnei{\Xt}}{\seta{\point}}&=
\defset{\V}{\Fclosed{\X}{\topology{}}}{\seta{\binary{\point}{\x}}\subseteq\V},
\label{thmindistinguishablityequiv1p2eq4}\\
\compl{\Fclosed{\X}{\topology{}}}
{\func{\cnei{\Xt}}{\seta{\point}}}&=
\defset{\V}{\Fclosed{\X}{\topology{}}}{\opair{\point}{\x}\subseteq\(\compl{\X}{\V}\)}.
\label{thmindistinguishablityequiv1p2eq5}
\end{align}
Moreover, it is clear that,
\begin{equation}\label{thmindistinguishablityequiv1p2eq6}
\Fclosed{\X}{\topology{}}=\func{\cnei{\Xt}}{\seta{\point}}\cup
\(\compl{\Fclosed{\X}{\topology{}}}{\func{\cnei{\Xt}}{\seta{\point}}}\).
\end{equation}
\Ref{thmindistinguishablityequiv1p2eq4},
\Ref{thmindistinguishablityequiv1p2eq5},
and
\Ref{thmindistinguishablityequiv1p2eq6}
yield,
\begin{equation}
\Fclosed{\X}{\topology{}}=
\defset{\V}{\Fclosed{\X}{\topology{}}}{\opair{\point}{\x}\subseteq\(\compl{\X}{\V}\)}\cup
\defset{\V}{\Fclosed{\X}{\topology{}}}{\opair{\point}{\x}\subseteq\V},
\end{equation}
and thus according to \refdef{deffamilyofclosedsets},
\begin{equation}
\topology{}=\defset{\U}{\topology{}}{\seta{\binary{\point}{\x}}\subseteq\U}\cup
\defset{\U}{\topology{}}{\seta{\binary{\point}{\x}}\subseteq\(\compl{\X}{\U}\)},
\end{equation}
and hence according to \refthm{thmindistinguishablityequiv0},
\begin{equation}
\opair{\point}{\x}\in\InDist{\Xt}.
\end{equation}
\endp
\end{itemize}
\endthm
%%%%%%%%%%%%%%%%%%%%%%%%%%%%%%%%%%%%%%%%%%%%%%%%%%%%%%%%%%%%%%%%%%%%%%%%%%%%%%%%%%%%%%%%%%%%%%%%%%%%%%%%%%%%%%%%%%%
\theorem\label{thmindistinguishablityequiv2}
$\Xt=\opair{\X}{\topology{}}$
is taken as a topological space.
\begin{align}
\InDist{\Xt}
&=\defset{\opair{\point}{\x}}{\Cprod{\X}{\X}}
{\intersection{\func{\nei{\Xt}}{\seta{\point}}}=
\intersection{\func{\nei{\Xt}}{\seta{\x}}}}.
\end{align}
\prooff
\begin{itemize}
\item[${\textbf{\textsf{p1}}}$]
According to \refdef{defindistinguishablepoints},
\begin{equation}
\Foreach{\opair{\point}{\x}}{\InDist{\Xt}}
\func{\nei{\Xt}}{\seta{\point}}=
\func{\nei{\Xt}}{\seta{\x}},
\end{equation}
and hence,
\begin{equation}
\Foreach{\opair{\point}{\x}}{\InDist{\Xt}}
\intersection{\func{\nei{\Xt}}{\seta{\point}}}=
\intersection{\func{\nei{\Xt}}{\seta{\x}}}.
\end{equation}
\endp
\end{itemize}
\begin{itemize}
\item[${\textbf{\textsf{p2}}}$]
$\opair{\point}{\x}$
is taken as such an arbitrary element of $\Cprod{\X}{\X}$ that,
\begin{equation}
\intersection{\func{\nei{\Xt}}{\seta{\point}}}=
\intersection{\func{\nei{\Xt}}{\seta{\x}}}.
\end{equation}
Then considering that,
\begin{equation}
\point\in\intersection{\func{\nei{\Xt}}{\seta{\point}}},
\end{equation}
it is clear that,
\begin{equation}
\point\in\intersection{\func{\nei{\Xt}}{\seta{\x}}},
\end{equation}
and hence,
\begin{equation}
\Foreach{\U}{\func{\nei{\Xt}}{\seta{\x}}}
\point\in\U.
\end{equation}
Therefore, according to \refdef{defnbdclassofsets},
\begin{equation}
\Foreach{\U}{\func{\nei{\Xt}}{\seta{\x}}}
\U\in\func{\nei{\Xt}}{\seta{\point}},
\end{equation}
which means,
\begin{equation}
\func{\nei{\Xt}}{\seta{\x}}\subseteq
\func{\nei{\Xt}}{\seta{\point}}.
\end{equation}
Similarly, it can be seen that,
\begin{equation}
\func{\nei{\Xt}}{\seta{\point}}\subseteq
\func{\nei{\Xt}}{\seta{\x}}.
\end{equation}
Therefore,
\begin{equation}
\func{\nei{\Xt}}{\seta{\point}}=
\func{\nei{\Xt}}{\seta{\x}},
\end{equation}
which according to \refdef{defindistinguishablepoints}, means,
\begin{equation}
\opair{\point}{\x}\in\InDist{\X}{\X}.
\end{equation}
\endp
\end{itemize}
\endthm
%%%%%%%%%%%%%%%%%%%%%%%%%%%%%%%%%%%%%%%%%%%%%%%%%%%%%%%%%%%%%%%%%%%%%%%%%%%%%%%%%%%%%%%%%%%%%%%%%%%%%%%%%%%%%%%%%%%
\theorem\label{thmindistinguishablityequiv3}
$\Xt=\opair{\X}{\topology{}}$
is taken as a topological space.
\begin{align}
\InDist{\Xt}
&=\defset{\opair{\point}{\x}}{\Cprod{\X}{\X}}
{\func{\Cl{\Xt}}{\seta{\point}}=\func{\Cl{\Xt}}{\seta{\x}}}.
\end{align}
\prooff
\begin{itemize}
\item[${\textbf{\textsf{p1}}}$]
According to \refthm{thmindistinguishablityequiv1},
\begin{equation}
\Foreach{\opair{\point}{\x}}{\InDist{\Xt}}
\func{\cnei{\Xt}}{\seta{\point}}=
\func{\cnei{\Xt}}{\seta{\x}},
\end{equation}
and hence according to \refdef{defcnbdclassofsets} and \refdef{defclosureofset},
\begin{align}
\Foreach{\opair{\point}{\x}}{\InDist{\Xt}}
\func{\Cl{\Xt}}{\seta{\point}}&=
\intersection{\func{\cnei{\Xt}}{\seta{\point}}}\cr
&=\intersection{\func{\cnei{\Xt}}{\seta{\x}}}\cr
&=\func{\Cl{\Xt}}{\seta{\x}}.
\end{align}
\endp
\end{itemize}
\begin{itemize}
\item[${\textbf{\textsf{p2}}}$]
$\opair{\point}{\x}$
is taken as such an arbitrary element of $\Cprod{\X}{\X}$ that
\begin{equation}
\func{\Cl{\Xt}}{\seta{\point}}=
\func{\Cl{\Xt}}{\seta{\x}}.
\end{equation}
Then considering that,
\begin{equation}
\point\in\func{\Cl{\Xt}}{\seta{\point}},
\end{equation}
and according to \refdef{defcnbdclassofsets} and \refdef{defclosureofset}, it is clear that,
\begin{equation}
\point\in\intersection{\func{\cnei{\Xt}}{\seta{\x}}},
\end{equation}
and hence,
\begin{equation}
\Foreach{\U}{\func{\cnei{\Xt}}{\seta{\x}}}
\point\in\U.
\end{equation}
Therefore, according to \refdef{defcnbdclassofsets},
\begin{equation}
\Foreach{\U}{\func{\cnei{\Xt}}{\seta{\x}}}
\U\in\func{\cnei{\Xt}}{\seta{\point}},
\end{equation}
which means,
\begin{equation}
\func{\cnei{\Xt}}{\seta{\x}}\subseteq
\func{\cnei{\Xt}}{\seta{\point}}.
\end{equation}
Similarly, it can be seen that,
\begin{equation}
\func{\cnei{\Xt}}{\seta{\point}}\subseteq
\func{\cnei{\Xt}}{\seta{\x}}.
\end{equation}
Therefore,
\begin{equation}
\func{\cnei{\Xt}}{\seta{\point}}=
\func{\cnei{\Xt}}{\seta{\x}},
\end{equation}
which according to \refthm{thmindistinguishablityequiv1}, means,
\begin{equation}
\opair{\point}{\x}\in\InDist{\X}{\X}.
\end{equation}
\endp
\end{itemize}
\endthm
%%%%%%%%%%%%%%%%%%%%%%%%%%%%%%%%%%%%%%%%%%%%%%%%%%%%%%%%%%%%%%%%%%%%%%%%%%%%%%%%%%%%%%%%%%%%%%%%%%%%%%%%%%%%%%%%%%%
\theorem\label{thmdistinguishablityequiv}
$\Xt=\opair{\X}{\topology{}}$
is taken as a topological space.
\begin{align}
\Dist{\Xt}&=
\defset{\opair{\point}{\x}}{\Cprod{\X}{\X}}
{\[\func{\cnei{\Xt}}{\seta{\point}}\neq
\func{\cnei{\Xt}}{\seta{\x}}\]}\cr
&=\defset{\opair{\point}{\x}}{\Cprod{\X}{\X}}
{\[\intersection{\func{\nei{\Xt}}{\seta{\point}}}\neq
\intersection{\func{\nei{\Xt}}{\seta{\x}}}\]}\cr
&=\defset{\opair{\point}{\x}}{\Cprod{\X}{\X}}
{\[\func{\Cl{\Xt}}{\seta{\point}}\neq
\func{\Cl{\Xt}}{\seta{\x}}\]}
\end{align}
\prooff
According to \refthm{thm(in)distinguishablepoints},
\refthm{thmindistinguishablityequiv1},
\refthm{thmindistinguishablityequiv2}, and \refthm{thmindistinguishablityequiv3}, it is obvious.
\endthm
%%%%%%%%%%%%%%%%%%%%%%%%%%%%%%%%%%%%%%%%%%%%%%%%%%%%%%%%%%%%%%%%%%%%%%%%%%%%%%%%%%%%%%%%%%%%%%%%%%%%%%%%%%%%%%%%%%%
%%%%%%%%%%%%%%%%%%%%%%%%%%%%%%%%%%%%%%%%%%%%%%%%%%%%%%%%%%%%%%%%%%%%%%%%%%%%%%%%%%%%%%%%%%%%%%%%%%%%%%%%%%%%%%%%%%%
%%%%%%%%%%%%%%%%%%%%%%%%%%%%%%%%%%%%%%%%%%%%%%%%%%%%%%%%%%%%%%%%%%%%%%%%%%%%%%%%%%%%%%%%%%%%%%%%%%%%%%%%%%%%%%%%%%%
%%%%%%%%%%%%%%%%%%%%%%%%%%%%%%%%%%%%%%%%%%%%%%%%%%%%%%%%%%%%%%%%%%%%%%%%%%%%%%%%%%%%%%%%%%%%%%%%%%%%%%%%%%%%%%%%%%%
\section{
$\sepA{0}$ (Kolmogorov) Axiom
}
\definition\label{defT0space}
$\Xt=\opair{\X}{\topology{}}$
is taken as a topological space.
$\Xt$
is referred to as a $\quotl$$\sepA{0}$ topological space$\quotr$, or a $\quotl$Kolmogorov space$\quotr$
iff
\begin{equation}
\InDist{\Xt}=\iCprod{\X}{\X}.
\end{equation}
In other words,
$\Xt$
is called a $\quotl$$\sepA{0}$ topological space$\quotr$
iff there exists no pair of distinct points of $\Xt$ that are indistinguishable.
\endef
%%%%%%%%%%%%%%%%%%%%%%%%%%%%%%%%%%%%%%%%%%%%%%%%%%%%%%%%%%%%%%%%%%%%%%%%%%%%%%%%%%%%%%%%%%%%%%%%%%%%%%%%%%%%%%%%%%%
\theorem\label{thmT0space1}
$\Xt=\opair{\X}{\topology{}}$
is taken as a topological space.
$\Xt$ is a $\sepA{0}$ topological space if and only if every pair of $\psCprod{\X}{\X}$
(every pair of distinct points of $\Xt$) is a pair of partially-distinguishable points of $\Xt$. That is,
$\Xt$
is a $\sepA{0}$ topological space if and only if
\begin{equation}
\Dist{\Xt}=\psCprod{\X}{\X}.
\end{equation}
\prooff
According to \refthm{thmpseudoCartesianproductscompl} and \refthm{thm(in)distinguishablepoints},
it is clear.
\endthm
%%%%%%%%%%%%%%%%%%%%%%%%%%%%%%%%%%%%%%%%%%%%%%%%%%%%%%%%%%%%%%%%%%%%%%%%%%%%%%%%%%%%%%%%%%%%%%%%%%%%%%%%%%%%%%%%%%%
\definition\label{defsetofT0topologies}
$\X$
is taken as a set.
The set of all topologies $\topology{}$ on $\X$ tsuch hat
$\opair{\X}{\topology{}}$ is a $\sepA{0}$ topological space, will be denoted by $\CtopsK{\X}$:
\begin{equation}
\CtopsK{\X}:=\defset{\topology{}}{\Ctops{\X}}
{\InDist{\opair{\X}{\topology{}}}=\iCprod{\X}{\X}}.
\end{equation}
\endef
%%%%%%%%%%%%%%%%%%%%%%%%%%%%%%%%%%%%%%%%%%%%%%%%%%%%%%%%%%%%%%%%%%%%%%%%%%%%%%%%%%%%%%%%%%%%%%%%%%%%%%%%%%%%%%%%%%%
%%%%%%%%%%%%%%%%%%%%%%%%%%%%%%%%%%%%%%%%%%%%%%%%%%%%%%%%%%%%%%%%%%%%%%%%%%%%%%%%%%%%%%%%%%%%%%%%%%%%%%%%%%%%%%%%%%%
\theorem\label{thmmultipointedindiscretespaceisnotT0}
$\X$
is taken as a set. If
$\CarD{\X}\geq 2$
($\X$ possesses at least two points),
then the indiscrete topological space $\opair{\X}{\seta{\binary{\empty}{\X}}}$ is not a $\sepA{0}$ topological space.
\proof
It is assumed that,
\begin{equation}\label{thmmultipointedindiscretespaceisnotT0peq1}
\CarD{\X}\geq 2.
\end{equation}
Then it is trivial that,
\begin{equation}\label{thmmultipointedindiscretespaceisnotT0peq2}
\psCprod{\X}{\X}\neq\empty.
\end{equation}
Moreover, it is clear that,
\begin{equation}\label{thmmultipointedindiscretespaceisnotT0peq3}
\Foreach{\point}{\X}
\func{\nei{\opair{\X}{\seta{\binary{\empty}{\X}}}}}{\seta{\point}}=\seta{\X},
\end{equation}
and hence according to \refdef{defdistinguishablepoints},
\begin{equation}\label{thmmultipointedindiscretespaceisnotT0peq4}
\Dist{\opair{\X}{\seta{\binary{\empty}{\X}}}}=\empty.
\end{equation}
\Ref{thmmultipointedindiscretespaceisnotT0peq2} and
\Ref{thmmultipointedindiscretespaceisnotT0peq4}
imply,
\begin{equation}
\Dist{\opair{\X}{\seta{\binary{\empty}{\X}}}}
\neq\psCprod{\X}{\X},
\end{equation}
which according to \refthm{thmT0space1}, means,
$\opair{\X}{\seta{\binary{\empty}{\X}}}$ is not a $\sepA{0}$ topological space.
\endthm
%%%%%%%%%%%%%%%%%%%%%%%%%%%%%%%%%%%%%%%%%%%%%%%%%%%%%%%%%%%%%%%%%%%%%%%%%%%%%%%%%%%%%%%%%%%%%%%%%%%%%%%%%%%%%%%%%%%
\theorem\label{thmT0space2}
$\Xt=\opair{\X}{\topology{}}$
is taken as a topological space.
$\Xt$ is a $\sepA{0}$ topological space if and only if,
\begin{equation}
\psCprod{\X}{\X}=\defset{\opair{\point}{\x}}{\Cprod{\X}{\X}}
{\[\func{\Cl{\Xt}}{\seta{\point}}\neq
\func{\Cl{\Xt}}{\seta{\x}}\]}.
\end{equation}
\proof
According to \refthm{thmdistinguishablityequiv} and \refthm{thmT0space1}, it is clear.
\endthm
%%%%%%%%%%%%%%%%%%%%%%%%%%%%%%%%%%%%%%%%%%%%%%%%%%%%%%%%%%%%%%%%%%%%%%%%%%%%%%%%%%%%%%%%%%%%%%%%%%%%%%%%%%%%%%%%%%%
%%%%%%%%%%%%%%%%%%%%%%%%%%%%%%%%%%%%%%%%%%%%%%%%%%%%%%%%%%%%%%%%%%%%%%%%%%%%%%%%%%%%%%%%%%%%%%%%%%%%%%%%%%%%%%%%%%%
%%%%%%%%%%%%%%%%%%%%%%%%%%%%%%%%%%%%%%%%%%%%%%%%%%%%%%%%%%%%%%%%%%%%%%%%%%%%%%%%%%%%%%%%%%%%%%%%%%%%%%%%%%%%%%%%%%%
%%%%%%%%%%%%%%%%%%%%%%%%%%%%%%%%%%%%%%%%%%%%%%%%%%%%%%%%%%%%%%%%%%%%%%%%%%%%%%%%%%%%%%%%%%%%%%%%%%%%%%%%%%%%%%%%%%%
%%%%%%%%%%%%%%%%%%%%%%%%%%%%%%%%%%%%%%%%%%%%%%%%%%%%%%%%%%%%%%%%%%%%%%%%%%%%%%%%%%%%%%%%%%%%%%%%%%%%%%%%%%%%%%%%%%%
\section{
$\sepA{1}$ Axiom
}
\definition\label{defT1space}
$\Xt=\opair{\X}{\topology{}}$
is taken as a topological space.
$\Xt$
is referred to as a $\quotl$$\sepA{1}$ topological space$\quotr$ iff
\begin{equation}
\Sep{\Xt}=\psCprod{\X}{\X}.
\end{equation}
\endef
%%%%%%%%%%%%%%%%%%%%%%%%%%%%%%%%%%%%%%%%%%%%%%%%%%%%%%%%%%%%%%%%%%%%%%%%%%%%%%%%%%%%%%%%%%%%%%%%%%%%%%%%%%%%%%%%%%%
\definition\label{defsetofT1topologies}
$\X$
is taken as a set.
The set of all topologies $\topology{}$ in $\X$ such that $\opair{\X}{\topology{}}$ is a $\sepA{1}$
topological space, will be denoted by $\CtopsF{\X}$:
\begin{equation}
\CtopsF{\X}:=\defset{\topology{}}{\Ctops{\X}}
{\Sep{\opair{\X}{\topology{}}}=\psCprod{\X}{\X}}.
\end{equation}
\endef
%%%%%%%%%%%%%%%%%%%%%%%%%%%%%%%%%%%%%%%%%%%%%%%%%%%%%%%%%%%%%%%%%%%%%%%%%%%%%%%%%%%%%%%%%%%%%%%%%%%%%%%%%%%%%%%%%%%
\theorem\label{thmsetofT1topologiesisnonemty}
$\X$
is taken as a set.
\begin{equation}
\CtopsF{\X}\neq\empty.
\end{equation}
\proof
According to \refthm{thmdiscreteTspaceisT1}  and \refdef{defsetofT1topologies},
\begin{equation}
\CSs{\X}\in\CtopsF{\X}.
\end{equation}
\endthm
%%%%%%%%%%%%%%%%%%%%%%%%%%%%%%%%%%%%%%%%%%%%%%%%%%%%%%%%%%%%%%%%%%%%%%%%%%%%%%%%%%%%%%%%%%%%%%%%%%%%%%%%%%%%%%%%%%%
\theorem\label{thmT1spaceequiv0}
$\Xt=\opair{\X}{\topology{}}$
is taken as a topological space.
$\Xt$ is a $\sepA{1}$ topological space if and only if for every $\opair{\point}{\x}$ in $\psCprod{\X}{\X}$,
there exists a neighborhood $\U$ of $\seta{\point}$ in $\Xt$, and a neighborhood $\V$ of $\seta{\x}$ in $\Xt$,
such that neither $\U$ contains $\x$, nor $\V$ contains $\point$. That is,
\begin{align}
\bigg(\topology{}\in\CtopsF{\X}\bigg)\thenn
\(\Foreach{\opair{\point}{\x}}{\psCprod{\X}{\X}}
\bigg[\Exists{\U}{\func{\nei{\Xt}}{\seta{\point}}}
\x\notin\U\bigg]\).
\end{align}
\proof
According to \refdef{defnbdclassofsets},
\refdef{defseparatedpoints},
and
\refdef{defT1space},
It is clear.
\endthm
%%%%%%%%%%%%%%%%%%%%%%%%%%%%%%%%%%%%%%%%%%%%%%%%%%%%%%%%%%%%%%%%%%%%%%%%%%%%%%%%%%%%%%%%%%%%%%%%%%%%%%%%%%%%%%%%%%%
\theorem\label{thmT1equivs}
$\Xt=\opair{\X}{\topology{}}$
is taken as a topological space.
The following propositions are equivalent.
\begin{itemize}
\item[$\[{\mathfrak{Eq}{\mathbf{0}}}\]$]
$\Xt$
is a $\sepA{1}$ topological space.
\item[$\[{\mathfrak{Eq}{\mathbf{1}}}\]$]
Every singleton subset of $\X$ is a closed set of $\Xt$:
\begin{equation}
\Foreach{\point}{\X}
\seta{\point}\in\Fclosed{\X}{\topology{}}.
\end{equation}
\item[$\[{\mathfrak{Eq}{\mathbf{2}}}\]$]
Every finite subset of $\X$ is a closed set of $\Xt$:
\begin{equation}
\defset{\asubset}{\X}{\CarD{\asubset}\in\Zp}\subseteq
\Fclosed{\X}{\topology{}}.
\end{equation}
\item[$\[{\mathfrak{Eq}{\mathbf{3}}}\]$]
Every subset of $\X$
equals the intersection of all of its neighborhoods in $\Xt$:
\begin{equation}
\Foreach{\asubset}{\CSs{\X}}
\asubset=\intersection{\func{\nei{\Xt}}{\asubset}}.
\end{equation}
\end{itemize}
\prooff\\
$\[{\mathfrak{Eq}{\mathbf{0}}}\leftrightarrow{\mathfrak{Eq}{\mathbf{1}}}\]$\\
It is clear that,
\begin{gather}
\Foreach{\point}{\X}
\(\Foreach{\x}{\(\compl{\X}{\seta{\point}}\)}
\bigg[\Exists{\V}{\func{\nei{\Xt}}{\seta{\x}}}
\point\notin\V\bigg]\)\cr
\vthenn\cr
\Foreach{\point}{\X}
\(\Foreach{\x}{\(\compl{\X}{\seta{\point}}\)}
\bigg[\Exists{\V}{\func{\nei{\Xt}}{\seta{\x}}}
\V\subseteq\(\compl{\X}{\seta{\point}}\)\bigg]\).
\end{gather}
According to \refdef{definteriorpoint} and \refthm{thmintofsetissetofintpoints},
\begin{gather}
\Foreach{\point}{\X}
\(\Foreach{\x}{\(\compl{\X}{\seta{\point}}\)}
\bigg[\Exists{\V}{\func{\nei{\Xt}}{\seta{\x}}}
\V\subseteq\(\compl{\X}{\seta{\point}}\)\bigg]\)\cr
\vthenn\cr
\Foreach{\point}{\X}
\bigg[\Foreach{\x}{\(\compl{\X}{\seta{\point}}\)}
\x\in\func{\Int{\Xt}}{\compl{\X}{\seta{\point}}}\bigg].
\end{gather}
According to \refcor{corintofset0},
\begin{gather}
\Foreach{\point}{\X}
\bigg[\Foreach{\x}{\(\compl{\X}{\seta{\point}}\)}
\x\in\func{\Int{\Xt}}{\compl{\X}{\seta{\point}}}\bigg]\cr
\vthenn\cr
\Foreach{\point}{\X}
\bigg[\func{\Int{\Xt}}{\compl{\X}{\seta{\point}}}=
\(\compl{\X}{\seta{\point}}\)\bigg].
\end{gather}
According to \refthm{thmintofopenset},
\begin{gather}
\Foreach{\point}{\X}
\bigg[\func{\Int{\Xt}}{\compl{\X}{\seta{\point}}}=
\(\compl{\X}{\seta{\point}}\)\bigg]\cr
\vthenn\cr
\Foreach{\point}{\X}
\bigg[
\(\compl{\X}{\seta{\point}}\)\in\topology{}\bigg].
\end{gather}
According to \refdef{deffamilyofclosedsets},
\begin{equation}
\(\Foreach{\point}{\X}
\bigg[
\(\compl{\X}{\seta{\point}}\)\in\topology{}\bigg]\)
\thenn
\bigg(\Foreach{\point}{\X}
\seta{\point}\in\Fclosed{\X}{\topology{}}\bigg).
\end{equation}
Therefore,
\begin{gather}
\Foreach{\point}{\X}
\(\Foreach{\x}{\(\compl{\X}{\seta{\point}}\)}
\bigg[\Exists{\V}{\func{\nei{\Xt}}{\seta{\x}}}
\point\notin\V\bigg]\)\cr
\vthenn\cr
\bigg(\Foreach{\point}{\X}
\seta{\point}\in\Fclosed{\X}{\topology{}}\bigg).
\end{gather}
According to \refthm{thmT1spaceequiv0},
this means $\Xt$ is a $\sepA{1}$ topological space if and only if,
\begin{equation*}
\Foreach{\point}{\X}
\seta{\point}\in\Fclosed{\X}{\topology{}}.
\end{equation*}
$\[{\mathfrak{Eq}{\mathbf{1}}}\leftrightarrow{\mathfrak{Eq}{\mathbf{2}}}\]$\\
According to \refthm{thmclosedsets}, it is clear that,
\begin{align}
\bigg(\Foreach{\point}{\X}
\seta{\point}\in\Fclosed{\X}{\topology{}}\bigg)
\thenn
\bigg(\defset{\asubset}{\X}{\CarD{\asubset}\in\Zp}\subseteq
\Fclosed{\X}{\topology{}}\bigg).
\end{align}
$\[{\mathfrak{Eq}{\mathbf{0}}}\leftrightarrow{\mathfrak{Eq}{\mathbf{3}}}\]$
\begin{itemize}
\item[${\textbf{\textsf{p1}}}$]
It is assumed that
$\Xt$ is a $\sepA{1}$ topological space.
Then according to \refthm{thmT1spaceequiv0},
\begin{equation}\label{thmT1equivs03p1eq1}
\Foreach{\opair{\point}{\x}}{\psCprod{\X}{\X}}
\bigg[\Exists{\U}{\func{\nei{\Xt}}{\seta{\point}}}
\x\notin\U\bigg].
\end{equation}
According to \Ref{thmT1equivs03p1eq1},
\begin{align}
\Foreach{\asubset}{\CSs{\X}}
\(\Foreach{\point}{\asubset}
\bigg[\Foreach{\x}{\(\compl{\X}{\asubset}\)}
\x\notin\intersection{\func{\nei{\Xt}}{\seta{\point}}}\bigg]\).
\end{align}
Thus considering that,
\begin{equation}
\Foreach{\asubset}{\CSs{\X}}
\asubset\subseteq\intersection{\func{\nei{\Xt}}{\asubset}},
\end{equation}
it is clear that,
\begin{equation*}
\Foreach{\asubset}{\CSs{\X}}
\intersection{\func{\nei{\Xt}}{\asubset}}=\asubset.
\end{equation*}
\endp
\end{itemize}
\begin{itemize}
\item[${\textbf{\textsf{p2}}}$]
It is assumed that,
\begin{equation}
\Foreach{\asubset}{\CSs{\X}}
\intersection{\func{\nei{\Xt}}{\asubset}}=\asubset.
\end{equation}
Then,
\begin{equation}
\Foreach{\point}{\X}
\intersection{\func{\nei{\Xt}}{\seta{\point}}}=
\seta{\point}.
\end{equation}
Therefore,
\begin{equation}
\Foreach{\opair{\point}{\x}}{\psCprod{\X}{\X}}
\x\notin\intersection{\func{\nei{\Xt}}{\seta{\point}}},
\end{equation}
which means,
\begin{equation}
\Foreach{\opair{\point}{\x}}{\psCprod{\X}{\X}}
\Exists{\U}{\func{\nei{\Xt}}{\seta{\point}}}
\x\notin\U,
\end{equation}
and hence according to \refthm{thmT1spaceequiv0},
$\Xt$ is a $\sepA{1}$ topological space.
\endp
\end{itemize}
\endthm
%%%%%%%%%%%%%%%%%%%%%%%%%%%%%%%%%%%%%%%%%%%%%%%%%%%%%%%%%%%%%%%%%%%%%%%%%%%%%%%%%%%%%%%%%%%%%%%%%%%%%%%%%%%%%%%%%%%
\theorem\label{thmT1ishereditary}
$\Xt=\opair{\X}{\topology{}}$
is taken as a topological space.
If $\Xt$ is a $\sepA{1}$ topological space, then for every $\asubset$ in $\CSs{\X}$,
$\opair{\asubset}{\stopology{\topology{}}{\asubset}}$
is a $\sepA{1}$ topological space:
\begin{equation}
\bigg(\topology{}\in\CtopsF{\X}\bigg)\then
\bigg[\Foreach{\asubset}{\CSs{\X}}
\stopology{\topology{}}{\asubset}\in\CtopsF{\asubset}\bigg].
\end{equation}
\proof
It is assumed that,
$\Xt$ is a $\sepA{1}$ topological space.
Then according to \refthm{thmT1spaceequiv0},
\begin{equation}\label{thmT1ishereditarypeq1}
\Foreach{\opair{\point}{\x}}{\psCprod{\X}{\X}}
\bigg[\Exists{\U}{\func{\nei{\Xt}}{\seta{\point}}}
\x\notin\U\bigg].
\end{equation}
\begin{itemize}
\item[${\textbf{\textsf{p1}}}$]
$\asubset$
is taken as an arbitrary subset of $\X$. According to \refdef{defnbdclassofsets} and
\refdef{defsubspacetopology1},
\begin{align}\label{thmT1ishereditaryp1eq1}
\Foreach{\point}{\asubset}
\func{\nei{\opair{\asubset}{\stopology{\topology{}}{\asubset}}}}{\seta{\point}}&=
\defset{\V}{\stopology{\topology{}}{\asubset}}{\point\in\V}\cr
&=\defset{\V}{\CSs{\asubset}}
{\bigg[\Exists{\U}{\func{\nei{\Xt}}{\seta{\point}}}\V=\asubset\cap\U\bigg]}.\cr
&{}
\end{align}
\begin{itemize}
\item[${\textbf{\textsf{p1-1}}}$]
$\opair{\point}{\x}$
is taken as an arbitrary element of $\psCprod{\asubset}{\asubset}$. According to \Ref{thmT1ishereditarypeq1},
\begin{equation}
\Existsis{\U}{\func{\nei{\Xt}}{\seta{\point}}}
\x\notin\U.
\end{equation}
Therefore, according to \Ref{thmT1ishereditaryp1eq1},
\begin{equation}
\Existsis{\(\asubset\cap\U\)}
{\func{\nei{\opair{\asubset}{\stopology{\topology{}}{\asubset}}}}{\seta{\point}}}
\x\not\in\(\asubset\cap\U\).
\end{equation}
\endp
\end{itemize}
Therefore,
\begin{align}
\Foreach{\opair{\point}{\x}}{\psCprod{\asubset}{\asubset}}
\bigg[\Exists{\(\V\)}
{\func{\nei{\opair{\asubset}{\stopology{\topology{}}{\asubset}}}}{\seta{\point}}}
\x\not\in\V\bigg],
\end{align}
which according to \refthm{thmT1spaceequiv0}, means,
$\opair{\asubset}{\stopology{\topology{}}{\asubset}}$ is a $\sepA{1}$ topological space.
\endp
\end{itemize}
\endthm
%%%%%%%%%%%%%%%%%%%%%%%%%%%%%%%%%%%%%%%%%%%%%%%%%%%%%%%%%%%%%%%%%%%%%%%%%%%%%%%%%%%%%%%%%%%%%%%%%%%%%%%%%%%%%%%%%%%
\theorem\label{thmdiscreteTspaceisT1}
$\X$
is taken as a set.
The discrete topological space $\opair{\X}{\CSs{\X}}$
is a $\sepA{1}$ topological space.
\proof
Considering that,
\begin{equation}
\Foreach{\point}{\X}
\seta{\point}\in\func{\nei{\opair{\X}{\CSs{\X}}}}{\seta{\point}},
\end{equation}
and
\begin{equation}
\Foreach{\opair{\point}{\x}}{\psCprod{\X}{\X}}
\seta{\point}\notin\func{\nei{\opair{\X}{\CSs{\X}}}}{\seta{\x}},
\end{equation}
it is clear that,
\begin{align}
&\Foreach{\opair{\point}{\x}}{\psCprod{\X}{\X}}\cr
&\bigg[
\func{\nei{\opair{\X}{\CSs{\X}}}}{\seta{\point}}\nsubseteq
\func{\nei{\opair{\X}{\CSs{\X}}}}{\seta{\x}},~
\func{\nei{\opair{\X}{\CSs{\X}}}}{\seta{\x}}\nsubseteq
\func{\nei{\opair{\X}{\CSs{\X}}}}{\seta{\point}}
\bigg],\cr
&{}
\end{align}
and hence according to \refdef{defseparatedpoints} and \refthm{thmseparatedpointsaredifferent},
\begin{equation}
\Sep{\opair{\X}{\CSs{\X}}}=\psCprod{\X}{\X},
\end{equation}
which according to \refdef{defT1space}, means
$\opair{\X}{\CSs{\X}}$ is a $\sepA{1}$ topological space.
\endthm
%%%%%%%%%%%%%%%%%%%%%%%%%%%%%%%%%%%%%%%%%%%%%%%%%%%%%%%%%%%%%%%%%%%%%%%%%%%%%%%%%%%%%%%%%%%%%%%%%%%%%%%%%%%%%%%%%%%
\theorem\label{thmsingletonTspaceisT1}
$\x$
is taken as a set.
The singleton topological space $\singletonTS{\x}$
is a $\sepA{1}$ topological space.
\proof
Considering that
$\singletonTS{\x}$ is a discrete topological space,
(\refthm{thmsingletontopology0}),
according to \refthm{thmdiscreteTspaceisT1}, it is clear that $\singletonTS{\x}$ is a $\sepA{1}$ topological space.
\endthm
%%%%%%%%%%%%%%%%%%%%%%%%%%%%%%%%%%%%%%%%%%%%%%%%%%%%%%%%%%%%%%%%%%%%%%%%%%%%%%%%%%%%%%%%%%%%%%%%%%%%%%%%%%%%%%%%%%%
\theorem\label{thmcoarsestT1topologyofset}
$\X$ is taken as a set.
\begin{equation}
\[\intersection{\CtopsF{\X}}\]\in\CtopsF{\X}.
\end{equation}
This means, the partially-ordered set $\opair{\CtopsF{\X}}{\subseteq}$ has a minimum.
In other words, there is a topologyy $\topology{}$ on $\X$ such that
$\opair{\X}{\topology{}}$ is a $\sepA{1}$ topological space, and there is no topology on $\X$ coarser
than $\topology{}$ that turns $\X$ into a $\sepA{1}$ topological space.
\proof
According to \refdef{defsetofT1topologies} and
\refthm{thmsetofT1topologiesisnonemty},
$\CtopsF{\X}$
is a non-empty subset of $\Ctops{\X}$:
\begin{equation}
\CtopsF{\X}\in\[\compl{\CSs{\Ctops{\X}}}{\seta{\empty}}\],
\end{equation}
and hence according to \refthm{thmintersectionoftopologies},
\begin{equation}
\[\intersection{\CtopsF{\X}}\]\in\Ctops{\X}.
\end{equation}
\begin{itemize}
\item[${\textbf{\textsf{p1}}}$]
$\point$
is taken as an arbitrary element of $\X$. According to \refthm{thmT1equivs}, for every $\topology{}$ in $\CtopsF{\X}$,
$\seta{\point}$ is a closed set of the topological space $\opair{\X}{\topology{}}$:
\begin{equation}
\Foreach{\topology{}}{\CtopsF{\X}}
\seta{\point}\in\Fclosed{\X}{\topology{}},
\end{equation}
which according to \refdef{deffamilyofclosedsets}, means,
\begin{equation}
\Foreach{\topology{}}{\CtopsF{\X}}
\(\compl{\X}{\seta{\point}}\)\in\topology{}.
\end{equation}
Therefore,
\begin{equation}
\(\compl{\X}{\seta{\point}}\)\in\[\intersection{\CtopsF{\X}}\],
\end{equation}
which according to \refdef{deffamilyofclosedsets}, means $\seta{\point}$ is a closed set of the
topological space $\opair{\X}{\intersection{\CtopsF{\X}}}$:
\begin{equation}
\seta{\point}\in\Fclosed{\X}{\intersection{\CtopsF{\X}}}.
\end{equation}
\endp
\end{itemize}
Therefore it is clear that every singleton subset of $\X$
is a closed set of the topological space $\opair{\X}{\intersection{\CtopsF{\X}}}$:
\begin{equation}
\Foreach{\point}{\X}
\seta{\point}\in\Fclosed{\X}{\intersection{\CtopsF{\X}}},
\end{equation}
and hence according to \refthm{thmT1equivs},
$\opair{\X}{\intersection{\CtopsF{\X}}}$ is a $\sepA{1}$ topological space.
\endthm
%%%%%%%%%%%%%%%%%%%%%%%%%%%%%%%%%%%%%%%%%%%%%%%%%%%%%%%%%%%%%%%%%%%%%%%%%%%%%%%%%%%%%%%%%%%%%%%%%%%%%%%%%%%%%%%%%%%
\theorem\label{thmcontinuousmapfromindiscretespacetoT1space}
$\X$ is taken as a non-empty set, and $\Yt=\opair{\Y}{\topology{\Y}}$ as a $\sepA{1}$
topological space.
Every $\cf$ in $\CF{\opair{\X}{\seta{\binary{\empty}{\X}}}}{\Yt}$
(every continuous map from the indiscrete topological space $\opair{\X}{\seta{\binary{\empty}{\X}}}$
to the topological space $\Yt$) is a constant map. That is,
\begin{align}
\Foreach{\cf}{\CF{\opair{\X}{\seta{\binary{\empty}{\X}}}}{\Yt}}
\bigg[\Exists{\y_{0}}{\Y}\func{\image{\cf}}{\X}=\seta{\y_{0}}\bigg].
\end{align}
\proof
$\cf$ is taken as an arbitrary element of $\CF{\opair{\X}{\seta{\binary{\empty}{\X}}}}{\Yt}$.
Then according to \refthm{thmcontiniuityequiv1},
the pre-image under $\cf$ of every closed set of $\Yt$
is a closed set of $\opair{\X}{\seta{\binary{\empty}{\X}}}$:
\begin{equation}\label{thmcontinuousmapfromindiscretespacetoT1spacepeq1}
\Foreach{\U}{\Fclosed{\Y}{\topology{\Y}}}
\func{\pimage{\cf}}{\U}\in\Fclosed{\X}{\seta{\binary{\empty}{\X}}}.
\end{equation}
Therefore, considering that
\begin{equation}\label{thmcontinuousmapfromindiscretespacetoT1spacepeq2}
\Fclosed{\X}{\seta{\binary{\empty}{\X}}}=\seta{\binary{\empty}{\X}},
\end{equation}
it is clear that,
\begin{equation}\label{thmcontinuousmapfromindiscretespacetoT1spacepeq3}
\Foreach{\U}{\Fclosed{\Y}{\topology{\Y}}}
\func{\pimage{\cf}}{\U}\in\seta{\binary{\empty}{\X}}.
\end{equation}
Moreover, considering that $\Yt$ is a $\sepA{1}$ topological space, according to \refthm{thmT1equivs},
every singleton subset of $\Y$ is a closed set of $\Yt$:
\begin{equation}\label{thmcontinuousmapfromindiscretespacetoT1spacepeq4}
\Foreach{\y}{\Y}
\seta{\y}\in\Fclosed{\Y}{\topology{\Y}}.
\end{equation}
\Ref{thmcontinuousmapfromindiscretespacetoT1spacepeq3} and \Ref{thmcontinuousmapfromindiscretespacetoT1spacepeq4}
imply that the pre-image under $\cf$ of every singleton subset of $\Y$ must be either $\empty$ or $\X$:
\begin{equation}
\Foreach{\y}{\Y}
\func{\pimage{\cf}}{\seta{\y}}\in\seta{\binary{\empty}{\X}}.
\end{equation}
There it is clear that there is an element $\y_{0}$ of $\Y$ such that the pre-image under $\cf$ of $\seta{\y_{0}}$
equals $\X$:
\begin{equation}
\Exists{\y_{0}}{\Y}
\func{\pimage{\cf}}{\seta{\y_{0}}}=\X,
\end{equation}
which means,
\begin{equation}
\Exists{\y_{0}}{\Y}
\func{\image{\cf}}{\X}=\seta{\y_{0}}.
\end{equation}
\endthm
%%%%%%%%%%%%%%%%%%%%%%%%%%%%%%%%%%%%%%%%%%%%%%%%%%%%%%%%%%%%%%%%%%%%%%%%%%%%%%%%%%%%%%%%%%%%%%%%%%%%%%%%%%%%%%%%%%%%%
\theorem\label{thminjectivecontinuousfunctiontoT1space}
$\Xt=\opair{\X}{\topology{\X}}$
is taken as a topological space, and $\Yt=\opair{\Y}{\topology{\Y}}$ as a $\sepA{1}$ topological space.
If there exists an injective continuous map from $\Xt$ to $\Yt$, then
$\Xt$ is also a $\sepA{1}$ topological space.
\begin{equation}
\[\CF{\Xt}{\Yt}\cap\InF{\X}{\Y}\]\neq\empty
\then
\topology{\X}\in\CtopsF{\X}.
\end{equation}
\prooff
It is assumed that,
\begin{equation}\label{thminjectivecontinuousfunctiontoT1spacepeq1}
\CF{\Xt}{\Yt}\cap\InF{\X}{\Y}\neq\empty.
\end{equation}
Then,
\begin{equation}\label{thminjectivecontinuousfunctiontoT1spacepeq2}
\Existsis{\cf}{\InF{\X}{\Y}}
\cf\in\CF{\Xt}{\Yt}.
\end{equation}
Considering that,
$\Yt$ is a $\sepA{1}$ topological space, according to \refthm{thmT1equivs},
\begin{equation}\label{thminjectivecontinuousfunctiontoT1spacepeq3}
\Foreach{\y}{\Y}
\seta{\y}\in\Fclosed{\Y}{\topology{\Y}}.
\end{equation}
Moreover, considering that $\cf$ is a continuous map from $\Xt$ to $\Yt$, according to \refthm{thmcontiniuityequiv1},
\begin{equation}\label{thminjectivecontinuousfunctiontoT1spacepeq4}
\Foreach{\U}{\Fclosed{\Y}{\topology{\Y}}}
\func{\pimage{\cf}}{\U}\in\Fclosed{\X}{\topology{\X}}.
\end{equation}
\Ref{thminjectivecontinuousfunctiontoT1spacepeq3} and
\Ref{thminjectivecontinuousfunctiontoT1spacepeq4}
imply that,
\begin{equation}\label{thminjectivecontinuousfunctiontoT1spacepeq5}
\Foreach{\y}{\Y}
\func{\pimage{\cf}}{\seta{\y}}\in\Fclosed{\X}{\topology{\X}}.
\end{equation}
Moreover, considering the injectivity of $\cf$,
\begin{equation}\label{thminjectivecontinuousfunctiontoT1spacepeq6}
\Foreach{\x}{\X}
\func{\pimage{\cf}}{\seta{\func{\cf}{\x}}}=\seta{\x}.
\end{equation}
\Ref{thminjectivecontinuousfunctiontoT1spacepeq5} and
\Ref{thminjectivecontinuousfunctiontoT1spacepeq6}
yiekd,
\begin{equation}
\Foreach{\x}{\X}
\seta{\x}\in\Fclosed{\X}{\topology{\X}},
\end{equation}
which according to \refthm{thmT1equivs}, means,
\begin{equation}
\topology{\X}\in\CtopsF{\X}.
\end{equation}
\endthm
%%%%%%%%%%%%%%%%%%%%%%%%%%%%%%%%%%%%%%%%%%%%%%%%%%%%%%%%%%%%%%%%%%%%%%%%%%%%%%%%%%%%%%%%%%%%%%%%%%%%%%%%%%%%%%%%%
\theorem\label{thmT1isstrongerthanT0}
$\X$
is taken as a set.
\begin{equation}
\CtopsK{\X}\supseteq\CtopsF{\X}.
\end{equation}
\proof
According to
\refthm{thmdistinguishablepointsaredifferent},
\refthm{thmlibsepdis},
\refthm{thmT0space1},
\refdef{defsetofT0topologies},
\refdef{defT1space}, and
\refdef{defsetofT1topologies},
it is clear.
\endthm
%%%%%%%%%%%%%%%%%%%%%%%%%%%%%%%%%%%%%%%%%%%%%%%%%%%%%%%%%%%%%%%%%%%%%%%%%%%%%%%%%%%%%%%%%%%%%%%%%%%%%%%%%%%%%%%%%%%%%%%%%%%%%%%%%%%%%%%%%%%%%%%%%%%%%%
%%%%%%%%%%%%%%%%%%%%%%%%%%%%%%%%%%%%%%%%%%%%%%%%%%%%%%%%%%%%%%%%%%%%%%%%%%%%%%%%%%%%%%%%%%%%%%%%%%%%%%%%%%%%%%%%%%%%%%%%%%%%%%%%%%%%%%%%%%%%%%%%%%%%%%
%%%%%%%%%%%%%%%%%%%%%%%%%%%%%%%%%%%%%%%%%%%%%%%%%%%%%%%%%%%%%%%%%%%%%%%%%%%%%%%%%%%%%%%%%%%%%%%%%%%%%%%%%%%%%%%%%%%%%%%%%%%%%%%%%%%%%%%%%%%%%%%%%%%%%%
%%%%%%%%%%%%%%%%%%%%%%%%%%%%%%%%%%%%%%%%%%%%%%%%%%%%%%%%%%%%%%%%%%%%%%%%%%%%%%%%%%%%%%%%%%%%%%%%%%%%%%%%%%%%%%%%%%%%%%%%%%%%%%%%%%%%%%%%%%%%%%%%%%%%%%
%%%%%%%%%%%%%%%%%%%%%%%%%%%%%%%%%%%%%%%%%%%%%%%%%%%%%%%%%%%%%%%%%%%%%%%%%%%%%%%%%%%%%%%%%%%%%%%%%%%%%%%%%%%%%%%%%%%%%%%%%%%%%%%%%%%%%%%%%%%%%%%%%%%%%%
%%%%%%%%%%%%%%%%%%%%%%%%%%%%%%%%%%%%%%%%%%%%%%%%%%%%%%%%%%%%%%%%%%%%%%%%%%%%%%%%%%%%%%%%%%%%%%%%%%%%%%%%%%%%%%%%%%%%%%%%%%%%%%%%%%%%%%%%%%%%%%%%%%%%%%
%%%%%%%%%%%%%%%%%%%%%%%%%%%%%%%%%%%%%%%%%%%%%%%%%%%%%%%%%%%%%%%%%%%%%%%%%%%%%%%%%%%%%%%%%%%%%%%%%%%%%%%%%%%%%%%%%%%%%%%%%%%%%%%%%%%%%%%%%%%%%%%%%%%%%%
%%%%%%%%%%%%%%%%%%%%%%%%%%%%%%%%%%%%%%%%%%%%%%%%%%%%%%%%%%%%%%%%%%%%%%%%%%%%%%%%%%%%%%%%%%%%%%%%%%%%%%%%%%%%%%%%%%%%%%%%%%%%%%%%%%%%%%%%%%%%%%%%%%%%%%
%%%%%%%%%%%%%%%%%%%%%%%%%%%%%%%%%%%%%%%%%%%%%%%%%%%%%%%%%%%%%%%%%%%%%%%%%%%%%%%%%%%%%%%%%%%%%%%%%%%%%%%%%%%%%%%%%%%%%%%%%%%%%%%%%%%%%%%%%%%%%%%%%%%%%%
\section{
$\sepA{2}$ (Hausdorff) Axiom
}
\definition\label{defT2space}
$\Xt=\opair{\X}{\topology{}}$
is taken as a topological space.
$\Xt$ is referred to as a $\quotl$$\sepA{2}$ topological space$\quotr$, or a
$\quotl$Hausdorff space$\quotr$ iff
\begin{equation}
\Lib{\Xt}=\psCprod{\X}{\X}.
\end{equation}
\endef
%%%%%%%%%%%%%%%%%%%%%%%%%%%%%%%%%%%%%%%%%%%%%%%%%%%%%%%%%%%%%%%%%%%%%%%%%%%%%%%%%%%%%%%%%%%%%%%%%%%%%%%%%%%%%%%%%%%
\definition\label{defsetofT2topologies}
$\X$
is taken as a set.
The set of all topologies $\topology{}$ on $\X$ such that $\opair{\X}{\topology{}}$ is a $\sepA{2}$
topological space, will be denoted by $\CtopsH{\X}$:
\begin{equation}
\CtopsH{\X}:=\defset{\topology{}}{\Ctops{\X}}
{\Lib{\opair{\X}{\topology{}}}=\psCprod{\X}{\X}}.
\end{equation}
\endef
%%%%%%%%%%%%%%%%%%%%%%%%%%%%%%%%%%%%%%%%%%%%%%%%%%%%%%%%%%%%%%%%%%%%%%%%%%%%%%%%%%%%%%%%%%%%%%%%%%%%%%%%%%%%%%%%%%%
\theorem\label{thmHausdorffequiv1}
$\Xt=\opair{\X}{\topology{}}$
is taken as a topological space.
$\Xt$ is a $\sepA{2}$ topological space if and only if every singleton subset of $\X$
equals the intersection of all closed neighborhoods of it in $\Xt$:
\begin{equation}
\bigg[\topology{}\in\CtopsH{\X}\bigg]
\thenn
\bigg[\Foreach{\point}{\X}
\seta{\point}=\intersection{\func{\cnei{\Xt}}{\seta{\point}}}\bigg].
\end{equation}
%%%%%%%%%%%%%%%%%%%%%%%%%%%%%%%%%%%%%%%%%%%%%%%%%%%%%%%%%%%%%%%%%%%%%%%%%%%%%%%%%%%%%%%%%%%%%%%%%%%%%%%%%%%%%%%%%%%
%%%%%%%%%%%%%%%%%%%%%%%%%%%%%%%%%%%%%%%%%%%%%%%%%%%%%%%%%%%%%%%%%%%%%%%%%%%%%%%%%%%%%%%%%%%%%%%%%%%%%%%%%%%%%%%%%%%%%%%%%%%%%%%%%%%%%%%%%%%%%%%%%%%%%%
%%%%%%%%%%%%%%%%%%%%%%%%%%%%%%%%%%%%%%%%%%%%%%%%%%%%%%%%%%%%%%%%%%%%%%%%%%%%%%%%%%%%%%%%%%%%%%%%%%%%%%%%%%%%%%%%%%%%%%%%%%%%%%%%%%%%%%%%%%%%%%%%%%%%%%
%%%%%%%%%%%%%%%%%%%%%%%%%%%%%%%%%%%%%%%%%%%%%%%%%%%%%%%%%%%%%%%%%%%%%%%%%%%%%%%%%%%%%%%%%%%%%%%%%%%%%%%%%%%%%%%%%%%%%%%%%%%%%%%%%%%%%%%%%%%%%%%%%%%%%%
%%%%%%%%%%%%%%%%%%%%%%%%%%%%%%%%%%%%%%%%%%%%%%%%%%%%%%%%%%%%%%%%%%%%%%%%%%%%%%%%%%%%%%%%%%%%%%%%%%%%%%%%%%%%%%%%%%%%%%%%%%%%%%%%%%%%%%%%%%%%%%%%%%%%%%
%%%%%%%%%%%%%%%%%%%%%%%%%%%%%%%%%%%%%%%%%%%%%%%%%%%%%%%%%%%%%%%%%%%%%%%%%%%%%%%%%%%%%%%%%%%%%%%%%%%%%%%%%%%%%%%%%%%%%%%%%%%%%%%%%%%%%%%%%%%%%%%%%%%%%%
%%%%%%%%%%%%%%%%%%%%%%%%%%%%%%%%%%%%%%%%%%%%%%%%%%%%%%%%%%%%%%%%%%%%%%%%%%%%%%%%%%%%%%%%%%%%%%%%%%%%%%%%%%%%%%%%%%%%%%%%%%%%%%%%%%%%%%%%%%%%%%%%%%%%%%
%%%%%%%%%%%%%%%%%%%%%%%%%%%%%%%%%%%%%%%%%%%%%%%%%%%%%%%%%%%%%%%%%%%%%%%%%%%%%%%%%%%%%%%%%%%%%%%%%%%%%%%%%%%%%%%%%%%%%%%%%%%%%%%%%%%%%%%%%%%%%%%%%%%%%%
%%%%%%%%%%%%%%%%%%%%%%%%%%%%%%%%%%%%%%%%%%%%%%%%%%%%%%%%%%%%%%%%%%%%%%%%%%%%%%%%%%%%%%%%%%%%%%%%%%%%%%%%%%%%%%%%%%%%%%%%%%%%%%%%%%%%%%%%%%%%%%%%%%%%%%
%%%%%%%%%%%%%%%%%%%%%%%%%%%%%%%%%%%%%%%%%%%%%%%%%%%%%%%%%%%%%%%%%%%%%%%%%%%%%%%%%%%%%%%%%%%%%%%%%%%%%%%%%%%%%%%%%%%%%%%%%%%%%%%%%%%%%%%%%%%%%%%%%%%%%%
\section{
$\sepA{3}$ Axiom, and Regular Spaces
}
\definition\label{defT3space}
$\Xt=\opair{\X}{\topology{}}$
is taken as a topological space.
$\Xt$
is referred to as a $\sepA{3}$ topological space$\quotr$
iff every closed set of $\Xt$ and every point of $\Xt$ in the complement of
that closed set possess at least one neighborhood in $\Xt$ each, that do not intersect each other.
That is, $\Xt$ is called a $\sepA{3}$ topological space iff
\begin{equation}
\Foreach{\asubset}{\Fclosed{\X}{\topology{}}}\Foreach{\point}{\(\compl{\X}{\asubset}\)}
\bigg[\Exists{\opair{\U}{\V}}{\bigg(\Cprod{\func{\nei{\Xt}}{\asubset}}{\func{\nei{\Xt}}{\seta{\point}}}\bigg)}
\U\cap\V=\empty\bigg].
\end{equation}
\endef
%%%%%%%%%%%%%%%%%%%%%%%%%%%%%%%%%%%%%%%%%%%%%%%%%%%%%%%%%%%%%%%%%%%%%%%%%%%%%%%%%%%%%%%%%%%%%%%%%%%%%%%%%%%%%%%%%%%
\definition\label{defsetofT3topologies}
$\X$ is taken as a set.
The set of all topologies $\topology{}$ on $\X$ such that $\opair{\X}{\topology{}}$ is a $\sepA{3}$
topological space, will be denoted by $\sCtops{3}{\X}$:
\begin{align}
&\sCtops{3}{\X}:=\cr
&\defset{\topology{}}{\Ctops{\X}}
{\Foreach{\asubset}{\Fclosed{\X}{\topology{}}}\Foreach{\point}{\(\compl{\X}{\asubset}\)}
\bigg[\Exists{\opair{\U}{\V}}{\bigg(\Cprod{\func{\nei{\Xt}}{\asubset}}{\func{\nei{\Xt}}{\seta{\point}}}\bigg)}
\U\cap\V=\empty\bigg]}.\cr
&{}
\end{align}
\endef
%%%%%%%%%%%%%%%%%%%%%%%%%%%%%%%%%%%%%%%%%%%%%%%%%%%%%%%%%%%%%%%%%%%%%%%%%%%%%%%%%%%%%%%%%%%%%%%%%%%%%%%%%%%%%%%%%%%
\theorem\label{thmT3spaceequiv0}
$\Xt=\opair{\X}{\topology{}}$
is taken as a topological space.
$\Xt$ is a $\sepA{3}$ topological space if and only if for every $\point$ in $\X$,
every neighborhood of $\seta{\point}$ in $\Xt$ includes the closure of at least one neighborhood of $\seta{\point}$
in $\Xt$. That is,
\begin{align}
\(\topology{}\in\sCtops{3}{\X}\)\thenn
\[\Foreach{\point}{\X}
\Foreach{\U}{\func{\nei{\Xt}}{\seta{\point}}}
\bigg(\Exists{\V}{\func{\nei{\Xt}}{\seta{\point}}}
\func{\Cl{\Xt}}{\V}\subseteq\U\bigg)\].
\end{align}
\prooff
\begin{itemize}
\item[${\textbf{\textsf{p1}}}$]
It is assumed that $\topology{}\in\sCtops{3}{\X}$.
Then according to \refdef{defsetofT3topologies},
\begin{equation}\label{thmT3spaceequiv0p1eq1}
\Foreach{\asubset}{\Fclosed{\X}{\topology{}}}\Foreach{\point}{\(\compl{\X}{\asubset}\)}
\bigg[\Exists{\opair{\U}{\V}}{\(\Cprod{\func{\nei{\Xt}}{\asubset}}{\func{\nei{\Xt}}{\seta{\point}}}\)}
\U\cap\V=\empty\bigg].
\end{equation}
\begin{itemize}
\item[${\textbf{\textsf{p1-1}}}$]
$\point$ is taken as an arbitrary element of $\X$, and $\U$ as an arbitrary element of $\func{\nei{\Xt}}{\seta{\point}}$.
Then according to \refdef{defnbdclassofsets} and \refdef{deffamilyofclosedsets}, it is clear that $\(\compl{\X}{\U}\)$
is a closed set of $\Xt$ not containing $\point$.
\begin{align}
\(\compl{\X}{\U}\)&\in\Fclosed{\X}{\topology{}},
\label{thmT3spaceequiv0p11eq1}\\
\point&\in\[\compl{\X}{\(\compl{\X}{\U}\)}\].
\label{thmT3spaceequiv0p11eq2}
\end{align}
Therefore, according to \Ref{thmT3spaceequiv0p1eq1}, it is inferred that
there exists a neighborhood $\U_0$ of $\(\compl{\X}{\U}\)$ in $\Xt$, and a neighborhood $\V$ of $\seta{\point}$
in $\Xt$ that does not intersect each other.
\begin{equation}
\Existsis{\opair{\U_0}{\V}}
{\big(\Cprod{\func{\nei{\Xt}}{\compl{\X}{\U}}}{\func{\nei{\Xt}}{\seta{\point}}}\big)}
\U_0\cap\V=\empty.
\end{equation}
Considering that $\U_0$ is an open set that does not intersect $\V$, and according to \refthm{thmopensetsintersectingclosure},
it becomes clear that $\U_0$ does not intersect the closure of $\V$ in $\Xt$ too.
\begin{equation}
\U_0\cap\func{\Cl{\Xt}}{\V}=\empty.
\end{equation}
Therefore,
\begin{equation}
\func{\Cl{\Xt}}{\V}\subseteq\(\compl{\X}{\U_0}\).
\end{equation}
Moreover, considering that $\(\compl{\X}{\U}\)\subseteq\U_0$, it is clear that,
\begin{equation}
\(\compl{\X}{\U_0}\)\subseteq\U.
\end{equation}
Therefore,
\begin{equation}
\func{\Cl{\Xt}}{\V}\subseteq\U.
\end{equation}
\endp
\end{itemize}
Therefore,
\begin{equation*}
\[\Foreach{\point}{\X}
\Foreach{\U}{\func{\nei{\Xt}}{\seta{\point}}}
\bigg(\Exists{\V}{\func{\nei{\Xt}}{\seta{\point}}}
\func{\Cl{\Xt}}{\V}\subseteq\U\bigg)\].
\end{equation*}
\endp
\end{itemize}
\begin{itemize}
\item[${\textbf{\textsf{p2}}}$]
It is assumed that,
\begin{equation}\label{thmT3spaceequiv0p2eq1}
\[\Foreach{\point}{\X}
\Foreach{\U}{\func{\nei{\Xt}}{\seta{\point}}}
\bigg(\Exists{\V}{\func{\nei{\Xt}}{\seta{\point}}}
\func{\Cl{\Xt}}{\V}\subseteq\U\bigg)\].
\end{equation}
\begin{itemize}
\item[${\textbf{\textsf{p2-1}}}$]
$\asubset$ is taken as an arbitrary element of $\Fclosed{\X}{\topology{}}$, and $\point$
as an arbitrary element of $\(\compl{\X}{\asubset}\)$. According to \refdef{deffamilyofclosedsets}, it is clear that,
\begin{equation}\label{thmT3spaceequiv0p21eq1}
\(\compl{\X}{\asubset}\)\in\topology{},
\end{equation}
and hence according to \refdef{defnbdclassofsets},
$\(\compl{\X}{\asubset}\)$ is a neighborhood of $\seta{\point}$ in $\Xt$:
\begin{equation}\label{thmT3spaceequiv0p21eq2}
\(\compl{\X}{\asubset}\)\in\func{\nei{\Xt}}{\seta{\point}}.
\end{equation}
Therefore, according to \Ref{thmT3spaceequiv0p2eq1},
\begin{equation}
\Existsis{\V_0}{\func{\nei{\Xt}}{\seta{\point}}}
\func{\Cl{\Xt}}{\V_0}\subseteq\(\compl{\X}{\asubset}\).
\end{equation}
Therefore,
\begin{equation}
\asubset\subseteq\[\compl{\X}{\func{\Cl{\Xt}}{\V_0}}\].
\end{equation}
Moreover, according to \refthm{thmclosureofclosedset},
$\func{\Cl{\Xt}}{\V_0}$ is a closed set of $\Xt$, and hence according to \refdef{deffamilyofclosedsets},
$\[\compl{\X}{\func{\Cl{\Xt}}{\V_0}}\]$ is an open set of $\Xt$.
\begin{equation}
\[\compl{\X}{\func{\Cl{\Xt}}{\V_0}}\]\in\topology{}.
\end{equation}
Therefore, according to \refdef{defnbdclassofsets},
$\[\compl{\X}{\func{\Cl{\Xt}}{\V_0}}\]$ is a neighborhood of $\asubset$ in $\Xt$.
\begin{equation}
\[\compl{\X}{\func{\Cl{\Xt}}{\V_0}}\]\in\func{\nei{\Xt}}{\asubset}.
\end{equation}
Moreover, considering that $\V_0\subseteq\func{\Cl{\Xt}}{\V_0}$, it is clear that,
\begin{equation}
\V_0\cap\[\compl{\X}{\func{\Cl{\Xt}}{\V_0}}\]=\empty.
\end{equation}
Therefore,
\begin{equation}
\Existsis{\opair{\[\compl{\X}{\func{\Cl{\Xt}}{\V_0}}\]}{\V_0}}
{\big(\Cprod{\func{\nei{\Xt}}{\asubset}}{\func{\nei{\Xt}}{\seta{\point}}}\big)}
\V_0\cap\[\compl{\X}{\func{\Cl{\Xt}}{\V_0}}\]=\empty.
\end{equation}
\endp
\end{itemize}
Therefore,
\begin{equation}
\Foreach{\asubset}{\Fclosed{\X}{\topology{}}}\Foreach{\point}{\(\compl{\X}{\asubset}\)}
\bigg[\Exists{\opair{\U}{\V}}{\(\Cprod{\func{\nei{\Xt}}{\asubset}}{\func{\nei{\Xt}}{\seta{\point}}}\)}
\U\cap\V=\empty\bigg],
\end{equation}
which according to \refdef{defsetofT3topologies}, means,
\begin{equation}
\topology{}\in\sCtops{3}{\X}.
\end{equation}
\endp
\end{itemize}
\endthm
%%%%%%%%%%%%%%%%%%%%%%%%%%%%%%%%%%%%%%%%%%%%%%%%%%%%%%%%%%%%%%%%%%%%%%%%%%%%%%%%%%%%%%%%%%%%%%%%%%%%%%%%%%%%%%%%%%%
\theorem\label{thmT3ishereditary}
$\Xt=\opair{\X}{\topology{}}$
is taken as a topological space. If $\Xt$ is a $\sepA{3}$ topological space, then for every $\asubset$
in $\CSs{\X}$,
$\opair{\asubset}{\stopology{\topology{}}{\asubset}}$ is a $\sepA{3}$ topological space:
\begin{equation}
\bigg(\topology{}\in\sCtops{3}{\X}\bigg)\then
\bigg[\Foreach{\asubset}{\CSs{\X}}
\stopology{\topology{}}{\asubset}\in\sCtops{3}{\asubset}\bigg].
\end{equation}
\proof
It is assumed that $\Xt$ is a $\sepA{3}$ topological space. Then according to \refdef{defT3space},
\begin{equation}\label{thmT3ishereditarypeq1}
\Foreach{\asubset}{\Fclosed{\X}{\topology{}}}\Foreach{\point}{\(\compl{\X}{\asubset}\)}
\bigg[\Exists{\opair{\U}{\V}}{\(\Cprod{\func{\nei{\Xt}}{\asubset}}{\func{\nei{\Xt}}{\seta{\point}}}\)}
\U\cap\V=\empty\bigg].
\end{equation}
\begin{itemize}
\item[${\textbf{\textsf{p1}}}$]
$\asubset$
is taken as a subset of $\X$.
\begin{itemize}
\item[${\textbf{\textsf{p1-1}}}$]
$\bsubset$ is taken as an element of $\Fclosed{\asubset}{\stopology{\topology{}}{\asubset}}$, and $\point$
as an element of $\(\compl{\asubset}{\bsubset}\)$. According to \refthm{thmsubspaceclosedsets},
\begin{equation}\label{thmT3ishereditaryp11eq1}
\Existsis{\p{\bsubset}}{\Fclosed{\X}{\topology{}}}
\bsubset=\asubset\cap\p{\bsubset}.
\end{equation}
It is clear that,
\begin{equation}\label{thmT3ishereditaryp11eq2}
\point\in\(\compl{\X}{\p{\bsubset}}\).
\end{equation}
Therefore, according to \Ref{thmT3ishereditarypeq1},
\begin{equation}\label{thmT3ishereditaryp11eq3}
\Existsis{\opair{\U}{\V}}{\(\Cprod{\func{\nei{\Xt}}{\p{\bsubset}}}{\func{\nei{\Xt}}{\seta{\point}}}\)}
\U\cap\V=\empty
\end{equation}
Considering that $\U$ is a neighborhood of $\p{\bsubset}$ in $\Xt$, according to \refdef{defnbdclassofsets},
\refdef{defsubspacetopology1}, and \Ref{thmT3ishereditaryp11eq1}, it is clear that
$\(\asubset\cap\U\)$ is a neighborhood of $\bsubset$ in the topological space
$\opair{\asubset}{\stopology{\topology{}}{\asubset}}$.
\begin{equation}
\(\asubset\cap\U\)\in\func{\nei{\opair{\asubset}{\stopology{\topology{}}{\asubset}}}}{\bsubset}.
\end{equation}
Similarly, according to \refdef{defnbdclassofsets} and \refdef{defsubspacetopology1}, it is clear that
$\(\asubset\cap\V\)$ is a neighborhood of $\seta{\point}$ in the topological space
$\opair{\asubset}{\stopology{\topology{}}{\asubset}}$.
\begin{equation}
\(\asubset\cap\V\)\in\func{\nei{\opair{\asubset}{\stopology{\topology{}}{\asubset}}}}{\seta{\point}}.
\end{equation}
Moreover, according to \Ref{thmT3ishereditaryp11eq3},
\begin{equation}
\(\asubset\cap\U\)\cap\(\asubset\cap\V\)=\empty.
\end{equation}
\endp
\end{itemize}
Therefore,
\begin{align}
\Foreach{\bsubset}{\Fclosed{\asubset}{\stopology{\topology{}}{\asubset}}}\Foreach{\point}{\(\compl{\asubset}{\bsubset}\)}
\bigg[\Exists{\opair{\U}{\V}}{\(\Cprod{\func{\nei{\opair{\asubset}{\stopology{\topology{}}{\asubset}}}}{\bsubset}}
{\func{\nei{\opair{\asubset}{\stopology{\topology{}}{\asubset}}}}{\seta{\point}}}\)}
\U\cap\V=\empty\bigg],
\end{align}
which according to \refdef{defsetofT3topologies}, means,
\begin{equation}
\stopology{\topology{}}{\asubset}\in\sCtops{3}{\asubset}.
\end{equation}
\endp
\end{itemize}
\endthm
%%%%%%%%%%%%%%%%%%%%%%%%%%%%%%%%%%%%%%%%%%%%%%%%%%%%%%%%%%%%%%%%%%%%%%%%%%%%%%%%%%%%%%%%%%%%%%%%%%%%%%%%%%%%%%%%%%%
\definition\label{defregularspace}
$\Xt=\opair{\X}{\topology{}}$
is taken as a topological space.
$\Xt$ is referred to as a $\quotl$regular topological space$\quotr$ iff
$\Xt$ is simultaneously a $\sepA{2}$ and a $\sepA{3}$ topological space. That is,
$\Xt$ is called a regular topological space iff
\begin{equation}
\topology{}\in\(\sCtops{2}{\X}\cap\sCtops{3}{\X}\).
\end{equation}
\endef
%%%%%%%%%%%%%%%%%%%%%%%%%%%%%%%%%%%%%%%%%%%%%%%%%%%%%%%%%%%%%%%%%%%%%%%%%%%%%%%%%%%%%%%%%%%%%%%%%%%%%%%%%%%%%%%%%%%
%%%%%%%%%%%%%%%%%%%%%%%%%%%%%%%%%%%%%%%%%%%%%%%%%%%%%%%%%%%%%%%%%%%%%%%%%%%%%%%%%%%%%%%%%%%%%%%%%%%%%%%%%%%%%%%%%%%%%%%%%%%%%%%%%%%%%%%%%%%%%%%%%%%%%%
%%%%%%%%%%%%%%%%%%%%%%%%%%%%%%%%%%%%%%%%%%%%%%%%%%%%%%%%%%%%%%%%%%%%%%%%%%%%%%%%%%%%%%%%%%%%%%%%%%%%%%%%%%%%%%%%%%%%%%%%%%%%%%%%%%%%%%%%%%%%%%%%%%%%%%
%%%%%%%%%%%%%%%%%%%%%%%%%%%%%%%%%%%%%%%%%%%%%%%%%%%%%%%%%%%%%%%%%%%%%%%%%%%%%%%%%%%%%%%%%%%%%%%%%%%%%%%%%%%%%%%%%%%%%%%%%%%%%%%%%%%%%%%%%%%%%%%%%%%%%%
%%%%%%%%%%%%%%%%%%%%%%%%%%%%%%%%%%%%%%%%%%%%%%%%%%%%%%%%%%%%%%%%%%%%%%%%%%%%%%%%%%%%%%%%%%%%%%%%%%%%%%%%%%%%%%%%%%%%%%%%%%%%%%%%%%%%%%%%%%%%%%%%%%%%%%
%%%%%%%%%%%%%%%%%%%%%%%%%%%%%%%%%%%%%%%%%%%%%%%%%%%%%%%%%%%%%%%%%%%%%%%%%%%%%%%%%%%%%%%%%%%%%%%%%%%%%%%%%%%%%%%%%%%%%%%%%%%%%%%%%%%%%%%%%%%%%%%%%%%%%%
%%%%%%%%%%%%%%%%%%%%%%%%%%%%%%%%%%%%%%%%%%%%%%%%%%%%%%%%%%%%%%%%%%%%%%%%%%%%%%%%%%%%%%%%%%%%%%%%%%%%%%%%%%%%%%%%%%%%%%%%%%%%%%%%%%%%%%%%%%%%%%%%%%%%%%
%%%%%%%%%%%%%%%%%%%%%%%%%%%%%%%%%%%%%%%%%%%%%%%%%%%%%%%%%%%%%%%%%%%%%%%%%%%%%%%%%%%%%%%%%%%%%%%%%%%%%%%%%%%%%%%%%%%%%%%%%%%%%%%%%%%%%%%%%%%%%%%%%%%%%%
%%%%%%%%%%%%%%%%%%%%%%%%%%%%%%%%%%%%%%%%%%%%%%%%%%%%%%%%%%%%%%%%%%%%%%%%%%%%%%%%%%%%%%%%%%%%%%%%%%%%%%%%%%%%%%%%%%%%%%%%%%%%%%%%%%%%%%%%%%%%%%%%%%%%%%
%%%%%%%%%%%%%%%%%%%%%%%%%%%%%%%%%%%%%%%%%%%%%%%%%%%%%%%%%%%%%%%%%%%%%%%%%%%%%%%%%%%%%%%%%%%%%%%%%%%%%%%%%%%%%%%%%%%%%%%%%%%%%%%%%%%%%%%%%%%%%%%%%%%%%%
\section{
$\sepA{4}$ Axiom, and Normal Spaces
}
\definition\label{defT4space}
$\Xt=\opair{\X}{\topology{}}$
is taken as a topological space.
$\Xt$ is referred to as a $\quotl$ $\sepA{4}$ topological space$\quotr$ iff
any pair of non-intersection closed sets of $\Xt$ possess neighborhoods in $\Xt$
that are disjoint. That is $\Xt$ is called a $\sepA{4}$ topological space iff
\begin{align}
&\Foreach{\opair{\asubset_1}{\asubset_2}}{\[\Cprod{\Fclosed{\X}{\topology{}}}{\Fclosed{\X}{\topology{}}}\]}\cr
&\bigg[\bigg(\asubset_1\cap\asubset_2=\empty\bigg)\then\Exists{\opair{\U_1}{\U_2}}{\bigg(\Cprod{\func{\nei{\Xt}}{\asubset_1}}{\func{\nei{\Xt}}{\asubset_2}}\bigg)}
\U_1\cap\U_2=\empty\bigg].
\end{align}
\endef
%%%%%%%%%%%%%%%%%%%%%%%%%%%%%%%%%%%%%%%%%%%%%%%%%%%%%%%%%%%%%%%%%%%%%%%%%%%%%%%%%%%%%%%%%%%%%%%%%%%%%%%%%%%%%%%%%%%
\definition\label{defsetofT4topologies}
$\X$ is taken as a set.
The set of all topologies $\topology{}$ on $\X$ such that $\opair{\X}{\topology{}}$
is a $\sepA{4}$ topological space, will be denoted by $\sCtops{4}{\X}$:
\begin{align}
&\sCtops{4}{\X}:=\left\{\topology{}\in\Ctops{\X}\thickspace:\thickspace
\bigg[\Foreach{\opair{\asubset_1}{\asubset_2}}{\Fclosed{\X}{\topology{}}^2}\right.\cr
&\left.\bigg(\asubset_1\cap\asubset_2=\empty\bigg)\then\Exists{\opair{\U_1}{\U_2}}{\bigg(\Cprod{\func{\nei{\Xt}}{\asubset_1}}{\func{\nei{\Xt}}{\asubset_2}}\bigg)}
\U_1\cap\U_2=\empty\bigg]\right\}.
\end{align}
\endef
%%%%%%%%%%%%%%%%%%%%%%%%%%%%%%%%%%%%%%%%%%%%%%%%%%%%%%%%%%%%%%%%%%%%%%%%%%%%%%%%%%%%%%%%%%%%%%%%%%%%%%%%%%%%%%%%%%%
\theorem\label{thmnormalspaceequiv0}
$\Xt=\opair{\X}{\topology{}}$
is taken as a topological space.
$\Xt$ is a $\sepA{4}$ topological space if and only if for every $\asubset$ in $\Fclosed{\X}{\topology{}}$,
every neighborhood of $\asubset$ in $\Xt$ includes the closure of at least one neighborhood of $\asubset$
in $\Xt$. That is,
\begin{align}\label{thmnormalspaceequiv0p1eq1}
\(\topology{}\in\sCtops{4}{\X}\)\thenn
\[\Foreach{\asubset}{\Fclosed{\X}{\topology{}}}
\Foreach{\U}{\func{\nei{\Xt}}{\asubset}}
\bigg(\Exists{\V}{\func{\nei{\Xt}}{\asubset}}
\func{\Cl{\Xt}}{\V}\subseteq\U\bigg)\].
\end{align}
\proof
\begin{itemize}
\item[${\textbf{\textsf{p1}}}$]
It is assumed that $\topology{}\in\sCtops{4}{\X}$.
Then according to \refdef{defsetofT4topologies},
\begin{align}
&\Foreach{\opair{\asubset_1}{\asubset_2}}{\[\Cprod{\Fclosed{\X}{\topology{}}}{\Fclosed{\X}{\topology{}}}\]}\cr
&\bigg[\bigg(\asubset_1\cap\asubset_2=\empty\bigg)\then\Exists{\opair{\U_1}{\U_2}}{\bigg(\Cprod{\func{\nei{\Xt}}{\asubset_1}}{\func{\nei{\Xt}}{\asubset_2}}\bigg)}
\U_1\cap\U_2=\empty\bigg].
\end{align}
\begin{itemize}
\item[${\textbf{\textsf{p1-1}}}$]
$\asubset$ is taken as an arbitrary element of $\Fclosed{\X}{\topology{}}$, and $\U$
as an arbitrary element of $\func{\nei{\Xt}}{\asubset}$. Then according to \refdef{defnbdclassofsets}
and \refdef{deffamilyofclosedsets}, it is clear that $\(\compl{\X}{\U}\)$
is a closed set of $\Xt$ that does not include $\asubset$.
\begin{align}
\(\compl{\X}{\U}\)&\in\Fclosed{\X}{\topology{}},
\label{thmnormalspaceequiv0p11eq1}\\
\asubset\cap\(\compl{\X}{\U}\)&=\empty.
\label{thmnormalspaceequiv0p11eq2}
\end{align}
Therefore, according to \Ref{thmnormalspaceequiv0p1eq1},
there exists disjoint neighborhoods of $\(\compl{\X}{\U}\)$ and $\asubset$ in $\Xt$:
\begin{equation}
\Existsis{\opair{\U_0}{\V}}
{\big[\Cprod{\func{\nei{\Xt}}{\compl{\X}{\U}}}{\func{\nei{\Xt}}{\asubset}}\big]}
\U_0\cap\V=\empty.
\end{equation}
Considering that $\U_0$ is an open set that does not intersect $\V$, and according to \refthm{thmopensetsintersectingclosure},
it is clear that $\U_0$ does not intersect the closure of $\V$ in $\Xt$:
\begin{equation}
\U_0\cap\func{\Cl{\Xt}}{\V}=\empty.
\end{equation}
Therefore,
\begin{equation}
\func{\Cl{\Xt}}{\V}\subseteq\(\compl{\X}{\U_0}\).
\end{equation}
Moreover, considering that $\(\compl{\X}{\U}\)\subseteq\U_0$, it is clear that,
\begin{equation}
\(\compl{\X}{\U_0}\)\subseteq\U.
\end{equation}
Therefore,
\begin{equation}
\func{\Cl{\Xt}}{\V}\subseteq\U.
\end{equation}
\endp
\end{itemize}
Therefore,
\begin{equation*}
\[\Foreach{\asubset}{\Fclosed{\X}{\topology{}}}
\Foreach{\U}{\func{\nei{\Xt}}{\asubset}}
\bigg(\Exists{\V}{\func{\nei{\Xt}}{\asubset}}
\func{\Cl{\Xt}}{\V}\subseteq\U\bigg)\].
\end{equation*}
\endp
\end{itemize}
\begin{itemize}
\item[${\textbf{\textsf{p2}}}$]
It is assumed that,
\begin{equation}\label{thmnormalspaceequiv0p2eq1}
\[\Foreach{\asubset}{\Fclosed{\X}{\topology{}}}
\Foreach{\U}{\func{\nei{\Xt}}{\asubset}}
\bigg(\Exists{\V}{\func{\nei{\Xt}}{\asubset}}
\func{\Cl{\Xt}}{\V}\subseteq\U\bigg)\].
\end{equation}
\begin{itemize}
\item[${\textbf{\textsf{p2-1}}}$]
$\asubset_1$ and $\asubset_2$ are taken as such elements of $\Fclosed{\X}{\topology{}}$ that,
\begin{equation}\label{thmnormalspaceequiv0p21eq1}
\asubset_1\cap\asubset_2=\empty.
\end{equation}
Considering that $\asubset_1$ is a closed set of $\Xt$, and according to \refdef{deffamilyofclosedsets},
it is clear that,
\begin{equation}\label{thmnormalspaceequiv0p21eq2}
\(\compl{\X}{\asubset_1}\)\in\topology{},
\end{equation}
and hence according to \refdef{defnbdclassofsets} and \Ref{thmnormalspaceequiv0p21eq1},
$\(\compl{\X}{\asubset_1}\)$ is a neighborhood of $\asubset_2$ in $\Xt$.
\begin{equation}\label{thmnormalspaceequiv0p21eq3}
\(\compl{\X}{\asubset}\)\in\func{\nei{\Xt}}{\asubset_2}.
\end{equation}
Therefore, according to \Ref{thmnormalspaceequiv0p2eq1},
\begin{equation}
\Existsis{\V_0}{\func{\nei{\Xt}}{\asubset_2}}
\func{\Cl{\Xt}}{\V_0}\subseteq\(\compl{\X}{\asubset_1}\).
\end{equation}
Therefore,
\begin{equation}
\asubset_1\subseteq\[\compl{\X}{\func{\Cl{\Xt}}{\V_0}}\].
\end{equation}
Moreover, according to \refthm{thmclosureofclosedset},
$\func{\Cl{\Xt}}{\V_0}$ is a closed set of $\Xt$, and thus according to \refdef{deffamilyofclosedsets},
$\[\compl{\X}{\func{\Cl{\Xt}}{\V_0}}\]$ is an open set of $\Xt$.
\begin{equation}
\[\compl{\X}{\func{\Cl{\Xt}}{\V_0}}\]\in\topology{}.
\end{equation}
Therefore, according to \refdef{defnbdclassofsets},
$\[\compl{\X}{\func{\Cl{\Xt}}{\V_0}}\]$
is a neighborhood of $\asubset_1$ in $\Xt$.
\begin{equation}
\[\compl{\X}{\func{\Cl{\Xt}}{\V_0}}\]\in\func{\nei{\Xt}}{\asubset_1}.
\end{equation}
Moreover, considering that $\V_0\subseteq\func{\Cl{\Xt}}{\V_0}$, it is clear that,
\begin{equation}
\V_0\cap\[\compl{\X}{\func{\Cl{\Xt}}{\V_0}}\]=\empty.
\end{equation}
Therefore,
\begin{equation}
\Existsis{\opair{\[\compl{\X}{\func{\Cl{\Xt}}{\V_0}}\]}{\V_0}}
{\big(\Cprod{\func{\nei{\Xt}}{\asubset_1}}{\func{\nei{\Xt}}{\asubset_2}}\big)}
\V_0\cap\[\compl{\X}{\func{\Cl{\Xt}}{\V_0}}\]=\empty.
\end{equation}
\endp
\end{itemize}
Therefore,
\begin{align}
&\Foreach{\opair{\asubset_1}{\asubset_2}}{\[\Cprod{\Fclosed{\X}{\topology{}}}{\Fclosed{\X}{\topology{}}}\]}\cr
&\bigg[\bigg(\asubset_1\cap\asubset_2=\empty\bigg)\then\Exists{\opair{\U_1}{\U_2}}{\bigg(\Cprod{\func{\nei{\Xt}}{\asubset_1}}{\func{\nei{\Xt}}{\asubset_2}}\bigg)}
\U_1\cap\U_2=\empty\bigg],
\end{align}
which according to \refdef{defsetofT4topologies}, means,
\begin{equation}
\topology{}\in\sCtops{4}{\X}.
\end{equation}
\endp
\end{itemize}
\endthm
%%%%%%%%%%%%%%%%%%%%%%%%%%%%%%%%%%%%%%%%%%%%%%%%%%%%%%%%%%%%%%%%%%%%%%%%%%%%%%%%%%%%%%%%%%%%%%%%%%%%%%%%%%%%%%%%%%%
\definition\label{defnormalspace}
$\Xt=\opair{\X}{\topology{}}$
is taken as a topological space.
$\Xt$ is referred to as a $\quotl$normal topological space$\quotr$ iff $\Xt$
is simultaneously  a $\sepA{2}$ and a $\sepA{4}$ topological space. That is,
$\Xt$ is called a normal topological space iff
\begin{equation}
\topology{}\in\(\sCtops{2}{\X}\cap\sCtops{4}{\X}\).
\end{equation}
\endef
%%%%%%%%%%%%%%%%%%%%%%%%%%%%%%%%%%%%%%%%%%%%%%%%%%%%%%%%%%%%%%%%%%%%%%%%%%%%%%%%%%%%%%%%%%%%%%%%%%%%%%%%%%%%%%%%%%%
\theorem\label{thmclosedsetofnormalspaceisnormal}
$\Xt=\opair{\X}{\topology{}}$ is taken as a topological space. If
$\Xt$ is a normal topological space, then for every $\asubset$ in $\Fclosed{\X}{\topology{}}$,
$\opair{\asubset}{\stopology{\topology{}}{\asubset}}$ is a normal topological space:
\begin{equation}
\bigg[\topology{}\in\(\sCtops{2}{\X}\cap\sCtops{4}{\X}\)\bigg]\then
\bigg[\Foreach{\asubset}{\Fclosed{\X}{\topology{}}}
\stopology{\topology{}}{\asubset}\in\(\sCtops{2}{\X}\cap\sCtops{4}{\X}\)\bigg].
\end{equation}
\endthm
\chapteR{
Compactness
}
\thispagestyle{fancy}
\section{
Open Covering, and Relative Open Covering
}
\definition\label{defrelopencover}
$\Xt=\opair{\X}{\topology{}}$
is taken as a topological space.
The function $\ocov{\Xt}$ is defined as,
\begin{align}
&\ocov{\Xt}\indef\Func{\CSs{\X}}{\CSs{\CSs{\CSs{\X}}}},\cr
&\Foreach{\asubset}{\CSs{\X}}
\func{\ocov{\Xt}}{\asubset}\eqdef
\defset{\acover}{\CSs{\CSs{\X}}}{\[\(\union{\acover}\)\supseteq\asubset,~\acover\subseteq\topology{}\]}.\cr
&{}
\end{align}
The function $\focov{\Xt}$ is defined as,
\begin{align}
&\focov{\Xt}\indef\Func{\CSs{\X}}{\CSs{\CSs{\CSs{\X}}}},\cr
&\Foreach{\asubset}{\CSs{\X}}
\func{\focov{\Xt}}{\asubset}\eqdef
\defset{\acover}{\func{\ocov{\Xt}}{\asubset}}{\CarD{\acover}\in\Zpz}.\cr
&{}
\end{align}
\begin{itemize}
\item
For every $\asubset$ in $\CSs{\X}$, every element of $\func{\ocov{\Xt}}{\asubset}$
is referred to as a $\quotl$relative open covering of $\asubset$ in the topological space $\Xt$$\quotr$.
\item
For every $\asubset$ in $\CSs{\X}$, every element of $\func{\focov{\Xt}}{\asubset}$ is called a
$\quotl$finite relative open covering of $\asubset$ in the topological space $\Xt$$\quotr$.
\item
Every element of $\func{\ocov{\Xt}}{\X}$
is referred to as an $\quotl$open covering of the topological space $\Xt$$\quotr$, and,
\begin{equation}
\Ocov{\Xt}:=\func{\ocov{\Xt}}{\X}.
\end{equation}
\item
Every element of $\func{\focov{\Xt}}{\X}$ is referred to as a $\quotl$finite open covering of
the topological space $\Xt$$\quotr$, and,
\begin{equation}
\fOcov{\Xt}:=\func{\focov{\Xt}}{\X}.
\end{equation}
\end{itemize}
\endef
%%%%%%%%%%%%%%%%%%%%%%%%%%%%%%%%%%%%%%%%%%%%%%%%%%%%%%%%%%%%%%%%%%%%%%%%%%%%%%%%%%%%%%%%%%%%%%%%%%%%%%%%%%%%%%%%%%%
\corollary\label{corunionofopencoverisaneiofset}
$\Xt=\opair{\X}{\topology{}}$
is taken as a topological space, and $\asubset$ as a subset of $\X$.
\begin{itemize}
\item
Every finite open covering of $\Xt$ is an open covering of $\Xt$.
\begin{equation}
\fOcov{\Xt}\subseteq\Ocov{\Xt}.
\end{equation}
\item
Every finite open covering of $\asubset$ in $\Xt$ is an open covering of $\asubset$ in $\Xt$:
\begin{equation}
\func{\focov{\Xt}}{\asubset}\subseteq\func{\ocov{\Xt}}{\asubset}.
\end{equation}
\item
The union of any open covering of $\asubset$ in $\Xt$ is a neighborhood of $\asubset$ in $\Xt$.
\begin{equation}
\Foreach{\acover}{\func{\ocov{\Xt}}{\asubset}}
\(\union{\acover}\)\in\func{\nei{\Xt}}{\asubset}.
\end{equation}
\end{itemize}
\endcor
%%%%%%%%%%%%%%%%%%%%%%%%%%%%%%%%%%%%%%%%%%%%%%%%%%%%%%%%%%%%%%%%%%%%%%%%%%%%%%%%%%%%%%%%%%%%%%%%%%%%%%%%%%%%%%%%%%%
\section{
Compact Topological Spaces
}
\definition\label{defcompactness}
$\Xt=\opair{\X}{\topology{}}$
is taken as a topological space.
$\Xt$ is referred to as a $\quotl$compact topological space$\quotr$ iff
every open covering of $\Xt$ includes at least one finite open covering of $\Xt$. That is,
$\Xt$ is called a compact topological space iff
\begin{equation}
\Foreach{\acover}{\Ocov{\Xt}}
\defset{\p{\acover}}{\fOcov{\Xt}}
{\p{\acover}\subseteq\acover}\neq\empty.
\end{equation}
\begin{itemize}
\item
It is said that $\quotl$$\X$ is compact relative to the topology $\topology{}$$\quotr$ iff
$\opair{\X}{\topology{}}$ is a compact topological space.
\item
$\Xt$ is referred to as a $\quotl$non-compact topological space$\quotr$ iff $\Xt$ is not a compact topological space.
That is, $\Xt$ is called a non-compact topological space iff
\begin{equation*}
\Exists{\acover}{\Ocov{\Xt}}
\defset{\p{\acover}}{\fOcov{\Xt}}
{\p{\acover}\subseteq\acover}=\empty.
\end{equation*}
\end{itemize}
\endef
%%%%%%%%%%%%%%%%%%%%%%%%%%%%%%%%%%%%%%%%%%%%%%%%%%%%%%%%%%%%%%%%%%%%%%%%%%%%%%%%%%%%%%%%%%%%%%%%%%%%%%%%%%%%%%%%%%%
\definition\label{defcompactsets}
$\Xt=\opair{\X}{\topology{}}$
is taken as a topological space.
The set of all subsets of $\X$ that are compact relative to the topology induced from $\topology{}$,
will be denoted by $\compacts{\Xt}$. That is,
\begin{align}
&\compacts{\Xt}:=\cr
&\defset{\asubset}{\CSs{\X}}
{\bigg[\Foreach{\acover}{\Ocov{\opair{\asubset}{\stopology{\topology{}}{\asubset}}}}
\defset{\p{\acover}}{\fOcov{\opair{\asubset}{\stopology{\topology{}}{\asubset}}}}
{\p{\acover}\subseteq\acover}\neq\empty\bigg]}.\cr
&{}
\end{align}
\begin{itemize}
\item
For every $\asubset$ in $\CSs{\X}$,
It is said that $\quotl$$\asubset$ is compact in the topological space $\Xt$$\quotr$, or
$\asubset$ is called a $\quotl$compact set of the topological space $\Xt$$\quotr$ iff
$\opair{\asubset}{\stopology{\topology{}}{\asubset}}$
is a compact topological space (or $\asubset\in\compacts{\Xt}$).
\end{itemize}
\endef
%%%%%%%%%%%%%%%%%%%%%%%%%%%%%%%%%%%%%%%%%%%%%%%%%%%%%%%%%%%%%%%%%%%%%%%%%%%%%%%%%%%%%%%%%%%%%%%%%%%%%%%%%%%%%%%%%%%
\theorem\label{thmcompactsetequiv0}
$\Xt=\opair{\X}{\topology{}}$
is taken as a topological space.
\begin{equation}
\Foreach{\asubset}{\CSs{\X}}
\bigg(
\asubset\in\compacts{\Xt}\thenn
\asubset\in\compacts{\opair{\asubset}{\stopology{\topology{}}{\asubset}}}
\bigg).
\end{equation}
\proof
According to \refdef{defcompactsets}, it is clear.
\endthm
%%%%%%%%%%%%%%%%%%%%%%%%%%%%%%%%%%%%%%%%%%%%%%%%%%%%%%%%%%%%%%%%%%%%%%%%%%%%%%%%%%%%%%%%%%%%%%%%%%%%%%%%%%%%%%%%%%%
\theorem\label{thmcompactsubspaceofcompactset}
$\Xt=\opair{\X}{\topology{}}$
is taken as a topological space, and $\Y$ as a subset of $\X$. For every $\asubset$ in $\CSs{\X}$,
$\asubset$ is a compact set of $\Xt$ if and only if $\asubset$ is a compact set of
$\opair{\Y}{\stopology{\topology{}}{\Y}}$. That is,
\begin{equation}
\compacts{\opair{\Y}{\stopology{\topology{}}{\Y}}}=
\CSs{\Y}\cap\compacts{\Xt}.
\end{equation}
\prooff
$\asubset$ is taken as an arbitrary subset of $\Y$. According to \refdef{defcompactsets},
\refthm{thmsubspacetopologytransitivity}, and
\refthm{thmcompactsetequiv0},
\begin{align}
\compacts{\opair{\Y}{\stopology{\topology{}}{\Y}}}&=
\defset{\asubset}{\CSs{\Y}}{\asubset\in\compacts{\opair{\asubset}{\stopology{\stopology{\topology{}}{\Y}}{\asubset}}}}\cr
&=\defset{\asubset}{\CSs{\Y}}{\asubset\in\compacts{\opair{\asubset}{\stopology{\topology{}}{\asubset}}}}\cr
&=\defset{\asubset}{\CSs{\Y}}{\asubset\in\compacts{\Xt}}\cr
&=\CSs{\Y}\cap\compacts{\Xt}.
\end{align}
\endthm
%%%%%%%%%%%%%%%%%%%%%%%%%%%%%%%%%%%%%%%%%%%%%%%%%%%%%%%%%%%%%%%%%%%%%%%%%%%%%%%%%%%%%%%%%%%%%%%%%%%%%%%%%%%%%%%%%%%
\theorem\label{thmcompactsets}
$\Xt=\opair{\X}{\topology{}}$
is taken as a topological space.
For every $\asubset$ in $\CSs{\X}$,
$\asubset$ is a compact set of $\Xt$ if and only if every open covering of $\asubset$ in $\Xt$
includes at least one one finite open covering of $\asubset$ in $\Xt$. That is,
\begin{align}
\compacts{\Xt}=\defset{\asubset}{\CSs{\X}}
{\bigg[\Foreach{\acover}{\func{\ocov{\Xt}}{\asubset}}
\defset{\p{\acover}}{\func{\focov{\Xt}}{\asubset}}
{\p{\acover}\subseteq\acover}\neq\empty\bigg]}.
\end{align}
\prooff
\begin{itemize}
\item[${\textbf{\textsf{p1}}}$]
$\asubset$ is taken as an arbitrary element of $\compacts{\Xt}$. Then according to \refdef{defcompactsets},
\begin{equation}\label{thmcompactsetsp1eq1}
\Foreach{\acover}{\Ocov{\opair{\asubset}{\stopology{\topology{}}{\asubset}}}}
\defset{\p{\acover}}{\fOcov{\opair{\asubset}{\stopology{\topology{}}{\asubset}}}}
{\p{\acover}\subseteq\acover}\neq\empty.
\end{equation}
\begin{itemize}
\item[${\textbf{\textsf{p1-1}}}$]
$\acover$ is taken as an arbitrary element of $\func{\ocov{\Xt}}{\asubset}$. Then according to \refdef{defrelopencover},
\begin{align}
\(\union{\acover}\)&\supseteq\asubset,
\label{thmcompactsetsp11eq1}\\
\acover&\subseteq\topology{}.
\label{thmcompactsetsp11eq2}
\end{align}
The set $\acover_{\asubset}$ is defined as,
\begin{equation}\label{thmcompactsetsp11eq3}
\acover_{\asubset}:=\defSet{\asubset\cap\U}{\U\in\acover}.
\end{equation}
\Ref{thmcompactsetsp11eq1} and
\Ref{thmcompactsetsp11eq3}
imply that,
\begin{equation}\label{thmcompactsetsp11eq4}
\union{\acover_{\asubset}}=\asubset.
\end{equation}
According to \refdef{defsubspacetopology1},
\Ref{thmcompactsetsp11eq2}, and
\Ref{thmcompactsetsp11eq3},
\begin{equation}\label{thmcompactsetsp11eq5}
\acover_{\asubset}\subseteq\stopology{\topology{}}{\asubset}.
\end{equation}
According to \refdef{defrelopencover}, \Ref{thmcompactsetsp11eq4}, and
\Ref{thmcompactsetsp11eq5}, it is inferred that $\acover_{\asubset}$
is an open covering of $\opair{\asubset}{\stopology{\topology{}}{\asubset}}$:
\begin{equation}\label{thmcompactsetsp11eq6}
\acover_{\asubset}\in\Ocov{\opair{\asubset}{\stopology{\topology{}}{\asubset}}}.
\end{equation}
Therefore, according to \Ref{thmcompactsetsp1eq1},
\begin{equation}\label{thmcompactsetsp11eq7}
\Existsis{\p{\acover_{\asubset}}}{\fOcov{\opair{\asubset}{\stopology{\topology{}}{\asubset}}}}
\p{\acover_{\asubset}}\subseteq\acover_{\asubset}.
\end{equation}
Therefore according to \Ref{thmcompactsetsp11eq3},
\begin{equation}\label{thmcompactsetsp11eq8}
\Existssis{\p{\acover}}{\acover}
\bigg[\p{\acover_{\asubset}}=\defSet{\asubset\cap\U}{\U\in\p{\acover}},~
\CarD{\p{\acover}}=\CarD{\p{\acover_{\asubset}}}\bigg].
\end{equation}
Thus considering that $\p{\acover_{\asubset}}$ is a finite set,
$\p{\acover}$ is also a finite set:
\begin{equation}\label{thmcompactsetsp11eq9}
\CarD{\p{\acover}}\in\Zpz.
\end{equation}
Moreover, considering that $\acover$ is a collection of open sets of $\Xt$, and
$\p{\acover}\subseteq\acover$, it is clear that $\p{\acover}$
is also a collection of open sets of $\Xt$:
\begin{equation}\label{thmcompactsetsp11eq10}
\p{\acover}\subseteq\topology{}.
\end{equation}
Furthermore, considering that $\p{\acover_\asubset}$ is an open covering of $\Xt$, and according to \Ref{thmcompactsetsp11eq8},
it os clear that,
\begin{align}\label{thmcompactsetsp11eq11}
\asubset\cap\(\union{\p{\acover}}\)&=
\union{\p{\acover_{\asubset}}}\cr
&=\asubset.
\end{align}
Therefore,
\begin{equation}\label{thmcompactsetsp11eq12}
\(\union{\p{\acover}}\)\supseteq\asubset.
\end{equation}
According to \refdef{defrelopencover},
\Ref{thmcompactsetsp11eq9},
\Ref{thmcompactsetsp11eq10}, and \Ref{thmcompactsetsp11eq12},
\begin{equation}\label{thmcompactsetsp11eq13}
\p{\acover}\in\func{\focov{\Xt}}{\asubset}.
\end{equation}
Therefore according to \Ref{thmcompactsetsp11eq8} and \Ref{thmcompactsetsp11eq13}, it is clear that,
\begin{equation}
\Existsis{\p{\acover}}{\func{\focov{\Xt}}{\asubset}}
\p{\acover}\subseteq\acover.
\end{equation}
\endp
\end{itemize}
Therefore,
\begin{equation*}
\Foreach{\acover}{\func{\ocov{\Xt}}{\asubset}}
\defset{\p{\acover}}{\func{\focov{\Xt}}{\asubset}}
{\p{\acover}\subseteq\acover}\neq\empty.
\end{equation*}
\endp
\end{itemize}
\begin{itemize}
\item[${\textbf{\textsf{p2}}}$]
$\asubset$
is taken as such an arbitrary element of $\CSs{\X}$ that,
\begin{equation}\label{thmcompactsetsp2eq1}
\Foreach{\acover}{\func{\ocov{\Xt}}{\asubset}}
\defset{\p{\acover}}{\func{\focov{\Xt}}{\asubset}}
{\p{\acover}\subseteq\acover}\neq\empty.
\end{equation}
\begin{itemize}
\item[${\textbf{\textsf{p2-1}}}$]
$\acover_\asubset$
is taken as an arbitrary element of $\Ocov{\opair{\asubset}{\stopology{\topology{}}{\asubset}}}$
Then according to \refdef{defrelopencover},
\begin{align}
\(\union{\acover_\asubset}\)&=\asubset,
\label{thmcompactsetsp21eq1}\\
\acover_\asubset&\subseteq\stopology{\topology{}}{\asubset}.
\label{thmcompactsetsp21eq2}
\end{align}
According to \refdef{defsubspacetopology1} and
\Ref{thmcompactsetsp21eq2},
\begin{equation}\label{thmcompactsetsp21eq3}
\Existssis{\acover}{\topology{}}
\acover_{\asubset}=\defSet{\asubset\cap\U}{\U\in\acover}.
\end{equation}
Therefore,
\begin{align}\label{thmcompactsetsp21eq4}
\(\union{\acover_\asubset}\)=\asubset\cap\(\union{\acover}\).
\end{align}
\Ref{thmcompactsetsp21eq1} and
\Ref{thmcompactsetsp21eq3}
imply that,
\begin{equation}\label{thmcompactsetsp21eq5}
\(\union{\acover}\)\supseteq\asubset.
\end{equation}
According to \refdef{defrelopencover}, \Ref{thmcompactsetsp21eq3}, and
\Ref{thmcompactsetsp21eq5},
\begin{equation}\label{thmcompactsetsp21eq6}
\acover\in\func{\ocov{\Xt}}{\asubset}.
\end{equation}
\Ref{thmcompactsetsp2eq1} and \Ref{thmcompactsetsp21eq6} imply that,
\begin{equation}\label{thmcompactsetsp21eq7}
\Existsis{\p{\acover}}{\func{\focov{\Xt}}{\asubset}}
\p{\acover}\subseteq\acover.
\end{equation}
According to \refdef{defrelopencover}, and considering that $\acover\subseteq\topology{}$,
\begin{align}
\(\union{\p{\acover}}\)&\supseteq\asubset,
\label{thmcompactsetsp21eq8}\\
\p{\acover}&\subseteq\topology{},
\label{thmcompactsetsp21eq9}\\
\CarD{\p{\acover}}&\in\Zpz.
\label{thmcompactsetsp21eq10}
\end{align}
The set
$\p{\acover_\asubset}$ is defined as,
\begin{equation}\label{thmcompactsetsp21eq11}
\p{\acover_\asubset}:=\defSet{\asubset\cap\U}{\U\in\p{\acover}}.
\end{equation}
\Ref{thmcompactsetsp21eq8} and
\Ref{thmcompactsetsp21eq11} yield,
\begin{equation}\label{thmcompactsetsp21eq12}
\(\union{\p{\acover_\asubset}}\)=\asubset.
\end{equation}
\Ref{thmcompactsetsp21eq9} and
\Ref{thmcompactsetsp21eq11} yield,
\begin{equation}\label{thmcompactsetsp21eq13}
\p{\acover_\asubset}\subseteq\stopology{\topology{}}{\asubset}.
\end{equation}
\Ref{thmcompactsetsp21eq10} and
\Ref{thmcompactsetsp21eq11} imply taht,
\begin{equation}\label{thmcompactsetsp21eq14}
\CarD{\p{\acover_\asubset}}\in\Zpz.
\end{equation}
According to \refdef{defrelopencover},
\Ref{thmcompactsetsp21eq12},
\Ref{thmcompactsetsp21eq13}, and \Ref{thmcompactsetsp21eq14},
$\p{\acover_\asubset}$
is a finite open overing of $\opair{\asubset}{\stopology{\topology{}}{\asubset}}$:
\begin{equation}\label{thmcompactsetsp21eq15}
\p{\acover_\asubset}\in\fOcov{\opair{\asubset}{\stopology{\topology{}}{\asubset}}}.
\end{equation}
Moreover, considering that $\p{\acover}\subseteq\acover$,
\Ref{thmcompactsetsp21eq3} and
\Ref{thmcompactsetsp21eq11} imply that,
\begin{equation}
\p{\acover_\asubset}\subseteq\acover_\asubset.
\end{equation}
Therefore,
\begin{equation}
\Existsis{\p{\acover_\asubset}}{\fOcov{\opair{\asubset}{\stopology{\topology{}}{\asubset}}}}
\p{\acover_\asubset}\subseteq\acover_\asubset.
\end{equation}
\endp
\end{itemize}
Therefore,
\begin{equation}
\Foreach{\acover_\asubset}{\Ocov{\opair{\asubset}{\stopology{\topology{}}{\asubset}}}}
\defset{\p{\acover_\asubset}}{\fOcov{\opair{\asubset}{\stopology{\topology{}}{\asubset}}}}
{\p{\acover_\asubset}\subseteq\acover_\asubset}\neq\empty.
\end{equation}
According to \refdef{defcompactness}, this means,
$\opair{\asubset}{\stopology{\topology{}}{\asubset}}$ is compact.
\begin{equation*}
\asubset\in\compacts{\Xt}.
\end{equation*}
\endp
\end{itemize}
\endthm
%%%%%%%%%%%%%%%%%%%%%%%%%%%%%%%%%%%%%%%%%%%%%%%%%%%%%%%%%%%%%%%%%%%%%%%%%%%%%%%%%%%%%%%%%%%%%%%%%%%%%%%%%%%%%%%%%%%
\theorem\label{thmclosedsetofcompactspaceiscompact}
$\Xt=\opair{\X}{\topology{}}$
is taken as a topological space. If $\Xt$ is compact, then every closed set of $\Xt$
is a compact set of $\Xt$. That is,
\footnote{In the precise language of set theory:\\
$\defset{\topology{}}{\Ctops{\X}}{\X\in\compacts{\opair{\X}{\topology{}}}}\subseteq
\defset{\topology{}}{\Ctops{\X}}
{\Fclosed{\X}{\topology{}}\subseteq\compacts{\opair{\X}{\topology{}}}}.$}
\begin{equation}
\bigg(\X\in\compacts{\Xt}\bigg)\then
\bigg(\Fclosed{\X}{\topology{}}\subseteq\compacts{\Xt}\bigg).
\end{equation}
\prooff
It is assumed that $\Xt$ is compact:
\begin{equation}\label{thmclosedsetofcompactspaceiscompactpeq1}
\X\in\compacts{\Xt}.
\end{equation}
Then according to \refdef{defcompactness},
\begin{equation}\label{thmclosedsetofcompactspaceiscompactpeq2}
\Foreach{\acover}{\Ocov{\Xt}}
\defset{\p{\acover}}{\fOcov{\Xt}}
{\p{\acover}\subseteq\acover}\neq\empty.
\end{equation}
\begin{itemize}
\item[${\textbf{\textsf{p1}}}$]
$\asubset$
is taken as an element of $\Fclosed{\X}{\topology{}}$
(a closed set of $\Xt$).
Then according to \refdef{deffamilyofclosedsets},
\begin{equation}\label{thmclosedsetofcompactspaceiscompactp1eq1}
\(\compl{\X}{\asubset}\)\in\topology{}.
\end{equation}
\begin{itemize}
\item[${\textbf{\textsf{p1-1}}}$]
$\acover$
is taken as an element of $\func{\ocov{\Xt}}{\asubset}$
(an open covering of $\asubset$ in $\Xt$).
Considering that $\func{\ocov{\Xt}}{\asubset}$ is a collection of open sets of $\Xt$,
by adding $\(\compl{\Xt}{\asubset}\)$ to this set, a new collection of open sets of $\Xt$ is formed.
In other words, according to \refdef{defrelopencover} and \Ref{thmclosedsetofcompactspaceiscompactp1eq1},
\begin{equation}\label{thmclosedsetofcompactspaceiscompactp11eq1}
\bigg(\acover\cup\seta{\compl{\X}{\asubset}}\bigg)
\subseteq\topology{}.
\end{equation}
Moreover, according to \refdef{defrelopencover},
$\acover$ covers the set $\asubset$:
\begin{equation}\label{thmclosedsetofcompactspaceiscompactp11eq2}
\(\union{\acover}\)\supseteq\asubset.
\end{equation}
Therefore it is clear that,
$\bigg(\acover\cup\seta{\compl{\X}{\asubset}}\bigg)$
also covers $\X$:
\begin{equation}\label{thmclosedsetofcompactspaceiscompactp11eq3}
\[\union{\bigg(\acover\cup\seta{\compl{\X}{\asubset}}\bigg)}\]=\X.
\end{equation}
According to \refdef{defrelopencover},
\Ref{thmclosedsetofcompactspaceiscompactp11eq1}, and
\Ref{thmclosedsetofcompactspaceiscompactp11eq3},
$\bigg(\acover\cup\seta{\compl{\X}{\asubset}}\bigg)$
is an open covering of $\Xt$:
\begin{equation}\label{thmclosedsetofcompactspaceiscompactp11eq4}
\bigg(\acover\cup\seta{\compl{\X}{\asubset}}\bigg)\in\Ocov{\Xt}.
\end{equation}
\Ref{thmclosedsetofcompactspaceiscompactpeq2} and
\Ref{thmclosedsetofcompactspaceiscompactp11eq4}
imply that,
\begin{equation}\label{thmclosedsetofcompactspaceiscompactp11eq5}
\Existsis{\p{\acover}}{\fOcov{\Xt}}
\p{\acover}\subseteq
\bigg(\acover\cup\seta{\compl{\X}{\asubset}}\bigg),
\end{equation}
and hence according to \refdef{defrelopencover},
\begin{align}
\(\union{\p{\acover}}\)&=\X,
\label{thmclosedsetofcompactspaceiscompactp11eq6}\\
\p{\acover}&\subseteq\topology{},
\label{thmclosedsetofcompactspaceiscompactp11eq7}\\
\CarD{\p{\acover}}&\in\Zpz.
\label{thmclosedsetofcompactspaceiscompactp11eq8}
\end{align}
According to \Ref{thmclosedsetofcompactspaceiscompactp11eq5}, it is clear that,
\begin{equation}\label{thmclosedsetofcompactspaceiscompactp11eq9}
\(\compl{\p{\acover}}{\seta{\compl{\X}{\asubset}}}\)\subseteq\acover.
\end{equation}
According to \Ref{thmclosedsetofcompactspaceiscompactp11eq6}, since $\p{\acover}$ covers
$\X$, $\(\compl{\p{\acover}}{\seta{\compl{\X}{\asubset}}}\)$ must cover the set $\asubset$. Because,
$\(\compl{\X}{\asubset}\)$ has no role in covering $\asubset$.
\begin{equation}\label{thmclosedsetofcompactspaceiscompactp11eq10}
\[\union{\(\compl{\p{\acover}}{\seta{\compl{\X}{\asubset}}}\)}\]
\supseteq\asubset.
\end{equation}
Acording to \Ref{thmclosedsetofcompactspaceiscompactp11eq7}, it is clear that
$\(\compl{\p{\acover}}{\seta{\compl{\X}{\asubset}}}\)$
is a collection of open sets of $\Xt$.
\begin{equation}\label{thmclosedsetofcompactspaceiscompactp11eq11}
\(\compl{\p{\acover}}{\seta{\compl{\X}{\asubset}}}\)\subseteq\topology{}
\end{equation}
It is also clear according to \Ref{thmclosedsetofcompactspaceiscompactp11eq8} that
$\(\compl{\p{\acover}}{\seta{\compl{\X}{\asubset}}}\)$ is a finite set.
\begin{equation}\label{thmclosedsetofcompactspaceiscompactp11eq12}
\CarD{\compl{\p{\acover}}{\seta{\compl{\X}{\asubset}}}}\in\Zpz.
\end{equation}
According to \refdef{defrelopencover},
\Ref{thmclosedsetofcompactspaceiscompactp11eq10},
\Ref{thmclosedsetofcompactspaceiscompactp11eq11},
and \Ref{thmclosedsetofcompactspaceiscompactp11eq12},
$\(\compl{\p{\acover}}{\seta{\compl{\X}{\asubset}}}\)$
is a finite open covering of $\asubset$ in $\Xt$.
\begin{equation}\label{thmclosedsetofcompactspaceiscompactp11eq13}
\(\compl{\p{\acover}}{\seta{\compl{\X}{\asubset}}}\)\in
\func{\focov{\Xt}}{\asubset}.
\end{equation}
Thus according to \Ref{thmclosedsetofcompactspaceiscompactp11eq9} and \Ref{thmclosedsetofcompactspaceiscompactp11eq10},
\begin{equation}
\Existsis{\(\compl{\p{\acover}}{\seta{\compl{\X}{\asubset}}}\)}
{\func{\focov{\Xt}}{\asubset}}
\(\compl{\p{\acover}}{\seta{\compl{\X}{\asubset}}}\)\subseteq\acover.
\end{equation}
\endp
\end{itemize}
Therefore,
\begin{equation}
\Foreach{\acover}{\func{\ocov{\Xt}}{\asubset}}
\defset{\p{\acover}}{\func{\focov{\Xt}}{\asubset}}
{\p{\acover}\subseteq\acover}\neq\empty.
\end{equation}
According to \refthm{thmcompactsets}, this means
$\opair{\asubset}{\stopology{\topology{}}{\asubset}}$ is compact.
\begin{equation}
\asubset\in\compacts{\Xt}.
\end{equation}
\endp
\end{itemize}
\endthm
%%%%%%%%%%%%%%%%%%%%%%%%%%%%%%%%%%%%%%%%%%%%%%%%%%%%%%%%%%%%%%%%%%%%%%%%%%%%%%%%%%%%%%%%%%%%%%%%%%%%%%%%%%%%%%%%%%%
\theorem\label{thmunionofcompactsets}
$\Xt=\opair{\X}{\topology{}}$ is taken as a topological space.
The union of every pair of compact sets of $\Xt$ is a compact set of $\Xt$. That is,
\begin{equation}
\Foreach{\opair{\acompactset{1}}{\acompactset{2}}}
{\Cprod{\compacts{\Xt}}{\compacts{\Xt}}}
\(\acompactset{1}\cup\acompactset{2}\)\in\compacts{\Xt}.
\end{equation}
\prooff
Each $\acompactset{1}$ and $\acompactset{2}$ is taken as an arbitrary element of $\compacts{\Xt}$.
\begin{itemize}
\item[${\textbf{\textsf{p1}}}$]
$\acover$
is taken as an arbitrary element of $\func{\ocov{\Xt}}{\acompactset{1}\cup\acompactset{2}}$
(an open covering of $\asubset$ in $\Xt$).
It is clear that $\acover$ is simultaneously an open covering of $\acompactset{1}$ in $\Xt$,
and an open covering of $\acompactset{2}$ in $\Xt$.
\begin{align}
\acover&\in\func{\ocov{\Xt}}{\acompactset{1}},
\label{thmunionofcompactsetspeq1}\\
\acover&\in\func{\ocov{\Xt}}{\acompactset{2}}.
\label{thmunionofcompactsetspeq2}
\end{align}
Therefore, considering that each $\acompactset{1}$ and $\acompactset{2}$ is a compact set of $\Xt$,
and according to \refthm{thmcompactsets},
\begin{align}
&\Existsis{\acover_1}{\func{\focov{\Xt}}{\acompactset{1}}}
\acover_1\subseteq\acover,
\label{thmunionofcompactsetspeq3}\\
&\Existsis{\acover_2}{\func{\focov{\Xt}}{\acompactset{2}}}
\acover_2\subseteq\acover.
\label{thmunionofcompactsetspeq4}
\end{align}
According to \refdef{defrelopencover},
\begin{align}
\(\union{\acover_1}\)&\supseteq\acompactset{1},
\label{thmunionofcompactsetspeq5}\\
\(\union{\acover_2}\)&\supseteq\acompactset{2},
\label{thmunionofcompactsetspeq6}
\end{align}
and hence,
\begin{equation}\label{thmunionofcompactsetspeq7}
\(\union{\(\acover_1\cup\acover_2\)}\)\supseteq\(\acompactset{1}\cup\acompactset{2}\).
\end{equation}
Moreover, according to \refdef{defrelopencover},
each $\acover_1$ and $\acover_2$ is a finite collection of open sets of $\Xt$. Thus it is clear that $\acover_1\cup\acover_2$
is also a finite collection of open sets of $\Xt$.
\begin{align}
\(\acover_1\cup\acover_2\)&\subseteq\topology{},
\label{thmunionofcompactsetspeq8}\\
\CarD{\acover_1\cup\acover_2}&\in\Zpz.
\label{thmunionofcompactsetspeq9}
\end{align}
According to \refdef{defrelopencover},
\Ref{thmunionofcompactsetspeq7},
\Ref{thmunionofcompactsetspeq8}, and
\Ref{thmunionofcompactsetspeq9},
$\(\acover_1\cup\acover_2\)$ is a finite open covering of $\acompactset{1}\cup\acompactset{2}$ in $\Xt$:
\begin{equation}\label{thmunionofcompactsetspeq10}
\(\acover_1\cup\acover_2\)\in\func{\focov{\Xt}}{\acompactset{1}\cup\acompactset{2}}.
\end{equation}
Moreover, considering that each $\acover_1$ and $\acover_2$ is a subset of $\acover$, it is clear that,
\begin{equation}\label{thmunionofcompactsetspeq11}
\(\acover_1\cup\acover_2\)\subseteq\acover.
\end{equation}
Thus there exists a finite open covering of $\acompactset{1}\cup\acompactset{2}$ in $\Xt$
which is a subset of $\acover$.
\begin{equation}
\Existsis{\(\acover_1\cup\acover_2\)}
{\func{\focov{\Xt}}{\acompactset{1}\cup\acompactset{2}}}
\(\acover_1\cup\acover_2\)\subseteq\acover.
\end{equation}
\endp
\end{itemize}
Therefore,
\begin{equation}
\Foreach{\acover}{\func{\ocov{\Xt}}{\acompactset{1}\cup\acompactset{2}}}
\defset{\p{\acover}}{\func{\focov{\Xt}}{\acompactset{1}\cup\acompactset{2}}}
{\p{\acover}\subseteq\acover}\neq\empty,
\end{equation}
and hence according to \refdef{defcompactsets},
$\acompactset{1}\cup\acompactset{2}$ is a compact set of $\Xt$.
\begin{equation}
\(\acompactset{1}\cup\acompactset{2}\)\in\compacts{\Xt}.
\end{equation}
\endthm
%%%%%%%%%%%%%%%%%%%%%%%%%%%%%%%%%%%%%%%%%%%%%%%%%%%%%%%%%%%%%%%%%%%%%%%%%%%%%%%%%%%%%%%%%%%%%%%%%%%%%%%%%%%%%%%%%%%
\theorem\label{thmintersectionofclosedcompactsets}
$\Xt=\opair{\X}{\topology{}}$ is taken as a topological space.
The intersection of any colection of non-empty, closed, and compact sets of $\Xt$
is a closed and compact set of $\Xt$. That is,
\begin{equation}
\Foreach{\somecompacts{}}{\[\compl{\CSs{\compacts{\Xt}\cap\Fclosed{\X}{\topology{}}}}{\seta{\empty}}\]}
\(\intersection{\somecompacts{}}\)\in\(\compacts{\Xt}\cap\Fclosed{\X}{\topology{}}\).
\end{equation}
\prooff
$\(\compacts{\Xt}\cap\Fclosed{\X}{\topology{}}\)$
is taken as an arbitrary element of $\[\compl{\CSs{\compacts{\Xt}}}{\seta{\empty}}\]$
(a non-empty collection of closed and compact sets of $\Xt$). According to \refthm{thmclosedsets},
$\(\intersection{\somecompacts{}}\)$ is a closed set of $\Xt$.
\begin{equation}\label{thmintersectionofclosedcompactsetspeq1}
\(\intersection{\somecompacts{}}\)\in\Fclosed{\X}{\topology{}}.
\end{equation}
Since $\somecompacts{}$ is non-empty, it has at least an element like $\acompactset{}$. It is trivial that,
\begin{equation}\label{thmintersectionofclosedcompactsetspeq2}
\(\intersection{\somecompacts{}}\)\subseteq\acompactset{}.
\end{equation}
According to \refthm{thmsubspaceclosedsets},
\Ref{thmintersectionofclosedcompactsetspeq1}, and
\Ref{thmintersectionofclosedcompactsetspeq2}, it is clear that
$\(\intersection{\somecompacts{}}\)$
is a closed set of the topological space $\opair{\acompactset{}}{\stopology{\topology{}}{\acompactset{}}}$.s
\begin{equation}\label{thmintersectionofclosedcompactsetspeq3}
\(\intersection{\somecompacts{}}\)\in\Fclosed{\acompactset{}}
{\stopology{\topology{}}{\acompactset{}}}.
\end{equation}
Considering that,
\begin{equation}\label{thmintersectionofclosedcompactsetspeq4}
\acompactset{}\in\compacts{\Xt},
\end{equation}
and according to \refdef{defcompactsets}, it is clear that
$\opair{\acompactset{}}{\stopology{\topology{}}{\acompactset{}}}$ is a compact topological space. In other words,
\begin{equation}\label{thmintersectionofclosedcompactsetspeq5}
\acompactset{}\in\compacts{\opair{\acompactset{}}{\stopology{\topology{}}{\acompactset{}}}}.
\end{equation}
Thus according to \Ref{thmintersectionofclosedcompactsetspeq3} and
\Ref{thmintersectionofclosedcompactsetspeq5},
$\(\intersection{\somecompacts{}}\)$
is a closed set of the compact topological space $\opair{\acompactset{}}{\stopology{\topology{}}{\acompactset{}}}$,
and hence according to \refthm{thmclosedsetofcompactspaceiscompact},
$\(\intersection{\somecompacts{}}\)$, it is a compact set of the topological space
$\opair{\acompactset{}}{\stopology{\topology{}}{\acompactset{}}}$.
\begin{equation}
\(\intersection{\somecompacts{}}\)\in
\compacts{\opair{\acompactset{}}{\stopology{\topology{}}{\acompactset{}}}}.
\end{equation}
Therefore, considering that $\acompactset{}$ is a compact set of $\Xt$, and
$\(\intersection{\somecompacts{}}\)$ is a compact set of
$\opair{\acompactset{}}{\stopology{\topology{}}{\acompactset{}}}$, and according to
\refthm{thmcompactsubspaceofcompactset}, it becomes clear that $\(\intersection{\somecompacts{}}\)$
is a compact set of $\Xt$.
\begin{equation}
\(\intersection{\somecompacts{}}\)\in\compacts{\Xt}.
\end{equation}
\endthm
%%%%%%%%%%%%%%%%%%%%%%%%%%%%%%%%%%%%%%%%%%%%%%%%%%%%%%%%%%%%%%%%%%%%%%%%%%%%%%%%%%%%%%%%%%%%%%%%%%%%%%%%%%%%%%%%%%%
\subsection{
Compactness in Hausdorff Spaces
}
%%%%%%%%%%%%%%%%%%%%%%%%%%%%%%%%%%%%%%%%%%%%%%%%%%%%%%%%%%%%%%%%%%%%%%%%%%%%%%%%%%%%%%%%%%%%%%%%%%%%%%%%%%%%%%%%%%%
\theorem\label{thmcompactsetofT2spacepoint}
$\Xt=\opair{\X}{\topology{}}$
is taken as a $\sepA{2}$ (Hausdorff) space.
Every compact set of $\Xt$ and every point of $\Xt$ out of it
possess disjoint neighborhoods in $\Xt$. That is,
\begin{align}
\Foreach{\acompactset{}}{\(\compl{\compacts{\Xt}}{\seta{\empty}}\)}
\Foreach{\point}{\(\compl{\X}{\acompactset{}}\)}
\bigg(\Exists{\opair{\U}{\V}}
{\Cprod{\func{\nei{\Xt}}{\seta{\point}}}{\func{\nei{\Xt}}{\acompactset{}}}}
\U\cap\V=\empty\bigg).
\end{align}
\prooff
$\asubset$
is taken as an arbitrary element of $\compl{\compacts{\Xt}}{\seta{\empty}}$
(a non-empty compact set of $\Xt$).
\begin{itemize}
\item[${\textbf{\textsf{p1}}}$]
$\point$ is taken as an arbitrary point of $\(\compl{\X}{\asubset}\)$. Considering that,
\begin{equation}\label{thmcompactsetofT2spacepointp1eq1}
\Foreach{\x}{\asubset}
\opair{\point}{\x}\in\psCprod{\X}{\X},
\end{equation}
and according to \refdef{defT2space}, it is clear that there exists at least a pair of functions
$\opair{u}{v}$, each from $\asubset$ to the collection of all open sets of $\Xt$ ($\topology{}$),
such that for every $\x$ in $\asubset$, $\func{u}{\x}$ is a neighborhood of $\seta{\point}$ in $\Xt$, and
$\func{v}{\x}$ is a neighborhood of $\seta{\x}$ in $\Xt$, and
$\func{u}{\x}$ and $\func{v}{\x}$ are disjoint.
\begin{align}\label{thmcompactsetofT2spacepointp1eq2}
&\Existsis{\opair{u}{v}}
{\Cprod{\Func{\asubset}{\topology{}}}{\Func{\asubset}{\topology{}}}}\cr
&\[\Foreach{\x}{\asubset}
\bigg(\func{u}{x}\in\func{\nei{\Xt}}{\seta{\point}},~
\func{v}{x}\in\func{\nei{\Xt}}{\seta{\x}},~
\func{u}{\x}\cap\func{v}{\x}=\empty\bigg)\].
\end{align}
It is clear that,
\begin{align}
\(\union{\func{\image{v}}{\asubset}}\)&\supseteq\asubset,
\label{thmcompactsetofT2spacepointp1eq3}\\
\func{\image{v}}{\asubset}&\subseteq\topology{}.
\label{thmcompactsetofT2spacepointp1eq4}
\end{align}
Therefore, according to \refdef{defrelopencover},
$\func{\image{v}}{\asubset}$ is an open covering of $\asubset$ in $\Xt$:
\begin{equation}\label{thmcompactsetofT2spacepointp1eq5}
\func{\image{v}}{\asubset}\in\func{\ocov{\Xt}}{\asubset}.
\end{equation}
Considering that $\asubset$ is a compact set of $\Xt$, and according to \refthm{thmcompactsets},
there exists a finite open covering of $\asubset$ in $\Xt$ which is a subset of $\func{\image{v}}{\asubset}$.
\begin{equation}\label{thmcompactsetofT2spacepointp1eq6}
\Existsis{\acover}{\func{\focov{\Xt}}{\asubset}}
\acover\subseteq\func{\image{v}}{\asubset}.
\end{equation}
Therefore, according to \refdef{defrelopencover} and non-emptiness of $\asubset$,
\begin{align}
\(\union{\acover}\)&\supseteq\asubset,
\label{thmcompactsetofT2spacepointp1eq7}\\
\CarD{\acover}&\in\Zp.
\label{thmcompactsetofT2spacepointp1eq8}
\end{align}
Therefore there exists a finite subset $\asubset_{0}$ of $\asubset$,
with the same number of elements as that of $\acover$, and
$\func{\image{v}}{\asubset_0}=\acover$.
\begin{equation}\label{thmcompactsetofT2spacepointp1eq9}
\Existssis{\asubset_0}{\asubset}
\bigg(\func{\image{v}}{\asubset_0}=\acover,~
\CarD{\asubset_0}=\CarD{\acover}\bigg).
\end{equation}
$\U$ is defined as the intrsection of all sets $\func{u}{x}$ when $\x$ ranges over $\asubset_0$:
\begin{align}\label{thmcompactsetofT2spacepointp1eq10}
\U:=\(\intersection{\func{\image{u}}{\asubset_0}}\)
\end{align}
Therefore, considering that $\func{\image{u}}{\asubset_0}$ is a finite collection of open sets of $\Xt$,
it is clear that,
\begin{equation}\label{thmcompactsetofT2spacepointp1eq11}
\U\in\topology{}.
\end{equation}
Moreover, considering that $\func{u}{\x}$ contains $\point$, it is clear that,
\begin{equation}\label{thmcompactsetofT2spacepointp1eq12}
\point\in\U.
\end{equation}
Thus $\U$ is a neighborhood of $\seta{\point}$.
\begin{equation}\label{thmcompactsetofT2spacepointp1eq13}
\U\in\func{\nei{\Xt}}{\seta{\point}}.
\end{equation}
Moreover, according to \Ref{thmcompactsetofT2spacepointp1eq6} and
\refcor{corunionofopencoverisaneiofset},
$\(\union{\acover}\)$ is a neighborhood of $\asubset$ in $\Xt$.
\begin{equation}
\(\union{\acover}\)\in\func{\nei{\Xt}}{\asubset}.
\end{equation}
Furthermore, according to \Ref{thmcompactsetofT2spacepointp1eq10},
\begin{equation}\label{thmcompactsetofT2spacepointp1eq14}
\Foreach{\x}{\asubset_0}
\U\subseteq\func{u}{\x},
\end{equation}
and according to \Ref{thmcompactsetofT2spacepointp1eq2},
\begin{equation}\label{thmcompactsetofT2spacepointp1eq15}
\Foreach{\x}{\asubset_0}
\func{u}{\x}\cap\func{v}{\x}=\empty.
\end{equation}
Therefore,
\begin{equation}\label{thmcompactsetofT2spacepointp1eq16}
\Foreach{\x}{\asubset_0}
\U\cap\func{v}{\x}=\empty,
\end{equation}
and hence,
\begin{equation}\label{thmcompactsetofT2spacepointp1eq17}
\U\cap\(\union{\func{\image{v}}{\asubset_0}}\)=\empty,
\end{equation}
and thus according to \Ref{thmcompactsetofT2spacepointp1eq9},
\begin{equation}\label{thmcompactsetofT2spacepointp1eq18}
\U\cap\(\union{\acover}\)=\empty.
\end{equation}
Therefore, according to \Ref{thmcompactsetofT2spacepointp1eq13},
\Ref{thmcompactsetofT2spacepointp1eq14}, and
\Ref{thmcompactsetofT2spacepointp1eq18},
\begin{equation}
\Existsis{\opair{\U}{\(\union{\acover}\)}}
{\Cprod{\func{\nei{\Xt}}{\seta{\point}}}{\func{\nei{\Xt}}{\asubset}}}
\U\cap\(\union{\acover}\)=\empty
\end{equation}
\endp
\end{itemize}
\endthm
%%%%%%%%%%%%%%%%%%%%%%%%%%%%%%%%%%%%%%%%%%%%%%%%%%%%%%%%%%%%%%%%%%%%%%%%%%%%%%%%%%%%%%%%%%%%%%%%%%%%%%%%%%%%%%%%%%%
\theorem\label{thmcompactsetofT2spaceisclosed}
$\Xt=\opair{\X}{\topology{}}$
is taken as a $\sepA{2}$ (Hausdorff) space. Every compact set of $\Xt$
is a closed set of $\Xt$. That is,
\begin{equation}
\compacts{\Xt}\subseteq\Fclosed{\X}{\topology{}}.
\end{equation}
\prooff
It is known that $\empty$ is both a compact and closed set of $\Xt$.
So we consider non-empty sets.\\
$\acompactset{}$ is taken as an arbitrary element of $\compl{\compacts{\Xt}}{\seta{\empty}}$
(a non-empty compact set of $\Xt$). According to \refthm{thmcompactsetofT2spacepoint},
\begin{equation}
\Foreach{\point}{\(\compl{\X}{\acompactset{}}\)}
\bigg(\Exists{\opair{\U}{\V}}
{\Cprod{\func{\nei{\Xt}}{\seta{\point}}}{\func{\nei{\Xt}}{\acompactset{}}}}
\U\cap\V=\empty\bigg).
\end{equation}
Therefore, it is clear that,
\begin{equation}
\Foreach{\point}{\(\compl{\X}{\acompactset{}}\)}
\Exists{\U}{\func{\nei{\Xt}}{\seta{\point}}}
\U\cap\acompactset{}=\empty.
\end{equation}
or equivalently,
\begin{equation}
\Foreach{\point}{\(\compl{\X}{\acompactset{}}\)}
\Exists{\U}{\func{\nei{\Xt}}{\seta{\point}}}
\U\subseteq\(\compl{\X}{\acompactset{}}\).
\end{equation}
Thus, according to \refdef{definteriorpoint} and \refthm{thmintofsetissetofintpoints},
\begin{equation}
\Foreach{\point}{\(\compl{\X}{\acompactset{}}\)}
\point\in\func{\Int{\Xt}}{\compl{\X}{\acompactset{}}},
\end{equation}
\endp
and hence knowing that the interior of every set is a subset of that set (\refcor{corintofset0}),
it becomes clear that,
\begin{equation}
\(\compl{\X}{\acompactset{}}\)=\func{\Int{\Xt}}{\compl{\X}{\acompactset{}}},
\end{equation}
and thus according to \refthm{thmintofopenset},
\begin{equation}
\(\compl{\X}{\acompactset{}}\)\in\topology{},
\end{equation}
and hence according to \refdef{deffamilyofclosedsets},
\begin{equation}
\acompactset{}\in\Fclosed{\X}{\topology{}}.
\end{equation}
\endthm
%%%%%%%%%%%%%%%%%%%%%%%%%%%%%%%%%%%%%%%%%%%%%%%%%%%%%%%%%%%%%%%%%%%%%%%%%%%%%%%%%%%%%%%%%%%%%%%%%%%%%%%%%%%%%%%%%%%
\theorem\label{thmdisjoinctcompactsetsofT2space}
$\Xt=\opair{\X}{\topology{}}$
is taken as a $\sepA{2}$ (Hausdorff) space.
Every pair of non-empty compact and disjoint sets of $\Xt$ have disjoint neighborhoods in $\Xt$. That is,
\begin{align}
&\Foreach{\opair{\acompactset{1}}{\acompactset{2}}}
{\Cprod{\(\compl{\compacts{\Xt}}{\seta{\empty}}\)}{\(\compl{\compacts{\Xt}}{\seta{\empty}}\)}}\cr
&\[\bigg(\acompactset{1}\cap\acompactset{2}=\empty\bigg)
\then
\bigg(\Exists{\opair{\U}{\V}}{\Cprod{\func{\nei{\Xt}}{\acompactset{1}}}{\func{\nei{\Xt}}{\acompactset{2}}}}
\U\cap\V=\empty\bigg)\].
\end{align}
\prooff
Each
$\acompactset{1}$
and
$\acompactset{2}$
is taken as an element of $\(\compl{\compacts{\Xt}}{\seta{\empty}}\)$
(a non-empty compact set of $\Xt$) such that,
\begin{equation}\label{thmdisjoinctcompactsetsofT2spacepeq1}
\acompactset{1}\cap\acompactset{2}=\empty.
\end{equation}
Therefore,
\begin{equation}\label{thmdisjoinctcompactsetsofT2spacepeq2}
\Foreach{\point}{\acompactset{1}}
\point\in\(\compl{\X}{\acompactset{2}}\).
\end{equation}
Therefore according to \refthm{thmcompactsetofT2spacepoint},
there exists a pair of functions
$u,~\Theta$,
each from $\acompactset{1}$
to the collection of all open sets of $\Xt$ ($\topology{}$), such that for every $\x$ in $\acompactset{1}$,
$\func{u}{\x}$ and $\func{\Theta}{\x}$ are disjoint neighborhoods of
$\seta{\x}$ and $\acompactset{2}$ in $\Xt$ , respectively:
\begin{align}\label{thmdisjoinctcompactsetsofT2spacepeq3}
&\Existsis{\opair{u}{\Theta}}
{\Cprod{\Func{\acompactset{1}}{\topology{}}}{\Func{\acompactset{1}}{\topology{}}}}\cr
&\[\Foreach{\x}{\acompactset{1}}
\bigg(\func{u}{x}\in\func{\nei{\Xt}}{\seta{\x}},
\func{\Theta}{x}\in\func{\nei{\Xt}}{\acompactset{2}},
\func{u}{\x}\cap\func{\Theta}{\x}=\empty\bigg)\].
\end{align}
It is clear that,
\begin{align}
\(\union{\func{\image{u}}{\acompactset{1}}}\)&\supseteq\acompactset{1},
\label{thmdisjoinctcompactsetsofT2spacepeq4}\\
\func{\image{u}}{\acompactset{1}}&\subseteq\topology{}.
\label{thmdisjoinctcompactsetsofT2spacepeq5}
\end{align}
Therefore, according to \refdef{defrelopencover},
$\func{\image{u}}{\acompactset{1}}$ is an open covering of $\acompactset{1}$ in $\Xt$.
\begin{equation}\label{thmdisjoinctcompactsetsofT2spacepeq6}
\func{\image{u}}{\acompactset{1}}\in\func{\ocov{\Xt}}{\acompactset{1}}.
\end{equation}
Considering that $\acompactset{1}$ is a compact set of $\Xt$, and according to \refthm{thmcompactsets},
there exists a finite open covering of $\acompactset{1}$ in $\Xt$ which is a subset of $\func{\image{u}}{\acompactset{1}}$.
\begin{equation}\label{thmdisjoinctcompactsetsofT2spacepeq7}
\Existsis{\acover}{\func{\focov{\Xt}}{\acompactset{1}}}
\acover\subseteq\func{\image{u}}{\acompactset{1}}.
\end{equation}
Therefore according to \refdef{defrelopencover},
and the non-emptiness of $\acompactset{1}$,
\begin{align}
\(\union{\acover}\)&\supseteq\acompactset{1},
\label{thmdisjoinctcompactsetsofT2spacepeq8}\\
\CarD{\acover}&\in\Zp.
\label{thmdisjoinctcompactsetsofT2spacepeq9}
\end{align}
Thus there exists a finite subset $\asubset$ of $\acompactset{1}$,
with the same number of elements as that of $\acover$, and,
$\func{\image{u}}{\asubset}=\acover$.
\begin{equation}\label{thmdisjoinctcompactsetsofT2spacepeq10}
\Existssis{\asubset}{\acompactset{1}}
\bigg(\func{\image{u}}{\asubset}=\acover,~
\CarD{\asubset}=\CarD{\acover}\bigg).
\end{equation}
$\U$
is defined as the intersection of all $\func{\Theta}{x}$ when $\x$ ranges over $\asubset$.
\begin{align}\label{thmdisjoinctcompactsetsofT2spacepeq11}
\U:=\(\intersection{\func{\image{\Theta}}{\asubset}}\)
\end{align}
Therefore, considering that $\func{\image{\Theta}}{\asubset}$ is a finite collection of open sets of $\Xt$,
it is clear that,
\begin{equation}\label{thmdisjoinctcompactsetsofT2spacepeq12}
\U\in\topology{}.
\end{equation}
Moreover, since evry $\func{\Theta}{\x}$ includes $\acompactset{2}$, it is trivial that,
\begin{equation}\label{thmdisjoinctcompactsetsofT2spacepeq13}
\acompactset{2}\subseteq\U.
\end{equation}
Thus, $\U$ is a neighborhood of $\acompactset{2}$.
\begin{equation}\label{thmdisjoinctcompactsetsofT2spacepeq14}
\U\in\func{\nei{\Xt}}{\acompactset{2}}.
\end{equation}
Moreover, according to \Ref{thmdisjoinctcompactsetsofT2spacepeq7} and \refcor{corunionofopencoverisaneiofset},
$\(\union{\acover}\)$ is a neighborhood of $\acompactset{1}$ in $\Xt$.
\begin{equation}\label{thmdisjoinctcompactsetsofT2spacepeq15}
\(\union{\acover}\)\in\func{\nei{\Xt}}{\acompactset{1}}.
\end{equation}
Furthermore, according to \Ref{thmdisjoinctcompactsetsofT2spacepeq11},
\begin{equation}\label{thmdisjoinctcompactsetsofT2spacepeq16}
\Foreach{\x}{\asubset}
\U\subseteq\func{\Theta}{\x},
\end{equation}
and according to \Ref{thmdisjoinctcompactsetsofT2spacepeq3},
\begin{equation}\label{thmdisjoinctcompactsetsofT2spacepeq17}
\Foreach{\x}{\asubset}
\func{\Theta}{\x}\cap\func{u}{\x}=\empty.
\end{equation}
Therefore,
\begin{equation}\label{thmdisjoinctcompactsetsofT2spacepeq18}
\Foreach{\x}{\asubset}
\U\cap\func{u}{\x}=\empty,
\end{equation}
and hence,
\begin{equation}\label{thmdisjoinctcompactsetsofT2spacepeq19}
\U\cap\(\union{\func{\image{u}}{\asubset}}\)=\empty,
\end{equation}
and hence according to \Ref{thmdisjoinctcompactsetsofT2spacepeq10},
\begin{equation}\label{thmdisjoinctcompactsetsofT2spacepeq20}
\U\cap\(\union{\acover}\)=\empty.
\end{equation}
Therefore according to \Ref{thmdisjoinctcompactsetsofT2spacepeq14},
\Ref{thmdisjoinctcompactsetsofT2spacepeq15}, and
\Ref{thmdisjoinctcompactsetsofT2spacepeq20}
it is clear that,
\begin{equation}
\Existsis{\opair{\(\union{\acover}\)}{\U}}
{\Cprod{\func{\nei{\Xt}}{\acompactset{1}}}{\func{\nei{\Xt}}{\acompactset{2}}}}
\(\union{\acover}\)\cap\U=\empty
\end{equation}
\endthm
%%%%%%%%%%%%%%%%%%%%%%%%%%%%%%%%%%%%%%%%%%%%%%%%%%%%%%%%%%%%%%%%%%%%%%%%%%%%%%%%%%%%%%%%%%%%%%%%%%%%%%%%%%%%%%%%%%%
\theorem\label{thmdisjoinctcompactsetsofT2space}
$\Xt=\opair{\X}{\topology{}}$
is taken as a $\sepA{2}$ (Hausdorff) space.
the intersection of a non-empty collection of compact sets of $\Xt$ is a compact set of $\Xt$. That is,
\begin{equation}
\Foreach{\somecompacts{}}{\[\compl{\CSs{\compacts{\Xt}}}{\seta{\empty}}\]}
\(\intersection{\somecompacts{}}\)\in\compacts{\Xt}.
\end{equation}
\prooff
$\somecompacts{}$
is taken as an arbitrary element of $\[\compl{\CSs{\compacts{\Xt}}}{\seta{\empty}}\]$
(a non-empty collection of compact sets of the Hausdorff space $\Xt$). According to \refthm{thmcompactsetofT2spaceisclosed},
$\somecompacts{}$
is a collection of closed sets of $\Xt$.
\begin{equation}\label{thmdisjoinctcompactsetsofT2spacepeq1}
\somecompacts{}\subseteq\Fclosed{\X}{\topology{}},
\end{equation}
and hence,
\begin{equation}
\somecompacts{}\in
\[\compl{\CSs{\compacts{\Xt}\cap\Fclosed{\X}{\topology{}}}}{\seta{\empty}}\],
\end{equation}
and hence according to \refthm{thmintersectionofclosedcompactsets},
\begin{equation}
\(\intersection{\somecompacts{}}\)\in\compacts{\Xt}.
\end{equation}
\endthm
%%%%%%%%%%%%%%%%%%%%%%%%%%%%%%%%%%%%%%%%%%%%%%%%%%%%%%%%%%%%%%%%%%%%%%%%%%%%%%%%%%%%%%%%%%%%%%%%%%%%%%%%%%%%%%%%%%%
\theorem
Each $\Xt_1$ and $\Xt_2$ is taken as a topological space.
If $\Xt_1$ and $\Xt_2$ are compact, then so is $\topprod{\Xt_1}{\Xt_2}$.
\proof
It is left as an exercise.
\endthm
%%%%%%%%%%%%%%%%%%%%%%%%%%%%%%%%%%%%%%%%%%%%%%%%%%%%%%%%%%%%%%%%%%%%%%%%%%%%%%%%%%%%%%%%%%%%%%%%%%%%%%%%%%%%%%%%%%%%%%%%%%%%%%%%
\theorem
$\Xt=\opair{\X}{\topology{\X}}$ is taken as a topological space, and $\eqrel{}$ as an equivalence relation on $\X$.
If $\Xt$ is compact, then so is $\topq{\Xt}{\eqrel{}}$.
\proof
It is left as an exercise.
\endthm
%%%%%%%%%%%%%%%%%%%%%%%%%%%%%%%%%%%%%%%%%%%%%%%%%%%%%%%%%%%%%%%%%%%%%%%%%%%%%%%%%%%%%%%%%%%%%%%%%%%%%%%%%%%%%%%%%%%
\subsection{
Compactness as a Topological Property
}
\theorem\label{thmcontinuousimageofcompactspace}
Each $\Xt=\opair{\X}{\topology{\X}}$ and
$\Yt=\opair{\Y}{\topology{\Y}}$
is taken as a topological space. If $\Xt$ is compact, then the image of $\X$ under every continuous map
from $\Xt$ to $\Yt$ is a compact set of $\Yt$. That is,
\begin{equation}
\(\X\in\compacts{\Xt}\)\then
\[\Foreach{\cf}{\CF{\Xt}{\Yt}}
\func{\image{\cf}}{\X}\in\compacts{\Yt}\].
\end{equation}
\prooff
It is assumed that $\Xt$ is compact, and $\cf$is taken as an arbitrary element of $\CF{\Xt}{\Yt}$.
$\cg$
is defined as the codomain-restriction of $\cf$ to $\func{\image{\cf}}{\X}$:
\begin{equation}
\cg:=\func{\rescd{\cf}}{\func{\image{\cf}}{\X}}.
\end{equation}
It is trivial that the image of $\cg$ equals $\func{\image{\cf}}{\X}$,
and hence $\cg$ is a surjective function from $\X$ to $\func{\image{\cf}}{\X}$.
\begin{equation}
\cg\in\surFunc{\X}{\func{\image{\cf}}{\X}}.
\end{equation}
Moreover, according to \refcor{correstrictionofcontinuousfunction}, $\cg$ is a continuous map from $\Xt$ to
$\opair{\func{\image{\cf}}{\X}}{\stopology{\topology{\Y}}{\func{\image{\cf}}{\X}}}$:
\begin{equation}
\cg\in\CF{\Xt}{\opair{\func{\image{\cf}}{\X}}{\stopology{\topology{\Y}}{\func{\image{\cf}}{\X}}}},
\end{equation}
and hence according to, \refdef{defcontinuousfunction},
\begin{equation}\label{thmcontinuousimageofcompactspacepeq1}
\Foreach{\U}{\stopology{\topology{\Y}}{\func{\image{\cf}}{\X}}}
\func{\pimage{\cg}}{\U}\in\topology{\X}.
\end{equation}
\begin{itemize}
\item[${\textbf{\textsf{p1}}}$]
$\acover_{\Y}$
is taken as an arbitrary element of
$\Ocov{\opair{\func{\image{\cf}}{\X}}{\stopology{\topology{\Y}}{\func{\image{\cf}}{\X}}}}$
(an open covering of
$\opair{\func{\image{\cf}}{\X}}{\stopology{\topology{\Y}}{\func{\image{\cf}}{\X}}}$)
Then according to \refdef{defrelopencover},
\begin{align}
\(\union{\acover_{\Y}}\)&=\func{\image{\cf}}{\X},
\label{thmcontinuousimageofcompactspacep1eq1}\\
\acover_{\Y}&\subseteq\stopology{\topology{\Y}}{\func{\image{\cf}}{\X}}.
\label{thmcontinuousimageofcompactspacep1eq2}
\end{align}
The set $\acover_{\X}$ is defined as,
\begin{equation}\label{thmcontinuousimageofcompactspacep1eq3}
\acover_{\X}:=\defSet{\func{\pimage{\cg}}{\U}}{\U\in\acover_{\Y}}.
\end{equation}
According to \Ref{thmcontinuousimageofcompactspacep1eq1}, and considering the surjectivity of $\cg$,
\begin{align}\label{thmcontinuousimageofcompactspacep1eq4}
\(\union{\acover_{\X}}\)&=\Union{\U}{\acover_{\Y}}{\func{\pimage{\cg}}{\U}}\cr
%&=\func{\pimage{\cf}}{\Union{\U}{\acover}{\U}}\cr
&=\func{\pimage{\cg}}{\union{\acover_{\Y}}}\cr
&=\X
\end{align}
Moreover, according to \Ref{thmcontinuousimageofcompactspacepeq1},
\Ref{thmcontinuousimageofcompactspacep1eq2}, and
\Ref{thmcontinuousimageofcompactspacep1eq3},
$\acover_{\X}$
is a collection of open sets of $\Xt$.
\begin{equation}\label{thmcontinuousimageofcompactspacep1eq5}
\acover_{\X}\subseteq\topology{\X}.
\end{equation}
According to \refdef{defrelopencover},
\Ref{thmcontinuousimageofcompactspacep1eq4}, and
\Ref{thmcontinuousimageofcompactspacep1eq5},
$\acover_{\X}$
is an open covering of $\Xt$:
\begin{equation}\label{thmcontinuousimageofcompactspacep1eq6}
\acover_{\X}\in\Ocov{\Xt}.
\end{equation}
Therefore, considering that $\Xt$ is compact, and according to \refdef{defcompactness},
it is inferred that there exists a finite open covering of $\Xt$ which is a subset of $\acover_{\X}$:
\begin{equation}\label{thmcontinuousimageofcompactspacep1eq7}
\Existsis{\p{\acover_{\X}}}{\fOcov{\Xt}}
\p{\acover_{\X}}\subseteq\acover_{\X}.
\end{equation}
According to \refdef{defrelopencover},
\begin{align}
\(\union{\p{\acover_{\X}}}\)&=\X,
\label{thmcontinuousimageofcompactspacep1eq8}\\
\CarD{\p{\acover_{\X}}}&\in\Zpz.
\label{thmcontinuousimageofcompactspacep1eq9}
\end{align}
The set $\p{\acover_{\Y}}$ is defined as,
\begin{equation}\label{thmcontinuousimageofcompactspacep1eq10}
\p{\acover_{\Y}}:=\defSet{\func{\image{\cg}}{\V}}{\V\in\p{\acover_{\X}}}.
\end{equation}
Considering that $\p{\acover_{\X}}\subseteq\acover_{\X}$, and $\cg$
is a surjective function from $\X$ to $\func{\image{\cf}}{\X}$, and according to
\Ref{thmcontinuousimageofcompactspacep1eq3} and
\Ref{thmcontinuousimageofcompactspacep1eq10},
\begin{equation}\label{thmcontinuousimageofcompactspacep1eq11}
\p{\acover_{\Y}}\subseteq\acover_{\Y},
\end{equation}
and hence according to \Ref{thmcontinuousimageofcompactspacep1eq2},
\begin{equation}\label{thmcontinuousimageofcompactspacep1eq12}
\p{\acover_{\Y}}\subseteq\stopology{\topology{\Y}}{\func{\image{\cf}}{\X}}.
\end{equation}
The surjectivity of $\cg$,
\Ref{thmcontinuousimageofcompactspacep1eq8} and
\Ref{thmcontinuousimageofcompactspacep1eq10}
imply that,
\begin{equation}\label{thmcontinuousimageofcompactspacep1eq13}
\(\union{\p{\acover_{\Y}}}\)=\func{\image{\cf}}{\X}.
\end{equation}
Moreover, the finiteness of $\p{\acover_{\X}}$ implies the finiteness of $\p{\acover_{\Y}}$:
\begin{equation}\label{thmcontinuousimageofcompactspacep1eq14}
\CarD{\p{\acover_{\Y}}}\in\Zpz.
\end{equation}
According to \refdef{defrelopencover},
\Ref{thmcontinuousimageofcompactspacep1eq11},
\Ref{thmcontinuousimageofcompactspacep1eq12},
\Ref{thmcontinuousimageofcompactspacep1eq13}, and
\Ref{thmcontinuousimageofcompactspacep1eq14},
it is inferred that $\p{\acover_{\Y}}$
is a finite open covering of
$\opair{\func{\image{\cf}}{\X}}{\stopology{\topology{\Y}}{\func{\image{\cf}}{\X}}}$,
and a subset of $\acover_{\Y}$.
\begin{equation}
\Existsis{\p{\acover_{\Y}}}{\fOcov{\opair{\func{\image{\cf}}{\X}}{\stopology{\topology{\Y}}{\func{\image{\cf}}{\X}}}}}
\p{\acover_{\Y}}\subseteq\acover_{\Y}.
\end{equation}
Therefore,
\begin{align}
&\Foreach{\acover_{\Y}}{\Ocov{\opair{\func{\image{\cf}}{\X}}{\stopology{\topology{\Y}}{\func{\image{\cf}}{\X}}}}}\cr
&\defset{\p{\acover_{\Y}}}{\fOcov{\opair{\func{\image{\cf}}{\X}}{\stopology{\topology{\Y}}{\func{\image{\cf}}{\X}}}}}
{\p{\acover_{\Y}}\subseteq\acover_{\Y}}\neq\empty,
\end{align}
and hence according to \refdef{defcompactness},
$\opair{\func{\image{\cf}}{\X}}{\stopology{\topology{\Y}}{\func{\image{\cf}}{\X}}}$
is compact:
\begin{equation}
\func{\image{\cf}}{\X}\in
\compacts{\opair{\func{\image{\cf}}{\X}}{\stopology{\topology{\Y}}{\func{\image{\cf}}{\X}}}}.
\end{equation}
Equivalently, according to \refdef{defcompactsets},
$\func{\image{\cf}}{\X}$ is a compact set of $\Yt$:
\begin{equation}
\func{\image{\cf}}{\X}\in
\compacts{\Yt}.
%\compacts{\opair{\func{\image{\cf}}{\X}}{\stopology{\topology{\Y}}{\func{\image{\cf}}{\X}}}}.
\end{equation}
\end{itemize}
\endthm
%%%%%%%%%%%%%%%%%%%%%%%%%%%%%%%%%%%%%%%%%%%%%%%%%%%%%%%%%%%%%%%%%%%%%%%%%%%%%%%%%%%%%%%%%%%%%%%%%%%%%%%%%%%%%%%%%%%
\theorem\label{thmcontinuousimageofcompactset}
Each
$\Xt=\opair{\X}{\topology{\X}}$
and
$\Yt=\opair{\Y}{\topology{\Y}}$
is taken as a topological space.
For every continuous map $\cf$ from $\Xt$ to $\Yt$, the image under $\cf$ of
every compact set of $\Xt$ is a compact set of $\Yt$:
\begin{equation}
\Foreach{\cf}{\CF{\Xt}{\Yt}}
\[\Foreach{\acompactset{}}{\compacts{\Xt}}
\bigg(\func{\image{\cf}}{\acompactset{}}\in\compacts{\Yt}\bigg)\].
\end{equation}
\prooff
$\cf$
is taken as an arbitrary element of $\CF{\Xt}{\Yt}$, and
$\acompactset{}$ as an arbitrary element of $\compacts{\Xt}$. According to \refthm{thmcontiniuityequiv2},
$\func{\resd{\cf}}{\acompactset{}}$
is a continuous map from the topological space
$\opair{\acompactset{}}{\stopology{\topology{\X}}{\acompactset{}}}$ to $\Yt$:
\begin{equation}
\func{\resd{\cf}}{\acompactset{}}
\in\CF{\opair{\acompactset{}}{\stopology{\topology{\X}}{\acompactset{}}}}{\Yt}.
\end{equation}
Moreover, considering that $\acompactset{}\in\compacts{\Xt}$, and according to \refthm{thmcompactsetequiv0},
\begin{equation}
\acompactset{}\in\compacts{\opair{\acompactset{}}{\stopology{\topology{\X}}{\acompactset{}}}}.
\end{equation}
Therefore according to \refthm{thmcontinuousimageofcompactspace},
\begin{align}
\func{\image{\[\func{\resd{\cf}}{\acompactset{}}\]}}{\acompactset{}}\in
\compacts{\Yt}.
\end{align}
Therefore considering that,
\begin{equation}
\func{\image{\[\func{\resd{\cf}}{\acompactset{}}\]}}{\acompactset{}}=
\func{\image{\cf}}{\acompactset{}},
\end{equation}
it becomes clear that,
\begin{equation}
\func{\image{\cf}}{\acompactset{}}\in\compacts{\Yt}.
\end{equation}
\endthm
%%%%%%%%%%%%%%%%%%%%%%%%%%%%%%%%%%%%%%%%%%%%%%%%%%%%%%%%%%%%%%%%%%%%%%%%%%%%%%%%%%%%%%%%%%%%%%%%%%%%%%%%%%%%%%%%%%%%%%%%%%%%%%%%
\theorem\label{thmcompactnessisatopologicalproperty}
Each
$\Xt=\opair{\X}{\topology{\X}}$
and
$\Yt=\opair{\Y}{\topology{\Y}}$
is taken as a topological space. If $\Xt$ and $\Yt$
are homeomorphic, then $\Xt$ is compact if and only if $\Yt$ is compact. That is,
\begin{equation}
\bigg(\homeomorphic{\Xt}{\Yt}\bigg)\then
\bigg(\X\in\compacts{\Xt}\thenn\Y\in\compacts{\Yt}\bigg).
\end{equation}
$\caution$
In other words, compactness is a topological property.
\proof
It is assumed taht $\Xt$ and $\Yt$ are homeomorphic:
\begin{equation}\label{thmcompactnessisatopologicalpropertypeq1}
\homeomorphic{\Xt}{\Yt}.
\end{equation}
According to \refdef{defhomeomorphic},
this means there exists at least one homeomorphism from $\Xt$ to $\Yt$:
\begin{equation}\label{thmcompactnessisatopologicalpropertypeq2}
\HOF{\Xt}{\Yt}\neq\empty.
\end{equation}
$\hf$
is taken as an element of $\HOF{\Xt}{\Yt}$
($\HOF{\Xt}{\Yt}=\seta{\binary{\hf}{\ldots}}$).
According to \refdef{defhomeomorphism},
\begin{align}
\hf&\in\IF{\X}{\Y},
\label{thmcompactnessisatopologicalpropertypeq3}\\
\hf&\in\CF{\Xt}{\Yt},
\label{thmcompactnessisatopologicalpropertypeq4}\\
\finv{\hf}&\in\CF{\Yt}{\Xt}.
\label{thmcompactnessisatopologicalpropertypeq5}
\end{align}
According to \Ref{thmcompactnessisatopologicalpropertypeq3}
(bijectivity of $\hf$),
it is clear that,
\begin{align}
\func{\image{\hf}}{\X}&=\Y,
\label{thmcompactnessisatopologicalpropertypeq6}\\
\func{\image{\[\finv{\hf}\]}}{\Y}&=\X.
\label{thmcompactnessisatopologicalpropertypeq7}
\end{align}
According to \refthm{thmconnectednessandcontinuousfunctions},
\Ref{thmcompactnessisatopologicalpropertypeq4},
\Ref{thmcompactnessisatopologicalpropertypeq5},
\Ref{thmcompactnessisatopologicalpropertypeq6}, and
\Ref{thmcompactnessisatopologicalpropertypeq7},
\begin{equation}
\bigg(\X\in\compacts{\Xt}\thenn\Y\in\compacts{\Yt}\bigg).
\end{equation}
\endthm
%%%%%%%%%%%%%%%%%%%%%%%%%%%%%%%%%%%%%%%%%%%%%%%%%%%%%%%%%%%%%%%%%%%%%%%%%%%%%%%%%%%%%%%%%%%%%%%%%%%%%%%%%%%%%%%%%%%%%%%%%%%%%%%%
%%%%%%%%%%%%%%%%%%%%%%%%%%%%%%%%%%%%%%%%%%%%%%%%%%%%%%%%%%%%%%%%%%%%%%%%%%%%%%%%%%%%%%%%%%%%%%%%%%%%%%%%%%%%%%%%%%%%%%%%%%%%%%%%
\subsection{
Embedding of a Compact Space in a Hausdorff Space
}
\theorem\label{thmcontinuousmapfromcompacttohausdorffisclosed}
Each
$\Xt=\opair{\X}{\topology{\X}}$
and
$\Yt=\opair{\Y}{\topology{\Y}}$
is taken as a topological space. If $\Xt$ is a compact topological space, and $\Yt$
is a Hausdorff space, then every continuous map from $\Xt$ to $\Yt$
is a closed map from $\Xt$ to $\Yt$:
\begin{equation}
\(\X\in\compacts{\Xt},~\topology{\Y}\in\CtopsH{\X}\)\then
\(\CF{\Xt}{\Yt}\subseteq\CM{\Xt}{\Yt}\).
\end{equation}
\prooff
It is assumed that $\Xt$ is compact, and $\Yt$ is a $\sepA{2}$ space.
$\cf$ is taken as an element of $\CF{\Xt}{\Yt}$.
\begin{itemize}
\item[${\textbf{\textsf{p1}}}$]
$\asubset$
is taken as an arbitrary element of $\Fclosed{\X}{\topology{\X}}$
(a closed set of $\Xt$).
Then considering that $\Xt$ is compact, according to \refthm{thmclosedsetofcompactspaceiscompact},
$\asubset$ is a compact set of $\Xt$:
\begin{equation}
\asubset\in\compacts{\Xt}.
\end{equation}
Thus according to \refthm{thmcontinuousimageofcompactset}, and
considering that $\cf$ is a continuous map from $\Xt$ to $\Yt$,
$\func{\image{\cf}}{\asubset}$
is a compact set of $\Yt$.
\begin{equation}
\func{\image{\cf}}{\asubset}\in\compacts{\Yt}.
\end{equation}
Therefore, considering that $\Yt$ is a $\sepA{2}$ space, and according to \refthm{thmcompactsetofT2spaceisclosed},
it becomes clear that $\func{\image{\cf}}{\asubset}$ is a closed set of $\Yt$:
\begin{equation}
\func{\image{\cf}}{\asubset}\in\Fclosed{\Y}{\topology{\Y}}.
\end{equation}
\endp
\end{itemize}
Therefore,
\begin{equation}
\Foreach{\asubset}{\Fclosed{\X}{\topology{\X}}}
\func{\image{\cf}}{\asubset}\in\Fclosed{\Y}{\topology{\Y}},
\end{equation}
and hence according to \refdef{defclosedmap},
$\cf$ is a closed map from $\Xt$ to $\Yt$:
\begin{equation}
\cf\in\CM{\Xt}{\Yt}.
\end{equation}
\endthm
%%%%%%%%%%%%%%%%%%%%%%%%%%%%%%%%%%%%%%%%%%%%%%%%%%%%%%%%%%%%%%%%%%%%%%%%%%%%%%%%%%%%%%%%%%%%%%%%%%%%%%%%%%%%%%%%%%%%%%%%%%%%%%%%
\theorem\label{thmcontinuousbijectionfromcompacttohausforffishomeomorphism}
$\Xt=\opair{\X}{\topology{\X}}$
is taken as a compact topological space, and
$\Yt=\opair{\Y}{\topology{\Y}}$ as a $\sepA{2}$ (Hausdorff) space.
Every continuous and bijective map from $\Xt$ to $\Yt$ is a homeomorphism from $\Xt$ to $\Yt$.
\begin{equation}
\(\CF{\Xt}{\Yt}\cap\IF{\X}{\Y}\)=\HOF{\Xt}{\Yt}.
\end{equation}
\prooff
According to \refthm{thmcontinuousmapfromcompacttohausdorffisclosed},
$\CF{\Xt}{\Yt}\subseteq\CM{\Xt}{\Yt}$,
and hence according to \refthm{thmhomeomorphismsisclosedandcontinuousmap},
\begin{align}
\(\CF{\Xt}{\Yt}\cap\IF{\X}{\Y}\)&=
\(\CF{\Xt}{\Yt}\cap\CM{\Xt}{\Yt}\cap\IF{\Xt}{\Yt}\)\cr
&=\HOF{\Xt}{\Yt}.
\end{align}
\endthm
%%%%%%%%%%%%%%%%%%%%%%%%%%%%%%%%%%%%%%%%%%%%%%%%%%%%%%%%%%%%%%%%%%%%%%%%%%%%%%%%%%%%%%%%%%%%%%%%%%%%%%%%%%%%%%%%%%%%%%%%%%%%%%%%
\theorem\label{thmcontinuousinjectionfromcompacttohausforffisembedding}
$\Xt=\opair{\X}{\topology{\X}}$
is taken as a compact topological space, and $\Yt=\opair{\Y}{\topology{\Y}}$
as a $\sepA{2}$ (Hausdorff) space.
Every continuous and injective map from $\Xt$ to $\Yt$ is an embedding of $\Xt$ in $\Yt$.
\begin{equation}
\(\CF{\Xt}{\Yt}\cap\InF{\Xt}{\Yt}\)=\EM{\Xt}{\Yt}.
\end{equation}
\prooff
$\cf$
is taken as an element of $\(\CF{\Xt}{\Yt}\cap\InF{\Xt}{\Yt}\)$
(a continuous and injective map from $\Xt$ to $\Yt$). Then according to \refcor{correstrictionofcontinuousfunction},
the restriction of $\cf$ to $\func{\image{\cf}}{\X}$ is a continuous map from $\Xt$ to
$\opair{\func{\image{\cf}}{\X}}{\stopology{\topology{\Y}}{\func{\image{\cf}}{\X}}}$:
\begin{equation}
\func{\rescd{\cf}}{\func{\image{\cf}}{\X}}\in
\CF{\Xt}{\opair{\func{\image{\cf}}{\X}}{\stopology{\topology{\Y}}{\func{\image{\cf}}{\X}}}}.
\end{equation}
Moreover, the injectivity of $\cf$ clearly implies the injectivity of $\func{\image{\cf}}{\X}$.
\begin{equation}
\func{\rescd{\cf}}{\func{\image{\cf}}{\X}}\in\InF{\X}{\func{\image{\cf}}{\X}}.
\end{equation}
Moreover, it is trivial that $\func{\rescd{\cf}}{\func{\image{\cf}}{\X}}$ is a surjective function.
\begin{equation}
\func{\rescd{\cf}}{\func{\image{\cf}}{\X}}\in\surFunc{\X}{\func{\image{\cf}}{\X}}.
\end{equation}
Therefore,
$\func{\rescd{\cf}}{\func{\image{\cf}}{\X}}$
is an bijective and continuous map from $\Xt$ to
$\opair{\func{\image{\cf}}{\X}}{\stopology{\topology{\Y}}{\func{\image{\cf}}{\X}}}$:
\begin{equation}
\func{\rescd{\cf}}{\func{\image{\cf}}{\X}}
\in\[\CF{\Xt}{\opair{\func{\image{\cf}}{\X}}{\stopology{\topology{\Y}}{\func{\image{\cf}}{\X}}}}
\cap\IF{\X}{\func{\image{\cf}}{\X}}\].
\end{equation}
Considering that $\Yt$ is a Hausdorff space, it is clear that,
$\opair{\func{\image{\cf}}{\X}}{\stopology{\topology{\Y}}{\func{\image{\cf}}{\X}}}$
is also a Hausdorff space, and hence according to \refthm{thmcontinuousbijectionfromcompacttohausforffishomeomorphism},
\begin{equation}
\[\CF{\Xt}{\opair{\func{\image{\cf}}{\X}}{\stopology{\topology{\Y}}{\func{\image{\cf}}{\X}}}}
\cap\IF{\X}{\func{\image{\cf}}{\X}}\]=
\HOF{\Xt}{\opair{\func{\image{\cf}}{\X}}{\stopology{\topology{\Y}}{\func{\image{\cf}}{\X}}}}.
\end{equation}
Therefore,
\begin{equation}
\func{\rescd{\cf}}{\func{\image{\cf}}{\X}}\in
\HOF{\Xt}{\opair{\func{\image{\cf}}{\X}}{\stopology{\topology{\Y}}{\func{\image{\cf}}{\X}}}}.
\end{equation}
Therefore according to \refdef{defembedding}, $\cf$ is an embedding of $\Xt$ in $\Yt$:
\begin{equation}
\cf\in\EM{\Xt}{\Yt}.
\end{equation}
\endthm
%%%%%%%%%%%%%%%%%%%%%%%%%%%%%%%%%%%%%%%%%%%%%%%%%%%%%%%%%%%%%%%%%%%%%%%%%%%%%%%%%%%%%%%%%%%%%%%%%%%%%%%%%%%%%%%%%%%%%%%%%%%%%%%%
%%%%%%%%%%%%%%%%%%%%%%%%%%%%%%%%%%%%%%%%%%%%%%%%%%%%%%%%%%%%%%%%%%%%%%%%%%%%%%%%%%%%%%%%%%%%%%%%%%%%%%%%%%%%%%%%%%%%%%%%%%%%%%%%
%%%%%%%%%%%%%%%%%%%%%%%%%%%%%%%%%%%%%%%%%%%%%%%%%%%%%%%%%%%%%%%%%%%%%%%%%%%%%%%%%%%%%%%%%%%%%%%%%%%%%%%%%%%%%%%%%%%%%%%%%%%%%%%%
%%%%%%%%%%%%%%%%%%%%%%%%%%%%%%%%%%%%%%%%%%%%%%%%%%%%%%%%%%%%%%%%%%%%%%%%%%%%%%%%%%%%%%%%%%%%%%%%%%%%%%%%%%%%%%%%%%%%%%%%%%%%%%%%
\section{
One-Point Compactification (Alexandroff Extension)
}
\definition\label{defalexandrofftopology}
$\X$
is taken as a set.
$\alexT{\X}$
is defined as,
\begin{align}
&\alexT{\X}\indef\Func{\Ctops{\X}}{\CSs{\X\cup\seta{\X}}},\cr\cr
&\Foreach{\topology{}}{\Ctops{\X}}\cr
&\func{\alexT{\X}}{\topology{}}\eqdef
\topology{}\cup\defset{\asubset}{\CSs{\X\cup\seta{\X}}}
{\[\Exists{\acompactset{}}{\bigg(\compacts{\opair{\X}{\topology{}}}\cap\Fclosed{\X}{\topology{}}\bigg)}
\asubset=\(\compl{\X}{\acompactset{}}\)\cup\seta{\X}\]}.\cr
&{}
\end{align}
Equivalently, according to \refdef{deffamilyofclosedsets},
\begin{align}
&\Foreach{\topology{}}{\Ctops{\X}}\cr
&\func{\alexT{\X}}{\topology{}}\eqdef\
\topology{}\cup\defset{\asubset}{\CSs{\X\cup\seta{\X}}}
{\[\Exists{\U}{\topology{}}
\bigg(\asubset=\U\cup\seta{\X},~\(\compl{\X}{\U}\)\in\compacts{\opair{\X}{\topology{}}}\bigg)\]}.\cr
&{}
\end{align}
In other words, $\alexT{\X}$
is defined as the function from the collection of all topologies on $\X$ ($\Ctops{\X}$)
to the power set of $\X\cup\seta{\X}$, so that for every topology $\topology{}$ on $\X$,
$\func{\alexT{\X}}{\topology{}}$
equals the union of $\topology{}$ and the collection of all subsets of $\X\cup\seta{\X}$
which are the union of $\seta{\X}$ and the complement of a compact and closed set of $\opair{\X}{\topology{}}$:
\begin{align}
\Foreach{\topology{}}{\Ctops{\X}}
\func{\alexT{\X}}{\topology{}}&\eqdef
\topology{}\cup\defSet{\(\compl{\X}{\acompactset{}}\)\cup\seta{\X}}
{\acompactset{}\in\bigg(\compacts{\opair{\X}{\topology{}}}\cap\Fclosed{\X}{\topology{}}\bigg)}\cr
&=\topology{}\cup\defSet{\U\cup\seta{\X}}{\[\U\in\topology{},~\(\compl{\X}{\U}\)\in\compacts{\opair{\X}{\topology{}}}\]}.
\end{align}
\endef
%%%%%%%%%%%%%%%%%%%%%%%%%%%%%%%%%%%%%%%%%%%%%%%%%%%%%%%%%%%%%%%%%%%%%%%%%%%%%%%%%%%%%%%%%%%%%%%%%%%%%%%%%%%%%%%%%%%%%%%%%%%%%%%%
\theorem\label{thmalexandrofftopology}
$\X$
is taken as a set.
For every topology $\topology{}$ on $\X$,
$\func{\alexT{\X}}{\topology{}}$
is a topology on $\X\cup\seta{\X}$:
\begin{equation}
\Foreach{\topology{}}{\Ctops{\X}}
\func{\alexT{\X}}{\topology{}}\in\Ctops{\X\cup\seta{\X}}.
\end{equation}
In other words, for every $\topology{}$ in $\Ctops{\X}$,
$\opair{\X\cup\seta{\X}}{\func{\alexT{\X}}{\topology{}}}$
is a topological space.
\proof
$\topology{}$
is taken as an arbitrary element of $\Ctops{\X}$
(a topology on $\X$). Then according to \refdef{deftopologicalspace},
\begin{align}
&\empty\in\topology{},
\label{thmalexandrofftopologypeq1}\\
&\X\in\topology{},
\label{thmalexandrofftopologypeq2}\\
&\Foreach{\opair{\U_{1}}{\U_{2}}}{\(\topology{}\times\topology{}\)}\U_{1}\cap\U_{2}\in\topology{},
\label{thmalexandrofftopologypeq3}\\
&\Foreachs{\sC}{\topology{}}
\(\union{\sC}\)\in\topology{}.
\label{thmalexandrofftopologypeq4}
\end{align}
According to \refdef{defalexandrofftopology} and
\Ref{thmalexandrofftopologypeq1},
\begin{equation}\label{thmalexandrofftopologypeq5}
\empty\in\func{\alexT{\X}}{\topology{}}.
\end{equation}
According to \refdef{defalexandrofftopology} and
\Ref{thmalexandrofftopologypeq2}, and considering that
\begin{equation}\label{thmalexandrofftopologypeq6}
\compl{\X}{\X}=\empty\in\compacts{\Xt},
\end{equation}
it becomes clear that,
\begin{equation}\label{thmalexandrofftopologypeq7}
\(\X\cap\seta{\X}\)\in\func{\alexT{\X}}{\topology{}}.
\end{equation}
The set $\tau$ is defined as,
\begin{equation}\label{thmalexandrofftopologypeq8}
\tau:=
\defSet{\U\cup\seta{\X}}{\[\U\in\topology{},~\(\compl{\X}{\U}\)\in\compacts{\opair{\X}{\topology{}}}\]}.
\end{equation}
\begin{itemize}
\item[${\textbf{\textsf{p1}}}$]
Each $\asubset_1$ and $\asubset_2$
is taken as an arbitrary element of $\func{\alexT{\X}}{\topology{}}$. Then according to \refdef{defalexandrofftopology},
\begin{align}\label{thmalexandrofftopologyp1eq1}
&\OR{\(\asubset_1\in\topology{}\)}{\(\asubset_1\in\tau\)},\\
&\OR{\(\asubset_2\in\topology{}\)}{\(\asubset_2\in\tau\)}.
\end{align}
According to \Ref{thmalexandrofftopologypeq3},
\begin{equation}\label{thmalexandrofftopologyp1eq2}
\bigg(\asubset_1\in\topology{},~\asubset_2\in\topology{}\bigg)\then
\bigg(\asubset_1\cap\asubset_2\in\topology{}\bigg).
\end{equation}
\begin{itemize}
\item[${\textbf{\textsf{p1-1}}}$]
It is assumed that $\asubset_1\in\tau$ and $\asubset_2\in\tau$.
Then,
\begin{align}
&\Existsis{\U_1}{\topology{}}\asubset_1=\U_1\cup\seta{\X},~\(\compl{\X}{\U_1}\)\in\compacts{\opair{\X}{\topology{}}},\\
&\Existsis{\U_2}{\topology{}}\asubset_2=\U_2\cup\seta{\X},~\(\compl{\X}{\U_2}\)\in\compacts{\opair{\X}{\topology{}}}.
\end{align}
Therefore, according to \refthm{thmunionofcompactsets},
\begin{equation}
\asubset_1\cap\asubset_2=\(\U_1\cap\U_2\)\cap\seta{\X},
\end{equation}
and
\begin{align}
\compl{\X}{\(\U_1\cap\U_2\)}&=
\(\compl{\X}{\U_1}\)\cup\(\compl{\X}{\U_2}\)\cr
&\in\compacts{\opair{\X}{\topology{}}}.
\end{align}
Moreover according to \Ref{thmalexandrofftopologypeq3},
\begin{equation}
\(\U_1\cap\U_2\)\in\topology{}.
\end{equation}
Therefore,
\begin{equation}
\(\asubset_1\cap\asubset_2\)\in\tau.
\end{equation}
\endp
\end{itemize}
Thus,
\begin{equation}\label{thmalexandrofftopologyp1eq3}
\bigg(\asubset_1\in\topology{},~\asubset_2\in\tau\bigg)\then
\bigg(\asubset_1\cap\asubset_2\in\tau\bigg).
\end{equation}
Moreover, it is clear that,
\begin{align}\label{thmalexandrofftopologyp1eq4}
\bigg(\asubset_1\in\topology{},~\asubset_2\in\tau\bigg)&\then
\bigg(\asubset_1\cap\asubset_2\in\topology{}\bigg),\\
\bigg(\asubset_1\in\tau,~\asubset_2\in\topology{}\bigg)&\then
\bigg(\asubset_1\cap\asubset_2\in\topology{}\bigg).
\end{align}
Therefore, considering that
$\func{\alexT{\X}}{\topology{}}=\topology{}\cup\tau$,
\begin{equation}
\(\asubset_1\cap\asubset_2\)\in\func{\alexT{\X}}{\topology{}}.
\end{equation}
\endp
\end{itemize}
Therefore,
\begin{equation}\label{thmalexandrofftopologypeq9}
\Foreach{\opair{\asubset_1}{\asubset_2}}
{\(\Cprod{\func{\alexT{\X}}{\topology{}}}{\func{\alexT{\X}}{\topology{}}}\)}
\(\asubset_1\cap\asubset_2\)\in\func{\alexT{\X}}{\topology{}}.
\end{equation}
\begin{itemize}
\item[${\textbf{\textsf{p2}}}$]
$\sCi$
is taken as a subset of $\func{\alexT{\X}}{\topology{}}$. Considering that
$\func{\alexT{\X}}{\topology{}}=\topology{}\cup\tau$, it is clear that,
\begin{equation}\label{thmalexandrofftopologyp2eq1}
\sCi=\(\sCi\cap\topology{}\)\cup\(\sCi\cap\tau\),
\end{equation}
and hence,
\begin{equation}\label{thmalexandrofftopologyp2eq2}
\(\union{\sCi}\)=\[\union{\(\sCi\cap\topology{}\)}\]\cup
\[\union{\(\sCi\cap\tau\)}\].
\end{equation}
According to \Ref{thmalexandrofftopologypeq4},
\begin{equation}\label{thmalexandrofftopologyp2eq3}
\[\union{\(\sCi\cap\topology{}\)}\]\in\topology{}.
\end{equation}
According to \Ref{thmalexandrofftopologypeq1}, it is clear that,
\begin{align}\label{thmalexandrofftopologyp2eq4}
\(\sCi\cap\tau\)=\empty
\then
\[\union{\(\sCi\cap\tau\)}\]\in\topology{}.
\end{align}
\begin{itemize}
\item[${\textbf{\textsf{p2-1}}}$]
It is assumed that $\(\sCi\cap\tau\)\neq\empty$.
Then,
\begin{align}\label{thmalexandrofftopologyp21eq1}
\[\union{\(\sCi\cap\tau\)}\]&=\seta{\X}\cup
\(\union{\defset{\U}{\topology{}}{\(\U\cup{\X}\)\in\(\sCi\cap\tau\)}}\).
\end{align}
According to \Ref{thmalexandrofftopologypeq4}, it is clear that,
\begin{equation}\label{thmalexandrofftopologyp21eq2}
\(\union{\defset{\U}{\topology{}}{\(\U\cup{\X}\)\in\(\sCi\cap\tau\)}}\)\in\topology{}.
\end{equation}
Considering that $\(\sCi\cap\tau\)\subseteq\tau$, it is clear that,
\begin{equation}\label{thmalexandrofftopologyp21eq3}
\Foreach{\V}{\defset{\U}{\topology{}}{\(\U\cup{\X}\)\in\(\sCi\cap\tau\)}}
\(\compl{\X}{\V}\)\in\(\compacts{\opair{\X}{\topology{}}}\cap\Fclosed{\X}{\topology{}}\),
\end{equation}
and hence according to \refthm{thmintersectionofclosedcompactsets},
and by defining
\begin{equation}\label{thmalexandrofftopologyp21eq4}
O:=\defset{\U}{\topology{}}{\(\U\cup{\X}\)\in\(\sCi\cap\tau\)},
\end{equation}
it is clear that,
\begin{align}\label{thmalexandrofftopologyp21eq5}
\compl{\X}{\(\union{\defset{\U}{\topology{}}{\(\U\cup{\X}\)\in\(\sCi\cap\tau\)}}\)}&=
\Intersection{\V}{O}{\(\compl{\X}{\V}\)}\cr
&\in\compacts{\opair{\X}{\topology{}}}.
\end{align}
\Ref{thmalexandrofftopologyp21eq1},
\Ref{thmalexandrofftopologyp21eq2}, and
\Ref{thmalexandrofftopologyp21eq5}
imply that,
\begin{equation}
\[\union{\(\sCi\cap\tau\)}\]\in\tau.
\end{equation}
\endp
\end{itemize}
Therefore,
\begin{equation}\label{thmalexandrofftopologyp2eq5}
\(\sCi\cap\tau\)\neq\empty\then
\[\union{\(\sCi\cap\tau\)}\]\in\tau.
\end{equation}
Considering that
$\func{\alexT{\X}}{\topology{}}=\topology{}\cup\tau$, and according to
\Ref{thmalexandrofftopologyp2eq2},
\Ref{thmalexandrofftopologyp2eq3},
\Ref{thmalexandrofftopologyp2eq4}, and
\Ref{thmalexandrofftopologyp2eq5},
it is inferred that,
\begin{equation}
\(\union{\sCi}\)\in\func{\alexT{\X}}{\topology{}}.
\end{equation}
\endp
\end{itemize}
Therefore,
\begin{equation}\label{thmalexandrofftopologypeq10}
\Foreachs{\sCi}{\func{\alexT{\X}}{\topology{}}}
\(\union{\sCi}\)\in\func{\alexT{\X}}{\topology{}}.
\end{equation}
According to \refdef{deftopologicalspace},
\Ref{thmalexandrofftopologypeq5},
\Ref{thmalexandrofftopologypeq7},
\Ref{thmalexandrofftopologypeq9}, and
\Ref{thmalexandrofftopologypeq10},
\begin{equation}
\func{\alexT{\X}}{\topology{}}\in\Ctops{\X\cup\seta{\X}}.
\end{equation}
\endthm
%%%%%%%%%%%%%%%%%%%%%%%%%%%%%%%%%%%%%%%%%%%%%%%%%%%%%%%%%%%%%%%%%%%%%%%%%%%%%%%%%%%%%%%%%%%%%%%%%%%%%%%%%%%%%%%%%%%%%%%%%%%%%%%%
\definition\label{defalexandroffextension}
$\X$
is taken as a set.
For every $\topology{}$ in $\Ctops{\X}$,
the topological space\\
$\opair{\X\cup\seta{\X}}{\func{\alexT{\X}}{\topology{}}}$
is referred to as the $\quotl$Alexandroff extension of the topological space $\opair{\X}{\topology{}}$$\quotr$, and
\begin{equation}
\Foreach{\topology{}}{\Ctops{\X}}
\alexTS{\opair{\X}{\topology{}}}:=
\opair{\X\cup\seta{\X}}{\func{\alexT{\X}}{\topology{}}}.
\end{equation}
\endef
%%%%%%%%%%%%%%%%%%%%%%%%%%%%%%%%%%%%%%%%%%%%%%%%%%%%%%%%%%%%%%%%%%%%%%%%%%%%%%%%%%%%%%%%%%%%%%%%%%%%%%%%%%%%%%%%%%%%%%%%%%%%%%%%
\definition\label{thmalexandrofftopologysubspace}
$\X$
is taken as a set. For every $\topology{}$ in $\Ctops{\X}$,
$\opair{\X}{\topology{}}$ is a topological subspace of
$\opair{\x\cup\seta{\X}}{\func{\alexT{\X}}{\topology{}}}$. That is,
\begin{equation}
\Foreach{\topology{}}{\Ctops{\X}}
\stopology{\func{\alexT{\X}}{\topology{}}}{\X}=\topology{}.
\end{equation}
\prooff
According to \refdef{defalexandrofftopology}, it is clear.
\endthm
%%%%%%%%%%%%%%%%%%%%%%%%%%%%%%%%%%%%%%%%%%%%%%%%%%%%%%%%%%%%%%%%%%%%%%%%%%%%%%%%%%%%%%%%%%%%%%%%%%%%%%%%%%%%%%%%%%%%%%%%%%%%%%%%
\theorem\label{thmalexandrofftopologyiscompact}
$\Xt=\opair{\X}{\topology{}}$
is taken as a topological space.
The topological space\\ $\opair{\X\cup\seta{\X}}{\func{\alexT{\X}}{\topology{}}}$ is compact.
\begin{equation}
\(\X\cup\seta{\X}\)\in\compacts{\opair{\X\cup\seta{\X}}{\func{\alexT{\X}}{\topology{}}}}.
\end{equation}
\prooff
\begin{itemize}
\item[${\textbf{\textsf{p1}}}$]
$\acover$
is taken as an element of $\Ocov{\opair{\X\cup\seta{\X}}{\func{\alexT{\X}}{\topology{}}}}$. Then according to
\refdef{defrelopencover},
\begin{align}
\(\union{\acover}\)&=\X\cup\seta{\X},
\label{thmalexandrofftopologyiscompactpeq1}\\
\acover&\subseteq\func{\alexT{\X}}{\topology{}}.
\label{thmalexandrofftopologyiscompactpeq2}
\end{align}
According to \Ref{thmalexandrofftopologyiscompactpeq1},
\begin{equation}\label{thmalexandrofftopologyiscompactpeq3}
\Existsis{\V}{\acover}
\X\in\V.
\end{equation}
According to \Ref{thmalexandrofftopologyiscompactpeq2} and
\Ref{thmalexandrofftopologyiscompactpeq3},
$\V$
is a neighborhood of $\seta{\X}$ in
$\opair{\X\cup\seta{\X}}{\func{\alexT{\X}}{\topology{}}}$:
\begin{equation}\label{thmalexandrofftopologyiscompactpeq4}
\X\in\V\in\func{\alexT{\X}}{\topology{}}.
\end{equation}
Thus according to \refdef{defalexandrofftopology},
\begin{equation}\label{thmalexandrofftopologyiscompactpeq5}
\Existsis{\U}{\topology{}}
\V=\U\cup\seta{\X},~\(\compl{\X}{\U}\)\in\compacts{\Xt}.
\end{equation}
Thus considering that $\(\compl{\X}{\U}\)$ is a compact set of $\Xt$, and according to
\refthm{thmalexandrofftopologysubspace},
and the fact that $\Xt$ is a topological subspace of $\opair{\X\cup\seta{\X}}{\func{\alexT{\X}}{\topology{}}}$,
based on \refthm{thmcompactsubspaceofcompactset} it is inferred that
$\(\compl{\X}{\U}\)$
is a compact set of the topological space
$\opair{\X\cup\seta{\X}}{\func{\alexT{\X}}{\topology{}}}$:
\begin{equation}\label{thmalexandrofftopologyiscompactpeq6}
\(\compl{\X}{\U}\)\in
\compacts{\opair{\X\cup\seta{\X}}{\func{\alexT{\X}}{\topology{}}}}.
\end{equation}
Moreover, considering that $\acover$ is an open covering of
$\opair{\X\cup\seta{\X}}{\func{\alexT{\X}}{\topology{}}}$, and according to \refdef{defrelopencover},
it is clear that it is also an open covering of $\(\compl{\X}{\U}\)$ in
$\opair{\X\cup\seta{\X}}{\func{\alexT{\X}}{\topology{}}}$.
\begin{equation}\label{thmalexandrofftopologyiscompactpeq7}
\acover\in\func{\ocov{\opair{\X\cup\seta{\X}}{\func{\alexT{\X}}{\topology{}}}}}{\compl{\X}{\U}}.
\end{equation}
According to \refthm{thmcompactsets},
\Ref{thmalexandrofftopologyiscompactpeq6}, and
\Ref{thmalexandrofftopologyiscompactpeq7},
it is inferred that there exists a finite open covering of $\(\compl{\X}{\U}\)$ in the topological space
$\opair{\X\cup\seta{\X}}{\func{\alexT{\X}}{\topology{}}}$ which is a subset of $\acover$.
\begin{equation}\label{thmalexandrofftopologyiscompactpeq8}
\Existsis{\p{\acover}}{\func{\focov{\opair{\X\cup\seta{\X}}{\func{\alexT{\X}}{\topology{}}}}}{\compl{\X}{\U}}}
\p{\acover}\subseteq\acover.
\end{equation}
Therefore according to \refdef{defrelopencover},
$\p{\acover}$
is a finite collection of open sets of
$\opair{\X\cup\seta{\X}}{\func{\alexT{\X}}{\topology{}}}$
that covers $\(\compl{\X}{\U}\)$.
\begin{align}
\p{\acover}&\subseteq\func{\alexT{\X}}{\topology{}},
\label{thmalexandrofftopologyiscompactpeq9}\\
\(\union{\p{\acover}}\)&\supseteq\(\compl{\X}{\U}\),
\label{thmalexandrofftopologyiscompactpeq10}\\
\CarD{\p{\acover}}&\in\Zpz.
\label{thmalexandrofftopologyiscompactpeq11}
\end{align}
\Ref{thmalexandrofftopologyiscompactpeq5},
\Ref{thmalexandrofftopologyiscompactpeq9},
and
\Ref{thmalexandrofftopologyiscompactpeq10}
imply that,
\begin{equation}\label{thmalexandrofftopologyiscompactpeq12}
\[\union{\(\p{\acover}\cup\seta{\V}\)}\]=\X\cup\seta{\X}.
\end{equation}
Moreover, considering that $\p{\acover}\subseteq\acover$, and
$\V\in\acover$, it is clear that,
\begin{equation}\label{thmalexandrofftopologyiscompactpeq13}
\(\p{\acover}\cup\seta{\V}\)\subseteq\acover,
\end{equation}
and hence according to \Ref{thmalexandrofftopologyiscompactpeq2},
\begin{equation}\label{thmalexandrofftopologyiscompactpeq14}
\(\p{\acover}\cup\seta{\V}\)\subseteq\func{\alexT{\X}}{\topology{}}.
\end{equation}
Moreover, according to \Ref{thmalexandrofftopologyiscompactpeq11},
it is trivial that,
\begin{equation}\label{thmalexandrofftopologyiscompactpeq15}
\CarD{\p{\acover}\cup\seta{\V}}=\CarD{\p{\acover}}+1\in\Zpz.
\end{equation}
According to \refdef{defrelopencover},
\Ref{thmalexandrofftopologyiscompactpeq12},
\Ref{thmalexandrofftopologyiscompactpeq13},
\Ref{thmalexandrofftopologyiscompactpeq14}, and
\Ref{thmalexandrofftopologyiscompactpeq15},
it is inferred that $\p{\acover}\cup\seta{\V}$ is a finite open covering of
$\opair{\X\cup\seta{\X}}{\func{\alexT{\X}}{\topology{}}}$ and a subset of $\acover$.
\begin{equation}
\Existsis{\(\p{\acover}\cup\seta{\V}\)}
{\fOcov{\opair{\X\cup\seta{\X}}{\func{\alexT{\X}}{\topology{}}}}}
\(\p{\acover}\cup\seta{\V}\)\subseteq\acover.
\end{equation}
\endp
\end{itemize}
Therefore,
\begin{equation}
\Foreach{\acover}{\Ocov{\opair{\X\cup\seta{\X}}{\func{\alexT{\X}}{\topology{}}}}}
\Exists{\acover_0}{\fOcov{\opair{\X\cup\seta{\X}}{\func{\alexT{\X}}{\topology{}}}}}
\acover_0\subseteq\acover,
\end{equation}
and hence according to \refdef{defcompactness},
$\opair{\X\cup\seta{\X}}{\func{\alexT{\X}}{\topology{}}}$ is a compact topological space.
Equivalently, according to \refdef{defcompactsets},
\begin{equation}
\(\X\cup\seta{\X}\)\in\compacts{\opair{\X\cup\seta{\X}}{\func{\alexT{\X}}{\topology{}}}}.
\end{equation}
\endthm
%%%%%%%%%%%%%%%%%%%%%%%%%%%%%%%%%%%%%%%%%%%%%%%%%%%%%%%%%%%%%%%%%%%%%%%%%%%%%%%%%%%%%%%%%%%%%%%%%%%%%%%%%%%%%%%%%%%%%%%%%%%%%%%%
\theorem\label{thmalexandrofftopologyiscompact}
$\Xt=\opair{\X}{\topology{}}$
is taken as a topological space.
The inclusion map of $\X$ into $\X\cup\seta{\X}$ is an embedding of $\Xt$ in $\alexTS{\Xt}$:
\begin{equation}
\(\incf{\X}{\X\cup\seta{\X}}\)\in
\EM{\opair{\X}{\topology{}}}{\opair{\X\cup\seta{\X}}{\func{\alexT{\X}}{\topology{}}}}
\end{equation}
\prooff
Considering that the codomain restriction of the function $\incf{\X}{\X\cup\seta{\X}}$
to its image equals the identity map on $\X$, and according to \refthm{thmtrivialhomeomorphismofaspace},
\refdef{defembedding}, and
\refthm{thmalexandrofftopologysubspace},
it is trivial.
\endthm
%%%%%%%%%%%%%%%%%%%%%%%%%%%%%%%%%%%%%%%%%%%%%%%%%%%%%%%%%%%%%%%%%%%%%%%%%%%%%%%%%%%%%%%%%%%%%%%%%%%%%%%%%%%%%%%%%%%%%%%%%%%%%%%%
%%%%%%%%%%%%%%%%%%%%%%%%%%%%%%%%%%%%%%%%%%%%%%%%%%%%%%%%%%%%%%%%%%%%%%%%%%%%%%%%%%%%%%%%%%%%%%%%%%%%%%%%%%%%%%%%%%%%%%%%%%%%%%%%
%%%%%%%%%%%%%%%%%%%%%%%%%%%%%%%%%%%%%%%%%%%%%%%%%%%%%%%%%%%%%%%%%%%%%%%%%%%%%%%%%%%%%%%%%%%%%%%%%%%%%%%%%%%%%%%%%%%%%%%%%%%%%%%%
\section{
Local Compactness
}
\definition\label{deflocalcompactness}
$\Xt=\opair{\X}{\topology{}}$
is taken as a topological space.
$\Xt$
is referred to as a $\quotl$locally compact topological space$\quotr$ iff
every point  $\x$ of $\Xt$ has a neighborhood in $\Xt$ that is a subset of
at least one compact set of $\Xt$. That is,
$\Xt$ is called a locally compact topological space, iff
\begin{equation}
\Foreach{\x}{\X}
\[\Exists{\opair{\U}{\acompactset{}}}{\(\Cprod{\func{\nei{\Xt}}{\seta{\x}}}{\compacts{\Xt}}\)}
\U\subseteq\acompactset{}\].
\end{equation}
\endef
\chapteR{A Taste of Topology in Functional Analysis}
\thispagestyle{fancy}
\section{Limit of a Topological Map}
Each $\Xt=\opair{\X}{\topology{\X}}$ and $\Yt=\opair{\Y}{\topology{\Y}}$ is
taken as a topological-space. $\asubset$ is taken as a subset of $\X$,
$\point$ an element of $\Lim{\Xt}{\asubset}$ (a limit-point of $\asubset$ in $\Xt$), and
$\cf$ is an element of $\Func{\asubset}{\Y}$. Since $\point$ is a limit-point of $\asubset$ in $\Xt$,
\begin{equation}\label{limitpointproperty}
\Foreach{\U}{\func{\nei{\Xt}}{\seta{\point}}}
\compl{\[\U\cap\asubset\]}{\seta{\point}}\neq\empty.
\end{equation}
%%%%%%%%%%%%%%%%%%%%%%%%%%%%%%%%%%%%%%%%%%%%%%%%%%%%%%%%%%%%%%%%%%%%%%%%%%%%%%%%%%%%%%%%%%%%%%%%%%%%%%%%%%
\definition\label{deflimitoftopologicalmap}
%and $\p{\point}$ an element of $\Y$.
\begin{align}
\Limits{\Xt}{\Yt}{\cf}{\point}:=
\defset{\p{\point}}{\Y}{\bigg[
\Foreach{\p{\U}}{\func{\nei{\Yt}}{\seta{\p{\point}}}}
\Exists{\U}{\func{\nei{\Xt}}{\seta{\point}}}
\func{\image{\cf}}{\compl{\[\U\cap\asubset\]}{\seta{\point}}}\subseteq\p{\U}\bigg]}.
\end{align}
Each element of $\Limits{\Xt}{\Yt}{\cf}{\point}$ is referred to as a
$\quotl$$\opair{\Xt}{\Yt}$-limit of $\cf$ at $\point$$\quotr$, or simply as a
$\quotl$limit of $\cf$ at $\point$$\quotr$.
\endef
%%%%%%%%%%%%%%%%%%%%%%%%%%%%%%%%%%%%%%%%%%%%%%%%%%%%%%%%%%%%%%%%%%%%%%%%%%%%%%%%%%%%%%%%%%%%%%%%%%%%%%%%%%
\theorem\label{thmuniquelimit}
If $\Yt$ is a Hausdorff-space, then
\begin{equation}
\Card{\Limits{\Xt}{\Yt}{\cf}{\point}}\leq1.
\end{equation}
In other words, if $\Yt$ is a Hausdorff-space, and if there exists a limit of $\cf$ at $\point$,
then this limit is unique.
\proof
It is assumed that $\Yt$ is a Hasudorrf-space. This means, for every pair of distinct points
of $\Yt$, there exists a pair of disjoint open sets of $\Yt$, each being a neighborhood of
one of those points.
\begin{equation}\label{thmuniquelimitpeq1}
\Foreach{\opair{\point_1}{\point_2}}{\psCprod{\Y}{\Y}}
\Existsis{\opair{\U_1}{\U_2}}{\Cprod{\func{\nei{\Xt}}{\seta{\point_1}}}{\func{\nei{\Xt}}{\seta{\point_2}}}}
\U_1\cap\U_2=\empty.
\end{equation}
\begin{itemize}
\item[${\textbf{\textsf{p1}}}$]
Each $\pp$ and $\y$ is taken as an element of $\Limits{\Xt}{\Yt}{\cf}{\point}$.
According to \refdef{deflimitoftopologicalmap},
\begin{align}
\Foreach{\p{\U}}{\func{\nei{\Yt}}{\seta{\p{\point}}}}
&\Exists{\U}{\func{\nei{\Xt}}{\seta{\point}}}
\func{\image{\cf}}{\compl{\[\U\cap\asubset\]}{\seta{\point}}}\subseteq\p{\U},
\label{thmuniquelimitp11eq1}\\
\Foreach{\p{\U}}{\func{\nei{\Yt}}{\seta{\y}}}
&\Exists{\U}{\func{\nei{\Xt}}{\seta{\point}}}
\func{\image{\cf}}{\compl{\[\U\cap\asubset\]}{\seta{\point}}}\subseteq\p{\U}.
\label{thmuniquelimitp11eq2}
\end{align}
\begin{itemize}
\item[${\textbf{\textsf{p1-1}}}$]
$\p{\U_1}$ and $\p{\U_2}$ are taken to be arbitrary elements of $\func{\nei{\Xt}}{\seta{\point}}$ and
$\func{\nei{\Xt}}{\seta{\y}}$, respectively. Then, according to \Ref{thmuniquelimitp11eq1}, and
\Ref{thmuniquelimitp11eq2},
\begin{align}
&\Existsis{\U_1}{\func{\nei{\Xt}}{\seta{\point}}}
\func{\image{\cf}}{\compl{\[\U_1\cap\asubset\]}{\seta{\point}}}\subseteq\p{\U_1},
\label{thmuniquelimitp111eq1}\\
&\Existsis{\U_2}{\func{\nei{\Xt}}{\seta{\point}}}
\func{\image{\cf}}{\compl{\[\U_2\cap\asubset\]}{\seta{\point}}}\subseteq\p{\U_2}.
\label{thmuniquelimitp111eq2}
\end{align}
Since $U_1\in\topology{\X}$, and $\U_2\in\topology{\X}$, it is obvious that,
\begin{equation}\label{thmuniquelimitp111eq3}
\(\U_1\cap\U_2\)\in\topology{\X}.
\end{equation}
On the other hand, it is clear that $\point\in\(\U_1\cap\U_2\)$. Therefore, $\(\U_1\cap\U_2\)$
is a neighborhood of $\point$ in $\Xt$.
\begin{equation}\label{thmuniquelimitp111eq4}
\(\U_1\cap\U_2\)\in\func{\nei{\Xt}}{\seta{\point}}.
\end{equation}\label{thmuniquelimitp111eq5}
Thus, according to \Ref{limitpointproperty},
\begin{equation}\label{thmuniquelimitp111eq6}
\compl{\[\(\U_1\cap\U_2\)\cap\asubset\]}{\seta{\point}}\neq\empty,
\end{equation}
or equivalently,
\begin{equation}\label{thmuniquelimitp111eq7}
\func{\image{\cf}}{\compl{\[\(\U_1\cap\U_2\)\cap\asubset\]}{\seta{\point}}}\neq\empty.
\end{equation}
On the other hand, it is obvious that,
\begin{align}
&\compl{\[\(\U_1\cap\U_2\)\cap\asubset\]}{\seta{\point}}\subseteq
\compl{\[\U_1\cap\asubset\]}{\seta{\point}},
\label{thmuniquelimitp111eq8}\\
&\compl{\[\(\U_1\cap\U_2\)\cap\asubset\]}{\seta{\point}}\subseteq
\compl{\[\U_2\cap\asubset\]}{\seta{\point}},
\label{thmuniquelimitp111eq9}
\end{align}
and hence,
\begin{align}
&\func{\image{\cf}}{\compl{\[\(\U_1\cap\U_2\)\cap\asubset\]}{\seta{\point}}}\subseteq
\func{\image{\cf}}{\compl{\[\U_1\cap\asubset\]}{\seta{\point}}},
\label{thmuniquelimitp111eq10}\\
&\func{\image{\cf}}{\compl{\[\(\U_1\cap\U_2\)\cap\asubset\]}{\seta{\point}}}\subseteq
\func{\image{\cf}}{\compl{\[\U_2\cap\asubset\]}{\seta{\point}}},
\label{thmuniquelimitp111eq11}
\end{align}
and followingly, according to \Ref{thmuniquelimitp111eq1}, and \Ref{thmuniquelimitp111eq2},
\begin{equation}\label{thmuniquelimitp111eq12}
\func{\image{\cf}}{\compl{\[\(\U_1\cap\U_2\)\cap\asubset\]}{\seta{\point}}}\subseteq
\(\p{\U_1}\cap\p{\U_2}\).
\end{equation}
\Ref{thmuniquelimitp111eq7} and \Ref{thmuniquelimitp111eq12} clearly imply,
\begin{equation}
\p{\U_1}\cap\p{\U_2}\neq\empty.
\end{equation}
\endp
\end{itemize}
So, it is seen that,
\begin{equation}
\Foreach{\opair{\p{\U_1}}{\p{\U_2}}}
{\Cprod{\func{\nei{\Xt}}{\seta{\point}}}{\func{\nei{\Xt}}{\seta{\y}}}}
\p{\U_1}\cap\p{\U_2}\neq\empty,
\end{equation}
and thus, according to \Ref{thmuniquelimitpeq1},
\begin{equation}
\y=\pp.
\end{equation}
\endp
\end{itemize}
So, it is seen that,
\begin{equation}
\Foreach{\opair{\pp}{\y}}{\Cprod{\Limits{\Xt}{\Yt}{\cf}{\point}}{\Limits{\Xt}{\Yt}{\cf}{\point}}}
\y=\pp,
\end{equation}
which means that either $\Limits{\Xt}{\Yt}{\cf}{\point}$ is empty or it has only one element.
%\end{itemize}
\endthm
%%%%%%%%%%%%%%%%%%%%%%%%%%%%%%%%%%%%%%%%%%%%%%%%%%%%%%%%%%%%%%%%%%%%%%%%%%%%%%%%%%%%%%%%%%%%%%%%%%%%%%%%%%%%%%%%%%%%%%%%%%%%%%%%%%%%%%%%%%%%%%%
%%%%%%%%%%%%%%%%%%%%%%%%%%%%%%%%%%%%%%%%%%%%%%%%%%%%%%%%%%%%%%%%%%%%%%%%%%%%%%%%%%%%%%%%%%%%%%%%%%%%%%%%%%%%%%%%%%%%%%%%%%%%%%%%%%%%%%%%%%%%%%%
%%%%%%%%%%%%%%%%%%%%%%%%%%%%%%%%%%%%%%%%%%%%%%%%%%%%%%%%%%%%%%%%%%%%%%%%%%%%%%%%%%%%%%%%%%%%%%%%%%%%%%%%%%%%%%%%%%%%%%%%%%%%%%%%%%%%%%%%%%%%%%%
%%%%%%%%%%%%%%%%%%%%%%%%%%%%%%%%%%%%%%%%%%%%%%%%%%%%%%%%%%%%%%%%%%%%%%%%%%%%%%%%%%%%%%%%%%%%%%%%%%%%%%%%%%%%%%%%%%%%%%%%%%%%%%%%%%%%%%%%%%%%%%%
%%%%%%%%%%%%%%%%%%%%%%%%%%%%%%%%%%%%%%%%%%%%%%%%%%%%%%%%%%%%%%%%%%%%%%%%%%%%%%%%%%%%%%%%%%%%%%%%%%%%%%%%%%%%%%%%%%%%%%%%%%%%%%%%%%%%%%%%%%%%%%%
%%%%%%%%%%%%%%%%%%%%%%%%%%%%%%%%%%%%%%%%%%%%%%%%%%%%%%%%%%%%%%%%%%%%%%%%%%%%%%%%%%%%%%%%%%%%%%%%%%%%%%%%%%%%%%%%%%%%%%%%%%%%%%%%%%%%%%%%%%%%%%%
\noindent
\textit{$\F=\triple{\f}{\fsum}{\fpro}$ is taken as one of the fields, $\R$ or $\C$.}
\section{Structure of a Normed Vector Space}
%%%%%%%%%%%%%%%%%%%%%%%%%%%%%%%%%%%%%%%%%%%%%%%%%%%%%%%%%%%%%%%%%%%%%%%%%%%%%%%%%%%%%%%%%%%%%%%%%%%%%%%%%%
\definition
$\V$ is taken as a set. The set of all vector-spaces having the set $\V$ as its' vectors,
and $\F$ as the field of scalars, is denoted by $\vecspaces{\V}{\F}$.
\begin{align}
&\vecspaces{\V}{\F}:=\cr
&\defset{\tuple{\V}{\vsum{}}{\spro{}}{\F}}
{\Cprod{\seta{\V}}{\Cprod{\Func{\Cprod{\V}{\V}}{\V}}{\Cprod{\Func{\Cprod{\f}{\V}}{\V}}{\seta{\F}}}}}
{\bigg[\tuple{\V}{\vsum{}}{\spro{}}{\F}~is~a~V.S\bigg]}.\cr
&{}
\end{align}
\endef
%%%%%%%%%%%%%%%%%%%%%%%%%%%%%%%%%%%%%%%%%%%%%%%%%%%%%%%%%%%%%%%%%%%%%%%%%%%%%%%%%%%%%%%%%%%%%%%%%%%%%%%%%%
\definition\label{defnorm}
$\V$ is taken as a non-empty set, and
$\VS=\tuple{\V}{\vsum{}}{\spro{}}{\F}$ is taken as an element of $\vecspaces{\V}{\F}$.
\begin{itemize}
\item
$\anorm{}$ is taken as an element of $\Func{\V}{\Rpz}$. $\anorm{}$ is called a
$\quotl$norm on the vector-space $\VS$$\quotr$ iff it possesses these properties.
\begin{itemize}
\item[${\textbf{\textsf{Nr1}}}$]\quad
$\Foreach{\v}{\V}\norm{\v}{}=0\thenn\v=0$.
\item[${\textbf{\textsf{Nr2}}}$]\quad
$\Foreach{\opair{\x}{\v}}{\Cprod{\f}{\V}}
\norm{\x\spro{}\v}{}=\abs{\x}\norm{\v}{}$.
\item[${\textbf{\textsf{Nr3}}}$]\quad
$\Foreach{\opair{\v}{\w}}{\Cprod{\V}{\V}}
\norm{\v+\w}{}\leq\norm{\v}{}+\norm{\w}{}$.
\end{itemize}
\item
The set of all norms on $\VS$ is denoted by $\norms{\VS}$.
\item
For each $\anorm{}$ in $\norms{\VS}$, $\opair{\VS}{\anorm{}}$
is referred to as a normed-space.
\end{itemize}
\endef
%%%%%%%%%%%%%%%%%%%%%%%%%%%%%%%%%%%%%%%%%%%%%%%%%%%%%%%%%%%%%%%%%%%%%%%%%%%%%%%%%%%%%%%%%%%%%%%%%%%%%%%%%%
\definition
Each $\V_1$ and $\V_2$ is taken as a non-empty set, and $\VS_1=\tuple{\V_1}{\vsum{1}}{\spro{1}}{\F}$ and
$\VS_2=\tuple{\V_2}{\vsum{2}}{\spro{2}}{\F}$ are taken as
elements of $\vecspaces{\V_1}{\F}$ and $\vecspaces{\V_2}{\F}$, respectively. $\anorm{1}$ and
$\anorm{2}$ are taken as elements of $\norms{\VS_1}$ and $\norms{\VS_2}$, respectively.
\begin{align}
&\anorm{\normsum{1}{2}}\indef\Func{\Cprod{\V_1}{\V_2}}{\Rpz},\cr
&\Foreach{\opair{\v}{\w}}{\Cprod{\V_1}{\V_2}}
\norm{\opair{\v}{\w}}{\normsum{1}{2}}\eqdef\norm{\v}{1}+\norm{\w}{2}.
\end{align}
\endef
%%%%%%%%%%%%%%%%%%%%%%%%%%%%%%%%%%%%%%%%%%%%%%%%%%%%%%%%%%%%%%%%%%%%%%%%%%%%%%%%%%%%%%%%%%%%%%%%%%%%%%%%%%
\theorem\label{thmcanonicalnormofproductspace}
Each $\V_1$ and $\V_2$ is taken as a non-empty set, and $\VS_1=\tuple{\V_1}{\vsum{1}}{\spro{1}}{\F}$ and
$\VS_2=\tuple{\V_2}{\vsum{2}}{\spro{2}}{\F}$ are taken as
elements of $\vecspaces{\V_1}{\F}$ and $\vecspaces{\V_2}{\F}$, respectively. $\anorm{1}$ and
$\anorm{2}$ are taken as elements of $\norms{\VS_1}$ and $\norms{\VS_2}$, respectively.
$\anorm{\normsum{1}{2}}$ is a norm on the vector-space $\Cprod{\VS_1}{\VS_2}$.
\begin{equation}
\anorm{\normsum{1}{2}}\in\norms{\Cprod{\VS_1}{\VS_2}}.
\end{equation}
\proof
It is clear according to \refdef{defnorm}.
\endthm
%%%%%%%%%%%%%%%%%%%%%%%%%%%%%%%%%%%%%%%%%%%%%%%%%%%%%%%%%%%%%%%%%%%%%%%%%%%%%%%%%%%%%%%%%%%%%%%%%%%%%%%%%%
\definition\label{defsetnorms}
$\V$ is taken as a non-empty set.
\begin{align}
\setnorms{\F}{\V}&:=\Union{\VS}{\vecspaces{\V}{\F}}{\norms{\VS}}\cr
&=\defset{\anorm{}}{\Func{\V}{\Rpz}}{\[\Exists{\VS}{\vecspaces{\V}{\F}}\anorm{}\in\norms{\VS}\]}.
\end{align}
\endef
%%%%%%%%%%%%%%%%%%%%%%%%%%%%%%%%%%%%%%%%%%%%%%%%%%%%%%%%%%%%%%%%%%%%%%%%%%%%%%%%%%%%%%%%%%%%%%%%%%%%%%%%%%
\definition\label{defnormdistance}
$\V$ is taken as a non-empty set.
\begin{align}
&\normdist{\F}{\V}\indef\Func{\setnorms{\F}{\V}}{\Func{\Cprod{\V}{\V}}{\Rpz}},\cr
&\Foreach{\anorm{}}{\setnorms{\F}{\V}}
\[\Foreach{\opair{\v}{\w}}{\Cprod{\V}{\V}}
\func{\func{\normdist{\F}{\V}}{\anorm{}}}{\binary{\v}{\w}}\eqdef\norm{\v-\w}{}\].
\end{align}
For each norm $\anorm{}$ in $\setnorms{\F}{\V}$, $\func{\normdist{\F}{\V}}{\anorm{}}$ is referred to as
the $\quotl$$\anorm{}$-induced distance on $\V$$\quotr$.
\endef
%%%%%%%%%%%%%%%%%%%%%%%%%%%%%%%%%%%%%%%%%%%%%%%%%%%%%%%%%%%%%%%%%%%%%%%%%%%%%%%%%%%%%%%%%%%%%%%%%%%%%%%%%%
\corollary\label{cornormdistanceisametric}
$\V$ is taken as a non-empty set.
For every $\anorm{}$ in $\setnorms{\F}{\V}$, $\opair{\V}{\func{\normdist{\F}{\V}}{\anorm{}}}$
is a metric-space.
\begin{align}
&\Foreach{\anorm{}}{\setnorms{\F}{\V}}\cr
&\quad\Foreach{\opair{\v}{\w}}{\Cprod{\V}{\V}}
\func{\func{\normdist{\F}{\VS}}{\anorm{}}}{\binary{\v}{\w}}=0\then\v=\w,\cr
&\quad\Foreach{\opair{\v}{\w}}{\Cprod{\V}{\V}}
\func{\func{\normdist{\F}{\V}}{\anorm{}}}{\binary{\v}{\w}}=
\func{\func{\normdist{\F}{\V}}{\anorm{}}}{\binary{\w}{\v}},\cr
&\quad\Foreach{\triple{\u}{\v}{\w}}{\V^{3}}
\func{\func{\normdist{\F}{\V}}{\anorm{}}}{\binary{\u}{\w}}\leq
\func{\func{\normdist{\F}{\V}}{\anorm{}}}{\binary{\u}{\v}}+
\func{\func{\normdist{\F}{\V}}{\anorm{}}}{\binary{\v}{\w}}.
\end{align}
\endcor
%%%%%%%%%%%%%%%%%%%%%%%%%%%%%%%%%%%%%%%%%%%%%%%%%%%%%%%%%%%%%%%%%%%%%%%%%%%%%%%%%%%%%%%%%%%%%%%%%%%%%%%%%%
\definition\label{deftopologyinducedbynorm}
$\V$ is taken as a non-empty set.
\begin{itemize}
\item
\begin{align}
&\normtop{\F}{\V}\indef\Func{\setnorms{\F}{\V}}{\Ctops{\V}},\cr
&\Foreach{\anorm{}}{\setnorms{\F}{\V}}\func{\normtop{\F}{\V}}{\anorm{}}\eqdef
\func{\metrictop{\V}}{\func{\normdist{\F}{\V}}{\anorm{}}},
\end{align}
where $\func{\metrictop{\V}}{\func{\normdist{\F}{\V}}{\anorm{}}}$ denotes the topology
induced by the metric $\func{\normdist{\F}{\V}}{\anorm{}}$ on the set $\V$.
For every $\anorm{}$ in $\setnorms{\F}{\V}$, $\func{\normtop{\F}{\VS}}{\anorm{}}$ is called the
$\quotl$canonical topology on $\V$ induced by the norm $\anorm{}$$\quotr$.
\item
For every $\anorm{}$ in $\setnorms{\F}{\V}$, the topological-space $\opair{\V}{\func{\normtop{\F}{\V}}{\anorm{}}}$
is denoted by $\normedtopspace{\anorm{}}$, and is referred to as the
$\quotl$canonical topological-space induced by the norm $\anorm{}$$\quotr$.
\end{itemize}
\endef
%%%%%%%%%%%%%%%%%%%%%%%%%%%%%%%%%%%%%%%%%%%%%%%%%%%%%%%%%%%%%%%%%%%%%%%%%%%%%%%%%%%%%%%%%%%%%%%%%%%%%%%%%%
\definition\label{defnormedspaceball}
$\V$ is taken as a non-empty set, and $\VS=\tuple{\V}{\vsum{}}{\spro{}}{\F}$ is taken as
an element of $\vecspaces{\V}{\F}$, and $\anorm{}$ an element of $\norms{\VS}$.
\begin{align}
&\ball{\opair{\VS}{\anorm{}}}\indef\Func{\Cprod{\V}{\Rp}}{\CSs{\V}},\cr
&\Foreach{\opair{\point}{r}}{\Cprod{\V}{\Rp}}
\func{\ball{\opair{\VS}{\anorm{}}}}{\binary{\point}{r}}\eqdef
\func{\mball{\V}{\func{\normdist{\F}{\V}}{\anorm{}}}}{\binary{\point}{r}},
\end{align}
and
\begin{align}
&\cball{\opair{\VS}{\anorm{}}}\indef\Func{\Cprod{\V}{\Rp}}{\CSs{\V}},\cr
&\Foreach{\opair{\point}{r}}{\Cprod{\V}{\Rp}}
\func{\cball{\opair{\VS}{\anorm{}}}}{\binary{\point}{r}}\eqdef
\func{\mcball{\V}{\func{\normdist{\F}{\V}}{\anorm{}}}}{\binary{\point}{r}}.
\end{align}
In other words, for every $\point$ in $\V$, and every positive real-number $r$,
$\func{\ball{\opair{\VS}{\anorm{}}}}{\binary{\point}{r}}$ is defined to be the $r$-ball
centered at $\point$ in the metric-space $\opair{\V}{\func{\normdist{\F}{\V}}{\anorm{}}}$,
and
$\func{\cball{\opair{\VS}{\anorm{}}}}{\binary{\point}{r}}$ is defined to be the $r$-closed-ball
centered at $\point$ in the metric-space $\opair{\V}{\func{\normdist{\F}{\V}}{\anorm{}}}$.
for every $\point$ in $\V$, and every positive real-number $r$,
$\func{\ball{\opair{\VS}{\anorm{}}}}{\binary{\point}{r}}$ is referred to as the
$\quotl$$r$-ball
centered at $\point$ in the normed-space $\opair{\VS}{\anorm{}}$$\quotr$,
and
$\func{\cball{\opair{\VS}{\anorm{}}}}{\binary{\point}{r}}$ is referred to as the
$\quotl$$r$-closed-ball
centered at $\point$ in the normed-space $\opair{\VS}{\anorm{}}$$\quotr$.
\endef
%%%%%%%%%%%%%%%%%%%%%%%%%%%%%%%%%%%%%%%%%%%%%%%%%%%%%%%%%%%%%%%%%%%%%%%%%%%%%%%%%%%%%%%%%%%%%%%%%%%%%%%%%%
\corollary\label{cornormedspaceball}
$\V$ is taken as a non-empty set, and $\VS=\tuple{\V}{\vsum{}}{\spro{}}{\F}$ is taken as
an element of $\vecspaces{\V}{\F}$, and $\anorm{}$ an element of $\norms{\VS}$.
\begin{itemize}
\item
\begin{align}
\Foreach{\opair{\point}{r}}{\Cprod{\V}{\Rp}}
\begin{cases}
\func{\ball{\opair{\VS}{\anorm{}}}}{\binary{\point}{r}}=
\defset{\v}{\V}{\norm{\v-\point}{}<r},\cr
\func{\cball{\opair{\VS}{\anorm{}}}}{\binary{\point}{r}}=
\defset{\v}{\V}{\norm{\v-\point}{}\leq r}.
\end{cases}
\end{align}
\item
$\defSet{\func{\ball{\opair{\VS}{\anorm{}}}}{\binary{\point}{r}}}{\point\in\V,~r\in\Rp}$
is a base for the topological space $\opair{\V}{\func{\normtop{\F}{\V}}{\anorm{}}}$.
\end{itemize}
\endcor
%%%%%%%%%%%%%%%%%%%%%%%%%%%%%%%%%%%%%%%%%%%%%%%%%%%%%%%%%%%%%%%%%%%%%%%%%%%%%%%%%%%%%%%%%%%%%%%%%%%%%%%%%%
\definition\label{defequivalentnorms}
$\V$ is taken as a non-empty set, and $\VS=\tuple{\V}{\vsum{}}{\spro{}}{\F}$ is taken as
an element of $\vecspaces{\V}{\F}$.
\begin{equation}
\normeqR{\VS}:=\defset{\opair{\anorm{1}}{\anorm{2}}}{\Cprod{\norms{\VS}}{\norms{\VS}}}
{\bigg[\func{\normtop{\F}{\V}}{\anorm{1}}=\func{\normtop{\F}{\V}}{\anorm{2}}\bigg]}.
\end{equation}
For each $\opair{\anorm{1}}{\anorm{2}}$ in $\normeqR{\VS}$, it is said that
$\quotl$$\anorm{1}$is equivalent to $\anorm{2}$$\quotr$.
\endef
%%%%%%%%%%%%%%%%%%%%%%%%%%%%%%%%%%%%%%%%%%%%%%%%%%%%%%%%%%%%%%%%%%%%%%%%%%%%%%%%%%%%%%%%%%%%%%%%%%%%%%%%%%
\corollary
$\V$ is taken as a non-empty set, and $\VS=\tuple{\V}{\vsum{}}{\spro{}}{\F}$ is taken as
an element of $\vecspaces{\V}{\F}$. $\normeqR{\VS}$ is an equivalence-relation on
the set of all norms on $\VS$, $\norms{\VS}$.
\endcor
%%%%%%%%%%%%%%%%%%%%%%%%%%%%%%%%%%%%%%%%%%%%%%%%%%%%%%%%%%%%%%%%%%%%%%%%%%%%%%%%%%%%%%%%%%%%%%%%%%%%%%%%%%
\theorem\label{thmNSnormequivalence}
$\V$ is taken as a non-empty set, and $\VS=\tuple{\V}{\vsum{}}{\spro{}}{\F}$ is taken as
an element of $\vecspaces{\V}{\F}$.
\begin{align}
%&\Foreach{\opair{\anorm{1}}{\anorm{2}}}{\Cprod{\norms{\VS}}{\norms{\VS}}}\cr
&\normeqR{\VS}=\cr
&\defset{\opair{\anorm{1}}{\anorm{2}}}{\norms{\VS}^2}
{\bigg[\Exists{\opair{m}{M}}{\(\Rp\)^2}\bigg(\Foreach{\v}{\V}m\norm{\v}{2}\leq
\norm{\v}{1}\leq M\norm{\v}{2}\bigg)\bigg]}.\cr
&{}
\end{align}
\proof

\endthm
%%%%%%%%%%%%%%%%%%%%%%%%%%%%%%%%%%%%%%%%%%%%%%%%%%%%%%%%%%%%%%%%%%%%%%%%%%%%%%%%%%%%%%%%%%%%%%%%%%%%%%%%%%
\theorem\label{thmnormsoffiniteVSareequivalent}
$\V$ is taken as a non-empty set, and $\VS=\tuple{\V}{\vsum{}}{\spro{}}{\F}$ is taken as
an element of $\vecspaces{\V}{\F}$. If $\VS$ is of finite dimension, then all norms on
$\VS$ are equivalent.
\begin{equation}
\func{\Vdim{}}{\VS}\in\N\then\normeqR{\VS}=\Cprod{\norms{\VS}}{\norms{\VS}}.
\end{equation}
\proof

\endthm
%%%%%%%%%%%%%%%%%%%%%%%%%%%%%%%%%%%%%%%%%%%%%%%%%%%%%%%%%%%%%%%%%%%%%%%%%%%%%%%%%%%%%%%%%%%%%%%%%%%%%%%%%%
\definition\label{defabsnorm}
$\VV{}$ is taken as a non-empty set, and $\VVS{}=\tuple{\VV{}}{\vsum{}}{\spro{}}{\F}$ is taken as
an element of $\vecspaces{\VV{}}{\F}$, such that $\func{\Vdim{}}{\VVS{}}\in\Zp$. $\vbase{}$ is taken
as an ordered-basis of the vector-space $\VVS{}$.
\begin{align}
&\anorm{\triple{1}{\VVS{}}{\vbase{}}}\indef\Func{\VV{}}{\Rpz},\cr
&\Foreach{\vv{}}{\VV{}}
\norm{\vv{}}{\triple{1}{\VVS{}}{\vbase{}}}\eqdef\ssum{k=1}{\func{\Vdim{}}{\VVS{}}}{\abs{\func{\[\com{\VVS{}}{\vv{}}{\vbase{}}\]}{k}}},
\end{align}
where, $\func{\[\com{\VVS{}}{\vv{}}{\vbase{}}\]}{k}$ denotes the $k$-th component of
the expansion of $\vv{}$ with respect to the ordered-basis $\vbase{}$.
\endef
%%%%%%%%%%%%%%%%%%%%%%%%%%%%%%%%%%%%%%%%%%%%%%%%%%%%%%%%%%%%%%%%%%%%%%%%%%%%%%%%%%%%%%%%%%%%%%%%%%%%%%%%%%
\theorem\label{thmabsnorm}
$\VV{}$ is taken as a non-empty set, and $\VVS{}=\tuple{\VV{}}{\vsum{}}{\spro{}}{\F}$ is taken as
an element of $\vecspaces{\VV{}}{\F}$, such that $\func{\Vdim{}}{\VVS{}}\in\Zp$. $\vbase{}$ is taken
as an ordered-basis of the vector-space $\VVS{}$. $\anorm{\triple{1}{\VVS{}}{\vbase{}}}$ is a norm on
$\VVS{}$.
\begin{equation}
\anorm{\triple{1}{\VVS{}}{\vbase{}}}\in\norms{\VVS{}}.
\end{equation}
\proof
It is trivial.
\endthm
%%%%%%%%%%%%%%%%%%%%%%%%%%%%%%%%%%%%%%%%%%%%%%%%%%%%%%%%%%%%%%%%%%%%%%%%%%%%%%%%%%%%%%%%%%%%%%%%%%%%%%%%%%
\theorem\label{thmabsnorm1}
$\VV{}$ is taken as a non-empty set, and $\VVS{}=\tuple{\VV{}}{\vsum{}}{\spro{}}{\F}$ is taken as
an element of $\vecspaces{\VV{}}{\F}$, such that $\func{\Vdim{}}{\VVS{}}\in\Zp$. $\vbase{}$ is taken
as an ordered-basis of the vector-space $\VVS{}$.
\begin{align}
\Foreach{\anorm{}}{\norms{\VVS{}}}
\Exists{M}{\Rp}
\bigg[\Foreach{\vv{}}{\VV{}}
\ssum{k=1}{\func{\Vdim{}}{\VVS{}}}{\abs{\func{\com{\VVS{}}{\vv{}}{\vbase{}}}{k}}}\leq
M\norm{\vv{}}{}\bigg].
\end{align}
\proof
According to \refthm{thmNSnormequivalence}, \refthm{thmnormsoffiniteVSareequivalent},
\refdef{defabsnorm}, and \refthm{thmabsnorm}, it is clear.
\endthm
%%%%%%%%%%%%%%%%%%%%%%%%%%%%%%%%%%%%%%%%%%%%%%%%%%%%%%%%%%%%%%%%%%%%%%%%%%%%%%%%%%%%%%%%%%%%%%%%%%%%%%%%%%
%%%%%%%%%%%%%%%%%%%%%%%%%%%%%%%%%%%%%%%%%%%%%%%%%%%%%%%%%%%%%%%%%%%%%%%%%%%%%%%%%%%%%%%%%%%%%%%%%%%%%%%%%%
%%%%%%%%%%%%%%%%%%%%%%%%%%%%%%%%%%%%%%%%%%%%%%%%%%%%%%%%%%%%%%%%%%%%%%%%%%%%%%%%%%%%%%%%%%%%%%%%%%%%%%%%%%
\subsection{Sequences and Convergence}
\definition\label{defconvergence}
$\V$ is taken as a non-empty set, and $\VS=\tuple{\V}{\vsum{}}{\spro{}}{\F}$ is taken as
an element of $\vecspaces{\V}{\F}$, and $\anorm{}$ an element of $\norms{\VS}$. $\seq{}$ is
taken as an element of $\Func{\N}{\V}$ (a sequence in $\V$). Also, $\NVS{}:=\opair{\VS}{\anorm{}}$.
\begin{itemize}
\item
$\v$ is taken as an element of $\V$. By definition, $\quotl$$\seq{}$ converges to $\v$
in the normed vector-space $\opair{\V}{\anorm{}}$$\quotr$, iff,
\begin{equation}
\Foreach{\varepsilon}{\Rp}
\Exists{N}{\N}\defSet{\func{\seq{}}{n}}{n\geq N}\subseteq\func{\ball{\opair{\VS}{\anorm{}}}}{\binary{\v}{\varepsilon}}.
\end{equation}
This property of $\v$ is also denoted by $\quotl$$\displaystyle\slimit{\NVS{}}{n}{\seq{n}}{\v}$$\quotr$.
\item
The set of all vectors of $\VS$ that $\seq{}$ vonverges to them in $\opair{\VS}{\anorm{}}$ is
denoted by $\conv{\opair{\VS}{\anorm{}}}{\seq{}}$. 	
\begin{equation}
\conv{\NVS{}}{\seq{}}:=
\defset{\v}{\V}
{\bigg[\Foreach{\varepsilon}{\Rp}
\Exists{N}{\N}\defSet{\func{\seq{}}{n}}{n\geq N}\subseteq
\func{\ball{\opair{\VS}{\anorm{}}}}{\binary{\v}{\varepsilon}}\bigg]}.
\end{equation}
\item
$\seq{}$ is defined to be a $\quotl$convergent sequence in the normed vector-space $\NVS{}$$\quotr$, iff
$\conv{\NVS{}}{\seq{}}\neq\empty$.
\end{itemize}
\endef
%%%%%%%%%%%%%%%%%%%%%%%%%%%%%%%%%%%%%%%%%%%%%%%%%%%%%%%%%%%%%%%%%%%%%%%%%%%%%%%%%%%%%%%%%%%%%%%%%%%%%%%%%%
\corollary\label{corconvergenceuniqueness}
$\V$ is taken as a non-empty set, and $\VS=\tuple{\V}{\vsum{}}{\spro{}}{\F}$ is taken as
an element of $\vecspaces{\V}{\F}$, and $\anorm{}$ an element of $\norms{\VS}$. $\seq{}$ is
taken as an element of $\Func{\N}{\V}$ (a sequence in $\V$). Also, $\NVS{}:=\opair{\VS}{\anorm{}}$.
\begin{equation}
\Card{\conv{\opair{\VS}{\anorm{}}}{\seq{}}}\leq 1.
\end{equation}
This means, if there exists a vector $\v$ of $\VS$ such that
$\displaystyle\slimit{\NVS{}}{n}{\seq{n}}{\v}$, then this is the only vector with this property.
\endcor
%%%%%%%%%%%%%%%%%%%%%%%%%%%%%%%%%%%%%%%%%%%%%%%%%%%%%%%%%%%%%%%%%%%%%%%%%%%%%%%%%%%%%%%%%%%%%%%%%%%%%%%%%%
\corollary\label{corconvergenceproperties1}
$\V$ is taken as a non-empty set, and $\VS=\tuple{\V}{\vsum{}}{\spro{}}{\F}$ is taken as
an element of $\vecspaces{\V}{\F}$, and $\anorm{}$ an element of $\norms{\VS}$. Each $\seq{}$ and $\seqq{}$
is taken as an element of $\Func{\N}{\V}$ (a sequence in $\V$), and $\x$ an element of $\ff$.
Also, $\NVS{}:=\opair{\VS}{\anorm{}}$.
\begin{align}
&\Foreach{\opair{\v}{\w}}{\Cprod{\conv{\NVS{}}{\seq{}}}{\conv{\NVS{}}{\seqq{}}}}
\v\vsum{}\w\in\conv{\NVS{}}{\seq{}+\seqq{}}.\\
&\Foreach{\v}{\conv{\NVS{}}{\seq{}}}\x\spro{}\v\in\conv{\NVS{}}{\x\spro{}\seq{}}.\\
&\Foreach{\v}{\conv{\NVS{}}{\seq{}}}\slimit{}{n}{\norm{\seq{n}-\v}{}}{0}.
\end{align}
\endcor
%%%%%%%%%%%%%%%%%%%%%%%%%%%%%%%%%%%%%%%%%%%%%%%%%%%%%%%%%%%%%%%%%%%%%%%%%%%%%%%%%%%%%%%%%%%%%%%%%%%%%%%%%%
\definition\label{defsetofallconvergentsequences}
$\V$ is taken as a non-empty set, and $\VS=\tuple{\V}{\vsum{}}{\spro{}}{\F}$ is taken as
an element of $\vecspaces{\V}{\F}$, and $\anorm{}$ an element of $\norms{\VS}$.
Also, $\NVS{}:=\opair{\VS}{\anorm{}}$. The set of all convergent sequences in $\NVS{}$ is
denoted by $\Con{\NVS{}}$.
\begin{equation}
\Con{\NVS{}}:=\defset{\seq{}}{\Func{\N}{\V}}{\conv{\NVS{}}{\seq{}}\neq\empty}.
\end{equation}
\endef
%%%%%%%%%%%%%%%%%%%%%%%%%%%%%%%%%%%%%%%%%%%%%%%%%%%%%%%%%%%%%%%%%%%%%%%%%%%%%%%%%%%%%%%%%%%%%%%%%%%%%%%%%%
\definition\label{defcauchysequence}
$\V$ is taken as a non-empty set, and $\VS=\tuple{\V}{\vsum{}}{\spro{}}{\F}$ is taken as
an element of $\vecspaces{\V}{\F}$, and $\anorm{}$ an element of $\norms{\VS}$.
Also, $\NVS{}:=\opair{\VS}{\anorm{}}$.
\begin{itemize}
\item
$\seq{}$ is taken as an element of $\Func{\N}{\V}$. $\seq{}$ is called a
$\quotl$Cauchy sequence in the normed vector-space $\NVS{}$$\quotr$, iff,
\begin{equation}
\Foreach{\varepsilon}{\Rp}
\Exists{N}{\N}\defSet{\norm{\func{\seq{}}{n}-\func{\seq{}}{m}}{}}{n\geq N,~m\geq N}
\subseteq\linterval{0}{\varepsilon}.
\end{equation}
\item
The set of all Cauchy sequences in $\NVS{}$ is denoted by $\Cauchy{\NVS{}}$.
\begin{align}
&\Cauchy{\NVS{}}:=\cr
&\defset{\seq{}}{\Func{\N}{\V}}
{\bigg[
\Foreach{\varepsilon}{\Rp}
\Exists{N}{\N}\defSet{\norm{\func{\seq{}}{n}-\func{\seq{}}{m}}{}}{n\geq N,~m\geq N}
\subseteq\linterval{0}{\varepsilon}\bigg]}.\cr
&{}
\end{align}
\end{itemize}
\endef
%%%%%%%%%%%%%%%%%%%%%%%%%%%%%%%%%%%%%%%%%%%%%%%%%%%%%%%%%%%%%%%%%%%%%%%%%%%%%%%%%%%%%%%%%%%%%%%%%%%%%%%%%%
\corollary\label{coreveryconvergentsequenceisCauchy}
$\V$ is taken as a non-empty set, and $\VS=\tuple{\V}{\vsum{}}{\spro{}}{\F}$ is taken as
an element of $\vecspaces{\V}{\F}$, and $\anorm{}$ an element of $\norms{\VS}$.
Also, $\NVS{}:=\opair{\VS}{\anorm{}}$.
Every convergent sequence in $\NVS{}$ is a Cauchy sequence in $\NVS{}$.
\begin{equation}
\Con{\NVS{}}\subseteq\Cauchy{\NVS{}}.
\end{equation}
\endcor
%%%%%%%%%%%%%%%%%%%%%%%%%%%%%%%%%%%%%%%%%%%%%%%%%%%%%%%%%%%%%%%%%%%%%%%%%%%%%%%%%%%%%%%%%%%%%%%%%%%%%%%%%%
%%%%%%%%%%%%%%%%%%%%%%%%%%%%%%%%%%%%%%%%%%%%%%%%%%%%%%%%%%%%%%%%%%%%%%%%%%%%%%%%%%%%%%%%%%%%%%%%%%%%%%%%%%
%%%%%%%%%%%%%%%%%%%%%%%%%%%%%%%%%%%%%%%%%%%%%%%%%%%%%%%%%%%%%%%%%%%%%%%%%%%%%%%%%%%%%%%%%%%%%%%%%%%%%%%%%%
\subsection{Banach Spaces}
\definition\label{defBanachspace}
$\V$ is taken as a non-empty set, and $\VS=\tuple{\V}{\vsum{}}{\spro{}}{\F}$ is taken as
an element of $\vecspaces{\V}{\F}$.
\begin{itemize}
\item
$\anorm{}$ is taken as an element of an element of $\norms{\VS}$.
Also, $\NVS{}:=\opair{\VS}{\anorm{}}$.
$\NVS{}$ is called a $\quotl$Banach-space$\quotr$, iff
\begin{equation}
\Cauchy{\NVS{}}=\Con{\NVS{}},
\end{equation}
which means, every Cauchy sequence in $\NVS{}$ is a convergent sequence in $\NVS{}$.
\item
The set of all norms on $\VS{}$ that make it a Banach-space is denoted by
$\Banachnorms{\VS{}}$.
\begin{equation}
\Banachnorms{\VS{}}:=\defset{\anorm{}}{\norms{\VS{}}}
{\Cauchy{\opair{\VS{}}{\anorm{}}}=\Con{\opair{\VS{}}{\anorm{}}}}.
\end{equation}
\end{itemize}
\endef
%%%%%%%%%%%%%%%%%%%%%%%%%%%%%%%%%%%%%%%%%%%%%%%%%%%%%%%%%%%%%%%%%%%%%%%%%%%%%%%%%%%%%%%%%%%%%%%%%%%%%%%%%%%%%%%%%%
\theorem\label{thmfinitedimensionalNVSisBanach}
$\V$ is taken as a non-empty set, and $\VS=\tuple{\V}{\vsum{}}{\spro{}}{\F}$ is taken as
an element of $\vecspaces{\V}{\F}$, and $\anorm{}$ an element of $\norms{\VS}$.
If $\func{\Vdim{}}{\VS}\in\N$, then $\opair{\VS}{\anorm{}}$ is a Banach-space.
\proof

\endthm
%%%%%%%%%%%%%%%%%%%%%%%%%%%%%%%%%%%%%%%%%%%%%%%%%%%%%%%%%%%%%%%%%%%%%%%%%%%%%%%%%%%%%%%%%%%%%%%%%%%%%%%%%%%%%%%%%%
%%%%%%%%%%%%%%%%%%%%%%%%%%%%%%%%%%%%%%%%%%%%%%%%%%%%%%%%%%%%%%%%%%%%%%%%%%%%%%%%%%%%%%%%%%%%%%%%%%%%%%%%%%%%%%%%%%
%%%%%%%%%%%%%%%%%%%%%%%%%%%%%%%%%%%%%%%%%%%%%%%%%%%%%%%%%%%%%%%%%%%%%%%%%%%%%%%%%%%%%%%%%%%%%%%%%%%%%%%%%%%%%%%%%%
\section{Linear-Continuous Maps}
\definition\label{defcontinuity}
Each $\V_1$ and $\V_2$ is taken as a non-empty set, and $\VS_1=\tuple{\V_1}{\vsum{1}}{\spro{1}}{\F}$ and
$\VS_2=\tuple{\V_2}{\vsum{2}}{\spro{2}}{\F}$ are taken as
elements of $\vecspaces{\V_1}{\F}$ and $\vecspaces{\V_2}{\F}$, respectively. $\anorm{1}$ and
$\anorm{2}$ are taken as elements of $\norms{\VS_1}$ and $\norms{\VS_2}$, respectively.
\begin{itemize}
\item
\begin{align}
\CF{\opair{\VS_1}{\anorm{1}}}{\opair{\VS_2}{\anorm{2}}}&:=
\CF{\opair{\V_1}{\func{\normtop{\F}{\V_1}}{\anorm{1}}}}{\opair{\V_2}{\func{\normtop{\F}{\V_2}}{\anorm{2}}}}\cr
&=\CF{\normedtopspace{\anorm{1}}}{\normedtopspace{\anorm{2}}}.
\end{align}
By definition, Each element of $\CF{\opair{\VS_1}{\anorm{1}}}{\opair{\VS_2}{\anorm{2}}}$ is
referred to as a $\quotl$continuous-map from the normed-space $\opair{\VS_1}{\anorm{1}}$
to the normed-space $\opair{\VS_2}{\anorm{2}}$$\quotr$,\\
or simply a $\quotl$continuous-map$\quotr$.
\item
$\cf$ is taken as an element of $\Func{\V_1}{\V_2}$, and $\point$ an element of $\V_1$.
$\cf$ is defined to be a $\quotl$continuous map from $\opair{\VS_1}{\anorm{1}}$ to
$\opair{\VS_2}{\anorm{2}}$ at $\point$$\quotr$ iff $\cf$ is a continuous map from the topological space
$\opair{\V_1}{\func{\normtop{\F}{\V}}{\anorm{1}}}$ to the topological space
$\opair{\V_2}{\func{\normtop{\F}{\V}}{\anorm{2}}}$ at $\point$. It can also simply be said that
$\quotl$$\cf$ is continuous at $\point$$\quotr$.
\end{itemize}
\endef
%%%%%%%%%%%%%%%%%%%%%%%%%%%%%%%%%%%%%%%%%%%%%%%%%%%%%%%%%%%%%%%%%%%%%%%%%%%%%%%%%%%%%%%%%%%%%%%%%%%%%%%%%%
\corollary\label{corcontinuityequivs}
Each $\V_1$ and $\V_2$ is taken as a non-empty set, and $\VS_1=\tuple{\V_1}{\vsum{1}}{\spro{1}}{\F}$ and
$\VS_2=\tuple{\V_2}{\vsum{2}}{\spro{2}}{\F}$ are taken as
elements of $\vecspaces{\V_1}{\F}$ and $\vecspaces{\V_2}{\F}$, respectively. $\anorm{1}$ and
$\anorm{2}$ are taken as elements of $\norms{\VS_1}$ and $\norms{\VS_2}$, respectively.
$\cf$ is taken as an element of $\Func{\V_1}{\V_2}$.
\begin{itemize}
\item
For every $\point$ in $\V_1$, $\cf$ is continuous at $\point$ if-and-only-if
\begin{equation}
\Foreach{\varepsilon}{\Rp}
\Exists{\delta}{\Rp}\func{\image{\cf}}{\func{\ball{\opair{\VS_1}{\anorm{1}}}}{\binary{\point}{\delta}}}\subseteq
\func{\ball{\opair{\VS_2}{\anorm{2}}}}{\binary{\func{\cf}{\point}}{\varepsilon}}.
\end{equation}
\item
$\cf$ is continuous if-and-only-if it is continuous at every $\point$ in $\V_1$.
\begin{gather}
\cf\in\CF{\opair{\VS_1}{\anorm{1}}}{\opair{\VS_2}{\anorm{2}}}\cr
\vthenn\cr
\Foreach{\point}{\V_1}\Foreach{\varepsilon}{\Rp}
\Exists{\delta}{\Rp}\func{\image{\cf}}{\func{\ball{\opair{\VS_1}{\anorm{1}}}}{\binary{\point}{\delta}}}\subseteq
\func{\ball{\opair{\VS_2}{\anorm{2}}}}{\binary{\func{\cf}{\point}}{\varepsilon}}.
\end{gather}
\end{itemize}
\endcor
%%%%%%%%%%%%%%%%%%%%%%%%%%%%%%%%%%%%%%%%%%%%%%%%%%%%%%%%%%%%%%%%%%%%%%%%%%%%%%%%%%%%%%%%%%%%%%%%%%%%%%%%%%
\definition\label{deflinearcontinuousmap}
Each $\V_1$ and $\V_2$ is taken as a non-empty set, and $\VS_1=\tuple{\V_1}{\vsum{1}}{\spro{1}}{\F}$ and
$\VS_2=\tuple{\V_2}{\vsum{2}}{\spro{2}}{\F}$ are taken as
elements of $\vecspaces{\V_1}{\F}$ and $\vecspaces{\V_2}{\F}$, respectively. $\anorm{1}$ and
$\anorm{2}$ are taken as elements of $\norms{\VS_1}$ and $\norms{\VS_2}$, respectively.
\begin{equation}
\LinC{\opair{\VS_1}{\anorm{1}}}{\opair{\VS_2}{\anorm{2}}}:=
\Lin{\VS_1}{\VS_2}\cap\CF{\opair{\VS_1}{\anorm{1}}}{\opair{\VS_2}{\anorm{2}}}.
\end{equation}
Every element of $\LinC{\opair{\VS_1}{\anorm{1}}}{\opair{\VS_2}{\anorm{2}}}$ is
called a $\quotl$linear-continuous map from the normed-space $\opair{\VS_1}{\anorm{1}}$
to the normed-space $\opair{\VS_2}{\anorm{2}}$$\quotr$.\\
In other words, a map $\cf:\V_1\to\V_2$ is called a linear-continuous map from
$\opair{\VS_1}{\anorm{1}}$ to $\opair{\VS_2}{\anorm{2}}$ iff
it is both a linear map from $\VS_1$ to $\VS_2$, and a continuous map from
$\opair{\VS_1}{\anorm{1}}$ to $\opair{\VS_2}{\anorm{2}}$.
\endef
%%%%%%%%%%%%%%%%%%%%%%%%%%%%%%%%%%%%%%%%%%%%%%%%%%%%%%%%%%%%%%%%%%%%%%%%%%%%%%%%%%%%%%%%%%%%%%%%%%%%%%%%%%
\theorem\label{thmlinearcontinuousmapequivs}
Each $\V_1$ and $\V_2$ is taken as a non-empty set, and $\VS_1=\tuple{\V_1}{\vsum{1}}{\spro{1}}{\F}$ and
$\VS_2=\tuple{\V_2}{\vsum{2}}{\spro{2}}{\F}$ are taken as
elements of $\vecspaces{\V_1}{\F}$ and $\vecspaces{\V_2}{\F}$, respectively. $\anorm{1}$ and
$\anorm{2}$ are taken as elements of $\norms{\VS_1}$ and $\norms{\VS_2}$, respectively.
$\lcf$ is taken as an element of $\Lin{\VS_1}{\VS_2}$.
These propositions are equivalent.
\begin{itemize}
\item[${\propos{1}}.$]\quad
$\lcf\in\LinC{\opair{\VS_1}{\anorm{1}}}{\opair{\VS_2}{\anorm{2}}}$.
\item[${\propos{2}}.$]\quad
$\lcf$ is continuous at $\zerov{\VS_1}$.
\item[${\propos{3}}.$]\quad
$\Exists{K}{\Rp}\Foreach{\v}{\defset{\w}{\V_1}{\norm{\w}{1}\leq 1}}\norm{\func{\lcf}{\v}}{2}\leq K$.
\item[${\propos{4}}.$]\quad
$\Exists{K}{\Rp}\Foreach{\v}{\V_1}\norm{\func{\lcf}{\v}}{2}\leq K\norm{\v}{1}$.
\end{itemize}
\proof\\
$\[{\propos{1}\rightarrow\propos{2}}\]$
According to \refcor{corcontinuityequivs}, it is obvious.\\
$\[{\propos{2}\rightarrow\propos{3}}\]$
It is assumed that $\lcf$ is continuous at $\zerov{\VS_1}$. So, according to \refcor{corcontinuityequivs},
\begin{equation}\label{thmlinearcontinuousmapequivsp2p3eq1}
\Foreach{\varepsilon}{\Rp}
\Exists{\delta}{\Rp}\func{\image{\lcf}}{\func{\ball{\opair{\VS_1}{\anorm{1}}}}{\binary{\zerov{\VS_1}}{\delta}}}\subseteq
\func{\ball{\opair{\VS_2}{\anorm{2}}}}{\binary{\func{\lcf}{\zerov{\VS_1}}}{\varepsilon}}.
\end{equation}
Thus, according to \refcor{cornormedspaceball}, and considering the fact that
$\func{\lcf}{\zerov{\VS_1}}=\zerov{\VS_2}$, it is clear that,
\begin{align}\label{thmlinearcontinuousmapequivsp2p3eq2}
\Foreach{\varepsilon}{\Rp}
\Exists{\delta}{\Rp}\func{\image{\lcf}}{\defset{\w}{\V_1}{\norm{\w}{1}\leq\delta}}\subseteq
\defset{\w}{\V_2}{\norm{\w}{2}\leq\varepsilon}.
\end{align}
So,
\begin{equation}\label{thmlinearcontinuousmapequivsp2p3eq3}
\Existsis{\delta_1}{\Rp}\func{\image{\lcf}}{\defset{\w}{\V_1}{\norm{\w}{1}\leq\delta_1}}\subseteq
\defset{\w}{\V_2}{\norm{\w}{2}\leq 1},
\end{equation}
and equivalently,
\begin{equation}\label{thmlinearcontinuousmapequivsp2p3eq4}
\Foreach{\v}{\defset{\w}{\V_1}{\norm{\w}{1}\leq\delta_1}}\norm{\func{\lcf}{\v}}{2}\leq 1.
\end{equation}
On the other hand, according to \refdef{defnorm}, it is clear that,
\begin{equation}\label{thmlinearcontinuousmapequivsp2p3eq5}
\Foreach{\v}{\defset{\w}{\V_1}{\norm{\w}{1}\leq 1}}
\delta_1\v\in{\defset{\w}{\V_1}{\norm{\w}{1}\leq\delta_1}},
\end{equation}
and hence, according to \Ref{thmlinearcontinuousmapequivsp2p3eq4},
\begin{equation}\label{thmlinearcontinuousmapequivsp2p3eq6}
\Foreach{\v}{\defset{\w}{\V_1}{\norm{\w}{1}\leq 1}}
\norm{\func{\lcf}{\delta_1\v}}{2}\leq 1.
\end{equation}
On the other hand, according to \refdef{defnorm} and linearity of $\lcf$,
\begin{equation}\label{thmlinearcontinuousmapequivsp2p3eq7}
\Foreach{\v}{\defset{\w}{\V_1}{\norm{\w}{1}\leq 1}}
\norm{\func{\lcf}{\delta_1\v}}{2}=\delta_1\norm{\func{\lcf}{\v}}{2}.
\end{equation}
\Ref{thmlinearcontinuousmapequivsp2p3eq6} and \Ref{thmlinearcontinuousmapequivsp2p3eq7}
imply,
\begin{equation}
\Existsis{\delta_1^{-1}}{\Rp}\Foreach{\v}{\defset{\w}{\V_1}{\norm{\w}{1}\leq 1}}\norm{\func{\lcf}{\v}}{2}\leq\delta_1^{-1}.
\end{equation}
$\[{\propos{3}\rightarrow\propos{4}}\]$
It is assumed that,
\begin{equation}
\Existsis{K}{\Rp}\Foreach{\v}{\defset{\w}{\V_1}{\norm{\w}{1}\leq 1}}
\norm{\func{\lcf}{\v}}{2}\leq K.
\end{equation}
Then, since
\begin{equation}
\Foreach{\v}{\compl{\V_1}{\seta{\zerov{\VS_1}}}}
\norm{\frac{1}{\norm{\v}{1}}\v}{1}=1,
\end{equation}
it is clear that,
\begin{equation}
\Foreach{\v}{\compl{\V_1}{\seta{\zerov{\VS_1}}}}
\norm{\func{\lcf}{\frac{1}{\norm{\v}{1}}\v}}{2}\leq K,
\end{equation}
and hence, according to \refdef{defnorm}, and linearity of $\lcf$,
\begin{equation}
\Foreach{\v}{\compl{\V_1}{\seta{\zerov{\VS_1}}}}
\norm{\func{\lcf}{\v}}{2}\leq K\norm{\v}{1}.
\end{equation}
Additionally, since $\norm{\zerov{\VS_2}}{2}=\norm{\func{\lcf}{\zerov{\VS_1}}}{2}=\norm{\zerov{\VS_1}}{1}=0$,
it is obvious that,
\begin{equation}
\norm{\func{\lcf}{\zerov{\VS_1}}}{2}\leq K\norm{\zerov{\VS_1}}{1}.
\end{equation}
So,
\begin{equation}
\Existsis{K}{\Rp}\Foreach{\v}{\V_1}\norm{\func{\lcf}{\v}}{2}\leq K\norm{\v}{1}.
\end{equation}
$\[{\propos{4}\rightarrow\propos{1}}\]$
It is assumed that,
\begin{equation}
\Existsis{K}{\Rp}\Foreach{\v}{\V_1}\norm{\func{\lcf}{\v}}{2}\leq K\norm{\v}{1}.
\end{equation}
Then, according to linearity of $\lcf$,
\begin{align}
&\Foreach{\point}{\V_1}\cr
&\Foreach{\varepsilon}{\Rp}
\Existsis{K^{-1}\varepsilon}{\Rp}
\func{\image{\lcf}}{\func{\ball{\opair{\VS_1}{\anorm{1}}}}{\binary{\point}{K^{-1}\varepsilon}}}\subseteq
\func{\ball{\opair{\VS_2}{\anorm{2}}}}{\binary{\func{\lcf}{\point}}{\varepsilon}},
\end{align}
and hence, according to \refcor{corcontinuityequivs},
\begin{align}
\lcf\in\LinC{\opair{\VS_1}{\anorm{1}}}{\opair{\VS_2}{\anorm{2}}}.
\end{align}
\endthm
%%%%%%%%%%%%%%%%%%%%%%%%%%%%%%%%%%%%%%%%%%%%%%%%%%%%%%%%%%%%%%%%%%%%%%%%%%%%%%%%%%%%%%%%%%%%%%%%%%%%%%%%%%
\corollary\label{corcompositionoflinearcontinuousmaps}
Each $\NVS{1}=\opair{\VS_1}{\anorm{1}}$, $\NVS{2}=\opair{\VS_2}{\anorm{2}}$, and
$\NVS{3}:=\opair{\VS_3}{\anorm{3}}$ is taken to be a normed vector-space.
\begin{align}
\Foreach{\opair{\lcf}{\p{\lcf}}}{\Cprod{\LinC{\NVS{1}}{\NVS{2}}}
{\LinC{\NVS{2}}{\NVS{3}}}}\cmp{\p{\lcf}}{\lcf}\in\LinC{\NVS{1}}{\NVS{3}}.
\end{align}
\endcor
%%%%%%%%%%%%%%%%%%%%%%%%%%%%%%%%%%%%%%%%%%%%%%%%%%%%%%%%%%%%%%%%%%%%%%%%%%%%%%%%%%%%%%%%%%%%%%%%%%%%%%%%%%
\corollary\label{coridentityisalinearcontinuousmap}
$\V$ is taken as a non-empty set, and $\VS=\tuple{\V}{\vsum{}}{\spro{}}{\F}$ is taken as
an element of $\vecspaces{\V}{\F}$. $\anorm{}$ is taken as an element of $\norms{\VS}$.
\begin{equation}
\identity{\V}\in\LinC{\opair{\VS}{\anorm{}}}{\opair{\VS}{\anorm{}}}.
\end{equation}
\endcor
%%%%%%%%%%%%%%%%%%%%%%%%%%%%%%%%%%%%%%%%%%%%%%%%%%%%%%%%%%%%%%%%%%%%%%%%%%%%%%%%%%%%%%%%%%%%%%%%%%%%%%%%%%
\theorem\label{thmlinearcontinuousmapsspace}
Each $\V_1$ and $\V_2$ is taken as a non-empty set, and $\VS_1=\tuple{\V_1}{\vsum{1}}{\spro{1}}{\F}$ and
$\VS_2=\tuple{\V_2}{\vsum{2}}{\spro{2}}{\F}$ are taken as
elements of $\vecspaces{\V_1}{\F}$ and $\vecspaces{\V_2}{\F}$, respectively. $\anorm{1}$ and
$\anorm{2}$ are taken as elements of $\norms{\VS_1}$ and $\norms{\VS_2}$, respectively.
\begin{align}
\tuple{\LinC{\opair{\VS_1}{\anorm{1}}}{\opair{\VS_2}{\anorm{2}}}}{\vsum{lc}}{\spro{lc}}{\F}\in
\vecspaces{\LinC{\opair{\VS_1}{\anorm{1}}}{\opair{\VS_2}{\anorm{2}}}}{\F},
\end{align}
where,
\begin{align}
\vsum{lc}&:=\func{\resd{\vsum{l}}}{\Cprod{\LinC{\opair{\VS_1}{\anorm{1}}}{\opair{\VS_2}{\anorm{2}}}}
{\LinC{\opair{\VS_1}{\anorm{1}}}{\opair{\VS_2}{\anorm{2}}}}},\\
\spro{lc}&:=\func{\resd{\spro{l}}}{\Cprod{\f}{\LinC{\opair{\VS_1}{\anorm{1}}}{\opair{\VS_2}{\anorm{2}}}}},
\end{align}
$\vsum{l}$ and $\spro{l}$ being the addition and scalar-multiplication operations of
the vector-space $\VLin{\VS_1}{\VS_2}$.
\proof
$\LinC{\opair{\VS_1}{\anorm{1}}}{\opair{\VS_2}{\anorm{2}}}$ is a non-empty subset of
$\Lin{\VS_1}{\VS_2}$, since the constant map that maps every vector of
$\VS_1$ to $\zerov{\VS_2}$ is a linear-continuous map from $\opair{\VS_1}{\anorm{1}}$
to $\opair{\VS_2}{\anorm{2}}$.
\begin{itemize}
\item[${\textbf{\textsf{p1}}}$]
Each $\lcf$ and $\p{\lcf}$ is taken as an element of
$\LinC{\opair{\VS_1}{\anorm{1}}}{\opair{\VS_2}{\anorm{2}}}$. So, according to
\refdef{deflinearcontinuousmap}, $\lcf$ and $\p{\lcf}$ belong to $\Lin{\VS_1}{\VS_2}$,
and hence,
\begin{equation}
\lcf\vsum{l}\p{\lcf}\in\Lin{\VS_1}{\VS_2}.
\end{equation}
On the other hand, according to \refthm{thmlinearcontinuousmapequivs},
\begin{align}
&\Existsis{K}{\Rp}\Foreach{\v}{\V_1}\norm{\func{\lcf}{\v}}{2}\leq K\norm{\v}{1},\\
&\Existsis{\p{K}}{\Rp}\Foreach{\v}{\V_1}\norm{\func{\p{\lcf}}{\v}}{2}\leq\p{K}\norm{\v}{1}.
\end{align}
Therefore, according to \refdef{defnorm},
\begin{align}
\Existsis{K+\p{K}}{\Rp}
\Foreach{\v}{\V_1}\norm{\func{\(\lcf\vsum{l}\p{\lcf}\)}{\v}}{2}&=
\norm{\func{\lcf}{\v}\vsum{2}\func{\p{\lcf}}{\v}}{2}\cr
&\leq\norm{\func{\lcf}{\v}}{2}+\norm{\func{\p{\lcf}}{\v}}{2}\cr
&\leq K\norm{\v}{1}+\p{K}\norm{\v}{1}\cr
&=\(K+\p{K}\)\norm{\v}{1}.
\end{align}
So, according to \refthm{thmlinearcontinuousmapequivs},
\begin{equation}
\(\lcf\vsum{l}\p{\lcf}\)\in\LinC{\opair{\VS_1}{\anorm{1}}}{\opair{\VS_2}{\anorm{2}}}.
\end{equation}
\endp
\end{itemize}
\begin{itemize}
\item[${\textbf{\textsf{p2}}}$]
$\lcf$ is taken as an element of
$\LinC{\opair{\VS_1}{\anorm{1}}}{\opair{\VS_2}{\anorm{2}}}$, and $\x$ an element of $\f$.
So, according to \refdef{deflinearcontinuousmap}, $\lcf$ belongs to $\Lin{\VS_1}{\VS_1}$,
and hence,
\begin{equation}
\(\x\spro{l}\lcf\)\in\Lin{\VS_1}{\VS_2}.
\end{equation}
On the other hand, according to \refthm{thmlinearcontinuousmapequivs},
\begin{equation}
\Existsis{K}{\Rp}\Foreach{\v}{\V_1}\norm{\func{\lcf}{\v}}{2}\leq K\norm{\v}{1},
\end{equation}
and thus, according to \refdef{defnorm},
\begin{align}
\Existsis{\(\abs{\x}K+1\)}{\Rp}
\Foreach{\v}{\V_1}\norm{\func{\(\x\spro{l}\lcf\)}{\v}}{2}&=
\norm{\x\spro{2}\(\func{\lcf}{\v}\)}{2}\cr
&=\abs{\x}\norm{\func{\lcf}{\v}}{2}\cr
&\leq\(\abs{\x}K\)\norm{\v}{1}\cr
&\leq\(\abs{\x}K+1\)\norm{\v}{1}.
\end{align}
So, according to \refthm{thmlinearcontinuousmapequivs},
\begin{equation}
\(\x\spro{l}\lcf\)\in\LinC{\opair{\VS_1}{\anorm{1}}}{\opair{\VS_2}{\anorm{2}}}.
\end{equation}
\endp
\end{itemize}
\endthm
%%%%%%%%%%%%%%%%%%%%%%%%%%%%%%%%%%%%%%%%%%%%%%%%%%%%%%%%%%%%%%%%%%%%%%%%%%%%%%%%%%%%%%%%%%%%%%%%%%%%%%%%%%
\theorem\label{thmLmapswithfinitedimensionaldomainisLC}
Each $\V_1$ and $\V_2$ is taken as a non-empty set, and $\VS_1=\tuple{\V_1}{\vsum{1}}{\spro{1}}{\F}$ and
$\VS_2=\tuple{\V_2}{\vsum{2}}{\spro{2}}{\F}$ are taken as
elements of $\vecspaces{\V_1}{\F}$ and $\vecspaces{\V_2}{\F}$, respectively. $\anorm{1}$ and
$\anorm{2}$ are taken as elements of $\norms{\VS_1}$ and $\norms{\VS_2}$, respectively. If
$\VS_1$ is of finite dimension, then every linear map from $\VS_1$ to $\VS_2$ is a
linear-continuous map from $\opair{\VS_1}{\anorm{1}}$ to $\opair{\VS_2}{\anorm{2}}$.
\begin{equation}
\func{\Vdim{}}{\VS_1}\in\N\then
\LinC{\opair{\VS_1}{\anorm{1}}}{\opair{\VS_2}{\anorm{2}}}=
\Lin{\VS_1}{\VS_2}.
\end{equation}
\proof

\endthm
%%%%%%%%%%%%%%%%%%%%%%%%%%%%%%%%%%%%%%%%%%%%%%%%%%%%%%%%%%%%%%%%%%%%%%%%%%%%%%%%%%%%%%%%%%%%%%%%%%%%%%%%%%
\definition\label{deflinearcontinuousmapsspace}
Each $\V_1$ and $\V_2$ is taken as a non-empty set, and $\VS_1=\tuple{\V_1}{\vsum{1}}{\spro{1}}{\F}$ and
$\VS_2=\tuple{\V_2}{\vsum{2}}{\spro{2}}{\F}$ are taken as
elements of $\vecspaces{\V_1}{\F}$ and $\vecspaces{\V_2}{\F}$, respectively. $\anorm{1}$ and
$\anorm{2}$ are taken as elements of $\norms{\VS_1}$ and $\norms{\VS_2}$, respectively. Also,
$\NVS{1}:=\opair{\VS_1}{\anorm{1}}$, $\NVS{2}:=\opair{\VS_2}{\anorm{2}}$.
\begin{itemize}
\item
\begin{equation}
\VLinC{\NVS{1}}{\NVS{2}}:=
\tuple{\LinC{\NVS{1}}{\NVS{2}}}{\vsum{lc}}{\spro{lc}}{\F}.
\end{equation}
$\VLinC{\NVS{1}}{\NVS{2}}$ is called the
$\quotl$vector-space of linear-continuous maps from $\NVS{1}$
to $\NVS{2}$.
\item
\begin{align}
&\VLinCnorm{\NVS{1}}{\NVS{2}}\indef\Func{\LinC{\NVS{1}}{\NVS{2}}}{\Rpz},\cr
&\Foreach{\lcf}{\LinC{\NVS{1}}{\NVS{2}}}\VLinCnormf{\lcf}{\NVS{1}}{\NVS{2}}\eqdef
\sup\defSet{\norm{\func{\lcf}{\v}}{2}}{\norm{\v}{1}\leq 1}.
\end{align}
\end{itemize}
\endef
%%%%%%%%%%%%%%%%%%%%%%%%%%%%%%%%%%%%%%%%%%%%%%%%%%%%%%%%%%%%%%%%%%%%%%%%%%%%%%%%%%%%%%%%%%%%%%%%%%%%%%%%%%
\corollary\label{corlinearcontinuousmapnorm}
Each $\NVS{1}=\opair{\VS_1}{\anorm{1}}$, and $\NVS{2}=\opair{\VS_2}{\anorm{2}}$
is taken to be a normed vector-space.
$\VLinCnorm{\NVS{1}}{\NVS{2}}$ is a norm on $\VLinC{\NVS{1}}{\NVS{2}}$.
\begin{equation}
\VLinCnorm{\NVS{1}}{\NVS{2}}\in\norms{\VLinC{\NVS{1}}{\NVS{2}}}.
\end{equation}
\endcor
%%%%%%%%%%%%%%%%%%%%%%%%%%%%%%%%%%%%%%%%%%%%%%%%%%%%%%%%%%%%%%%%%%%%%%%%%%%%%%%%%%%%%%%%%%%%%%%%%%%%%%%%%%
\definition\label{defnormedlinearcontinuousmapspace}
Each $\V_1$ and $\V_2$ is taken as a non-empty set, and $\VS_1=\tuple{\V_1}{\vsum{1}}{\spro{1}}{\F}$ and
$\VS_2=\tuple{\V_2}{\vsum{2}}{\spro{2}}{\F}$ are taken as
elements of $\vecspaces{\V_1}{\F}$ and $\vecspaces{\V_2}{\F}$, respectively. $\anorm{1}$ and
$\anorm{2}$ are taken as elements of $\norms{\VS_1}$ and $\norms{\VS_2}$, respectively. Also,
$\NVS{1}:=\opair{\VS_1}{\anorm{1}}$, $\NVS{2}:=\opair{\VS_2}{\anorm{2}}$.
\begin{equation}
\NVLinC{\NVS{1}}{\NVS{2}}:=\opair{\VLinC{\NVS{1}}{\NVS{2}}}{\VLinCnorm{\NVS{1}}{\NVS{2}}}.
\end{equation}
$\NVLinC{\NVS{1}}{\NVS{2}}$ is called the
$\quotl$normed vector-space of linear-continuous maps from
$\NVS{1}$ to $\NVS{2}$$\quotr$.
\endef
%%%%%%%%%%%%%%%%%%%%%%%%%%%%%%%%%%%%%%%%%%%%%%%%%%%%%%%%%%%%%%%%%%%%%%%%%%%%%%%%%%%%%%%%%%%%%%%%%%%%%%%%%%
\theorem\label{thmlcm1}
Each $\V_1$ and $\V_2$ is taken as a non-empty set, and $\VS_1=\tuple{\V_1}{\vsum{1}}{\spro{1}}{\F}$ and
$\VS_2=\tuple{\V_2}{\vsum{2}}{\spro{2}}{\F}$ are taken as
elements of $\vecspaces{\V_1}{\F}$ and $\vecspaces{\V_2}{\F}$, respectively. $\anorm{1}$ and
$\anorm{2}$ are taken as elements of $\norms{\VS_1}$ and $\norms{\VS_2}$, respectively. Also,
$\NVS{1}:=\opair{\VS_1}{\anorm{1}}$, $\NVS{2}:=\opair{\VS_2}{\anorm{2}}$.
\begin{equation}
\Foreach{\lcf}{\LinC{\NVS{1}}{\NVS{2}}}
\bigg[\Foreach{\v}{\V_1}
\norm{\func{\lcf}{\v}}{2}\leq\VLinCnormf{\lcf}{\NVS{1}}{\NVS{2}}\norm{\v}{1}\bigg].
\end{equation}
\proof
\begin{itemize}
\item[${\textbf{\textsf{p1}}}$]
$\lcf$ is taken as an arbitrary element of $\LinC{\NVS{1}}{\NVS{2}}$, and $\v$ an arbitrary element of\\
$\compl{\V_1}{\seta{\zerov{\VS_1}}}$.
According to \refdef{defnorm},
\begin{align}
\norm{\func{\lcf}{\v}}{2}&=\norm{\func{\lcf}{\norm{\v}{1}\(\frac{1}{\norm{\v}{1}}\v\)}}{2}\cr
&=\norm{\v}{1}\norm{\func{\lcf}{\frac{1}{\norm{\v}{1}}\v}}{2}.
\end{align}
On the other hand, considering that,
\begin{equation}
\norm{\frac{1}{\norm{\v}{1}}\v}{1}=1,
\end{equation}
according to \refdef{deflinearcontinuousmapsspace},
\begin{equation}
\norm{\func{\lcf}{\frac{1}{\norm{\v}{1}}\v}}{2}\leq
\VLinCnormf{\lcf}{\NVS{1}}{\NVS{2}}.
\end{equation}
Thus,
\begin{equation}
\norm{\func{\lcf}{\v}}{2}\leq\VLinCnormf{\lcf}{\NVS{1}}{\NVS{2}}\norm{\v}{1}.
\end{equation}
\endp
\end{itemize}
Additionally, it is obvious that,
\begin{equation}
0=\norm{\func{\lcf}{\zerov{\VS_1}}}{2}
\leq\VLinCnormf{\lcf}{\NVS{1}}{\NVS{2}}\norm{\zerov{\VS_1}}{1}=0.
\end{equation}
So,
\begin{equation}
\Foreach{\v}{\V_1}
\norm{\func{\lcf}{\v}}{2}\leq\VLinCnormf{\lcf}{\NVS{1}}{\NVS{2}}\norm{\v}{1}.
\end{equation}
\endthm
%%%%%%%%%%%%%%%%%%%%%%%%%%%%%%%%%%%%%%%%%%%%%%%%%%%%%%%%%%%%%%%%%%%%%%%%%%%%%%%%%%%%%%%%%%%%%%%%%%%%%%%%%%
\corollary\label{corlcm2}
Each $\V_1$ and $\V_2$ is taken as a non-empty set, and $\VS_1=\tuple{\V_1}{\vsum{1}}{\spro{1}}{\F}$ and
$\VS_2=\tuple{\V_2}{\vsum{2}}{\spro{2}}{\F}$ are taken as
elements of $\vecspaces{\V_1}{\F}$ and $\vecspaces{\V_2}{\F}$, respectively. $\anorm{1}$ and
$\anorm{2}$ are taken as elements of $\norms{\VS_1}$ and $\norms{\VS_2}$, respectively. Also,
$\NVS{1}:=\opair{\VS_1}{\anorm{1}}$, $\NVS{2}:=\opair{\VS_2}{\anorm{2}}$.
\begin{equation}
\Foreach{\lcf}{\LinC{\NVS{1}}{\NVS{2}}}
\VLinCnormf{\lcf}{\NVS{1}}{\NVS{2}}=
\min\defset{M}{\Rpz}{\bigg[\Foreach{\v}{\V_1}
\norm{\func{\lcf}{\v}}{2}\leq M\norm{\v}{1}\bigg]}.
\end{equation}
\endcor
%%%%%%%%%%%%%%%%%%%%%%%%%%%%%%%%%%%%%%%%%%%%%%%%%%%%%%%%%%%%%%%%%%%%%%%%%%%%%%%%%%%%%%%%%%%%%%%%%%%%%%%%%%
\theorem
Each $\V_1$, $\V_2$, and $\V_3$ is taken as a non-empty set, and $\VS_1=\tuple{\V_1}{\vsum{1}}{\spro{3}}{\F}$,
$\VS_2=\tuple{\V_2}{\vsum{2}}{\spro{2}}{\F}$, and $\VS_3=\tuple{\V_3}{\vsum{3}}{\spro{3}}{\F}$ are taken as
elements of $\vecspaces{\V_1}{\F}$, $\vecspaces{\V_2}{\F}$, and $\vecspaces{\V_3}{\F}$ respectively.
$\anorm{1}$, $\anorm{2}$, and $\anorm{3}$ are taken as elements of $\norms{\VS_1}$, $\norms{\VS_2}$,
and $\norms{\VS_3}$, respectively. Also,
$\NVS{1}:=\opair{\VS_1}{\anorm{1}}$, $\NVS{2}:=\opair{\VS_2}{\anorm{2}}$, and
$\NVS{3}:=\opair{\VS_3}{\anorm{3}}$.
\begin{align}
\Foreach{\opair{\lcf}{\p{\lcf}}}{\Cprod{\LinC{\NVS{1}}{\NVS{2}}}{\LinC{\NVS{2}}{\NVS{3}}}}
\VLinCnormf{\cmp{\p{\lcf}}{\lcf}}{\NVS{1}}{\NVS{3}}\leq
\VLinCnormf{\p{\lcf}}{\NVS{2}}{\NVS{3}}\VLinCnormf{\lcf}{\NVS{1}}{\NVS{2}}.
\end{align}
\proof
$\lcf$ and $\p{\lcf}$ are taken as arbitrary elements of
$\LinC{\NVS{1}}{\NVS{2}}$ and $\LinC{\NVS{2}}{\NVS{3}}$, respectively.
According to \refthm{thmlcm1},
\begin{align}
\Foreach{\v}{\V_1}
\norm{\func{\(\cmp{\p{\lcf}}{\lcf}\)}{\v}}{3}&=
\norm{\func{\p{\lcf}}{\func{\lcf}{\v}}}{3}\cr
&\leq\VLinCnormf{\p{\lcf}}{\NVS{2}}{\NVS{3}}
\norm{\func{\lcf}{\v}}{2}\cr
&\leq\(\VLinCnormf{\p{\lcf}}{\NVS{2}}{\NVS{3}}
\VLinCnormf{\lcf}{\NVS{1}}{\NVS{2}}\)\norm{\v}{1}.
\end{align}
So, according to \refcor{corlcm2},
\begin{equation}
\VLinCnormf{\cmp{\p{\lcf}}{\lcf}}{\NVS{1}}{\NVS{3}}\leq
\(\VLinCnormf{\p{\lcf}}{\NVS{2}}{\NVS{3}}
\VLinCnormf{\lcf}{\NVS{1}}{\NVS{2}}\).
\end{equation}
\endthm
%%%%%%%%%%%%%%%%%%%%%%%%%%%%%%%%%%%%%%%%%%%%%%%%%%%%%%%%%%%%%%%%%%%%%%%%%%%%%%%%%%%%%%%%%%%%%%%%%%%%%%%%%%
\theorem\label{thmthmLCcodomainBanach}
Each $\V_1$ and $\V_2$ is taken as a non-empty set, and $\VS_1=\tuple{\V_1}{\vsum{1}}{\spro{1}}{\F}$ and
$\VS_2=\tuple{\V_2}{\vsum{2}}{\spro{2}}{\F}$ are taken as
elements of $\vecspaces{\V_1}{\F}$ and $\vecspaces{\V_2}{\F}$, respectively. $\anorm{1}$ and
$\anorm{2}$ are taken as elements of $\norms{\VS_1}$ and $\norms{\VS_2}$, respectively. Also,
$\NVS{1}:=\opair{\VS_1}{\anorm{1}}$, $\NVS{2}:=\opair{\VS_2}{\anorm{2}}$.
If $\NVS{2}$ is a Banach-space, then $\NVLinC{\NVS{1}}{\NVS{2}}$ is also a Banach-space.
\proof

\endthm
%%%%%%%%%%%%%%%%%%%%%%%%%%%%%%%%%%%%%%%%%%%%%%%%%%%%%%%%%%%%%%%%%%%%%%%%%%%%%%%%%%%%%%%%%%%%%%%%%%%%%%%%%%
%%%%%%%%%%%%%%%%%%%%%%%%%%%%%%%%%%%%%%%%%%%%%%%%%%%%%%%%%%%%%%%%%%%%%%%%%%%%%%%%%%%%%%%%%%%%%%%%%%%%%%%%%%
%%%%%%%%%%%%%%%%%%%%%%%%%%%%%%%%%%%%%%%%%%%%%%%%%%%%%%%%%%%%%%%%%%%%%%%%%%%%%%%%%%%%%%%%%%%%%%%%%%%%%%%%%%
\subsection{Dual of a Normed Vector Space}
\theorem
$\tuple{\f}{\fsum}{\fpro}{\F}$ is a vector-space having the set $\f$ as its' vectors,
and $\F$ as its' field of scalars.
\begin{equation}
\tuple{\f}{\fsum}{\fpro}{\F}\in\vecspaces{\f}{\F}.
\end{equation}
Additionally, $\anorm{\ff}\in\norms{\tuple{\f}{\fsum}{\fpro}{\F}}$, where,
\begin{equation}
\Foreach{\x}{\f}\norm{\x}{\ff}:=\abs{\x}.
\end{equation}
\endthm
%%%%%%%%%%%%%%%%%%%%%%%%%%%%%%%%%%%%%%%%%%%%%%%%%%%%%%%%%%%%%%%%%%%%%%%%%%%%%%%%%%%%%%%%%%%%%%%%%%%%%%%%%%
\definition
$\FF$ is called the $\quotl$canonical normed vector-space of $\F$$\quotr$, where,
\begin{equation}
\FF:=\tuple{\f}{\fsum}{\fpro}{\F}.
\end{equation}
$\opair{\FF}{\anorm{\ff}}$ is called the $\quotl$canonical normed vector-space of $\F$$\quotr$,
and it is denoted by $\NF$.
\endef
%%%%%%%%%%%%%%%%%%%%%%%%%%%%%%%%%%%%%%%%%%%%%%%%%%%%%%%%%%%%%%%%%%%%%%%%%%%%%%%%%%%%%%%%%%%%%%%%%%%%%%%%%%
\theorem\label{thmNFisaBanachspace}
$\NF$ is a Banach-space.
\endthm
%%%%%%%%%%%%%%%%%%%%%%%%%%%%%%%%%%%%%%%%%%%%%%%%%%%%%%%%%%%%%%%%%%%%%%%%%%%%%%%%%%%%%%%%%%%%%%%%%%%%%%%%%%
\definition\label{defdualofnormedspace}
$\V$ is taken as a non-empty set, and $\VS=\tuple{\V}{\vsum{}}{\spro{}}{\F}$ is taken as
an element of $\vecspaces{\V}{\F}$. $\anorm{}$ is taken as an element of $\norms{\VS}$. Also,
$\NVS{}:=\opair{\VS}{\anorm{}}$.
\begin{equation}
\dualNVS{\NVS{}}:=\NVLinC{\NVS{}}{\NF}.
\end{equation}
$\dualNVS{\NVS{}}$ is referred to as the $\quotl$dual-space of ${\NVS{}}$$\quotr$.
\endef
%%%%%%%%%%%%%%%%%%%%%%%%%%%%%%%%%%%%%%%%%%%%%%%%%%%%%%%%%%%%%%%%%%%%%%%%%%%%%%%%%%%%%%%%%%%%%%%%%%%%%%%%%%
\theorem\label{thmdualspaceofNVSisaBanachspace}
$\V$ is taken as a non-empty set, and $\VS=\tuple{\V}{\vsum{}}{\spro{}}{\F}$ is taken as
an element of $\vecspaces{\V}{\F}$. $\anorm{}$ is taken as an element of $\norms{\VS}$. Also,
$\NVS{}:=\opair{\VS}{\anorm{}}$. $\dualNVS{\NVS{}}$ is a Banach-space.
\proof
According to \refthm{thmthmLCcodomainBanach}, and \refthm{thmNFisaBanachspace},
it is clear.
\endthm
%%%%%%%%%%%%%%%%%%%%%%%%%%%%%%%%%%%%%%%%%%%%%%%%%%%%%%%%%%%%%%%%%%%%%%%%%%%%%%%%%%%%%%%%%%%%%%%%%%%%%%%%%%
%%%%%%%%%%%%%%%%%%%%%%%%%%%%%%%%%%%%%%%%%%%%%%%%%%%%%%%%%%%%%%%%%%%%%%%%%%%%%%%%%%%%%%%%%%%%%%%%%%%%%%%%%%
%%%%%%%%%%%%%%%%%%%%%%%%%%%%%%%%%%%%%%%%%%%%%%%%%%%%%%%%%%%%%%%%%%%%%%%%%%%%%%%%%%%%%%%%%%%%%%%%%%%%%%%%%%
\subsection{Isomorphisms Between Normed Vector Spaces}
\definition\label{defnormedisomorphism}
Each $\V_1$ and $\V_2$ is taken as a non-empty set, and $\VS_1=\tuple{\V_1}{\vsum{1}}{\spro{1}}{\F}$ and
$\VS_2=\tuple{\V_2}{\vsum{2}}{\spro{2}}{\F}$ are taken as
elements of $\vecspaces{\V_1}{\F}$ and $\vecspaces{\V_2}{\F}$, respectively. $\anorm{1}$ and
$\anorm{2}$ are taken as elements of $\norms{\VS_1}$ and $\norms{\VS_2}$, respectively. Also,
$\NVS{1}:=\opair{\VS_1}{\anorm{1}}$, $\NVS{2}:=\opair{\VS_2}{\anorm{2}}$.
\begin{itemize}
\item
$\lcf$ is taken as an
element of $\Func{\V_1}{\V_2}$. $\lcf$ is called an $\quotl$isomorphism from the
normed vector-space $\NVS{1}$ to the normed vector-space $\NVS{2}$$\quotr$, iff it posesses these
properties.
\begin{itemize}
\item[${\textbf{\textsf{NIso1}}}$]\quad
$\lcf\in\IF{\V_1}{\V_2}$.
\item[${\textbf{\textsf{NIso2}}}$]\quad
$\lcf\in\LinC{\NVS{1}}{\NVS{2}}$.
\item[${\textbf{\textsf{NIso2}}}$]\quad
$\finv{\lcf}\in\LinC{\NVS{2}}{\NVS{2}}$.
\end{itemize}
\item
The set of all isomorphisms from $\NVS{1}$ to $\NVS{2}$ is denoted by $\NVSiso{\NVS{1}}{\NVS{2}}$.
\begin{equation}
\NVSiso{\NVS{1}}{\NVS{2}}:=\defset{\lcf}{\IF{\V_1}{\V_2}\cap\LinC{\NVS{1}}{\NVS{2}}}
{\finv{\lcf}\in\LinC{\NVS{2}}{\NVS{1}}}.
\end{equation}
\end{itemize}
\endef
%%%%%%%%%%%%%%%%%%%%%%%%%%%%%%%%%%%%%%%%%%%%%%%%%%%%%%%%%%%%%%%%%%%%%%%%%%%%%%%%%%%%%%%%%%%%%%%%%%%%%%%%%%
\corollary\label{cornormedisomorphismequiv1}
Each $\V_1$ and $\V_2$ is taken as a non-empty set, and $\VS_1=\tuple{\V_1}{\vsum{1}}{\spro{1}}{\F}$ and
$\VS_2=\tuple{\V_2}{\vsum{2}}{\spro{2}}{\F}$ are taken as
elements of $\vecspaces{\V_1}{\F}$ and $\vecspaces{\V_2}{\F}$, respectively. $\anorm{1}$ and
$\anorm{2}$ are taken as elements of $\norms{\VS_1}$ and $\norms{\VS_2}$, respectively. Also,
$\NVS{1}:=\opair{\VS_1}{\anorm{1}}$, $\NVS{2}:=\opair{\VS_2}{\anorm{2}}$.
\begin{equation}
\NVSiso{\NVS{1}}{\NVS{2}}:=\Lin{\VS_1}{\VS_2}\cap\HOF{\opair{\V_1}{\func{\normtop{\F}{\V_1}}{\anorm{1}}}}
{\opair{\V_2}{\func{\normtop{\F}{\V_2}}{\anorm{2}}}}.
\end{equation}
\endcor
%%%%%%%%%%%%%%%%%%%%%%%%%%%%%%%%%%%%%%%%%%%%%%%%%%%%%%%%%%%%%%%%%%%%%%%%%%%%%%%%%%%%%%%%%%%%%%%%%%%%%%%%%%
\theorem
Each $\V_1$ and $\V_2$ is taken as a non-empty set, and $\VS_1=\tuple{\V_1}{\vsum{1}}{\spro{1}}{\F}$ and
$\VS_2=\tuple{\V_2}{\vsum{2}}{\spro{2}}{\F}$ are taken as
elements of $\vecspaces{\V_1}{\F}$ and $\vecspaces{\V_2}{\F}$, respectively. $\anorm{1}$ and
$\anorm{2}$ are taken as elements of $\norms{\VS_1}$ and $\norms{\VS_2}$, respectively. Also,
$\NVS{1}:=\opair{\VS_1}{\anorm{1}}$, $\NVS{2}:=\opair{\VS_2}{\anorm{2}}$. If both $\NVS{1}$ and
$\NVS{2}$ are Banach-spaces, then
\begin{equation}
\NVSiso{\NVS{1}}{\NVS{2}}=
\IF{\V_1}{\V_2}\cap\LinC{\NVS{1}}{\NVS{2}}.
\end{equation}
\endthm
%%%%%%%%%%%%%%%%%%%%%%%%%%%%%%%%%%%%%%%%%%%%%%%%%%%%%%%%%%%%%%%%%%%%%%%%%%%%%%%%%%%%%%%%%%%%%%%%%%%%%%%%%%%%%%%%%%%%%%%%%%%%%%%%%%%%%%%%%%%%
%%%%%%%%%%%%%%%%%%%%%%%%%%%%%%%%%%%%%%%%%%%%%%%%%%%%%%%%%%%%%%%%%%%%%%%%%%%%%%%%%%%%%%%%%%%%%%%%%%%%%%%%%%%%%%%%%%%%%%%%%%%%%%%%%%%%%%%%%%%%
%%%%%%%%%%%%%%%%%%%%%%%%%%%%%%%%%%%%%%%%%%%%%%%%%%%%%%%%%%%%%%%%%%%%%%%%%%%%%%%%%%%%%%%%%%%%%%%%%%%%%%%%%%%%%%%%%%%%%%%%%%%%%%%%%%%%%%%%%%%%
%%%%%%%%%%%%%%%%%%%%%%%%%%%%%%%%%%%%%%%%%%%%%%%%%%%%%%%%%%%%%%%%%%%%%%%%%%%%%%%%%%%%%%%%%%%%%%%%%%%%%%%%%%%%%%%%%%%%%%%%%%%%%%%%%%%%%%%%%%%%
%%%%%%%%%%%%%%%%%%%%%%%%%%%%%%%%%%%%%%%%%%%%%%%%%%%%%%%%%%%%%%%%%%%%%%%%%%%%%%%%%%%%%%%%%%%%%%%%%%%%%%%%%%%%%%%%%%%%%%%%%%%%%%%%%%%%%%%%%%%%
\section{Multilinear-Continuous Maps}
%%%%%%%%%%%%%%%%%%%%%%%%%%%%%%%%%%%%%%%%%%%%%%%%%%%%%%%%%%%%%%%%%%%%%%%%%%%%%%%%%%%%%%%%%%%%%%%%%%%%%%%%%%
\definition
$\collection{}$ is taken as a set, such that $\empty\notin\collection{}$.
\begin{equation}
\Cproduct{\collection{}}:=\defset{\alphaa{}}{\Func{\collection{}}{\union{\collection{}}}}
{\Foreach{\aset{}}{\collection{}}\func{\alphaa{}}{\aset{}}\in\aset{}}.
\end{equation}
$\Cproduct{\VV{}}$ is referred to as the $\quotl$Cartesian-product of $\collection{}$.
\footnote{Axiom of choice states that $\Cproduct{\collection{}}\neq\empty$.}
\endef
%%%%%%%%%%%%%%%%%%%%%%%%%%%%%%%%%%%%%%%%%%%%%%%%%%%%%%%%%%%%%%%%%%%%%%%%%%%%%%%%%%%%%%%%%%%%%%%%%%%%%%%%%%
\definition\label{defCartesianproduct}
Each $\collection{}$ and $\index$ is taken as a set, such that $\empty\notin\collection{}$.
$\VV{}$ is taken as an element of
$\IF{\index}{\collection{}}$ (a bijective map from $\index$ to $\collection{}$).
\begin{equation}
\Cproduct{\VV{}}:=\defset{\alphaa{}}{\Func{\index}{\union{\collection{}}}}
{\Foreach{i}{\index}\func{\alphaa{}}{i}\in\func{\VV{}}{i}}.
\end{equation}
$\Cproduct{\VV{}}$ is referred to as the $\quotl$Cartesian-product of $\VV{}$$\quotr$.
For clarity, the symbol $\sindex{\VV{i}}{i}{\findex{}}$ is allowed to stand for $\VV{}$, and
for each $i$ in $\index$, $\VV{i}$ stands for $\func{\VV{}}{i}$.
Also, for every $\alphaa{}$ in $\Cproduct{\VV{}}$ and
for each $i$ in $\index$, $\alphaa{i}$ stands for $\func{\alphaa{}}{i}$.
\endef
%%%%%%%%%%%%%%%%%%%%%%%%%%%%%%%%%%%%%%%%%%%%%%%%%%%%%%%%%%%%%%%%%%%%%%%%%%%%%%%%%%%%%%%%%%%%%%%%%%%%%%%%%%
\definition\label{defCartesianproduct}
Each $\collection{}$ and $\index$ is taken as a set, such that $\empty\notin\collection{}$.
$\VV{}$ is taken as an element of
$\IF{\index}{\collection{}}$ (a bijective map from $\index$ to $\collection{}$).
\begin{equation}
\Foreachs{I}{\index}
\sindex{\VV{i}}{i}{I}:=
\Cproduct{\func{\rescd{\func{\resd{\VV{}}}{I}}}{\func{\image{\VV{}}}{I}}}.
\end{equation}
Also,
\begin{equation}
\Foreach{k}{\index}
\seta{\VV{k}}_{k}:=\sindex{\VV{i}}{i}{\seta{k}}.
\end{equation}
\endef
%%%%%%%%%%%%%%%%%%%%%%%%%%%%%%%%%%%%%%%%%%%%%%%%%%%%%%%%%%%%%%%%%%%%%%%%%%%%%%%%%%%%%%%%%%%%%%%%%%%%%%%%%%
\definition
For each $n$ in $\N$, $\findex{n}:=\seta{\suc{1}{n}}$.
\endef
%%%%%%%%%%%%%%%%%%%%%%%%%%%%%%%%%%%%%%%%%%%%%%%%%%%%%%%%%%%%%%%%%%%%%%%%%%%%%%%%%%%%%%%%%%%%%%%%%%%%%%%%%%
\definition
$\collection{}$ is taken as a set, such that $\empty\notin\collection{}$.
$n$ is taken as an element of $\Zp$.
$\VV{}$ is taken as an element of
$\IF{\findex{n}}{\collection{}}$ (a bijective map from $\findex{n}$ to $\collection{}$).
\begin{equation}
\Times{\VV{1}}{\VV{n}}:=\Cproduct{\VV{}},
\end{equation}
and,
\begin{equation}
\Foreach{\alphaa{}}{\Cproduct{\VV{}}}
\mtuple{\alphaa{1}}{\alphaa{n}}:=\alphaa{}.
\end{equation}
\endef
%%%%%%%%%%%%%%%%%%%%%%%%%%%%%%%%%%%%%%%%%%%%%%%%%%%%%%%%%%%%%%%%%%%%%%%%%%%%%%%%%%%%%%%%%%%%%%%%%%%%%%%%%%
\definition\label{defproductvectorspaceinjection}
$\sindex{\VV{i}}{i}{\findex{n}}$
is taken to be an indexed-collection of non-empty sets, for some $n$ in $\N$.
For each $i$ in $\seta{\suc{1}{n}}$, $\VVS{i}=\tuple{\VV{i}}{\vsum{i}}{\spro{i}}{\F}$
is taken as an element of $\vecspaces{\VV{i}}{\F}$.
For simplicity, $\indexedcollection{V}$ stands for
$\mtuple{\VVS{1}}{\VVS{n}}$.
\begin{align}
&\Vprodinj{\indexedcollection{V}}\indef
\Func{\Union{I}{\CSs{\seta{\suc{1}{n}}}}{\Cproduct{\sindex{\VV{i}}{i}{I}}}}
{\Cproduct{\sindex{\VV{i}}{i}{\findex{n}}}},\cr
&\Foreachs{I}{\seta{\suc{1}{n}}}
\Foreach{\alphaa{}}{\sindex{\VV{i}}{i}{I}}
\func{\[\func{\Vprodinj{\indexedcollection{V}}}{\alphaa{}}\]}{k}\eqdef
\begin{cases}
\func{\alphaa{}}{k},\, & k\in I,\cr
\zerov{\VVS{k}},\, & k\in\(\compl{\findex{n}}{I}\).
\end{cases}
\end{align}
\endef
%%%%%%%%%%%%%%%%%%%%%%%%%%%%%%%%%%%%%%%%%%%%%%%%%%%%%%%%%%%%%%%%%%%%%%%%%%%%%%%%%%%%%%%%%%%%%%%%%%%%%%%%%%
\definition\label{defmultilinearmap}
%$\seta{\suc{\VV{1}}{\VV{n+1}}}$
$\sindex{\VV{i}}{i}{\findex{n}}$
is taken to be an indexed-collection of non-empty sets, for some $n$ in $\N$.
For each $i$ in $\seta{\suc{1}{n}}$, $\VVS{i}=\tuple{\VV{i}}{\vsum{i}}{\spro{i}}{\F}$
is taken as an element of $\vecspaces{\VV{i}}{\F}$.
$\WW{}$ is taken as a non-empty set, and $\WS{}=\tuple{\WW{}}{\vsum{w}}{\spro{w}}{\F}$
an element of $\vecspaces{\WW{}}{\F}$. For convenience, $\indexedcollection{V}$ stands for
$\mtuple{\VVS{1}}{\VVS{n}}$.
Also, $\VVS{1}\times\ldots\times\VVS{n}=
\tuple{\displaystyle\Cproduct{\sindex{\VV{i}}{i}{\findex{n}}}}{\vsum{}}{\spro{}}{\F}$,
where $\VVS{1}\times\ldots\times\VVS{n}$ stands for the product vector-space of the considered
vector-spaces.
%and $\anorm{i}$ an element of $\norms{\VVS{i}}$.
\begin{itemize}
\item
$\ml{}$ is taken as an element of $\displaystyle\Func{\Cproduct{\sindex{\VV{i}}{i}{\findex{n}}}}{\WW{}}$.
$\ml{}$ is referred to as a $\quotl$multilinear map from $\mtuple{\VVS{1}}{\VVS{n}}$ to $\WS{}$$\quotr$,
iff,
\begin{align}
&\Foreach{k}{\findex{n}}
\Foreach{\alphaa{}}{\Cproduct{\sindex{\VV{i}}{i}{\(\compl{\findex{n}}{\seta{k}}\)}}}
\Foreach{\triple{\x}{\v}{\w}}{\Cprod{\f}{\Cprod{\Cproduct{\seta{\VV{k}}_{k}}}{\Cproduct{\seta{\VV{k}}_{k}}}}}\cr
&\begin{aligned}
\func{\ml{}}{\func{\Vprodinj{\indexedcollection{V}}}{\alphaa{}}+
\func{\Vprodinj{\indexedcollection{V}}}{\x\spro{i}\v\vsum{i}\w}}=&\quad
\x\(\spro{w}\)\func{\ml{}}{\func{\Vprodinj{\indexedcollection{V}}}{\alphaa{}}+
\func{\Vprodinj{\indexedcollection{V}}}{\v}}\cr
&\(\vsum{w}\)
\func{\ml{}}{\func{\Vprodinj{\indexedcollection{V}}}{\alphaa{}}+
\func{\Vprodinj{\indexedcollection{V}}}{\w}}.\cr
\end{aligned}\cr
&{}
\end{align}
In other words, $\ml{}$ is defined to be a multilinear map from $\mtuple{\VVS{1}}{\VVS{n}}$ to $\WS{}$,
iff, for each $k$ in $\findex{n}$, and for every $u_i\in\VV{i}$ such that $i\in\findex{n}$ and $i\neq k$,
\begin{align}
&\Foreach{\triple{\x}{\v}{\w}}{\Cprod{\f}{\Cprod{\VV{k}}{\VV{k}}}}\cr
&\begin{aligned}
\func{\ml{}}{\binary{\suc{u_1}{u_{k-1}}}{\binary{\x\spro{k}\v\vsum{k}\w}{\suc{u_{k+1}}{u_n}}}}=&\quad
\x\spro{w}\func{\ml{}}{\binary{\suc{u_1}{u_{k-1}}}{\binary{\v}{\suc{u_{k+1}}{u_n}}}}\cr
&\vsum{w}
\func{\ml{}}{\binary{\suc{u_1}{u_{k-1}}}{\binary{\w}{\suc{u_{k+1}}{u_n}}}}.
\end{aligned}\cr
&{}
\end{align}
\item
The set of all multilinear maps from $\mtuple{\VVS{1}}{\VVS{n}}$
to $\WS{}$ is denoted by
$\MLin{\suc{\VVS{1}}{\VVS{n}}}{\WS{}}$.
\end{itemize}
\endef
%%%%%%%%%%%%%%%%%%%%%%%%%%%%%%%%%%%%%%%%%%%%%%%%%%%%%%%%%%%%%%%%%%%%%%%%%%%%%%%%%%%%%%%%%%%%%%%%%%%%%%%%%%
\theorem\label{thmmultilinearmapsspace}
$\sindex{\VV{i}}{i}{\findex{n}}$
is taken to be an indexed-collection of non-empty sets, for some $n$ in $\N$.
For each $i$ in $\seta{\suc{1}{n}}$, $\VVS{i}=\tuple{\VV{i}}{\vsum{i}}{\spro{i}}{\F}$
is taken as an element of $\vecspaces{\VV{i}}{\F}$.
$\WW{}$ is taken as a non-empty set, and $\WS{}=\tuple{\WW{}}{\vsum{w}}{\spro{w}}{\F}$
an element of $\vecspaces{\WW{}}{\F}$. For convenience, $\indexedcollection{V}$ stands for
$\mtuple{\VVS{1}}{\VVS{n}}$.
Also, $\VVS{1}\times\ldots\times\VVS{n}=
\tuple{\displaystyle\Cproduct{\sindex{\VV{i}}{i}{\findex{n}}}}{\vsum{}}{\spro{}}{\F}$,
where $\VVS{1}\times\ldots\times\VVS{n}$ stands for the product vector-space of the considered
vector-spaces.
\begin{align}
\tuple{\MLin{\suc{\VVS{1}}{\VVS{n}}}{\WS{}}}{\vsum{ml}}{\spro{ml}}{\F}\in
\vecspaces{\MLin{\suc{\VVS{1}}{\VVS{n}}}{\WS{}}}{\F},
\end{align}
where,
\begin{align}
\vsum{ml}&:=\func{\resd{\vsum{f}}}{\Cprod{\MLin{\suc{\NVS{1}}{\NVS{n}}}{\NWS{}}}{\MLin{\suc{\NVS{1}}{\NVS{n}}}{\NWS{}}}},\\
\spro{ml}&:=\func{\resd{\spro{f}}}{\Cprod{\f}{\MLin{\suc{\NVS{1}}{\NVS{n}}}{\NWS{}}}},
\end{align}
%\Cprod{\Func{\Times{\VV{1}}{\VV{n}}}{\W{}}}{\Func{\Times{\VV{1}}{\VV{n}}}{\W{}}}
$\vsum{f}$ and $\spro{f}$ being the addition and scalar-multiplication operations of
the vector-space of all functions from $\Times{\VV{1}}{\VV{n}}$ to the vector-space $\WS{}$, that is
$\VFunc{\Times{\VV{1}}{\VV{n}}}{\WS{}}$.
\proof
$\MLin{\suc{\VVS{1}}{\VVS{n}}}{\WS{}}$ is a non-empty subset of
$\VFunc{\Times{\VV{1}}{\VV{n}}}{\WS{}}$, since the constant map that maps every element of
$\Times{\VV{1}}{\VV{n}}$ to $\zerov{\WS{}}$ is a multilinear map from $\mtuple{\VVS{1}}{\VVS{n}}$
to $\WS{}$.
%It can be easily verified that addition of two arbitrary linear-continuous maps is a continuous map,
%and any scalar-multiplication of an arbitrary linear-continuous map is again a continuous map.
%So, $\triple{\LinC{\opair{\VS_1}{\anorm{1}}}{\opair{\VS_2}{\anorm{2}}}}{\vsum{lc}}{\spro{lc}}$
%is a linear-subspace of $\VLin{\VS_1}{\VS_2}$.
\begin{itemize}
\item[${\textbf{\textsf{p1}}}$]
Each $\ml{1}$ and $\ml{2}$ is taken as an element of
$\MLin{\suc{\VVS{1}}{\VVS{n}}}{\WS{}}$, and $c$ an element of $\f$. According to \refdef{defmultilinearmap},
\begin{align}
&\Foreach{k}{\findex{n}}
\Foreach{\alphaa{}}{\Cproduct{\sindex{\VV{i}}{i}{\(\compl{\findex{n}}{\seta{k}}\)}}}
\Foreach{\triple{\x}{\v}{\w}}{\Cprod{\f}{\Cprod{\Cproduct{\seta{\VV{k}}_{k}}}{\Cproduct{\seta{\VV{k}}_{k}}}}}\cr
&\begin{aligned}
&\quad\func{\[\(c\spro{f}\ml{1}\)\vsum{f}\ml{2}\]}{\func{\Vprodinj{\indexedcollection{V}}}{\alphaa{}}+
\func{\Vprodinj{\indexedcollection{V}}}{\x\spro{i}\v\vsum{i}\w}}\cr
&=c\spro{w}\func{\ml{1}}{\func{\Vprodinj{\indexedcollection{V}}}{\alphaa{}}+
\func{\Vprodinj{\indexedcollection{V}}}{\x\spro{i}\v\vsum{i}\w}}\vsum{w}
\func{\ml{2}}{\func{\Vprodinj{\indexedcollection{V}}}{\alphaa{}}+
\func{\Vprodinj{\indexedcollection{V}}}{\x\spro{i}\v\vsum{i}\w}}\cr
&=\quad c\spro{w}\bigg[\x\spro{w}\func{\ml{1}}{\func{\Vprodinj{\indexedcollection{V}}}{\alphaa{}}+
\func{\Vprodinj{\indexedcollection{V}}}{\v}}
\vsum{w}
\func{\ml{1}}{\func{\Vprodinj{\indexedcollection{V}}}{\alphaa{}}+
\func{\Vprodinj{\indexedcollection{V}}}{\w}}\bigg]\cr
&\quad\vsum{w}\bigg[\x\spro{w}\func{\ml{2}}{\func{\Vprodinj{\indexedcollection{V}}}{\alphaa{}}+
\func{\Vprodinj{\indexedcollection{V}}}{\v}}
\vsum{w}
\func{\ml{2}}{\func{\Vprodinj{\indexedcollection{V}}}{\alphaa{}}+
\func{\Vprodinj{\indexedcollection{V}}}{\w}}\bigg]\cr
&=\quad\x\spro{w}\bigg[c\spro{w}\func{\ml{1}}{\func{\Vprodinj{\indexedcollection{V}}}{\alphaa{}}+
\func{\Vprodinj{\indexedcollection{V}}}{\v}}\vsum{w}\func{\ml{2}}{\func{\Vprodinj{\indexedcollection{V}}}{\alphaa{}}+
\func{\Vprodinj{\indexedcollection{V}}}{\v}}\bigg]\cr
&\quad\vsum{w}\bigg[c\spro{w}\func{\ml{1}}{\func{\Vprodinj{\indexedcollection{V}}}{\alphaa{}}+
\func{\Vprodinj{\indexedcollection{V}}}{\w}}\vsum{w}\func{\ml{2}}{\func{\Vprodinj{\indexedcollection{V}}}{\alphaa{}}+
\func{\Vprodinj{\indexedcollection{V}}}{\w}}\bigg]\cr
&=\x\spro{w}\func{\[\(c\spro{f}\ml{1}\)\vsum{f}\ml{2}\]}{\func{\Vprodinj{\indexedcollection{V}}}{\alphaa{}}+
\func{\Vprodinj{\indexedcollection{V}}}{\v}}
\vsum{w}
\func{\[\(c\spro{f}\ml{1}\)\vsum{f}\ml{2}\]}{\func{\Vprodinj{\indexedcollection{V}}}{\alphaa{}}+
\func{\Vprodinj{\indexedcollection{V}}}{\w}}.
\end{aligned}\cr
&{}
\end{align}
So, according to \refdef{defmultilinearmap},
\begin{equation}
\[\(c\spro{f}\ml{1}\)\vsum{f}\ml{2}\]\in\MLin{\suc{\VVS{1}}{\VVS{n}}}{\WS{}}.
\end{equation}
\endp
\end{itemize}
\endthm
%%%%%%%%%%%%%%%%%%%%%%%%%%%%%%%%%%%%%%%%%%%%%%%%%%%%%%%%%%%%%%%%%%%%%%%%%%%%%%%%%%%%%%%%%%%%%%%%%%%%%%%%%%
\definition\label{defproducttopologyofNVSs}
$\sindex{\VV{i}}{i}{\findex{n}}$
is taken to be an indexed-collection of non-empty sets, for some $n$ in $\N$.
\begin{align}
&\prodnormtop{\F}{\Cproduct{\sindex{\VV{i}}{i}{\findex{n}}}}\indef
\Func{\Times{\setnorms{\F}{\VV{1}}}{\setnorms{\F}{\VV{n}}}}
{\Ctops{\Cproduct{\sindex{\VV{i}}{i}{\findex{n}}}}},\cr
&\begin{aligned}
&\Foreach{\mtuple{\anorm{1}}{\anorm{n}}}{\Times{\setnorms{\F}{\VV{1}}}{\setnorms{\F}{\VV{n}}}}\cr
&\func{\prodnormtop{\F}{\Cproduct{\sindex{\VV{i}}{i}{\findex{n}}}}}{\mult{\anorm{1}}{\anorm{n}}}\eqdef
\func{\prodtop{\Cproduct{\sindex{\VV{i}}{i}{\findex{n}}}}}
{\mult{\func{\normtop{\F}{\VV{1}}}{\anorm{1}}}{\func{\normtop{\F}{\VV{2}}}{\anorm{n}}}}.
\end{aligned}
\end{align}
$\func{\prodtop{\Cproduct{\sindex{\VV{i}}{i}{\findex{n}}}}}
{\mult{\func{\normtop{\F}{\VV{1}}}{\anorm{1}}}{\func{\normtop{\F}{\VV{2}}}{\anorm{n}}}}$
denotes the topology of the product topological space
$\Times{\opair{\VV{1}}{\func{\normtop{\F}{\VV{1}}}{\anorm{1}}}}
{\opair{\VV{1}}{\func{\normtop{\F}{\VV{n}}}{\anorm{n}}}}$.\\
So, having the vector-spaces $\VVS{1}\in\vecspaces{\VV{1}}{\F}$,$\cdots$,$\VVS{n}\in\vecspaces{\VV{n}}{\F}$,
and the norms $\anorm{1}$,$\cdots$,$\anorm{n}$ on the vector-spaces $\VVS{1}$,$\ldots$,$\VVS{n}$,
respectively,
$\func{\prodnormtop{\F}{\Cproduct{\sindex{\VV{i}}{i}{\findex{n}}}}}{\mult{\anorm{1}}{\anorm{n}}}$
is the product topology of the induced topologies from $\anorm{i}$s on $\VV{i}$s.
The topological space\\ $\opair{\Times{\VV{1}}{\VV{n}}}
{\func{\prodnormtop{\F}{\Cproduct{\sindex{\VV{i}}{i}{\findex{n}}}}}{\mult{\anorm{1}}{\anorm{n}}}}$
is called the $\quotl$canonical topological space induced by the norms $\anorm{i}$s$\quotr$,
and denoted briefly by $\prodnormedtopspace{\suc{\anorm{1}}{\anorm{n}}}$.
\endef
%%%%%%%%%%%%%%%%%%%%%%%%%%%%%%%%%%%%%%%%%%%%%%%%%%%%%%%%%%%%%%%%%%%%%%%%%%%%%%%%%%%%%%%%%%%%%%%%%%%%%%%%%%
\definition\label{defmulticontinuousmap}
$\sindex{\VV{i}}{i}{\findex{n}}$
is taken to be an indexed-collection of non-empty sets, for some $n$ in $\N$.
For each $i$ in $\seta{\suc{1}{n}}$, $\VVS{i}=\tuple{\VV{i}}{\vsum{i}}{\spro{i}}{\F}$
is taken as an element of $\vecspaces{\VV{i}}{\F}$, and $\anorm{i}$
an element of $\norms{\VVS{i}}$.
$\WW{}$ is taken as a non-empty set.
$\WS{}=\tuple{\WW{}}{\vsum{w}}{\spro{w}}{\F}$ is taken as
an element of $\vecspaces{\WW{}}{\F}$, and $\anorm{}$ an element of $\norms{\WS{}}$.
Also, for each $i$ in $\findex{n}$, $\NVS{i}:=\opair{\VVS{i}}{\anorm{i}}$, and
$\NWS{}:=\opair{\WS{}}{\anorm{}}$.
\begin{itemize}
\item
\begin{align}
&\MCF{\suc{\NVS{1}}{\NVS{n}}}{\NWS{}}:=\cr
&\CF{\opair{\Cproduct{\sindex{\VV{i}}{i}{\findex{n}}}}
{\func{\prodnormtop{\F}{\Cproduct{\sindex{\VV{i}}{i}{\findex{n}}}}}{\mult{\anorm{1}}{\anorm{n}}}}}
{\opair{\WW{}}{\func{\normtop{\F}{\WW{}}}{\anorm{}}}}.
\end{align}
Each element of $\MCF{\suc{\NVS{1}}{\NVS{n}}}{\NWS{}}$ is referred to as a
$\quotl$multi-continuous map from the $n$-tuple of normed-spaces
$\mtuple{\NVS{1}}{\NVS{n}}$ to $\NWS{}$$\quotr$.\\
So, a map from $\Times{\VV{1}}{\VV{n}}$ to the normed-space
$\WW{}$ is called a multi-continuous map
from $\mtuple{\NVS{1}}{\NVS{n}}$ to $\NWS{}$, iff it is a continuous map
from the canonical topological space induced by the norms $\anorm{i}$s
to the canonical topological space induced by the norm $\anorm{}$.
\item
$\cf$ is taken as an element of $\Func{\Times{\VV{1}}{\VV{n}}}{\WW{}}$, and $\point$ an element of
$\Times{\VV{1}}{\VV{n}}$.
$\cf$ is defined to be a $\quotl$multi-continuous map from $\mtuple{\NVS{1}}{\NVS{n}}$ to
$\NWS{}$ at $\point$$\quotr$ iff $\cf$ is a continuous map from the topological space
$\opair{\Cproduct{\sindex{\VV{i}}{i}{\findex{n}}}}
{\func{\prodnormtop{\F}{\Cproduct{\sindex{\VV{i}}{i}{\findex{n}}}}}{\mult{\anorm{1}}{\anorm{n}}}}$ to the
topological space $\opair{\WW{}}{\func{\normtop{\F}{\WW{}}}{\anorm{}}}$ at $\point$.
It can also simply be said that $\quotl$$\cf$ is multi-continuous at $\point$$\quotr$.
\end{itemize}
\endef
%%%%%%%%%%%%%%%%%%%%%%%%%%%%%%%%%%%%%%%%%%%%%%%%%%%%%%%%%%%%%%%%%%%%%%%%%%%%%%%%%%%%%%%%%%%%%%%%%%%%%%%%%%
\corollary\label{cormulticontinuityequivs}
$\sindex{\VV{i}}{i}{\findex{n}}$
is taken to be an indexed-collection of non-empty sets, for some $n$ in $\N$.
For each $i$ in $\seta{\suc{1}{n}}$, $\VVS{i}=\tuple{\VV{i}}{\vsum{i}}{\spro{i}}{\F}$
is taken as an element of $\vecspaces{\VV{i}}{\F}$, and $\anorm{i}$
an element of $\norms{\VVS{i}}$.
$\WW{}$ is taken as a non-empty set.
$\WS{}=\tuple{\WW{}}{\vsum{w}}{\spro{w}}{\F}$ is taken as
an element of $\vecspaces{\WW{}}{\F}$, and $\anorm{}$ an element of $\norms{\WS{}}$.
Also, for each $i$ in $\findex{n}$, $\NVS{i}:=\opair{\VVS{i}}{\anorm{i}}$, and
$\NWS{}:=\opair{\WS{}}{\anorm{}}$. $\cf$ is taken as an element of
$\Func{\Times{\VV{1}}{\VV{n}}}{\WW{}}$.
\begin{itemize}
\item
For every $\point$ in $\Times{\VV{1}}{\VV{n}}$, $\cf$ is multi-continuous at $\point$ if-and-only-if,
\begin{equation}
\Foreach{\varepsilon}{\Rp}
\Exists{\delta}{\Rp}\func{\image{\cf}}{\Times{\func{\ball{\NVS{1}}}{\binary{\pnt{1}}{\delta}}}
{\func{\ball{\NVS{n}}}{\binary{\pnt{n}}{\delta}}}}\subseteq
\func{\ball{\NWS{}}}{\binary{\func{\cf}{\pnt{}}}{\varepsilon}}.
\end{equation}
\item
$\cf$ is multi-continuous if-and-only-if it is multi-continuous at every $\pnt{}$ in
$\Times{\VV{1}}{\VV{n}}$.
\begin{gather}
\cf\in\MCF{\suc{\NVS{1}}{\NVS{n}}}{\NWS{}}\cr
\vthenn\cr
\[\begin{aligned}
&\Foreach{\pnt{}}{\Times{\VV{1}}{\VV{n}}}\cr
&\Foreach{\varepsilon}{\Rp}
\Exists{\delta}{\Rp}\func{\image{\cf}}{\Times{\func{\ball{\NVS{1}}}{\binary{\pnt{1}}{\delta}}}
{\func{\ball{{\NVS{n}}}}{\binary{\pnt{n}}{\delta}}}}\subseteq
\func{\ball{\NWS{}}}{\binary{\func{\cf}{\point}}{\varepsilon}}
\end{aligned}\].
\end{gather}
\end{itemize}
\endcor
%%%%%%%%%%%%%%%%%%%%%%%%%%%%%%%%%%%%%%%%%%%%%%%%%%%%%%%%%%%%%%%%%%%%%%%%%%%%%%%%%%%%%%%%%%%%%%%%%%%%%%%%%%
\definition\label{defmultilinearcontinuousmap}
$\sindex{\VV{i}}{i}{\findex{n}}$
is taken to be an indexed-collection of non-empty sets, for some $n$ in $\N$.
For each $i$ in $\seta{\suc{1}{n}}$, $\VVS{i}=\tuple{\VV{i}}{\vsum{i}}{\spro{i}}{\F}$
is taken as an element of $\vecspaces{\VV{i}}{\F}$, and $\anorm{i}$
an element of $\norms{\VVS{i}}$.
$\WW{}$ is taken as a non-empty set.
$\WS{}=\tuple{\WW{}}{\vsum{w}}{\spro{w}}{\F}$ is taken as
an element of $\vecspaces{\WW{}}{\F}$, and $\anorm{}$ an element of $\norms{\WS{}}$.
Also, for each $i$ in $\findex{n}$, $\NVS{i}:=\opair{\VVS{i}}{\anorm{i}}$, and
$\NWS{}:=\opair{\WS{}}{\anorm{}}$.
\begin{equation}
\MLC{\suc{\NVS{1}}{\NVS{n}}}{\NWS{}}:=
\MLin{\suc{\VVS{1}}{\VVS{n}}}{\WS{}}\cap
\MCF{\suc{\NVS{1}}{\NVS{2}}}{\NWS{}}.
\end{equation}
Every element of $\MLC{\suc{\NVS{1}}{\NVS{n}}}{\NWS{}}$ is
called a $\quotl$multilinear-continuous map from the $n$-tuple of normed-spaces
$\mtuple{\NVS{1}}{\NVS{n}}$ to the normed-space $\NWS{}$$\quotr$.\\
In other words, a map $\cf:\Times{\VV{1}}{\VV{n}}\to\WW{}$ is called a multilinear-continuous map from
$\mtuple{\NVS{1}}{\NVS{n}}$ to $\NWS{}$, iff
it is both a multilinear map from $\mtuple{\VVS{1}}{\VVS{n}}$ to $\WS{}$, and a continuous map from
$\mtuple{\NVS{1}}{\NVS{n}}$ to $\NWS{}$.
\endef
%%%%%%%%%%%%%%%%%%%%%%%%%%%%%%%%%%%%%%%%%%%%%%%%%%%%%%%%%%%%%%%%%%%%%%%%%%%%%%%%%%%%%%%%%%%%%%%%%%%%%%%%%%
\theorem\label{thmmultilinearcontinuousmapequivs}
$\sindex{\VV{i}}{i}{\findex{n}}$
is taken to be an indexed-collection of non-empty sets, for some $n$ in $\N$.
For each $i$ in $\seta{\suc{1}{n}}$, $\VVS{i}=\tuple{\VV{i}}{\vsum{i}}{\spro{i}}{\F}$
is taken as an element of $\vecspaces{\VV{i}}{\F}$, and $\anorm{i}$
an element of $\norms{\VVS{i}}$.
$\WW{}$ is taken as a non-empty set.
$\WS{}=\tuple{\WW{}}{\vsum{w}}{\spro{w}}{\F}$ is taken as
an element of $\vecspaces{\WW{}}{\F}$, and $\anorm{}$ an element of $\norms{\WS{}}$.
Also, for each $i$ in $\findex{n}$, $\NVS{i}:=\opair{\VVS{i}}{\anorm{i}}$, and
$\NWS{}:=\opair{\WS{}}{\anorm{}}$. $\ml{}$ is taken as an element of $\MLin{\suc{\VVS{1}}{\VVS{n}}}{\WS{}}$.
These propositions are equivalent.
\begin{itemize}
\item[${\propos{1}}.$]\quad
$\ml{}\in\MLC{\suc{\NVS{1}}{\NVS{n}}}{\NWS{}}$.
\item[${\propos{2}}.$]\quad
$\ml{}$ is multiliniear-continuous at $\mtuple{\zerov{\VVS{1}}}{\zerov{\VVS{n}}}$.
\item[${\propos{3}}.$]\quad
$\Exists{K}{\Rp}\Foreach{\vv{}}{\defset{\p{\vv{}}}{\Times{\VV{1}}{\VV{n}}}{\[\Foreach{i}{\findex{n}}\norm{\p{\vv{i}}}{i}\leq 1\]}}\norm{\func{\ml{}}{\vv{}}}{}\leq K$.
\item[${\propos{4}}.$]\quad
$\Exists{K}{\Rp}\Foreach{\vv{}}{\Times{\VV{1}}{\VV{n}}}\norm{\func{\ml{}}{\vv{}}}{}\leq K{\norm{\vv{1}}{1}}\ldots{\norm{\vv{n}}{n}}$.
\end{itemize}
\proof\\
$\[{\propos{1}\rightarrow\propos{2}}\]$
According to \refcor{cormulticontinuityequivs}, it is obvious.\\
$\[{\propos{2}\rightarrow\propos{3}}\]$
It is assumed that $\ml{}$ is continuous at $\mtuple{\zerov{\VVS{1}}}{\zerov{\VVS{n}}}$.
So, according to \refcor{cormulticontinuityequivs},
\begin{equation}\label{thmmultilinearcontinuousmapequivsp2p3eq1}
\Foreach{\varepsilon}{\Rp}
\Exists{\delta}{\Rp}\func{\image{\ml{}}}{\Times{\func{\ball{\NVS{1}}}{\binary{\zerov{\VVS{1}}}{\delta}}}
{\func{\ball{\NVS{n}}}{\binary{\zerov{\VVS{n}}}{\delta}}}}\subseteq
\func{\ball{\NWS{}}}{\binary{\func{\ml{}}{\suc{\zerov{\VVS{1}}}{\zerov{\VVS{n}}}}}{\varepsilon}}.
\end{equation}
Thus, according to \refcor{cornormedspaceball}, and considering the fact that
$\func{\ml{}}{\suc{\zerov{\VVS{1}}}{\zerov{\VVS{n}}}}=\zerov{\WS{}}$, it is clear that,
\begin{align}\label{thmmultilinearcontinuousmapequivsp2p3eq2}
&\Foreach{\varepsilon}{\Rp}\cr
&\Exists{\delta}{\Rp}\func{\image{\ml{}}}{\defset{\p{\vv{}}}{\Times{\VV{1}}{\VV{n}}}{\[\Foreach{i}{\findex{n}}\norm{\p{\vv{i}}}{i}\leq\delta\]}}\subseteq
\defset{\ww{}}{\W{}}{\norm{\ww{}}{}\leq\varepsilon}.
\end{align}
So,
\begin{equation}\label{thmmultilinearcontinuousmapequivsp2p3eq3}
\Existsis{\delta_1}{\Rp}\func{\image{\ml{}}}{\defset{\p{\vv{}}}{\Times{\VV{1}}{\VV{n}}}{\[\Foreach{i}{\findex{n}}\norm{\p{\vv{i}}}{i}\leq\delta_1\]}}\subseteq
\defset{\ww{}}{\W}{\norm{\ww{}}{}\leq 1},
\end{equation}
and equivalently,
\begin{equation}\label{thmmultilinearcontinuousmapequivsp2p3eq4}
\Foreach{\vv{}}{\defset{\p{\vv{}}}{\Times{\VV{1}}{\VV{n}}}{\[\Foreach{i}{\findex{n}}\norm{\p{\vv{i}}}{i}\leq\delta_1\]}}\norm{\func{\ml{}}{\vv{}}}{}\leq 1.
\end{equation}
On the other hand, according to \refdef{defnorm}, it is clear that,
\begin{align}\label{thmmultilinearcontinuousmapequivsp2p3eq5}
&\Foreach{\v}{\defset{\p{\vv{}}}{\Times{\VV{1}}{\VV{n}}}{\[\Foreach{i}{\findex{n}}\norm{\p{\vv{i}}}{i}\leq 1\]}}\cr
&\delta_1\v\in\defset{\p{\vv{}}}{\Times{\VV{1}}{\VV{n}}}{\[\Foreach{i}{\findex{n}}\norm{\p{\vv{i}}}{i}\leq\delta_1\]},
\end{align}
and hence, according to \Ref{thmmultilinearcontinuousmapequivsp2p3eq4},
\begin{align}\label{thmmultilinearcontinuousmapequivsp2p3eq6}
\Foreach{\v}{\defset{\p{\vv{}}}{\Times{\VV{1}}{\VV{n}}}{\[\Foreach{i}{\findex{n}}\norm{\p{\vv{i}}}{i}\leq 1\]}}
\norm{\func{\ml{}}{\delta_1\vv{}}}{}\leq 1.
\end{align}
On the other hand, according to \refdef{defnorm} and multilinearity of $\ml{}$,
\begin{equation}\label{thmmultilinearcontinuousmapequivsp2p3eq7}
\Foreach{\v}{\defset{\p{\vv{}}}{\Times{\VV{1}}{\VV{n}}}{\[\Foreach{i}{\findex{n}}\norm{\p{\vv{i}}}{i}\leq 1\]}}
\norm{\func{\ml{}}{\delta_1\vv{}}}{}=\delta_1^{n}\norm{\func{\ml{}}{\vv{}}}{}.
\end{equation}
\Ref{thmmultilinearcontinuousmapequivsp2p3eq6} and \Ref{thmmultilinearcontinuousmapequivsp2p3eq7}
imply,
\begin{equation}
\Existsis{\delta_1^{-n}}{\Rp}
\Foreach{\v}{\defset{\p{\vv{}}}{\Times{\VV{1}}{\VV{n}}}{\[\Foreach{i}{\findex{n}}\norm{\p{\vv{i}}}{i}\leq 1\]}}
\norm{\func{\ml{}}{\vv{}}}{}\leq\delta_1^{-n}.
\end{equation}
%%%%%%%%%%%%%%%%%%%%
$\[{\propos{3}\rightarrow\propos{4}}\]$
It is assumed that,
\begin{equation}
\Exists{K}{\Rp}\Foreach{\vv{}}{\defset{\p{\vv{}}}{\Times{\VV{1}}{\VV{n}}}{\[\Foreach{i}{\findex{n}}
\norm{\p{\vv{i}}}{i}\leq 1\]}}\norm{\func{\ml{}}{\vv{}}}{}\leq K.
\end{equation}
Then, since
\begin{equation}
\Foreach{i}{\findex{n}}
\Foreach{\vv{}}{\compl{\VV{i}}{\seta{\zerov{\VVS{i}}}}}
\norm{\frac{1}{\norm{\vv{}}{i}}\vv{}}{i}=1,
\end{equation}
it is clear that,
\begin{equation}
\Foreach{\vv{}}{\compl{\Times{\VV{1}}{\VV{n}}}{\seta{\mtuple{\zerov{\VVS{1}}}{\zerov{\VVS{n}}}}}}
\norm{\func{\ml{}}{\suc{\frac{1}{\norm{\vv{1}}{1}}\vv{1}}{\frac{1}{\norm{\vv{n}}{n}}\vv{n}}}}{}\leq K,
\end{equation}
and hence, according to \refdef{defnorm}, and multilinearity of $\ml{}$,
\begin{equation}
\Foreach{\vv{}}{\compl{\Times{\VV{1}}{\VV{n}}}{\seta{\mtuple{\zerov{\VVS{1}}}{\zerov{\VVS{n}}}}}}
\norm{\func{\ml{}}{\vv{}}}{}\leq K{\norm{\vv{1}}{1}}\ldots{\norm{\vv{n}}{n}}.
\end{equation}
Additionally, since $\norm{\zerov{\WS{}}}{}=
\norm{\func{\ml{}}{\suc{\zerov{\VVS{1}}}{\zerov{\VVS{n}}}}}{2}=\norm{\zerov{\VVS{1}}}{1}=
\ldots=\norm{\zerov{\VVS{n}}}{n}=0$,
it is obvious that,
\begin{equation}
\norm{\func{\ml{}}{\suc{\zerov{\VVS{1}}}{\zerov{\VVS{n}}}}}{2}\leq
K{\norm{\zerov{\VVS{1}}}{1}}\ldots{\norm{\zerov{\VVS{n}}}{n}}.
\end{equation}
So,
\begin{equation}
\Existsis{K}{\Rp}
\Foreach{\vv{}}{\Times{\VV{1}}{\VV{n}}}
\norm{\func{\ml{}}{\vv{}}}{}\leq K{\norm{\vv{1}}{1}}\ldots{\norm{\vv{n}}{n}}.
\end{equation}
%%%%%%%%%%%%%%%%%%%%%%%%%%%%%%%%%%%
$\[{\propos{4}\rightarrow\propos{1}}\]$
It is assumed that,
\begin{equation}
\Existsis{K}{\Rp}
\Foreach{\vv{}}{\Times{\VV{1}}{\VV{n}}}
\norm{\func{\ml{}}{\vv{}}}{}\leq K{\norm{\vv{1}}{1}}\ldots{\norm{\vv{n}}{n}}.
\end{equation}
\begin{itemize}
\item[${\textbf{\textsf{p1}}}$]
$\pnt{}$ is taken as an arbitrary element of $\Times{\VV{1}}{\VV{n}}$,
and $\varepsilon$ an arbitrary element of $\Rp$.
$\xi$ is defined as,
\begin{equation}
\xi:=\(\sum_{i=1}^{n}\norm{\pnt{i}}{i}\)+1,
\end{equation}
and $\delta$ as,
\begin{equation}
\delta:=\min{\seta{\binary{1}{\frac{\varepsilon}{nK\xi^{n-1}}}}}.
\end{equation}
So, according to \refdef{defnorm},
\begin{align}
&\Foreach{\vv{}}{\defset{\p{\vv{}}}{\Times{\VV{1}}{\VV{n}}}{\[\Foreach{i}{\findex{n}}\norm{\p{\vv{i}}-\pnt{i}}{i}<\delta\]}}\cr
&\begin{aligned}
&\quad\norm{\func{\ml{}}{\vv{}}-\func{\ml{}}{\pnt{}}}{}\cr
&=\norm{\func{\ml{}}{\suc{\vv{1}}{\vv{n}}}-\func{\ml{}}{\suc{\pnt{1}}{\pnt{n}}}}{}\cr
&=\norm{\func{\ml{}}{\binary{\vv{1}-\pnt{1}}{\suc{\vv{2}}{\vv{n}}}}+
\func{\ml{}}{\binary{\pnt{1}}{\suc{\vv{2}-\pnt{2}}{\vv{n}}}}+\ldots+
\func{\ml{}}{\binary{\suc{\pnt{1}}{\pnt{n_1}}}{\vv{n}-\pnt{n}}}}{}\cr
&\leq\quad\norm{\func{\ml{}}{\binary{\vv{1}-\pnt{1}}{\suc{\vv{2}}{\vv{n}}}}}{}+
\norm{\func{\ml{}}{\binary{\pnt{1}}{\binary{\vv{2}-\pnt{2}}{\suc{\vv{3}}{\vv{n}}}}}}{}+\ldots\cr
&\quad+
\norm{\func{\ml{}}{\binary{\suc{\pnt{1}}{\pnt{n-1}}}{\vv{n}-\pnt{n}}}}{}\cr
&\leq\quad K\norm{\vv{1}-\pnt{1}}{1}{\norm{\vv{2}}{2}}\ldots{\norm{\vv{n}}{n}}+
K\norm{\pnt{1}}{1}{\norm{\vv{2}-\pnt{2}}{2}}{\norm{\vv{3}}{3}}\ldots{\norm{\vv{n}}{n}}+\ldots\cr
&\quad+
K{\norm{\pnt{1}}{1}}\ldots{\norm{\pnt{n-1}}{n-1}}\norm{\vv{n}-\pnt{n}}{n}\cr
&<\(nK\xi^{n-1}\)\delta\leq\varepsilon.
\end{aligned}
\end{align}
This means,
\begin{align}
\func{\image{\ml{}}}{\Times{\func{\ball{\NVS{1}}}{\binary{\pnt{1}}{\delta}}}
{\func{\ball{\NVS{n}}}{\binary{\pnt{n}}{\delta}}}}\subseteq
\func{\ball{\NWS{}}}{\binary{\func{\ml{}}{\pnt{}}}{\varepsilon}}.
\end{align}
\endp
\end{itemize}
So,
\begin{align}
&\Foreach{\pnt{}}{\Times{\VV{1}}{\VV{n}}}\cr
&\Foreach{\varepsilon}{\Rp}
\Exists{\delta}{\Rp}\func{\image{\cf}}{\Times{\func{\ball{\NVS{1}}}{\binary{\pnt{1}}{\delta}}}
{\func{\ball{{\NVS{n}}}}{\binary{\pnt{n}}{\delta}}}}\subseteq
\func{\ball{\NWS{}}}{\binary{\func{\cf}{\point}}{\varepsilon}},
\end{align}
which according to \refcor{cormulticontinuityequivs}, means,
\begin{equation}
\ml{}\in\MLC{\suc{\NVS{1}}{\NVS{n}}}{\NWS{}}.
\end{equation}
\endthm
%%%%%%%%%%%%%%%%%%%%%%%%%%%%%%%%%%%%%%%%%%%%%%%%%%%%%%%%%%%%%%%%%%%%%%%%%%%%%%%%%%%%%%%%%%%%%%%%%%%%%%%%%%
\theorem\label{thmmultilinearcontinuousmapsspace}
$\sindex{\VV{i}}{i}{\findex{n}}$
is taken to be an indexed-collection of non-empty sets, for some $n$ in $\N$.
For each $i$ in $\seta{\suc{1}{n}}$, $\VVS{i}=\tuple{\VV{i}}{\vsum{i}}{\spro{i}}{\F}$
is taken as an element of $\vecspaces{\VV{i}}{\F}$, and $\anorm{i}$
an element of $\norms{\VVS{i}}$.
$\WW{}$ is taken as a non-empty set.
$\WS{}=\tuple{\WW{}}{\vsum{w}}{\spro{w}}{\F}$ is taken as
an element of $\vecspaces{\WW{}}{\F}$, and $\anorm{}$ an element of $\norms{\WS{}}$.
Also, for each $i$ in $\findex{n}$, $\NVS{i}:=\opair{\VVS{i}}{\anorm{i}}$, and
$\NWS{}:=\opair{\WS{}}{\anorm{}}$.
\begin{align}
\tuple{\MLC{\suc{\NVS{1}}{\NVS{n}}}{\NWS{}}}{\vsum{mlc}}{\spro{mlc}}{\F}\in
\vecspaces{\MLC{\suc{\NVS{1}}{\NVS{n}}}{\NWS{}}}{\F},
\end{align}
where,
\begin{align}
\vsum{mlc}&:=\func{\resd{\vsum{ml}}}{\Cprod{\MLC{\suc{\NVS{1}}{\NVS{n}}}{\NWS{}}}{\MLC{\suc{\NVS{1}}{\NVS{n}}}{\NWS{}}}},\\
\spro{mlc}&:=\func{\resd{\spro{ml}}}{\Cprod{\f}{\MLC{\suc{\NVS{1}}{\NVS{n}}}{\NWS{}}}},
\end{align}
$\vsum{ml}$ and $\spro{ml}$ being the addition and scalar-multiplication operations of
the vector-space of all multilinear maps from $\suc{\VVS{1}}{\VVS{n}}$ to the vector-space $\WS{}$.
\proof
$\MLC{\suc{\NVS{1}}{\NVS{n}}}{\NWS{}}$ is a non-empty subset of
$\VFunc{\Times{\VV{1}}{\VV{n}}}{\WS{}}$, since the constant map that maps every element of
$\Times{\VV{1}}{\VV{n}}$ to $\zerov{\WS{}}$ is a multilinear-continuous map from $\mtuple{\NVS{1}}{\NVS{n}}$
to $\NWS{}$.
\begin{itemize}
\item[${\textbf{\textsf{p1}}}$]
Each $\ml{1}$ and $\ml{2}$ is taken as an element of
$\MLC{\suc{\NVS{1}}{\NVS{n}}}{\NWS{}}$. So, according to
\refdef{defmultilinearcontinuousmap}, $\ml{1}$ and $\ml{2}$ belong to $\MLin{\suc{\VVS{1}}{\VVS{n}}}{\WS{}}$,
and hence according to \refthm{thmmultilinearmapsspace},
\begin{equation}
\ml{1}\vsum{ml}\ml{2}\in\MLin{\suc{\VVS{1}}{\VS{n}}}{\WS{}}.
\end{equation}
On the other hand, according to \refthm{thmmultilinearcontinuousmapequivs},
\begin{align}
&\Exists{K}{\Rp}\Foreach{\vv{}}{\Times{\VV{1}}{\VV{n}}}\norm{\func{\ml{1}}{\vv{}}}{}\leq K{\norm{\vv{1}}{1}}\ldots{\norm{\vv{n}}{n}},\\
&\Exists{\p{K}}{\Rp}\Foreach{\vv{}}{\Times{\VV{1}}{\VV{n}}}\norm{\func{\ml{2}}{\vv{}}}{}\leq\p{K}{\norm{\vv{1}}{1}}\ldots{\norm{\vv{n}}{n}}.
\end{align}
Therefore, according to \refdef{defnorm},
\begin{align}
&\Existsis{K+\p{K}}{\Rp}
\Foreach{\vv{}}{\Times{\VV{1}}{\VV{n}}}\cr
&\begin{aligned}
\norm{\func{\(\ml{1}\vsum{ml}\ml{2}\)}{\v}}{}&=
\norm{\func{\ml{1}}{\v}\vsum{w}\func{\ml{2}}{\v}}{}\cr
&\leq\norm{\func{\ml{1}}{\v}}{}+\norm{\func{\ml{2}}{\v}}{}\cr
&\leq K{\norm{\vv{1}}{1}}\ldots{\norm{\vv{n}}{n}}+
\p{K}{\norm{\vv{1}}{1}}\ldots{\norm{\vv{n}}{n}}\cr
&=\(K+\p{K}\){\norm{\vv{1}}{1}}\ldots{\norm{\vv{n}}{n}}.\cr
\end{aligned}\cr
&{}
\end{align}
So, according to \refthm{thmmultilinearcontinuousmapequivs},
\begin{equation}
\(\ml{1}\vsum{ml}\ml{2}\)\in\MLC{\suc{\NVS{1}}{\NVS{n}}}{\NWS{}}.
\end{equation}
\endp
\end{itemize}
\begin{itemize}
\item[${\textbf{\textsf{p2}}}$]
$\ml{}$ is taken as an element of
$\MLC{\suc{\NVS{1}}{\NVS{n}}}{\NWS{}}$, and $\x$ an element of $\f$.
So, according to \refdef{defmultilinearcontinuousmap}, $\ml{}$ belongs to $\MLin{\suc{\VVS{1}}{\VVS{n}}}{\WS{}}$,
and hence according to \refthm{thmmultilinearmapsspace},
\begin{equation}
\(\x\spro{ml}\ml{}\)\in\MLin{\suc{\VVS{1}}{\VVS{n}}}{\NWS{}}.
\end{equation}
On the other hand, according to \refthm{thmmultilinearcontinuousmapequivs},
\begin{equation}
\Exists{K}{\Rp}\Foreach{\vv{}}{\Times{\VV{1}}{\VV{n}}}\norm{\func{\ml{}}{\vv{}}}{}\leq K{\norm{\vv{1}}{1}}\ldots{\norm{\vv{n}}{n}},
\end{equation}
and thus, according to \refdef{defnorm},
\begin{align}
\Existsis{\(\abs{\x}K+1\)}{\Rp}
\Foreach{\vv{}}{\Times{\VV{1}}{\VV{n}}}
\norm{\func{\(\x\spro{ml}\ml{}\)}{\v}}{}&=
\norm{\x\spro{2}\(\func{\ml{}}{\v}\)}{}\cr
&=\abs{\x}\norm{\func{\ml{}}{\v}}{}\cr
&\leq\(\abs{\x}K\)\norm{\vv{1}}{1}\ldots\norm{\vv{n}}{n}\cr
&\leq\(\abs{\x}K+1\)\norm{\vv{1}}{1}\ldots\norm{\vv{n}}{n}.\cr
&{}
\end{align}
So, according to \refthm{thmmultilinearcontinuousmapequivs},
\begin{equation}
\(\x\spro{ml}\ml{}\)\in\MLC{\suc{\NVS{1}}{\NVS{n}}}{\NWS{}}.
\end{equation}
\endp
\end{itemize}
\endthm
%%%%%%%%%%%%%%%%%%%%%%%%%%%%%%%%%%%%%%%%%%%%%%%%%%%%%%%%%%%%%%%%%%%%%%%%%%%%%%%%%%%%%%%%%%%%%%%%%%%%%%%%%%
\definition\label{defmultilinearcontinuousmapsnorm}
$\sindex{\VV{i}}{i}{\findex{n}}$
is taken to be an indexed-collection of non-empty sets, for some $n$ in $\N$.
For each $i$ in $\seta{\suc{1}{n}}$, $\VVS{i}=\tuple{\VV{i}}{\vsum{i}}{\spro{i}}{\F}$
is taken as an element of $\vecspaces{\VV{i}}{\F}$, and $\anorm{i}$
an element of $\norms{\VVS{i}}$.
$\WW{}$ is taken as a non-empty set.
$\WS{}=\tuple{\WW{}}{\vsum{w}}{\spro{w}}{\F}$ is taken as
an element of $\vecspaces{\WW{}}{\F}$, and $\anorm{}$ an element of $\norms{\WS{}}$.
Also, for each $i$ in $\findex{n}$, $\NVS{i}:=\opair{\VVS{i}}{\anorm{i}}$, and
$\NWS{}:=\opair{\WS{}}{\anorm{}}$.
\begin{itemize}
\item
\begin{equation}
\VMLC{\suc{\NVS{1}}{\NVS{n}}}{\NWS{}}:=
\tuple{\MLC{\suc{\NVS{1}}{\NVS{n}}}{\NWS{}}}{\vsum{mlc}}{\spro{mlc}}{\F}.
\end{equation}
$\VMLC{\suc{\NVS{1}}{\NVS{n}}}{\NWS{}}$ is called the
$\quotl$vector-space of multilinear-continuous maps from\\
$\mtuple{\NVS{1}}{\NVS{n}}$ to $\NWS{}$$\quotr$.
\item
\begin{align}
&\VMLCnorm{\mtuple{\NVS{1}}{\NVS{2}}}{\NWS{}}\indef\Func{\MLC{\suc{\NVS{1}}{\NVS{n}}}{\NWS{}}}{\Rpz},\cr\cr
&\Foreach{\ml{}}{\MLC{\suc{\NVS{1}}{\NVS{n}}}{\NWS{}}}\cr
&\VMLCnormf{\ml{}}{\mtuple{\NVS{1}}{\NVS{n}}}{\NWS{}}\eqdef
\sup\defSet{\norm{\func{\ml{}}{\v}}{}}
{\v\in\Times{\func{\cball{\NVS{1}}}{\binary{\zerov{\VVS{1}}}{1}}}{\func{\cball{\NVS{n}}}{\binary{\zerov{\VVS{n}}}{n}}}}.\cr
&{}
\end{align}
\end{itemize}
\endef
%%%%%%%%%%%%%%%%%%%%%%%%%%%%%%%%%%%%%%%%%%%%%%%%%%%%%%%%%%%%%%%%%%%%%%%%%%%%%%%%%%%%%%%%%%%%%%%%%%%%%%%%%%
\theorem\label{thmmultilinearcontinuousmapsnorm}
$\sindex{\VV{i}}{i}{\findex{n}}$
is taken to be an indexed-collection of non-empty sets, for some $n$ in $\N$.
For each $i$ in $\seta{\suc{1}{n}}$, $\VVS{i}=\tuple{\VV{i}}{\vsum{i}}{\spro{i}}{\F}$
is taken as an element of $\vecspaces{\VV{i}}{\F}$, and $\anorm{i}$
an element of $\norms{\VVS{i}}$.
$\WW{}$ is taken as a non-empty set.
$\WS{}=\tuple{\WW{}}{\vsum{w}}{\spro{w}}{\F}$ is taken as
an element of $\vecspaces{\WW{}}{\F}$, and $\anorm{}$ an element of $\norms{\WS{}}$.
Also, for each $i$ in $\findex{n}$, $\NVS{i}:=\opair{\VVS{i}}{\anorm{i}}$, and
$\NWS{}:=\opair{\WS{}}{\anorm{}}$.
Also, $\VMLC{\suc{\NVS{1}}{\NVS{n}}}{\NWS{}}=\tuple{\MLC{\suc{\NVS{1}}{\NVS{n}}}{\NWS{}}}{\vsum{mlc}}{\spro{mlc}}{\F}$.
$\VMLCnorm{\mtuple{\NVS{1}}{\NVS{2}}}{\NWS{}}$ is a norm on $\VMLC{\suc{\NVS{1}}{\NVS{n}}}{\NWS{}}$.
\begin{equation}
\VMLCnorm{\mtuple{\NVS{1}}{\NVS{2}}}{\NWS{}}\in\norms{\VMLC{\suc{\NVS{1}}{\NVS{n}}}{\NWS{}}}.
\end{equation}
\proof
\begin{itemize}
\item[${\textbf{\textsf{p1}}}$]
$\ml{}$ is taken as an arbitrary element of $\MLC{\suc{\NVS{1}}{\NVS{n}}}{\NWS{}}$, and it is assumed that
\begin{equation}\label{thmmultilinearcontinuousmapsnormp1eq1}
\VMLCnormf{\ml{}}{\mtuple{\NVS{1}}{\NVS{n}}}{\NWS{}}=0.
\end{equation}
Then according to \refdef{defmultilinearcontinuousmapsnorm},
\begin{equation}\label{thmmultilinearcontinuousmapsnormp1eq2}
\sup\defSet{\norm{\func{\ml{}}{\v}}{}}
{\v\in\Times{\func{\cball{\NVS{1}}}{\binary{\zerov{\VVS{1}}}{1}}}{\func{\cball{\NVS{n}}}{\binary{\zerov{\VVS{n}}}{n}}}}=0,
\end{equation}
which clearly implies that,
\begin{equation}\label{thmmultilinearcontinuousmapsnormp1eq3}
\Foreach{\vv{}}{\Times{\func{\cball{\NVS{1}}}{\binary{\zerov{\VVS{1}}}{1}}}{\func{\cball{\NVS{n}}}{\binary{\zerov{\VVS{n}}}{n}}}}
\norm{\func{\ml{}}{\vv{}}}{}=0.
\end{equation}
\begin{itemize}
\item[${\textbf{\textsf{p1-1}}}$]
$\vv{}$ is taken as an arbitrary element of
$\Times{\(\compl{\VV{1}}{\seta{\zerov{\VVS{1}}}}\)}{\(\compl{\VV{n}}{\seta{\zerov{\VVS{n}}}}\)}$.
According to multilinearity of $\ml{}$,
\begin{equation}\label{thmmultilinearcontinuousmapsnormp11eq1}
\func{\ml{}}{\vv{}}=
\(\norm{\vv{1}}{1}\ldots\norm{\vv{n}}{n}\)\spro{w}\func{\ml{}}{\suc{\frac{1}{\norm{\vv{1}}{1}}\spro{1}\vv{1}}{\frac{1}{\norm{\vv{n}}{n}}\spro{n}\vv{n}}}.
\end{equation}
Additionally, since,
\begin{equation}\label{thmmultilinearcontinuousmapsnormp11eq2}
\mtuple{\frac{1}{\norm{\vv{1}}{1}}\spro{1}\vv{1}}{\frac{1}{\norm{\vv{n}}{n}}\spro{n}\vv{n}}\in
\Times{\func{\cball{\NVS{1}}}{\binary{\zerov{\VVS{1}}}{1}}}{\func{\cball{\NVS{n}}}{\binary{\zerov{\VVS{n}}}{n}}},
\end{equation}
according to \Ref{thmmultilinearcontinuousmapsnormp1eq3},
\begin{equation}\label{thmmultilinearcontinuousmapsnormp11eq3}
\norm{\func{\ml{}}{\suc{\frac{1}{\norm{\vv{1}}{1}}\spro{1}\vv{1}}{\frac{1}{\norm{\vv{n}}{n}}\spro{n}\vv{n}}}}{}=0,
\end{equation}
and hence according to \Ref{thmmultilinearcontinuousmapsnormp11eq1}, and \refdef{defnorm},
\begin{align}
\norm{\func{\ml{}}{\vv{}}}{}=0.
\end{align}
\endp
\end{itemize}
So,
\begin{equation}
\Foreach{\vv{}}{\Times{\(\compl{\VV{1}}{\seta{\zerov{\VVS{1}}}}\)}{\(\compl{\VV{n}}{\seta{\zerov{\VVS{n}}}}\)}}
\norm{\func{\ml{}}{\vv{}}}{}=0.
\end{equation}
Additionally, it is known that,
\begin{align}
\Foreach{\vv{}}{\[\compl{\Times{\VV{1}}{\VV{n}}}
{\bigg(\Times{\(\compl{\VV{1}}{\seta{\zerov{\VVS{1}}}}\)}{\(\compl{\VV{n}}{\seta{\zerov{\VVS{n}}}}\)}\bigg)}\]}
\norm{\func{\ml{}}{\vv{}}}{}=0.
\end{align}
Therefore,
\begin{equation}
\Foreach{\vv{}}{\Times{\VV{1}}{\VV{n}}}
\norm{\func{\ml{}}{\vv{}}}{}=0,
\end{equation}
and according to \refdef{defnorm},
\begin{equation}
\Foreach{\vv{}}{\Times{\VV{1}}{\VV{n}}}
\func{\ml{}}{\vv{}}=\zerov{\WS{}}.
\end{equation}
This means that $\ml{}$ is the zero element of the vector-space $\VMLC{\suc{\NVS{1}}{\NVS{n}}}{\NWS{}}$.
\begin{equation}
\ml{}=\zerov{\VFunc{\suc{\VV{1}}{\VV{n}}}{\WS{}}}.
\end{equation}
\endp
\end{itemize}
So,
\begin{equation}\label{thmmultilinearcontinuousmapsnormpeq1}
\Foreach{\ml{}}{\MLC{\suc{\NVS{1}}{\NVS{n}}}{\NWS{}}}
\bigg[\VMLCnormf{\ml{}}{\mtuple{\NVS{1}}{\NVS{n}}}{\NWS{}}=0\bigg]\then
\bigg[\ml{}=\zerov{\VFunc{\suc{\VV{1}}{\VV{n}}}{\WS{}}}\bigg].
\end{equation}
Additionally, it is obvious that,
\begin{equation}\label{thmmultilinearcontinuousmapsnormpeq2}
\VMLCnormf{\zerov{\VFunc{\suc{\VV{1}}{\VV{n}}}{\WS{}}}}{\mtuple{\NVS{1}}{\NVS{n}}}{\NWS{}}=0.
\end{equation}
\begin{itemize}
\item[${\textbf{\textsf{p2}}}$]
$\ml{}$ is taken as an arbitrary element of $\MLC{\suc{\NVS{1}}{\NVS{n}}}{\NWS{}}$, and $\x$ an element of
$\f$. According to \refdef{defmultilinearcontinuousmapsnorm}, and \refdef{defnorm},
\begin{align}
\VMLCnormf{\(\x\spro{mlc}\ml{}\)}{\mtuple{\NVS{1}}{\NVS{n}}}{\NWS{}}&=
\sup\defSet{\norm{\func{\(\x\spro{mlc}\ml{}\)}{\v}}{}}
{\v\in\Times{\func{\cball{\NVS{1}}}{\binary{\zerov{\VVS{1}}}{1}}}{\func{\cball{\NVS{n}}}{\binary{\zerov{\VVS{n}}}{n}}}}\cr
&=\sup\defSet{\abs{\x}\norm{\func{\ml{}}{\v}}{}}
{\v\in\Times{\func{\cball{\NVS{1}}}{\binary{\zerov{\VVS{1}}}{1}}}{\func{\cball{\NVS{n}}}{\binary{\zerov{\VVS{n}}}{n}}}}\cr
&=\abs{\x}\sup\defSet{\norm{\func{\ml{}}{\v}}{}}
{\v\in\Times{\func{\cball{\NVS{1}}}{\binary{\zerov{\VVS{1}}}{1}}}{\func{\cball{\NVS{n}}}{\binary{\zerov{\VVS{n}}}{n}}}}\cr
&=\abs{\x}\VMLCnormf{\ml{}}{\mtuple{\NVS{1}}{\NVS{n}}}{\NWS{}}.
\end{align}
\endp
\end{itemize}
\begin{itemize}
\item[${\textbf{\textsf{p3}}}$]
Each $\ml{1}$ and $\ml{2}$ is taken as an arbitrary element of $\MLC{\suc{\NVS{1}}{\NVS{n}}}{\NWS{}}$.
According to \refdef{defnorm},
\begin{align}
&\quad\VMLCnormf{\ml{1}\vsum{mlc}\ml{2}}{\mtuple{\NVS{1}}{\NVS{n}}}{\NWS{}}\cr
&=\sup\defSet{\norm{\func{\(\ml{1}\vsum{mlc}\ml{2}\)}{\v}}{}}
{\v\in\Times{\func{\cball{\NVS{1}}}{\binary{\zerov{\VVS{1}}}{1}}}{\func{\cball{\NVS{n}}}{\binary{\zerov{\VVS{n}}}{n}}}}\cr
&=\sup\defSet{\norm{\func{\ml{1}}{\vv{}}\vsum{w}\func{\ml{2}}{\vv{}}}{}}
{\v\in\Times{\func{\cball{\NVS{1}}}{\binary{\zerov{\VVS{1}}}{1}}}{\func{\cball{\NVS{n}}}{\binary{\zerov{\VVS{n}}}{n}}}}\cr
&\leq\sup\defSet{\norm{\func{\ml{1}}{\vv{}}}{}\vsum{}\norm{\func{\ml{2}}{\vv{}}}{}}
{\v\in\Times{\func{\cball{\NVS{1}}}{\binary{\zerov{\VVS{1}}}{1}}}{\func{\cball{\NVS{n}}}{\binary{\zerov{\VVS{n}}}{n}}}}\cr
&\leq\quad\sup\defSet{\norm{\func{\ml{1}}{\v}}{}}
{\v\in\Times{\func{\cball{\NVS{1}}}{\binary{\zerov{\VVS{1}}}{1}}}{\func{\cball{\NVS{n}}}{\binary{\zerov{\VVS{n}}}{n}}}}\cr
&\quad+\sup\defSet{\norm{\func{\ml{2}}{\v}}{}}
{\v\in\Times{\func{\cball{\NVS{1}}}{\binary{\zerov{\VVS{1}}}{1}}}{\func{\cball{\NVS{n}}}{\binary{\zerov{\VVS{n}}}{n}}}}\cr
&=\VMLCnormf{\ml{1}}{\mtuple{\NVS{1}}{\NVS{n}}}{\NWS{}}+\VMLCnormf{\ml{2}}{\mtuple{\NVS{1}}{\NVS{n}}}{\NWS{}}.
\end{align}
\endp
\end{itemize}
So, according to \refdef{defnorm},
$\VMLCnorm{\mtuple{\NVS{1}}{\NVS{2}}}{\NWS{}}$ is a norm on $\VMLC{\suc{\NVS{1}}{\NVS{n}}}{\NWS{}}$.
\endthm
%%%%%%%%%%%%%%%%%%%%%%%%%%%%%%%%%%%%%%%%%%%%%%%%%%%%%%%%%%%%%%%%%%%%%%%%%%%%%%%%%%%%%%%%%%%%%%%%%%%%%%%%%%
\definition\label{defnormedspaceofmultilinearcontinuousmaps}
$\sindex{\VV{i}}{i}{\findex{n}}$
is taken to be an indexed-collection of non-empty sets, for some $n$ in $\N$.
For each $i$ in $\seta{\suc{1}{n}}$, $\VVS{i}=\tuple{\VV{i}}{\vsum{i}}{\spro{i}}{\F}$
is taken as an element of $\vecspaces{\VV{i}}{\F}$, and $\anorm{i}$
an element of $\norms{\VVS{i}}$.
$\WW{}$ is taken as a non-empty set.
$\WS{}=\tuple{\WW{}}{\vsum{w}}{\spro{w}}{\F}$ is taken as
an element of $\vecspaces{\WW{}}{\F}$, and $\anorm{}$ an element of $\norms{\WS{}}$.
Also, for each $i$ in $\findex{n}$, $\NVS{i}:=\opair{\VVS{i}}{\anorm{i}}$, and
$\NWS{}:=\opair{\WS{}}{\anorm{}}$.
\begin{equation}
\NVMLC{\suc{\NVS{1}}{\NVS{n}}}{\NWS{}}:=
\opair{\VMLC{\suc{\NVS{1}}{\NVS{n}}}{\NWS{}}}{\VMLCnorm{\mtuple{\NVS{1}}{\NVS{2}}}{\NWS{}}}.
\end{equation}
$\NVMLC{\suc{\NVS{1}}{\NVS{n}}}{\NWS{}}$ is called the
$\quotl$normed vector-space of multilinear-continuous maps from\\
$\mtuple{\NVS{1}}{\NVS{n}}$ to $\NWS{}$$\quotr$.
\endef
%%%%%%%%%%%%%%%%%%%%%%%%%%%%%%%%%%%%%%%%%%%%%%%%%%%%%%%%%%%%%%%%%%%%%%%%%%%%%%%%%%%%%%%%%%%%%%%%%%%%%%%%%%
\theorem\label{thmmultilinearcontinuousmapsnormproperty1}
$\sindex{\VV{i}}{i}{\findex{n}}$
is taken to be an indexed-collection of non-empty sets, for some $n$ in $\N$.
For each $i$ in $\seta{\suc{1}{n}}$, $\VVS{i}=\tuple{\VV{i}}{\vsum{i}}{\spro{i}}{\F}$
is taken as an element of $\vecspaces{\VV{i}}{\F}$, and $\anorm{i}$
an element of $\norms{\VVS{i}}$.
$\WW{}$ is taken as a non-empty set.
$\WS{}=\tuple{\WW{}}{\vsum{w}}{\spro{w}}{\F}$ is taken as
an element of $\vecspaces{\WW{}}{\F}$, and $\anorm{}$ an element of $\norms{\WS{}}$.
Also, for each $i$ in $\findex{n}$, $\NVS{i}:=\opair{\VVS{i}}{\anorm{i}}$, and
$\NWS{}:=\opair{\WS{}}{\anorm{}}$.
\begin{align}
&\Foreach{\ml{}}{\MLC{\suc{\NVS{1}}{\NVS{n}}}{\NWS{}}}\cr
&\Foreach{\vv{}}{\Times{\VV{1}}{\VV{n}}}
\norm{\func{\ml{}}{\vv{}}}{}\leq
\VMLCnormf{\ml{}}{\mtuple{\NVS{1}}{\NVS{n}}}{\NWS{}}
\norm{\vv{1}}{1}\ldots\norm{\vv{n}}{n}.
\end{align}
\proof
\begin{itemize}
\item[${\textbf{\textsf{p1}}}$]
$\ml{}$ is taken as an arbitrary element of $\MLC{\suc{\NVS{1}}{\NVS{n}}}{\NWS{}}$,
and $\vv{}$ an arbitrary element of
$\Times{\(\compl{\VV{1}}{\seta{\zerov{\VVS{1}}}}\)}{\(\compl{\VV{n}}{\seta{\zerov{\VVS{n}}}}\)}$.
According to \refdef{defnorm}, and multilinearity of $\ml{}$,
\begin{align}
\norm{\func{\ml{}}{\v}}{}&=
\norm{\norm{\vv{1}}{1}\ldots\norm{\vv{n}}{n}\spro{w}\func{\ml{}}{\suc{\frac{1}{\norm{\vv{1}}{1}}\spro{1}\vv{1}}{\frac{1}{\norm{\vv{n}}{n}}\spro{n}\vv{n}}}}{}\cr
&=\norm{\vv{1}}{1}\ldots\norm{\vv{n}}{n}\norm{\func{\ml{}}{\suc{\frac{1}{\norm{\vv{1}}{1}}\spro{1}\vv{1}}{\frac{1}{\norm{\vv{n}}{n}}\spro{n}\vv{n}}}}{}\cr.
\end{align}
Additionally, since,
\begin{equation}\label{thmmultilinearcontinuousmapsnormp11eq2}
\mtuple{\frac{1}{\norm{\vv{1}}{1}}\spro{1}\vv{1}}{\frac{1}{\norm{\vv{n}}{n}}\spro{n}\vv{n}}\in
\Times{\func{\cball{\NVS{1}}}{\binary{\zerov{\VVS{1}}}{1}}}{\func{\cball{\NVS{n}}}{\binary{\zerov{\VVS{n}}}{n}}},
\end{equation}
according to \refdef{defmultilinearcontinuousmapsnorm},
\begin{equation}
\norm{\func{\ml{}}{\suc{\frac{1}{\norm{\vv{1}}{1}}\spro{1}\vv{1}}{\frac{1}{\norm{\vv{n}}{n}}\spro{n}\vv{n}}}}{}\leq
\VMLCnormf{\ml{}}{\mtuple{\NVS{1}}{\NVS{n}}}{\NWS{}}.
\end{equation}
Thus,
\begin{equation}
\norm{\func{\ml{}}{\v}}{}\leq
\VMLCnormf{\ml{}}{\mtuple{\NVS{1}}{\NVS{n}}}{\NWS{}}
\norm{\vv{1}}{1}\ldots\norm{\vv{n}}{n}.
\end{equation}
\endp
\end{itemize}
Additionally, it is obvious that,
\begin{align}
&\Foreach{\vv{}}{\[\compl{\Times{\VV{1}}{\VV{n}}}
{\bigg(\Times{\(\compl{\VV{1}}{\seta{\zerov{\VVS{1}}}}\)}{\(\compl{\VV{n}}{\seta{\zerov{\VVS{n}}}}\)}\bigg)}\]}\cr
&0=\norm{\func{\ml{}}{\vv{}}}{}
\leq\VMLCnormf{\ml{}}{\mtuple{\NVS{1}}{\NVS{n}}}{\NWS{}}
\norm{\vv{1}}{1}\ldots\norm{\vv{n}}{n}=0.
\end{align}
So,
\begin{equation}
\Foreach{\vv{}}{\Times{\VV{1}}{\VV{n}}}
\norm{\func{\ml{}}{\vv{}}}{}\leq
\VMLCnormf{\ml{}}{\mtuple{\NVS{1}}{\NVS{n}}}{\NWS{}}
\norm{\vv{1}}{1}\ldots\norm{\vv{n}}{n}.
\end{equation}
\endthm
%%%%%%%%%%%%%%%%%%%%%%%%%%%%%%%%%%%%%%%%%%%%%%%%%%%%%%%%%%%%%%%%%%%%%%%%%%%%%%%%%%%%%%%%%%%%%%%%%%%%%%%%%%
\theorem\label{thmmultilinearcontinuousmapsnormproperty2}
$\sindex{\VV{i}}{i}{\findex{n}}$
is taken to be an indexed-collection of non-empty sets, for some $n$ in $\N$.
For each $i$ in $\seta{\suc{1}{n}}$, $\VVS{i}=\tuple{\VV{i}}{\vsum{i}}{\spro{i}}{\F}$
is taken as an element of $\vecspaces{\VV{i}}{\F}$, and $\anorm{i}$
an element of $\norms{\VVS{i}}$.
$\WW{}$ is taken as a non-empty set.
$\WS{}=\tuple{\WW{}}{\vsum{w}}{\spro{w}}{\F}$ is taken as
an element of $\vecspaces{\WW{}}{\F}$, and $\anorm{}$ an element of $\norms{\WS{}}$.
Also, for each $i$ in $\findex{n}$, $\NVS{i}:=\opair{\VVS{i}}{\anorm{i}}$, and
$\NWS{}:=\opair{\WS{}}{\anorm{}}$.
\begin{align}
&\Foreach{\ml{}}{\MLC{\suc{\NVS{1}}{\NVS{n}}}{\NWS{}}}\cr
&\VMLCnormf{\ml{}}{\mtuple{\NVS{1}}{\NVS{n}}}{\NWS{}}=\min
\defset{M}{\Rpz}
{\bigg[\Foreach{\vv{}}{\Times{\VV{1}}{\VV{n}}}
\norm{\func{\ml{}}{\vv{}}}{}\leq
M\norm{\vv{1}}{1}\ldots\norm{\vv{n}}{n}\bigg]}.\cr
&{}
\end{align}
\proof
$\ml{}$ is taken as an element of $\MLC{\suc{\NVS{1}}{\NVS{n}}}{\NWS{}}$.
According to \refthm{thmmultilinearcontinuousmapsnorm}
and \refthm{thmmultilinearcontinuousmapsnormproperty1},
\begin{equation}\label{thmmultilinearcontinuousmapsnormproperty2peq1}
\VMLCnormf{\ml{}}{\mtuple{\NVS{1}}{\NVS{n}}}{\NWS{}}\in
\defset{M}{\Rpz}
{\bigg[\Foreach{\vv{}}{\Times{\VV{1}}{\VV{n}}}
\norm{\func{\ml{}}{\vv{}}}{}\leq
M\norm{\vv{1}}{1}\ldots\norm{\vv{n}}{n}\bigg]}.
\end{equation}
\begin{itemize}
\item[${\textbf{\textsf{p1}}}$]
$M$ is taken as an element of\\
$\defset{\p{M}}{\Rpz}{\bigg[\Foreach{\vv{}}{\Times{\VV{1}}{\VV{n}}}\norm{\func{\ml{}}{\vv{}}}{}\leq
\p{M}\norm{\vv{1}}{1}\ldots\norm{\vv{n}}{n}\bigg]}$.
Then it is clear that,
\begin{align}
\Foreach{\vv{}}
{\Times{\func{\cball{\NVS{1}}}{\binary{\zerov{\VVS{1}}}{1}}}{\func{\cball{\NVS{n}}}{\binary{\zerov{\VVS{n}}}{n}}}}
\norm{\func{\ml{}}{\vv{}}}{}&\leq
M\norm{\vv{1}}{1}\ldots\norm{\vv{n}}{n}\cr
&\leq M,
\end{align}
and hence according to \refdef{defmultilinearcontinuousmapsnorm},
\begin{equation}
\VMLCnormf{\ml{}}{\mtuple{\NVS{1}}{\NVS{n}}}{\NWS{}}\leq M.
\end{equation}
\endp
\end{itemize}
So,
\begin{align}\label{thmmultilinearcontinuousmapsnormproperty2peq2}
&\Foreach{M}{\defset{\p{M}}{\Rpz}{\bigg[\Foreach{\vv{}}{\Times{\VV{1}}{\VV{n}}}\norm{\func{\ml{}}{\vv{}}}{}\leq
\p{M}\norm{\vv{1}}{1}\ldots\norm{\vv{n}}{n}\bigg]}}\cr
&\VMLCnormf{\ml{}}{\mtuple{\NVS{1}}{\NVS{n}}}{\NWS{}}\leq M.
\end{align}
\Ref{thmmultilinearcontinuousmapsnormproperty2peq1} and \Ref{thmmultilinearcontinuousmapsnormproperty2peq2}
imply,
\begin{equation}
\VMLCnormf{\ml{}}{\mtuple{\NVS{1}}{\NVS{n}}}{\NWS{}}=\min
\defset{M}{\Rpz}
{\bigg[\Foreach{\vv{}}{\Times{\VV{1}}{\VV{n}}}
\norm{\func{\ml{}}{\vv{}}}{}\leq
M\norm{\vv{1}}{1}\ldots\norm{\vv{n}}{n}\bigg]}.
\end{equation}
\endthm
%%%%%%%%%%%%%%%%%%%%%%%%%%%%%%%%%%%%%%%%%%%%%%%%%%%%%%%%%%%%%%%%%%%%%%%%%%%%%%%%%%%%%%%%%%%%%%%%%%%%%%%%%%
\theorem\label{thmBanachMLC}
$\sindex{\VV{i}}{i}{\findex{n}}$
is taken to be an indexed-collection of non-empty sets, for some $n$ in $\N$.
For each $i$ in $\seta{\suc{1}{n}}$, $\VVS{i}=\tuple{\VV{i}}{\vsum{i}}{\spro{i}}{\F}$
is taken as an element of $\vecspaces{\VV{i}}{\F}$, and $\anorm{i}$
an element of $\norms{\VVS{i}}$.
$\WW{}$ is taken as a non-empty set.
$\WS{}=\tuple{\WW{}}{\vsum{w}}{\spro{w}}{\F}$ is taken as
an element of $\vecspaces{\WW{}}{\F}$, and $\anorm{}$ an element of $\norms{\WS{}}$.
Also, for each $i$ in $\findex{n}$, $\NVS{i}:=\opair{\VVS{i}}{\anorm{i}}$, and
$\NWS{}:=\opair{\WS{}}{\anorm{}}$.
If $\NWS{}$ is a Banach-space, then $\NVMLC{\suc{\NVS{1}}{\NVS{n}}}{\NWS{}}$ is also a Banach-space.
\begin{equation}
\anorm{}\in\Banachnorms{\WS{}}\then
\VMLCnorm{\mtuple{\NVS{1}}{\NVS{2}}}{\NWS{}}\in
\Banachnorms{\VMLC{\suc{\NVS{1}}{\NVS{n}}}{\NWS{}}}.
\end{equation}
\proof

\endthm
%%%%%%%%%%%%%%%%%%%%%%%%%%%%%%%%%%%%%%%%%%%%%%%%%%%%%%%%%%%%%%%%%%%%%%%%%%%%%%%%%%%%%%%%%%%%%%%%%%%%%%%%%%
\theorem\label{thmmultilinearmapsonfinitedimensionalspacesarecontinuous}
$\sindex{\VV{i}}{i}{\findex{n}}$
is taken to be an indexed-collection of non-empty sets, for some $n$ in $\N$.
For each $i$ in $\seta{\suc{1}{n}}$, $\VVS{i}=\tuple{\VV{i}}{\vsum{i}}{\spro{i}}{\F}$
is taken as an element of $\vecspaces{\VV{i}}{\F}$, and $\anorm{i}$
an element of $\norms{\VVS{i}}$.
$\WW{}$ is taken as a non-empty set.
$\WS{}=\tuple{\WW{}}{\vsum{w}}{\spro{w}}{\F}$ is taken as
an element of $\vecspaces{\WW{}}{\F}$, and $\anorm{}$ an element of $\norms{\WS{}}$.
Also, for each $i$ in $\findex{n}$, $\NVS{i}:=\opair{\VVS{i}}{\anorm{i}}$, and
$\NWS{}:=\opair{\WS{}}{\anorm{}}$. If for each $i$ in $\findex{n}$, $\VVS{i}$ is a
finite-dimensional vector-space, then every multilinear map from $\mtuple{\VVS{1}}{\VVS{n}}$
to $\WS{}$ is a continuous (and hence a multilinear-continuous) map from $\mtuple{\NVS{1}}{\NVS{n}}$
to $\NWS{}$.
\begin{align}
\big[\Foreach{i}{\findex{n}}\func{\Vdim{}}{\VVS{i}}\in\Zp\big]
\then
\MLin{\suc{\VVS{1}}{\VVS{n}}}{\WS{}}=\MLC{\suc{\NVS{1}}{\NVS{n}}}{\NWS{}}.
\end{align}
\proof
$\ml{}$ is taken as an element of $\MLin{\suc{\VVS{1}}{\VVS{n}}}{\WS{}}$.\\
For every $i$ in $\findex{n}$, $\vbase{i}$ is taken as an ordered-basis of the
vector-space $\VVS{i}$, and
\begin{equation}
\Foreach{i}{\findex{n}}
\vdim{i}:=\func{\Vdim{}}{\VVS{i}}.
\end{equation}
So,
\begin{equation}
\Foreach{i}{\findex{n}}
\Foreach{\vv{}}{\VVS{i}}
\vv{}=\ssum{k=1}{\vdim{i}}{\func{\[\com{\VVS{i}}{\vv{}}{\vbase{i}}\]}{k}\spro{i}\func{\vbase{i}}{k}},
\end{equation}
where, $\func{\[\com{\VVS{i}}{\vv{}}{\vbase{i}}\]}{k}$ denotes the $k$-th component of
the expansion of $\vv{}$ with respect to the ordered-basis $\vbase{i}$.
According to multilinearity of $\ml{}$,
\begin{align}
&\Foreach{\vv{}}{\(\Cproduct{\sindex{\VV{i}}{i}{\findex{n}}}\)}\cr
&\func{\ml{}}{\vv{}}=\sum_{k_1=1}^{\vdim{1}}\ldots\sum_{k_n=1}^{\vdim{n}}
\bigg(\func{\[\com{\VVS{1}}{\vv{1}}{\vbase{1}}\]}{k_1}\ldots
\func{\[\com{\VVS{n}}{\vv{n}}{\vbase{n}}\]}{k_n}\bigg)\spro{w}
\func{\ml{}}{\suc{\func{\vbase{1}}{k_1}}{\func{\vbase{n}}{k_n}}}.\cr
&{}
\end{align}
So, according to \refdef{defnorm},
\begin{align}
&\Foreach{\vv{}}{\(\Cproduct{\sindex{\VV{i}}{i}{\findex{n}}}\)}\cr
&\norm{\func{\ml{}}{\vv{}}}{}\leq\sum_{k_1=1}^{\vdim{1}}\ldots\sum_{k_n=1}^{\vdim{n}}
\bigg(\abs{\func{\[\com{\VVS{1}}{\vv{1}}{\vbase{1}}\]}{k_1}}\ldots
\abs{\func{\[\com{\VVS{n}}{\vv{n}}{\vbase{n}}\]}{k_n}}\bigg)
\norm{\func{\ml{}}{\suc{\func{\vbase{1}}{k_1}}{\func{\vbase{n}}{k_n}}}}{},\cr
&{}
\end{align}
and hence by defining,
\begin{equation}
C:=1+\max\defSet{\norm{\func{\ml{}}{\suc{\func{\vbase{1}}{k_1}}{\func{\vbase{n}}{k_n}}}}{}}
{k_1\in\findex{\vdim{1}},\ldots,k_n\in\findex{\vdim{n}}},
\end{equation}
it can be seen that,
\begin{align}
&\Foreach{\vv{}}{\(\Cproduct{\sindex{\VV{i}}{i}{\findex{n}}}\)}\cr
&\norm{\func{\ml{}}{\vv{}}}{}\leq C\sum_{k_1=1}^{\vdim{1}}\ldots\sum_{k_n=1}^{\vdim{n}}
\bigg(\abs{\func{\[\com{\VVS{1}}{\vv{1}}{\vbase{1}}\]}{k_1}}\ldots
\abs{\func{\[\com{\VVS{n}}{\vv{n}}{\vbase{n}}\]}{k_n}}\bigg),
\end{align}
and thus,
\begin{align}
&\Foreach{\vv{}}{\(\Cproduct{\sindex{\VV{i}}{i}{\findex{n}}}\)}\cr
&\norm{\func{\ml{}}{\vv{}}}{}\leq C\bigg(\ssum{k_1=1}{\vdim{1}}
{\abs{\func{\[\com{\VVS{1}}{\vv{1}}{\vbase{1}}\]}{k_1}}}\bigg)\ldots
\bigg(\ssum{k_n=1}{\vdim{n}}
{\abs{\func{\[\com{\VVS{n}}{\vv{n}}{\vbase{n}}\]}{k_n}}}\bigg).
\end{align}
Additionally, according to \refthm{thmabsnorm1},
\begin{equation}
%&\Foreach{\vv{}}{\(\Cproduct{\sindex{\VV{i}}{i}{\findex{n}}}\)}\cr
\Foreach{i}{\findex{n}}
\Existsis{M_i}{\Rp}\bigg[\Foreach{\uu{}}{\VV{i}}
\bigg(\ssum{k_i=1}{\vdim{i}}
{\abs{\func{\[\com{\VVS{i}}{\uu{}}{\vbase{i}}\]}{k_i}}}\bigg)\leq M_i\norm{\uu{}}{i}\bigg].
\end{equation}
Therefore,
\begin{align}
\Existsis{C\(M_1\ldots M_n\)}{\Rp}
\Foreach{\vv{}}{\(\Cproduct{\sindex{\VV{i}}{i}{\findex{n}}}\)}
\norm{\func{\ml{}}{\vv{}}}{}\leq \[C\(M_1\ldots M_n\)\]
\norm{\vv{1}}{1}\ldots\norm{\vv{n}}{n},
\end{align}
and hence according to \refthm{thmmultilinearcontinuousmapequivs},
\begin{equation}
\ml{}\in\MLC{\suc{\NVS{1}}{\NVS{n}}}{\NWS{}}.
\end{equation}
\endthm
\newpage
\Bibliography{}
\renewcommand{\addcontentsline}[3]{}

\let\addcontentsline\oldaddcontentsline
%%%%%%%%%%%%%%%%%%%%%%%%%%%%%%%%%%%%%%%%%%%%%%%%%%%%%%%%%%%%%%%%%%%%%%%%%%%%%%%%%%%%%%%%%%%%%%%%%%%%%%%

\end{document}